\documentclass {msuphddissertation}

\usepackage {amsmath, amssymb, amsthm}
\usepackage {mathrsfs}
\usepackage {graphicx}
\usepackage{hyperref}
\usepackage {tikz}
\usepackage{tikz-cd}
\usepackage[shortlabels]{enumitem}
\usetikzlibrary {snakes, calc}

\usepackage {libertine}
\usepackage [T1]{fontenc}
\usepackage[section]{placeins}
\usepackage{physics}
\usepackage{mathdots}

\theoremstyle {plain}
\newtheorem {thm} {Theorem}[chapter]

\newtheorem {lem} [thm]{Lemma}
\newtheorem {cor} [thm]{Corollary}
\newtheorem {prop} [thm]{Proposition}

\theoremstyle {definition}
\newtheorem {defn} [thm]{Definition}
\newtheorem {ex} [thm]{Example}

\newtheorem {rmk} [thm]{Remark}

\author {Andrew Claussen} 
\title {Expansion Posets for Polygon Cluster Algebras} 

\unit{Mathematics --- Doctor of Philosophy} 

\begin{document}

\maketitlepage

\begin {abstract} 
Define an \textit{expansion poset} to be the poset of monomials of a cluster variable attached to an arc in a polygon, where each monomial is represented by the corresponding combinatorial object from some fixed combinatorial cluster expansion formula. We introduce an involution on several of the interrelated combinatorial objects and constructions associated to type $A$ surface cluster algebras, including certain classes of arcs, triangulations, and distributive lattices. We use these involutions to formulate a dual version of skein relations for arcs, and dual versions of three existing expansion posets. In particular, this leads to two new cluster expansion formulas, and recovers the lattice path expansion of Propp et al. We provide an explicit, structure-preserving poset isomorphism between an expansion poset and its dual version from the dual arc. We also show that an expansion poset and its dual version constructed from the same arc are dual in the sense of distributive lattices.

We show that any expansion poset is isomorphic to a closed interval in one of the lattices $L(m,n)$ of Young diagrams contained in an $m \times n$ grid, and that any $L(m,n)$ has a covering by such intervals. In particular, this implies that any expansion poset is isomorphic to an interval in Young's lattice.

We give two formulas for the rank function of any lattice path expansion poset, and prove that this rank function is unimodal whenever the underlying snake graph is built from at most four maximal straight segments. This gives a partial solution to a recent conjecture by Ovsienko and Morier-Genoud. We also characterize which expansion posets have symmetric rank generating functions, based on the shape of the underlying snake graph.

We show that the support of any type $A$ cluster variable is the orbit of a groupoid. This implies that any such cluster variable can be reconstructed from any one of its monomials. 

Finally, in work joint with Nicholas Ovenhouse, we partially generalize $T$-paths to configurations of affine flags, and prove that a $T$-path expansion analogous to the type $A$ case holds when the initial seed is from a fan triangulation. We finish by describing the structure of two types of expansion posets in this context.
\end {abstract}

\begin {acknowledgment}
First and foremost, I would like to thank my advisor Misha Shapiro. Certainly, this thesis would not have been possible without his tireless guidance, endless patience, and unwavering personal and academic support. 

I would like to thank the other members of my committee Jonathan Hall, Rajesh Kulkarni, and Bruce Sagan, both for their inspiring courses and lectures, and for their continual support and encouragement throughout the years. 

A special thanks goes to my friend and collaborator Nicholas Ovenhouse. In particular, I am grateful for both the innumerable comments and corrections that he offered during the writing of this thesis, and also for his willingness to have our joint work included here.

I would like to thank my fellow graduate students for their friendship and support, especially Nick Ovenhouse, Blake Icabone, Duff Baker-Jarvis, and Sami Merhi. 

Likewise, thanks to everyone I met through the advanced track program at MSU, especially Daniel Diroff, Bar Roytman, Christopher Klerkx, Rebecca Sodervick, Jonathan Jonker, and Katrina Suchoski. 

I am grateful to Tsvetanka Sendova, Russel Schwab, and Jeanne Wald for their commitment to my professional development, and for all the opportunities they have presented me with during my time at MSU.

Indeed, thanks to all the faculty and staff at MSU with whom I have had the pleasure of working with. 

I would like to thank Mark Naber of Monroe County Community College for inspiring me to further my mathematical education. 

Finally, I am thankful to all my wonderful family and friends for the love and support they have given me throughout this endeavor. I'd like to give special thanks to my dad (Rusty), Marshia, Aaron and Autumn, Chris, Jake, Graham, and the Minney family. 
\end {acknowledgment}

\begin{dedication}
To my parents, Kay Marie Claussen and Irvin Boyd Claussen.
\end{dedication}

\TOC 

\LOT 

\LOF 

\newpage

\pagenumbering{arabic}

\begin {doublespace}


\chapter{Introduction}

\textit{Cluster algebras} are a class of inherently combinatorial commutative algebras that were defined by Fomin and Zelevinsky in \cite{fomin2002cluster}. The definition of cluster algebras was motivated by observations in representation theory. Since then, cluster algebra structures have been recognized and studied in various other fields of mathematics, such as decorated Teichm\"uller theory and Poisson geometry (see \cite{gekhtman2005cluster} and \cite{fomin2006cluster}, \cite{fomin2018cluster}), higher Teichm\"uller theory \cite{fock2006moduli}, rings of invariants (see \cite{fomin2012tensor} and \cite{fomin2014webs}),  elementary number theory \cite{ccanakcci2018cluster} including diophantine equations \cite{rabideau2018continued}, and knot theory \cite{lee2019cluster}, just to name a few.

Each cluster algebra has a distinguished set of generators, called the \textit{cluster variables}. These generators are grouped into overlapping subsets, called \textit{clusters}, all having the same cardinality, called the \textit{rank} of the cluster algebra.

A \textit{seed} of a cluster algebra is a triple consisting of a cluster, a \textit{coefficient tuple}, and a skew-symmetrizable matrix called an \textit{exchange matrix}. In any rank $n$ cluster algebra, each seed can be \textit{mutated} in direction $i$ for any $1 \leq i \leq n$ to produce $n$ more seeds. By construction, seed mutation is an involution. Clusters in adjacent seeds differing by a mutation in direction $i$ are equal, except that the $i^{th}$ variables in each differ from one another by what is called an \textit{exchange relation}. This means that their product is a certain binomial sum whose form is governed by the exchange matrices. Any cluster algebra can be computed by fixing an \textit{initial seed} and iterating seed mutation to produce all the cluster variables.

A \textit{Laurent polynomial} is a polynomial with negative exponents allowed, i.e., any Laurent polynomial can be written as $\sum_{I \in \mathbb{Z}^{n}} a_I X^{I}$. One of the first fundamental results of cluster algebra theory is that any cluster variable can be written as a Laurent polynomial with respect to any cluster. This is called the \textit{Laurent phenomenon}.

We will work within the subclass of cluster algebras in which exchange matrices and matrix mutation may be replaced by \textit{quivers} and \textit{quiver mutation}, respectively. These cluster algebras are called \textit{cluster algebras from quivers}, or \textit{skew-symmetric cluster algebras of geometric type}

A large subclass of cluster algebras are the \textit{cluster algebras from surfaces} \cite{felikson2012cluster}. In such a cluster algebra, the cluster variables are in bijection with certain curves in the surface called \textit{tagged arcs}, seeds are in bijection with \textit{tagged triangulations} of the surface, and mutation corresponds to a \textit{flip} of a tagged arc in a tagged triangulation.

The combinatorics of the subclass of cluster algebras from a disc with $n+3$ marked points on the boundary (i.e., a polygon) are the main focus of this paper. These cluster algebras are examples of \textit{cluster algebras of finite type $A_n$}. As the notion of tagged triangulations and tagged arcs is unnecessary in this restricted level of generality, we will henceforth only speak of \textit{(ordinary) arcs}. Any seed in a cluster algebra of finite type $A_n$ can be modeled by a triangulated $(n+3)$-gon, made up of $n$ triangulating arcs (which we also call \textit{internal diagonals}, or just \textit{diagonals}), and $n+3$ segments that make up the boundary of the polygon, called \textit{boundary segments}. Mutation corresponds to a flip of one of the $n$ internal diagonals in the triangulation. Figure \ref{fig:flip_intro} below shows one of the two possible flips inside a triangulated polygon. 

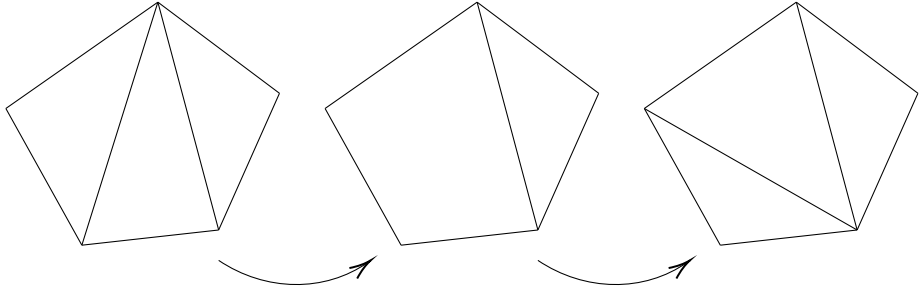
\begin {figure}[h!]
    \centering
    \caption{A flip inside a triangulated pentagon}
    \label{fig:flip_intro}
    \begin{tikzpicture}[x=0.75pt,y=0.75pt,yscale=-1,xscale=1]

\draw    (20,93.67) -- (96.67,40) ;
\draw    (158,86) -- (96.67,40) ;
\draw    (58.33,162.67) -- (20,93.67) ;
\draw    (127.33,155) -- (158,86) ;
\draw    (58.33,162.67) -- (127.33,155) ;

\draw    (181,93.67) -- (257.67,40) ;
\draw    (319,86) -- (257.67,40) ;
\draw    (219.33,162.67) -- (181,93.67) ;
\draw    (288.33,155) -- (319,86) ;
\draw    (219.33,162.67) -- (288.33,155) ;

\draw    (342,93.67) -- (418.67,40) ;
\draw    (480,86) -- (418.67,40) ;
\draw    (380.33,162.67) -- (342,93.67) ;
\draw    (449.33,155) -- (480,86) ;
\draw    (380.33,162.67) -- (449.33,155) ;

\draw [color={rgb, 255:red, 0; green, 0; blue, 0 }  ,draw opacity=1 ]   (96.67,40) -- (58.33,162.67) ;
\draw    (96.67,40) -- (127.33,155) ;
\draw    (257.67,40) -- (288.33,155) ;
\draw    (418.67,40) -- (449.33,155) ;
\draw [color={rgb, 255:red, 0; green, 0; blue, 0 }  ,draw opacity=1 ]   (342,93.67) -- (449.33,155) ;
\draw    (127.33,170.33) .. controls (151,186.26) and (180.93,186.58) .. (202.68,171.29) ;
\draw [shift={(204,170.33)}, rotate = 503.36] [color={rgb, 255:red, 0; green, 0; blue, 0 }  ][line width=0.75]    (10.93,-3.29) .. controls (6.95,-1.4) and (3.31,-0.3) .. (0,0) .. controls (3.31,0.3) and (6.95,1.4) .. (10.93,3.29)   ;
\draw    (288.33,170.33) .. controls (312,186.26) and (341.93,186.58) .. (363.68,171.29) ;
\draw [shift={(365,170.33)}, rotate = 503.36] [color={rgb, 255:red, 0; green, 0; blue, 0 }  ][line width=0.75]    (10.93,-3.29) .. controls (6.95,-1.4) and (3.31,-0.3) .. (0,0) .. controls (3.31,0.3) and (6.95,1.4) .. (10.93,3.29)   ;

\end{tikzpicture}
  
\end {figure}

Numerous \textit{combinatorial cluster expansion formulas} for surface type and type $A_n$ cluster algebras have been developed in recent years. Each such expansion gives an explicit formula to compute any cluster variable by writing it as a weighted sum over a certain class of combinatorial objects.  We recall three such expansion formulas from the literature - \textit{perfect matchings of snake graphs}, \textit{perfect matchings of angles}, and \textit{$T$-paths}. (see \cite{musiker2011positivity} and \cite{musiker2013bases}, \cite{yurikusa2018combinatorial} and \cite{yurikusa2019cluster}, and \cite{schiffler2008cluster}, \cite{schiffler2009cluster}, and \cite{schiffler2010cluster} respectively). 

The monomials of any cluster variable can be naturally arranged into a poset \cite{fomin2007cluster}. We can equip each such poset with the additional node structure it inherits from some combinatorial expansion formula. We call any such poset an \textit{expansion poset} (in fact, each of these posets is isomorphic to a distributive lattice \cite{musiker2013bases}).

In this paper, inspired by the snake graph involution introduced in \cite{propp2005combinatorics} and the involution ``Jimm'' in \cite{uludaug2015jimm}, we define a new notion of \textit{duality}, a certain involution on several of the combinatorial objects found in type A cluster theory, including for instance triangulations of an $n$-gon, and a certain class of distributive lattices. Each such object is parameterized by a binary word $w$, so that duality between objects is controlled by duality $w \longleftrightarrow w^{*}$ on the underlying words.  These constructions yield a dual version of skein relations for arcs. Furthermore, we use this duality to produce an equivalent yet combinatorially distinct version (built from the dual arc in the dual triangulation) of each of the three expansion posets mentioned above.

Namely, we observe that there is an explicit poset isomorphism (respecting additional node structure) between the perfect matching expansion poset $\mathbb{P}_{w}$, and the \textit{dual lattice path expansion poset} $\mathbb{L}_{w^{*}}$ associated to the dual word. This duality, minus poset structures, was given in \cite{propp2005combinatorics} in the context of frieze patterns. We define two more expansion posets, new to the best of our knowledge, which we call \textit{lattice paths of angles} $\mathbb{B}_w$ and \textit{$S$-paths} $\mathbb{S}_w$, respectively. We observe that there is an explicit structure-preserving poset isomorphism between the expansion posets $\mathbb{A}_w$ and $\mathbb{B}_{w^{*}},$ and likewise between the $T$-path expansion poset $\mathbb{T}_w$ and $\mathbb{S}_{w^{*}}.$

The three isomorphisms just indicated are all written in terms of either snake graph duality, or traingulation duality. Here, we call any such isomorphism an \textit{expansion duality}. The horizontal maps in the figure below represent expansion dualities. The maps comprising the two triangles are written either in terms of the \textit{folding/unfolding maps} from \cite{musiker2010cluster}, or the \textit{angle projection} maps defined in \cite{yurikusa2019cluster} (or some combination of the two).

\begin {figure}[h!]
    \centering
    \caption{Expansion duality}
    \label{fig:exp_duality}
    \begin{tikzcd}
                                                & \mathbb{P}_w \arrow[ldd] \arrow[r] \arrow[d]        & \mathbb{L}_{w}^{*} \arrow[ldd] \arrow[l] \arrow[d] \\
                                                & \mathbb{A}_w \arrow[u] \arrow[r] \arrow[ld]         & \mathbb{B}_{w}^{*} \arrow[u] \arrow[l] \arrow[ld]  \\
\mathbb{T}_{w} \arrow[ru] \arrow[r] \arrow[ruu] & \mathbb{S}_{w}^{*} \arrow[ru] \arrow[l] \arrow[ruu] &                                                   
\end{tikzcd}
\end {figure}
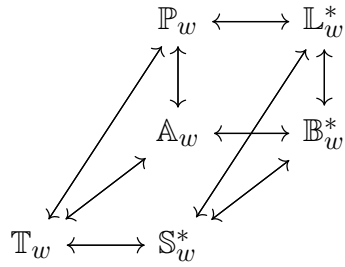

Furthermore, we show that the two isomorphism classes of expansion posets attached to any arc are dual in the sense of distributive lattices mentioned above.

A \textit{partition} of a positive number $m$ is a weakly decreasing sequence $(\lambda_1 , \lambda_2 , \dots , \lambda_{l} )$ such that $m = \lambda_1 + \lambda_2 + \dots + \lambda_l.$  A \textit{Young diagram} is a visualization of a partition of a positive integer by rows of boxes (see \cite{sagan2013symmetric}). \textit{Young's lattice} is the infinite lattice whose nodes are Young diagrams, which are ordered by inclusion (see \cite{sagan2013symmetric} and \cite{stanley1988differential}). Young's lattice possesses a collection of finite sublattices, typically called $L(m,n)$, whose nodes are all those Young diagrams that fit within an $m \times n$ rectangular grid, and whose rank function is the $q$-binomial coefficient ${n+m \brack m}_q$. We give a groupoid structure on the set of all snake graphs which reconstructs the posets $L(m,n)$. That is, each orbit of this groupoid can be given the structure of a poset, which is isomorphic to one of the lattices $L(m,n)$, and furthermore (the isomorphism class of) every $L(m,n)$ can be obtained in this way. We show that each expansion poset considered above is isomorphically embedded as a closed interval in some $L(m,n)$, and moreover that each $L(m,n)$ has a covering by such embeddings.

We provide two closed formulas for the rank function $\mathbb{L}_w (q)$ of any lattice path expansion poset $\mathbb{L}_w$ of lattice paths on the snake graph $G_w$. The first is written as sums of products of \textit{hook snake graphs}, and its terms are parameterized by the nodes in a Boolean lattice. The second formula we introduce is a refinement of the first, and is written in terms of $q$-numbers corresponding to lengths of the maximal straight segments that $G_w$ is built from. The latter formula is obtained by considering the snake graph $G_w$ itself as a lattice path on another snake graph. In fact, this construction makes $G_w$ into the central node in a \textit{Fibonacci cube}, an interconnection topology defined in \cite{hsu1993fibonacci}.

Recently, a conjecture was made in \cite{morier2018q} that is equivalent to asking if the coefficients of the rank generating function of any expansion poset are \textit{unimodal}, meaning that they weakly increase, and then weakly decrease. By again studying lattice path expansions, we show this conjecture to be true for snake graphs built from at most four maximal straight segments. We also characterize, based on the shape of the underlying snake graph, which expansion posets have palindromic, or \textit{symmetric}, coefficients.

Next, we interpret the support of any cluster variable $x_w$ in a type $A$ surface cluster algebra as an orbit of a certain groupoid. It follows that any two cluster variables, written with respect to the same initial seed, have disjoint supports. Thus, any cluster variable $x_w$ is completely determined by any one of its Laurent monomials.

Lastly, in work joint with Nicholas Ovenhouse, we partially generalize the $T$-path expansion mentioned above to configurations of affine flags (see \cite{fock2006moduli}). We describe in two special cases the poset structure on the Laurent monomials from an expansion in this context.

\chapter {Cluster Algebras}

We begin this chapter by defining the subclass of cluster algebras called \textit{cluster algebras from quivers}, or \textit{skew-symmetric cluster algebras of geometric type}. Next, we introduce \textit{cluster algebras from surfaces}. Finally, we describe \textit{surface cluster algebras of type $A$}, which is the class of cluster algebras we study in the rest of the paper.

\section {Cluster Algebras from Quivers}

\begin {defn}
   
   A \textit{quiver} is a $4$-tuple $Q = (Q_0 , Q_1 , s , t)$, where $Q_0$ is a set of \textit{nodes} or \textit{vertices}, $Q_1$ is a multiset of \textit{arrows} whose elements are from $Q_0 \times Q_0$, and the \textit{source} and \textit{target} functions $s,t: Q_1 \longrightarrow Q_0$ are the first and second projection, respectively. A quiver $Q$ contains a \textit{loop} if there exists a node $v \in Q_0$ and an arrow $e \in Q_1$ such that $v = s(e) = t(e)$. We say $Q$ contains a \textit{$2$-cycle} if there exists two distinct nodes $v_1 , v_2 \in Q_0$ and two arrows $e_1 , e_2 \in Q_1$ such that $v_1 = s(e_1) = t(e_2)$ and $v_2 = s(e_2) = t(e_1).$ A quiver $Q$ contains a \textit{$3$-cycle} if there exist three distinct nodes $v_1 , v_2 , v_3 \in Q_0 $ and three arrows $e_1 , e_2 , e_3 \in Q_1$ such that $v_1 = s(e_1) = t(e_3)$, $v_2 = s(e_2) = t(e_1)$, and $v_3 = s(e_3) = t(e_2)$. We say $Q$ is \textit{finite} if both $Q_0$ and $Q_1$ are finite sets. In this case, we label the vertices of $Q$ by $1,2,... , \# \{ \text{nodes in } Q \} .$

\end {defn}

 From now on, we only consider finite quivers (unless stated otherwise) that do not contain any loops or $2$-cycles.

\begin {defn}
   Let $Q = (Q_0 , Q_1)$ be a quiver with vertices $Q_0 = \{ 1,2,...,n \}$. We define \textit{quiver mutation at vertex $k$} to be the the map on quivers $\mu_k$ whose image $\mu_k (Q)$ is defined by the following procedure:

\begin{enumerate}

\item For each pair of arrows $i \rightarrow k \rightarrow j$, create a new directed edge $i \rightarrow j$ (note that prohibiting $2$-cycles implies that $i \neq j$).

\item Reverse the orientation of each arrow incident to $k$

\item Delete any and all $2$-cycles created in step 1. 

\end{enumerate}
\end {defn}

\begin {ex}
   In Figure \ref{fig:quiver_mutation}, we illustrate the three-step process in Definition 2.2 by mutating the left-most quiver at vertex $1$ to produce the right-most quiver.
   
   
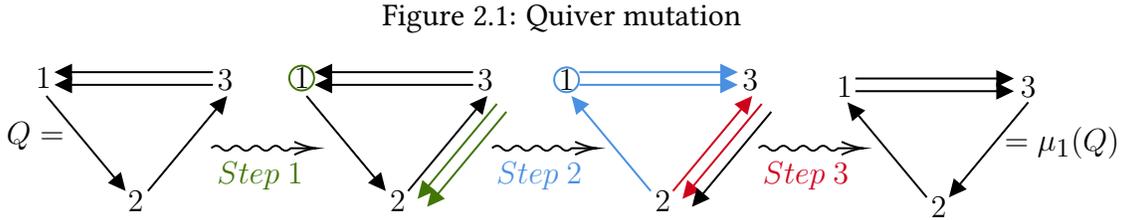
\begin {figure}[h!]
    \centering
    \caption{Quiver mutation}
    \label{fig:quiver_mutation}

\tikzset{every picture/.style={line width=0.75pt}} 

\begin{tikzpicture}[x=0.75pt,y=0.75pt,yscale=-1,xscale=1]

\draw    (220.8,192.8) -- (143.8,192.8) ;
\draw [shift={(140.8,192.8)}, rotate = 360] [fill={rgb, 255:red, 0; green, 0; blue, 0 }  ][line width=0.08]  [draw opacity=0] (8.93,-4.29) -- (0,0) -- (8.93,4.29) -- cycle    ;
\draw    (220.8,186.4) -- (143.8,186.4) ;
\draw [shift={(140.8,186.4)}, rotate = 360] [fill={rgb, 255:red, 0; green, 0; blue, 0 }  ][line width=0.08]  [draw opacity=0] (8.93,-4.29) -- (0,0) -- (8.93,4.29) -- cycle    ;
\draw    (136.8,198.4) -- (174.88,244.1) ;
\draw [shift={(176.8,246.4)}, rotate = 230.19] [fill={rgb, 255:red, 0; green, 0; blue, 0 }  ][line width=0.08]  [draw opacity=0] (8.93,-4.29) -- (0,0) -- (8.93,4.29) -- cycle    ;
\draw    (226.08,200.7) -- (188,246.4) ;
\draw [shift={(228,198.4)}, rotate = 129.81] [fill={rgb, 255:red, 0; green, 0; blue, 0 }  ][line width=0.08]  [draw opacity=0] (8.93,-4.29) -- (0,0) -- (8.93,4.29) -- cycle    ;
\draw    (352,192.8) -- (275,192.8) ;
\draw [shift={(272,192.8)}, rotate = 360] [fill={rgb, 255:red, 0; green, 0; blue, 0 }  ][line width=0.08]  [draw opacity=0] (8.93,-4.29) -- (0,0) -- (8.93,4.29) -- cycle    ;
\draw    (352,186.4) -- (275,186.4) ;
\draw [shift={(272,186.4)}, rotate = 360] [fill={rgb, 255:red, 0; green, 0; blue, 0 }  ][line width=0.08]  [draw opacity=0] (8.93,-4.29) -- (0,0) -- (8.93,4.29) -- cycle    ;
\draw    (268,198.4) -- (306.08,244.1) ;
\draw [shift={(308,246.4)}, rotate = 230.19] [fill={rgb, 255:red, 0; green, 0; blue, 0 }  ][line width=0.08]  [draw opacity=0] (8.93,-4.29) -- (0,0) -- (8.93,4.29) -- cycle    ;
\draw    (357.28,200.7) -- (319.2,246.4) ;
\draw [shift={(359.2,198.4)}, rotate = 129.81] [fill={rgb, 255:red, 0; green, 0; blue, 0 }  ][line width=0.08]  [draw opacity=0] (8.93,-4.29) -- (0,0) -- (8.93,4.29) -- cycle    ;
\draw [color={rgb, 255:red, 65; green, 117; blue, 5 }  ,draw opacity=1 ]   (364,202.4) -- (325.92,248.1) ;
\draw [shift={(324,250.4)}, rotate = 309.81] [fill={rgb, 255:red, 65; green, 117; blue, 5 }  ,fill opacity=1 ][line width=0.08]  [draw opacity=0] (8.93,-4.29) -- (0,0) -- (8.93,4.29) -- cycle    ;
\draw [color={rgb, 255:red, 65; green, 117; blue, 5 }  ,draw opacity=1 ]   (368.8,206.4) -- (330.72,252.1) ;
\draw [shift={(328.8,254.4)}, rotate = 309.81] [fill={rgb, 255:red, 65; green, 117; blue, 5 }  ,fill opacity=1 ][line width=0.08]  [draw opacity=0] (8.93,-4.29) -- (0,0) -- (8.93,4.29) -- cycle    ;
\draw [color={rgb, 255:red, 74; green, 144; blue, 226 }  ,draw opacity=1 ]   (482.6,192.8) -- (405.6,192.8) ;
\draw [shift={(485.6,192.8)}, rotate = 180] [fill={rgb, 255:red, 74; green, 144; blue, 226 }  ,fill opacity=1 ][line width=0.08]  [draw opacity=0] (8.93,-4.29) -- (0,0) -- (8.93,4.29) -- cycle    ;
\draw [color={rgb, 255:red, 74; green, 144; blue, 226 }  ,draw opacity=1 ]   (482.6,186.4) -- (405.6,186.4) ;
\draw [shift={(485.6,186.4)}, rotate = 180] [fill={rgb, 255:red, 74; green, 144; blue, 226 }  ,fill opacity=1 ][line width=0.08]  [draw opacity=0] (8.93,-4.29) -- (0,0) -- (8.93,4.29) -- cycle    ;
\draw [color={rgb, 255:red, 74; green, 144; blue, 226 }  ,draw opacity=1 ]   (403.52,200.7) -- (441.6,246.4) ;
\draw [shift={(401.6,198.4)}, rotate = 50.19] [fill={rgb, 255:red, 74; green, 144; blue, 226 }  ,fill opacity=1 ][line width=0.08]  [draw opacity=0] (8.93,-4.29) -- (0,0) -- (8.93,4.29) -- cycle    ;
\draw [color={rgb, 255:red, 208; green, 2; blue, 27 }  ,draw opacity=1 ]   (490.88,200.7) -- (452.8,246.4) ;
\draw [shift={(492.8,198.4)}, rotate = 129.81] [fill={rgb, 255:red, 208; green, 2; blue, 27 }  ,fill opacity=1 ][line width=0.08]  [draw opacity=0] (8.93,-4.29) -- (0,0) -- (8.93,4.29) -- cycle    ;
\draw [color={rgb, 255:red, 208; green, 2; blue, 27 }  ,draw opacity=1 ]   (497.6,202.4) -- (459.52,248.1) ;
\draw [shift={(457.6,250.4)}, rotate = 309.81] [fill={rgb, 255:red, 208; green, 2; blue, 27 }  ,fill opacity=1 ][line width=0.08]  [draw opacity=0] (8.93,-4.29) -- (0,0) -- (8.93,4.29) -- cycle    ;
\draw [color={rgb, 255:red, 0; green, 0; blue, 0 }  ,draw opacity=1 ]   (502.4,206.4) -- (464.32,252.1) ;
\draw [shift={(462.4,254.4)}, rotate = 309.81] [fill={rgb, 255:red, 0; green, 0; blue, 0 }  ,fill opacity=1 ][line width=0.08]  [draw opacity=0] (8.93,-4.29) -- (0,0) -- (8.93,4.29) -- cycle    ;
\draw    (361.6,224.8) .. controls (363.27,223.13) and (364.93,223.13) .. (366.6,224.8) .. controls (368.27,226.47) and (369.93,226.47) .. (371.6,224.8) .. controls (373.27,223.13) and (374.93,223.13) .. (376.6,224.8) .. controls (378.27,226.47) and (379.93,226.47) .. (381.6,224.8) .. controls (383.27,223.13) and (384.93,223.13) .. (386.6,224.8) .. controls (388.27,226.47) and (389.93,226.47) .. (391.6,224.8) .. controls (393.27,223.13) and (394.93,223.13) .. (396.6,224.8) .. controls (398.27,226.47) and (399.93,226.47) .. (401.6,224.8) -- (404.4,224.8) -- (412.4,224.8) ;
\draw [shift={(414.4,224.8)}, rotate = 180] [color={rgb, 255:red, 0; green, 0; blue, 0 }  ][line width=0.75]    (10.93,-3.29) .. controls (6.95,-1.4) and (3.31,-0.3) .. (0,0) .. controls (3.31,0.3) and (6.95,1.4) .. (10.93,3.29)   ;
\draw    (220,224.8) .. controls (221.67,223.13) and (223.33,223.13) .. (225,224.8) .. controls (226.67,226.47) and (228.33,226.47) .. (230,224.8) .. controls (231.67,223.13) and (233.33,223.13) .. (235,224.8) .. controls (236.67,226.47) and (238.33,226.47) .. (240,224.8) .. controls (241.67,223.13) and (243.33,223.13) .. (245,224.8) .. controls (246.67,226.47) and (248.33,226.47) .. (250,224.8) .. controls (251.67,223.13) and (253.33,223.13) .. (255,224.8) .. controls (256.67,226.47) and (258.33,226.47) .. (260,224.8) -- (262.8,224.8) -- (270.8,224.8) ;
\draw [shift={(272.8,224.8)}, rotate = 180] [color={rgb, 255:red, 0; green, 0; blue, 0 }  ][line width=0.75]    (10.93,-3.29) .. controls (6.95,-1.4) and (3.31,-0.3) .. (0,0) .. controls (3.31,0.3) and (6.95,1.4) .. (10.93,3.29)   ;
\draw    (496,224.8) .. controls (497.67,223.13) and (499.33,223.13) .. (501,224.8) .. controls (502.67,226.47) and (504.33,226.47) .. (506,224.8) .. controls (507.67,223.13) and (509.33,223.13) .. (511,224.8) .. controls (512.67,226.47) and (514.33,226.47) .. (516,224.8) .. controls (517.67,223.13) and (519.33,223.13) .. (521,224.8) .. controls (522.67,226.47) and (524.33,226.47) .. (526,224.8) .. controls (527.67,223.13) and (529.33,223.13) .. (531,224.8) .. controls (532.67,226.47) and (534.33,226.47) .. (536,224.8) -- (538.8,224.8) -- (546.8,224.8) ;
\draw [shift={(548.8,224.8)}, rotate = 180] [color={rgb, 255:red, 0; green, 0; blue, 0 }  ][line width=0.75]    (10.93,-3.29) .. controls (6.95,-1.4) and (3.31,-0.3) .. (0,0) .. controls (3.31,0.3) and (6.95,1.4) .. (10.93,3.29)   ;
\draw    (621.8,196) -- (544.8,196) ;
\draw [shift={(624.8,196)}, rotate = 180] [fill={rgb, 255:red, 0; green, 0; blue, 0 }  ][line width=0.08]  [draw opacity=0] (8.93,-4.29) -- (0,0) -- (8.93,4.29) -- cycle    ;
\draw    (621.8,189.6) -- (544.8,189.6) ;
\draw [shift={(624.8,189.6)}, rotate = 180] [fill={rgb, 255:red, 0; green, 0; blue, 0 }  ][line width=0.08]  [draw opacity=0] (8.93,-4.29) -- (0,0) -- (8.93,4.29) -- cycle    ;
\draw    (542.72,203.9) -- (580.8,249.6) ;
\draw [shift={(540.8,201.6)}, rotate = 50.19] [fill={rgb, 255:red, 0; green, 0; blue, 0 }  ][line width=0.08]  [draw opacity=0] (8.93,-4.29) -- (0,0) -- (8.93,4.29) -- cycle    ;
\draw    (632,201.6) -- (593.92,247.3) ;
\draw [shift={(592,249.6)}, rotate = 309.81] [fill={rgb, 255:red, 0; green, 0; blue, 0 }  ][line width=0.08]  [draw opacity=0] (8.93,-4.29) -- (0,0) -- (8.93,4.29) -- cycle    ;
\draw  [color={rgb, 255:red, 65; green, 117; blue, 5 }  ,draw opacity=1 ] (259.3,189.4) .. controls (259.3,185.95) and (262.1,183.15) .. (265.55,183.15) .. controls (269,183.15) and (271.8,185.95) .. (271.8,189.4) .. controls (271.8,192.85) and (269,195.65) .. (265.55,195.65) .. controls (262.1,195.65) and (259.3,192.85) .. (259.3,189.4) -- cycle ;
\draw  [color={rgb, 255:red, 74; green, 144; blue, 226 }  ,draw opacity=1 ] (392.9,190.25) .. controls (392.9,186.8) and (395.7,184) .. (399.15,184) .. controls (402.6,184) and (405.4,186.8) .. (405.4,190.25) .. controls (405.4,193.7) and (402.6,196.5) .. (399.15,196.5) .. controls (395.7,196.5) and (392.9,193.7) .. (392.9,190.25) -- cycle ;

\draw (130,182.7) node [anchor=north west][inner sep=0.75pt]    {$1$};
\draw (221.6,183.7) node [anchor=north west][inner sep=0.75pt]    {$3$};
\draw (176.4,244.7) node [anchor=north west][inner sep=0.75pt]    {$2$};
\draw (260.2,182.7) node [anchor=north west][inner sep=0.75pt]    {$\textcolor[rgb]{0,0,0}{1}$};
\draw (352.8,183.7) node [anchor=north west][inner sep=0.75pt]    {$3$};
\draw (308.6,244.7) node [anchor=north west][inner sep=0.75pt]    {$2$};
\draw (393.8,183.7) node [anchor=north west][inner sep=0.75pt]    {$1$};
\draw (486.4,183.7) node [anchor=north west][inner sep=0.75pt]    {$3$};
\draw (442.2,244.7) node [anchor=north west][inner sep=0.75pt]    {$2$};
\draw (534,186.9) node [anchor=north west][inner sep=0.75pt]    {$1$};
\draw (626.6,186.9) node [anchor=north west][inner sep=0.75pt]    {$3$};
\draw (581.4,247.9) node [anchor=north west][inner sep=0.75pt]    {$2$};
\draw (221.5,230.1) node [anchor=north west][inner sep=0.75pt]    {$\textcolor[rgb]{0.25,0.46,0.02}{Step\ 1}$};
\draw (362.3,230.1) node [anchor=north west][inner sep=0.75pt]    {$\textcolor[rgb]{0.29,0.56,0.89}{Step\ 2}$};
\draw (496.7,230.1) node [anchor=north west][inner sep=0.75pt]    {$\textcolor[rgb]{0.82,0.01,0.11}{Step\ 3}$};
\draw (618.6,212.6) node [anchor=north west][inner sep=0.75pt]    {$=\mu _{1}( Q)$};
\draw (115.2,209.9) node [anchor=north west][inner sep=0.75pt]    {$Q=$};

\end{tikzpicture}

\end {figure}

\end {ex}

One can check that quiver mutation at vertex $k$ is an involution, so that mutation equivalence is indeed an equivalence relation.

\begin {defn}
    Fix $n \leq m$ and let $Q = (Q_0 , Q_1)$ be a quiver with vertices $ Q_0 = \{ 1,2,...,n , n+1 , ... , m \}$. Partition $Q_0$ into the two sets $Q_{0}^{mutable} = \{  1,2,...,n \}$ and $Q_{0}^{frozen} = \{ n+1 , n+2 , ... , m \},$ which we call the \textit{mutable vertices} and the \textit{frozen vertices}, respectively. Fix an \textit{ambient field} $F \cong \mathbb{Q} (f_1 , f_2 , ... , f_n , f_{n+1} ,  ... , f_{m}).$ By associating to each vertex $i$ in $Q_0$ the indeterminate $f_i$, we obtain a \textit{seed} in $F$, denoted $\big( (f_1 , f_2 , ... , f_m ), Q \big).$ The variables $f_1 , f_2 , ... , f_n$ are called \textit{cluster variables} (or \textit{mutable variables}), and $f_{n+1} , ... , f_{m}$ are called \textit{frozen variables}. The $n$-tuple $(f_1 , f_2 , ... , f_n)$ is called the \textit{cluster} of the seed $\big( (f_1 , f_2 , ... , f_m ), Q \big)$, and the $m$-tuple $(f_1 , f_2 , ... , f_{m})$ is the \textit{extended cluster} of the seed $\big( (f_1 , f_2 , ... , f_m ), Q \big)$. 
\end {defn}

\begin {defn}

    Consider the seed $\big( (f_1 , f_2 , ... , f_m ), Q \big).$ For $1 \leq k \leq n$, we define \textit{seed mutation at variable $k$} to be the map on seeds $\mu_k$ whose image $\mu_{k} \big( (f_1 , f_2 , ... , f_m ), Q \big) = \big( (\widetilde{f_{1}} , \widetilde{f_{2}} , ... , \widetilde{f_{m}} ) , \mu_k (Q) \big)$ is defined as follows:

    \begin{itemize}
    \item If $j \neq k,$ then $\widetilde{f}_j = f_j.$
    
    \item If $j = k$, then $f_k$ and $\widetilde{f_k}$ are related by the following \textit{exchange relation}: 
    
    $$f_k \widetilde{f_k} =  \prod_{s \rightarrow k } f_s + \prod_{k \rightarrow t } f_t .$$
    \end{itemize}

\end {defn}

\begin {defn} 

    Consider the \textit{initial seed} $\big( (x_1 , x_2 , ... , x_m ) , Q \big)$ with \textit{initial cluster} $(x_1 , x_2 , ... , x_n )$ and \textit{initial cluster variables} $x_1 , x_2 , ... , x_n$. Let $\mathcal{S}$ be the set of all seeds mutation equivalent to $\big( (x_1 , x_2 , ... , x_m ) , Q \big),$ and let $\mathcal{X}$ be the union of all cluster variables in all seeds in $\mathcal{S}.$ Let $R = \mathbb{Z} [ x_{n+1} , ... , x_{m}].$ The \textit{cluster algebra $\mathcal{A} = \mathcal{A} \big( (x_1 , x_2 , ... , x_m ) , Q \big)$ from the quiver $Q$} is the $R$-algebra generated by $\mathcal{X}.$ The \textit{rank} of the cluster algebra $\mathcal{A}$ is the cardinality $n$ of any of its clusters. 

\end {defn}

Below we state the Laurent Phenomenon in the restricted generality of cluster algebras from quivers.

\begin {thm} (Theorem 3.1 in \cite{fomin2002cluster})
Let $\mathcal{A}$ be the cluster algebra from the quiver $Q$ with arbitrary initial seed $\big( (x_1 , x_2 , ... , x_m ), Q \big)$. Then any cluster variable in $\mathcal{A}$ can be expressed as a Laurent polynomial in the variables $x_1 , x_2 , ... , x_n$, with coefficients in $\mathbb{Z}$.  

\end {thm}

\section {Cluster Algebras from Surfaces}

See \cite{fomin2006cluster}, \cite{fock2006moduli} for details on the topics presented in this section. In particular, we remind the reader that tagged arcs are required for the general discussion but will not be mentioned here.

\begin {defn}
   An oriented surface $\Sigma$ with nonempty boundary is called a \textit{marked surface} if $\Sigma$ comes equipped with a collection of finitely many marked points on its boundary.
\end {defn}

\begin {defn}
    An \textit{arc} is any curve inside $\Sigma$ with endpoints at marked points, considered up to isotopy relative the set of marked points on the boundary of $\Sigma$, and such that the relative interior of this curve is disjoint from the boundary of $\Sigma$. Any curve beginning and ending at distinct marked points which lies entirely within the boundary of $\Sigma$ and does not contain any marked points in its interior is called a \textit{boundary segment}.
\end {defn}

\begin {defn}
    Two arcs in $\Sigma$ are called \textit{compatible} if they have isotopy representatives that do not intersect, except possibly at endpoints. An \textit{ideal triangulation} $\Delta$ is a maximal collection of distinct pairwise compatible arcs, along with all boundary segments. The arcs of a triangulation cut the surface into \textit{ideal triangles}. A \textit{flip} (or \textit{Whitehead move}) of an ideal triangulation at an arc $\gamma$ is the process of removing $\gamma$ from $\Delta$ and replacing it with the unique arc $\widetilde{\gamma}$ that gives another ideal triangulation of $\Sigma$ (see Figure \ref{fig:flip}).
\end {defn}

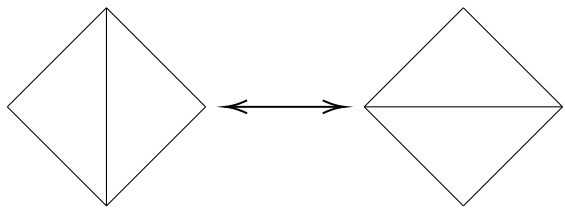
\begin {figure}[h!]
\centering
\caption {Flip inside a quadrilateral}
\label {fig:flip}
\begin{tikzpicture}[x=0.75pt,y=0.75pt,yscale=-1,xscale=1]

\draw    (200,90) -- (250,140) ;
\draw    (150,140) -- (200,90) ;
\draw    (150,140) -- (200,190) ;
\draw    (200,190) -- (250,140) ;
\draw    (200,90) -- (200,190) ;
\draw    (330,140) -- (380,90) ;
\draw    (380,90) -- (430,140) ;
\draw    (330,140) -- (380,190) ;
\draw    (430,140) -- (380,190) ;
\draw    (330,140) -- (430,140) ;
\draw [line width=0.75]    (262,140) -- (318,140) ;
\draw [shift={(320,140)}, rotate = 180] [color={rgb, 255:red, 0; green, 0; blue, 0 }  ][line width=0.75]    (10.93,-3.29) .. controls (6.95,-1.4) and (3.31,-0.3) .. (0,0) .. controls (3.31,0.3) and (6.95,1.4) .. (10.93,3.29)   ;
\draw [shift={(260,140)}, rotate = 0] [color={rgb, 255:red, 0; green, 0; blue, 0 }  ][line width=0.75]    (10.93,-3.29) .. controls (6.95,-1.4) and (3.31,-0.3) .. (0,0) .. controls (3.31,0.3) and (6.95,1.4) .. (10.93,3.29)   ;

\end{tikzpicture}

\end {figure}

It is known that any two ideal triangulations of $\Sigma$ are connected by a sequence of flips.

\begin{defn}

    Let $\Sigma$ be a marked surface. We now construct a cluster algebra $\mathcal{A}(\Sigma)$, called the \textit{cluster algebra from the surface $\Sigma$}, that depends only on $\Sigma$. To do this, we construct a quiver whose nodes are in one-to-one correspondence with the collection of arcs and boundary segments of some ideal triangulation $\Delta$ of $\Sigma$, and whose arrows form clockwise $3$-cycles inside each ideal triangle. This construction does not depend on the choice of ideal triangulation.

\end{defn}

The correspondences below follow from the construction given in the previous definition. $$\text{initial cluster variables of } \mathcal{A}(\Sigma) \longleftrightarrow \text{ arcs in } \Delta$$ $$\text{non-initial cluster variables of } \mathcal{A}(\Sigma) \longleftrightarrow \text{ arcs in } \Sigma \text{ that are not in } \Delta$$ $$\text{frozen variables of } \mathcal{A}(\Sigma) \longleftrightarrow \text{ boundary segments of }  \Sigma$$ $$\text{seeds of } \mathcal{A}(\Sigma) \longleftrightarrow \text{ triangulations } \Delta \text{ of } \Sigma$$ $$\text{ seed mutations in  } \mathcal{A}(\Sigma) \longleftrightarrow \text{ flips of arcs in } \Delta$$

\section {Cluster Algebras of Finite Type $A_n$}

\begin {defn}
  A \textit{cluster algebra of (finite) type $A_n$} is any cluster algebra from a quiver $Q = (Q_0 , Q_1 , s , t)$ such that the induced subquiver on the mutable vertices $Q_{0}^{mutable}$ is mutation equivalent to some orientation of a type $A_n$ Dynkin diagram.
\end {defn}

Any cluster algebra from a surface $\mathcal{A} (\Sigma)$ such that $\Sigma$ is an \textit{$(n+3)$-gon}, i.e. $\Sigma$ is a closed disc with $n+3$ marked points on its boundary, is a rank $n$ cluster algebra of type $A_n$. As is the case for general surfaces, cluster variables correspond to arcs in the $(n+3)$-gon, frozen variables correspond to the boundary segments, seeds correspond to triangulations, and seed mutations correspond to flips of arcs.

\begin {ex}
   Figure \ref{fig:seed} below shows one seed of a surface cluster algebra of type $A_3.$

   \begin {figure}[h!] \label{}
    \centering
    \caption{One seed for $A_3$.}
     \label{fig:seed}

\begin{tikzpicture}[x=0.75pt,y=0.75pt,yscale=-1,xscale=1]

\draw    (546.08,88.51) -- (620.92,163.34) ;
\draw    (620.92,163.34) -- (546.08,238.18) ;
\draw    (433.83,88.51) -- (546.08,88.51) ;
\draw    (433.83,238.18) -- (546.08,238.18) ;
\draw    (433.83,88.51) -- (359,163.34) ;
\draw    (359,163.34) -- (433.83,238.18) ;

\draw    (433.83,88.51) -- (433.83,238.18) ;
\draw    (546.08,88.51) -- (546.08,238.18) ;
\draw    (546.08,88.51) -- (433.83,238.18) ;

\draw (440,160.5) node    {$1$};
\draw (503,160.5) node    {$2$};
\draw (553,160.5) node    {$3$};
\draw (491,79.5) node    {$4$};
\draw (492.38,246.49) node    {$5$};
\draw (392.22,208.34) node    {$6$};
\draw (392.83,119.11) node    {$7$};
\draw (588.55,207.26) node    {$8$};
\draw (588.16,117.69) node    {$9$};

\end{tikzpicture}

\end{figure}
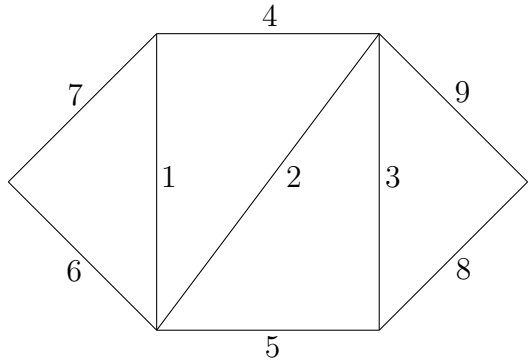

\end{ex}

We will use this seed as the basis for a running example to be followed throughout the rest of this text.

\chapter {Combinatorial Constructions}

In this chapter, we recall several interrelated objects and constructions naturally occuring in type $A$ cluster combinatorics. The objects found in this chapter are parameterized by binary words.

\section {Words}

\begin{defn}
    A \textit{(binary) word of length $n-1$} is a finite string formed from $n-1$ choices of elements from the set $ \{ a,b \}$. A word of length $n-1$ will be denoted $w = w_1 w_2 ... w_{n-1} $ and its \textit{length} is $l(w) = n-1$.
\end{defn}

The word 
$$w = \underbrace{aa \cdots a}_{k_1 \text{times}} \underbrace{bb \cdots b}_{k_2 \text{times}} \cdots$$

will be abbreviated $w = a^{k_1} b^{k_2} \cdots.$ 

As mentioned above, the constructions that follow are parameterized by the words $w$.

\begin{defn}
    A word $w$ is \textit{straight} if the only letter occurring in $w$ is $a$, or the only letter occurring in $w$ is $b$. Conversely, a word is \textit{zigzag} if neither of the substrings $aa$ nor $bb$ occur in $w$.
\end{defn}

\begin{ex}
    We fix the word $w=ab$ for our concrete running example to be followed throughout this section. Note that the word $w=ab$ is zigzag. 
\end{ex}

\section {Type $A_n$ Dynkin Quivers}

A \textit{type $A_n$ Dynkin diagram} is an undirected graph with $n$ nodes $1,2,...,n$ and one edge between each pair of consecutive nodes $i$ and $i+1$ for $1 \leq i \leq n-1.$ We picture any type $A_n$ Dynkin diagram as shown in Figure \ref{fig:A_n}. 

\begin {figure}[h!]
    \centering
    \caption{Type $A_n$ Dynkin diagram}
    \label{fig:A_n}
    \begin{tikzcd}
1 \arrow[r, no head] & 2 \arrow[r, no head] & \cdots \arrow[r, no head] & n
\end{tikzcd}
\end {figure}

Order the nodes $1,2,...,n$ and edges $\overline{i,i+1}$ of any type $A_n$ Dynkin diagram by $1 < 2 < ... < n$ and $\overline{1,2} < \overline{2,3} < ... < \overline{n-1,n}$, respectively. An \textit{orientation} of a type $A_n$ Dynkin diagram is a choice of orientation for each of the $n-1$ edges of $A_n$.

\begin {defn}
   A \textit{type $A_n$ Dynkin quiver} is a quiver that is mutation equivalent to an orientation of a type $A_n$ Dynkin diagram.
\end {defn}

Figure \ref{fig:generic_dynkin_quiver} shows one of the $2^{n-1}$ possible orientations of a Dynkin diagram of type $A_n,$ each of which is an example of a type $A_n$ Dynkin quiver (since it is mutation equivalent to itself).

\begin {figure}[h!]
    \centering
    \caption{Orientation of a type $A_n$ Dynkin diagram}
    \label{fig:generic_dynkin_quiver}
    \begin{tikzcd}
1 \arrow[r] & 2 \arrow[r] & \cdots \arrow[r] & n
\end{tikzcd}
\end {figure}

\begin {defn}
    Let $w = w_1 w_2 ... w_{n-1}$ be a word of length $l(w) = n-1$, and $A_n$ the Dynkin diagram of type $A_n$ with nodes labeled $1,2,...,n$. The \textit{type $A_n$ Dynkin quiver $A_w$ associated to $w$} is defined by mapping each $w_i$ to the $i^{th}$ edge $\overline{i,i+1}$ of the Dynkin diagram $A_{n}$ and declaring that any edge labeled by $a$ becomes oriented $i \longleftarrow i+1$, and that any edge labeled by $b$ becomes oriented $i \longrightarrow i+1$.
\end {defn}

\begin{ex}
    The word $w=ab$ gives the Dynkin quiver $A_w$ shown in Figure \ref{fig:A_ab_2}. 
\end{ex}

\begin {figure}[h!]
    \centering
    \caption{The Dynkin quiver $A_{ab}$}
    \label{fig:A_ab_2}

    \begin{tikzcd}
1 & 2 \arrow[l, "a", swap] \arrow[r, "b"] & 3
\end{tikzcd}

\end {figure}

\begin {rmk}

    If the word $w$ is straight then the edges in the Dynkin quiver $A_w$ are all oriented in the same direction. Conversely, if the word $w$ is zigzag then the edges in $A_w$ alternate in orientation. 
    
\end {rmk}

\section {Posets}

\bigskip

The posets we define here are called \textit{piecewise-linear posets} in \cite{baileycluster}, and \textit{zig-zag chains} in \cite{knauer2018lattice}.

\begin {defn}
    Define the \textit{poset $C_w$ associated to $w$} to be the Hasse diagram of a poset whose underlying graph is the Dynkin diagram $A_n$ and covering relations are $i \lessdot j$ in $C_w$ iff $i \rightarrow j$ in $A_w$. We visualize the edges of $C_w$ as taking unit diagonal steps upwards.
\end {defn}

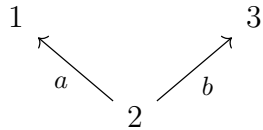
\begin {figure}[h!]
    \centering
    \caption{The poset $C_{ab}$}
    \label{fig:C_ab}
    
  \begin{tikzcd}
1 && 3 \\
& 2 \arrow[ul, "a"] \arrow[ur, "b", swap]
\end{tikzcd}

\end {figure}

\begin {ex}
   Figure \ref{fig:C_ab} shows the poset $C_{ab}.$ 
\end {ex}

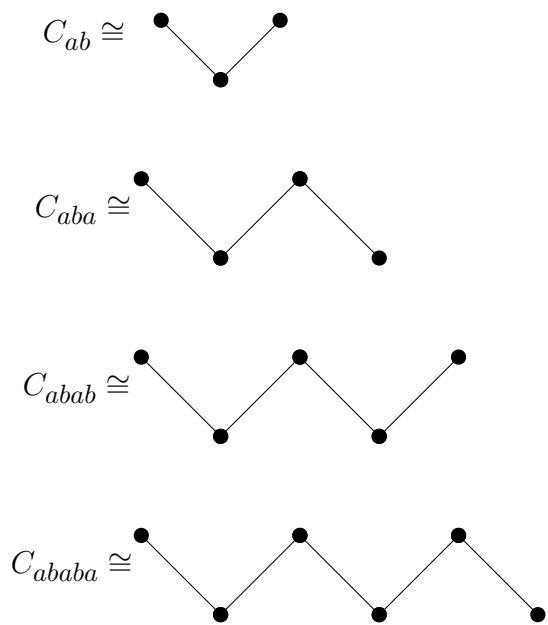
\begin {figure}[h!]
    \centering
    \caption{Four examples of fence posets $C_w$}
    \label{fig:fence}
    \begin{tikzpicture}[x=0.75pt,y=0.75pt,yscale=-1,xscale=1]

\draw    (290,40) -- (260,70) ;
\draw [shift={(260,70)}, rotate = 135] [color={rgb, 255:red, 0; green, 0; blue, 0 }  ][fill={rgb, 255:red, 0; green, 0; blue, 0 }  ][line width=0.75]      (0, 0) circle [x radius= 3.35, y radius= 3.35]   ;
\draw [shift={(290,40)}, rotate = 135] [color={rgb, 255:red, 0; green, 0; blue, 0 }  ][fill={rgb, 255:red, 0; green, 0; blue, 0 }  ][line width=0.75]      (0, 0) circle [x radius= 3.35, y radius= 3.35]   ;
\draw    (260,70) -- (230,40) ;
\draw [shift={(230,40)}, rotate = 225] [color={rgb, 255:red, 0; green, 0; blue, 0 }  ][fill={rgb, 255:red, 0; green, 0; blue, 0 }  ][line width=0.75]      (0, 0) circle [x radius= 3.35, y radius= 3.35]   ;
\draw [shift={(260,70)}, rotate = 225] [color={rgb, 255:red, 0; green, 0; blue, 0 }  ][fill={rgb, 255:red, 0; green, 0; blue, 0 }  ][line width=0.75]      (0, 0) circle [x radius= 3.35, y radius= 3.35]   ;
\draw    (300,120) -- (260,160) ;
\draw [shift={(260,160)}, rotate = 135] [color={rgb, 255:red, 0; green, 0; blue, 0 }  ][fill={rgb, 255:red, 0; green, 0; blue, 0 }  ][line width=0.75]      (0, 0) circle [x radius= 3.35, y radius= 3.35]   ;
\draw [shift={(300,120)}, rotate = 135] [color={rgb, 255:red, 0; green, 0; blue, 0 }  ][fill={rgb, 255:red, 0; green, 0; blue, 0 }  ][line width=0.75]      (0, 0) circle [x radius= 3.35, y radius= 3.35]   ;
\draw    (260,160) -- (220,120) ;
\draw [shift={(220,120)}, rotate = 225] [color={rgb, 255:red, 0; green, 0; blue, 0 }  ][fill={rgb, 255:red, 0; green, 0; blue, 0 }  ][line width=0.75]      (0, 0) circle [x radius= 3.35, y radius= 3.35]   ;
\draw [shift={(260,160)}, rotate = 225] [color={rgb, 255:red, 0; green, 0; blue, 0 }  ][fill={rgb, 255:red, 0; green, 0; blue, 0 }  ][line width=0.75]      (0, 0) circle [x radius= 3.35, y radius= 3.35]   ;
\draw    (380,210) -- (340,250) ;
\draw [shift={(340,250)}, rotate = 135] [color={rgb, 255:red, 0; green, 0; blue, 0 }  ][fill={rgb, 255:red, 0; green, 0; blue, 0 }  ][line width=0.75]      (0, 0) circle [x radius= 3.35, y radius= 3.35]   ;
\draw [shift={(380,210)}, rotate = 135] [color={rgb, 255:red, 0; green, 0; blue, 0 }  ][fill={rgb, 255:red, 0; green, 0; blue, 0 }  ][line width=0.75]      (0, 0) circle [x radius= 3.35, y radius= 3.35]   ;
\draw    (340,160) -- (300,120) ;
\draw [shift={(300,120)}, rotate = 225] [color={rgb, 255:red, 0; green, 0; blue, 0 }  ][fill={rgb, 255:red, 0; green, 0; blue, 0 }  ][line width=0.75]      (0, 0) circle [x radius= 3.35, y radius= 3.35]   ;
\draw [shift={(340,160)}, rotate = 225] [color={rgb, 255:red, 0; green, 0; blue, 0 }  ][fill={rgb, 255:red, 0; green, 0; blue, 0 }  ][line width=0.75]      (0, 0) circle [x radius= 3.35, y radius= 3.35]   ;
\draw    (300,210) -- (260,250) ;
\draw [shift={(260,250)}, rotate = 135] [color={rgb, 255:red, 0; green, 0; blue, 0 }  ][fill={rgb, 255:red, 0; green, 0; blue, 0 }  ][line width=0.75]      (0, 0) circle [x radius= 3.35, y radius= 3.35]   ;
\draw [shift={(300,210)}, rotate = 135] [color={rgb, 255:red, 0; green, 0; blue, 0 }  ][fill={rgb, 255:red, 0; green, 0; blue, 0 }  ][line width=0.75]      (0, 0) circle [x radius= 3.35, y radius= 3.35]   ;
\draw    (260,250) -- (220,210) ;
\draw [shift={(220,210)}, rotate = 225] [color={rgb, 255:red, 0; green, 0; blue, 0 }  ][fill={rgb, 255:red, 0; green, 0; blue, 0 }  ][line width=0.75]      (0, 0) circle [x radius= 3.35, y radius= 3.35]   ;
\draw [shift={(260,250)}, rotate = 225] [color={rgb, 255:red, 0; green, 0; blue, 0 }  ][fill={rgb, 255:red, 0; green, 0; blue, 0 }  ][line width=0.75]      (0, 0) circle [x radius= 3.35, y radius= 3.35]   ;
\draw    (340,250) -- (300,210) ;
\draw [shift={(300,210)}, rotate = 225] [color={rgb, 255:red, 0; green, 0; blue, 0 }  ][fill={rgb, 255:red, 0; green, 0; blue, 0 }  ][line width=0.75]      (0, 0) circle [x radius= 3.35, y radius= 3.35]   ;
\draw [shift={(340,250)}, rotate = 225] [color={rgb, 255:red, 0; green, 0; blue, 0 }  ][fill={rgb, 255:red, 0; green, 0; blue, 0 }  ][line width=0.75]      (0, 0) circle [x radius= 3.35, y radius= 3.35]   ;
\draw    (380,300) -- (340,340) ;
\draw [shift={(340,340)}, rotate = 135] [color={rgb, 255:red, 0; green, 0; blue, 0 }  ][fill={rgb, 255:red, 0; green, 0; blue, 0 }  ][line width=0.75]      (0, 0) circle [x radius= 3.35, y radius= 3.35]   ;
\draw [shift={(380,300)}, rotate = 135] [color={rgb, 255:red, 0; green, 0; blue, 0 }  ][fill={rgb, 255:red, 0; green, 0; blue, 0 }  ][line width=0.75]      (0, 0) circle [x radius= 3.35, y radius= 3.35]   ;
\draw    (300,300) -- (260,340) ;
\draw [shift={(260,340)}, rotate = 135] [color={rgb, 255:red, 0; green, 0; blue, 0 }  ][fill={rgb, 255:red, 0; green, 0; blue, 0 }  ][line width=0.75]      (0, 0) circle [x radius= 3.35, y radius= 3.35]   ;
\draw [shift={(300,300)}, rotate = 135] [color={rgb, 255:red, 0; green, 0; blue, 0 }  ][fill={rgb, 255:red, 0; green, 0; blue, 0 }  ][line width=0.75]      (0, 0) circle [x radius= 3.35, y radius= 3.35]   ;
\draw    (260,340) -- (220,300) ;
\draw [shift={(220,300)}, rotate = 225] [color={rgb, 255:red, 0; green, 0; blue, 0 }  ][fill={rgb, 255:red, 0; green, 0; blue, 0 }  ][line width=0.75]      (0, 0) circle [x radius= 3.35, y radius= 3.35]   ;
\draw [shift={(260,340)}, rotate = 225] [color={rgb, 255:red, 0; green, 0; blue, 0 }  ][fill={rgb, 255:red, 0; green, 0; blue, 0 }  ][line width=0.75]      (0, 0) circle [x radius= 3.35, y radius= 3.35]   ;
\draw    (340,340) -- (300,300) ;
\draw [shift={(300,300)}, rotate = 225] [color={rgb, 255:red, 0; green, 0; blue, 0 }  ][fill={rgb, 255:red, 0; green, 0; blue, 0 }  ][line width=0.75]      (0, 0) circle [x radius= 3.35, y radius= 3.35]   ;
\draw [shift={(340,340)}, rotate = 225] [color={rgb, 255:red, 0; green, 0; blue, 0 }  ][fill={rgb, 255:red, 0; green, 0; blue, 0 }  ][line width=0.75]      (0, 0) circle [x radius= 3.35, y radius= 3.35]   ;
\draw    (420,340) -- (380,300) ;
\draw [shift={(380,300)}, rotate = 225] [color={rgb, 255:red, 0; green, 0; blue, 0 }  ][fill={rgb, 255:red, 0; green, 0; blue, 0 }  ][line width=0.75]      (0, 0) circle [x radius= 3.35, y radius= 3.35]   ;
\draw [shift={(420,340)}, rotate = 225] [color={rgb, 255:red, 0; green, 0; blue, 0 }  ][fill={rgb, 255:red, 0; green, 0; blue, 0 }  ][line width=0.75]      (0, 0) circle [x radius= 3.35, y radius= 3.35]   ;

\draw (169,39) node [anchor=north west][inner sep=0.75pt]    {$C_{ab} \cong $};
\draw (165,127) node [anchor=north west][inner sep=0.75pt]    {$C_{aba} \cong $};
\draw (159,217) node [anchor=north west][inner sep=0.75pt]    {$C_{abab} \cong $};
\draw (153,307) node [anchor=north west][inner sep=0.75pt]    {$C_{ababa} \cong $};

\end{tikzpicture}
  
\end {figure}

\begin{rmk}
     If the word $w$ of length $l(w) = n-1$ is straight then the poset $C_w$ is a linear chain with $n$ elements and $n-1$ edges. In this case, the covering relations are $$1 > 2 > 3 > \cdots > n$$ if $w = a^{n-1},$ or $$1 < 2 < 3 < \cdots < n$$
    if $w = b^{n-1}.$ Conversely, if the word $w$ is zigzag, then the poset $C_w$ is a \textit{fence} or \textit{zigzag poset} (see \cite{munarini2002rank}) with $n$ elements and $n-1$ edges, with covering relations $$1 > 2 < 3 > \cdots $$ if $w = ababa \cdots,$ or $$1 < 2 > 3 < \cdots$$ if $w = babab \cdots$. See Figure \ref{fig:fence} for four examples of a fence poset from a word $w$. 
\end{rmk}
    
\section {Triangulations}

Form a new quiver $Q_w$, containing $A_w$ as a complete subquiver, by adding $n+3$ frozen nodes and $2n+4$ directed edges to and from $A_w$ as follows:

\begin{itemize}
    \item For each edge $\overline{i,i+1}$ in $A_w$, introduce a node $n+i$ and two directed edges between $n+i$ and the endpoints $i$ and $i+1$ of $\overline{i,i+1}$ such that a clockwise $3$-cycle is formed. See Figure \ref{fig:3cycle_1}.

\begin {figure}[h!]
    \centering
    \caption{Clockwise $3$-cycles created from edges in $A_w$}
    \label{fig:3cycle_1}
     \begin{tikzcd}
             & i+n \arrow[rd] &                &           &              &                &                \\
\cdots i \arrow[ru] &                & i+1 \cdots \arrow[ll] & \text{or} & \cdots i \arrow[rr] &                & i+1 \cdots \arrow[ld] \\
             &                &                &           &              & i+n \arrow[lu] &               
\end{tikzcd}
  
\end {figure}
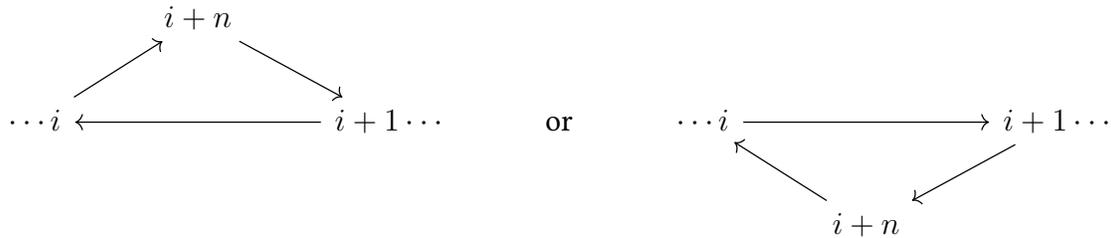

    \item Add two nodes $2n$ and $2n+1$ that form a clockwise $3$-cycle with the first node $1$ of $A_w$.

    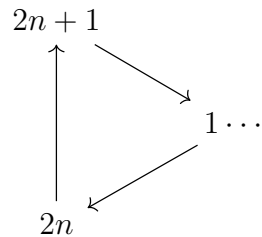
\begin {figure}[h!]
    \centering
    \caption{Leftmost $3$-cycle}
    \label{fig:3cycle_2}
    \begin{tikzcd}
2n+1 \arrow[rd] &              \\
                & 1 \arrow[ld] \cdots \\
2n \arrow[uu]   &             
\end{tikzcd}
  
\end {figure}

    \item Add two nodes $2n+2$ and $2n+3$ that form a clockwise $3$-cycle with the last node $n$ of $A_w$. If $l(w)$ is odd, or $l(w)$ is even and $w$ ends in $a^2$ or $b^2$, form the following $3$-cycle:

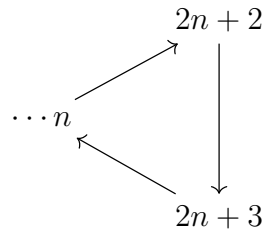
\begin {figure}[h!]
    \centering
    \caption{Rightmost $3$-cycle, if $l(w)$ is odd}
    \label{fig:3cycle_3}
    \begin{tikzcd}
             & 2n+2 \arrow[dd] \\
\cdots n \arrow[ru] &                 \\
             & 2n+3 \arrow[lu]
\end{tikzcd}
  
\end {figure}

If $l(w)$ is even and $w$ ends in $ba$ or $ab$, form the following $3$-cycle instead:

     \begin {figure}[h!]
    \centering
    \caption{Rightmost $3$-cycle, if $l(w)$ is even}
    \label{fig:3cycle_4}
    \begin{tikzcd}
             & 2n+3 \arrow[dd] \\
\cdots n \arrow[ru] &                 \\
             & 2n+2 \arrow[lu]
\end{tikzcd}
  
\end{figure}
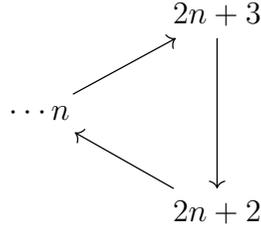

\end{itemize}

\newpage

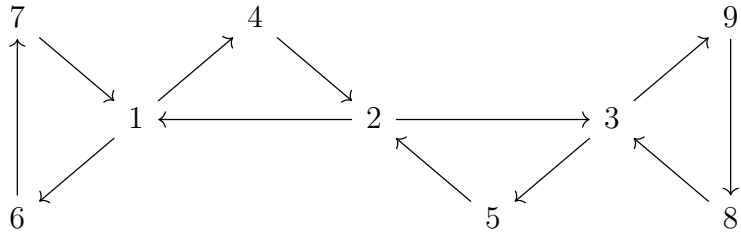
\begin {figure}[h!]
    \centering
    \caption{The quiver $Q_{ab}$}
    \label{fig:Q_ab}
    
  \begin{tikzcd}
7 \arrow[dr] && 4 \arrow[dr] &&&& 9 \arrow[dd] \\
& 1 \arrow[dl] \arrow[ur] && 2 \arrow[ll] \arrow[rr] && 3 \arrow[dl] \arrow[ur] \\
6 \arrow[uu] &&&& 5 \arrow[ul] && 8 \arrow[ul]
\end{tikzcd}
\end {figure}

\begin{ex}
     Figure \ref{fig:Q_ab} shows the quiver $Q_{ab}$. 
\end{ex}

\begin{defn} \label{Delta_w}
    The quiver $Q_w$ induces the \textit{triangulation $\Delta_w$ associated to $w$ of the $(n+3)$-gon $\Sigma$} whose elements are the internal diagonals labeled by $1 , 2 , \dots , n$ and boundary edges labeled by $n+1 , n+2 , \dots , 2n+3$.
\end{defn}

The ordering of nodes in $A_w \hookrightarrow Q_w$ induces an ordering of the internal diagonals $$\delta_{1} < \delta_{2} < \dots < \delta_{n}$$ of $\Delta_w$.

For $1 \leq i \leq n-1$, let $\Delta_i$ be the unique triangle cut out by $\Delta_w$ such that the two internal diagonals $\delta_{i}$ and $\delta_{i+1}$ are sides of $\Delta_i$. Let $\Delta_0$ be the unique triangle with sides consisting of two boundary segments and the internal diagonal $\delta_{1}$, and $\Delta_{n}$ the unique triangle with sides consisting of two boundary segments and the internal diagonal $\delta_{n}$. The ordering of the internal diagonals induces an ordering $\Delta_0 < \Delta_1 < \cdots < \Delta_n$ on the triangles $\Delta_i$.

By construction, the pair of consecutive triangles $\Delta_{i-1}$ and $\Delta_i$ each has precisely one edge labeled $i$, and there are no other common labels among their edges. We use the notation $\Sigma_{w} = [\Delta_0 , \Delta_1 , \dots , \Delta_n]$ to indicate that the surface $\Sigma$ with the triangulation $\Delta_w$ can be built by successively gluing $\Delta_{i+1}$ to $\Delta_{i}$ along the edge labeled $i+1$ for each $i$.

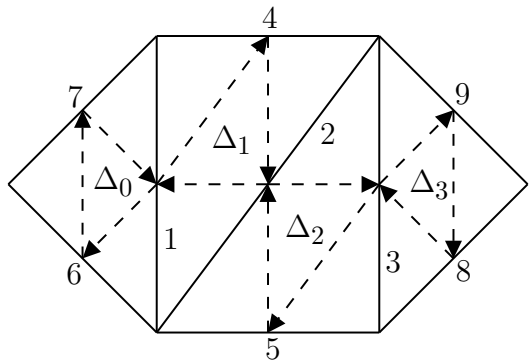
\begin {figure}[h!]
    \centering
    \caption{The triangulated polygon $\Sigma_{ab}$}
    \label{fig:Sigma_ab}

\tikzset{every picture/.style={line width=0.75pt}} 

\begin{tikzpicture}[x=0.75pt,y=0.75pt,yscale=-1,xscale=1]

\draw    (247.08,72.02) -- (321.92,146.85) ;
\draw    (321.92,146.85) -- (247.08,221.69) ;
\draw    (134.83,72.02) -- (247.08,72.02) ;
\draw    (134.83,221.69) -- (247.08,221.69) ;
\draw    (134.83,72.02) -- (60,146.85) ;
\draw    (60,146.85) -- (134.83,221.69) ;

\draw    (134.83,72.02) -- (134.83,221.69) ;
\draw    (247.08,72.02) -- (247.08,221.69) ;
\draw    (247.08,72.02) -- (134.83,221.69) ;

\draw  [dash pattern={on 4.5pt off 4.5pt}]  (97.42,109.43) -- (132.71,144.73) ;
\draw [shift={(134.83,146.85)}, rotate = 225] [fill={rgb, 255:red, 0; green, 0; blue, 0 }  ][line width=0.08]  [draw opacity=0] (8.93,-4.29) -- (0,0) -- (8.93,4.29) -- cycle    ;
\draw  [dash pattern={on 4.5pt off 4.5pt}]  (134.83,146.85) -- (99.54,182.15) ;
\draw [shift={(97.42,184.27)}, rotate = 315] [fill={rgb, 255:red, 0; green, 0; blue, 0 }  ][line width=0.08]  [draw opacity=0] (8.93,-4.29) -- (0,0) -- (8.93,4.29) -- cycle    ;
\draw  [dash pattern={on 4.5pt off 4.5pt}]  (97.42,184.27) -- (97.42,112.43) ;
\draw [shift={(97.42,109.43)}, rotate = 450] [fill={rgb, 255:red, 0; green, 0; blue, 0 }  ][line width=0.08]  [draw opacity=0] (8.93,-4.29) -- (0,0) -- (8.93,4.29) -- cycle    ;
\draw  [dash pattern={on 4.5pt off 4.5pt}]  (190.96,72.02) -- (190.96,143.85) ;
\draw [shift={(190.96,146.85)}, rotate = 270] [fill={rgb, 255:red, 0; green, 0; blue, 0 }  ][line width=0.08]  [draw opacity=0] (8.93,-4.29) -- (0,0) -- (8.93,4.29) -- cycle    ;
\draw  [dash pattern={on 4.5pt off 4.5pt}]  (190.96,146.85) -- (137.83,146.85) ;
\draw [shift={(134.83,146.85)}, rotate = 360] [fill={rgb, 255:red, 0; green, 0; blue, 0 }  ][line width=0.08]  [draw opacity=0] (8.93,-4.29) -- (0,0) -- (8.93,4.29) -- cycle    ;
\draw  [dash pattern={on 4.5pt off 4.5pt}]  (134.83,146.85) -- (189.16,74.42) ;
\draw [shift={(190.96,72.02)}, rotate = 486.87] [fill={rgb, 255:red, 0; green, 0; blue, 0 }  ][line width=0.08]  [draw opacity=0] (8.93,-4.29) -- (0,0) -- (8.93,4.29) -- cycle    ;
\draw  [dash pattern={on 4.5pt off 4.5pt}]  (190.96,146.85) -- (244.08,146.85) ;
\draw [shift={(247.08,146.85)}, rotate = 180] [fill={rgb, 255:red, 0; green, 0; blue, 0 }  ][line width=0.08]  [draw opacity=0] (8.93,-4.29) -- (0,0) -- (8.93,4.29) -- cycle    ;
\draw  [dash pattern={on 4.5pt off 4.5pt}]  (247.08,146.85) -- (192.76,219.29) ;
\draw [shift={(190.96,221.69)}, rotate = 306.87] [fill={rgb, 255:red, 0; green, 0; blue, 0 }  ][line width=0.08]  [draw opacity=0] (8.93,-4.29) -- (0,0) -- (8.93,4.29) -- cycle    ;
\draw  [dash pattern={on 4.5pt off 4.5pt}]  (190.96,221.69) -- (190.96,149.85) ;
\draw [shift={(190.96,146.85)}, rotate = 450] [fill={rgb, 255:red, 0; green, 0; blue, 0 }  ][line width=0.08]  [draw opacity=0] (8.93,-4.29) -- (0,0) -- (8.93,4.29) -- cycle    ;
\draw  [dash pattern={on 4.5pt off 4.5pt}]  (247.08,146.85) -- (282.38,111.56) ;
\draw [shift={(284.5,109.43)}, rotate = 495] [fill={rgb, 255:red, 0; green, 0; blue, 0 }  ][line width=0.08]  [draw opacity=0] (8.93,-4.29) -- (0,0) -- (8.93,4.29) -- cycle    ;
\draw  [dash pattern={on 4.5pt off 4.5pt}]  (284.5,109.43) -- (284.5,181.27) ;
\draw [shift={(284.5,184.27)}, rotate = 270] [fill={rgb, 255:red, 0; green, 0; blue, 0 }  ][line width=0.08]  [draw opacity=0] (8.93,-4.29) -- (0,0) -- (8.93,4.29) -- cycle    ;
\draw  [dash pattern={on 4.5pt off 4.5pt}]  (284.5,184.27) -- (249.21,148.97) ;
\draw [shift={(247.08,146.85)}, rotate = 405] [fill={rgb, 255:red, 0; green, 0; blue, 0 }  ][line width=0.08]  [draw opacity=0] (8.93,-4.29) -- (0,0) -- (8.93,4.29) -- cycle    ;

\draw (141.98,174.58) node    {$1$};
\draw (221.1,121.41) node    {$2$};
\draw (254.16,186.08) node    {$3$};
\draw (192,63.01) node    {$4$};
\draw (193.38,230) node    {$5$};
\draw (93.22,191.85) node    {$6$};
\draw (93.83,102.62) node    {$7$};
\draw (289.55,190.77) node    {$8$};
\draw (289.16,101.2) node    {$9$};
\draw (272.5,147) node    {$\Delta _{3}$};
\draw (209.83,169.67) node    {$\Delta _{2}$};
\draw (172.83,124.5) node    {$\Delta _{1}$};
\draw (113,145.5) node    {$\Delta _{0}$};

\end{tikzpicture}

\end {figure}

 \begin{ex} \label{Sigma_ab}
      Figure \ref{fig:Sigma_ab} shows the triangulated polygon $\Sigma_{ab}$, with both edge labels $i$ and triangles $\Delta_i$ indicated (along with the quiver $Q_{ab},$ its edges pictured here with dashed arrows).
 \end{ex}

The quiver $Q_w$ induces a type $A_n$ cluster algebra $\mathcal{A} (\Sigma)_w$ with initial extended cluster equal to $(x_1 , x_2 , \dots , x_n , x_{n+1} , \dots , x_{2n+3}).$ In the next section we assign to any word $w$ a cluster variable in the associated cluster algebra.

\begin{defn} \label{fan_def}
    We say that $\Delta_w$ is a \textit{fan triangulation} if there exists some vertex $v$ of $\Sigma$ such that each internal diagonal $\delta_{1}, \delta_{2}, \dots , \delta_{n}$ has $v$ as one of its endpoints. We say that $\Delta_w$ is a \textit{zigzag triangulation} if no three internal diagonals share a common endpoint. 
\end{defn}

\begin{rmk} \label{striaght_zigzag_triangulation}
    If $w$ is straight, then $\Delta_w$ is a fan triangulation. Conversely, if $w$ is zigzag, then $\Delta_w$ is zigzag triangulation.  
\end{rmk}

\begin{ex}
    The triangulation $\Delta_{ab}$ shown in Example \ref{Sigma_ab} is a zigzag triangulation. Figure \ref{fig:fan} below shows one example of a fan triangulation. This particular triangulation will be encountered again later, starting in Chapter \ref{dual_const_section}.
\end{ex}

\begin {figure}[h!]
    \centering
    \caption{A fan triangulation}
    \label{fig:fan}
    \begin{tikzpicture}[x=0.75pt,y=0.75pt,yscale=-1,xscale=1]

\draw    (302.25,104) -- (344.5,146.25) ;
\draw    (260,146.25) -- (302.25,104) ;
\draw    (344.5,230.75) -- (344.5,146.25) ;
\draw    (260,230.75) -- (260,146.25) ;
\draw    (260,230.75) -- (302.25,273) ;
\draw    (302.25,273) -- (344.5,230.75) ;
\draw    (260,230.75) -- (302.25,104) ;
\draw    (344.5,230.75) -- (302.25,104) ;
\draw    (302.25,273) -- (302.25,104) ;

\end{tikzpicture}

\end {figure}
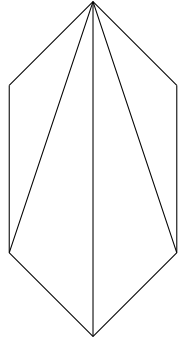

\section{Arcs, Cluster Variables, and Resolutions}

\begin{defn}
    Let $a$ be the vertex of $\Delta_0$ that is not an endpoint of edge $\delta_{1}$, and let $b$ be the vertex of $\Delta_n$ that is not an endpoint of edge $\delta_{n}$. Consider the oriented arc $\gamma_w = \gamma_{a \rightarrow b}$ in $\Sigma_w$ with initial vertex $a$ and terminal vertex $b$. We call the arc $\gamma_w$ the \textit{arc in $\Sigma_w$ associated to $w$}, and the resulting cluster variable $x_{w}$ the \textit{cluster variable associated to $w$}.
\end{defn}

\begin{ex}
    Figure \ref{fig:gamma_ab} below shows the arc $\gamma_{ab}$ and the associated cluster variable $x_{ab}$, parameterized by the word $w=ab$. 
\end{ex}

\begin {figure}[h!]
    \centering
    \caption{The arc $\gamma_{ab}$ inside $\Sigma_{ab}$, and the associated cluster variable $x_{ab}$}
    \label{fig:gamma_ab}
\begin{tikzpicture}[x=0.75pt,y=0.75pt,yscale=-1,xscale=1]

\draw    (466.08,98.02) -- (540.92,172.85) ;
\draw    (540.92,172.85) -- (466.08,247.69) ;
\draw    (353.83,98.02) -- (466.08,98.02) ;
\draw    (353.83,247.69) -- (466.08,247.69) ;
\draw    (353.83,98.02) -- (279,172.85) ;
\draw    (279,172.85) -- (353.83,247.69) ;

\draw    (353.83,98.02) -- (353.83,247.69) ;
\draw    (466.08,98.02) -- (466.08,247.69) ;
\draw    (466.08,98.02) -- (353.83,247.69) ;

\draw [color={rgb, 255:red, 0; green, 0; blue, 0 }  ,draw opacity=1 ]   (279,172.85) -- (540.92,172.85) ;
\draw [shift={(540.92,172.85)}, rotate = 0] [color={rgb, 255:red, 0; green, 0; blue, 0 }  ,draw opacity=1 ][fill={rgb, 255:red, 0; green, 0; blue, 0 }  ,fill opacity=1 ][line width=0.75]      (0, 0) circle [x radius= 3.35, y radius= 3.35]   ;
\draw [shift={(279,172.85)}, rotate = 0] [color={rgb, 255:red, 0; green, 0; blue, 0 }  ,draw opacity=1 ][fill={rgb, 255:red, 0; green, 0; blue, 0 }  ,fill opacity=1 ][line width=0.75]      (0, 0) circle [x radius= 3.35, y radius= 3.35]   ;
\draw  [dash pattern={on 0.84pt off 2.51pt}]  (260.33,50.67) .. controls (252.26,151.29) and (316.26,139.29) .. (330.26,169.29) ;

\draw (268,172) node  [color={rgb, 255:red, 0; green, 0; blue, 0 }  ,opacity=1 ]  {$a$};
\draw (550,171) node  [color={rgb, 255:red, 0; green, 0; blue, 0 }  ,opacity=1 ]  {$b$};
\draw (440,180) node  [color={rgb, 255:red, 0; green, 0; blue, 0 }  ,opacity=1 ]  {$\gamma _{ab}$};
\draw (398,38) node  [color={rgb, 255:red, 0; green, 0; blue, 0 }  ,opacity=1 ]  {$x_{ab} =\frac{x^{2}_{2} x_{7} x_{8} +x_{2} x_{5} x_{7} x_{9} +x_{4} x_{5} x_{6} x_{9} +x_{2} x_{4} x_{6} x_{8} +x_{1} x_{3} x_{6} x_{9}}{x_{1} x_{2} x_{3}}$};
\draw (358.98,146.58) node    {$1$};
\draw (440.1,146.41) node    {$2$};
\draw (472.16,206.08) node    {$3$};
\draw (412,90.01) node    {$4$};
\draw (413.38,256) node    {$5$};
\draw (313.22,217.85) node    {$6$};
\draw (313.83,127.62) node    {$7$};
\draw (509.55,215.77) node    {$8$};
\draw (508.16,127.2) node    {$9$};

\end{tikzpicture}

\end {figure}
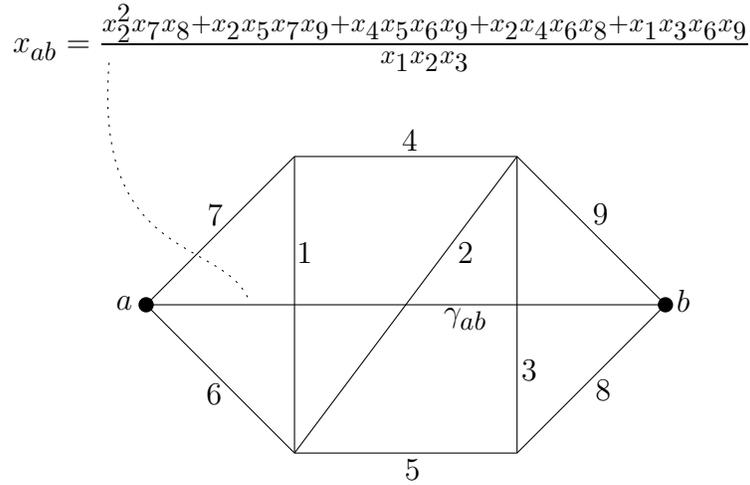

Any cluster variable $x_{w}$ can be written as $$x_{w} = \frac{f(x_1 , x_2 , \dots x_n)}{x_1 x_2 \dots x_n},$$ where $f$ is a polynomial with coefficients from $\mathbb{Z} [x_{n+1} , x_{n+2} , \dots , x_{2n+3}]$. The first goal of the remainder of this section is to explain the resolution process given in \cite{gekhtman2003cluster} used to compute the monomials in $f$, and hence the cluster variable $x_{w}$. The second goal is to define the set $\text{Res}(w)$ of resolutions associated to $w$, and the set $\text{Tree}(w)$ of resolution trees associated to $w$.

Fix $w$ and consider the arc $\gamma_w$ inside the triangulated polygon $\Sigma_w$. Recall the diagonals of $\Delta_w$ are $\delta_{1}, \delta_{2}, \dots , \delta_{n}$ and that by construction $\gamma_w$ crosses each of these $n$ internal diagonals, creating $n$ intersection points in $\Sigma_w$. Call these intersection points $p_i = \gamma_w \cap \delta_{i}.$

To resolve the intersection point $p_i$, we choose a small neighborhood of $p_i$ and replace it with a pair of nonintersecting smooth curves in one of two different ways, shown below in Figure \ref{fig:smoothing}.

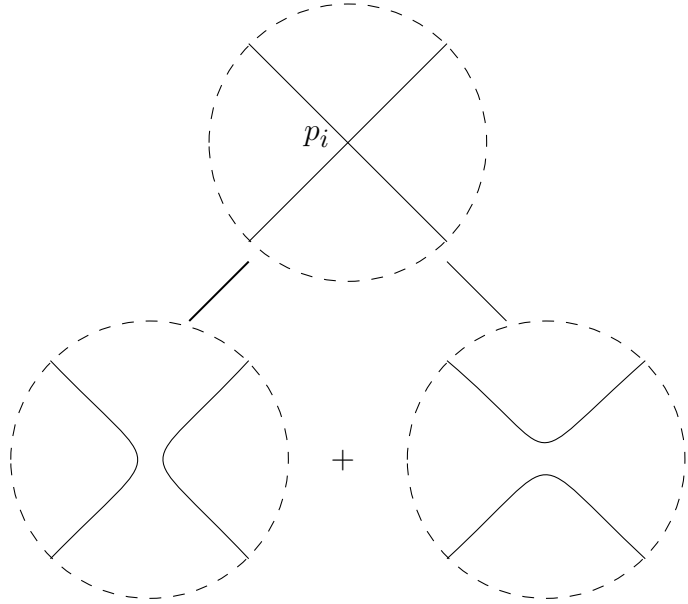
\begin {figure}[h!]
    \centering
    \caption{Resolution of the intersection point $p_i$}
    \label{fig:smoothing}
    \begin{tikzpicture}[x=0.75pt,y=0.75pt,yscale=-1,xscale=1]

\draw  [dash pattern={on 4.5pt off 4.5pt}] (150,290) .. controls (150,251.34) and (181.34,220) .. (220,220) .. controls (258.66,220) and (290,251.34) .. (290,290) .. controls (290,328.66) and (258.66,360) .. (220,360) .. controls (181.34,360) and (150,328.66) .. (150,290) -- cycle ;
\draw  [dash pattern={on 4.5pt off 4.5pt}] (350,290) .. controls (350,251.34) and (381.34,220) .. (420,220) .. controls (458.66,220) and (490,251.34) .. (490,290) .. controls (490,328.66) and (458.66,360) .. (420,360) .. controls (381.34,360) and (350,328.66) .. (350,290) -- cycle ;
\draw  [dash pattern={on 4.5pt off 4.5pt}] (250,130) .. controls (250,91.34) and (281.34,60) .. (320,60) .. controls (358.66,60) and (390,91.34) .. (390,130) .. controls (390,168.66) and (358.66,200) .. (320,200) .. controls (281.34,200) and (250,168.66) .. (250,130) -- cycle ;
\draw    (270,80) -- (370,180) ;
\draw    (170,240) .. controls (228,300.47) and (229.33,279.13) .. (170,340) ;
\draw    (270,240) .. controls (213.33,301.13) and (211.33,277.8) .. (270,340) ;
\draw    (270,180) -- (370,80) ;
\draw    (370,340) .. controls (430,284.47) and (408.67,282.47) .. (470,340) ;
\draw    (370,240) .. controls (431.5,292.2) and (404.5,298.2) .. (470,240) ;
\draw [line width=0.75]    (270,190) -- (240,220) ;
\draw    (370,190) -- (400,220) ;

\draw (317.5,290) node    {$+$};
\draw (304.5,127) node    {$p_{i}$};

\end{tikzpicture}

\end {figure}

Choose one resolution out of the two above for each intersection point $p_i$; this results in a collection of $n+1$ nonintersecting curves in $\Sigma.$ Note that closed curves based at some $v \in \Sigma$ can occur, but that this process never leads to a closed curve that is not attached to some vertex of $\Sigma$.

\begin{defn}
    The set of \textit{resolutions $\text{Res}(w)$ associated to $w$} consists of those diagrams that can be obtained from resolving each $p_i$ in one of the two possible ways, in some chosen order.
\end{defn}

Each element in $\text{Res}(w)$ is weighted as follows. First, replace each arc in $r \in \text{Res}(w)$ with distinct endpoints with the arc or boundary segment from $\Delta_w$ that it is isotopic to. Let $E(r)$ be the collection of arcs and boundary segments from $\Delta_w$ produced from the resolution $r$, along with $\varnothing$ if any closed loops are present. Define the \textit{weight} of any resolution $r$ to be $x_{r} = \prod_{j \in E(r)} x_{j},$ where $x_{\varnothing} = 0.$

\begin{prop} (Proposition 2.1 in \cite{gekhtman2003cluster})
    Fix the word $w$. Consider the arc $\gamma_w$ in the triangulated polygon $\Sigma_w$ triangulated by $\Delta_w$. Then any internal diagonals obtained by a resolution belong to $\Delta_w$. The cluster variable $x_{w}$ is equal to 
    $$x_{w} = \frac{1}{x_1 x_2 \dots x_n} \sum_{r \in \text{Res}(w)} x_{r}.$$
\end{prop}

We now describe how to produce a \textit{resolution tree from $w$}. Each node of such a tree is a diagram of arcs inside the $(n+3)$-gon $\Sigma$, and is weighted by the product of cluster variables associated to those arcs (or zero if there is a closed loop in the diagram). The root of any resolution tree from $w$ is the diagram consisting of the arc $\gamma_w$ inside $\Sigma_w$. Choosing an intersection point $p_i$ to resolve at creates two children of this root (see Figure \ref{fig:smoothing}). If we continue along in this way (choosing an intersection point to resolve at in each child, etc.), and halt whenever we have resolved every intersection point, a binary tree (with additional node structure) is produced.

\begin{defn}
    The set of \textit{resolution trees $\text{Tree}(w)$ associated to $w$} is the set whose elements are the binary resolution trees from $w$ as described above.
\end{defn}

Note that $\text{Res}(w)$ is equal to the union of the leaves of the trees in $\text{Tree}(w).$

\begin{ex}
    The figure below shows one element of $\text{Tree}(w)$ for the word $w=ab.$ 
\end{ex}

\begin {figure}[h!]
    \centering
    \caption{One element of $\text{Tree}(ab)$}
    \label{fig:tree_ab}
    


\end {figure}

Although our construction of arcs seems restrictive, the next lemma shows that there is in fact no loss of generality.

\begin{lem}
    Any cluster variable associated to an arc in a polygon can be computed as $x_{w}$ for some word $w$. 
\end{lem}

\begin{proof}
    Let $\nu$ be an arc in the triangulated polygon $\Sigma_{\Delta}$, triangulated by $\Delta$. If $\nu$ crosses every internal diagonal in $\Delta$, then we are done. Otherwise, we work inside the triangulated subpolygon $\Sigma_{\Delta}^{\prime}$ in $\Sigma_{\Delta}$ obtained by deleting from $\Delta$ any edge that is not an edge of some triangle crossed by $\nu$. Furthermore, we ``freeze'' any edge $\eta$ bordering a deleted triangle, meaning we disallow this arc to be flipped, so that also $x_{\eta}$ cannot be mutated.

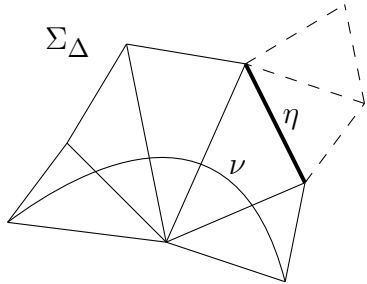
\begin {figure}[h!]
    \centering
    \caption{An arc in a triangulated subpolygon}
    \label{fig:subpolygon}
    \begin{tikzpicture}[x=0.75pt,y=0.75pt,yscale=-1,xscale=1]

\draw    (130,200) -- (160,150) ;
\draw    (130,200) -- (180,250) ;
\draw    (180,250) -- (160,150) ;
\draw    (160,150) -- (220,160) ;
\draw    (180,250) -- (220,160) ;
\draw [line width=1.5]    (250,220) -- (220,160) ;
\draw    (250,220) -- (180,250) ;
\draw    (240,270) -- (250,220) ;
\draw    (240,270) -- (180,250) ;
\draw  [dash pattern={on 4.5pt off 4.5pt}]  (280,180) -- (220,160) ;
\draw  [dash pattern={on 4.5pt off 4.5pt}]  (250,220) -- (280,180) ;
\draw    (130,200) -- (100,240) ;
\draw    (180,250) -- (100,240) ;
\draw [color={rgb, 255:red, 0; green, 0; blue, 0 }  ,draw opacity=1 ]   (100,240) .. controls (140,210) and (215,172.75) .. (240,270) ;
\draw  [dash pattern={on 4.5pt off 4.5pt}]  (280,180) -- (270,130) ;
\draw  [dash pattern={on 4.5pt off 4.5pt}]  (220,160) -- (270,130) ;

\draw (216,212) node    {$\nu $};
\draw (130.5,149) node    {$\Sigma _{\Delta }$};
\draw (243,188) node    {$\eta $};
\end{tikzpicture}
\end{figure}

The result now follows by noting that any arc obtained by resolving $\nu$ will be contained entirely within the triangulated subpolygon just mentioned. 

\end{proof}

\section{Snake Graphs}

\begin{defn}
    The \textit{snake graph $G_w$ associated to $w$} is the labeled planar graph recursively defined by the procedure given below.

\begin{enumerate}
    \item Choose an orientation-preserving embedding of the triangulated square $[\Delta_0 , \Delta_1 ]$ into the discrete plane $\mathbb{Z}^2$ such that its image $\widetilde{T}_1$ is a triangulated unit square  with vertices $(0,0),(1,0),(0,1),$ and $(1,1)$ in $\mathbb{Z}^2$, and such that the point $a \in \Delta_0$ maps to the point $(0,0)$. Note that the (line spanned by the) image of the triangulating edge will have slope $-1$.
    
    \item Choose an orientation-reversing map of $[\Delta_{1} , \Delta_{2} ]$ into $\mathbb{Z}^2$ such that its image $\widetilde{T}_2$ is a triangulated unit square (again, with triangulating edge having slope $-1$) glued to $\widetilde{T}_1$ along the unique edge in each $\widetilde{T_i}$ labeled $j \in \{ n+1 , ... , 2n+3 \}.$ Note that if the intersection point of the diagonals $\delta_{1}$ and $\delta_{2}$ is to the left (resp. right) of $\gamma$, then $\widetilde{T_2}$ is the triangulated square directly to the right of (resp. above) $\widetilde{T}_1$.
    
    \item Continue this process, using orientation-preserving maps for $i$ odd and orientation-reversing maps for $i$ even, to get the graph $\widetilde{G_w}$, built from triangulated unit squares in $\mathbb{Z}^2$ (with all triangulating edges having slope $-1$) glued either above or to the right of the previous square. Each $\widetilde{T_i}$ will be called a \textit{tile} of $\widetilde{G_w}$. The triangulating edge of each $\widetilde{T_i}$ is called the \textit{diagonal} of $\widetilde{T_i}.$
    
    \item The \textit{snake graph} $G_w$ is the graph in $\mathbb{Z}^2$ gotten by deleting each diagonal from each tile in $\widetilde{G_w}$. 
\end{enumerate}
\end{defn}

Let $T_i$ be the tile $\widetilde{T_i}$ after its diagonal has been removed. We will call $T_i$ a \textit{tile} of $G_w$. We will often refer to the corners (SW, SE, NE, NW) and edges (S,E,N,W) of a tile $T_i$ as indicated in the next figure.

\begin {figure}[h!]
    \centering
    \caption{Shorthand to describe the corners and edges of a tile $T_i$}
    \label{fig:tile}
    \begin{tikzpicture}[x=0.75pt,y=0.75pt,yscale=-1,xscale=1]

\draw    (203.99,129.85) -- (203.99,220.94) ;
\draw    (295.08,129.85) -- (295.08,220.94) ;
\draw    (203.99,129.85) -- (295.08,129.85) ;
\draw    (203.99,220.94) -- (295.08,220.94) ;

\draw (189.57,117.7) node    {$\text{NW}$};
\draw (311.02,117.7) node    {$\text{NE}$};
\draw (308.75,230.05) node    {$\text{SE}$};
\draw (184.25,230.05) node    {$\text{SW}$};
\draw (250.3,117.7) node    {$\text{N}$};
\draw (188.57,171.84) node    {$\text{W}$};
\draw (306.19,170.84) node    {$\text{E}$};
\draw (249.54,231.56) node    {$\text{S}$};
\draw (248.5,172) node    {$T_{i}$};

\end{tikzpicture}

\end{figure}
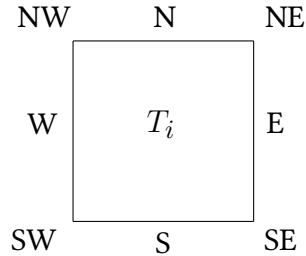

Order the tiles of $G_w$ by $T_1 < T_2 < \dots < T_n.$ A \textit{boundary edge} of $G_w$ is any edge of $G_w$ not occurring as a shared edge between any two of its consecutive tiles. Any edge that is not a boundary edge (i.e., each gluing edge) will be called an \textit{internal edge} of $G_w$. Let the internal edges of $G_w$ be labeled $e_1 , e_2 , \dots , e_{n-1}$, where $e_i$ is the gluing edge between tiles $T_i$ and $T_{i+1}$.

\begin{defn} \label{straight_zigzag_SG_def}
    A snake graph $G_w$ is called \textit{straight} if all of its tiles lie in a single row or column. A snake graph is called \textit{zigzag} if no three consecutive tiles are straight.
\end{defn}

\begin{ex}
    Below we illustrate the construction of the (straight) snake graph $G_{ab}$ associated to the (zigzag) word $w=ab.$

\begin {figure}[h!]
    \centering
    \caption{Construction of the snake graph $G_{ab}$}
    \label{fig:G_ab}
    \begin{tikzpicture}[x=0.75pt,y=0.75pt,yscale=-1,xscale=1]

\draw    (151.43,32.86) -- (200,81.43) ;
\draw    (200,81.43) -- (151.43,130) ;
\draw [color={rgb, 255:red, 0; green, 0; blue, 0 }  ,draw opacity=1 ][line width=1.5]    (78.57,32.86) -- (151.43,32.86) ;
\draw    (78.57,130) -- (151.43,130) ;
\draw [color={rgb, 255:red, 0; green, 0; blue, 0 }  ,draw opacity=1 ][line width=1.5]    (78.57,32.86) -- (30,81.43) ;
\draw [color={rgb, 255:red, 0; green, 0; blue, 0 }  ,draw opacity=1 ][line width=1.5]    (30,81.43) -- (78.57,130) ;
\draw [color={rgb, 255:red, 0; green, 0; blue, 0 }  ,draw opacity=1 ][line width=1.5]    (78.57,32.86) -- (78.57,130) ;
\draw    (151.43,32.86) -- (151.43,130) ;
\draw [color={rgb, 255:red, 0; green, 0; blue, 0 }  ,draw opacity=1 ][line width=1.5]    (151.43,32.86) -- (78.57,130) ;
\draw [color={rgb, 255:red, 0; green, 0; blue, 0 }  ,draw opacity=1 ] [dash pattern={on 0.84pt off 2.51pt}]  (30,81.43) -- (200,81.43) ;
\draw    (151.43,142.86) -- (200,191.43) ;
\draw    (200,191.43) -- (151.43,240) ;
\draw [color={rgb, 255:red, 0; green, 0; blue, 0 }  ,draw opacity=1 ][line width=1.5]    (78.57,142.86) -- (151.43,142.86) ;
\draw [color={rgb, 255:red, 0; green, 0; blue, 0 }  ,draw opacity=1 ][line width=1.5]    (78.57,240) -- (151.43,240) ;
\draw    (78.57,142.86) -- (30,191.43) ;
\draw    (30,191.43) -- (78.57,240) ;
\draw [color={rgb, 255:red, 0; green, 0; blue, 0 }  ,draw opacity=1 ][line width=1.5]    (78.57,142.86) -- (78.57,240) ;
\draw [color={rgb, 255:red, 0; green, 0; blue, 0 }  ,draw opacity=1 ][line width=1.5]    (151.43,142.86) -- (151.43,240) ;
\draw [color={rgb, 255:red, 0; green, 0; blue, 0 }  ,draw opacity=1 ][line width=1.5]    (151.43,142.86) -- (78.57,240) ;
\draw [color={rgb, 255:red, 0; green, 0; blue, 0 }  ,draw opacity=1 ] [dash pattern={on 0.84pt off 2.51pt}]  (30,191.43) -- (200,191.43) ;
\draw [color={rgb, 255:red, 0; green, 0; blue, 0 }  ,draw opacity=1 ][line width=1.5]    (151.43,252.86) -- (200,301.43) ;
\draw [color={rgb, 255:red, 0; green, 0; blue, 0 }  ,draw opacity=1 ][line width=1.5]    (200,301.43) -- (151.43,350) ;
\draw    (78.57,252.86) -- (151.43,252.86) ;
\draw [color={rgb, 255:red, 0; green, 0; blue, 0 }  ,draw opacity=1 ][line width=1.5]    (78.57,350) -- (151.43,350) ;
\draw    (78.57,252.86) -- (30,301.43) ;
\draw    (30,301.43) -- (78.57,350) ;
\draw    (78.57,252.86) -- (78.57,350) ;
\draw [color={rgb, 255:red, 0; green, 0; blue, 0 }  ,draw opacity=1 ][line width=1.5]    (151.43,252.86) -- (151.43,350) ;
\draw [color={rgb, 255:red, 0; green, 0; blue, 0 }  ,draw opacity=1 ][line width=1.5]    (151.43,252.86) -- (78.57,350) ;
\draw [color={rgb, 255:red, 0; green, 0; blue, 0 }  ,draw opacity=1 ] [dash pattern={on 0.84pt off 2.51pt}]  (30,301.43) -- (200,301.43) ;
\draw    (350,20) -- (350,120) ;
\draw    (410,80) -- (290,80) ;
\draw [color={rgb, 255:red, 0; green, 0; blue, 0 }  ,draw opacity=1 ][line width=1.5]    (350,60) -- (350,80) ;
\draw [color={rgb, 255:red, 0; green, 0; blue, 0 }  ,draw opacity=1 ][line width=1.5]    (370,60) -- (370,80) ;
\draw [color={rgb, 255:red, 0; green, 0; blue, 0 }  ,draw opacity=1 ][line width=1.5]    (370,80) -- (350,80) ;
\draw [color={rgb, 255:red, 0; green, 0; blue, 0 }  ,draw opacity=1 ][line width=1.5]    (370,60) -- (350,60) ;
\draw [color={rgb, 255:red, 0; green, 0; blue, 0 }  ,draw opacity=1 ]   (210,100) -- (298,100) ;
\draw [shift={(300,100)}, rotate = 180] [color={rgb, 255:red, 0; green, 0; blue, 0 }  ,draw opacity=1 ][line width=0.75]    (10.93,-3.29) .. controls (6.95,-1.4) and (3.31,-0.3) .. (0,0) .. controls (3.31,0.3) and (6.95,1.4) .. (10.93,3.29)   ;
\draw [shift={(210,100)}, rotate = 0] [color={rgb, 255:red, 0; green, 0; blue, 0 }  ,draw opacity=1 ][line width=0.75]      (0,-11.18) .. controls (-3.09,-11.18) and (-5.59,-8.68) .. (-5.59,-5.59) .. controls (-5.59,-2.5) and (-3.09,0) .. (0,0) ;
\draw    (350,130) -- (350,230) ;
\draw    (410,190) -- (290,190) ;
\draw [color={rgb, 255:red, 0; green, 0; blue, 0 }  ,draw opacity=1 ][line width=0.75]    (350,170) -- (350,190) ;
\draw [color={rgb, 255:red, 0; green, 0; blue, 0 }  ,draw opacity=1 ][line width=0.75]    (370,170) -- (370,190) ;
\draw [color={rgb, 255:red, 0; green, 0; blue, 0 }  ,draw opacity=1 ][line width=0.75]    (370,190) -- (350,190) ;
\draw [color={rgb, 255:red, 208; green, 2; blue, 27 }  ,draw opacity=1 ][line width=1.5]    (370,170) -- (350,170) ;
\draw [color={rgb, 255:red, 208; green, 2; blue, 27 }  ,draw opacity=1 ][line width=1.5]    (370,150) -- (350,150) ;
\draw [color={rgb, 255:red, 208; green, 2; blue, 27 }  ,draw opacity=1 ][line width=1.5]    (370,150) -- (370,170) ;
\draw [color={rgb, 255:red, 208; green, 2; blue, 27 }  ,draw opacity=1 ][line width=1.5]    (350,150) -- (350,170) ;
\draw    (350,240) -- (350,340) ;
\draw    (410,300) -- (290,300) ;
\draw [color={rgb, 255:red, 0; green, 0; blue, 0 }  ,draw opacity=1 ][line width=0.75]    (350,280) -- (350,300) ;
\draw [color={rgb, 255:red, 0; green, 0; blue, 0 }  ,draw opacity=1 ][line width=0.75]    (370,280) -- (370,300) ;
\draw [color={rgb, 255:red, 0; green, 0; blue, 0 }  ,draw opacity=1 ][line width=0.75]    (370,300) -- (350,300) ;
\draw [color={rgb, 255:red, 0; green, 0; blue, 0 }  ,draw opacity=1 ][line width=0.75]    (370,280) -- (350,280) ;
\draw [color={rgb, 255:red, 0; green, 0; blue, 0 }  ,draw opacity=1 ][line width=1.5]    (370,260) -- (350,260) ;
\draw [color={rgb, 255:red, 0; green, 0; blue, 0 }  ,draw opacity=1 ][line width=0.75]    (370,260) -- (370,280) ;
\draw [color={rgb, 255:red, 0; green, 0; blue, 0 }  ,draw opacity=1 ][line width=0.75]    (350,260) -- (350,280) ;
\draw [color={rgb, 255:red, 0; green, 0; blue, 0 }  ,draw opacity=1 ][line width=1.5]    (370,240) -- (350,240) ;
\draw [color={rgb, 255:red, 0; green, 0; blue, 0 }  ,draw opacity=1 ][line width=1.5]    (370,240) -- (370,260) ;
\draw [color={rgb, 255:red, 0; green, 0; blue, 0 }  ,draw opacity=1 ][line width=1.5]    (350,240) -- (350,260) ;
\draw [color={rgb, 255:red, 0; green, 0; blue, 0 }  ,draw opacity=1 ][line width=1.5]    (350,60) -- (370,80) ;
\draw [color={rgb, 255:red, 0; green, 0; blue, 0 }  ,draw opacity=1 ][line width=0.75]  [dash pattern={on 0.84pt off 2.51pt}]  (350,170) -- (370,190) ;
\draw [color={rgb, 255:red, 208; green, 2; blue, 27 }  ,draw opacity=1 ][line width=1.5]    (350,150) -- (370,170) ;
\draw [color={rgb, 255:red, 0; green, 0; blue, 0 }  ,draw opacity=1 ][line width=0.75]    (350,280) -- (370,300) ;
\draw [color={rgb, 255:red, 0; green, 0; blue, 0 }  ,draw opacity=1 ][line width=0.75]  [dash pattern={on 0.84pt off 2.51pt}]  (350,260) -- (370,280) ;
\draw [color={rgb, 255:red, 0; green, 0; blue, 0 }  ,draw opacity=1 ][line width=1.5]    (350,240) -- (370,260) ;
\draw [color={rgb, 255:red, 208; green, 2; blue, 27 }  ,draw opacity=1 ]   (210,210) -- (298,210) ;
\draw [shift={(300,210)}, rotate = 180] [color={rgb, 255:red, 208; green, 2; blue, 27 }  ,draw opacity=1 ][line width=0.75]    (10.93,-3.29) .. controls (6.95,-1.4) and (3.31,-0.3) .. (0,0) .. controls (3.31,0.3) and (6.95,1.4) .. (10.93,3.29)   ;
\draw [shift={(210,210)}, rotate = 0] [color={rgb, 255:red, 208; green, 2; blue, 27 }  ,draw opacity=1 ][line width=0.75]      (0,-11.18) .. controls (-3.09,-11.18) and (-5.59,-8.68) .. (-5.59,-5.59) .. controls (-5.59,-2.5) and (-3.09,0) .. (0,0) ;
\draw [color={rgb, 255:red, 0; green, 0; blue, 0 }  ,draw opacity=1 ]   (210,320) -- (298,320) ;
\draw [shift={(300,320)}, rotate = 180] [color={rgb, 255:red, 0; green, 0; blue, 0 }  ,draw opacity=1 ][line width=0.75]    (10.93,-3.29) .. controls (6.95,-1.4) and (3.31,-0.3) .. (0,0) .. controls (3.31,0.3) and (6.95,1.4) .. (10.93,3.29)   ;
\draw [shift={(210,320)}, rotate = 0] [color={rgb, 255:red, 0; green, 0; blue, 0 }  ,draw opacity=1 ][line width=0.75]      (0,-11.18) .. controls (-3.09,-11.18) and (-5.59,-8.68) .. (-5.59,-5.59) .. controls (-5.59,-2.5) and (-3.09,0) .. (0,0) ;
\draw    (520,150) -- (600,150) ;
\draw    (520,230) -- (600,230) ;
\draw    (600,150) -- (600,230) ;
\draw    (520,150) -- (520,230) ;
\draw    (520,230) -- (520,310) ;
\draw    (600,230) -- (600,310) ;
\draw    (520,310) -- (600,310) ;
\draw    (520,70) -- (600,70) ;
\draw    (600,150) -- (600,70) ;
\draw    (520,150) -- (520,70) ;
\draw    (400,252) .. controls (400.63,249.73) and (402.08,248.9) .. (404.35,249.53) .. controls (406.62,250.16) and (408.07,249.33) .. (408.7,247.06) .. controls (409.33,244.79) and (410.78,243.97) .. (413.05,244.6) .. controls (415.32,245.23) and (416.77,244.4) .. (417.39,242.13) .. controls (418.02,239.86) and (419.47,239.03) .. (421.74,239.66) .. controls (424.01,240.29) and (425.46,239.46) .. (426.09,237.19) .. controls (426.72,234.92) and (428.17,234.09) .. (430.44,234.72) .. controls (432.71,235.35) and (434.16,234.53) .. (434.79,232.26) .. controls (435.42,229.99) and (436.87,229.16) .. (439.14,229.79) .. controls (441.41,230.42) and (442.86,229.59) .. (443.48,227.32) .. controls (444.11,225.05) and (445.56,224.22) .. (447.83,224.85) .. controls (450.1,225.48) and (451.55,224.65) .. (452.18,222.38) .. controls (452.81,220.11) and (454.26,219.29) .. (456.53,219.92) .. controls (458.8,220.55) and (460.25,219.72) .. (460.88,217.45) .. controls (461.51,215.18) and (462.96,214.35) .. (465.23,214.98) -- (465.3,214.94) -- (472.26,210.99) ;
\draw [shift={(474,210)}, rotate = 510.42] [color={rgb, 255:red, 0; green, 0; blue, 0 }  ][line width=0.75]    (10.93,-3.29) .. controls (6.95,-1.4) and (3.31,-0.3) .. (0,0) .. controls (3.31,0.3) and (6.95,1.4) .. (10.93,3.29)   ;

\draw (245,87) node    {$Tile\ 1$};
\draw (117,100) node  [color={rgb, 255:red, 0; green, 0; blue, 0 }  ,opacity=1 ]  {$2$};
\draw (47,110) node  [color={rgb, 255:red, 0; green, 0; blue, 0 }  ,opacity=1 ]  {$6$};
\draw (47,50) node  [color={rgb, 255:red, 0; green, 0; blue, 0 }  ,opacity=1 ]  {$7$};
\draw (113,20) node  [color={rgb, 255:red, 0; green, 0; blue, 0 }  ,opacity=1 ]  {$4$};
\draw (73,180) node    {$1$};
\draw (163,180) node    {$3$};
\draw (113,150) node    {$4$};
\draw (117,230) node    {$5$};
\draw (113,290) node    {$2$};
\draw (183,270) node    {$9$};
\draw (113,360) node    {$5$};
\draw (183,330) node    {$8$};
\draw (360.03,175.75) node  [font=\scriptsize,color={rgb, 255:red, 208; green, 2; blue, 27 }  ,opacity=1 ]  {$4$};
\draw (361.14,143.06) node  [font=\scriptsize,color={rgb, 255:red, 208; green, 2; blue, 27 }  ,opacity=1 ]  {$5$};
\draw (376.2,160) node  [font=\scriptsize,color={rgb, 255:red, 208; green, 2; blue, 27 }  ,opacity=1 ]  {$3$};
\draw (346.17,159.14) node  [font=\scriptsize,color={rgb, 255:red, 208; green, 2; blue, 27 }  ,opacity=1 ]  {$1$};
\draw (375.95,69.1) node  [font=\scriptsize,color={rgb, 255:red, 0; green, 0; blue, 0 }  ,opacity=1 ]  {$2$};
\draw (361.46,85.84) node  [font=\scriptsize,color={rgb, 255:red, 0; green, 0; blue, 0 }  ,opacity=1 ]  {$6$};
\draw (346.5,69.1) node  [font=\scriptsize,color={rgb, 255:red, 0; green, 0; blue, 0 }  ,opacity=1 ]  {$7$};
\draw (361.86,53.14) node  [font=\scriptsize,color={rgb, 255:red, 0; green, 0; blue, 0 }  ,opacity=1 ]  {$4$};
\draw (361.3,265.94) node  [font=\scriptsize,color={rgb, 255:red, 0; green, 0; blue, 0 }  ,opacity=1 ]  {$5$};
\draw (375.83,249.51) node  [font=\scriptsize,color={rgb, 255:red, 0; green, 0; blue, 0 }  ,opacity=1 ]  {$8$};
\draw (360.94,233.83) node  [font=\scriptsize,color={rgb, 255:red, 0; green, 0; blue, 0 }  ,opacity=1 ]  {$9$};
\draw (345.94,249.8) node  [font=\scriptsize,color={rgb, 255:red, 0; green, 0; blue, 0 }  ,opacity=1 ]  {$2$};
\draw (87,70) node  [color={rgb, 255:red, 0; green, 0; blue, 0 }  ,opacity=1 ]  {$1$};
\draw (105,187) node    {$2$};
\draw (143,309) node  [color={rgb, 255:red, 65; green, 117; blue, 5 }  ,opacity=1 ]  {$\textcolor[rgb]{0,0,0}{3}$};
\draw (245,197) node    {$\textcolor[rgb]{0.82,0.01,0.11}{Tile\ 2}$};
\draw (245,307) node    {$Tile\ 3$};
\draw (558.33,318.67) node    {$6$};
\draw (513,272.67) node    {$7$};
\draw (608.33,271.33) node    {$2$};
\draw (559,237.33) node    {$4$};
\draw (514.67,189) node    {$1$};
\draw (607,188) node    {$3$};
\draw (559,158.67) node    {$5$};
\draw (514.33,108.67) node    {$2$};
\draw (607,108) node    {$8$};
\draw (559,60) node    {$9$};
\draw (461,260) node    {$ \begin{array}{l}
Remove\\
Diagonals
\end{array}$};

\end{tikzpicture}
  
\end {figure}

\end{ex}

\begin{rmk} \label{straight_zigzag_G}
    It follows easily from the constructions that a straight word $w$ yields a fan triangulation $\Delta_w$ (see Remark 3.16), which in turn results in a zigzag snake graph $G_w$. Conversely, a zigzag word $w$ yields a zigzag triangulation $\Delta_w$, which gives a straight snake graph $G_w$.
\end{rmk}

\begin{rmk}
    The definition we have given for building a snake graph from a triangulated surface is essentially a process called \textit{unfolding}. Conversely, any snake graph from a surface can be \textit{folded} back up to reconstruct the surface. Later, it will be convenient to use the folding/unfolding maps to relate certain expansions to others. For details see \cite{musiker2010cluster}. Figure \ref{fig:folding_unfolding} below illustrates the folding and unfolding maps.

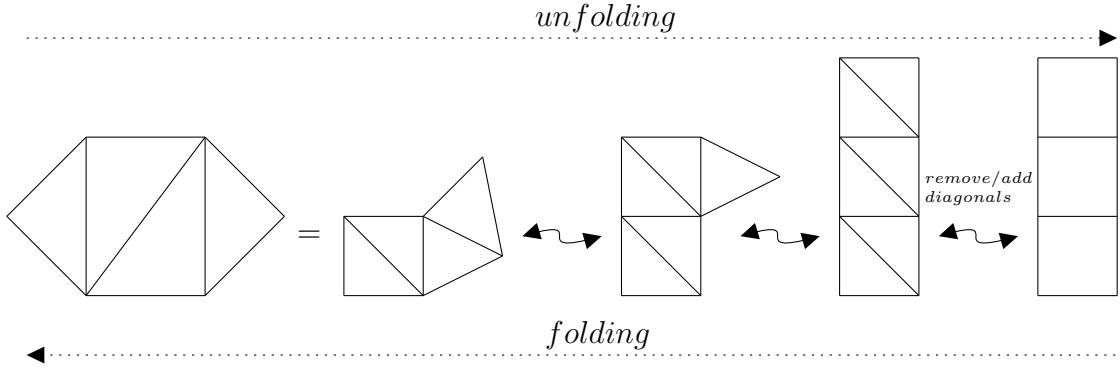
\begin {figure}[h!]
    \centering
    \caption{The folding and unfolding maps}
    \label{fig:folding_unfolding}
    \begin{tikzpicture}[x=0.75pt,y=0.75pt,yscale=-1,xscale=1]

\draw [color={rgb, 255:red, 0; green, 0; blue, 0 }  ,draw opacity=1 ]   (120,90) -- (160,130) ;
\draw [color={rgb, 255:red, 0; green, 0; blue, 0 }  ,draw opacity=1 ]   (60,90) -- (120,90) ;
\draw [color={rgb, 255:red, 0; green, 0; blue, 0 }  ,draw opacity=1 ]   (20,130) -- (60,90) ;
\draw [color={rgb, 255:red, 0; green, 0; blue, 0 }  ,draw opacity=1 ]   (20,130) -- (60,170) ;
\draw [color={rgb, 255:red, 0; green, 0; blue, 0 }  ,draw opacity=1 ]   (120,170) -- (160,130) ;
\draw [color={rgb, 255:red, 0; green, 0; blue, 0 }  ,draw opacity=1 ]   (60,170) -- (120,170) ;
\draw [color={rgb, 255:red, 0; green, 0; blue, 0 }  ,draw opacity=1 ]   (60,90) -- (60,170) ;
\draw [color={rgb, 255:red, 0; green, 0; blue, 0 }  ,draw opacity=1 ]   (60,170) -- (120,90) ;
\draw [color={rgb, 255:red, 0; green, 0; blue, 0 }  ,draw opacity=1 ]   (120,90) -- (120,170) ;
\draw [color={rgb, 255:red, 0; green, 0; blue, 0 }  ,draw opacity=1 ]   (190,170) -- (230,170) ;
\draw [color={rgb, 255:red, 0; green, 0; blue, 0 }  ,draw opacity=1 ]   (190,130) -- (190,170) ;
\draw [color={rgb, 255:red, 0; green, 0; blue, 0 }  ,draw opacity=1 ]   (190,130) -- (230,130) ;
\draw [color={rgb, 255:red, 0; green, 0; blue, 0 }  ,draw opacity=1 ]   (230,130) -- (230,170) ;
\draw [color={rgb, 255:red, 0; green, 0; blue, 0 }  ,draw opacity=1 ]   (190,130) -- (230,170) ;
\draw [color={rgb, 255:red, 0; green, 0; blue, 0 }  ,draw opacity=1 ]   (270,150) -- (230,170) ;
\draw [color={rgb, 255:red, 0; green, 0; blue, 0 }  ,draw opacity=1 ]   (230,130) -- (270,150) ;
\draw [color={rgb, 255:red, 0; green, 0; blue, 0 }  ,draw opacity=1 ]   (260,100) -- (270,150) ;
\draw [color={rgb, 255:red, 0; green, 0; blue, 0 }  ,draw opacity=1 ]   (260,100) -- (230,130) ;
\draw [color={rgb, 255:red, 0; green, 0; blue, 0 }  ,draw opacity=1 ]   (330,130) -- (330,170) ;
\draw [color={rgb, 255:red, 0; green, 0; blue, 0 }  ,draw opacity=1 ]   (283.27,139.05) .. controls (315.64,130.02) and (280.71,150.59) .. (317.65,140.65) ;
\draw [shift={(320,140)}, rotate = 524.31] [fill={rgb, 255:red, 0; green, 0; blue, 0 }  ,fill opacity=1 ][line width=0.08]  [draw opacity=0] (8.93,-4.29) -- (0,0) -- (8.93,4.29) -- cycle    ;
\draw [shift={(280,140)}, rotate = 343.39] [fill={rgb, 255:red, 0; green, 0; blue, 0 }  ,fill opacity=1 ][line width=0.08]  [draw opacity=0] (8.93,-4.29) -- (0,0) -- (8.93,4.29) -- cycle    ;
\draw [color={rgb, 255:red, 0; green, 0; blue, 0 }  ,draw opacity=1 ]   (393.27,139.05) .. controls (425.64,130.02) and (390.71,150.59) .. (427.65,140.65) ;
\draw [shift={(430,140)}, rotate = 524.31] [fill={rgb, 255:red, 0; green, 0; blue, 0 }  ,fill opacity=1 ][line width=0.08]  [draw opacity=0] (8.93,-4.29) -- (0,0) -- (8.93,4.29) -- cycle    ;
\draw [shift={(390,140)}, rotate = 343.39] [fill={rgb, 255:red, 0; green, 0; blue, 0 }  ,fill opacity=1 ][line width=0.08]  [draw opacity=0] (8.93,-4.29) -- (0,0) -- (8.93,4.29) -- cycle    ;
\draw [color={rgb, 255:red, 0; green, 0; blue, 0 }  ,draw opacity=1 ]   (330,130) -- (370,130) ;
\draw [color={rgb, 255:red, 0; green, 0; blue, 0 }  ,draw opacity=1 ]   (330,170) -- (370,170) ;
\draw [color={rgb, 255:red, 0; green, 0; blue, 0 }  ,draw opacity=1 ]   (370,130) -- (370,170) ;
\draw [color={rgb, 255:red, 0; green, 0; blue, 0 }  ,draw opacity=1 ]   (370,90) -- (370,130) ;
\draw [color={rgb, 255:red, 0; green, 0; blue, 0 }  ,draw opacity=1 ]   (330,90) -- (370,90) ;
\draw [color={rgb, 255:red, 0; green, 0; blue, 0 }  ,draw opacity=1 ]   (330,90) -- (330,130) ;
\draw [color={rgb, 255:red, 0; green, 0; blue, 0 }  ,draw opacity=1 ]   (330,90) -- (370,130) ;
\draw [color={rgb, 255:red, 0; green, 0; blue, 0 }  ,draw opacity=1 ]   (370,90) -- (410,110) ;
\draw [color={rgb, 255:red, 0; green, 0; blue, 0 }  ,draw opacity=1 ]   (370,130) -- (410,110) ;
\draw [color={rgb, 255:red, 0; green, 0; blue, 0 }  ,draw opacity=1 ]   (440,130) -- (440,170) ;
\draw [color={rgb, 255:red, 0; green, 0; blue, 0 }  ,draw opacity=1 ]   (440,130) -- (480,130) ;
\draw [color={rgb, 255:red, 0; green, 0; blue, 0 }  ,draw opacity=1 ]   (480,170) -- (480,130) ;
\draw [color={rgb, 255:red, 0; green, 0; blue, 0 }  ,draw opacity=1 ]   (440,130) -- (440,90) ;
\draw [color={rgb, 255:red, 0; green, 0; blue, 0 }  ,draw opacity=1 ]   (440,90) -- (480,90) ;
\draw [color={rgb, 255:red, 0; green, 0; blue, 0 }  ,draw opacity=1 ]   (480,130) -- (480,90) ;
\draw [color={rgb, 255:red, 0; green, 0; blue, 0 }  ,draw opacity=1 ]   (480,90) -- (480,50) ;
\draw [color={rgb, 255:red, 0; green, 0; blue, 0 }  ,draw opacity=1 ]   (440,90) -- (440,50) ;
\draw [color={rgb, 255:red, 0; green, 0; blue, 0 }  ,draw opacity=1 ]   (440,50) -- (480,50) ;

\draw [color={rgb, 255:red, 0; green, 0; blue, 0 }  ,draw opacity=1 ]   (440,170) -- (480,170) ;
\draw  [dash pattern={on 0.84pt off 2.51pt}]  (33,200) -- (580,200) ;
\draw [shift={(30,200)}, rotate = 0] [fill={rgb, 255:red, 0; green, 0; blue, 0 }  ][line width=0.08]  [draw opacity=0] (8.93,-4.29) -- (0,0) -- (8.93,4.29) -- cycle    ;
\draw  [dash pattern={on 0.84pt off 2.51pt}]  (30,40) -- (577,40) ;
\draw [shift={(580,40)}, rotate = 180] [fill={rgb, 255:red, 0; green, 0; blue, 0 }  ][line width=0.08]  [draw opacity=0] (8.93,-4.29) -- (0,0) -- (8.93,4.29) -- cycle    ;
\draw [color={rgb, 255:red, 0; green, 0; blue, 0 }  ,draw opacity=1 ]   (330,130) -- (370,170) ;
\draw [color={rgb, 255:red, 0; green, 0; blue, 0 }  ,draw opacity=1 ]   (440,130) -- (480,170) ;
\draw [color={rgb, 255:red, 0; green, 0; blue, 0 }  ,draw opacity=1 ]   (440,90) -- (480,130) ;
\draw [color={rgb, 255:red, 0; green, 0; blue, 0 }  ,draw opacity=1 ]   (440,50) -- (480,90) ;
\draw [color={rgb, 255:red, 0; green, 0; blue, 0 }  ,draw opacity=1 ]   (493.27,139.05) .. controls (525.64,130.02) and (490.71,150.59) .. (527.65,140.65) ;
\draw [shift={(530,140)}, rotate = 524.31] [fill={rgb, 255:red, 0; green, 0; blue, 0 }  ,fill opacity=1 ][line width=0.08]  [draw opacity=0] (8.93,-4.29) -- (0,0) -- (8.93,4.29) -- cycle    ;
\draw [shift={(490,140)}, rotate = 343.39] [fill={rgb, 255:red, 0; green, 0; blue, 0 }  ,fill opacity=1 ][line width=0.08]  [draw opacity=0] (8.93,-4.29) -- (0,0) -- (8.93,4.29) -- cycle    ;
\draw [color={rgb, 255:red, 0; green, 0; blue, 0 }  ,draw opacity=1 ]   (540,130) -- (540,170) ;
\draw [color={rgb, 255:red, 0; green, 0; blue, 0 }  ,draw opacity=1 ]   (540,130) -- (580,130) ;
\draw [color={rgb, 255:red, 0; green, 0; blue, 0 }  ,draw opacity=1 ]   (580,170) -- (580,130) ;
\draw [color={rgb, 255:red, 0; green, 0; blue, 0 }  ,draw opacity=1 ]   (540,130) -- (540,90) ;
\draw [color={rgb, 255:red, 0; green, 0; blue, 0 }  ,draw opacity=1 ]   (540,90) -- (580,90) ;
\draw [color={rgb, 255:red, 0; green, 0; blue, 0 }  ,draw opacity=1 ]   (580,130) -- (580,90) ;
\draw [color={rgb, 255:red, 0; green, 0; blue, 0 }  ,draw opacity=1 ]   (580,90) -- (580,50) ;
\draw [color={rgb, 255:red, 0; green, 0; blue, 0 }  ,draw opacity=1 ]   (540,90) -- (540,50) ;
\draw [color={rgb, 255:red, 0; green, 0; blue, 0 }  ,draw opacity=1 ]   (540,50) -- (580,50) ;

\draw [color={rgb, 255:red, 0; green, 0; blue, 0 }  ,draw opacity=1 ]   (540,170) -- (580,170) ;

\draw (172.5,140) node    {$=$};
\draw (50,142) node   [align=left] {};
\draw (318,190) node    {$folding$};
\draw (322,30) node    {$unfolding$};
\draw (510,116) node  [font=\tiny]  {$ \begin{array}{l}
remove/add\\
diagonals
\end{array}$};

\end{tikzpicture}
  
\end{figure}

\end{rmk}

\begin{defn}
    Fix a word $w$ and the snake graph $G_w$. Form a word $\text{sh}(G_w)$ of length $n-1$, called the \textit{shape of the snake graph $G_w$}, by letting $i$ run through $1,2, \dots ,n$ in the following rule: 
    If tile $T_{i+1}$ is glued to the right of tile $T_i$, write the letter $a$, and if tile $T_{i+1}$ is glued to the top of tile $T_i$, write the letter $b$
\end{defn}

\begin{ex}
    The shape of the snake graph for the word $w = ab$ is $\text{sh} (G_w) = bb$.
\end{ex}

\begin{rmk}
    If $w$ is straight then $G_w$ is zigzag (see Remark \ref{straight_zigzag_G}) and so $\text{sh}(G_w)$ is zigzag. Conversely, if $w$ is zigzag then $G_w$ is straight and hence $\text{sh}(G_w)$ is straight.
\end{rmk}

Definition \ref{sign_function} and Definition \ref{sign_sequence} below will be used in the next section to define the continued fraction associated to a word. 

Recall that $G_w$ has an embedding into $\mathbb{Z}^2$ such that the first tile $T_1$ of $G_w$ has vertices $(0,0),(1,0),(0,1)$, and $(1,1).$ Additionally, recall that each tile $T_i$ of $G_w$ is glued either above or to the right of the previous tile $T_{i-1}$. Informally, the snake graph $G_w$ ``goes up and to the right''.

\begin{defn} \label{sign_function} (see \cite{ccanakcci2018cluster})
    Fix the word $w$ and the associated snake graph $G_w$. Let $x$ and $y$ be the names of the coordinate functions on $\mathbb{Z}^2$. Note that the midpoint $m_e$ of each edge $e$ in $G_w$ lies on precisely one of the diagonal lines $y = x + (j+ \frac{1}{2})$ for $j \in \mathbb{Z}$. The \textit{sign function} $s$ on $G_w$ is the function on the edges $e$ of $G_w$ to the set $\{ - , + \}$ defined by
\[
  s (e) =
  \begin{cases}
                                   -, & \text{if $m_e$ lies on $y = x + (j+ \frac{1}{2})$ for $j$ even} \\
                                   +, & if \text{$m_e$ lies on $y = x + (j+ \frac{1}{2})$ for $j$ odd.}
  \end{cases}
\]
\end{defn}

\begin{ex} 
    The construction of the sign function $s$ on $G_{ab}$ is shown below in Figure \ref{fig:sign_ab}. 
\end{ex}  

\begin {figure}[h!]
    \centering
    \caption{The sign function $s$ on $G_{ab}$}
    \label{fig:sign_ab}
    \begin{tikzpicture}[x=0.75pt,y=0.75pt,yscale=-1,xscale=1]

\draw    (460,120) -- (460,160) ;

\draw    (460,120) -- (500,120) ;

\draw    (500,120) -- (500,160) ;

\draw    (460,160) -- (500,160) ;

\draw    (460,80) -- (460,120) ;

\draw    (460,80) -- (500,80) ;

\draw    (500,80) -- (500,120) ;

\draw    (460,160) -- (460,200) ;

\draw    (500,160) -- (500,200) ;

\draw    (460,200) -- (500,200) ;

\draw  [dash pattern={on 0.84pt off 2.51pt}]  (430,130) -- (510,50) ;

\draw [color={rgb, 255:red, 208; green, 2; blue, 27 }  ,draw opacity=1 ]   (470.5,195.5) -- (490.5,195.5) ;

\draw [color={rgb, 255:red, 208; green, 2; blue, 27 }  ,draw opacity=1 ]   (489.5,181.5) -- (509.5,181.5) ;

\draw [color={rgb, 255:red, 208; green, 2; blue, 27 }  ,draw opacity=1 ]   (469,156) -- (489,156) ;

\draw [color={rgb, 255:red, 208; green, 2; blue, 27 }  ,draw opacity=1 ]   (479,146) -- (479,166) ;

\draw [color={rgb, 255:red, 208; green, 2; blue, 27 }  ,draw opacity=1 ]   (447,179) -- (467,179) ;

\draw [color={rgb, 255:red, 208; green, 2; blue, 27 }  ,draw opacity=1 ]   (457,169) -- (457,189) ;

\draw [color={rgb, 255:red, 208; green, 2; blue, 27 }  ,draw opacity=1 ]   (493,140) -- (513,140) ;

\draw [color={rgb, 255:red, 208; green, 2; blue, 27 }  ,draw opacity=1 ]   (503,130) -- (503,150) ;

\draw [color={rgb, 255:red, 208; green, 2; blue, 27 }  ,draw opacity=1 ]   (451.5,140.5) -- (471.5,140.5) ;

\draw [color={rgb, 255:red, 208; green, 2; blue, 27 }  ,draw opacity=1 ]   (469.5,115.5) -- (489.5,115.5) ;

\draw [color={rgb, 255:red, 208; green, 2; blue, 27 }  ,draw opacity=1 ]   (489.5,99.5) -- (509.5,99.5) ;

\draw [color={rgb, 255:red, 208; green, 2; blue, 27 }  ,draw opacity=1 ]   (447,100) -- (467,100) ;

\draw [color={rgb, 255:red, 208; green, 2; blue, 27 }  ,draw opacity=1 ]   (457,90) -- (457,110) ;

\draw [color={rgb, 255:red, 208; green, 2; blue, 27 }  ,draw opacity=1 ]   (469,76) -- (489,76) ;

\draw [color={rgb, 255:red, 208; green, 2; blue, 27 }  ,draw opacity=1 ]   (479,66) -- (479,86) ;

\draw    (340,120) -- (340,160) ;

\draw    (340,120) -- (380,120) ;

\draw    (380,120) -- (380,160) ;

\draw    (340,160) -- (380,160) ;

\draw    (340,80) -- (340,120) ;

\draw    (340,80) -- (380,80) ;

\draw    (380,80) -- (380,120) ;

\draw    (340,160) -- (340,200) ;

\draw    (380,160) -- (380,200) ;

\draw    (340,200) -- (380,200) ;

\draw  [dash pattern={on 0.84pt off 2.51pt}]  (320,160) -- (400,80) ;

\draw [color={rgb, 255:red, 208; green, 2; blue, 27 }  ,draw opacity=1 ]   (350.5,195.5) -- (370.5,195.5) ;

\draw [color={rgb, 255:red, 208; green, 2; blue, 27 }  ,draw opacity=1 ]   (369.5,181.5) -- (389.5,181.5) ;

\draw [color={rgb, 255:red, 208; green, 2; blue, 27 }  ,draw opacity=1 ]   (349,156) -- (369,156) ;

\draw [color={rgb, 255:red, 208; green, 2; blue, 27 }  ,draw opacity=1 ]   (359,146) -- (359,166) ;

\draw [color={rgb, 255:red, 208; green, 2; blue, 27 }  ,draw opacity=1 ]   (327,179) -- (347,179) ;

\draw [color={rgb, 255:red, 208; green, 2; blue, 27 }  ,draw opacity=1 ]   (337,169) -- (337,189) ;

\draw [color={rgb, 255:red, 208; green, 2; blue, 27 }  ,draw opacity=1 ]   (373,140) -- (393,140) ;

\draw [color={rgb, 255:red, 208; green, 2; blue, 27 }  ,draw opacity=1 ]   (383,130) -- (383,150) ;

\draw [color={rgb, 255:red, 208; green, 2; blue, 27 }  ,draw opacity=1 ]   (331.5,140.5) -- (351.5,140.5) ;

\draw [color={rgb, 255:red, 208; green, 2; blue, 27 }  ,draw opacity=1 ]   (349.5,115.5) -- (369.5,115.5) ;

\draw [color={rgb, 255:red, 208; green, 2; blue, 27 }  ,draw opacity=1 ]   (369.5,99.5) -- (389.5,99.5) ;

\draw    (220,120) -- (220,160) ;

\draw    (220,120) -- (260,120) ;

\draw    (260,120) -- (260,160) ;

\draw    (220,160) -- (260,160) ;

\draw    (220,80) -- (220,120) ;

\draw    (220,80) -- (260,80) ;

\draw    (260,80) -- (260,120) ;

\draw    (220,160) -- (220,200) ;

\draw    (260,160) -- (260,200) ;

\draw    (220,200) -- (260,200) ;

\draw  [dash pattern={on 0.84pt off 2.51pt}]  (200,200) -- (280,120) ;

\draw [color={rgb, 255:red, 208; green, 2; blue, 27 }  ,draw opacity=1 ]   (230.5,195.5) -- (250.5,195.5) ;

\draw [color={rgb, 255:red, 208; green, 2; blue, 27 }  ,draw opacity=1 ]   (249.5,181.5) -- (269.5,181.5) ;

\draw [color={rgb, 255:red, 208; green, 2; blue, 27 }  ,draw opacity=1 ]   (229,156) -- (249,156) ;

\draw [color={rgb, 255:red, 208; green, 2; blue, 27 }  ,draw opacity=1 ]   (239,146) -- (239,166) ;

\draw [color={rgb, 255:red, 208; green, 2; blue, 27 }  ,draw opacity=1 ]   (207,179) -- (227,179) ;

\draw [color={rgb, 255:red, 208; green, 2; blue, 27 }  ,draw opacity=1 ]   (217,169) -- (217,189) ;

\draw [color={rgb, 255:red, 208; green, 2; blue, 27 }  ,draw opacity=1 ]   (253,140) -- (273,140) ;

\draw [color={rgb, 255:red, 208; green, 2; blue, 27 }  ,draw opacity=1 ]   (263,130) -- (263,150) ;

\draw    (100,120) -- (100,160) ;

\draw    (100,120) -- (140,120) ;

\draw    (140,120) -- (140,160) ;

\draw    (100,160) -- (140,160) ;

\draw    (100,80) -- (100,120) ;

\draw    (100,80) -- (140,80) ;

\draw    (140,80) -- (140,120) ;

\draw    (100,160) -- (100,200) ;

\draw    (140,160) -- (140,200) ;

\draw    (100,200) -- (140,200) ;

\draw  [dash pattern={on 0.84pt off 2.51pt}]  (90,230) -- (170,150) ;

\draw [color={rgb, 255:red, 208; green, 2; blue, 27 }  ,draw opacity=1 ]   (110.5,195.5) -- (130.5,195.5) ;

\draw [color={rgb, 255:red, 208; green, 2; blue, 27 }  ,draw opacity=1 ]   (129.5,181.5) -- (149.5,181.5) ;

\end{tikzpicture}
  
\end{figure}
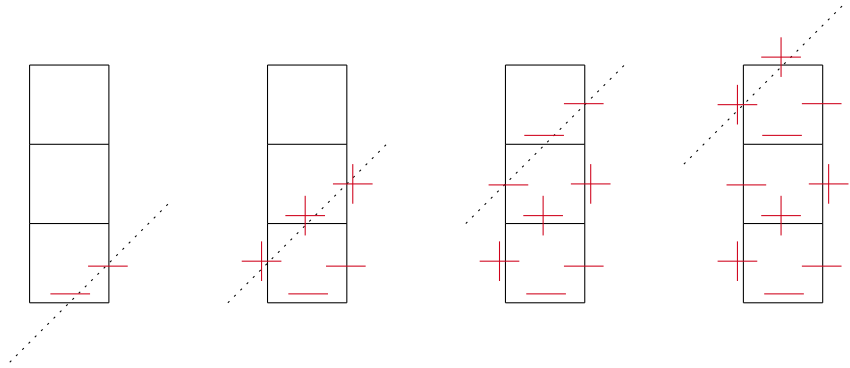

For $\epsilon \in \{ - , + \}$, define  

\[
  - \epsilon =
  \begin{cases}
                                   + & \text{if $\epsilon = -$} \\
                                   - & \text{if $\epsilon = +.$}
  \end{cases}
\]

Recall the internal edges of $G_w$ are labeled $e_1 , \dots , e_{n-1}.$ Let $e_0$ be the S edge of the first tile $T_1$ of $G_w$. Define $e_n$ to be either the N or E edge of tile $T_n$ according to whether the last three tiles of $G_w$ form a straight or zigzag snake graph; if the last three tiles are straight, the edge $e_n$ is across from $e_{n-1}$ (so that $s (e_n) = -s (e_{n-1})$), and if the last three tiles of $G_w$ instead form a zigzag snake graph then the edge $e_n$ is adjacent to $e_{n-1}$ (thus, $s (e_n) = s (e_{n-1})).$ If $\text{sh}(G_w) = bb,$ then we choose $e_2$ to be across from $e_{1}$. If $\text{sh}(G_w) = aa,$ then we choose $e_2$ to be adjacent to $e_{1}$. If $G_w$ has only one tile, then say $e_1$ is across from $e_0$.

\begin{defn} \label{sign_sequence} (see \cite{ccanakcci2018cluster})
    Fix the word $w$, the associated snake graph $G_w$, and the sign function $s$ on the edges of $G_w$. The \textit{sign sequence $\mathbf{s}_w$  associated to the word $w$} is the sequence $\mathbf{s}_w = (s(e_{0}) , s(e_{1}) , \dots , s(e_{n}))$.
\end{defn}
\begin{ex}
    For $w=ab$ the associated sign sequence is $\mathbf{s}_{ab} = (-,+,-,+).$

\begin {figure}[h!]
    \centering
    \caption{The sign sequence $\textbf{s}_{ab}$ on $G_{ab}$}
    \label{fig:s_ab}
      \begin{tikzpicture}[x=0.75pt,y=0.75pt,yscale=-1,xscale=1]

\draw    (460,120) -- (460,160) ;

\draw    (460,120) -- (500,120) ;

\draw    (500,120) -- (500,160) ;

\draw    (460,160) -- (500,160) ;

\draw    (460,80) -- (460,120) ;

\draw    (460,80) -- (500,80) ;

\draw    (500,80) -- (500,120) ;

\draw    (460,160) -- (460,200) ;

\draw    (500,160) -- (500,200) ;

\draw    (460,200) -- (500,200) ;

\draw [color={rgb, 255:red, 208; green, 2; blue, 27 }  ,draw opacity=1 ]   (470.5,195.5) -- (490.5,195.5) ;

\draw [color={rgb, 255:red, 208; green, 2; blue, 27 }  ,draw opacity=1 ]   (489.5,181.5) -- (509.5,181.5) ;

\draw [color={rgb, 255:red, 208; green, 2; blue, 27 }  ,draw opacity=1 ]   (469,156) -- (489,156) ;

\draw [color={rgb, 255:red, 208; green, 2; blue, 27 }  ,draw opacity=1 ]   (479,146) -- (479,166) ;

\draw [color={rgb, 255:red, 208; green, 2; blue, 27 }  ,draw opacity=1 ]   (447,179) -- (467,179) ;

\draw [color={rgb, 255:red, 208; green, 2; blue, 27 }  ,draw opacity=1 ]   (457,169) -- (457,189) ;

\draw [color={rgb, 255:red, 208; green, 2; blue, 27 }  ,draw opacity=1 ]   (493,140) -- (513,140) ;

\draw [color={rgb, 255:red, 208; green, 2; blue, 27 }  ,draw opacity=1 ]   (503,130) -- (503,150) ;

\draw [color={rgb, 255:red, 208; green, 2; blue, 27 }  ,draw opacity=1 ]   (451.5,140.5) -- (471.5,140.5) ;

\draw [color={rgb, 255:red, 208; green, 2; blue, 27 }  ,draw opacity=1 ]   (469.5,115.5) -- (489.5,115.5) ;

\draw [color={rgb, 255:red, 208; green, 2; blue, 27 }  ,draw opacity=1 ]   (489.5,99.5) -- (509.5,99.5) ;

\draw [color={rgb, 255:red, 208; green, 2; blue, 27 }  ,draw opacity=1 ]   (447,100) -- (467,100) ;

\draw [color={rgb, 255:red, 208; green, 2; blue, 27 }  ,draw opacity=1 ]   (457,90) -- (457,110) ;

\draw [color={rgb, 255:red, 208; green, 2; blue, 27 }  ,draw opacity=1 ]   (469,76) -- (489,76) ;

\draw [color={rgb, 255:red, 208; green, 2; blue, 27 }  ,draw opacity=1 ]   (479,66) -- (479,86) ;

\draw  [dash pattern={on 4.5pt off 4.5pt}] (464.6,196.8) .. controls (464.6,188.52) and (471.32,181.8) .. (479.6,181.8) .. controls (487.88,181.8) and (494.6,188.52) .. (494.6,196.8) .. controls (494.6,205.08) and (487.88,211.8) .. (479.6,211.8) .. controls (471.32,211.8) and (464.6,205.08) .. (464.6,196.8) -- cycle ;
\draw  [dash pattern={on 4.5pt off 4.5pt}] (463.8,156.4) .. controls (463.8,148.12) and (470.52,141.4) .. (478.8,141.4) .. controls (487.08,141.4) and (493.8,148.12) .. (493.8,156.4) .. controls (493.8,164.68) and (487.08,171.4) .. (478.8,171.4) .. controls (470.52,171.4) and (463.8,164.68) .. (463.8,156.4) -- cycle ;
\draw  [dash pattern={on 4.5pt off 4.5pt}] (464.4,116.2) .. controls (464.4,107.92) and (471.12,101.2) .. (479.4,101.2) .. controls (487.68,101.2) and (494.4,107.92) .. (494.4,116.2) .. controls (494.4,124.48) and (487.68,131.2) .. (479.4,131.2) .. controls (471.12,131.2) and (464.4,124.48) .. (464.4,116.2) -- cycle ;
\draw  [dash pattern={on 4.5pt off 4.5pt}] (463.8,76.4) .. controls (463.8,68.12) and (470.52,61.4) .. (478.8,61.4) .. controls (487.08,61.4) and (493.8,68.12) .. (493.8,76.4) .. controls (493.8,84.68) and (487.08,91.4) .. (478.8,91.4) .. controls (470.52,91.4) and (463.8,84.68) .. (463.8,76.4) -- cycle ;

\end{tikzpicture}
  
\end{figure}

\end{ex}

\section{Continued Fractions}

In \cite{ccanakcci2018cluster} (see also \cite{ccanakcci2020snake}, \cite{lee2019cluster}), connections between cluster algebras, continued fractions, and snake graphs were established. We review some of the basic definitions found there, and use them to define the continued fraction associated to a word $w$.

\begin{defn}
    Fix $a_1 \in \mathbb{Z}$ and for $2 \leq i \leq k$ fix the positive integers $a_i \in \mathbb{Z}.$ A \textit{finite continued fraction}, denoted by $[a_1 , a_2 , \dots , a_k]$ is an expression of the form 

\[
  a_1+\cfrac{1}{a_2+\cfrac{1}{a_3+\cfrac{1}{\ddots +\cfrac{1}{a_k}}}}
\]
\end{defn}

Say that the continued fraction $[a_1 , a_2 , \dots , a_k ]$ is \textit{positive} if $a_i > 0$ for every $i$. From now on we will only consider positive continued fractions.

\begin{defn}
    Fix a word $w$ of length $n-1$, along with the snake graph $G_w$. Recall the sign sequence $\mathbf{s}_w = (s(e_0) , \dots  , s(e_n))$ of $G_w$, and our convention that $s(e_0) = -.$ Let $\epsilon = -$. Define positive integers $a_1 , ... , a_k$ as indicated below.
    $$\mathbf{s}_w = (s(e_0) , \dots  , s(e_n)) = (\underbrace{\epsilon , \epsilon , \dots , \epsilon}_{a_1 \text{ times}} , \underbrace{-\epsilon , -\epsilon , \dots , -\epsilon}_{a_2 \text{ times}} , \dots , \underbrace{(-1)^k \epsilon , (-1)^k \epsilon , \dots (-1)^k \epsilon}_{a_k \text{ times}}).$$ Define the \textit{(finite, positive) continued fraction $\text{CF}(w)$ associated to $w$} to be $\text{CF}(w) = [a_1 , a_2 , \dots , a_k].$

\end{defn}

\begin{rmk}
     By Theorem A in \cite{ccanakcci2018cluster} (which we recall as Theorem \ref{CF_card_quotient_P} below), simplifying the continued fraction $\text{CF}(w)$ gives a rational number in lowest terms that is greater than $1$.
\end{rmk}

\begin{ex} \label{CF_ab}
    The continued fraction  $\text{CF}(ab)$ is equal to $\frac{5}{3}.$ Indeed, the associated sign sequence is $\mathbf{s}_{ab} = (-,+,-,+),$ so that 
    $$\text{CF}(ab) = [1,1,1,1] = 1 + \cfrac{1}{1+\cfrac{1}{1+\cfrac{1}{1}}} = \frac{5}{3}.$$ 
\end{ex}

\begin{rmk} \label{fib}
     Suppose the word $w$ is either straight or zigzag, and that its length is $l(w) = n-1$. Let the Fibonacci numbers be denoted by $F_1 = 1 , F_2 = 1 , F_3 = 2,$ etc. Consider the continued fraction $\text{CF}(w)$.

\begin{enumerate}[(a)]
    \item If $w$ is the straight word $w = a^{n-1}$ then $$\text{CF}(w) = [2,\underbrace{1,\dots ,1}_{n-1 \text{ times}}] = \frac{F_{n+2}}{F_n},$$ .
    \item If $w$ is the straight word $w = b^{n-1}$ then $$\text{CF}(w) = [\underbrace{1,1,1,\dots,1}_{n+1 \text{ times}}] = \frac{F_{n+2}}{F_{n+1}}.$$
    \item If $w$ is the zigzag word $w = bab \cdots$ then $\text{CF}(w) = [1,n] = \frac{n+1}{n}.$
    
    \item If $w$ is the zigzag word $w = aba \cdots$ then $\text{CF}(w) = [n] = \frac{n}{1} = n.$
\end{enumerate}

\end{rmk}

\section{Distributive Lattices}

\begin{defn}
    Let $D$ be a finite poset. The \textit{meet} of the elements $p \in D$ and $q \in D$ is the unique element denoted $p \wedge q \in D$ (if it exists) that satisfies: 
\begin{enumerate}[(a)]
    \item $p \wedge q < p$ and $p \wedge q < q$, and 
    
    \item if there exists some $r \in D$ such that $r < p$ and $r < q$, then $r < p \wedge q.$
\end{enumerate}

The \textit{join} of $p$ and $q$ is the unique element $p \vee q \in D$ (if it exists) that satisfies:
\begin{enumerate}[(a)]
    \item $p \vee q > p$ and $p \vee q > q$, and 
    
    \item if there exists some $r \in D$ such that $r > p$ and $r > q$, then $r > p \vee q.$
\end{enumerate}

We say $D$ is a \textit{lattice} if for any two elements $p , q \in D,$ both $p \wedge q$ and $p \vee q$ exist.

\end{defn}

\begin{rmk} \label{dist_min_max}
    It is easy to see that every finite lattice has a unique minimum element and a unique maximum element.
\end{rmk}

\begin{defn}
    Suppose that the finite poset $D$ is a lattice. We say $D$ is a \textit{distributive lattice} if for all $p,q,r \in D$ the following two distributive laws hold:

\begin{enumerate}[(a)]
    \item $p \wedge (q \vee r) = (p \wedge q) \vee (p \wedge r)$
    
    \item $p \vee (q \wedge r) = (p \vee q) \wedge (p \vee r).$
\end{enumerate}
\end{defn}

\begin{defn}
    Let $C$ be a finite poset. An \textit{order ideal} $I$ of $C$ is a subset of $C$ such that if $p \in I$ and $q \in C$ with $q < p$, then $q \in C$. We denote by $\mathcal{I}(C)$ the poset (ordered by inclusion) of order ideals of $C$.
\end{defn}

\begin{ex} \label{fib_cube_000}
    The poset of order ideals of a fence on $n$ vertices is called a $\textit{Fibonacci cube of order $n$},$ which we denote $\Gamma_n$. Alternatively, the Fibonacci cube of order $n$ may be defined as a graph, with vertices those binary words from $\{ 0,1 \}$ with $n$ bits that do not contain two consecutive instances of the bit $1$. There is an edge between two vertices if they differ by precisely one bit in the same position. For original references, see \cite{hoft1985fibonacci} and \cite{hsu1993fibonacci}. For a somewhat recent survey on Fibonacci cubes, see \cite{klavvzar2013structure}. Figure \ref{fig:fib_cubes} below shows the Fibonacci cubes which result from computing the poset of order ideals on each of the four zigzag posets shown in Figure \ref{fig:fence}.
\end{ex}

\begin {figure}[h!]
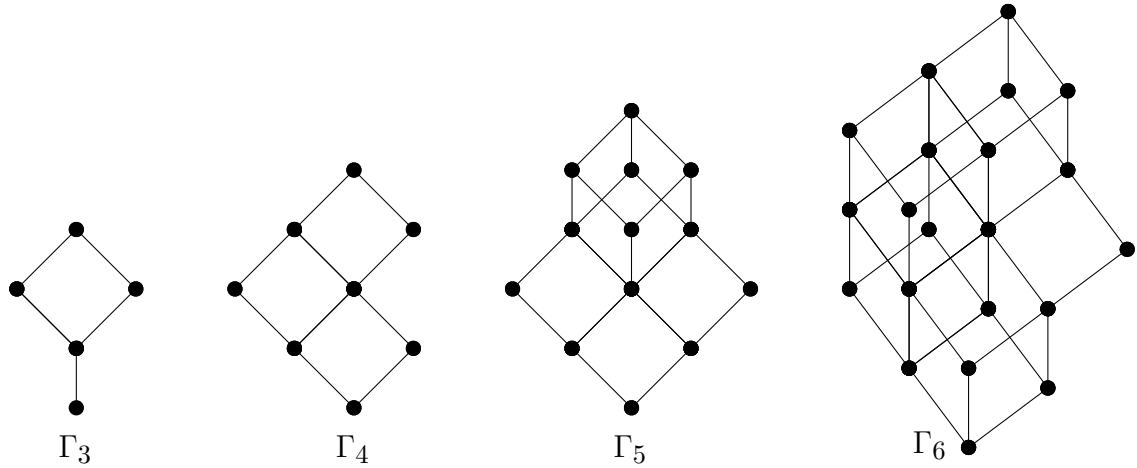

    \centering
    \caption{Four Fibonacci cubes}
    \label{fig:fib_cubes}


\end {figure}

\begin{defn}
    Let $D$ be a finite lattice. An element $r \in D$ is said to be \textit{join-irreducible} if $r$ is not the unique minimum element of $D$ (see Remark \ref{dist_min_max}) and there do not exist two elements $p,q \in D$ with $p < r , q < r,$ and $r = p \vee q$. We denote by $\mathcal{J}(D)$ the poset (with the induced order) of join-irreducible elements of $D$.
\end{defn}

\begin{thm} \label{birkhoff}
    (see \cite{birkhoff1937rings}) Let $D$ be a finite lattice. Let $C = \mathcal{J} (D)$ be the poset of join-irreducibles of $D$. Then $D$ is a distributive lattice if and only if $D$ is isomorphic to $\mathcal{I} (C).$
\end{thm}

By Theorem \ref{birkhoff}, $\mathcal{I}(C_w)$ is a distributive lattice for any word $w$. 

\begin{defn}
    The \textit{distributive lattice $D_w$ associated to $w$} is defined by $D_w = \mathcal{I}(C_w)$
\end{defn}

\begin{ex}
    The distributive lattice $D_{ab}$ is isomorphic to the Fibonacci cube $\Gamma_3$, shown as the leftmost poset in Figure \ref{fig:fib_cubes}.
\end{ex}

\chapter{Expansion Posets} \label{exp_poset_chapter}

Fix a word $w$ of length $l(w) = n-1$ and the associated objects defined in the previous chapter. The goals of this chapter are as follows:

\begin{enumerate}[(1)]
    \item Recall three known combinatorial interpretations of the terms in the Laurent expansion of any cluster variable $x_{w}$ parameterized by the arc $\gamma_w$.
    
    \item Equip each such representation with a poset structure. 
\end{enumerate}

\section{Perfect Matchings of Snake Graphs}

It is easy to see that any snake graph with $n$ tiles has an even number $2(n+1)$ of vertices.

\begin {defn}
   A \textit{perfect matching} $P$ of the snake graph $G_w$ is a choice of $n+1$ edges in $G_w$ such that each vertex of $G_w$ is the endpoint of exactly one edge $e$ in $P$.
\end {defn}

The \textit{weight} of the edge $e$ is the initial cluster variable $x_e$. The \textit{weight} of a perfect matching $P$ is defined to be the product of initial cluster variables $x_P = \prod_{e \in P} x_e.$ Let $\mathbb{P}_w$ be the set of all perfect matchings of the snake graph $G_w.$

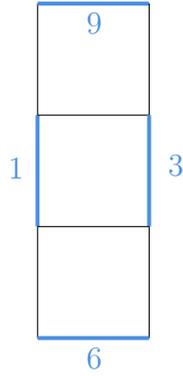
\begin {figure}[h!]
    \centering
    \caption{The perfect matching $P_{-}$ on $G_{ab}$}
    \label{fig:P_}
    \begin{tikzpicture}[x=0.75pt,y=0.75pt,yscale=-1,xscale=1]

\draw [color={rgb, 255:red, 74; green, 144; blue, 226 }  ,draw opacity=1 ][line width=1.5]    (306.98,157.38) -- (306.98,213.66) ;
\draw    (306.98,157.38) -- (363.25,157.38) ;
\draw [color={rgb, 255:red, 74; green, 144; blue, 226 }  ,draw opacity=1 ][line width=1.5]    (363.25,157.38) -- (363.25,213.66) ;
\draw    (306.98,213.66) -- (363.25,213.66) ;
\draw    (306.98,101.11) -- (306.98,157.38) ;
\draw [color={rgb, 255:red, 74; green, 144; blue, 226 }  ,draw opacity=1 ][line width=1.5]    (306.98,101.11) -- (363.25,101.11) ;
\draw    (363.25,101.11) -- (363.25,157.38) ;
\draw    (306.98,213.66) -- (306.98,269.93) ;
\draw    (363.25,213.66) -- (363.25,269.93) ;
\draw [color={rgb, 255:red, 74; green, 144; blue, 226 }  ,draw opacity=1 ][line width=1.5]    (306.98,269.93) -- (363.25,269.93) ;

\draw (335.65,280.17) node    {$\textcolor[rgb]{0.29,0.56,0.89}{6}$};
\draw (296.08,184.46) node    {$\textcolor[rgb]{0.29,0.56,0.89}{1}$};
\draw (376.8,182.88) node    {$\textcolor[rgb]{0.29,0.56,0.89}{3}$};
\draw (335.65,110.83) node    {$\textcolor[rgb]{0.29,0.56,0.89}{9}$};

\end{tikzpicture}
  
\end {figure}

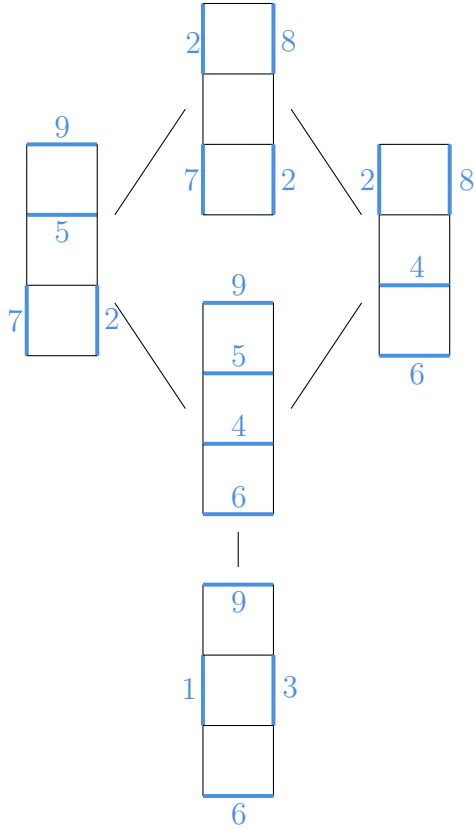
\begin {figure}[h!]
    \centering
    \caption{The poset $\mathbb{P}_{ab}$}
    \label{fig:P_ab}
    \begin{tikzpicture}[x=0.75pt,y=0.75pt,yscale=-1,xscale=1]

\draw [color={rgb, 255:red, 74; green, 144; blue, 226 }  ,draw opacity=1 ][line width=1.5]    (278.89,338.89) -- (278.89,374.44) ;
\draw    (278.89,338.89) -- (314.44,338.89) ;
\draw [color={rgb, 255:red, 74; green, 144; blue, 226 }  ,draw opacity=1 ][line width=1.5]    (314.44,338.89) -- (314.44,374.44) ;
\draw    (278.89,374.44) -- (314.44,374.44) ;
\draw    (278.89,303.33) -- (278.89,338.89) ;
\draw [color={rgb, 255:red, 74; green, 144; blue, 226 }  ,draw opacity=1 ][line width=1.5]    (278.89,303.33) -- (314.44,303.33) ;
\draw    (314.44,303.33) -- (314.44,338.89) ;
\draw    (278.89,374.44) -- (278.89,410) ;
\draw    (314.44,374.44) -- (314.44,410) ;
\draw [color={rgb, 255:red, 74; green, 144; blue, 226 }  ,draw opacity=1 ][line width=1.5]    (278.89,410) -- (314.44,410) ;
\draw    (278.89,196.67) -- (278.89,232.22) ;
\draw [color={rgb, 255:red, 74; green, 144; blue, 226 }  ,draw opacity=1 ][line width=1.5]    (278.89,196.67) -- (314.44,196.67) ;
\draw    (314.44,196.67) -- (314.44,232.22) ;
\draw [color={rgb, 255:red, 74; green, 144; blue, 226 }  ,draw opacity=1 ][line width=1.5]    (278.89,232.22) -- (314.44,232.22) ;
\draw    (278.89,161.11) -- (278.89,196.67) ;
\draw [color={rgb, 255:red, 74; green, 144; blue, 226 }  ,draw opacity=1 ][line width=1.5]    (278.89,161.11) -- (314.44,161.11) ;
\draw    (314.44,161.11) -- (314.44,196.67) ;
\draw    (278.89,232.22) -- (278.89,267.78) ;
\draw    (314.44,232.22) -- (314.44,267.78) ;
\draw [color={rgb, 255:red, 74; green, 144; blue, 226 }  ,draw opacity=1 ][line width=1.5]    (278.89,267.78) -- (314.44,267.78) ;
\draw    (367.78,116.67) -- (367.78,152.22) ;
\draw    (367.78,116.67) -- (403.33,116.67) ;
\draw    (403.33,116.67) -- (403.33,152.22) ;
\draw [color={rgb, 255:red, 74; green, 144; blue, 226 }  ,draw opacity=1 ][line width=1.5]    (367.78,152.22) -- (403.33,152.22) ;
\draw [color={rgb, 255:red, 74; green, 144; blue, 226 }  ,draw opacity=1 ][line width=1.5]    (367.78,81.11) -- (367.78,116.67) ;
\draw    (367.78,81.11) -- (403.33,81.11) ;
\draw [color={rgb, 255:red, 74; green, 144; blue, 226 }  ,draw opacity=1 ][line width=1.5]    (403.33,81.11) -- (403.33,116.67) ;
\draw    (367.78,152.22) -- (367.78,187.78) ;
\draw    (403.33,152.22) -- (403.33,187.78) ;
\draw [color={rgb, 255:red, 74; green, 144; blue, 226 }  ,draw opacity=1 ][line width=1.5]    (367.78,187.78) -- (403.33,187.78) ;
\draw    (190,116.67) -- (190,152.22) ;
\draw [color={rgb, 255:red, 74; green, 144; blue, 226 }  ,draw opacity=1 ][line width=1.5]    (190,116.67) -- (225.56,116.67) ;
\draw    (225.56,116.67) -- (225.56,152.22) ;
\draw    (190,152.22) -- (225.56,152.22) ;
\draw    (190,81.11) -- (190,116.67) ;
\draw [color={rgb, 255:red, 74; green, 144; blue, 226 }  ,draw opacity=1 ][line width=1.5]    (190,81.11) -- (225.56,81.11) ;
\draw    (225.56,81.11) -- (225.56,116.67) ;
\draw [color={rgb, 255:red, 74; green, 144; blue, 226 }  ,draw opacity=1 ][line width=1.5]    (190,152.22) -- (190,187.78) ;
\draw [color={rgb, 255:red, 74; green, 144; blue, 226 }  ,draw opacity=1 ][line width=1.5]    (225.56,152.22) -- (225.56,187.78) ;
\draw    (190,187.78) -- (225.56,187.78) ;
\draw    (278.89,45.56) -- (278.89,81.11) ;
\draw    (278.89,45.56) -- (314.44,45.56) ;
\draw    (314.44,45.56) -- (314.44,81.11) ;
\draw    (278.89,81.11) -- (314.44,81.11) ;
\draw [color={rgb, 255:red, 74; green, 144; blue, 226 }  ,draw opacity=1 ][line width=1.5]    (278.89,10) -- (278.89,45.56) ;
\draw    (278.89,10) -- (314.44,10) ;
\draw [color={rgb, 255:red, 74; green, 144; blue, 226 }  ,draw opacity=1 ][line width=1.5]    (314.44,10) -- (314.44,45.56) ;
\draw [color={rgb, 255:red, 74; green, 144; blue, 226 }  ,draw opacity=1 ][line width=1.5]    (278.89,81.11) -- (278.89,116.67) ;
\draw [color={rgb, 255:red, 74; green, 144; blue, 226 }  ,draw opacity=1 ][line width=1.5]    (314.44,81.11) -- (314.44,116.67) ;
\draw    (278.89,116.67) -- (314.44,116.67) ;
\draw    (296.67,276.67) -- (296.67,294.44) ;
\draw    (234.44,161.11) -- (270,214.44) ;
\draw    (323.33,214.44) -- (358.89,161.11) ;
\draw    (234.44,116.67) -- (270,63.33) ;
\draw    (323.33,63.33) -- (358.89,116.67) ;

\draw (297,419) node    {$\textcolor[rgb]{0.29,0.56,0.89}{6}$};
\draw (272,356) node    {$\textcolor[rgb]{0.29,0.56,0.89}{1}$};
\draw (323,355) node    {$\textcolor[rgb]{0.29,0.56,0.89}{3}$};
\draw (297,312) node    {$\textcolor[rgb]{0.29,0.56,0.89}{9}$};
\draw (297,152) node    {$\textcolor[rgb]{0.29,0.56,0.89}{9}$};
\draw (208,72) node    {$\textcolor[rgb]{0.29,0.56,0.89}{9}$};
\draw (297,259) node    {$\textcolor[rgb]{0.29,0.56,0.89}{6}$};
\draw (387,197) node    {$\textcolor[rgb]{0.29,0.56,0.89}{6}$};
\draw (297,223) node    {$\textcolor[rgb]{0.29,0.56,0.89}{4}$};
\draw (387,143) node    {$\textcolor[rgb]{0.29,0.56,0.89}{4}$};
\draw (297,188) node    {$\textcolor[rgb]{0.29,0.56,0.89}{5}$};
\draw (208,125) node    {$\textcolor[rgb]{0.29,0.56,0.89}{5}$};
\draw (184,170) node    {$\textcolor[rgb]{0.29,0.56,0.89}{7}$};
\draw (233,169) node    {$\textcolor[rgb]{0.29,0.56,0.89}{2}$};
\draw (273,99) node    {$\textcolor[rgb]{0.29,0.56,0.89}{7}$};
\draw (322,99) node    {$\textcolor[rgb]{0.29,0.56,0.89}{2}$};
\draw (274,30) node    {$\textcolor[rgb]{0.29,0.56,0.89}{2}$};
\draw (322,29) node    {$\textcolor[rgb]{0.29,0.56,0.89}{8}$};
\draw (412,99) node    {$\textcolor[rgb]{0.29,0.56,0.89}{8}$};
\draw (362,99) node    {$\textcolor[rgb]{0.29,0.56,0.89}{2}$};

\end{tikzpicture}
  
\end {figure}

\begin {ex}
    Figure \ref{fig:P_} shows one perfect matching $P_{-}$ (see Definition \ref{min_P_def} and Figure \ref{fig:P_ab} below for the notation $P_{-}$) of the snake graph $G_{ab}$. The weight of this perfect matching is $x_1 x_3 x_6 x_9$. 
\end {ex}

\begin{thm}
    (Theorem 3.1 in \cite{musiker2010cluster}) Let $w$ be any word, and consider the set $\mathbb{P}_w$ of perfect matchings on the snake graph $G_w$. Then the cluster variable $x_{\gamma}$ can be written as the sum $$x_{\gamma} = \frac{1}{x_1 x_2 \dots x_n}  \sum_{P \in \mathbb{P}_w} x_P.$$ 
\end{thm}

\begin{ex}
   Figure \ref{fig:justpm} below shows the five perfect matchings on the snake graph $G_{ab}$. Note that summing over the weights of these perfect matchings gives the cluster variable $x_{ab}$ displayed in Figure \ref{fig:gamma_ab}.
\end{ex}

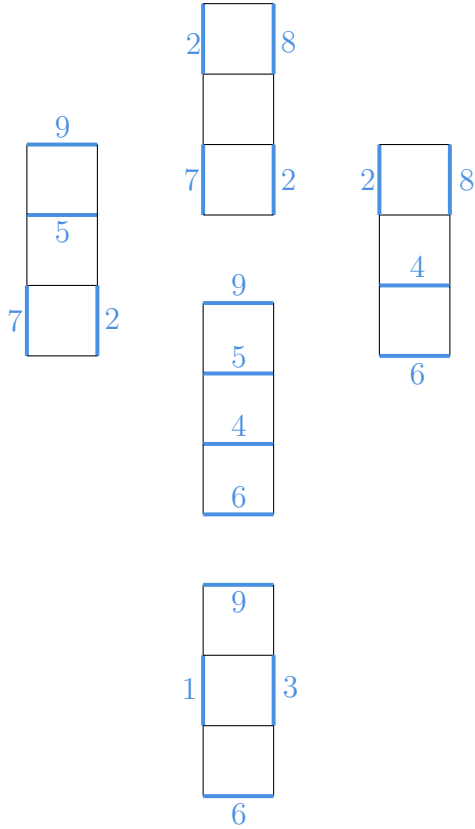
\begin {figure}[h!]
    \centering
    \caption{The five perfect matchings on $G_{ab}$}
    \label{fig:justpm}
    \begin{tikzpicture}[x=0.75pt,y=0.75pt,yscale=-1,xscale=1]

\draw [color={rgb, 255:red, 74; green, 144; blue, 226 }  ,draw opacity=1 ][line width=1.5]    (298.89,358.89) -- (298.89,394.44) ;
\draw    (298.89,358.89) -- (334.44,358.89) ;
\draw [color={rgb, 255:red, 74; green, 144; blue, 226 }  ,draw opacity=1 ][line width=1.5]    (334.44,358.89) -- (334.44,394.44) ;
\draw    (298.89,394.44) -- (334.44,394.44) ;
\draw    (298.89,323.33) -- (298.89,358.89) ;
\draw [color={rgb, 255:red, 74; green, 144; blue, 226 }  ,draw opacity=1 ][line width=1.5]    (298.89,323.33) -- (334.44,323.33) ;
\draw    (334.44,323.33) -- (334.44,358.89) ;
\draw    (298.89,394.44) -- (298.89,430) ;
\draw    (334.44,394.44) -- (334.44,430) ;
\draw [color={rgb, 255:red, 74; green, 144; blue, 226 }  ,draw opacity=1 ][line width=1.5]    (298.89,430) -- (334.44,430) ;
\draw    (298.89,216.67) -- (298.89,252.22) ;
\draw [color={rgb, 255:red, 74; green, 144; blue, 226 }  ,draw opacity=1 ][line width=1.5]    (298.89,216.67) -- (334.44,216.67) ;
\draw    (334.44,216.67) -- (334.44,252.22) ;
\draw [color={rgb, 255:red, 74; green, 144; blue, 226 }  ,draw opacity=1 ][line width=1.5]    (298.89,252.22) -- (334.44,252.22) ;
\draw    (298.89,181.11) -- (298.89,216.67) ;
\draw [color={rgb, 255:red, 74; green, 144; blue, 226 }  ,draw opacity=1 ][line width=1.5]    (298.89,181.11) -- (334.44,181.11) ;
\draw    (334.44,181.11) -- (334.44,216.67) ;
\draw    (298.89,252.22) -- (298.89,287.78) ;
\draw    (334.44,252.22) -- (334.44,287.78) ;
\draw [color={rgb, 255:red, 74; green, 144; blue, 226 }  ,draw opacity=1 ][line width=1.5]    (298.89,287.78) -- (334.44,287.78) ;
\draw    (387.78,136.67) -- (387.78,172.22) ;
\draw    (387.78,136.67) -- (423.33,136.67) ;
\draw    (423.33,136.67) -- (423.33,172.22) ;
\draw [color={rgb, 255:red, 74; green, 144; blue, 226 }  ,draw opacity=1 ][line width=1.5]    (387.78,172.22) -- (423.33,172.22) ;
\draw [color={rgb, 255:red, 74; green, 144; blue, 226 }  ,draw opacity=1 ][line width=1.5]    (387.78,101.11) -- (387.78,136.67) ;
\draw    (387.78,101.11) -- (423.33,101.11) ;
\draw [color={rgb, 255:red, 74; green, 144; blue, 226 }  ,draw opacity=1 ][line width=1.5]    (423.33,101.11) -- (423.33,136.67) ;
\draw    (387.78,172.22) -- (387.78,207.78) ;
\draw    (423.33,172.22) -- (423.33,207.78) ;
\draw [color={rgb, 255:red, 74; green, 144; blue, 226 }  ,draw opacity=1 ][line width=1.5]    (387.78,207.78) -- (423.33,207.78) ;
\draw    (210,136.67) -- (210,172.22) ;
\draw [color={rgb, 255:red, 74; green, 144; blue, 226 }  ,draw opacity=1 ][line width=1.5]    (210,136.67) -- (245.56,136.67) ;
\draw    (245.56,136.67) -- (245.56,172.22) ;
\draw    (210,172.22) -- (245.56,172.22) ;
\draw    (210,101.11) -- (210,136.67) ;
\draw [color={rgb, 255:red, 74; green, 144; blue, 226 }  ,draw opacity=1 ][line width=1.5]    (210,101.11) -- (245.56,101.11) ;
\draw    (245.56,101.11) -- (245.56,136.67) ;
\draw [color={rgb, 255:red, 74; green, 144; blue, 226 }  ,draw opacity=1 ][line width=1.5]    (210,172.22) -- (210,207.78) ;
\draw [color={rgb, 255:red, 74; green, 144; blue, 226 }  ,draw opacity=1 ][line width=1.5]    (245.56,172.22) -- (245.56,207.78) ;
\draw    (210,207.78) -- (245.56,207.78) ;
\draw    (298.89,65.56) -- (298.89,101.11) ;
\draw    (298.89,65.56) -- (334.44,65.56) ;
\draw    (334.44,65.56) -- (334.44,101.11) ;
\draw    (298.89,101.11) -- (334.44,101.11) ;
\draw [color={rgb, 255:red, 74; green, 144; blue, 226 }  ,draw opacity=1 ][line width=1.5]    (298.89,30) -- (298.89,65.56) ;
\draw    (298.89,30) -- (334.44,30) ;
\draw [color={rgb, 255:red, 74; green, 144; blue, 226 }  ,draw opacity=1 ][line width=1.5]    (334.44,30) -- (334.44,65.56) ;
\draw [color={rgb, 255:red, 74; green, 144; blue, 226 }  ,draw opacity=1 ][line width=1.5]    (298.89,101.11) -- (298.89,136.67) ;
\draw [color={rgb, 255:red, 74; green, 144; blue, 226 }  ,draw opacity=1 ][line width=1.5]    (334.44,101.11) -- (334.44,136.67) ;
\draw    (298.89,136.67) -- (334.44,136.67) ;

\draw (317,439) node    {$\textcolor[rgb]{0.29,0.56,0.89}{6}$};
\draw (292,376) node    {$\textcolor[rgb]{0.29,0.56,0.89}{1}$};
\draw (343,375) node    {$\textcolor[rgb]{0.29,0.56,0.89}{3}$};
\draw (317,332) node    {$\textcolor[rgb]{0.29,0.56,0.89}{9}$};
\draw (317,172) node    {$\textcolor[rgb]{0.29,0.56,0.89}{9}$};
\draw (228,92) node    {$\textcolor[rgb]{0.29,0.56,0.89}{9}$};
\draw (317,279) node    {$\textcolor[rgb]{0.29,0.56,0.89}{6}$};
\draw (407,217) node    {$\textcolor[rgb]{0.29,0.56,0.89}{6}$};
\draw (317,243) node    {$\textcolor[rgb]{0.29,0.56,0.89}{4}$};
\draw (407,163) node    {$\textcolor[rgb]{0.29,0.56,0.89}{4}$};
\draw (317,208) node    {$\textcolor[rgb]{0.29,0.56,0.89}{5}$};
\draw (228,145) node    {$\textcolor[rgb]{0.29,0.56,0.89}{5}$};
\draw (204,190) node    {$\textcolor[rgb]{0.29,0.56,0.89}{7}$};
\draw (253,189) node    {$\textcolor[rgb]{0.29,0.56,0.89}{2}$};
\draw (293,119) node    {$\textcolor[rgb]{0.29,0.56,0.89}{7}$};
\draw (342,119) node    {$\textcolor[rgb]{0.29,0.56,0.89}{2}$};
\draw (294,50) node    {$\textcolor[rgb]{0.29,0.56,0.89}{2}$};
\draw (342,49) node    {$\textcolor[rgb]{0.29,0.56,0.89}{8}$};
\draw (432,119) node    {$\textcolor[rgb]{0.29,0.56,0.89}{8}$};
\draw (382,119) node    {$\textcolor[rgb]{0.29,0.56,0.89}{2}$};

\end{tikzpicture}

\end {figure}

Consider the word $w$, the snake graph $G_w$, and the associated continued fraction $\text{CF}(w) = [a_1 , a_2 \dots a_k]$. Let $G_{w}^{a_1}$ be the subsnake graph of $G_w$ obtained from $G_w$ by deleting its first $a_1$ tiles. Denote the cardinality of the set $\mathbb{P}_w$ by $\abs{\mathbb{P}_w}$. Let $\mathbb{P}_{w}^{a_1}$ be the perfect matchings on $G_{w}^{a_1},$ and $\abs{\mathbb{P}_{w}^{a_1}}$ the cardinality of the set $\mathbb{P}_{w}^{a_1}$.

\newpage  

\begin{thm} \label{CF_card_quotient_P}
    (Theorem 3.4 in \cite{ccanakcci2018cluster}) For any word $w$, its associated continued fraction $\textit{CF}(w)$ is equal to the quotient of cardinalities
$$\textit{CF}(w) = \frac{\abs{\mathbb{P}_w}}{\abs{\mathbb{P}_{w}^{a_1}}},$$ and the fraction on the right hand side is reduced. 
\end{thm}

We now give the set $\mathbb{P}_w$ a poset structure.

Let $P \in \mathbb{P}_w$. A \textit{twist} of $P$ at tile $i$ is the local move that takes two edges in $P$ that are located opposite one another on tile $T_i$ of $G_w$ and replaces them with the remaining two edges of $T_i$.

Directly below is the local picture for the twist at a generic tile $T_i$.

\begin {figure}[h!]
    \centering
    \caption{Twist of a perfect matching at tile $T_i$}
    \label{fig:twist}

\tikzset{every picture/.style={line width=0.75pt}} 

\begin{tikzpicture}[x=0.75pt,y=0.75pt,yscale=-1,xscale=1]

\draw    (218,205) -- (218,245) ;
\draw [color={rgb, 255:red, 74; green, 144; blue, 226 }  ,draw opacity=1 ][line width=2.25]    (218,205) -- (258,205) ;
\draw [color={rgb, 255:red, 74; green, 144; blue, 226 }  ,draw opacity=1 ][line width=2.25]    (218,245) -- (258,245) ;
\draw    (258,205) -- (258,245) ;
\draw [color={rgb, 255:red, 74; green, 144; blue, 226 }  ,draw opacity=1 ][line width=2.25]    (318,205) -- (318,245) ;
\draw    (318,205) -- (358,205) ;
\draw    (318,245) -- (358,245) ;
\draw [color={rgb, 255:red, 74; green, 144; blue, 226 }  ,draw opacity=1 ][line width=2.25]    (358,205) -- (358,245) ;
\draw    (271,225) -- (305,225) ;
\draw [shift={(308,225)}, rotate = 180] [fill={rgb, 255:red, 0; green, 0; blue, 0 }  ][line width=0.08]  [draw opacity=0] (8.93,-4.29) -- (0,0) -- (8.93,4.29) -- cycle    ;
\draw [shift={(268,225)}, rotate = 0] [fill={rgb, 255:red, 0; green, 0; blue, 0 }  ][line width=0.08]  [draw opacity=0] (8.93,-4.29) -- (0,0) -- (8.93,4.29) -- cycle    ;

\draw (237.5,226) node    {$T_{i}$};
\draw (337.5,226) node    {$T_{i}$};

\end{tikzpicture}

\end {figure}
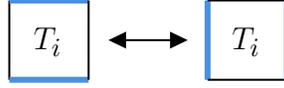

\newpage

An $\textit{up-twist}$ at tile $T_i$ is a twist that meets the \textit{twist-parity condition} (see Theorem 5.4 in \cite{musiker2013bases}):
    
    \begin{enumerate}[(1)]
        \item If $i$ is odd, the horizontal edges of $T_i$ are replaced with the vertical edges of $T_i$, or 
        
        \item If $i$ is even, the vertical edges of $T_i$ are replaced with the horizontal edges of $T_i$.
    \end{enumerate}

\begin{defn} \label{min_P_def}
     The \textit{minimal matching} $P_{-}$ of $\mathbb{P}_w$ is the unique perfect matching of $G_w$ such that every edge in $P_{-}$ is a boundary edge of $G_w$ and the S edge of tile $T_1$ is in $P_{-}.$ The \textit{maximal matching} $P_{+}$ of $\mathbb{P}_w$ is the unique perfect matching of $G_w$ such that every edge in $P_{+}$ is a boundary edge of $G_w$ and the S edge of tile $T_1$ is not in $P_{+}.$ 
\end{defn}

\begin{defn}
     Define a \textit{poset structure on $\mathbb{P}_w$} as follows. The unique minimal element of $\mathbb{P}_w$ is the minimal matching $P_{-}$, and the unique maximal element is $P_{+}$. A perfect matching $P_{2}$ covers a perfect matching $P_1$ if there exists a tile $T_i$ such that $P_2$ can be obtained from $P_1$ by performing a single up-twist of $P_1$ at $T_i$.
\end{defn}

\begin{ex}
   The poset $\mathbb{P}_{ab}$ of perfect matchings on the snake graph $G_{ab}$ is shown below in Figure \ref{fig:P_ab}. Note that $\mathbb{P}_{ab}$ is isomorphic to the Fibonacci cube $\Gamma_3$. The sum of the weights attached to each perfect matching gives the cluster variable $x_{ab}$ shown in Example \ref{fig:gamma_ab}.
\end{ex}

\newpage

\begin{rmk} \label{P_w_chain_or_fib}
    Let the length of $w$ be $l(w) = n-1$. 
 
 \begin{enumerate}[(a)]
     \item The poset of perfect matchings $\mathbb{P}_w$ on a zigzag snake graph $G_w$ is isomorphic to a linear chain with $n+1$ elements and $n$ edges. 
     
     \item The poset of perfect matchings $\mathbb{P}_w$ on a straight snake graph $G_w$ with $n$ tiles is isomorphic to the Fibonacci cube $\Gamma_n$. 
 \end{enumerate}

\end{rmk}

\begin{rmk}
    The two posets $\mathbb{A}_w$ and $\mathbb{T}_w$ defined below are each isomorphic to a linear chain when $w$ is straight, and are each isomorphic to a Fibonacci cube when $w$ is zigzag. In general, the three posets defined in this chapter are isomorphic to one another when they are parameterized by the same word. In Proposition \ref{bij_PAT}, we recall explicit isomorphisms between these posets which respects the additional node structure present in each.
\end{rmk}

\section {Perfect Matchings of Angles}

Fix a word $w$ of length $l(w) = n-1$. By the constructions above, this choice determines the arc $\gamma_w = \gamma_{a \rightarrow b}$ in the triangulated $(n+3)$-gon $\Sigma_w$. Recall the notation $\Sigma_w = [\Delta_0 , \Delta_1 , \dots , \Delta_n].$

\begin{defn}
     Any angle $o$ of the triangle $\Delta_i$ is \textit{incident} to its vertex. A \textit{perfect matching $\alpha$ of angles in $\Sigma_{w}$} is a selection of $n+1$ angles from the triangles $\Delta_0 , \Delta_1 , \dots , \Delta_n$ of $\Sigma_{w}$, one per triangle, such that 

\begin{enumerate}[(1)]

    \item Each angle is incident to an endpoint of one of the internal diagonals $\delta_{1} , \delta_{2}, \dots , \delta_{n}$.
    
    \item No two angles are incident to the same vertex of the polygon $\Sigma$. 
    
\end{enumerate}
\end{defn}

The \textit{weight} of each angle $o$ in $\Delta_i$ is the cluster variable $x_{o}$ associated to the edge of $\Delta_i$ opposite of $o$. The \textit{weight of $\alpha$} is defined to be the product of initial cluster variables $x_{\alpha} = \prod_{o \in \alpha} x_{o}.$ Let $\mathbb{A}_w$ be the set of all perfect matchings of angles in $\Sigma_{w}.$ 

\begin {figure}[h!]
    \centering
    \caption{The perfect matching of angles $A_{-}$ on $\Sigma_{ab}$.}
    \label{fig:A_}
    \begin{tikzpicture}[x=0.75pt,y=0.75pt,yscale=-1,xscale=1]

\draw    (254.5,329.38) -- (279.63,304.25) ;
\draw    (279.63,354.5) -- (254.5,329.38) ;
\draw    (329.88,354.5) -- (279.63,354.5) ;
\draw    (329.88,304.25) -- (279.63,304.25) ;
\draw    (329.88,354.5) -- (355,329.38) ;
\draw    (355,329.38) -- (329.88,304.25) ;
\draw    (279.63,354.5) -- (279.63,304.25) ;
\draw    (329.88,354.5) -- (329.88,304.25) ;
\draw    (279.63,354.5) -- (329.88,304.25) ;
\draw  [color={rgb, 255:red, 74; green, 144; blue, 226 }  ,draw opacity=1 ][fill={rgb, 255:red, 74; green, 144; blue, 226 }  ,fill opacity=1 ] (272.37,314.87) .. controls (272.37,313.42) and (273.56,312.24) .. (275.01,312.24) .. controls (276.47,312.24) and (277.65,313.42) .. (277.65,314.87) .. controls (277.65,316.33) and (276.47,317.51) .. (275.01,317.51) .. controls (273.56,317.51) and (272.37,316.33) .. (272.37,314.87) -- cycle ;
\draw  [color={rgb, 255:red, 74; green, 144; blue, 226 }  ,draw opacity=1 ][fill={rgb, 255:red, 74; green, 144; blue, 226 }  ,fill opacity=1 ] (316.16,308.77) .. controls (316.16,307.32) and (317.35,306.13) .. (318.8,306.13) .. controls (320.26,306.13) and (321.44,307.32) .. (321.44,308.77) .. controls (321.44,310.23) and (320.26,311.41) .. (318.8,311.41) .. controls (317.35,311.41) and (316.16,310.23) .. (316.16,308.77) -- cycle ;
\draw  [color={rgb, 255:red, 74; green, 144; blue, 226 }  ,draw opacity=1 ][fill={rgb, 255:red, 74; green, 144; blue, 226 }  ,fill opacity=1 ] (288.17,349.69) .. controls (288.17,348.23) and (289.35,347.05) .. (290.81,347.05) .. controls (292.26,347.05) and (293.44,348.23) .. (293.44,349.69) .. controls (293.44,351.15) and (292.26,352.33) .. (290.81,352.33) .. controls (289.35,352.33) and (288.17,351.15) .. (288.17,349.69) -- cycle ;
\draw  [color={rgb, 255:red, 74; green, 144; blue, 226 }  ,draw opacity=1 ][fill={rgb, 255:red, 74; green, 144; blue, 226 }  ,fill opacity=1 ] (331.6,343.23) .. controls (331.6,341.77) and (332.78,340.59) .. (334.24,340.59) .. controls (335.69,340.59) and (336.87,341.77) .. (336.87,343.23) .. controls (336.87,344.69) and (335.69,345.87) .. (334.24,345.87) .. controls (332.78,345.87) and (331.6,344.69) .. (331.6,343.23) -- cycle ;

\draw (262,345) node    {$\textcolor[rgb]{0.29,0.56,0.89}{6}$};
\draw (286,328) node    {$\textcolor[rgb]{0.29,0.56,0.89}{1}$};
\draw (325,328) node    {$\textcolor[rgb]{0.29,0.56,0.89}{3}$};
\draw (350,311) node    {$\textcolor[rgb]{0.29,0.56,0.89}{9}$};

\end{tikzpicture}

\end {figure}
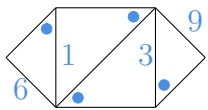

\begin{ex}
   Shown in Figure \ref{fig:A_} is one perfect matching of angles in $\Sigma_{ab}$ (in fact, it is the \textit{minimal matching} $A_{-}$, explained below). Selecting an angle to be included in a particular matching is visualized by placing a ball ``close'' to that angle inside the appropriate triangle. 
\end{ex}

Note that the weight $x_1 x_3 x_6 x_9$ of this element is the same as the weight of the snake graph perfect matching shown in Example \ref{fig:P_}.

\begin{thm}
    (Theorem 1.2 in \cite{yurikusa2019cluster}) Let $w$ be any word, and consider the set $\mathbb{A}_w$ of perfect matchings of angles on the triangulated surface $\Sigma_{w}.$ Then the cluster variable $x_{w}$ can be written as
    $$x_{w} = \frac{1}{x_1 x_2 \dots x_n} \sum_{\alpha \in \mathbb{A}_w} x_{\alpha}.$$
\end{thm}

\begin{defn}
     An angle is \textit{incident} to the two arcs in $\Delta_w$ which are sides of that angle. By a \textit{boundary angle} of $\Sigma_{w}$, we mean any angle of $\Sigma_{w}$ that is incident to exactly one boundary edge of $\Delta_w$. Any angle that is neither incident to $a$ or $b$ nor to a boundary angle will be called an \textit{internal angle}.
\end{defn}

We now give the set $\mathbb{A}_w$ a poset structure.

Let $A \in \mathbb{A}_w$. A \textit{twist} of $A$ at diagonal $\delta_{i}$ is the local move that takes two angles in $A$ incident to opposite vertices of the same internal diagonal $i$ and replaces them with the remaining two angles incident to $\delta_{i}$.

Directly below is the local picture for the twist at diagonal $\delta_{i}$.

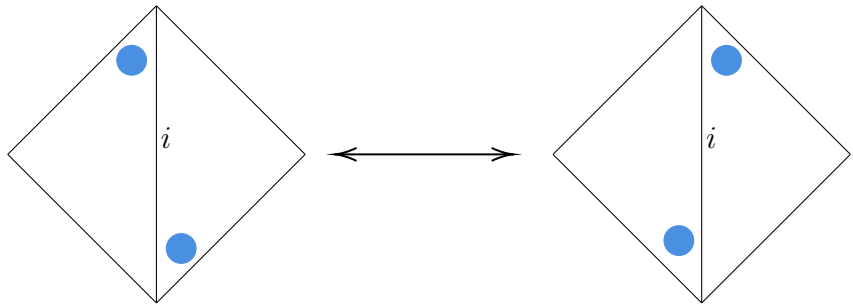
\begin {figure}[h!]
    \centering
    \caption{Twist of a perfect matching of angles at diagonal $\delta_{i}$}
    \label{fig:twist2}
    \begin{tikzpicture}[x=0.75pt,y=0.75pt,yscale=-1,xscale=1]

\draw  [color={rgb, 255:red, 74; green, 144; blue, 226 }  ,draw opacity=1 ][fill={rgb, 255:red, 74; green, 144; blue, 226 }  ,fill opacity=1 ] (304.98,242.49) .. controls (304.98,238.35) and (308.35,234.98) .. (312.49,234.98) .. controls (316.64,234.98) and (320,238.35) .. (320,242.49) .. controls (320,246.64) and (316.64,250) .. (312.49,250) .. controls (308.35,250) and (304.98,246.64) .. (304.98,242.49) -- cycle ;
\draw  [color={rgb, 255:red, 74; green, 144; blue, 226 }  ,draw opacity=1 ][fill={rgb, 255:red, 74; green, 144; blue, 226 }  ,fill opacity=1 ] (280,147.51) .. controls (280,143.36) and (283.36,140) .. (287.51,140) .. controls (291.65,140) and (295.02,143.36) .. (295.02,147.51) .. controls (295.02,151.65) and (291.65,155.02) .. (287.51,155.02) .. controls (283.36,155.02) and (280,151.65) .. (280,147.51) -- cycle ;
\draw    (300,120) -- (375,195) ;
\draw    (225,195) -- (300,120) ;
\draw    (225,195) -- (300,270) ;
\draw    (300,270) -- (375,195) ;
\draw    (300,120) -- (300,270) ;
\draw [line width=0.75]    (392,195) -- (478,195) ;
\draw [shift={(480,195)}, rotate = 180] [color={rgb, 255:red, 0; green, 0; blue, 0 }  ][line width=0.75]    (10.93,-3.29) .. controls (6.95,-1.4) and (3.31,-0.3) .. (0,0) .. controls (3.31,0.3) and (6.95,1.4) .. (10.93,3.29)   ;
\draw [shift={(390,195)}, rotate = 0] [color={rgb, 255:red, 0; green, 0; blue, 0 }  ][line width=0.75]    (10.93,-3.29) .. controls (6.95,-1.4) and (3.31,-0.3) .. (0,0) .. controls (3.31,0.3) and (6.95,1.4) .. (10.93,3.29)   ;
\draw  [color={rgb, 255:red, 74; green, 144; blue, 226 }  ,draw opacity=1 ][fill={rgb, 255:red, 74; green, 144; blue, 226 }  ,fill opacity=1 ] (580,147.51) .. controls (580,143.36) and (583.36,140) .. (587.51,140) .. controls (591.65,140) and (595.02,143.36) .. (595.02,147.51) .. controls (595.02,151.65) and (591.65,155.02) .. (587.51,155.02) .. controls (583.36,155.02) and (580,151.65) .. (580,147.51) -- cycle ;
\draw  [color={rgb, 255:red, 74; green, 144; blue, 226 }  ,draw opacity=1 ][fill={rgb, 255:red, 74; green, 144; blue, 226 }  ,fill opacity=1 ] (555.98,238.51) .. controls (555.98,234.36) and (559.35,231) .. (563.49,231) .. controls (567.64,231) and (571,234.36) .. (571,238.51) .. controls (571,242.65) and (567.64,246.02) .. (563.49,246.02) .. controls (559.35,246.02) and (555.98,242.65) .. (555.98,238.51) -- cycle ;
\draw    (575,120) -- (650,195) ;
\draw    (500,195) -- (575,120) ;
\draw    (500,195) -- (575,270) ;
\draw    (575,270) -- (650,195) ;
\draw    (575,120) -- (575,270) ;

\draw (301,180) node [anchor=north west][inner sep=0.75pt]    {$i$};
\draw (576,180) node [anchor=north west][inner sep=0.75pt]    {$i$};

\end{tikzpicture}
  
\end {figure}

The arc $\gamma_{a \rightarrow b}$ naturally partitions the set $\{\text{vertices of $\Sigma$}\} - \{ a,b \}$ into two sets; those vertices which are to the left of $\gamma_{a \rightarrow b}$, and those vertices of $\Sigma$ to the right of $\gamma_{a \rightarrow b}$. Let $l_i$ be the endpoint of $i$ to the left of $\gamma_{a \rightarrow b}$, and $r_i$ the endpoint to the right. A twist at $\delta_{i}$ is an \textit{up-twist} if the angle from $\Delta_{i}$ matched to $r_i$ is replaced with the angle matched to $l_i$, and the angle from $\Delta_{i-1}$ matched to $l_i$ is replaced with the angle matched to $r_i$.

\begin{defn}
     The \textit{minimal matching} $A_{-}$ of $\mathbb{A}_w$ is the unique perfect matching of angles in $\Sigma_{w}$ such that the boundary angle in $\Delta_0$ with boundary edge $\delta_{2n+1}$ is matched, and only boundary angles are used in $A_{-}$. The \textit{maximal matching $A_{+}$} is the unique perfect matching of angles in $\Sigma_w$ such that the boundary angle in $\Delta_0$ with boundary edge $\delta_{2n}$ is matched, and only boundary angles are used in $A_{+}$.
\end{defn}

\begin{defn}
     The \textit{poset structure on $\mathbb{A}_w$} is defined as follows. The unique minimal element of $\mathbb{A}_w$ is the minimal matching $A_{-}$, and the unique maximal element is $A_{+}$. A perfect matching of angles $A_2$ covers a perfect matching of angles $A_1$ if there exists a diagonal $\delta_{i}$ such that $A_2$ can be obtained from $A_1$ by performing a single up-twist of $A_1$ at $\delta_{i}$.
\end{defn}

\begin{rmk}
    Alternatively, $A_{-}$ can be defined by the following \textit{min-condition}, found in \cite{yurikusa2018combinatorial}. At each vertex $v$ of $\Sigma_{w}$, order the angles incident to $v$ in counterclockwise order around $v$. For each vertex $v$ of the triangulated polygon $\Sigma_{w}$ that is incident to at least one internal diagonal of $\Delta_w$, the angle $o \in A_{-}$ at $v$ is the first angle at $v$. Similarly, a \textit{max-condition} (replace ``counterclockwise'' with ``clockwise'' in the above) can be used to define $A_{+}$.
\end{rmk}

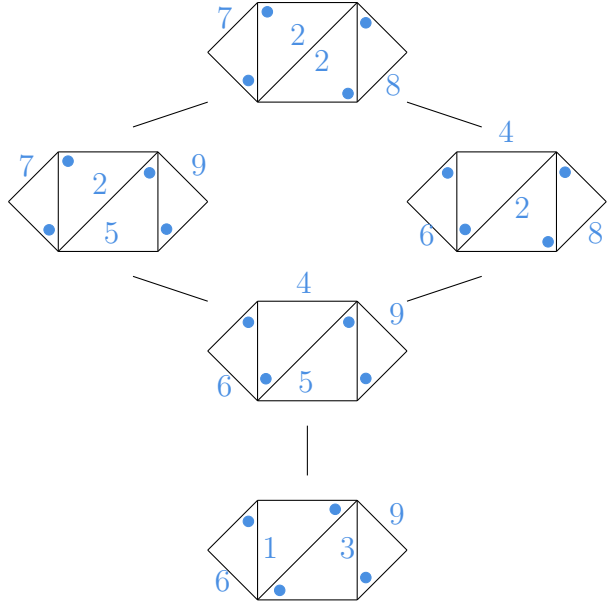
\begin {figure}[h!]
    \centering
    \caption{The poset $\mathbb{A}_{ab}$}
    \label{fig:A_ab}
    \begin{tikzpicture}[x=0.75pt,y=0.75pt,yscale=-1,xscale=1]

\draw    (254.5,329.38) -- (279.63,304.25) ;
\draw    (279.63,354.5) -- (254.5,329.38) ;
\draw    (329.88,354.5) -- (279.63,354.5) ;
\draw    (329.88,304.25) -- (279.63,304.25) ;
\draw    (329.88,354.5) -- (355,329.38) ;
\draw    (355,329.38) -- (329.88,304.25) ;
\draw    (279.63,354.5) -- (279.63,304.25) ;
\draw    (329.88,354.5) -- (329.88,304.25) ;
\draw    (279.63,354.5) -- (329.88,304.25) ;
\draw    (304.75,291.69) -- (304.75,266.56) ;
\draw    (355,203.75) -- (392.69,191.19) ;
\draw    (216.81,115.81) -- (254.5,103.25) ;
\draw    (355,103.25) -- (392.69,115.81) ;
\draw    (216.81,191.19) -- (254.5,203.75) ;
\draw  [color={rgb, 255:red, 74; green, 144; blue, 226 }  ,draw opacity=1 ][fill={rgb, 255:red, 74; green, 144; blue, 226 }  ,fill opacity=1 ] (272.37,314.87) .. controls (272.37,313.42) and (273.56,312.24) .. (275.01,312.24) .. controls (276.47,312.24) and (277.65,313.42) .. (277.65,314.87) .. controls (277.65,316.33) and (276.47,317.51) .. (275.01,317.51) .. controls (273.56,317.51) and (272.37,316.33) .. (272.37,314.87) -- cycle ;
\draw  [color={rgb, 255:red, 74; green, 144; blue, 226 }  ,draw opacity=1 ][fill={rgb, 255:red, 74; green, 144; blue, 226 }  ,fill opacity=1 ] (316.16,308.77) .. controls (316.16,307.32) and (317.35,306.13) .. (318.8,306.13) .. controls (320.26,306.13) and (321.44,307.32) .. (321.44,308.77) .. controls (321.44,310.23) and (320.26,311.41) .. (318.8,311.41) .. controls (317.35,311.41) and (316.16,310.23) .. (316.16,308.77) -- cycle ;
\draw  [color={rgb, 255:red, 74; green, 144; blue, 226 }  ,draw opacity=1 ][fill={rgb, 255:red, 74; green, 144; blue, 226 }  ,fill opacity=1 ] (288.17,349.69) .. controls (288.17,348.23) and (289.35,347.05) .. (290.81,347.05) .. controls (292.26,347.05) and (293.44,348.23) .. (293.44,349.69) .. controls (293.44,351.15) and (292.26,352.33) .. (290.81,352.33) .. controls (289.35,352.33) and (288.17,351.15) .. (288.17,349.69) -- cycle ;
\draw  [color={rgb, 255:red, 74; green, 144; blue, 226 }  ,draw opacity=1 ][fill={rgb, 255:red, 74; green, 144; blue, 226 }  ,fill opacity=1 ] (331.6,343.23) .. controls (331.6,341.77) and (332.78,340.59) .. (334.24,340.59) .. controls (335.69,340.59) and (336.87,341.77) .. (336.87,343.23) .. controls (336.87,344.69) and (335.69,345.87) .. (334.24,345.87) .. controls (332.78,345.87) and (331.6,344.69) .. (331.6,343.23) -- cycle ;
\draw    (254.5,228.88) -- (279.63,203.75) ;
\draw    (279.63,254) -- (254.5,228.88) ;
\draw    (329.88,254) -- (279.63,254) ;
\draw    (329.88,203.75) -- (279.63,203.75) ;
\draw    (329.88,254) -- (355,228.88) ;
\draw    (355,228.88) -- (329.88,203.75) ;
\draw    (279.63,254) -- (279.63,203.75) ;
\draw    (329.88,254) -- (329.88,203.75) ;
\draw    (279.63,254) -- (329.88,203.75) ;
\draw  [color={rgb, 255:red, 74; green, 144; blue, 226 }  ,draw opacity=1 ][fill={rgb, 255:red, 74; green, 144; blue, 226 }  ,fill opacity=1 ] (272.37,214.37) .. controls (272.37,212.92) and (273.56,211.74) .. (275.01,211.74) .. controls (276.47,211.74) and (277.65,212.92) .. (277.65,214.37) .. controls (277.65,215.83) and (276.47,217.01) .. (275.01,217.01) .. controls (273.56,217.01) and (272.37,215.83) .. (272.37,214.37) -- cycle ;
\draw  [color={rgb, 255:red, 74; green, 144; blue, 226 }  ,draw opacity=1 ][fill={rgb, 255:red, 74; green, 144; blue, 226 }  ,fill opacity=1 ] (281.2,242.88) .. controls (281.2,241.43) and (282.38,240.24) .. (283.83,240.24) .. controls (285.29,240.24) and (286.47,241.43) .. (286.47,242.88) .. controls (286.47,244.34) and (285.29,245.52) .. (283.83,245.52) .. controls (282.38,245.52) and (281.2,244.34) .. (281.2,242.88) -- cycle ;
\draw  [color={rgb, 255:red, 74; green, 144; blue, 226 }  ,draw opacity=1 ][fill={rgb, 255:red, 74; green, 144; blue, 226 }  ,fill opacity=1 ] (323.03,214.33) .. controls (323.03,212.87) and (324.21,211.69) .. (325.67,211.69) .. controls (327.12,211.69) and (328.3,212.87) .. (328.3,214.33) .. controls (328.3,215.79) and (327.12,216.97) .. (325.67,216.97) .. controls (324.21,216.97) and (323.03,215.79) .. (323.03,214.33) -- cycle ;
\draw  [color={rgb, 255:red, 74; green, 144; blue, 226 }  ,draw opacity=1 ][fill={rgb, 255:red, 74; green, 144; blue, 226 }  ,fill opacity=1 ] (331.6,242.73) .. controls (331.6,241.27) and (332.78,240.09) .. (334.24,240.09) .. controls (335.69,240.09) and (336.87,241.27) .. (336.87,242.73) .. controls (336.87,244.19) and (335.69,245.37) .. (334.24,245.37) .. controls (332.78,245.37) and (331.6,244.19) .. (331.6,242.73) -- cycle ;
\draw    (154,153.5) -- (179.13,128.38) ;
\draw    (179.13,178.63) -- (154,153.5) ;
\draw    (229.38,178.63) -- (179.13,178.63) ;
\draw    (229.38,128.38) -- (179.13,128.38) ;
\draw    (229.38,178.63) -- (254.5,153.5) ;
\draw    (254.5,153.5) -- (229.38,128.38) ;
\draw    (179.13,178.63) -- (179.13,128.38) ;
\draw    (229.38,178.63) -- (229.38,128.38) ;
\draw    (179.13,178.63) -- (229.38,128.38) ;
\draw  [color={rgb, 255:red, 74; green, 144; blue, 226 }  ,draw opacity=1 ][fill={rgb, 255:red, 74; green, 144; blue, 226 }  ,fill opacity=1 ] (171.87,167.71) .. controls (171.87,166.26) and (173.06,165.08) .. (174.51,165.08) .. controls (175.97,165.08) and (177.15,166.26) .. (177.15,167.71) .. controls (177.15,169.17) and (175.97,170.35) .. (174.51,170.35) .. controls (173.06,170.35) and (171.87,169.17) .. (171.87,167.71) -- cycle ;
\draw  [color={rgb, 255:red, 74; green, 144; blue, 226 }  ,draw opacity=1 ][fill={rgb, 255:red, 74; green, 144; blue, 226 }  ,fill opacity=1 ] (181.41,133.05) .. controls (181.41,131.59) and (182.59,130.41) .. (184.05,130.41) .. controls (185.51,130.41) and (186.69,131.59) .. (186.69,133.05) .. controls (186.69,134.51) and (185.51,135.69) .. (184.05,135.69) .. controls (182.59,135.69) and (181.41,134.51) .. (181.41,133.05) -- cycle ;
\draw  [color={rgb, 255:red, 74; green, 144; blue, 226 }  ,draw opacity=1 ][fill={rgb, 255:red, 74; green, 144; blue, 226 }  ,fill opacity=1 ] (222.53,138.95) .. controls (222.53,137.5) and (223.71,136.32) .. (225.17,136.32) .. controls (226.62,136.32) and (227.8,137.5) .. (227.8,138.95) .. controls (227.8,140.41) and (226.62,141.59) .. (225.17,141.59) .. controls (223.71,141.59) and (222.53,140.41) .. (222.53,138.95) -- cycle ;
\draw  [color={rgb, 255:red, 74; green, 144; blue, 226 }  ,draw opacity=1 ][fill={rgb, 255:red, 74; green, 144; blue, 226 }  ,fill opacity=1 ] (231.1,167.35) .. controls (231.1,165.9) and (232.28,164.72) .. (233.74,164.72) .. controls (235.19,164.72) and (236.37,165.9) .. (236.37,167.35) .. controls (236.37,168.81) and (235.19,169.99) .. (233.74,169.99) .. controls (232.28,169.99) and (231.1,168.81) .. (231.1,167.35) -- cycle ;
\draw    (355,153.5) -- (380.13,128.38) ;
\draw    (380.13,178.63) -- (355,153.5) ;
\draw    (430.38,178.63) -- (380.13,178.63) ;
\draw    (430.38,128.38) -- (380.13,128.38) ;
\draw    (430.38,178.63) -- (455.5,153.5) ;
\draw    (455.5,153.5) -- (430.38,128.38) ;
\draw    (380.13,178.63) -- (380.13,128.38) ;
\draw    (430.38,178.63) -- (430.38,128.38) ;
\draw    (380.13,178.63) -- (430.38,128.38) ;
\draw  [color={rgb, 255:red, 74; green, 144; blue, 226 }  ,draw opacity=1 ][fill={rgb, 255:red, 74; green, 144; blue, 226 }  ,fill opacity=1 ] (372.87,139) .. controls (372.87,137.54) and (374.06,136.36) .. (375.51,136.36) .. controls (376.97,136.36) and (378.15,137.54) .. (378.15,139) .. controls (378.15,140.46) and (376.97,141.64) .. (375.51,141.64) .. controls (374.06,141.64) and (372.87,140.46) .. (372.87,139) -- cycle ;
\draw  [color={rgb, 255:red, 74; green, 144; blue, 226 }  ,draw opacity=1 ][fill={rgb, 255:red, 74; green, 144; blue, 226 }  ,fill opacity=1 ] (381.7,167.51) .. controls (381.7,166.05) and (382.88,164.87) .. (384.33,164.87) .. controls (385.79,164.87) and (386.97,166.05) .. (386.97,167.51) .. controls (386.97,168.96) and (385.79,170.15) .. (384.33,170.15) .. controls (382.88,170.15) and (381.7,168.96) .. (381.7,167.51) -- cycle ;
\draw  [color={rgb, 255:red, 74; green, 144; blue, 226 }  ,draw opacity=1 ][fill={rgb, 255:red, 74; green, 144; blue, 226 }  ,fill opacity=1 ] (423.53,173.77) .. controls (423.53,172.31) and (424.71,171.13) .. (426.17,171.13) .. controls (427.62,171.13) and (428.8,172.31) .. (428.8,173.77) .. controls (428.8,175.23) and (427.62,176.41) .. (426.17,176.41) .. controls (424.71,176.41) and (423.53,175.23) .. (423.53,173.77) -- cycle ;
\draw  [color={rgb, 255:red, 74; green, 144; blue, 226 }  ,draw opacity=1 ][fill={rgb, 255:red, 74; green, 144; blue, 226 }  ,fill opacity=1 ] (432.1,138.64) .. controls (432.1,137.18) and (433.28,136) .. (434.74,136) .. controls (436.19,136) and (437.37,137.18) .. (437.37,138.64) .. controls (437.37,140.1) and (436.19,141.28) .. (434.74,141.28) .. controls (433.28,141.28) and (432.1,140.1) .. (432.1,138.64) -- cycle ;
\draw    (254.5,78.13) -- (279.63,53) ;
\draw    (279.63,103.25) -- (254.5,78.13) ;
\draw    (329.88,103.25) -- (279.63,103.25) ;
\draw    (329.88,53) -- (279.63,53) ;
\draw    (329.88,103.25) -- (355,78.13) ;
\draw    (355,78.13) -- (329.88,53) ;
\draw    (279.63,103.25) -- (279.63,53) ;
\draw    (329.88,103.25) -- (329.88,53) ;
\draw    (279.63,103.25) -- (329.88,53) ;
\draw  [color={rgb, 255:red, 74; green, 144; blue, 226 }  ,draw opacity=1 ][fill={rgb, 255:red, 74; green, 144; blue, 226 }  ,fill opacity=1 ] (272.37,92.34) .. controls (272.37,90.88) and (273.56,89.7) .. (275.01,89.7) .. controls (276.47,89.7) and (277.65,90.88) .. (277.65,92.34) .. controls (277.65,93.8) and (276.47,94.98) .. (275.01,94.98) .. controls (273.56,94.98) and (272.37,93.8) .. (272.37,92.34) -- cycle ;
\draw  [color={rgb, 255:red, 74; green, 144; blue, 226 }  ,draw opacity=1 ][fill={rgb, 255:red, 74; green, 144; blue, 226 }  ,fill opacity=1 ] (281.91,57.68) .. controls (281.91,56.22) and (283.09,55.04) .. (284.55,55.04) .. controls (286.01,55.04) and (287.19,56.22) .. (287.19,57.68) .. controls (287.19,59.13) and (286.01,60.31) .. (284.55,60.31) .. controls (283.09,60.31) and (281.91,59.13) .. (281.91,57.68) -- cycle ;
\draw  [color={rgb, 255:red, 74; green, 144; blue, 226 }  ,draw opacity=1 ][fill={rgb, 255:red, 74; green, 144; blue, 226 }  ,fill opacity=1 ] (322.62,98.91) .. controls (322.62,97.45) and (323.8,96.27) .. (325.25,96.27) .. controls (326.71,96.27) and (327.89,97.45) .. (327.89,98.91) .. controls (327.89,100.37) and (326.71,101.55) .. (325.25,101.55) .. controls (323.8,101.55) and (322.62,100.37) .. (322.62,98.91) -- cycle ;
\draw  [color={rgb, 255:red, 74; green, 144; blue, 226 }  ,draw opacity=1 ][fill={rgb, 255:red, 74; green, 144; blue, 226 }  ,fill opacity=1 ] (331.6,63.23) .. controls (331.6,61.77) and (332.78,60.59) .. (334.24,60.59) .. controls (335.69,60.59) and (336.87,61.77) .. (336.87,63.23) .. controls (336.87,64.68) and (335.69,65.86) .. (334.24,65.86) .. controls (332.78,65.86) and (331.6,64.68) .. (331.6,63.23) -- cycle ;

\draw (262,345) node    {$\textcolor[rgb]{0.29,0.56,0.89}{6}$};
\draw (286,328) node    {$\textcolor[rgb]{0.29,0.56,0.89}{1}$};
\draw (325,328) node    {$\textcolor[rgb]{0.29,0.56,0.89}{3}$};
\draw (350,311) node    {$\textcolor[rgb]{0.29,0.56,0.89}{9}$};
\draw (263.06,246.44) node    {$\textcolor[rgb]{0.29,0.56,0.89}{6}$};
\draw (350,210) node    {$\textcolor[rgb]{0.29,0.56,0.89}{9}$};
\draw (303.06,194.44) node    {$\textcolor[rgb]{0.29,0.56,0.89}{4}$};
\draw (304.06,244.44) node    {$\textcolor[rgb]{0.29,0.56,0.89}{5}$};
\draw (250,135) node    {$\textcolor[rgb]{0.29,0.56,0.89}{9}$};
\draw (206.06,169.44) node    {$\textcolor[rgb]{0.29,0.56,0.89}{5}$};
\draw (365.06,170.44) node    {$\textcolor[rgb]{0.29,0.56,0.89}{6}$};
\draw (405.06,118.44) node    {$\textcolor[rgb]{0.29,0.56,0.89}{4}$};
\draw (413.06,156.44) node    {$\textcolor[rgb]{0.29,0.56,0.89}{2}$};
\draw (450.06,169.44) node    {$\textcolor[rgb]{0.29,0.56,0.89}{8}$};
\draw (312.06,82.44) node    {$\textcolor[rgb]{0.29,0.56,0.89}{2}$};
\draw (348.06,94.44) node    {$\textcolor[rgb]{0.29,0.56,0.89}{8}$};
\draw (200.06,144.44) node    {$\textcolor[rgb]{0.29,0.56,0.89}{2}$};
\draw (163.06,135.44) node    {$\textcolor[rgb]{0.29,0.56,0.89}{7}$};
\draw (300.06,69.44) node    {$\textcolor[rgb]{0.29,0.56,0.89}{2}$};
\draw (263.06,60.44) node    {$\textcolor[rgb]{0.29,0.56,0.89}{7}$};

\end{tikzpicture}
  
\end {figure}

    Figure \ref{fig:A_ab} shows the poset of perfect matchings of angles in $\Sigma_{ab}$. This poset is isomorphic to the one given in Example \ref{fig:P_ab}, and one can check that the respective weights coincide as well.

\section{$T$-Paths} \label{T_path_section}

Recall that the internal diagonals of $\Delta_w$ are labeled $1,2,\dots , n$ and are ordered $\delta_{1} < \delta_{2} < \cdots < \delta_{n}$, and the notation for the intersection points $p_i = \gamma_w \cap \delta_{i}$.

\begin{defn}
     A \textit{$T$-path from $a$ to $b$}, denoted $T = (T_1 , T_2 , \dots , T_{l(T)} ),$ is an ordered selection of bicolored (either blue or red) edges from the triangulation $\Delta_w$ subject to the following conditions:

\begin{enumerate}

\item The edges in $T$ form a path from $a$ to $b$.

\item The number of edges $l(T)$ in $T$ is odd (we call $l(T)$ the \textit{length} of the $T$-path).

\item The odd-indexed edges in $T$ are all colored blue, and the even-indexed edges in $T$ are all colored red.

\item All edges in $T$ are distinct.

\item Every red edge in $T$ crosses $\gamma_{a \rightarrow b}.$

\item If $\delta_{i}$ and $\delta_{j}$ are two internal diagonals of the triangulation that $\gamma_{a \rightarrow b}$ crosses and $i<j$ then the crossing point of $\delta_{i}$ and $\gamma_{a \rightarrow b}$ is closer to $a$ than the crossing point of $\delta_{j}$ and $\gamma_{a \rightarrow b}$. 

\end{enumerate}
\end{defn}

For any $T \in \mathbb{T}_w$ let $B_{T}$ be the set of blue edges from $T$, and let $R_{T}$ be the red edges from $T$. Define the \textit{weight $x_T$ of $T$} by

$$x_T = \prod_{b \in B_{T}} x_{b} \prod_{r \in R_{T}} x^{-1}_{r}.$$ 

Let $\mathbb{T}_w$ be the set of all $T$-paths from $a$ to $b$.

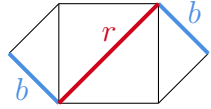
\begin {figure}[h!]
    \centering
    \caption{The $T$-path $T_{-}$ with edges from $\Delta_{ab}$}
    \label{fig:T_}
    \begin{tikzpicture}[x=0.75pt,y=0.75pt,yscale=-1,xscale=1]

\draw    (247.5,334.88) -- (272.63,309.75) ;
\draw [color={rgb, 255:red, 74; green, 144; blue, 226 }  ,draw opacity=1 ][line width=1.5]    (272.63,360) -- (247.5,334.88) ;
\draw    (322.88,360) -- (272.63,360) ;
\draw    (322.88,309.75) -- (272.63,309.75) ;
\draw    (322.88,360) -- (348,334.88) ;
\draw [color={rgb, 255:red, 74; green, 144; blue, 226 }  ,draw opacity=1 ][line width=1.5]    (348,334.88) -- (322.88,309.75) ;
\draw    (272.63,360) -- (272.63,309.75) ;
\draw    (322.88,360) -- (322.88,309.75) ;
\draw [color={rgb, 255:red, 208; green, 2; blue, 27 }  ,draw opacity=1 ][line width=1.5]    (272.63,360) -- (322.88,309.75) ;

\draw (254,353) node    {$\textcolor[rgb]{0.29,0.56,0.89}{b}$};
\draw (341,314) node    {$\textcolor[rgb]{0.29,0.56,0.89}{b}$};
\draw (298,325) node    {$\textcolor[rgb]{0.82,0.01,0.11}{r}$};

\end{tikzpicture}

\end {figure}

\begin{ex}
   Figure \ref{fig:T_} shows one $T$-path on the triangulation $\Sigma_{ab}$ (in fact, it is the \textit{minimal $T$-path} $T_{-}$, explained below). The weight of this $T$-path is $\frac{x_6 x_9}{x_2}.$ Note that multiplying this weight by $x_1 x_2 x_3$ gives $x_1 x_3 x_6 x_9,$ the weight of the elements in Example \ref{fig:P_} and Example \ref{fig:A_}. 
\end{ex}

\begin{thm}
    (Theorem 1.2 in \cite{schiffler2008cluster}) Let $w$ be any word, and consider the set $\mathbb{T}_w$ of $T$-paths on the triangulated surface $\Sigma_w$. Then the cluster variable $x_w$ can be written as 
$$x_{w} = \sum_{T \in \mathbb{T}_w } x_T.$$
\end{thm}

We now give the set $\mathbb{T}_w$ a poset structure.

Let $T \in \mathbb{T}_w$ and let $T = (T_1 , T_2 , \dots , T_{l(T)})$ be a $T$-path. Pick a red edge $T_r$ from $T$, which necessarily has as its underlying edge an internal diagonal of $\Delta_w$. The two triangles $\Delta_i$ and $\Delta_{i+1}$ from $\Sigma_w$ that are glued along the underlying edge of $T_r$ determine a triangulated quadrilateral $[\Delta_i , \Delta_{i+1}]$ with diagonal $T_r$.

Define a \textit{twist} of $T$ to be the local move that colors the four (non-triangulating) sides of the quadrilateral $[\Delta_i , \Delta_{i+1}]$ as follows: two edges of $[\Delta_i , \Delta_{i+1}]$ that are opposite one another are colored blue, the other two edges are colored red, and uncolored boundary edges in $\Delta_w$ are not allowed to be colored red. Here, if a red edge is colored blue then the colors cancel one another and that edge is not used in the resulting $T$-path, and similarly for if a blue edge is colored red \cite{gunawan2019cluster}.

Directly below is the local picture for the $T$-path twist at diagonal $\delta_{i}$. 

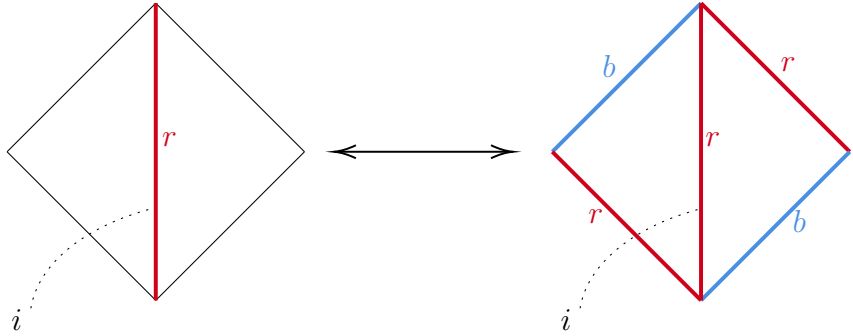
\begin {figure}[h!]
    \centering
    \caption{Twist of a $T$-path at diagonal $\delta_{i}$}
    \label{fig:twist3}
    
    \begin{tikzpicture}[x=0.75pt,y=0.75pt,yscale=-1,xscale=1]

\draw    (163,66) -- (238,141) ;
\draw    (88,141) -- (163,66) ;
\draw    (88,141) -- (163,216) ;
\draw    (163,216) -- (238,141) ;
\draw [color={rgb, 255:red, 208; green, 2; blue, 27 }  ,draw opacity=1 ][line width=1.5]    (163,66) -- (163,216) ;
\draw [line width=0.75]    (255,141) -- (341,141) ;
\draw [shift={(343,141)}, rotate = 180] [color={rgb, 255:red, 0; green, 0; blue, 0 }  ][line width=0.75]    (10.93,-3.29) .. controls (6.95,-1.4) and (3.31,-0.3) .. (0,0) .. controls (3.31,0.3) and (6.95,1.4) .. (10.93,3.29)   ;
\draw [shift={(253,141)}, rotate = 0] [color={rgb, 255:red, 0; green, 0; blue, 0 }  ][line width=0.75]    (10.93,-3.29) .. controls (6.95,-1.4) and (3.31,-0.3) .. (0,0) .. controls (3.31,0.3) and (6.95,1.4) .. (10.93,3.29)   ;
\draw [color={rgb, 255:red, 208; green, 2; blue, 27 }  ,draw opacity=1 ][line width=1.5]    (438,66) -- (513,141) ;
\draw [color={rgb, 255:red, 74; green, 144; blue, 226 }  ,draw opacity=1 ][line width=1.5]    (363,141) -- (438,66) ;
\draw [color={rgb, 255:red, 208; green, 2; blue, 27 }  ,draw opacity=1 ][line width=1.5]    (363,141) -- (438,216) ;
\draw [color={rgb, 255:red, 74; green, 144; blue, 226 }  ,draw opacity=1 ][line width=1.5]    (438,216) -- (513,141) ;
\draw [color={rgb, 255:red, 208; green, 2; blue, 27 }  ,draw opacity=1 ][line width=1.5]    (438,66) -- (438,216) ;
\draw  [dash pattern={on 0.84pt off 2.51pt}]  (100,220) .. controls (103.77,196.02) and (129.77,178.02) .. (160,170) ;
\draw  [dash pattern={on 0.84pt off 2.51pt}]  (377,220) .. controls (380.77,196.02) and (406.77,178.02) .. (437,170) ;

\draw (89,219) node [anchor=north west][inner sep=0.75pt]    {$i$};
\draw (380,169) node [anchor=north west][inner sep=0.75pt]    {$\textcolor[rgb]{0.82,0.01,0.11}{r}$};
\draw (477,92) node [anchor=north west][inner sep=0.75pt]    {$\textcolor[rgb]{0.82,0.01,0.11}{r}$};
\draw (483,169) node [anchor=north west][inner sep=0.75pt]    {$\textcolor[rgb]{0.29,0.56,0.89}{b}$};
\draw (387,91) node [anchor=north west][inner sep=0.75pt]    {$\textcolor[rgb]{0.29,0.56,0.89}{b}$};
\draw (165,130) node [anchor=north west][inner sep=0.75pt]    {$\textcolor[rgb]{0.82,0.01,0.11}{r}$};
\draw (439,130) node [anchor=north west][inner sep=0.75pt]    {$\textcolor[rgb]{0.82,0.01,0.11}{r}$};
\draw (366,219) node [anchor=north west][inner sep=0.75pt]    {$i$};

\end{tikzpicture}

\end {figure}

The arc $\gamma_{a \rightarrow b}$ naturally partitions the set $\{\text{vertices of $\Sigma$}\} - \{ a,b \}$ into two sets; those vertices which are to the left of $\gamma_{a \rightarrow b}$, and those vertices of $\Sigma$ to the right of $\gamma_{a \rightarrow b}$. Let $l_i$ be the endpoint of $\delta_{i}$ to the left of $\gamma_{a \rightarrow b}$, and $r_i$ the endpoint to the right of $\gamma_{a \rightarrow b}$. Let $u_{i-1}$ and $d_{i-1}$ be the two edges of $\Delta_{i-1}$ that are not equal to $\delta_{i}$, such that $u_{i-1}$ is adjacent to $l_i,$ and $d_{i-1}$ is adjacent to $r_i$. Let the two edges $u_{i}$ and $d_{i}$ of $\Delta_{i}$ be defined similarly. A twist at $\delta_{i}$ is an \textit{up-twist} if it colors $d_{i-1}$ and $u_i$ red, and colors $d_{i}$ and $u_{i-1}$ blue. 

Recall that the sides of the first triangle $\Delta_0$ of $\Sigma_{w}$ are labeled $1 \longrightarrow 2n \longrightarrow 2n+1 \longrightarrow 1$ in clockwise order. The \textit{minimal element} $T_{-}$ of $\mathbb{T}_w$ is the unique $T$-path that starts with the blue edge $\delta_{2n}$ and uses only internal diagonals for its edges, except for the first and last boundary edges. The \textit{maximal element} $T_{+}$ of $\mathbb{T}_w$ is the unique $T$-path from $a$ to $b$ that starts with the blue edge $\delta_{2n+1}$ and uses only internal diagonals for its edges, except for the first and last boundary edges.

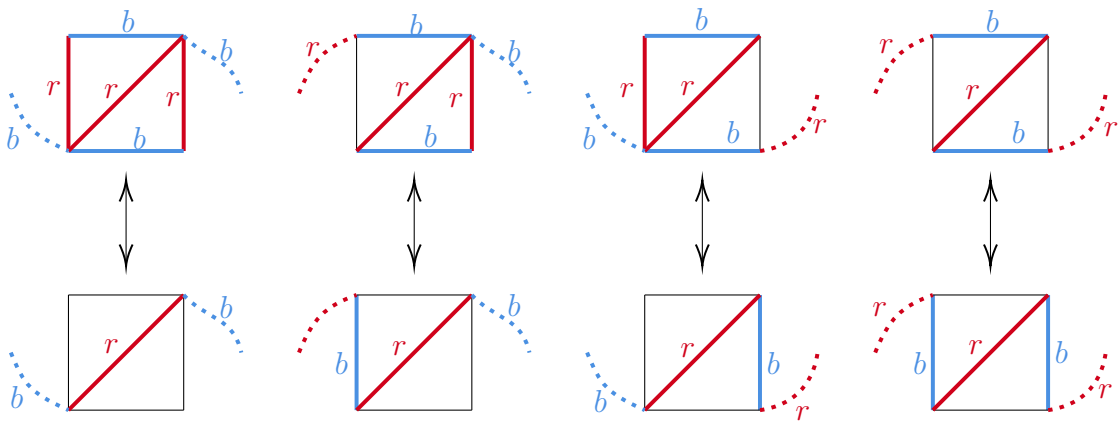
\begin {figure}[h!]
    \centering
    \caption{$T$-path twists at internal diagonals are always defined}
    \label{fig:Tpath_welldef}
    \begin{tikzpicture}[x=0.75pt,y=0.75pt,yscale=-1,xscale=1]

\draw    (74.05,224.09) -- (74.05,282.22) ;
\draw    (132.17,224.09) -- (132.17,282.22) ;
\draw    (74.05,224.09) -- (132.17,224.09) ;
\draw    (74.05,282.22) -- (132.17,282.22) ;
\draw [color={rgb, 255:red, 208; green, 2; blue, 27 }  ,draw opacity=1 ][line width=1.5]    (74.05,282.22) -- (132.17,224.09) ;
\draw [color={rgb, 255:red, 74; green, 144; blue, 226 }  ,draw opacity=1 ][line width=1.5]    (219.36,224.09) -- (219.36,282.22) ;
\draw    (277.49,224.09) -- (277.49,282.22) ;
\draw    (219.36,224.09) -- (277.49,224.09) ;
\draw    (219.36,282.22) -- (277.49,282.22) ;
\draw [color={rgb, 255:red, 208; green, 2; blue, 27 }  ,draw opacity=1 ][line width=1.5]    (219.36,282.22) -- (277.49,224.09) ;
\draw    (364.68,224.09) -- (364.68,282.22) ;
\draw [color={rgb, 255:red, 74; green, 144; blue, 226 }  ,draw opacity=1 ][line width=1.5]    (422.81,224.09) -- (422.81,282.22) ;
\draw    (364.68,224.09) -- (422.81,224.09) ;
\draw    (364.68,282.22) -- (422.81,282.22) ;
\draw [color={rgb, 255:red, 208; green, 2; blue, 27 }  ,draw opacity=1 ][line width=1.5]    (364.68,282.22) -- (422.81,224.09) ;
\draw [color={rgb, 255:red, 74; green, 144; blue, 226 }  ,draw opacity=1 ][line width=1.5]    (510,224.09) -- (510,282.22) ;
\draw [color={rgb, 255:red, 74; green, 144; blue, 226 }  ,draw opacity=1 ][line width=1.5]    (568.12,224.09) -- (568.12,282.22) ;
\draw    (510,224.09) -- (568.12,224.09) ;
\draw    (510,282.22) -- (568.12,282.22) ;
\draw [color={rgb, 255:red, 208; green, 2; blue, 27 }  ,draw opacity=1 ][line width=1.5]    (510,282.22) -- (568.12,224.09) ;
\draw [color={rgb, 255:red, 208; green, 2; blue, 27 }  ,draw opacity=1 ][line width=1.5]    (74.05,93.31) -- (74.05,151.43) ;
\draw [color={rgb, 255:red, 208; green, 2; blue, 27 }  ,draw opacity=1 ][line width=1.5]    (132.17,93.31) -- (132.17,151.43) ;
\draw [color={rgb, 255:red, 74; green, 144; blue, 226 }  ,draw opacity=1 ][line width=1.5]    (74.05,93.31) -- (132.17,93.31) ;
\draw [color={rgb, 255:red, 74; green, 144; blue, 226 }  ,draw opacity=1 ][line width=1.5]    (74.05,151.43) -- (132.17,151.43) ;
\draw [color={rgb, 255:red, 208; green, 2; blue, 27 }  ,draw opacity=1 ][line width=1.5]    (74.05,151.43) -- (132.17,93.31) ;
\draw    (219.36,93.31) -- (219.36,151.43) ;
\draw [color={rgb, 255:red, 208; green, 2; blue, 27 }  ,draw opacity=1 ][line width=1.5]    (277.49,93.31) -- (277.49,151.43) ;
\draw [color={rgb, 255:red, 74; green, 144; blue, 226 }  ,draw opacity=1 ][line width=1.5]    (219.36,93.31) -- (277.49,93.31) ;
\draw [color={rgb, 255:red, 74; green, 144; blue, 226 }  ,draw opacity=1 ][line width=1.5]    (219.36,151.43) -- (277.49,151.43) ;
\draw [color={rgb, 255:red, 208; green, 2; blue, 27 }  ,draw opacity=1 ][line width=1.5]    (219.36,151.43) -- (277.49,93.31) ;
\draw [color={rgb, 255:red, 208; green, 2; blue, 27 }  ,draw opacity=1 ][line width=1.5]    (364.68,93.31) -- (364.68,151.43) ;
\draw    (422.81,93.31) -- (422.81,151.43) ;
\draw [color={rgb, 255:red, 74; green, 144; blue, 226 }  ,draw opacity=1 ][line width=1.5]    (364.68,93.31) -- (422.81,93.31) ;
\draw [color={rgb, 255:red, 74; green, 144; blue, 226 }  ,draw opacity=1 ][line width=1.5]    (364.68,151.43) -- (422.81,151.43) ;
\draw [color={rgb, 255:red, 208; green, 2; blue, 27 }  ,draw opacity=1 ][line width=1.5]    (364.68,151.43) -- (422.81,93.31) ;
\draw    (510,93.31) -- (510,151.43) ;
\draw    (568.12,93.31) -- (568.12,151.43) ;
\draw [color={rgb, 255:red, 74; green, 144; blue, 226 }  ,draw opacity=1 ][line width=1.5]    (510,93.31) -- (568.12,93.31) ;
\draw [color={rgb, 255:red, 74; green, 144; blue, 226 }  ,draw opacity=1 ][line width=1.5]    (510,151.43) -- (568.12,151.43) ;
\draw [color={rgb, 255:red, 208; green, 2; blue, 27 }  ,draw opacity=1 ][line width=1.5]    (510,151.43) -- (568.12,93.31) ;
\draw [color={rgb, 255:red, 74; green, 144; blue, 226 }  ,draw opacity=1 ][line width=1.5]  [dash pattern={on 1.69pt off 2.76pt}]  (44.98,253.15) .. controls (51.77,269.86) and (55.16,273.74) .. (74.05,282.22) ;
\draw [color={rgb, 255:red, 74; green, 144; blue, 226 }  ,draw opacity=1 ][line width=1.5]  [dash pattern={on 1.69pt off 2.76pt}]  (335.62,253.15) .. controls (342.4,269.86) and (345.79,273.74) .. (364.68,282.22) ;
\draw [color={rgb, 255:red, 74; green, 144; blue, 226 }  ,draw opacity=1 ][line width=1.5]  [dash pattern={on 1.69pt off 2.76pt}]  (277.49,224.09) .. controls (289.12,236.81) and (302.19,234.63) .. (306.55,253.15) ;
\draw [color={rgb, 255:red, 74; green, 144; blue, 226 }  ,draw opacity=1 ][line width=1.5]  [dash pattern={on 1.69pt off 2.76pt}]  (132.17,224.09) .. controls (143.8,236.81) and (156.88,234.63) .. (161.24,253.15) ;
\draw [color={rgb, 255:red, 208; green, 2; blue, 27 }  ,draw opacity=1 ][line width=1.5]  [dash pattern={on 1.69pt off 2.76pt}]  (190.3,253.15) .. controls (196.84,236.81) and (203.38,229.54) .. (219.36,224.09) ;
\draw [color={rgb, 255:red, 208; green, 2; blue, 27 }  ,draw opacity=1 ][line width=1.5]  [dash pattern={on 1.69pt off 2.76pt}]  (480.93,253.15) .. controls (487.47,236.81) and (494.01,229.54) .. (510,224.09) ;
\draw [color={rgb, 255:red, 208; green, 2; blue, 27 }  ,draw opacity=1 ][line width=1.5]  [dash pattern={on 1.69pt off 2.76pt}]  (422.81,282.22) .. controls (440.25,279.67) and (449.69,268.77) .. (451.87,253.15) ;
\draw [color={rgb, 255:red, 208; green, 2; blue, 27 }  ,draw opacity=1 ][line width=1.5]  [dash pattern={on 1.69pt off 2.76pt}]  (568.12,282.22) .. controls (585.56,279.67) and (595.01,268.77) .. (597.19,253.15) ;
\draw [color={rgb, 255:red, 74; green, 144; blue, 226 }  ,draw opacity=1 ][line width=1.5]  [dash pattern={on 1.69pt off 2.76pt}]  (44.98,122.37) .. controls (51.77,139.08) and (55.16,142.95) .. (74.05,151.43) ;
\draw [color={rgb, 255:red, 74; green, 144; blue, 226 }  ,draw opacity=1 ][line width=1.5]  [dash pattern={on 1.69pt off 2.76pt}]  (132.17,93.31) .. controls (143.8,106.02) and (156.88,103.84) .. (161.24,122.37) ;
\draw [color={rgb, 255:red, 208; green, 2; blue, 27 }  ,draw opacity=1 ][line width=1.5]  [dash pattern={on 1.69pt off 2.76pt}]  (422.81,151.43) .. controls (440.25,148.89) and (449.69,137.99) .. (451.87,122.37) ;
\draw [color={rgb, 255:red, 208; green, 2; blue, 27 }  ,draw opacity=1 ][line width=1.5]  [dash pattern={on 1.69pt off 2.76pt}]  (190.3,122.37) .. controls (196.84,106.02) and (203.38,98.75) .. (219.36,93.31) ;
\draw [color={rgb, 255:red, 74; green, 144; blue, 226 }  ,draw opacity=1 ][line width=1.5]  [dash pattern={on 1.69pt off 2.76pt}]  (277.49,93.31) .. controls (289.12,106.02) and (302.19,103.84) .. (306.55,122.37) ;
\draw [color={rgb, 255:red, 74; green, 144; blue, 226 }  ,draw opacity=1 ][line width=1.5]  [dash pattern={on 1.69pt off 2.76pt}]  (335.62,122.37) .. controls (342.4,139.08) and (345.79,142.95) .. (364.68,151.43) ;
\draw [color={rgb, 255:red, 208; green, 2; blue, 27 }  ,draw opacity=1 ][line width=1.5]  [dash pattern={on 1.69pt off 2.76pt}]  (480.93,122.37) .. controls (487.47,106.02) and (494.01,98.75) .. (510,93.31) ;
\draw [color={rgb, 255:red, 208; green, 2; blue, 27 }  ,draw opacity=1 ][line width=1.5]  [dash pattern={on 1.69pt off 2.76pt}]  (568.12,151.43) .. controls (585.56,148.89) and (595.01,137.99) .. (597.19,122.37) ;
\draw    (103.11,167.96) -- (103.11,207.56) ;
\draw [shift={(103.11,209.56)}, rotate = 270] [color={rgb, 255:red, 0; green, 0; blue, 0 }  ][line width=0.75]    (10.93,-3.29) .. controls (6.95,-1.4) and (3.31,-0.3) .. (0,0) .. controls (3.31,0.3) and (6.95,1.4) .. (10.93,3.29)   ;
\draw [shift={(103.11,165.96)}, rotate = 90] [color={rgb, 255:red, 0; green, 0; blue, 0 }  ][line width=0.75]    (10.93,-3.29) .. controls (6.95,-1.4) and (3.31,-0.3) .. (0,0) .. controls (3.31,0.3) and (6.95,1.4) .. (10.93,3.29)   ;
\draw    (248.43,167.96) -- (248.43,207.56) ;
\draw [shift={(248.43,209.56)}, rotate = 270] [color={rgb, 255:red, 0; green, 0; blue, 0 }  ][line width=0.75]    (10.93,-3.29) .. controls (6.95,-1.4) and (3.31,-0.3) .. (0,0) .. controls (3.31,0.3) and (6.95,1.4) .. (10.93,3.29)   ;
\draw [shift={(248.43,165.96)}, rotate = 90] [color={rgb, 255:red, 0; green, 0; blue, 0 }  ][line width=0.75]    (10.93,-3.29) .. controls (6.95,-1.4) and (3.31,-0.3) .. (0,0) .. controls (3.31,0.3) and (6.95,1.4) .. (10.93,3.29)   ;
\draw    (393.74,167.96) -- (393.74,207.56) ;
\draw [shift={(393.74,209.56)}, rotate = 270] [color={rgb, 255:red, 0; green, 0; blue, 0 }  ][line width=0.75]    (10.93,-3.29) .. controls (6.95,-1.4) and (3.31,-0.3) .. (0,0) .. controls (3.31,0.3) and (6.95,1.4) .. (10.93,3.29)   ;
\draw [shift={(393.74,165.96)}, rotate = 90] [color={rgb, 255:red, 0; green, 0; blue, 0 }  ][line width=0.75]    (10.93,-3.29) .. controls (6.95,-1.4) and (3.31,-0.3) .. (0,0) .. controls (3.31,0.3) and (6.95,1.4) .. (10.93,3.29)   ;
\draw    (539.06,167.96) -- (539.06,207.56) ;
\draw [shift={(539.06,209.56)}, rotate = 270] [color={rgb, 255:red, 0; green, 0; blue, 0 }  ][line width=0.75]    (10.93,-3.29) .. controls (6.95,-1.4) and (3.31,-0.3) .. (0,0) .. controls (3.31,0.3) and (6.95,1.4) .. (10.93,3.29)   ;
\draw [shift={(539.06,165.96)}, rotate = 90] [color={rgb, 255:red, 0; green, 0; blue, 0 }  ][line width=0.75]    (10.93,-3.29) .. controls (6.95,-1.4) and (3.31,-0.3) .. (0,0) .. controls (3.31,0.3) and (6.95,1.4) .. (10.93,3.29)   ;

\draw (61.69,115.14) node [anchor=north west][inner sep=0.75pt]    {$\textcolor[rgb]{0.82,0.01,0.11}{r}$};
\draw (90.75,116.33) node [anchor=north west][inner sep=0.75pt]    {$\textcolor[rgb]{0.82,0.01,0.11}{r}$};
\draw (121.82,120.59) node [anchor=north west][inner sep=0.75pt]    {$\textcolor[rgb]{0.82,0.01,0.11}{r}$};
\draw (237.52,115.33) node [anchor=north west][inner sep=0.75pt]    {$\textcolor[rgb]{0.82,0.01,0.11}{r}$};
\draw (381.38,114.78) node [anchor=north west][inner sep=0.75pt]    {$\textcolor[rgb]{0.82,0.01,0.11}{r}$};
\draw (350.87,116.24) node [anchor=north west][inner sep=0.75pt]    {$\textcolor[rgb]{0.82,0.01,0.11}{r}$};
\draw (525.25,116.78) node [anchor=north west][inner sep=0.75pt]    {$\textcolor[rgb]{0.82,0.01,0.11}{r}$};
\draw (526.7,245.57) node [anchor=north west][inner sep=0.75pt]    {$\textcolor[rgb]{0.82,0.01,0.11}{r}$};
\draw (592.09,268.08) node [anchor=north west][inner sep=0.75pt]    {$\textcolor[rgb]{0.82,0.01,0.11}{r}$};
\draw (595,135.85) node [anchor=north west][inner sep=0.75pt]    {$\textcolor[rgb]{0.82,0.01,0.11}{r}$};
\draw (481.65,93.63) node [anchor=north west][inner sep=0.75pt]    {$\textcolor[rgb]{0.82,0.01,0.11}{r}$};
\draw (448.23,134.39) node [anchor=north west][inner sep=0.75pt]    {$\textcolor[rgb]{0.82,0.01,0.11}{r}$};
\draw (192.47,95.53) node [anchor=north west][inner sep=0.75pt]    {$\textcolor[rgb]{0.82,0.01,0.11}{r}$};
\draw (264.68,121.14) node [anchor=north west][inner sep=0.75pt]    {$\textcolor[rgb]{0.82,0.01,0.11}{r}$};
\draw (90.75,245.11) node [anchor=north west][inner sep=0.75pt]    {$\textcolor[rgb]{0.82,0.01,0.11}{r}$};
\draw (236.07,246.11) node [anchor=north west][inner sep=0.75pt]    {$\textcolor[rgb]{0.82,0.01,0.11}{r}$};
\draw (381.38,247.11) node [anchor=north west][inner sep=0.75pt]    {$\textcolor[rgb]{0.82,0.01,0.11}{r}$};
\draw (439.51,278.99) node [anchor=north west][inner sep=0.75pt]    {$\textcolor[rgb]{0.82,0.01,0.11}{r}$};
\draw (478.29,226.32) node [anchor=north west][inner sep=0.75pt]    {$\textcolor[rgb]{0.82,0.01,0.11}{r}$};
\draw (148.88,93.53) node [anchor=north west][inner sep=0.75pt]    {$\textcolor[rgb]{0.29,0.56,0.89}{b}$};
\draw (41.34,138.21) node [anchor=north west][inner sep=0.75pt]    {$\textcolor[rgb]{0.29,0.56,0.89}{b}$};
\draw (43.34,268.99) node [anchor=north west][inner sep=0.75pt]    {$\textcolor[rgb]{0.29,0.56,0.89}{b}$};
\draw (148.88,223.32) node [anchor=north west][inner sep=0.75pt]    {$\textcolor[rgb]{0.29,0.56,0.89}{b}$};
\draw (294.19,222.32) node [anchor=north west][inner sep=0.75pt]    {$\textcolor[rgb]{0.29,0.56,0.89}{b}$};
\draw (294.19,94.08) node [anchor=north west][inner sep=0.75pt]    {$\textcolor[rgb]{0.29,0.56,0.89}{b}$};
\draw (331.98,138.21) node [anchor=north west][inner sep=0.75pt]    {$\textcolor[rgb]{0.29,0.56,0.89}{b}$};
\draw (337.79,270.99) node [anchor=north west][inner sep=0.75pt]    {$\textcolor[rgb]{0.29,0.56,0.89}{b}$};
\draw (570.3,245.93) node [anchor=north west][inner sep=0.75pt]    {$\textcolor[rgb]{0.29,0.56,0.89}{b}$};
\draw (498.73,251.38) node [anchor=north west][inner sep=0.75pt]    {$\textcolor[rgb]{0.29,0.56,0.89}{b}$};
\draw (548.5,135.13) node [anchor=north west][inner sep=0.75pt]    {$\textcolor[rgb]{0.29,0.56,0.89}{b}$};
\draw (535.42,79) node [anchor=north west][inner sep=0.75pt]    {$\textcolor[rgb]{0.29,0.56,0.89}{b}$};
\draw (388.65,78) node [anchor=north west][inner sep=0.75pt]    {$\textcolor[rgb]{0.29,0.56,0.89}{b}$};
\draw (244.79,81) node [anchor=north west][inner sep=0.75pt]    {$\textcolor[rgb]{0.29,0.56,0.89}{b}$};
\draw (99.47,79) node [anchor=north west][inner sep=0.75pt]    {$\textcolor[rgb]{0.29,0.56,0.89}{b}$};
\draw (105.28,138.13) node [anchor=north west][inner sep=0.75pt]    {$\textcolor[rgb]{0.29,0.56,0.89}{b}$};
\draw (252.05,137.13) node [anchor=north west][inner sep=0.75pt]    {$\textcolor[rgb]{0.29,0.56,0.89}{b}$};
\draw (404.64,135.13) node [anchor=north west][inner sep=0.75pt]    {$\textcolor[rgb]{0.29,0.56,0.89}{b}$};
\draw (424.98,251.38) node [anchor=north west][inner sep=0.75pt]    {$\textcolor[rgb]{0.29,0.56,0.89}{b}$};
\draw (207.1,251.38) node [anchor=north west][inner sep=0.75pt]    {$\textcolor[rgb]{0.29,0.56,0.89}{b}$};

\end{tikzpicture}
  
\end {figure}

\begin{rmk}
    It is not obvious that a twist of a $T$-path is always defined. Figure \ref{fig:Tpath_welldef} shows the eight possibilities for how any $T$-path looks locally at an internal diagonal which is not the first or last. In each quadrilateral shown, the top two edges are boundary edges, the other three are internal diagonals, and the dotted lines can be either, as long as a red edge is not on the boundary. Each of the four arrows indicate when two $T$-paths are related by a twist. 
\end{rmk}

Thus, twists are well-defined for any diagonal that isn't the first or last. A similar analysis shows twists are always defined for the first and last diagonals, as well.

\begin{defn}
     The \textit{poset structure} on $\mathbb{T}_w$ is defined as follows. The unique minimal element of $\mathbb{T}_w$ is $T_{-}$, and the unique maximal element is $T_{+}$. A $T$-path $T_2$ covers a $T$-path $T_1$ if there exists a diagonal $\delta_{i}$ such that $T_2$ can be obtained from $T_1$ by performing a single up-twist of $T_1$ at diagonal $\delta_{i}$.
\end{defn}

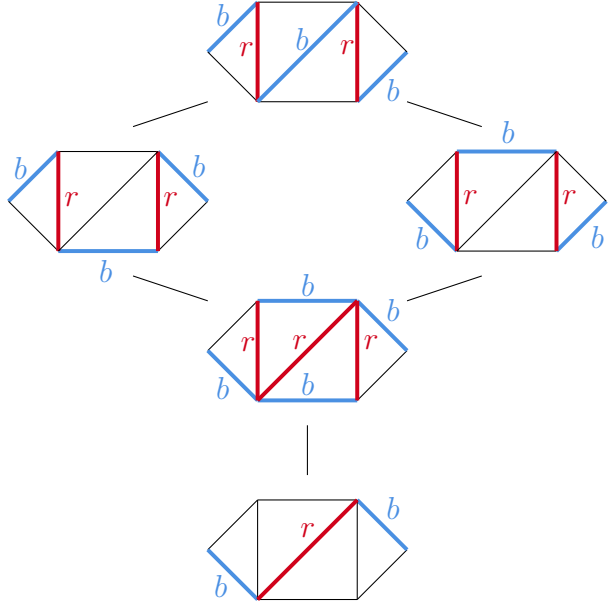
\begin {figure}[h!]
    \centering
    \caption{The poset $T_{ab}$}
    \label{fig:T_ab}
    \begin{tikzpicture}[x=0.75pt,y=0.75pt,yscale=-1,xscale=1]

\draw    (247.5,334.88) -- (272.63,309.75) ;
\draw [color={rgb, 255:red, 74; green, 144; blue, 226 }  ,draw opacity=1 ][line width=1.5]    (272.63,360) -- (247.5,334.88) ;
\draw    (322.88,360) -- (272.63,360) ;
\draw    (322.88,309.75) -- (272.63,309.75) ;
\draw    (322.88,360) -- (348,334.88) ;
\draw [color={rgb, 255:red, 74; green, 144; blue, 226 }  ,draw opacity=1 ][line width=1.5]    (348,334.88) -- (322.88,309.75) ;
\draw    (272.63,360) -- (272.63,309.75) ;
\draw    (322.88,360) -- (322.88,309.75) ;
\draw [color={rgb, 255:red, 208; green, 2; blue, 27 }  ,draw opacity=1 ][line width=1.5]    (272.63,360) -- (322.88,309.75) ;
\draw    (297.75,297.19) -- (297.75,272.06) ;
\draw    (348,209.25) -- (385.69,196.69) ;
\draw    (209.81,121.31) -- (247.5,108.75) ;
\draw    (348,108.75) -- (385.69,121.31) ;
\draw    (209.81,196.69) -- (247.5,209.25) ;
\draw    (247.5,234.38) -- (272.63,209.25) ;
\draw [color={rgb, 255:red, 74; green, 144; blue, 226 }  ,draw opacity=1 ][line width=1.5]    (272.63,259.5) -- (247.5,234.38) ;
\draw [color={rgb, 255:red, 74; green, 144; blue, 226 }  ,draw opacity=1 ][line width=1.5]    (322.88,259.5) -- (272.63,259.5) ;
\draw [color={rgb, 255:red, 74; green, 144; blue, 226 }  ,draw opacity=1 ][line width=1.5]    (322.88,209.25) -- (272.63,209.25) ;
\draw    (322.88,259.5) -- (348,234.38) ;
\draw [color={rgb, 255:red, 74; green, 144; blue, 226 }  ,draw opacity=1 ][line width=1.5]    (348,234.38) -- (322.88,209.25) ;
\draw [color={rgb, 255:red, 208; green, 2; blue, 27 }  ,draw opacity=1 ][line width=1.5]    (272.63,259.5) -- (272.63,209.25) ;
\draw [color={rgb, 255:red, 208; green, 2; blue, 27 }  ,draw opacity=1 ][line width=1.5]    (322.88,259.5) -- (322.88,209.25) ;
\draw [color={rgb, 255:red, 208; green, 2; blue, 27 }  ,draw opacity=1 ][line width=1.5]    (272.63,259.5) -- (322.88,209.25) ;
\draw [color={rgb, 255:red, 74; green, 144; blue, 226 }  ,draw opacity=1 ][line width=1.5]    (147,159) -- (172.13,133.88) ;
\draw    (172.13,184.13) -- (147,159) ;
\draw [color={rgb, 255:red, 74; green, 144; blue, 226 }  ,draw opacity=1 ][line width=1.5]    (222.38,184.13) -- (172.13,184.13) ;
\draw    (222.38,133.88) -- (172.13,133.88) ;
\draw    (222.38,184.13) -- (247.5,159) ;
\draw [color={rgb, 255:red, 74; green, 144; blue, 226 }  ,draw opacity=1 ][line width=1.5]    (247.5,159) -- (222.38,133.88) ;
\draw [color={rgb, 255:red, 208; green, 2; blue, 27 }  ,draw opacity=1 ][line width=1.5]    (172.13,184.13) -- (172.13,133.88) ;
\draw [color={rgb, 255:red, 208; green, 2; blue, 27 }  ,draw opacity=1 ][line width=1.5]    (222.38,184.13) -- (222.38,133.88) ;
\draw    (172.13,184.13) -- (222.38,133.88) ;
\draw    (348,159) -- (373.13,133.88) ;
\draw [color={rgb, 255:red, 74; green, 144; blue, 226 }  ,draw opacity=1 ][line width=1.5]    (373.13,184.13) -- (348,159) ;
\draw    (423.38,184.13) -- (373.13,184.13) ;
\draw [color={rgb, 255:red, 74; green, 144; blue, 226 }  ,draw opacity=1 ][line width=1.5]    (423.38,133.88) -- (373.13,133.88) ;
\draw [color={rgb, 255:red, 74; green, 144; blue, 226 }  ,draw opacity=1 ][line width=1.5]    (423.38,184.13) -- (448.5,159) ;
\draw    (448.5,159) -- (423.38,133.88) ;
\draw [color={rgb, 255:red, 208; green, 2; blue, 27 }  ,draw opacity=1 ][line width=1.5]    (373.13,184.13) -- (373.13,133.88) ;
\draw [color={rgb, 255:red, 208; green, 2; blue, 27 }  ,draw opacity=1 ][line width=1.5]    (423.38,184.13) -- (423.38,133.88) ;
\draw    (373.13,184.13) -- (423.38,133.88) ;
\draw [color={rgb, 255:red, 74; green, 144; blue, 226 }  ,draw opacity=1 ][line width=1.5]    (247.5,83.63) -- (272.63,58.5) ;
\draw    (272.63,108.75) -- (247.5,83.63) ;
\draw    (322.88,108.75) -- (272.63,108.75) ;
\draw    (322.88,58.5) -- (272.63,58.5) ;
\draw [color={rgb, 255:red, 74; green, 144; blue, 226 }  ,draw opacity=1 ][line width=1.5]    (322.88,108.75) -- (348,83.63) ;
\draw    (348,83.63) -- (322.88,58.5) ;
\draw [color={rgb, 255:red, 208; green, 2; blue, 27 }  ,draw opacity=1 ][line width=1.5]    (272.63,108.75) -- (272.63,58.5) ;
\draw [color={rgb, 255:red, 208; green, 2; blue, 27 }  ,draw opacity=1 ][line width=1.5]    (322.88,108.75) -- (322.88,58.5) ;
\draw [color={rgb, 255:red, 74; green, 144; blue, 226 }  ,draw opacity=1 ][line width=1.5]    (272.63,108.75) -- (322.88,58.5) ;

\draw (254,353) node    {$\textcolor[rgb]{0.29,0.56,0.89}{b}$};
\draw (341,314) node    {$\textcolor[rgb]{0.29,0.56,0.89}{b}$};
\draw (298,325) node    {$\textcolor[rgb]{0.82,0.01,0.11}{r}$};
\draw (294,230) node    {$\textcolor[rgb]{0.82,0.01,0.11}{r}$};
\draw (341,214) node    {$\textcolor[rgb]{0.29,0.56,0.89}{b}$};
\draw (255,253) node    {$\textcolor[rgb]{0.29,0.56,0.89}{b}$};
\draw (298,251) node    {$\textcolor[rgb]{0.29,0.56,0.89}{b}$};
\draw (298,201) node    {$\textcolor[rgb]{0.29,0.56,0.89}{b}$};
\draw (444,177) node    {$\textcolor[rgb]{0.29,0.56,0.89}{b}$};
\draw (399,124) node    {$\textcolor[rgb]{0.29,0.56,0.89}{b}$};
\draw (355.56,177.56) node    {$\textcolor[rgb]{0.29,0.56,0.89}{b}$};
\draw (380,157) node    {$\textcolor[rgb]{0.82,0.01,0.11}{r}$};
\draw (430,157) node    {$\textcolor[rgb]{0.82,0.01,0.11}{r}$};
\draw (295,77) node    {$\textcolor[rgb]{0.29,0.56,0.89}{b}$};
\draw (341,103) node    {$\textcolor[rgb]{0.29,0.56,0.89}{b}$};
\draw (179,158) node    {$\textcolor[rgb]{0.82,0.01,0.11}{r}$};
\draw (229,158) node    {$\textcolor[rgb]{0.82,0.01,0.11}{r}$};
\draw (267,82) node    {$\textcolor[rgb]{0.82,0.01,0.11}{r}$};
\draw (318,81) node    {$\textcolor[rgb]{0.82,0.01,0.11}{r}$};
\draw (255,65) node    {$\textcolor[rgb]{0.29,0.56,0.89}{b}$};
\draw (153,142) node    {$\textcolor[rgb]{0.29,0.56,0.89}{b}$};
\draw (196,194) node    {$\textcolor[rgb]{0.29,0.56,0.89}{b}$};
\draw (243,141) node    {$\textcolor[rgb]{0.29,0.56,0.89}{b}$};
\draw (268,231) node    {$\textcolor[rgb]{0.82,0.01,0.11}{r}$};
\draw (330,230) node    {$\textcolor[rgb]{0.82,0.01,0.11}{r}$};

\end{tikzpicture}
  
\end {figure}

\begin{ex}
   Figure \ref{fig:T_ab} shows the poset $\mathbb{T}_{ab}$ of $T$-paths from $a$ to $b$ associated to the word $w=ab$. 
\end{ex}

If we sum over the weights of this poset, we obtain
$$x_{ab} = \frac{x_6 x_9}{x_2} + \frac{x_4 x_5 x_6 x_9}{x_1 x_2 x_3} + \frac{x_5 x_7 x_9}{x_1 x_3} + \frac{x_4 x_6 x_8}{x_1 x_3} + \frac{x_2 x_7 x_8}{x_1 x_3}.$$
If we find a common denominator for the five terms in this sum, we see that this expression for $x_{ab}$ is equivalent to the one given in Example \ref{fig:gamma_ab}.

It is known that perfect matchings on $G_w$ are equivalent to both perfect matchings of angles on $\Sigma_w$ (see \cite{yurikusa2019cluster}) and $T$-paths on $\Sigma_w$ (see \cite{musiker2010cluster}).

\begin{prop} \label{bij_PAT}
    Fix the word $w$. Then there are poset isomorphisms  

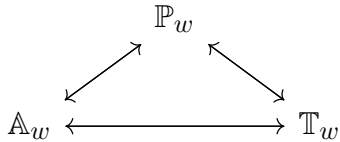
\begin{center}
\begin{tikzcd}
                                   & \mathbb{P}_w \arrow[ld] \arrow[rd] &                                    \\
\mathbb{A}_w \arrow[rr] \arrow[ru] &                                    & \mathbb{T}_w \arrow[lu] \arrow[ll]
\end{tikzcd}
\end{center}

which respect the additional node structure of each poset.
\end{prop}

We illustrate Proposition \ref{bij_PAT} with the three posets associated to the word $w=ab$ from our running example.  By Lemma 3.2 in \cite{yurikusa2019cluster} there is a bijection between the edges in $G_w$ and the angles in $\Sigma_{w},$ induced by identifying certain pairs of angles in $\widetilde{G_w}.$ The pairs of angles that are identified are those that are opposite one another in the quadrilateral determined by two consecutive tiles of $\widetilde{G_w}.$ Any pair of angles in $\widetilde{G_w}$ which have been identified correspond to a single internal angle in $\Sigma_w$.

\begin {figure}[h!]
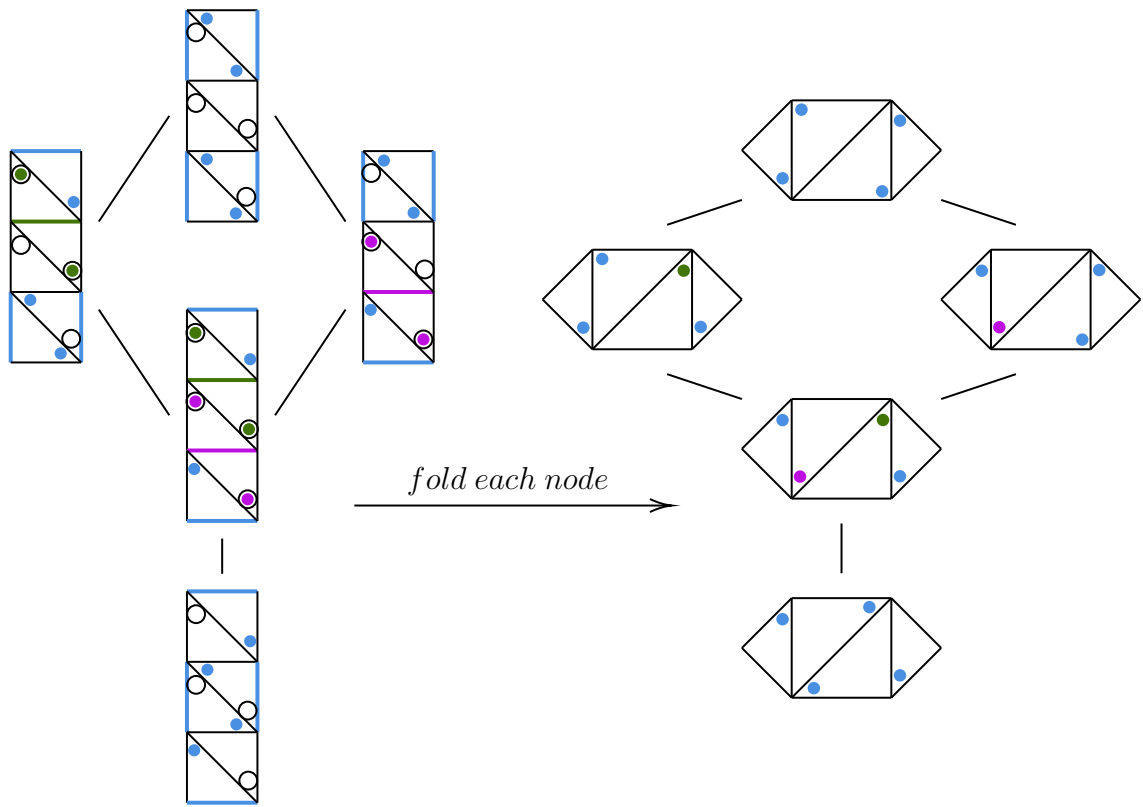

    \centering
    \caption{The map $\mathbb{P}_{ab} \longrightarrow \mathbb{A}_{ab}$ via angle identification and folding}
    \label{fig:P_to_A}

\tikzset{every picture/.style={line width=0.75pt}} 



\end {figure}

Thus, given a perfect matching $P \in \mathbb{P}_w$, we can associate to it a collection of angles in $\widetilde{G}_w$ and fold the result to obtain the corresponding perfect matching of angles $A$ in $\Sigma_w$. This assignment extends to a map of posets, since $\mathbb{A}_w$ can be obtained by starting with the poset $\mathbb{P}_w$ and folding every node. See Figure \ref{fig:P_to_A} for an illustration using the objects from our running example.

\newpage

The second map $\mathbb{A}_w \longrightarrow \mathbb{T}_w$ is defined by sending each angle in a matching $A$ to the edge in $\Delta_w$ it is opposite of, and then coloring each internal diagonal red.

\begin {figure}[h!]
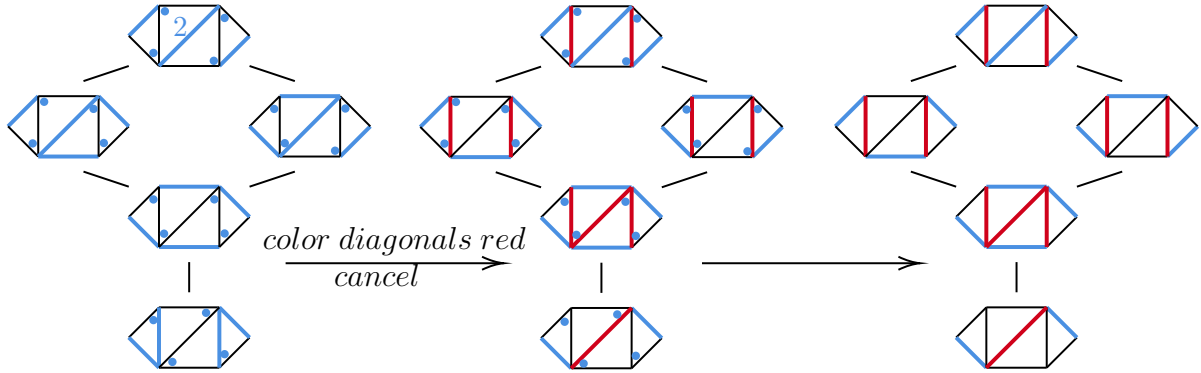

    \centering
    \caption{The map $\mathbb{A}_{ab} \longrightarrow \mathbb{T}_{ab}$ via diagonal coloring and cancellation}
    \label{fig:A_to_T}

\tikzset{every picture/.style={line width=0.75pt}} 

%

\end {figure}

\begin{ex}
   Figure \ref{fig:min_bijection} shows explicitly the application of the composition of bijections $\mathbb{P}_{ab} \longrightarrow \mathbb{A}_{ab} \longrightarrow \mathbb{T}_{ab} \longrightarrow \mathbb{P}_{ab}$ to the input $P_{-} \in \mathbb{P}_{ab}$.
\end{ex}

\begin {figure}[h!]
    \centering
    \caption{The composition $\mathbb{P}_{ab} \longrightarrow \mathbb{A}_{ab} \longrightarrow \mathbb{T}_{ab} \longrightarrow \mathbb{P}_{ab}$ restricted to minimal elements}
    \label{fig:min_bijection}
    %
  
\end {figure}

\chapter {Dual Combinatorial Constructions} \label{dual_const_section}

\section {Words}

For any $w_i \in \{ a,b \}$, let $w_{i}^{*} \in \{ a,b \}$ be the image of $w_i$ under the involution $a \longleftrightarrow b.$

\begin{defn}
     Let $w = w_1 w_2 \dots w_{n-1}$ be a word of length $n-1$. The \textit{dual word} $w^{*}$ is the word of length $n-1$ defined by $w^{*} = w_{1}^{*} w_2 w_{3}^{*} w_4 w_{5}^{*}  \cdots$
\end{defn}

\begin{ex}
   The dual of the word $w=ab$ is $w^{*} = a^{*}b = bb.$
\end{ex}

\begin{rmk}
    The dual of a straight word $w$ is a zigzag word $w^{*}$ and vice versa. 
\end{rmk}

\section{Type $A_n$ Dynkin Quivers}

\begin{defn}
     Let $A_w$ be the Dynkin quiver associated to the word $w$. The \textit{dual Dynkin quiver} $A_{w}^{*}$ is obtained by reversing the orientation of every other edge of $A_{w}$, starting with the first.
\end{defn}

The next result is clear from the definitions.

\begin{prop}
    For any word $w$, we have $A_{w}^{*} = A_{w^{*}}.$
\end{prop}

\begin {figure}[h!]
    \centering
    \caption{The quiver $A_{ab}$ and its dual $A_{bb}$}
    \label{fig:A_ab_dual}
    \begin{tikzcd}
1 & 2 \arrow[l, "a", swap] \arrow[r, "b"] & 3 && 1 \arrow[r, "b"] & 2 \arrow[r, "b"] & 3
\end{tikzcd}

\end {figure}
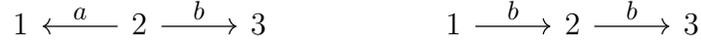

\begin{ex}
    In Figure \ref{fig:A_ab_dual}, the Dynkin quiver $A_{ab}$ is shown on the left, and the dual Dynkin quiver $A_{w}^{*} = A_{w^{*}} = A_{bb}$ is shown on the right. 
\end{ex}

\section {Posets}

\begin{defn}
     Recall the poset $C_w$ associated to the word $w$. Define the \textit{orientation} of an edge in $C_{w}$ to be either NW or NE, according to whether it is labeled by $a$ or $b$, respectively. The \textit{dual poset} $C_{w}^{*}$ is defined by changing the orientation of every other edge of $C_w$, starting with the first.
\end{defn}

\begin{prop}
     For any word $w$, we have $C_{w}^{*} = C_{w^{*}}.$
\end{prop}

\begin {figure}[h!]
    \centering
    \caption{The poset $C_{ab}$ and its dual $C_{bb}$}
    \label{fig:C_ab_dual}
    \begin{tikzcd}
1 && 3 &&&& 3  \\
& 2 \arrow[ul, "a"] \arrow[ur, "b", swap] &&&& 2 \arrow[ur, "b", swap] \\
&&&& 1 \arrow[ur, "b" , swap]
\end{tikzcd}
  
\end {figure}
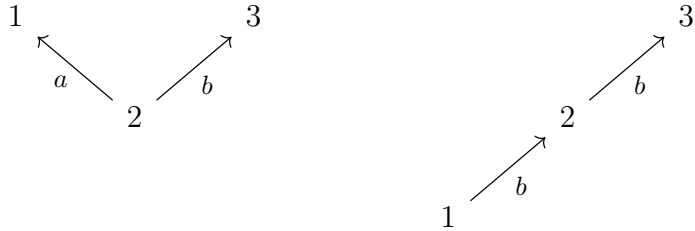

\begin{ex}
    The leftmost poset in Figure \ref{fig:C_ab_dual} is the poset $C_{ab}$ (see Figure \ref{fig:C_ab}), and on the right is the dual poset $C_{w}^{*} =C_{w^{*}} = C_{bb}$. 
\end{ex}

\begin{rmk}
    If $w$ is straight then $C_w$ is isomorphic to a linear chain, and the dual poset $C_{w}^{*}$ is isomorphic to a fence.
\end{rmk}

\section {Triangulations} 

Recall the notation $\Sigma_w = [\Delta_0 , \Delta_1 , \cdots , \Delta_n ]$, where $\Delta_i$ are the ideal triangles (with edge labels from Definition \ref{Delta_w}) cut out by the triangulation $\Delta_w$ of $\Sigma$. Let $\nabla_{i}$ be the edge-labeled triangle with the same positive integer labels as $\Delta_{i}$ but with opposite orientation. Define the \textit{triangle map} by the assignment $$\Sigma_w = [\Delta_0 , \Delta_1 , \dots , \Delta_n]  \mapsto [\Delta_0 , \nabla_1 , \Delta_2 , \nabla_3 , \Delta_4 ,  \dots].$$ 
Define the image of this map to be $\Sigma_{w}^{*}.$

\begin{defn} 
     Consider the triangulation $\Delta_w$ of $\Sigma$ associated to $w$. The \textit{dual triangulation $\Delta_{w}^{*}$ of $\Sigma$} is the triangulation of $\Sigma$ obtained by application of the triangle map $\Sigma_{w} \mapsto \Sigma_{w}^{*}$.
\end{defn}

\begin {figure}[h!]
    \centering
    \caption{The triangle map applied to $\Sigma_{ab}$ gives $\Sigma_{bb}$}
    \label{fig:triangle_map}

\tikzset{every picture/.style={line width=0.75pt}} 

\begin{tikzpicture}[x=0.75pt,y=0.75pt,yscale=-1,xscale=1]

\draw    (86,187.4) -- (126,187.4) ;
\draw    (126,187.4) -- (146,207.4) ;
\draw    (146,207.4) -- (126,227.4) ;
\draw    (86,187.4) -- (66,207.4) ;
\draw    (66,207.4) -- (86,227.4) ;
\draw    (86,227.4) -- (126,227.4) ;
\draw    (86,187.4) -- (86,227.4) ;
\draw    (126,187.4) -- (126,227.4) ;
\draw    (126,187.4) -- (86,227.4) ;
\draw    (206.6,186.9) -- (206.6,226.9) ;
\draw    (206.6,186.9) -- (186.6,206.9) ;
\draw    (186.6,206.9) -- (206.6,226.9) ;
\draw    (216.6,186.9) -- (216.6,226.9) ;
\draw    (256.6,186.9) -- (216.6,226.9) ;
\draw    (216.6,186.9) -- (256.6,186.9) ;
\draw    (266.6,186.9) -- (226.6,226.9) ;
\draw    (226.6,226.9) -- (266.6,226.9) ;
\draw    (266.6,186.9) -- (266.6,226.9) ;
\draw    (276.6,186.9) -- (276.6,226.9) ;
\draw    (296.6,206.9) -- (276.6,226.9) ;
\draw    (276.6,186.9) -- (296.6,206.9) ;
\draw  [dash pattern={on 0.84pt off 2.51pt}]  (170,230) .. controls (170.5,209) and (182.5,213) .. (194.07,206.74) ;
\draw  [dash pattern={on 0.84pt off 2.51pt}]  (213.41,244.07) .. controls (224.07,234.74) and (222.07,226.74) .. (222.07,214.07) ;
\draw  [dash pattern={on 0.84pt off 2.51pt}]  (262.74,239.41) .. controls (258.07,236.07) and (254.07,226.07) .. (255.41,216.74) ;
\draw  [dash pattern={on 0.84pt off 2.51pt}]  (280,220) .. controls (279,234.75) and (284.5,239.75) .. (290,250) ;
\draw    (379.5,188.1) -- (379.5,228.1) ;
\draw    (379.5,188.1) -- (359.5,208.1) ;
\draw    (359.5,208.1) -- (379.5,228.1) ;
\draw    (479.5,186.9) -- (439.5,226.9) ;
\draw    (439.5,226.9) -- (479.5,226.9) ;
\draw    (479.5,186.9) -- (479.5,226.9) ;
\draw    (389.5,186.9) -- (389.5,226.9) ;
\draw    (389.5,226.9) -- (429.5,226.9) ;
\draw    (389.5,186.9) -- (429.5,226.9) ;
\draw    (491,188) -- (491,228) ;
\draw    (511,208) -- (491,228) ;
\draw    (491,188) -- (511,208) ;
\draw  [dash pattern={on 0.84pt off 2.51pt}]  (502.81,204.58) .. controls (514.5,198.75) and (511,186.25) .. (510,180) ;
\draw  [dash pattern={on 0.84pt off 2.51pt}]  (399.65,205.98) .. controls (395.65,196.64) and (398.31,189.31) .. (406.31,183.98) ;
\draw    (579.5,169.4) -- (599.5,189.4) ;
\draw    (559.5,189.4) -- (579.5,169.4) ;
\draw    (599.5,229.4) -- (599.5,189.4) ;
\draw    (559.5,229.4) -- (559.5,189.4) ;
\draw    (559.5,229.4) -- (579.5,249.4) ;
\draw    (579.5,249.4) -- (599.5,229.4) ;
\draw    (559.5,229.4) -- (579.5,169.4) ;
\draw    (579.5,249.4) -- (579.5,169.4) ;
\draw    (599.5,229.4) -- (579.5,169.4) ;
\draw    (310,210) .. controls (311.67,208.33) and (313.33,208.33) .. (315,210) .. controls (316.67,211.67) and (318.33,211.67) .. (320,210) .. controls (321.67,208.33) and (323.33,208.33) .. (325,210) .. controls (326.67,211.67) and (328.33,211.67) .. (330,210) .. controls (331.67,208.33) and (333.33,208.33) .. (335,210) .. controls (336.67,211.67) and (338.33,211.67) .. (340,210) -- (348,210) ;
\draw [shift={(350,210)}, rotate = 180] [color={rgb, 255:red, 0; green, 0; blue, 0 }  ][line width=0.75]    (10.93,-3.29) .. controls (6.95,-1.4) and (3.31,-0.3) .. (0,0) .. controls (3.31,0.3) and (6.95,1.4) .. (10.93,3.29)   ;
\draw  [dash pattern={on 0.84pt off 2.51pt}]  (342.93,231) .. controls (343.43,210) and (355.43,214) .. (367,207.74) ;
\draw  [dash pattern={on 0.84pt off 2.51pt}]  (476.74,240.67) .. controls (472.07,237.33) and (468.07,227.33) .. (469.41,218) ;

\draw (82.6,207.8) node  [font=\scriptsize]  {$1$};
\draw (104.6,201.1) node  [font=\scriptsize]  {$2$};
\draw (121,209) node  [font=\scriptsize]  {$3$};
\draw (105.9,182) node  [font=\scriptsize]  {$4$};
\draw (107,233) node  [font=\scriptsize]  {$5$};
\draw (71,220.3) node  [font=\scriptsize]  {$6$};
\draw (73.4,193.8) node  [font=\scriptsize]  {$7$};
\draw (139.4,221) node  [font=\scriptsize]  {$8$};
\draw (139.4,191.4) node  [font=\scriptsize]  {$9$};
\draw (39.5,208.4) node    {$\Sigma _{w} =$};
\draw (195,191.3) node  [font=\scriptsize]  {$7$};
\draw (192.6,221.8) node  [font=\scriptsize]  {$6$};
\draw (236,180) node  [font=\scriptsize]  {$4$};
\draw (289.5,191.4) node  [font=\scriptsize]  {$9$};
\draw (290,222) node  [font=\scriptsize]  {$8$};
\draw (234.2,201.6) node  [font=\scriptsize]  {$2$};
\draw (248.7,213.6) node  [font=\scriptsize]  {$2$};
\draw (245.1,233) node  [font=\scriptsize]  {$5$};
\draw (262.1,209) node  [font=\scriptsize]  {$3$};
\draw (279.6,209) node  [font=\scriptsize]  {$3$};
\draw (203.7,206.3) node  [font=\scriptsize]  {$1$};
\draw (221.2,206.3) node  [font=\scriptsize]  {$1$};
\draw (168.5,234) node    {$\Delta _{0}$};
\draw (211.5,249) node    {$\Delta _{1}$};
\draw (267.5,246) node    {$\Delta _{2}$};
\draw (298.5,253) node    {$\Delta _{3}$};
\draw (367.9,192.5) node  [font=\scriptsize]  {$7$};
\draw (366.5,223) node  [font=\scriptsize]  {$6$};
\draw (375.6,207.5) node  [font=\scriptsize]  {$1$};
\draw (461.6,213.6) node  [font=\scriptsize]  {$2$};
\draw (463,233) node  [font=\scriptsize]  {$5$};
\draw (475,209) node  [font=\scriptsize]  {$3$};
\draw (393.83,209.23) node  [font=\scriptsize]  {$1$};
\draw (409.67,233) node  [font=\scriptsize]  {$4$};
\draw (410.93,199.4) node  [font=\scriptsize]  {$2$};
\draw (505.5,223) node  [font=\scriptsize]  {$9$};
\draw (505.5,192.5) node  [font=\scriptsize]  {$8$};
\draw (495,210.1) node  [font=\scriptsize]  {$3$};
\draw (411.5,174.9) node    {$\nabla _{1}$};
\draw (511.5,171) node    {$\nabla _{3}$};
\draw (594.5,176.4) node  [font=\scriptsize]  {$8$};
\draw (604,205.9) node  [font=\scriptsize]  {$9$};
\draw (592,193) node  [font=\scriptsize]  {$3$};
\draw (583.1,221.6) node  [font=\scriptsize]  {$2$};
\draw (566.67,194.57) node  [font=\scriptsize]  {$1$};
\draw (566.4,175.4) node  [font=\scriptsize]  {$7$};
\draw (555,203.9) node  [font=\scriptsize]  {$6$};
\draw (566.67,244) node  [font=\scriptsize]  {$4$};
\draw (594,243) node  [font=\scriptsize]  {$5$};
\draw (629.5,208.9) node    {$=\Sigma ^{*}_{w}$};
\draw (165,206) node    {$=$};
\draw (533,207) node    {$=$};
\draw (341.43,235) node    {$\Delta _{0}$};
\draw (481.5,247.26) node    {$\Delta _{2}$};
\draw (336,186) node    {$ \begin{array}{l}
triangle\\
map
\end{array}$};

\end{tikzpicture}

\end {figure}
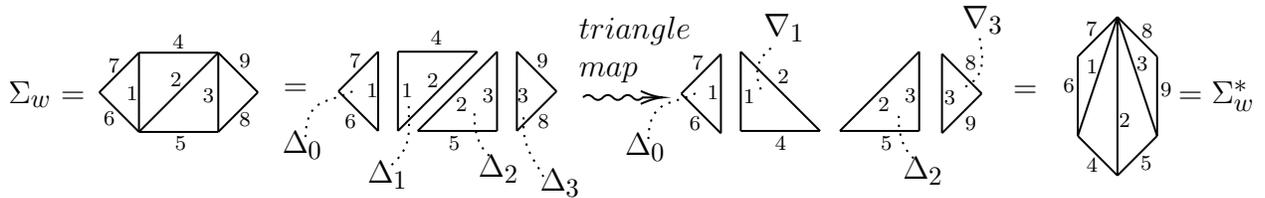

\begin{ex}
    Fix $w=ab$. Figure \ref{fig:triangle_map} shows how the dual triangulation $\Delta_{w}^{*} = \Delta_{ab}^{*} = \Delta_{bb}$ is built by applying the triangle map to $\Sigma_{ab}.$ 
\end{ex}

The next result again follows from the constructions given thus far.

\begin{prop}
     For any word $w,$ we have $\Sigma_{w}^{*} = \Sigma_{w^{*}}$ and $\Delta_{w}^{*} = \Delta_{w^{*}}.$ 
\end{prop}

\begin{rmk}
    If $w$ is straight then $\Delta_w$ is a fan triangulation and the dual $\Delta_{w}^{*}$ is a zigzag triangulation. Conversely, if $w$ is zigzag then $\Delta_w$ is a zigzag triangulation and the dual $\Delta_{w}^{*}$ is a fan triangulation.
\end{rmk}

Application of the triangle map to $\Sigma_w$ gives a new seed for a cluster algebra $\mathcal{A} (\Sigma)_{w^{*}}$ isomorphic to $\mathcal{A} (\Sigma)_w$. If the new initial variable $y_i$ is attached to the arc labeled $i$ in $\Sigma_{w}^{*}$, we relabel this arc with the variable $x_i$. This combinatorial relabeling is introduced so that when we compute cluster variables attached to arcs in the dual triangulation, the result is a Laurent monomial in the initial cluster variables from the original seed.

\section {Slides, Arcs, and Cluster Variables}

Fix a word $w$ and the associated arc $\gamma = \gamma_{a \rightarrow b}$ in $\Sigma_w$.  Recall the internal diagonals of $\Delta_w$ are $\delta_{1} , \delta_{2} , \dots , \delta_{n}$. For $1 \leq i \leq n$, let $p_i = \gamma_{w} \cap \delta_{i}$ be the intersection point of $\gamma_{w}$ with the $i^{th}$ internal diagonal $\delta_{i}$ of $\Delta_w$. Set $p_0 = a$ and $p_{n+1} = b$. Choose a point $m_i \in \text{Int}(\Delta_i) \cap \gamma$ for each $0 \leq i \leq n.$ Let $\gamma_{i}$ be the portion of the arc $\gamma_w$ strictly between $m_{i-1}$ and $m_i$. Let $\gamma_{i}^{left}$ be the portion of $\gamma_i$ strictly between $m_{i-1}$ and $p_i$, and $\gamma_{i}^{right}$ the portion of $\gamma_i$ strictly between $p_i$ and $m_{i}$.

\begin{defn}
     For $1 \leq i \leq n,$ define a \textit{slide of $\gamma_w$ at $p_i$} by the following two-step process. 

\begin{enumerate}[(1)]
    \item Perform the smooth isotopy that fixes $\gamma_w - \gamma_i$ and sends $p_i$ to one of the endpoints of $\delta_{i}$ such that the images of $\gamma_{i}^{left}$ and $\gamma_{i}^{right}$ do not intersect any arc, and have no self-intersections. 
    
    \item Delete the diagonal $\delta_{i}$. 
\end{enumerate}
\end{defn}

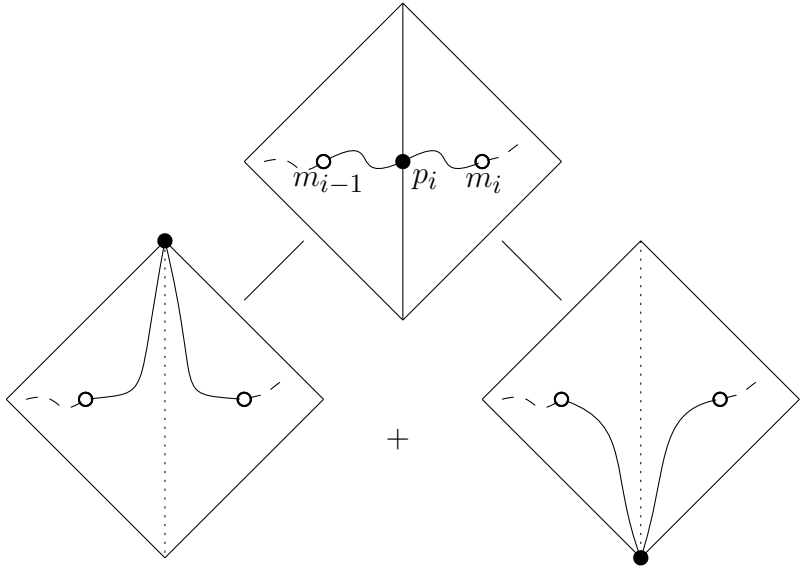
\begin {figure}[h!]
    \centering
    \caption{Dual resolution of the intersection point $p_i$}
    \label{fig:smoothing_dual}
    \begin{tikzpicture}[x=0.75pt,y=0.75pt,yscale=-1,xscale=1]

\draw    (320,10) -- (240,90) ;
\draw    (400,90) -- (320,170) ;
\draw    (320,10) -- (400,90) ;
\draw    (240,90) -- (320,170) ;
\draw    (200,130) -- (120,210) ;
\draw    (280,210) -- (200,290) ;
\draw    (200,130) -- (280,210) ;
\draw    (120,210) -- (200,290) ;
\draw    (440,130) -- (360,210) ;
\draw    (520,210) -- (440,290) ;
\draw    (440,130) -- (520,210) ;
\draw    (360,210) -- (440,290) ;
\draw    (320,10) -- (320,170) ;
\draw    (282.32,88.78) .. controls (311.16,74.02) and (290.51,103.16) .. (320,90) ;
\draw [shift={(280,90)}, rotate = 331.51] [color={rgb, 255:red, 0; green, 0; blue, 0 }  ][line width=0.75]      (0, 0) circle [x radius= 3.35, y radius= 3.35]   ;
\draw    (320,90) .. controls (351.61,72.85) and (330.63,102.28) .. (358.25,90.76) ;
\draw [shift={(360,90)}, rotate = 335.95] [color={rgb, 255:red, 0; green, 0; blue, 0 }  ][line width=0.75]      (0, 0) circle [x radius= 3.35, y radius= 3.35]   ;
\draw [shift={(320,90)}, rotate = 331.51] [color={rgb, 255:red, 0; green, 0; blue, 0 }  ][fill={rgb, 255:red, 0; green, 0; blue, 0 }  ][line width=0.75]      (0, 0) circle [x radius= 3.35, y radius= 3.35]   ;
\draw  [dash pattern={on 4.5pt off 4.5pt}]  (250,90) .. controls (269.06,84.69) and (263.21,100.8) .. (278.28,91.15) ;
\draw [shift={(280,90)}, rotate = 325.18] [color={rgb, 255:red, 0; green, 0; blue, 0 }  ][line width=0.75]      (0, 0) circle [x radius= 3.35, y radius= 3.35]   ;
\draw  [dash pattern={on 4.5pt off 4.5pt}]  (362.45,89.33) .. controls (378.29,85.14) and (363.53,91.46) .. (380,80) ;
\draw [shift={(360,90)}, rotate = 344.44] [color={rgb, 255:red, 0; green, 0; blue, 0 }  ][line width=0.75]      (0, 0) circle [x radius= 3.35, y radius= 3.35]   ;
\draw  [dash pattern={on 0.84pt off 2.51pt}]  (200,130) -- (200,290) ;
\draw    (162.51,209.72) .. controls (194.02,206.29) and (184.89,210.93) .. (200,130) ;
\draw [shift={(160,210)}, rotate = 353.39] [color={rgb, 255:red, 0; green, 0; blue, 0 }  ][line width=0.75]      (0, 0) circle [x radius= 3.35, y radius= 3.35]   ;
\draw    (200,130) .. controls (219.11,209.38) and (202.21,207.13) .. (237.76,209.83) ;
\draw [shift={(240,210)}, rotate = 4.46] [color={rgb, 255:red, 0; green, 0; blue, 0 }  ][line width=0.75]      (0, 0) circle [x radius= 3.35, y radius= 3.35]   ;
\draw [shift={(200,130)}, rotate = 76.46] [color={rgb, 255:red, 0; green, 0; blue, 0 }  ][fill={rgb, 255:red, 0; green, 0; blue, 0 }  ][line width=0.75]      (0, 0) circle [x radius= 3.35, y radius= 3.35]   ;
\draw  [dash pattern={on 4.5pt off 4.5pt}]  (130,210) .. controls (149.06,204.69) and (143.21,220.8) .. (158.28,211.15) ;
\draw [shift={(160,210)}, rotate = 325.18] [color={rgb, 255:red, 0; green, 0; blue, 0 }  ][line width=0.75]      (0, 0) circle [x radius= 3.35, y radius= 3.35]   ;
\draw  [dash pattern={on 4.5pt off 4.5pt}]  (242.45,209.33) .. controls (258.29,205.14) and (243.53,211.46) .. (260,200) ;
\draw [shift={(240,210)}, rotate = 344.44] [color={rgb, 255:red, 0; green, 0; blue, 0 }  ][line width=0.75]      (0, 0) circle [x radius= 3.35, y radius= 3.35]   ;
\draw  [dash pattern={on 0.84pt off 2.51pt}]  (440,130) -- (440,290) ;
\draw    (402.65,210.92) .. controls (435.87,223.21) and (423.91,247.1) .. (440,290) ;
\draw [shift={(400,210)}, rotate = 18.2] [color={rgb, 255:red, 0; green, 0; blue, 0 }  ][line width=0.75]      (0, 0) circle [x radius= 3.35, y radius= 3.35]   ;
\draw    (440,290) .. controls (454.21,246.88) and (446.81,216.25) .. (478.04,210.33) ;
\draw [shift={(480,210)}, rotate = 351.51] [color={rgb, 255:red, 0; green, 0; blue, 0 }  ][line width=0.75]      (0, 0) circle [x radius= 3.35, y radius= 3.35]   ;
\draw [shift={(440,290)}, rotate = 288.24] [color={rgb, 255:red, 0; green, 0; blue, 0 }  ][fill={rgb, 255:red, 0; green, 0; blue, 0 }  ][line width=0.75]      (0, 0) circle [x radius= 3.35, y radius= 3.35]   ;
\draw  [dash pattern={on 4.5pt off 4.5pt}]  (370,210) .. controls (389.06,204.69) and (383.21,220.8) .. (398.28,211.15) ;
\draw [shift={(400,210)}, rotate = 325.18] [color={rgb, 255:red, 0; green, 0; blue, 0 }  ][line width=0.75]      (0, 0) circle [x radius= 3.35, y radius= 3.35]   ;
\draw  [dash pattern={on 4.5pt off 4.5pt}]  (482.45,209.33) .. controls (498.29,205.14) and (483.53,211.46) .. (500,200) ;
\draw [shift={(480,210)}, rotate = 344.44] [color={rgb, 255:red, 0; green, 0; blue, 0 }  ][line width=0.75]      (0, 0) circle [x radius= 3.35, y radius= 3.35]   ;
\draw    (270,130) -- (240,160) ;
\draw    (370,130) -- (400,160) ;

\draw (317.5,230) node    {$+$};
\draw (331.5,99) node    {$p_{i}$};
\draw (282.5,101) node    {$m_{i-1}$};
\draw (361,101) node    {$m_{i}$};

\end{tikzpicture}
  
\end {figure}

Choosing one of the two possible slides at each $p_i$ results in a collection of nonintersecting curves in $\Sigma$. Note that closed curves based at some $v \in \Sigma$ are the only loops that can occur. Replace each curve with distinct endpoints by the arc or boundary segment from $\Delta_w$ with the same endpoints, and replace each closed curve with $\varnothing$. See Figure \ref{fig:smoothing_dual} below.

\begin{defn}
     The set of \textit{dual resolutions $\text{Res}(w)^{*}$ associated to $w$} has as its elements those diagrams that can be obtained from sliding each $p_i$ in one of the two possible ways, in some order. For $r^{*} \in \text{Res}(w)^{*}$, let $E(r^{*})$ be the collection of arcs and boundary segments from $\Delta_w$ produced from the resolution $r^{*}$, along with $\varnothing$ if any closed loops are present. Define the \textit{weight} of any dual resolution $r^{*}$ to be $x_{r^{*}} = \prod_{j \in E(r^{*})} x_{j},$ where $x_{\varnothing} = 0.$
\end{defn}

We now describe how to produce a \textit{dual resolution tree} from $w$. Each node of such a tree is a diagram of arcs inside the $(n+3)$-gon $\Sigma$, and is weighted by the product of cluster variables associated to those arcs, or zero if there is a closed loop in the diagram. The root of a dual resolution tree from $w$ is the diagram consisting of the arc $\gamma_w$ inside $\Sigma_w$. Choosing an intersection point $p_i$ to slide at creates two children of this root (see Figure \ref{fig:smoothing_dual}). Continuing in this way (choosing an intersection point to slide at in each child, etc.) and halting whenever either we create a loop or we have performed a slide at every intersection point, a binary tree (with additional node structure) is produced.

\begin{defn}
     The set of \textit{dual resolution trees $\text{Tree}(w)^{*}$ associated to $w$} is the set whose elements are dual slide trees from $w$ as described above.
\end{defn}

\begin{rmk}
    The set $\text{Res}(w)^{*}$ is equal to the union of the leaves of the trees in $\text{Tree}(w)^{*}.$
\end{rmk}

\begin {figure}[h!]
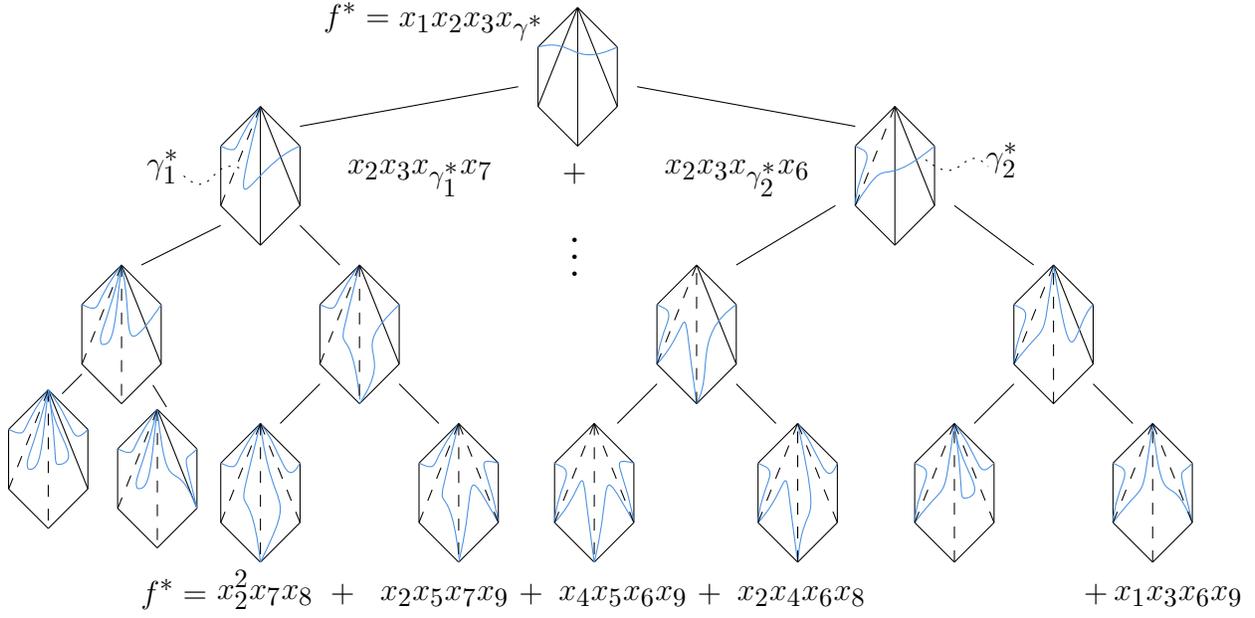

    \centering
    \caption{One element of $\text{Tree}(bb)^{*}$}
    \label{fig:tree_bb}


\end {figure}

\begin{ex}
    Fix $w=ab$. Figure \ref{fig:tree_bb} shows one element of $\text{Tree}(w^{*})^{*} = \text{Tree}(bb)^{*}$. Note that this tree is isomorphic to the element of $\text{Tree}(w) = \text{Tree}(ab)$ from Example \ref{fig:tree_ab}, and that the weights of the leaves here coincide with the weights of the leaves there.
\end{ex}

\begin{defn}
      Let $a^{*}$ and $b^{*}$ be the images of $a$ and $b$ under the triangle map $\Sigma_w \mapsto \Sigma_{w}^{*}$. The \textit{dual of the arc $\gamma$} is the oriented arc $\gamma_{w^{*}} = \gamma^{*}_{a^{*} \rightarrow b^{*}}$ from $a^{*}$ to $b^{*}$ inside the polygon $\Sigma_{w}^{*}$. The \textit{cluster variable $x_{w}^{*}$ dual to $x_{w}$} is defined by 
$$x_{w}^{*} = \frac{1}{x_1 x_2 ... x_n}\sum_{r^* \in \text{Res}(w)^{*}} x_{r^{*}}.$$ 
\end{defn}

\begin{rmk}
    We caution that in general the arc $\gamma_{w}^{*}$ is not equal to the arc $\gamma_{w}$. Furthermore, we do not yet know that $x_{w}^{*}$ is in fact a cluster variable in a cluster algebra; this is part (c) in Theorem \ref{duality_thm} below. 
\end{rmk}

\section {Snake Graphs}

The notion of a dual snake graph was introduced in \cite{propp2005combinatorics}. 

\begin{defn} \label{dual_snake}
      Fix the arc $\gamma$ in the triangulated polygon $\Sigma_{w}$ triangulated by $\Delta_{w}$. The \textit{dual snake graph $G_{w}^{*}$ associated to $w$} is the labeled planar graph recursively defined as follows: 

\begin{enumerate}[1.]
    \item Choose an orientation-preserving embedding of the triangulated square $[\Delta_0 , \nabla_1 ]$ into the discrete plane $\mathbb{Z}^2$ such that its image $\widetilde{T}_{1}^{*}$ is a triangulated unit square  with vertices $(0,0)$, $(1,0)$, $(0,1),$ and $(1,1)$ in $\mathbb{Z}^2$, and such that the point $a \in \Delta_0$ maps to the point $(0,0)$. Note that the (line spanned by the) image of the triangulating edge will have slope $-1$.
    
    \item Choose an orientation-preserving map of $[\Delta_{1} , \nabla_{2} ]$ into $\mathbb{Z}^2$ such that its image $\widetilde{T}_{2}^{*}$ is a triangulated unit square (again, with triangulating edge having slope $-1$) glued to $\widetilde{T}_{1}^{*}$ along the unique edge in each $\widetilde{T_{i}^{*}}$ labeled $j \in \{ n+1 , \dots , 2n+3 \}.$ Note that if the intersection point of the diagonals $\delta_{1}$ and $\delta_{2}$ is to the left (resp. right) of $\gamma_w$, then $\widetilde{T_{2}^{*}}$ is the triangulated square directly above (resp. to the right of) $\widetilde{T}_{1}^{*}$.
    
    \item Continue this process, using orientation-preserving maps for both $i$ odd and  even, to get the graph $\widetilde{G_{w}^{*}}$, built from triangulated unit squares in $\mathbb{Z}^2$ (with all triangulating edges having slope $-1$) glued either above or to the right of the previous square. Each $\widetilde{T_{i}^{*}}$ will be called a \textit{tile} of $\widetilde{G_{w}^{*}}$. The triangulating edge of each $\widetilde{T_{i}^{*}}$ is called the \textit{diagonal} of $\widetilde{T_{i}^{*}}.$
    
    \item The \textit{dual snake graph} $G_{w}^{*}$ is the graph in $\mathbb{Z}^2$ gotten by deleting each diagonal from each tile in $\widetilde{G_{w}^{*}}$. 
\end{enumerate}
\end{defn}

\begin {figure}[h!]
    \centering
    \caption{Construction of the dual snake graph $G_{bb}$}
    \label{fig:G_bb}
    \begin{tikzpicture}[x=0.75pt,y=0.75pt,yscale=-1,xscale=1]

\draw    (130.43,53.86) -- (179,102.43) ;
\draw    (179,102.43) -- (130.43,151) ;
\draw [color={rgb, 255:red, 0; green, 0; blue, 0 }  ,draw opacity=1 ][line width=1.5]    (57.57,53.86) -- (130.43,53.86) ;
\draw    (57.57,151) -- (130.43,151) ;
\draw [color={rgb, 255:red, 0; green, 0; blue, 0 }  ,draw opacity=1 ][line width=1.5]    (57.57,53.86) -- (9,102.43) ;
\draw [color={rgb, 255:red, 0; green, 0; blue, 0 }  ,draw opacity=1 ][line width=1.5]    (9,102.43) -- (57.57,151) ;
\draw [color={rgb, 255:red, 0; green, 0; blue, 0 }  ,draw opacity=1 ][line width=1.5]    (57.57,53.86) -- (57.57,151) ;
\draw    (130.43,53.86) -- (130.43,151) ;
\draw [color={rgb, 255:red, 0; green, 0; blue, 0 }  ,draw opacity=1 ][line width=1.5]    (130.43,53.86) -- (57.57,151) ;
\draw [color={rgb, 255:red, 0; green, 0; blue, 0 }  ,draw opacity=1 ] [dash pattern={on 0.84pt off 2.51pt}]  (9,102.43) -- (179,102.43) ;
\draw    (130.43,163.86) -- (179,212.43) ;
\draw    (179,212.43) -- (130.43,261) ;
\draw [color={rgb, 255:red, 0; green, 0; blue, 0 }  ,draw opacity=1 ][line width=1.5]    (57.57,163.86) -- (130.43,163.86) ;
\draw [color={rgb, 255:red, 0; green, 0; blue, 0 }  ,draw opacity=1 ][line width=1.5]    (57.57,261) -- (130.43,261) ;
\draw    (57.57,163.86) -- (9,212.43) ;
\draw    (9,212.43) -- (57.57,261) ;
\draw [color={rgb, 255:red, 0; green, 0; blue, 0 }  ,draw opacity=1 ][line width=1.5]    (57.57,163.86) -- (57.57,261) ;
\draw [color={rgb, 255:red, 0; green, 0; blue, 0 }  ,draw opacity=1 ][line width=1.5]    (130.43,163.86) -- (130.43,261) ;
\draw [color={rgb, 255:red, 0; green, 0; blue, 0 }  ,draw opacity=1 ][line width=1.5]    (130.43,163.86) -- (57.57,261) ;
\draw [color={rgb, 255:red, 0; green, 0; blue, 0 }  ,draw opacity=1 ] [dash pattern={on 0.84pt off 2.51pt}]  (9,212.43) -- (179,212.43) ;
\draw [color={rgb, 255:red, 0; green, 0; blue, 0 }  ,draw opacity=1 ][line width=1.5]    (130.43,273.86) -- (179,322.43) ;
\draw [color={rgb, 255:red, 0; green, 0; blue, 0 }  ,draw opacity=1 ][line width=1.5]    (179,322.43) -- (130.43,371) ;
\draw    (57.57,273.86) -- (130.43,273.86) ;
\draw [color={rgb, 255:red, 0; green, 0; blue, 0 }  ,draw opacity=1 ][line width=1.5]    (57.57,371) -- (130.43,371) ;
\draw    (57.57,273.86) -- (9,322.43) ;
\draw    (9,322.43) -- (57.57,371) ;
\draw    (57.57,273.86) -- (57.57,371) ;
\draw [color={rgb, 255:red, 0; green, 0; blue, 0 }  ,draw opacity=1 ][line width=1.5]    (130.43,273.86) -- (130.43,371) ;
\draw [color={rgb, 255:red, 0; green, 0; blue, 0 }  ,draw opacity=1 ][line width=1.5]    (130.43,273.86) -- (57.57,371) ;
\draw [color={rgb, 255:red, 0; green, 0; blue, 0 }  ,draw opacity=1 ] [dash pattern={on 0.84pt off 2.51pt}]  (9,322.43) -- (179,322.43) ;
\draw    (329,41) -- (329,141) ;
\draw    (389,101) -- (269,101) ;
\draw [color={rgb, 255:red, 0; green, 0; blue, 0 }  ,draw opacity=1 ][line width=1.5]    (329,81) -- (329,101) ;
\draw [color={rgb, 255:red, 208; green, 2; blue, 27 }  ,draw opacity=1 ][line width=1.5]    (349,81) -- (349,101) ;
\draw [color={rgb, 255:red, 0; green, 0; blue, 0 }  ,draw opacity=1 ][line width=1.5]    (349,101) -- (329,101) ;
\draw [color={rgb, 255:red, 208; green, 2; blue, 27 }  ,draw opacity=1 ][line width=1.5]    (349,81) -- (329,81) ;
\draw [color={rgb, 255:red, 208; green, 2; blue, 27 }  ,draw opacity=1 ]   (190,120) -- (278,120) ;
\draw [shift={(280,120)}, rotate = 180] [color={rgb, 255:red, 208; green, 2; blue, 27 }  ,draw opacity=1 ][line width=0.75]    (10.93,-3.29) .. controls (6.95,-1.4) and (3.31,-0.3) .. (0,0) .. controls (3.31,0.3) and (6.95,1.4) .. (10.93,3.29)   ;
\draw [shift={(190,120)}, rotate = 0] [color={rgb, 255:red, 208; green, 2; blue, 27 }  ,draw opacity=1 ][line width=0.75]      (0,-11.18) .. controls (-3.09,-11.18) and (-5.59,-8.68) .. (-5.59,-5.59) .. controls (-5.59,-2.5) and (-3.09,0) .. (0,0) ;
\draw    (329,151) -- (329,251) ;
\draw    (389,211) -- (269,211) ;
\draw [color={rgb, 255:red, 0; green, 0; blue, 0 }  ,draw opacity=1 ][line width=0.75]    (329,191) -- (329,211) ;
\draw [color={rgb, 255:red, 0; green, 0; blue, 0 }  ,draw opacity=1 ][line width=0.75]    (330,190) -- (350,190) ;
\draw [color={rgb, 255:red, 0; green, 0; blue, 0 }  ,draw opacity=1 ][line width=0.75]    (349,211) -- (329,211) ;
\draw [color={rgb, 255:red, 0; green, 0; blue, 0 }  ,draw opacity=1 ][line width=1.5]    (370,210) -- (350,210) ;
\draw [color={rgb, 255:red, 208; green, 2; blue, 27 }  ,draw opacity=1 ][line width=1.5]    (370,190) -- (350,190) ;
\draw [color={rgb, 255:red, 208; green, 2; blue, 27 }  ,draw opacity=1 ][line width=1.5]    (370,190) -- (370,210) ;
\draw [color={rgb, 255:red, 0; green, 0; blue, 0 }  ,draw opacity=1 ][line width=1.5]    (349,191) -- (349,211) ;
\draw    (329,261) -- (329,361) ;
\draw    (389,321) -- (269,321) ;
\draw [color={rgb, 255:red, 0; green, 0; blue, 0 }  ,draw opacity=1 ][line width=0.75]    (329,301) -- (329,321) ;
\draw [color={rgb, 255:red, 0; green, 0; blue, 0 }  ,draw opacity=1 ][line width=0.75]    (349,301) -- (349,321) ;
\draw [color={rgb, 255:red, 0; green, 0; blue, 0 }  ,draw opacity=1 ][line width=0.75]    (349,321) -- (329,321) ;
\draw [color={rgb, 255:red, 0; green, 0; blue, 0 }  ,draw opacity=1 ][line width=0.75]    (349,301) -- (329,301) ;
\draw [color={rgb, 255:red, 0; green, 0; blue, 0 }  ,draw opacity=1 ][line width=1.5]    (370,300) -- (350,300) ;
\draw [color={rgb, 255:red, 0; green, 0; blue, 0 }  ,draw opacity=1 ][line width=0.75]    (369,301) -- (369,321) ;
\draw [color={rgb, 255:red, 208; green, 2; blue, 27 }  ,draw opacity=1 ][line width=1.5]    (370,280) -- (350,280) ;
\draw [color={rgb, 255:red, 208; green, 2; blue, 27 }  ,draw opacity=1 ][line width=1.5]    (370,280) -- (370,300) ;
\draw [color={rgb, 255:red, 0; green, 0; blue, 0 }  ,draw opacity=1 ][line width=1.5]    (350,280) -- (350,300) ;
\draw [color={rgb, 255:red, 208; green, 2; blue, 27 }  ,draw opacity=1 ][line width=1.5]    (329,81) -- (349,101) ;
\draw [color={rgb, 255:red, 0; green, 0; blue, 0 }  ,draw opacity=1 ][line width=0.75]  [dash pattern={on 0.84pt off 2.51pt}]  (329,191) -- (349,211) ;
\draw [color={rgb, 255:red, 208; green, 2; blue, 27 }  ,draw opacity=1 ][line width=1.5]    (349,191) -- (369,211) ;
\draw [color={rgb, 255:red, 0; green, 0; blue, 0 }  ,draw opacity=1 ][line width=0.75]    (329,301) -- (349,321) ;
\draw [color={rgb, 255:red, 0; green, 0; blue, 0 }  ,draw opacity=1 ][line width=0.75]  [dash pattern={on 0.84pt off 2.51pt}]  (349,301) -- (369,321) ;
\draw [color={rgb, 255:red, 208; green, 2; blue, 27 }  ,draw opacity=1 ][line width=1.5]    (350,280) -- (370,300) ;
\draw [color={rgb, 255:red, 208; green, 2; blue, 27 }  ,draw opacity=1 ]   (189,231) -- (277,231) ;
\draw [shift={(279,231)}, rotate = 180] [color={rgb, 255:red, 208; green, 2; blue, 27 }  ,draw opacity=1 ][line width=0.75]    (10.93,-3.29) .. controls (6.95,-1.4) and (3.31,-0.3) .. (0,0) .. controls (3.31,0.3) and (6.95,1.4) .. (10.93,3.29)   ;
\draw [shift={(189,231)}, rotate = 0] [color={rgb, 255:red, 208; green, 2; blue, 27 }  ,draw opacity=1 ][line width=0.75]      (0,-11.18) .. controls (-3.09,-11.18) and (-5.59,-8.68) .. (-5.59,-5.59) .. controls (-5.59,-2.5) and (-3.09,0) .. (0,0) ;
\draw [color={rgb, 255:red, 208; green, 2; blue, 27 }  ,draw opacity=1 ]   (189,341) -- (277,341) ;
\draw [shift={(279,341)}, rotate = 180] [color={rgb, 255:red, 208; green, 2; blue, 27 }  ,draw opacity=1 ][line width=0.75]    (10.93,-3.29) .. controls (6.95,-1.4) and (3.31,-0.3) .. (0,0) .. controls (3.31,0.3) and (6.95,1.4) .. (10.93,3.29)   ;
\draw [shift={(189,341)}, rotate = 0] [color={rgb, 255:red, 208; green, 2; blue, 27 }  ,draw opacity=1 ][line width=0.75]      (0,-11.18) .. controls (-3.09,-11.18) and (-5.59,-8.68) .. (-5.59,-5.59) .. controls (-5.59,-2.5) and (-3.09,0) .. (0,0) ;
\draw    (386,292) .. controls (386.63,289.73) and (388.08,288.9) .. (390.35,289.53) .. controls (392.62,290.16) and (394.07,289.33) .. (394.7,287.06) .. controls (395.33,284.79) and (396.78,283.97) .. (399.05,284.6) .. controls (401.32,285.23) and (402.77,284.4) .. (403.39,282.13) .. controls (404.02,279.86) and (405.47,279.03) .. (407.74,279.66) .. controls (410.01,280.29) and (411.46,279.46) .. (412.09,277.19) .. controls (412.72,274.92) and (414.17,274.09) .. (416.44,274.72) .. controls (418.71,275.35) and (420.16,274.53) .. (420.79,272.26) .. controls (421.42,269.99) and (422.87,269.16) .. (425.14,269.79) .. controls (427.41,270.42) and (428.86,269.59) .. (429.48,267.32) .. controls (430.11,265.05) and (431.56,264.22) .. (433.83,264.85) .. controls (436.1,265.48) and (437.55,264.65) .. (438.18,262.38) .. controls (438.81,260.11) and (440.26,259.29) .. (442.53,259.92) .. controls (444.8,260.55) and (446.25,259.72) .. (446.88,257.45) .. controls (447.51,255.18) and (448.96,254.35) .. (451.23,254.98) -- (451.3,254.94) -- (458.26,250.99) ;
\draw [shift={(460,250)}, rotate = 510.42] [color={rgb, 255:red, 0; green, 0; blue, 0 }  ][line width=0.75]    (10.93,-3.29) .. controls (6.95,-1.4) and (3.31,-0.3) .. (0,0) .. controls (3.31,0.3) and (6.95,1.4) .. (10.93,3.29)   ;
\draw    (482.8,262.6) -- (562.8,262.6) ;
\draw    (482.8,262.6) -- (482.8,182.6) ;
\draw    (562.8,262.6) -- (562.8,182.6) ;
\draw    (482.8,182.6) -- (562.8,182.6) ;
\draw    (562.8,262.6) -- (642.8,262.6) ;
\draw    (642.8,262.6) -- (642.8,182.6) ;
\draw    (562.8,182.6) -- (642.8,182.6) ;
\draw    (562.8,182.6) -- (562.8,102.6) ;
\draw    (642.8,182.6) -- (642.8,102.6) ;
\draw    (562.8,102.6) -- (642.8,102.6) ;
\draw [color={rgb, 255:red, 0; green, 0; blue, 0 }  ,draw opacity=1 ][line width=0.75]    (349,301) -- (369,301) ;

\draw (225,107) node  [color={rgb, 255:red, 245; green, 166; blue, 35 }  ,opacity=1 ]  {$\textcolor[rgb]{0.82,0.01,0.11}{Tile\ 1}$};
\draw (96,121) node  [color={rgb, 255:red, 0; green, 0; blue, 0 }  ,opacity=1 ]  {$2$};
\draw (26,131) node  [color={rgb, 255:red, 0; green, 0; blue, 0 }  ,opacity=1 ]  {$6$};
\draw (26,71) node  [color={rgb, 255:red, 0; green, 0; blue, 0 }  ,opacity=1 ]  {$7$};
\draw (92,41) node  [color={rgb, 255:red, 0; green, 0; blue, 0 }  ,opacity=1 ]  {$4$};
\draw (52,201) node    {$1$};
\draw (142,201) node    {$3$};
\draw (92,171) node    {$4$};
\draw (96,251) node    {$5$};
\draw (92,311) node    {$2$};
\draw (162,291) node    {$9$};
\draw (92,381) node    {$5$};
\draw (162,351) node    {$8$};
\draw (343.67,200.33) node  [font=\scriptsize,color={rgb, 255:red, 0; green, 0; blue, 0 }  ,opacity=1 ]  {$4$};
\draw (360.33,183) node  [font=\scriptsize,color={rgb, 255:red, 208; green, 2; blue, 27 }  ,opacity=1 ]  {$5$};
\draw (376.33,200.33) node  [font=\scriptsize,color={rgb, 255:red, 208; green, 2; blue, 27 }  ,opacity=1 ]  {$3$};
\draw (360.33,217.67) node  [font=\scriptsize,color={rgb, 255:red, 0; green, 0; blue, 0 }  ,opacity=1 ]  {$1$};
\draw (340.66,74.1) node  [font=\scriptsize,color={rgb, 255:red, 208; green, 2; blue, 27 }  ,opacity=1 ]  {$2$};
\draw (340.46,106.84) node  [font=\scriptsize,color={rgb, 255:red, 0; green, 0; blue, 0 }  ,opacity=1 ]  {$6$};
\draw (323.5,90.1) node  [font=\scriptsize,color={rgb, 255:red, 0; green, 0; blue, 0 }  ,opacity=1 ]  {$7$};
\draw (355.43,90.14) node  [font=\scriptsize,color={rgb, 255:red, 208; green, 2; blue, 27 }  ,opacity=1 ]  {$4$};
\draw (359.71,305.64) node  [font=\scriptsize,color={rgb, 255:red, 0; green, 0; blue, 0 }  ,opacity=1 ]  {$5$};
\draw (360.46,273.75) node  [font=\scriptsize,color={rgb, 255:red, 208; green, 2; blue, 27 }  ,opacity=1 ]  {$8$};
\draw (374.64,288.96) node  [font=\scriptsize,color={rgb, 255:red, 208; green, 2; blue, 27 }  ,opacity=1 ]  {$9$};
\draw (345.61,289.29) node  [font=\scriptsize,color={rgb, 255:red, 0; green, 0; blue, 0 }  ,opacity=1 ]  {$2$};
\draw (66,91) node  [color={rgb, 255:red, 0; green, 0; blue, 0 }  ,opacity=1 ]  {$1$};
\draw (84,208) node    {$2$};
\draw (122,330) node  [color={rgb, 255:red, 65; green, 117; blue, 5 }  ,opacity=1 ]  {$\textcolor[rgb]{0,0,0}{3}$};
\draw (221.5,220) node  [color={rgb, 255:red, 245; green, 166; blue, 35 }  ,opacity=1 ]  {$\textcolor[rgb]{0.82,0.01,0.11}{Tile\ 2}$};
\draw (224,328) node  [color={rgb, 255:red, 245; green, 166; blue, 35 }  ,opacity=1 ]  {$\textcolor[rgb]{0.82,0.01,0.11}{Tile\ 3}$};
\draw (457,300) node    {$ \begin{array}{l}
Remove\\
Diagonals
\end{array}$};
\draw (525.6,270) node    {$6$};
\draw (477,223.6) node    {$7$};
\draw (524.2,172.8) node    {$2$};
\draw (569.6,222.6) node    {$4$};
\draw (602.8,270) node    {$1$};
\draw (650.6,223.8) node    {$3$};
\draw (604.4,172.6) node    {$5$};
\draw (557,142.6) node    {$2$};
\draw (605,93.4) node    {$8$};
\draw (650.8,139.6) node    {$9$};

\end{tikzpicture}
  
\end {figure}
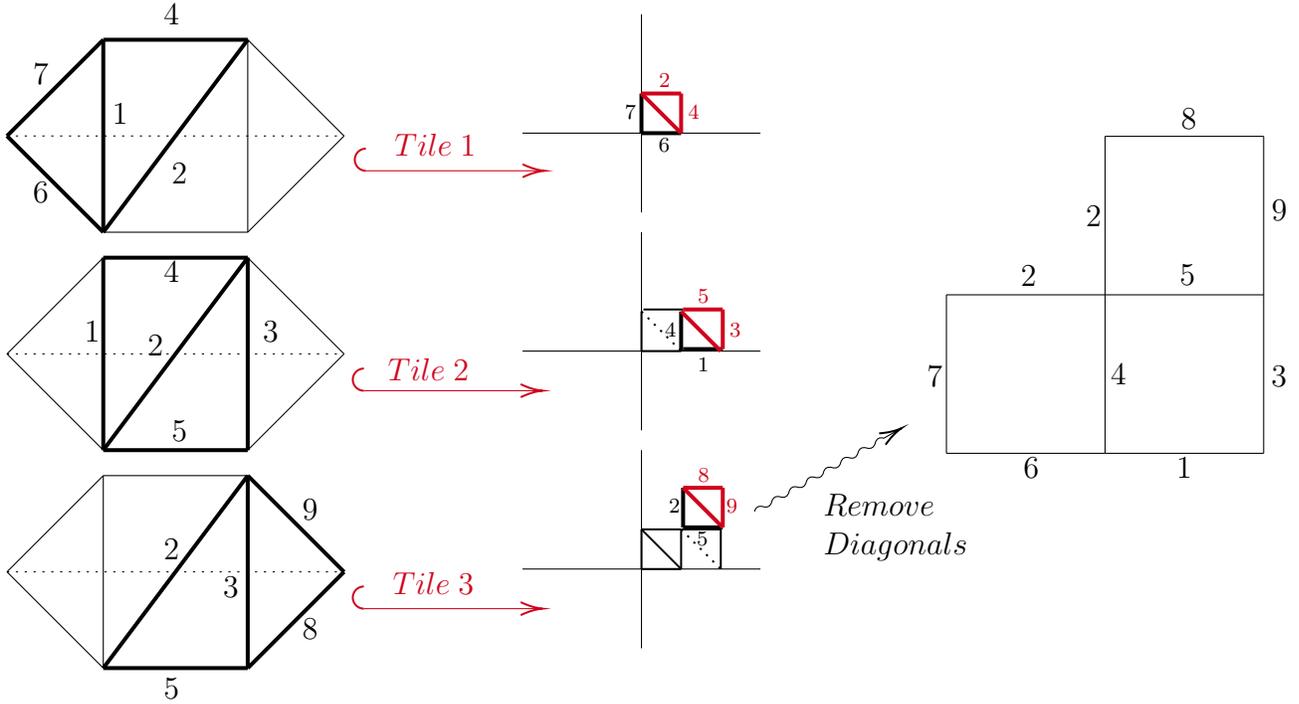

\begin{ex}
    Fix $w=ab$. Figure \ref{fig:G_bb} illustrates the construction of the dual snake graph $G_{w}^{*} = G_{bb}$ from the triangulation $\Delta_{ab}$.
\end{ex}

\newpage

We give now an explicit procedure $G_{w} \mapsto G_{w}^{*}$ for computing the dual snake graph, starting from $G_w$.

\begin{defn}
      Consider the snake graph $G_w$ with tiles $T_1 , T_2 , \dots , T_n$. Let the diagonal of $T_i$ be called $d_i$. The \textit{antidiagonal of tile $T_i$,} denoted $a_i$, is the line segment inside $T_i$ formed by joining the SW and NE corners of $T_i$.
\end{defn}

For any snake graph $H$ with $n$ tiles, the diagonal $d_i$ of $T_i$ gives two subgraphs $H_{i-1}$ and $H_i$ of $H$ that respectively consist of all vertices and edges of $H$ weakly below or weakly above the line spanned by $d_i$. Define $H^{T_i}$ to be the snake graph produced by reflecting $H_i$ about the line spanned by $a_i$ and regluing the image of $H_i$ to $H_{i-1}$. It is clear that the result of this operation is another snake graph. Write $(H^{T_i})^{T_j} = H^{T_{i} \circ T_{j}}.$

One can see from the constructions that we have $G_{w} \mapsto G_{w}^{T_1 \circ T_2 \circ \dots \circ T_n} = G_{w}^{*}.$ Namely, performing this composition of maps gives each tile of $G_w$ a half-twist (as in Definition \ref{dual_snake}), and furthermore the shape of $G_{w}^{T_1 \circ T_2 \circ \dots \circ T_n}$ coincides with the shape of $G_{w}^{*}.$

\begin {figure}[h!]
    \centering
    \caption{Transforming $G_{ab}$ into its dual $G_{bb}$}
    \label{fig:factorization}
    \begin{tikzpicture}[x=0.75pt,y=0.75pt,yscale=-1,xscale=1]

\draw    (20,240) -- (80,240) ;

\draw    (20,180) -- (20,240) ;

\draw    (20,180) -- (80,180) ;

\draw    (80,180) -- (80,240) ;

\draw    (80,120) -- (80,180) ;

\draw    (20,120) -- (20,180) ;

\draw    (20,120) -- (80,120) ;

\draw    (20,60) -- (20,120) ;

\draw    (80,60) -- (80,120) ;

\draw    (20,60) -- (80,60) ;

\draw    (120,240) -- (180,240) ;

\draw    (120,180) -- (120,240) ;

\draw    (120,180) -- (180,180) ;

\draw    (180,180) -- (180,240) ;

\draw    (180,240) -- (240,240) ;

\draw    (180,180) -- (240,180) ;

\draw    (240,180) -- (240,240) ;

\draw    (240,180) -- (300,180) ;

\draw    (240,240) -- (300,240) ;

\draw    (300,180) -- (300,240) ;

\draw    (340,240) -- (400,240) ;

\draw    (340,180) -- (340,240) ;

\draw    (340,180) -- (400,180) ;

\draw    (400,180) -- (400,240) ;

\draw    (400,240) -- (460,240) ;

\draw    (400,180) -- (460,180) ;

\draw    (460,180) -- (460,240) ;

\draw    (400,120) -- (400,180) ;

\draw    (460,120) -- (460,180) ;

\draw    (400,120) -- (460,120) ;

\draw    (500,240) -- (560,240) ;

\draw    (500,180) -- (500,240) ;

\draw    (500,180) -- (560,180) ;

\draw    (560,180) -- (560,240) ;

\draw    (560,240) -- (620,240) ;

\draw    (560,180) -- (620,180) ;

\draw    (620,180) -- (620,240) ;

\draw    (560,120) -- (560,180) ;

\draw    (620,120) -- (620,180) ;

\draw    (560,120) -- (620,120) ;

\draw  [dash pattern={on 0.84pt off 2.51pt}]  (0,160) -- (100,260) ;

\draw  [dash pattern={on 0.84pt off 2.51pt}]  (160,160) -- (260,260) ;

\draw  [dash pattern={on 0.84pt off 2.51pt}]  (380,100) -- (480,200) ;

\draw    (90,274.09) .. controls (118.07,293.99) and (156.33,296.23) .. (188.53,275.07) ;
\draw [shift={(190,274.09)}, rotate = 505.66] [color={rgb, 255:red, 0; green, 0; blue, 0 }  ][line width=0.75]    (10.93,-3.29) .. controls (6.95,-1.4) and (3.31,-0.3) .. (0,0) .. controls (3.31,0.3) and (6.95,1.4) .. (10.93,3.29)   ;

\draw    (270,274.09) .. controls (298.07,293.99) and (336.33,296.23) .. (368.53,275.07) ;
\draw [shift={(370,274.09)}, rotate = 505.66] [color={rgb, 255:red, 0; green, 0; blue, 0 }  ][line width=0.75]    (10.93,-3.29) .. controls (6.95,-1.4) and (3.31,-0.3) .. (0,0) .. controls (3.31,0.3) and (6.95,1.4) .. (10.93,3.29)   ;

\draw    (450,274.09) .. controls (478.07,293.99) and (516.33,296.23) .. (548.53,275.07) ;
\draw [shift={(550,274.09)}, rotate = 505.66] [color={rgb, 255:red, 0; green, 0; blue, 0 }  ][line width=0.75]    (10.93,-3.29) .. controls (6.95,-1.4) and (3.31,-0.3) .. (0,0) .. controls (3.31,0.3) and (6.95,1.4) .. (10.93,3.29)   ;

\draw (49,47) node    {$9$};
\draw (88,87) node    {$8$};
\draw (14,89) node    {$2$};
\draw (48,130) node    {$5$};
\draw (90,146) node    {$3$};
\draw (12,147) node    {$1$};
\draw (48,192) node    {$4$};
\draw (47,250) node    {$6$};
\draw (11,206) node    {$7$};
\draw (91,208) node    {$2$};
\draw (153,170) node    {$2$};
\draw (173,207) node    {$4$};
\draw (153,250) node    {$6$};
\draw (233,210) node    {$5$};
\draw (307,210) node    {$9$};
\draw (113,210) node    {$7$};
\draw (207,170) node    {$3$};
\draw (207,250) node    {$1$};
\draw (267,170) node    {$8$};
\draw (273,250) node    {$2$};
\draw (373,170) node    {$2$};
\draw (333,210) node    {$7$};
\draw (373,250) node    {$6$};
\draw (393,210) node    {$4$};
\draw (427,250) node    {$1$};
\draw (427,170) node    {$5$};
\draw (467,210) node    {$3$};
\draw (427,110) node    {$9$};
\draw (467,150) node    {$8$};
\draw (393,150) node    {$2$};
\draw (533,250) node    {$6$};
\draw (593,250) node    {$1$};
\draw (627,210) node    {$3$};
\draw (593,110) node    {$8$};
\draw (627,150) node    {$9$};
\draw (553,150) node    {$2$};
\draw (593,170) node    {$5$};
\draw (553,210) node    {$4$};
\draw (493,210) node    {$7$};
\draw (533,170) node    {$2$};
\draw (139,301) node    {$T_{1}$};
\draw (321,301) node    {$T_{2}$};
\draw (499,301) node    {$T_{3}$};

\end{tikzpicture}
  
\end {figure}
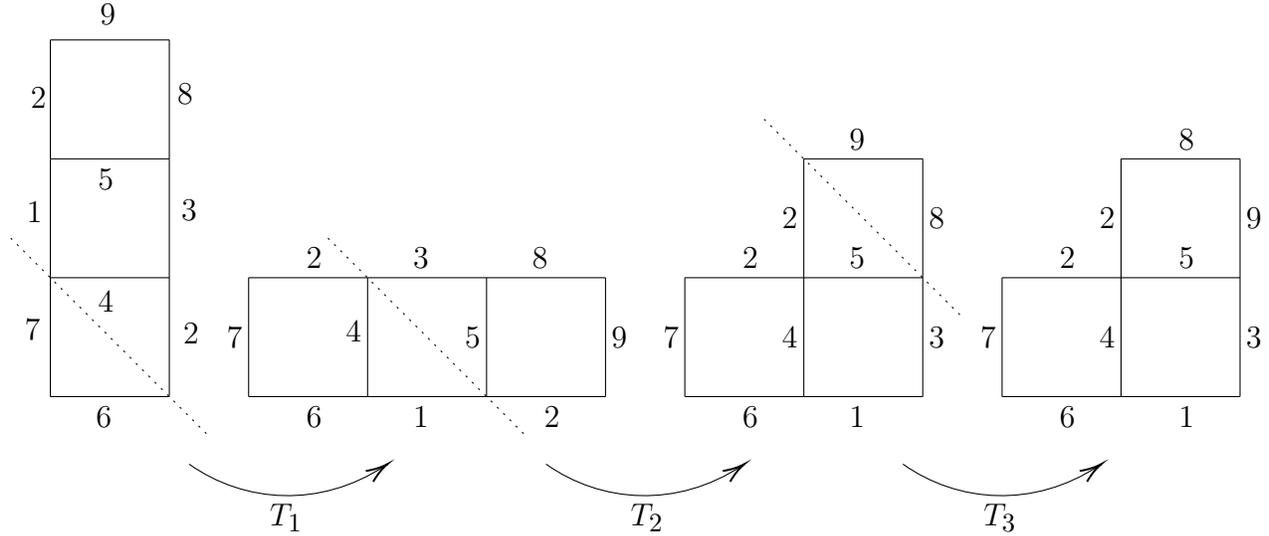

\begin{ex}
    We demonstrate in Figure \ref{fig:factorization} the factorization $G_{ab} \mapsto G_{ab}^{T_1 \circ T_2 \circ \dots \circ T_n} = G_{bb}$.
\end{ex}

The next result follows from the factorization just given.

\begin{prop} \label{shape_duality}
    For any word $w$, we have 

\begin{enumerate}[(a)]
\item $\textit{sh}(G_w) = w^{*}$  and $\textit{sh}(G_{w}^{*}) = w.$
    \item $\textit{sh}(G_w)^{*} = \textit{sh}(G_{w}^{*}).$ 
\item $G_{w}^{*} = G_{w^{*}}.$
\end{enumerate}
\end{prop}

\begin{defn}
      Fix the word $w$, and consider the sign sequence $\mathbf{s}_w = (s(e_0) , s(e_1) , \dots , s(e_n))$ on the snake graph $G_w$. The \textit{dual sign sequence} $s_{w}^{*}$ is defined by 
$$\mathbf{s}_{w}^{*} = (s(e_0) , -s(e_1) , s(e_2) , -s(e_3) , \dots).$$
\end{defn}

\begin{prop} \label{sign_dual}
    For any word $w$, we have $\mathbf{s}_{w}^{*} = \mathbf{s}_{w^{*}}.$
\end{prop}

\begin{proof}
    This follows immediately from Proposition \ref{shape_duality}.
\end{proof}

\begin{ex} 
     Below is the sign sequence $\mathbf{s}_{ab}$ and its dual $\mathbf{s}_{ab}^{*} = \mathbf{s}_{bb}.$
\end{ex}

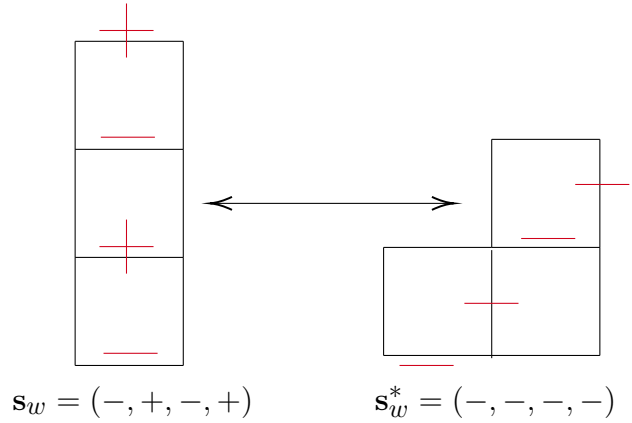
\begin {figure}[h!]
    \centering
    \caption{The sign sequence $\mathbf{s}_{ab}$ and its dual $\mathbf{s}_{bb}$}
    \label{fig:s_ab}
      \begin{tikzpicture}[x=0.75pt,y=0.75pt,yscale=-1,xscale=1]

\draw    (157,93.62) -- (157,148.16) ;
\draw    (157,93.62) -- (211.54,93.62) ;
\draw    (211.54,93.62) -- (211.54,148.16) ;
\draw    (157,148.16) -- (211.54,148.16) ;
\draw    (157,39.09) -- (157,93.62) ;
\draw    (157,39.09) -- (211.54,39.09) ;
\draw    (211.54,39.09) -- (211.54,93.62) ;
\draw    (157,148.16) -- (157,202.7) ;
\draw    (211.54,148.16) -- (211.54,202.7) ;
\draw    (157,202.7) -- (211.54,202.7) ;
\draw [color={rgb, 255:red, 208; green, 2; blue, 27 }  ,draw opacity=1 ]   (171.32,196.56) -- (198.58,196.56) ;
\draw [color={rgb, 255:red, 208; green, 2; blue, 27 }  ,draw opacity=1 ]   (169.27,142.71) -- (196.54,142.71) ;
\draw [color={rgb, 255:red, 208; green, 2; blue, 27 }  ,draw opacity=1 ]   (182.9,129.07) -- (182.9,156.34) ;

\draw [color={rgb, 255:red, 208; green, 2; blue, 27 }  ,draw opacity=1 ]   (169.95,87.49) -- (197.22,87.49) ;
\draw [color={rgb, 255:red, 208; green, 2; blue, 27 }  ,draw opacity=1 ]   (169.27,33.63) -- (196.54,33.63) ;
\draw [color={rgb, 255:red, 208; green, 2; blue, 27 }  ,draw opacity=1 ]   (182.9,20) -- (182.9,47.27) ;

\draw    (367.1,144.48) -- (367.1,199.02) ;
\draw    (367.1,143.12) -- (421.64,143.12) ;
\draw    (421.64,143.12) -- (421.64,197.65) ;
\draw    (367.1,88.58) -- (367.1,143.12) ;
\draw    (367.1,88.58) -- (421.64,88.58) ;
\draw    (421.64,88.58) -- (421.64,143.12) ;
\draw    (367.1,197.65) -- (421.64,197.65) ;
\draw    (312.57,197.65) -- (367.1,197.65) ;
\draw    (312.57,143.12) -- (312.57,197.65) ;
\draw    (312.57,143.12) -- (367.1,143.12) ;
\draw [color={rgb, 255:red, 208; green, 2; blue, 27 }  ,draw opacity=1 ]   (409.23,111.35) -- (436.5,111.35) ;
\draw [color={rgb, 255:red, 208; green, 2; blue, 27 }  ,draw opacity=1 ]   (320.61,202.7) -- (347.88,202.7) ;
\draw [color={rgb, 255:red, 208; green, 2; blue, 27 }  ,draw opacity=1 ]   (353.33,171.34) -- (380.6,171.34) ;
\draw [color={rgb, 255:red, 208; green, 2; blue, 27 }  ,draw opacity=1 ]   (381.96,138.62) -- (409.23,138.62) ;
\draw    (227.17,120.89) -- (345.88,120.89) ;
\draw [shift={(347.88,120.89)}, rotate = 180] [color={rgb, 255:red, 0; green, 0; blue, 0 }  ][line width=0.75]    (10.93,-3.29) .. controls (6.95,-1.4) and (3.31,-0.3) .. (0,0) .. controls (3.31,0.3) and (6.95,1.4) .. (10.93,3.29)   ;
\draw [shift={(225.17,120.89)}, rotate = 0] [color={rgb, 255:red, 0; green, 0; blue, 0 }  ][line width=0.75]    (10.93,-3.29) .. controls (6.95,-1.4) and (3.31,-0.3) .. (0,0) .. controls (3.31,0.3) and (6.95,1.4) .. (10.93,3.29)   ;

\draw (186,221) node    {$\mathbf{s}_{w} =( -,+,-,+)$};
\draw (369,221) node    {$\mathbf{s}^{*}_{w} =( -,-,-,-)$};

\end{tikzpicture}
 
\end {figure}

\section {Continued Fractions}

For more on the involution on continued fractions given in the next definition, see \cite{uludaug2015jimm}.

\begin{defn}
       Consider the finite positive continued fraction $[a_1 , a_2 , \dots , a_d]$. The \textit{dual continued fraction} $[a_1 , a_2 , \dots , a_d]^{*}$ is gotten from $[a_1 , a_2 , \dots , a_d]$ by first writing each $a_i$ as $1+1+...+1,$ substituting these expressions into their respective entries in the continued fraction, and applying the involution that exchanges the symbol ``,'' with the symbol ``+''.
\end{defn}

\begin{ex}
      The continued fraction associated to the word $w=ab$ is $\text{CF}(ab) = [1,1,1,1]$ (see Example \ref{CF_ab}). The dual continued fraction is computed as
$$\text{CF}(w)^{*} = \text{CF}(ab)^{*} = [1,1,1,1]^{*} = [1+1+1+1] = [4] = \text{CF}(bb) = \text{CF}(w^{*})$$

\end{ex}

\begin{rmk}
    The continued fractions in (a) from Remark \ref{fib} are dual to those in (c), and the same for (b) and (d).
\end{rmk}

\begin{prop}
    For any word $w$ we have $$\textit{CF}(w)^{*} = \textit{CF}(w^{*}).$$ 
\end{prop}

\begin{proof}
This follows immediately from Proposition \ref{sign_dual}.
\end{proof}

\section{Distributive Lattices} \label{dist_dual_sec}

\begin{defn} \label{dual_dist_def}
       The \textit{distributive lattice dual to $D_w$} is simply $D_{w^{*}} = \mathcal{I}(C_{w}^{*}).$
\end{defn}

\begin {figure}[h!]
    \centering
    \caption{The Fibonacci cube $\Gamma_3$ and its dual}
    \label{fig:fib_dual}
    \begin{tikzpicture}[x=0.75pt,y=0.75pt,yscale=-1,xscale=1]

\draw    (100,50) -- (130,80) ;
\draw [shift={(130,80)}, rotate = 45] [color={rgb, 255:red, 0; green, 0; blue, 0 }  ][fill={rgb, 255:red, 0; green, 0; blue, 0 }  ][line width=0.75]      (0, 0) circle [x radius= 3.35, y radius= 3.35]   ;
\draw [shift={(100,50)}, rotate = 45] [color={rgb, 255:red, 0; green, 0; blue, 0 }  ][fill={rgb, 255:red, 0; green, 0; blue, 0 }  ][line width=0.75]      (0, 0) circle [x radius= 3.35, y radius= 3.35]   ;
\draw    (70,80) -- (100,110) ;
\draw [shift={(100,110)}, rotate = 45] [color={rgb, 255:red, 0; green, 0; blue, 0 }  ][fill={rgb, 255:red, 0; green, 0; blue, 0 }  ][line width=0.75]      (0, 0) circle [x radius= 3.35, y radius= 3.35]   ;
\draw [shift={(70,80)}, rotate = 45] [color={rgb, 255:red, 0; green, 0; blue, 0 }  ][fill={rgb, 255:red, 0; green, 0; blue, 0 }  ][line width=0.75]      (0, 0) circle [x radius= 3.35, y radius= 3.35]   ;
\draw    (100,50) -- (70,80) ;
\draw [shift={(70,80)}, rotate = 135] [color={rgb, 255:red, 0; green, 0; blue, 0 }  ][fill={rgb, 255:red, 0; green, 0; blue, 0 }  ][line width=0.75]      (0, 0) circle [x radius= 3.35, y radius= 3.35]   ;
\draw [shift={(100,50)}, rotate = 135] [color={rgb, 255:red, 0; green, 0; blue, 0 }  ][fill={rgb, 255:red, 0; green, 0; blue, 0 }  ][line width=0.75]      (0, 0) circle [x radius= 3.35, y radius= 3.35]   ;
\draw    (130,80) -- (100,110) ;
\draw [shift={(100,110)}, rotate = 135] [color={rgb, 255:red, 0; green, 0; blue, 0 }  ][fill={rgb, 255:red, 0; green, 0; blue, 0 }  ][line width=0.75]      (0, 0) circle [x radius= 3.35, y radius= 3.35]   ;
\draw [shift={(130,80)}, rotate = 135] [color={rgb, 255:red, 0; green, 0; blue, 0 }  ][fill={rgb, 255:red, 0; green, 0; blue, 0 }  ][line width=0.75]      (0, 0) circle [x radius= 3.35, y radius= 3.35]   ;
\draw    (100,110) -- (100,140) ;
\draw [shift={(100,140)}, rotate = 90] [color={rgb, 255:red, 0; green, 0; blue, 0 }  ][fill={rgb, 255:red, 0; green, 0; blue, 0 }  ][line width=0.75]      (0, 0) circle [x radius= 3.35, y radius= 3.35]   ;
\draw [shift={(100,110)}, rotate = 90] [color={rgb, 255:red, 0; green, 0; blue, 0 }  ][fill={rgb, 255:red, 0; green, 0; blue, 0 }  ][line width=0.75]      (0, 0) circle [x radius= 3.35, y radius= 3.35]   ;
\draw    (210,110) -- (210,140) ;
\draw [shift={(210,140)}, rotate = 90] [color={rgb, 255:red, 0; green, 0; blue, 0 }  ][fill={rgb, 255:red, 0; green, 0; blue, 0 }  ][line width=0.75]      (0, 0) circle [x radius= 3.35, y radius= 3.35]   ;
\draw [shift={(210,110)}, rotate = 90] [color={rgb, 255:red, 0; green, 0; blue, 0 }  ][fill={rgb, 255:red, 0; green, 0; blue, 0 }  ][line width=0.75]      (0, 0) circle [x radius= 3.35, y radius= 3.35]   ;
\draw    (210,80) -- (210,110) ;
\draw [shift={(210,110)}, rotate = 90] [color={rgb, 255:red, 0; green, 0; blue, 0 }  ][fill={rgb, 255:red, 0; green, 0; blue, 0 }  ][line width=0.75]      (0, 0) circle [x radius= 3.35, y radius= 3.35]   ;
\draw [shift={(210,80)}, rotate = 90] [color={rgb, 255:red, 0; green, 0; blue, 0 }  ][fill={rgb, 255:red, 0; green, 0; blue, 0 }  ][line width=0.75]      (0, 0) circle [x radius= 3.35, y radius= 3.35]   ;
\draw    (210,50) -- (210,80) ;
\draw [shift={(210,80)}, rotate = 90] [color={rgb, 255:red, 0; green, 0; blue, 0 }  ][fill={rgb, 255:red, 0; green, 0; blue, 0 }  ][line width=0.75]      (0, 0) circle [x radius= 3.35, y radius= 3.35]   ;
\draw [shift={(210,50)}, rotate = 90] [color={rgb, 255:red, 0; green, 0; blue, 0 }  ][fill={rgb, 255:red, 0; green, 0; blue, 0 }  ][line width=0.75]      (0, 0) circle [x radius= 3.35, y radius= 3.35]   ;
\draw    (133,110) -- (187,110) ;
\draw [shift={(190,110)}, rotate = 180] [fill={rgb, 255:red, 0; green, 0; blue, 0 }  ][line width=0.08]  [draw opacity=0] (8.93,-4.29) -- (0,0) -- (8.93,4.29) -- cycle    ;
\draw [shift={(130,110)}, rotate = 0] [fill={rgb, 255:red, 0; green, 0; blue, 0 }  ][line width=0.08]  [draw opacity=0] (8.93,-4.29) -- (0,0) -- (8.93,4.29) -- cycle    ;

\draw (105,159) node    {$D_{ab}$};
\draw (215,159) node    {$D_{bb}$};

\end{tikzpicture}
  
\end {figure}
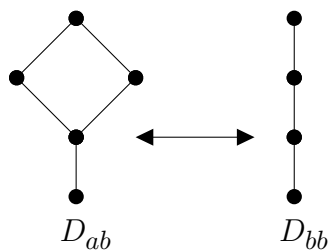

\begin{ex}
    The three expansion posets pictured in Examples \ref{fig:P_ab}, \ref{fig:A_ab}, and \ref{fig:T_ab} are all isomorphic to the same distributive lattice $D_{ab} \cong \Gamma_3$. The dual lattice $D_{bb}$ is a chain poset on $4$ vertices. See Figure \ref{fig:fib_dual}.
\end{ex}

\begin{rmk}
    In general, a chain poset with $n+1$ vertices is dual to the Fibonacci cube $\Gamma_{n}.$
\end{rmk}

\chapter{Dual Expansion Posets}

\section{Lattice Paths on Snake Graphs}

We recall here the \textit{lattice path expansion posets} from \cite{propp2005combinatorics}.

\begin{defn}
       A \textit{lattice path} in a snake graph $G_w$ with $n$ tiles is a choice of $n+1$ edges $L$ from $G_w$ which when concatenated form a path, taking only unit steps right or up, from the SW vertex of tile $T_1$ to the NE vertex of tile $T_n$.  The \textit{weight} $x_L$ of $L$ is defined to be the product of initial cluster variables $x_L = \prod_{l \in L} x_l.$ Let $\mathbb{L}_w$ be the set of all lattice paths of the snake graph $G_w.$
\end{defn}

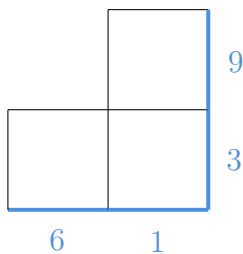
\begin {figure}[h!]
    \centering
    \caption{The lattice path $L_{-}$ on $G_{bb}$}
    \label{fig:L_}
    \begin{tikzpicture}[x=0.75pt,y=0.75pt,yscale=-1,xscale=1]

\draw [color={rgb, 255:red, 74; green, 144; blue, 226 }  ,draw opacity=1 ][line width=1.5]    (183.5,267.16) -- (234.02,267.16) ;
\draw    (183.5,216.65) -- (234.02,216.65) ;
\draw    (234.02,267.16) -- (234.02,216.65) ;
\draw    (183.5,267.16) -- (183.5,216.65) ;
\draw [color={rgb, 255:red, 74; green, 144; blue, 226 }  ,draw opacity=1 ][line width=1.5]    (284.54,267.16) -- (284.54,216.65) ;
\draw [color={rgb, 255:red, 74; green, 144; blue, 226 }  ,draw opacity=1 ][line width=1.5]    (234.02,267.16) -- (284.54,267.16) ;
\draw    (234.02,216.65) -- (284.54,216.65) ;
\draw    (234.02,166.13) -- (284.54,166.13) ;
\draw [color={rgb, 255:red, 74; green, 144; blue, 226 }  ,draw opacity=1 ][line width=1.5]    (284.54,216.65) -- (284.54,166.13) ;
\draw    (234.02,216.65) -- (234.02,166.13) ;

\draw (208.42,282.32) node    {$\textcolor[rgb]{0.29,0.56,0.89}{6}$};
\draw (258.94,283.16) node    {$\textcolor[rgb]{0.29,0.56,0.89}{1}$};
\draw (297.67,241.9) node    {$\textcolor[rgb]{0.29,0.56,0.89}{3}$};
\draw (298.51,192.23) node    {$\textcolor[rgb]{0.29,0.56,0.89}{9}$};

\end{tikzpicture}

\end {figure}

\begin{ex} 
    Figure \ref{fig:L_} shows one of the five lattice paths on the snake graph $G_{bb}$ dual to $G_{ab}.$ Note that the weight of this lattice path coincides with the weight of the the perfect matching shown in Figure \ref{fig:P_}.
\end{ex}

The next expansion formula is from \cite{propp2005combinatorics}.

\begin{thm} 
    Let $w$ be any word, and consider the set $\mathbb{L}_{w^{*}}$ of lattice paths on the dual snake graph $G_{w^{*}}$. Then the cluster variable $x_w$ can be written as
$$x_{w} = \frac{1}{x_1 x_2 \dots x_n}  \sum_{L \in \mathbb{L}_{w^{*}}} x_L.$$
\end{thm}

\begin{cor} \label{card_equal}
    For any word $w$, we have $\abs{\mathbb{P}_w} = \abs{\mathbb{L}_{w^{*}}}$ and $\abs{\mathbb{L}_w} = \abs{\mathbb{P}_{w^{*}}}.$
\end{cor}

Let $\text{CF}(w) = [a_1 , a_2 , \dots , a_k]$ and $\text{CF}(w^{*}) = [b_1 , b_2 , \dots , b_l].$ Let $G_{w}^{a_1}$ be the subsnake graph of $G_{w}$ obtained by deleting the first $a_1$ tiles of $G_w$ (see the discussion directly before Theorem \ref{CF_card_quotient_P}). Let $G_{w}^{b_1}$ be the subsnake graph of $G_{w}$ obtained by deleting the first $b_1$ tiles of $G_w$. Combining Corollary \ref{card_equal} with Theorem \ref{CF_card_quotient_P} gives the next result.

\begin{cor}
    For any word $w$, we have
$$\text{CF}(w) = \frac{\abs{\mathbb{P}_w}}{\abs{\mathbb{P}_{w}^{a_1}}} = \frac{\abs{\mathbb{L}_{w^{*}}}}{\abs{\mathbb{L}_{w^{*}}^{a_1}}}$$
and 
$$\text{CF}(w)^{*} = \frac{\abs{\mathbb{P}_{w^{*}}}}{\abs{\mathbb{P}_{w^{*}}^{b_1}}} = \frac{\abs{\mathbb{L}_w}}{\abs{\mathbb{L}_{w}^{b_1}}}.$$

\end{cor}

We now give the set $\mathbb{L}_w$ a poset structure.

A \textit{flip} of $L$ at tile $T_i$ is the local move that takes two edges of $L$ located on the same tile $T_i$ of $G_w$ that are incident to a common vertex of $T_i$ (so, necessarily the two edges in question are the S and E edges of $T_i$, or the W and N edges of tile $T_i$) and replaces them with the other two edges of $T_i$.

Directly below is the local picture for the flip at a generic tile $T_i$. 

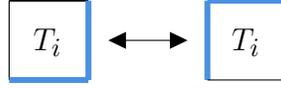
\begin {figure}[h!]
    \centering
    \caption{Flip of a lattice path at tile $T_i$}
    \label{fig:flip}
    \begin{tikzpicture}[x=0.75pt,y=0.75pt,yscale=-1,xscale=1]

\draw    (220,120) -- (220,160) ;
\draw [color={rgb, 255:red, 0; green, 0; blue, 0 }  ,draw opacity=1 ][line width=0.75]    (220,120) -- (260,120) ;
\draw [color={rgb, 255:red, 74; green, 144; blue, 226 }  ,draw opacity=1 ][line width=2.25]    (220,160) -- (260,160) ;
\draw [color={rgb, 255:red, 74; green, 144; blue, 226 }  ,draw opacity=1 ][line width=2.25]    (260,120) -- (260,160) ;
\draw [color={rgb, 255:red, 74; green, 144; blue, 226 }  ,draw opacity=1 ][line width=2.25]    (320,120) -- (320,160) ;
\draw [color={rgb, 255:red, 74; green, 144; blue, 226 }  ,draw opacity=1 ][line width=2.25]    (320,120) -- (360,120) ;
\draw    (320,160) -- (360,160) ;
\draw [color={rgb, 255:red, 0; green, 0; blue, 0 }  ,draw opacity=1 ][line width=0.75]    (360,120) -- (360,160) ;
\draw    (273,140) -- (307,140) ;
\draw [shift={(310,140)}, rotate = 180] [fill={rgb, 255:red, 0; green, 0; blue, 0 }  ][line width=0.08]  [draw opacity=0] (8.93,-4.29) -- (0,0) -- (8.93,4.29) -- cycle    ;
\draw [shift={(270,140)}, rotate = 0] [fill={rgb, 255:red, 0; green, 0; blue, 0 }  ][line width=0.08]  [draw opacity=0] (8.93,-4.29) -- (0,0) -- (8.93,4.29) -- cycle    ;

\draw (239.5,141) node    {$T_{i}$};
\draw (339.5,141) node    {$T_{i}$};

\end{tikzpicture}
  
\end {figure}

\newpage 

An \textit{up-flip} at tile $T_i$ is a flip that replaces the S and E edges of $T_i$ with the N and W edges of $T_i$.

\begin{defn} \label{L_min_def}
       The \textit{minimal element} $L_{-}$ of $\mathbb{L}_w$ is the unique lattice path of $G_w$ such that every edge in $L_{-}$ is a boundary edge of $G_w$ and the S edge of tile $T_1$ is in $L_{-}.$ The \textit{maximal element} $L_{+}$ of $\mathbb{L}_w$ is the unique lattice path of $G_w$ such that every edge in $L_{+}$ is a boundary edge of $G_w$ and the S edge of tile $T_1$ is not in $L_{-}.$
\end{defn}

\begin{defn}
       The \textit{poset structure on $\mathbb{L}_w$} is defined as follows. The unique minimal element of $\mathbb{L}_w$ is the minimal lattice path $L_{-},$ and the unique maximal element is $L_{+}.$ A lattice path $L_2$ covers a lattice path $L_1$ if there exists a tile $T_i$ such that $L_2$ can be obtained from $L_1$ by performing a single up-flip of $L_1$ at $T_i.$
\end{defn}

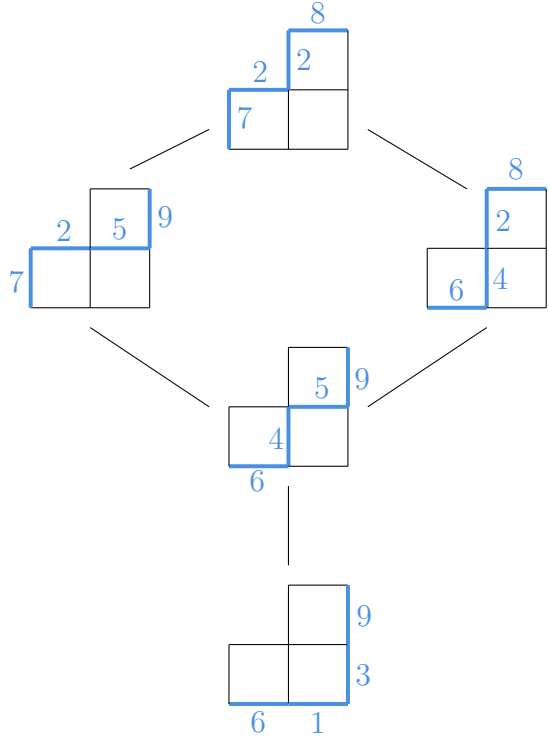
\begin {figure}[h!]
    \centering
    \caption{The poset $\mathbb{L}_{bb}$}
    \label{fig:L_bb}
    \begin{tikzpicture}[x=0.75pt,y=0.75pt,yscale=-1,xscale=1]

\draw [color={rgb, 255:red, 74; green, 144; blue, 226 }  ,draw opacity=1 ][line width=1.5]    (215,438) -- (245,438) ;
\draw    (215,408) -- (245,408) ;
\draw    (245,438) -- (245,408) ;
\draw    (215,438) -- (215,408) ;
\draw [color={rgb, 255:red, 74; green, 144; blue, 226 }  ,draw opacity=1 ][line width=1.5]    (275,438) -- (275,408) ;
\draw [color={rgb, 255:red, 74; green, 144; blue, 226 }  ,draw opacity=1 ][line width=1.5]    (245,438) -- (275,438) ;
\draw    (245,408) -- (275,408) ;
\draw    (245,378) -- (275,378) ;
\draw [color={rgb, 255:red, 74; green, 144; blue, 226 }  ,draw opacity=1 ][line width=1.5]    (275,408) -- (275,378) ;
\draw    (245,408) -- (245,378) ;
\draw [color={rgb, 255:red, 74; green, 144; blue, 226 }  ,draw opacity=1 ][line width=1.5]    (215,318) -- (245,318) ;
\draw    (215,288) -- (245,288) ;
\draw [color={rgb, 255:red, 74; green, 144; blue, 226 }  ,draw opacity=1 ][line width=1.5]    (245,318) -- (245,288) ;
\draw    (215,318) -- (215,288) ;
\draw    (275,318) -- (275,288) ;
\draw    (245,318) -- (275,318) ;
\draw [color={rgb, 255:red, 74; green, 144; blue, 226 }  ,draw opacity=1 ][line width=1.5]    (245,288) -- (275,288) ;
\draw    (245,258) -- (275,258) ;
\draw [color={rgb, 255:red, 74; green, 144; blue, 226 }  ,draw opacity=1 ][line width=1.5]    (275,288) -- (275,258) ;
\draw    (245,288) -- (245,258) ;
\draw [color={rgb, 255:red, 74; green, 144; blue, 226 }  ,draw opacity=1 ][line width=1.5]    (315,238) -- (345,238) ;
\draw    (315,208) -- (345,208) ;
\draw [color={rgb, 255:red, 74; green, 144; blue, 226 }  ,draw opacity=1 ][line width=1.5]    (345,238) -- (345,208) ;
\draw    (315,238) -- (315,208) ;
\draw    (375,238) -- (375,208) ;
\draw    (345,238) -- (375,238) ;
\draw    (345,208) -- (375,208) ;
\draw [color={rgb, 255:red, 74; green, 144; blue, 226 }  ,draw opacity=1 ][line width=1.5]    (345,178) -- (375,178) ;
\draw    (375,208) -- (375,178) ;
\draw [color={rgb, 255:red, 74; green, 144; blue, 226 }  ,draw opacity=1 ][line width=1.5]    (345,208) -- (345,178) ;
\draw    (115,238) -- (145,238) ;
\draw [color={rgb, 255:red, 74; green, 144; blue, 226 }  ,draw opacity=1 ][line width=1.5]    (115,208) -- (145,208) ;
\draw    (145,238) -- (145,208) ;
\draw [color={rgb, 255:red, 74; green, 144; blue, 226 }  ,draw opacity=1 ][line width=1.5]    (115,238) -- (115,208) ;
\draw    (175,238) -- (175,208) ;
\draw    (145,238) -- (175,238) ;
\draw [color={rgb, 255:red, 74; green, 144; blue, 226 }  ,draw opacity=1 ][line width=1.5]    (145,208) -- (175,208) ;
\draw    (145,178) -- (175,178) ;
\draw [color={rgb, 255:red, 74; green, 144; blue, 226 }  ,draw opacity=1 ][line width=1.5]    (175,208) -- (175,178) ;
\draw    (145,208) -- (145,178) ;
\draw    (215,158) -- (245,158) ;
\draw [color={rgb, 255:red, 74; green, 144; blue, 226 }  ,draw opacity=1 ][line width=1.5]    (215,128) -- (245,128) ;
\draw    (245,158) -- (245,128) ;
\draw [color={rgb, 255:red, 74; green, 144; blue, 226 }  ,draw opacity=1 ][line width=1.5]    (215,158) -- (215,128) ;
\draw    (275,158) -- (275,128) ;
\draw    (245,158) -- (275,158) ;
\draw    (245,128) -- (275,128) ;
\draw [color={rgb, 255:red, 74; green, 144; blue, 226 }  ,draw opacity=1 ][line width=1.5]    (245,98) -- (275,98) ;
\draw    (275,128) -- (275,98) ;
\draw [color={rgb, 255:red, 74; green, 144; blue, 226 }  ,draw opacity=1 ][line width=1.5]    (245,128) -- (245,98) ;
\draw    (145,248) -- (205,288) ;
\draw    (285,148) -- (335,178) ;
\draw    (165,168) -- (205,148) ;
\draw    (285,288) -- (345,248) ;
\draw    (245,328) -- (245,368) ;

\draw (229.8,447) node    {$\textcolor[rgb]{0.29,0.56,0.89}{6}$};
\draw (259.8,447.5) node    {$\textcolor[rgb]{0.29,0.56,0.89}{1}$};
\draw (282.8,423) node    {$\textcolor[rgb]{0.29,0.56,0.89}{3}$};
\draw (283.3,393.5) node    {$\textcolor[rgb]{0.29,0.56,0.89}{9}$};
\draw (282.3,274) node    {$\textcolor[rgb]{0.29,0.56,0.89}{9}$};
\draw (182.8,192) node    {$\textcolor[rgb]{0.29,0.56,0.89}{9}$};
\draw (229.3,325.5) node    {$\textcolor[rgb]{0.29,0.56,0.89}{6}$};
\draw (329.8,229) node    {$\textcolor[rgb]{0.29,0.56,0.89}{6}$};
\draw (238.8,304) node    {$\textcolor[rgb]{0.29,0.56,0.89}{4}$};
\draw (261.8,278.5) node    {$\textcolor[rgb]{0.29,0.56,0.89}{5}$};
\draw (159.8,198) node    {$\textcolor[rgb]{0.29,0.56,0.89}{5}$};
\draw (351.8,223.5) node    {$\textcolor[rgb]{0.29,0.56,0.89}{4}$};
\draw (131.8,199) node    {$\textcolor[rgb]{0.29,0.56,0.89}{2}$};
\draw (353.3,194) node    {$\textcolor[rgb]{0.29,0.56,0.89}{2}$};
\draw (107.8,225) node    {$\textcolor[rgb]{0.29,0.56,0.89}{7}$};
\draw (222.8,142) node    {$\textcolor[rgb]{0.29,0.56,0.89}{7}$};
\draw (230.8,119) node    {$\textcolor[rgb]{0.29,0.56,0.89}{2}$};
\draw (252.8,113) node    {$\textcolor[rgb]{0.29,0.56,0.89}{2}$};
\draw (259.8,89) node    {$\textcolor[rgb]{0.29,0.56,0.89}{8}$};
\draw (359.8,168) node    {$\textcolor[rgb]{0.29,0.56,0.89}{8}$};

\end{tikzpicture}
  
\end {figure}

\begin{ex}
    Fix $w=ab$. In Figure \ref{fig:L_bb}, we illustrate the poset $\mathbb{L}_{w^{*}} = \mathbb{L}_{bb}$ of lattice paths on the dual snake graph $G_{bb}$. Compare with Figures \ref{fig:P_ab}, \ref{fig:A_ab}, and \ref{fig:T_ab} from Chapter \ref{exp_poset_chapter}. 
\end{ex}

\newpage
\begin{rmk}
    If the word $w$ is straight, the poset $\mathbb{P}_w$ is isomorphic to a Fibonacci cube (see (b) in \ref{P_w_chain_or_fib}), and likewise $\mathbb{L}_{w^{*}}$ is isomorphic to the same Fibonacci cube. Conversely, if the word $w$ is zigzag, then both the posets $\mathbb{P}_w$ and $\mathbb{L}_{w^{*}}$ are isomorphic to the same linear chain. Similar remarks hold for the other expansion formulas. See Theorem \ref{duality_thm} below for details.
\end{rmk}

\section{Lattice Paths of Angles}

Say that a vertex $v$ of $\Sigma$ is \textit{incident} to any triangle that it is a vertex of.

\begin{defn}
       A \textit{lattice path of angles $\beta$} of $\Sigma_w$ is a selection of $n+1$ angles from the triangles $\Delta_0 , \Delta_1 , \dots , \Delta_n$ of $\Delta,$ one per triangle, such that 

\begin{enumerate}[(1)]
    \item Each angle is incident to an endpoint of one of the internal diagonals $\delta_{1} , \delta_{2}, \dots , \delta_{n}$.
    
    \item An even number of angles from $\beta$ are incident to each vertex of the polygon $\Sigma_w$ which is incident to an even number of triangles from $\Sigma_w$, and an odd number of angles from $\beta$ are incident to each vertex of the polygon $\Sigma_w$ which is incident to an odd number of triangles of $\Sigma_w$. 
    
\end{enumerate}

\end{defn}

Each angle $b$ in $\Delta_i$ can be assigned the cluster variable $x_b$ associated to the edge of $\Delta_i$ opposite of $b$. The \textit{weight of $\beta$} is defined to be the product of initial cluster variables $x_{\beta} = \prod_{b \in \beta} x_b.$ Let $\mathbb{B}_w$ be the set of all lattice paths of angles in $\Sigma_w.$

\begin{ex}
    Below we show one of the five lattice paths of angles on the dual triangulated surface $\Sigma_{ab}^{*} = \Sigma_{bb}$. 
\end{ex}

\begin {figure}[h!]
    \centering
    \caption{The lattice path of angles $B_{-}$ on $G_{bb}$}
    \label{fig:B_}
    
    \begin{tikzpicture}[x=0.75pt,y=0.75pt,yscale=-1,xscale=1]

\draw    (282.25,84) -- (324.5,126.25) ;
\draw    (240,126.25) -- (282.25,84) ;
\draw    (324.5,210.75) -- (324.5,126.25) ;
\draw    (240,210.75) -- (240,126.25) ;
\draw    (240,210.75) -- (282.25,253) ;
\draw    (282.25,253) -- (324.5,210.75) ;
\draw    (240,210.75) -- (282.25,84) ;
\draw    (324.5,210.75) -- (282.25,84) ;
\draw    (282.25,253) -- (282.25,84) ;
\draw  [color={rgb, 255:red, 74; green, 144; blue, 226 }  ,draw opacity=1 ][fill={rgb, 255:red, 74; green, 144; blue, 226 }  ,fill opacity=1 ] (256.61,117.04) .. controls (256.61,113.96) and (259.1,111.46) .. (262.18,111.46) .. controls (265.26,111.46) and (267.75,113.96) .. (267.75,117.04) .. controls (267.75,120.11) and (265.26,122.61) .. (262.18,122.61) .. controls (259.1,122.61) and (256.61,120.11) .. (256.61,117.04) -- cycle ;
\draw  [color={rgb, 255:red, 74; green, 144; blue, 226 }  ,draw opacity=1 ][fill={rgb, 255:red, 74; green, 144; blue, 226 }  ,fill opacity=1 ] (267.17,230.05) .. controls (267.17,226.98) and (269.67,224.48) .. (272.74,224.48) .. controls (275.82,224.48) and (278.32,226.98) .. (278.32,230.05) .. controls (278.32,233.13) and (275.82,235.63) .. (272.74,235.63) .. controls (269.67,235.63) and (267.17,233.13) .. (267.17,230.05) -- cycle ;
\draw  [color={rgb, 255:red, 74; green, 144; blue, 226 }  ,draw opacity=1 ][fill={rgb, 255:red, 74; green, 144; blue, 226 }  ,fill opacity=1 ] (285.89,230.82) .. controls (285.89,227.74) and (288.39,225.25) .. (291.46,225.25) .. controls (294.54,225.25) and (297.04,227.74) .. (297.04,230.82) .. controls (297.04,233.9) and (294.54,236.39) .. (291.46,236.39) .. controls (288.39,236.39) and (285.89,233.9) .. (285.89,230.82) -- cycle ;
\draw  [color={rgb, 255:red, 74; green, 144; blue, 226 }  ,draw opacity=1 ][fill={rgb, 255:red, 74; green, 144; blue, 226 }  ,fill opacity=1 ] (296.88,115.9) .. controls (296.88,112.82) and (299.37,110.33) .. (302.45,110.33) .. controls (305.53,110.33) and (308.02,112.82) .. (308.02,115.9) .. controls (308.02,118.98) and (305.53,121.47) .. (302.45,121.47) .. controls (299.37,121.47) and (296.88,118.98) .. (296.88,115.9) -- cycle ;

\draw (234,161) node    {$\textcolor[rgb]{0.29,0.56,0.89}{6}$};
\draw (265,161) node    {$\textcolor[rgb]{0.29,0.56,0.89}{1}$};
\draw (302,161) node    {$\textcolor[rgb]{0.29,0.56,0.89}{3}$};
\draw (332,160) node    {$\textcolor[rgb]{0.29,0.56,0.89}{9}$};
\end{tikzpicture}
  
\end {figure}
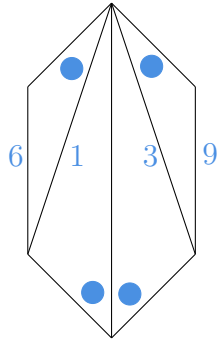

The next result is included as part of the statement of Theorem \ref{duality_thm} (e), and is proved there.

\begin{thm}
    Let $w$ be any word. Consider the set $\mathbb{B}_{w^{*}}$ of lattice paths of angles on the dual triangulated surface $\Sigma_{w}^{*}$. Then the cluster variable $x_w$ can be written as 
    $$x_{w} = \frac{1}{x_1 x_2 \dots x_n} \sum_{\beta \in \mathbb{B}^{*} }x_{\beta} .$$
\end{thm}

We now give the set $\mathbb{B}_w$ a poset structure.

A \textit{flip} of $B$ at diagonal $\delta_{i}$ is the local move that takes two angles in $B$ that are each incident to the same endpoint of the internal diagonal $\delta_{i}$ and replaces them with the remaining two angles incident to $\delta_{i}$.

Directly below is the flip at a generic internal diagonal $\delta_{i}$. 

\begin {figure}[h!]
    \centering
    \caption{Flip of a lattice path of angles at diagonal $\delta_{i}$}
    \label{fig:flip2}
    \begin{tikzpicture}[x=0.75pt,y=0.75pt,yscale=-1,xscale=1]

\draw  [color={rgb, 255:red, 74; green, 144; blue, 226 }  ,draw opacity=1 ][fill={rgb, 255:red, 74; green, 144; blue, 226 }  ,fill opacity=1 ] (162.98,86.49) .. controls (162.98,82.35) and (166.35,78.98) .. (170.49,78.98) .. controls (174.64,78.98) and (178,82.35) .. (178,86.49) .. controls (178,90.64) and (174.64,94) .. (170.49,94) .. controls (166.35,94) and (162.98,90.64) .. (162.98,86.49) -- cycle ;
\draw  [color={rgb, 255:red, 74; green, 144; blue, 226 }  ,draw opacity=1 ][fill={rgb, 255:red, 74; green, 144; blue, 226 }  ,fill opacity=1 ] (139,86.51) .. controls (139,82.36) and (142.36,79) .. (146.51,79) .. controls (150.65,79) and (154.02,82.36) .. (154.02,86.51) .. controls (154.02,90.65) and (150.65,94.02) .. (146.51,94.02) .. controls (142.36,94.02) and (139,90.65) .. (139,86.51) -- cycle ;
\draw    (158,60) -- (233,135) ;
\draw    (83,135) -- (158,60) ;
\draw    (83,135) -- (158,210) ;
\draw    (158,210) -- (233,135) ;
\draw    (158,60) -- (158,210) ;
\draw [line width=0.75]    (250,135) -- (336,135) ;
\draw [shift={(338,135)}, rotate = 180] [color={rgb, 255:red, 0; green, 0; blue, 0 }  ][line width=0.75]    (10.93,-3.29) .. controls (6.95,-1.4) and (3.31,-0.3) .. (0,0) .. controls (3.31,0.3) and (6.95,1.4) .. (10.93,3.29)   ;
\draw [shift={(248,135)}, rotate = 0] [color={rgb, 255:red, 0; green, 0; blue, 0 }  ][line width=0.75]    (10.93,-3.29) .. controls (6.95,-1.4) and (3.31,-0.3) .. (0,0) .. controls (3.31,0.3) and (6.95,1.4) .. (10.93,3.29)   ;
\draw  [color={rgb, 255:red, 74; green, 144; blue, 226 }  ,draw opacity=1 ][fill={rgb, 255:red, 74; green, 144; blue, 226 }  ,fill opacity=1 ] (438,183.51) .. controls (438,179.36) and (441.36,176) .. (445.51,176) .. controls (449.65,176) and (453.02,179.36) .. (453.02,183.51) .. controls (453.02,187.65) and (449.65,191.02) .. (445.51,191.02) .. controls (441.36,191.02) and (438,187.65) .. (438,183.51) -- cycle ;
\draw  [color={rgb, 255:red, 74; green, 144; blue, 226 }  ,draw opacity=1 ][fill={rgb, 255:red, 74; green, 144; blue, 226 }  ,fill opacity=1 ] (413.98,183.51) .. controls (413.98,179.36) and (417.35,176) .. (421.49,176) .. controls (425.64,176) and (429,179.36) .. (429,183.51) .. controls (429,187.65) and (425.64,191.02) .. (421.49,191.02) .. controls (417.35,191.02) and (413.98,187.65) .. (413.98,183.51) -- cycle ;
\draw    (433,60) -- (508,135) ;
\draw    (358,135) -- (433,60) ;
\draw    (358,135) -- (433,210) ;
\draw    (433,210) -- (508,135) ;
\draw    (433,60) -- (433,210) ;

\draw (159,120) node [anchor=north west][inner sep=0.75pt]    {$i$};
\draw (434,120) node [anchor=north west][inner sep=0.75pt]    {$i$};

\end{tikzpicture}
  
\end {figure}
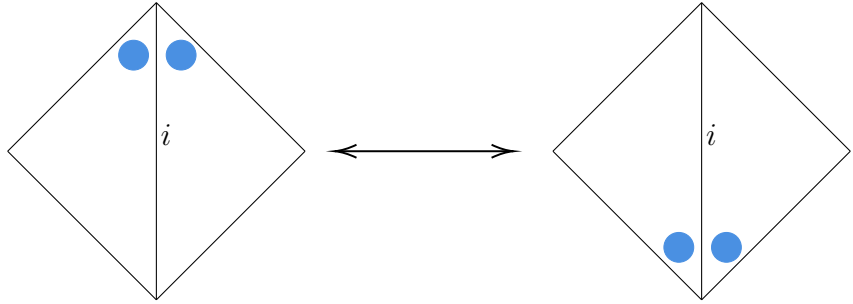

Let $l_i$ be the endpoint of $\delta_{i}$ to the left of $\gamma_{a \rightarrow b}$, and $r_i$ the endpoint to the right. An \textit{up-flip} of $\beta$ at diagonal $\delta_{i}$ is a flip that meets either of the following two conditions: 

\begin{enumerate}[(1)]
        \item $i$ is odd, and the two angles matched to $l_i$ are replaced with the two angles matched to $r_i$, or
        
        \item $i$ is even, and the two angles matched to $r_i$ are replaced with the two angles matched to $l_i$.
    \end{enumerate}

\begin{defn}
       The \textit{minimal element} $B_{-}$ of $\mathbb{B}_w$ is the unique lattice path of angles such that the boundary angle in $\Delta_0$ with boundary edge $\delta_{2n+1}$ is matched, and only boundary angles are used in $B_{-}$. The \textit{maximal element} $B_{+}$ of $\mathbb{B}_w$ is the unique lattice path of angles such that the angle the boundary angle in $\Delta_0$ with boundary edge $\delta_{2n}$ is matched, and only boundary angles are used in $B_{+}$.
\end{defn}

\begin{defn}
       The \textit{poset structure on $\mathbb{B}_w$} is defined as follows. The unique minimal element of $\mathbb{B}_w$ is the minimal lattice path of angles $B_{-},$ and the unique maximal element is $B_{+}.$ A lattice path of angles $B_2$ covers another lattice path of angles $B_1$ if there exists a diagonal $\delta_{i}$ such that $B_2$ can be obtained from $B_1$ by performing a single up-twist of $B_1$ at $\delta_{i}$.
\end{defn}

\begin{ex}
    Fix $w=ab.$ The poset of lattice paths of angles on the dual triangulated surface $\Sigma_{bb}$ is shown in Figure \ref{fig:B_bb} below.  
\end{ex}

\begin {figure}[h!]
    \centering
    \caption{The poset $\mathbb{B}_{bb}$}
    \label{fig:B_bb}
    
  \begin{tikzpicture}[x=0.75pt,y=0.75pt,yscale=-1,xscale=1]

\draw    (320,320) -- (340,340) ;
\draw    (300,340) -- (320,320) ;
\draw    (340,380) -- (340,340) ;
\draw    (300,380) -- (300,340) ;
\draw    (300,380) -- (320,400) ;
\draw    (320,400) -- (340,380) ;
\draw    (300,380) -- (320,320) ;
\draw    (340,380) -- (320,320) ;
\draw    (320,400) -- (320,320) ;
\draw    (320,310) -- (320,284.88) ;
\draw    (350,230) -- (400,200) ;
\draw    (240,200) -- (290,230) ;
\draw    (350,110) -- (400,140) ;
\draw    (240,140) -- (290,110) ;
\draw  [color={rgb, 255:red, 74; green, 144; blue, 226 }  ,draw opacity=1 ][fill={rgb, 255:red, 74; green, 144; blue, 226 }  ,fill opacity=1 ] (307.86,335.64) .. controls (307.86,334.18) and (309.04,333) .. (310.5,333) .. controls (311.96,333) and (313.14,334.18) .. (313.14,335.64) .. controls (313.14,337.1) and (311.96,338.28) .. (310.5,338.28) .. controls (309.04,338.28) and (307.86,337.1) .. (307.86,335.64) -- cycle ;
\draw  [color={rgb, 255:red, 74; green, 144; blue, 226 }  ,draw opacity=1 ][fill={rgb, 255:red, 74; green, 144; blue, 226 }  ,fill opacity=1 ] (312.86,389.14) .. controls (312.86,387.68) and (314.04,386.5) .. (315.5,386.5) .. controls (316.96,386.5) and (318.14,387.68) .. (318.14,389.14) .. controls (318.14,390.6) and (316.96,391.78) .. (315.5,391.78) .. controls (314.04,391.78) and (312.86,390.6) .. (312.86,389.14) -- cycle ;
\draw  [color={rgb, 255:red, 74; green, 144; blue, 226 }  ,draw opacity=1 ][fill={rgb, 255:red, 74; green, 144; blue, 226 }  ,fill opacity=1 ] (321.72,389.5) .. controls (321.72,388.04) and (322.9,386.86) .. (324.36,386.86) .. controls (325.82,386.86) and (327,388.04) .. (327,389.5) .. controls (327,390.96) and (325.82,392.14) .. (324.36,392.14) .. controls (322.9,392.14) and (321.72,390.96) .. (321.72,389.5) -- cycle ;
\draw  [color={rgb, 255:red, 74; green, 144; blue, 226 }  ,draw opacity=1 ][fill={rgb, 255:red, 74; green, 144; blue, 226 }  ,fill opacity=1 ] (326.92,335.1) .. controls (326.92,333.64) and (328.1,332.46) .. (329.56,332.46) .. controls (331.02,332.46) and (332.2,333.64) .. (332.2,335.1) .. controls (332.2,336.56) and (331.02,337.74) .. (329.56,337.74) .. controls (328.1,337.74) and (326.92,336.56) .. (326.92,335.1) -- cycle ;
\draw    (320,200) -- (340,220) ;
\draw    (300,220) -- (320,200) ;
\draw    (340,260) -- (340,220) ;
\draw    (300,260) -- (300,220) ;
\draw    (300,260) -- (320,280) ;
\draw    (320,280) -- (340,260) ;
\draw    (300,260) -- (320,200) ;
\draw    (340,260) -- (320,200) ;
\draw    (320,280) -- (320,200) ;

\draw  [color={rgb, 255:red, 74; green, 144; blue, 226 }  ,draw opacity=1 ][fill={rgb, 255:red, 74; green, 144; blue, 226 }  ,fill opacity=1 ] (307.86,215.64) .. controls (307.86,214.18) and (309.04,213) .. (310.5,213) .. controls (311.96,213) and (313.14,214.18) .. (313.14,215.64) .. controls (313.14,217.1) and (311.96,218.28) .. (310.5,218.28) .. controls (309.04,218.28) and (307.86,217.1) .. (307.86,215.64) -- cycle ;
\draw  [color={rgb, 255:red, 74; green, 144; blue, 226 }  ,draw opacity=1 ][fill={rgb, 255:red, 74; green, 144; blue, 226 }  ,fill opacity=1 ] (313.2,225.14) .. controls (313.2,223.68) and (314.38,222.5) .. (315.83,222.5) .. controls (317.29,222.5) and (318.47,223.68) .. (318.47,225.14) .. controls (318.47,226.6) and (317.29,227.78) .. (315.83,227.78) .. controls (314.38,227.78) and (313.2,226.6) .. (313.2,225.14) -- cycle ;
\draw  [color={rgb, 255:red, 74; green, 144; blue, 226 }  ,draw opacity=1 ][fill={rgb, 255:red, 74; green, 144; blue, 226 }  ,fill opacity=1 ] (321.39,225.5) .. controls (321.39,224.04) and (322.57,222.86) .. (324.03,222.86) .. controls (325.49,222.86) and (326.67,224.04) .. (326.67,225.5) .. controls (326.67,226.96) and (325.49,228.14) .. (324.03,228.14) .. controls (322.57,228.14) and (321.39,226.96) .. (321.39,225.5) -- cycle ;
\draw  [color={rgb, 255:red, 74; green, 144; blue, 226 }  ,draw opacity=1 ][fill={rgb, 255:red, 74; green, 144; blue, 226 }  ,fill opacity=1 ] (326.92,215.1) .. controls (326.92,213.64) and (328.1,212.46) .. (329.56,212.46) .. controls (331.02,212.46) and (332.2,213.64) .. (332.2,215.1) .. controls (332.2,216.56) and (331.02,217.74) .. (329.56,217.74) .. controls (328.1,217.74) and (326.92,216.56) .. (326.92,215.1) -- cycle ;
\draw    (220,130) -- (240,150) ;
\draw    (200,150) -- (220,130) ;
\draw    (240,190) -- (240,150) ;
\draw    (200,190) -- (200,150) ;
\draw    (200,190) -- (220,210) ;
\draw    (220,210) -- (240,190) ;
\draw    (200,190) -- (220,130) ;
\draw    (240,190) -- (220,130) ;
\draw    (220,210) -- (220,130) ;

\draw  [color={rgb, 255:red, 74; green, 144; blue, 226 }  ,draw opacity=1 ][fill={rgb, 255:red, 74; green, 144; blue, 226 }  ,fill opacity=1 ] (200.86,167.14) .. controls (200.86,165.68) and (202.04,164.5) .. (203.5,164.5) .. controls (204.96,164.5) and (206.14,165.68) .. (206.14,167.14) .. controls (206.14,168.6) and (204.96,169.78) .. (203.5,169.78) .. controls (202.04,169.78) and (200.86,168.6) .. (200.86,167.14) -- cycle ;
\draw  [color={rgb, 255:red, 74; green, 144; blue, 226 }  ,draw opacity=1 ][fill={rgb, 255:red, 74; green, 144; blue, 226 }  ,fill opacity=1 ] (203.2,187.14) .. controls (203.2,185.68) and (204.38,184.5) .. (205.83,184.5) .. controls (207.29,184.5) and (208.47,185.68) .. (208.47,187.14) .. controls (208.47,188.6) and (207.29,189.78) .. (205.83,189.78) .. controls (204.38,189.78) and (203.2,188.6) .. (203.2,187.14) -- cycle ;
\draw  [color={rgb, 255:red, 74; green, 144; blue, 226 }  ,draw opacity=1 ][fill={rgb, 255:red, 74; green, 144; blue, 226 }  ,fill opacity=1 ] (221.39,155.5) .. controls (221.39,154.04) and (222.57,152.86) .. (224.03,152.86) .. controls (225.49,152.86) and (226.67,154.04) .. (226.67,155.5) .. controls (226.67,156.96) and (225.49,158.14) .. (224.03,158.14) .. controls (222.57,158.14) and (221.39,156.96) .. (221.39,155.5) -- cycle ;
\draw  [color={rgb, 255:red, 74; green, 144; blue, 226 }  ,draw opacity=1 ][fill={rgb, 255:red, 74; green, 144; blue, 226 }  ,fill opacity=1 ] (226.92,145.1) .. controls (226.92,143.64) and (228.1,142.46) .. (229.56,142.46) .. controls (231.02,142.46) and (232.2,143.64) .. (232.2,145.1) .. controls (232.2,146.56) and (231.02,147.74) .. (229.56,147.74) .. controls (228.1,147.74) and (226.92,146.56) .. (226.92,145.1) -- cycle ;
\draw    (420,130) -- (440,150) ;
\draw    (400,150) -- (420,130) ;
\draw    (440,190) -- (440,150) ;
\draw    (400,190) -- (400,150) ;
\draw    (400,190) -- (420,210) ;
\draw    (420,210) -- (440,190) ;
\draw    (400,190) -- (420,130) ;
\draw    (440,190) -- (420,130) ;
\draw    (420,210) -- (420,130) ;

\draw  [color={rgb, 255:red, 74; green, 144; blue, 226 }  ,draw opacity=1 ][fill={rgb, 255:red, 74; green, 144; blue, 226 }  ,fill opacity=1 ] (407.86,145.64) .. controls (407.86,144.18) and (409.04,143) .. (410.5,143) .. controls (411.96,143) and (413.14,144.18) .. (413.14,145.64) .. controls (413.14,147.1) and (411.96,148.28) .. (410.5,148.28) .. controls (409.04,148.28) and (407.86,147.1) .. (407.86,145.64) -- cycle ;
\draw  [color={rgb, 255:red, 74; green, 144; blue, 226 }  ,draw opacity=1 ][fill={rgb, 255:red, 74; green, 144; blue, 226 }  ,fill opacity=1 ] (413.2,155.14) .. controls (413.2,153.68) and (414.38,152.5) .. (415.83,152.5) .. controls (417.29,152.5) and (418.47,153.68) .. (418.47,155.14) .. controls (418.47,156.6) and (417.29,157.78) .. (415.83,157.78) .. controls (414.38,157.78) and (413.2,156.6) .. (413.2,155.14) -- cycle ;
\draw  [color={rgb, 255:red, 74; green, 144; blue, 226 }  ,draw opacity=1 ][fill={rgb, 255:red, 74; green, 144; blue, 226 }  ,fill opacity=1 ] (431.89,187) .. controls (431.89,185.54) and (433.07,184.36) .. (434.53,184.36) .. controls (435.99,184.36) and (437.17,185.54) .. (437.17,187) .. controls (437.17,188.46) and (435.99,189.64) .. (434.53,189.64) .. controls (433.07,189.64) and (431.89,188.46) .. (431.89,187) -- cycle ;
\draw  [color={rgb, 255:red, 74; green, 144; blue, 226 }  ,draw opacity=1 ][fill={rgb, 255:red, 74; green, 144; blue, 226 }  ,fill opacity=1 ] (433.92,167.1) .. controls (433.92,165.64) and (435.1,164.46) .. (436.56,164.46) .. controls (438.02,164.46) and (439.2,165.64) .. (439.2,167.1) .. controls (439.2,168.56) and (438.02,169.74) .. (436.56,169.74) .. controls (435.1,169.74) and (433.92,168.56) .. (433.92,167.1) -- cycle ;
\draw    (320,60) -- (340,80) ;
\draw    (300,80) -- (320,60) ;
\draw    (340,120) -- (340,80) ;
\draw    (300,120) -- (300,80) ;
\draw    (300,120) -- (320,140) ;
\draw    (320,140) -- (340,120) ;
\draw    (300,120) -- (320,60) ;
\draw    (340,120) -- (320,60) ;
\draw    (320,140) -- (320,60) ;

\draw  [color={rgb, 255:red, 74; green, 144; blue, 226 }  ,draw opacity=1 ][fill={rgb, 255:red, 74; green, 144; blue, 226 }  ,fill opacity=1 ] (300.86,97.14) .. controls (300.86,95.68) and (302.04,94.5) .. (303.5,94.5) .. controls (304.96,94.5) and (306.14,95.68) .. (306.14,97.14) .. controls (306.14,98.6) and (304.96,99.78) .. (303.5,99.78) .. controls (302.04,99.78) and (300.86,98.6) .. (300.86,97.14) -- cycle ;
\draw  [color={rgb, 255:red, 74; green, 144; blue, 226 }  ,draw opacity=1 ][fill={rgb, 255:red, 74; green, 144; blue, 226 }  ,fill opacity=1 ] (303.2,117.14) .. controls (303.2,115.68) and (304.38,114.5) .. (305.83,114.5) .. controls (307.29,114.5) and (308.47,115.68) .. (308.47,117.14) .. controls (308.47,118.6) and (307.29,119.78) .. (305.83,119.78) .. controls (304.38,119.78) and (303.2,118.6) .. (303.2,117.14) -- cycle ;
\draw  [color={rgb, 255:red, 74; green, 144; blue, 226 }  ,draw opacity=1 ][fill={rgb, 255:red, 74; green, 144; blue, 226 }  ,fill opacity=1 ] (332.32,118.64) .. controls (332.32,117.18) and (333.5,116) .. (334.96,116) .. controls (336.42,116) and (337.6,117.18) .. (337.6,118.64) .. controls (337.6,120.1) and (336.42,121.28) .. (334.96,121.28) .. controls (333.5,121.28) and (332.32,120.1) .. (332.32,118.64) -- cycle ;
\draw  [color={rgb, 255:red, 74; green, 144; blue, 226 }  ,draw opacity=1 ][fill={rgb, 255:red, 74; green, 144; blue, 226 }  ,fill opacity=1 ] (333.32,96.3) .. controls (333.32,94.84) and (334.5,93.66) .. (335.96,93.66) .. controls (337.42,93.66) and (338.6,94.84) .. (338.6,96.3) .. controls (338.6,97.76) and (337.42,98.94) .. (335.96,98.94) .. controls (334.5,98.94) and (333.32,97.76) .. (333.32,96.3) -- cycle ;

\draw (295,354) node    {$\textcolor[rgb]{0.29,0.56,0.89}{6}$};
\draw (311,356) node    {$\textcolor[rgb]{0.29,0.56,0.89}{1}$};
\draw (327,356) node    {$\textcolor[rgb]{0.29,0.56,0.89}{3}$};
\draw (344,354) node    {$\textcolor[rgb]{0.29,0.56,0.89}{9}$};
\draw (295.06,237.44) node    {$\textcolor[rgb]{0.29,0.56,0.89}{6}$};
\draw (345,237) node    {$\textcolor[rgb]{0.29,0.56,0.89}{9}$};
\draw (306.06,272.44) node    {$\textcolor[rgb]{0.29,0.56,0.89}{4}$};
\draw (336.06,272.44) node    {$\textcolor[rgb]{0.29,0.56,0.89}{5}$};
\draw (245,171) node    {$\textcolor[rgb]{0.29,0.56,0.89}{9}$};
\draw (234.06,204.44) node    {$\textcolor[rgb]{0.29,0.56,0.89}{5}$};
\draw (215.06,173.44) node    {$\textcolor[rgb]{0.29,0.56,0.89}{2}$};
\draw (207.06,134.44) node    {$\textcolor[rgb]{0.29,0.56,0.89}{7}$};
\draw (395,165) node    {$\textcolor[rgb]{0.29,0.56,0.89}{6}$};
\draw (405.06,204.44) node    {$\textcolor[rgb]{0.29,0.56,0.89}{4}$};
\draw (424.06,170.44) node    {$\textcolor[rgb]{0.29,0.56,0.89}{2}$};
\draw (435.06,135.44) node    {$\textcolor[rgb]{0.29,0.56,0.89}{8}$};
\draw (326.06,107.44) node    {$\textcolor[rgb]{0.29,0.56,0.89}{2}$};
\draw (335.06,65.44) node    {$\textcolor[rgb]{0.29,0.56,0.89}{8}$};
\draw (315.06,107.44) node    {$\textcolor[rgb]{0.29,0.56,0.89}{2}$};
\draw (307.06,64.44) node    {$\textcolor[rgb]{0.29,0.56,0.89}{7}$};

\end{tikzpicture}
\end {figure}
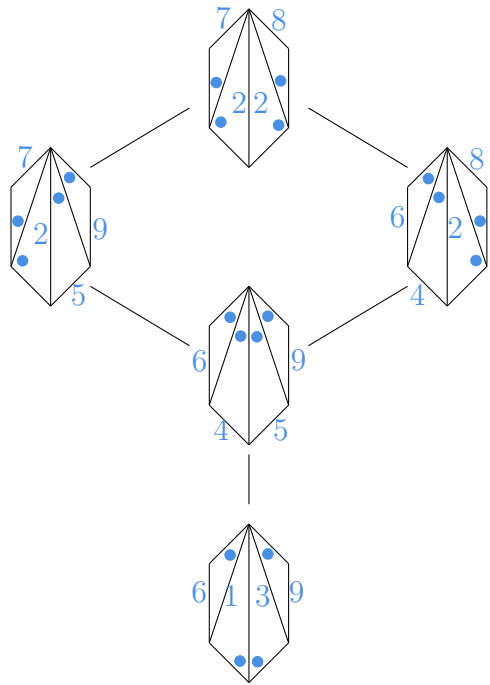

\section{$S$-paths}

\begin{defn}
       An \textit{$S$-path from $a$ to $b$} is an ordered (uncolored) selection $S$ of edges from the triangulation $\Delta_w$ subject to the following conditions: 

\begin{enumerate}[(1)]
    \item The edges in $S$ form a path from $a$ to $b$
    
    \item The number of edges in $S$ is $n+1$
    
    \item At least one edge from each triangle $\Delta_i$ in $\Sigma_w$ occurs in $S$.

\end{enumerate}
\end{defn}

Define the \textit{weight $x_S$ of $S$} by $x_S = \prod_{s \in S} x_s.$ Let $\mathbb{S}_w$ be the set of all $S$-paths from $a$ to $b$.

Again, the proof of the next result is given in Theorem \ref{duality_thm} below.

\begin{thm}
    Let $w$ be any word. Consider the set $\mathbb{S}_{w^{*}}$ of $S$-paths from $a^{*}$ to $b^{*}$ with edges taken from the dual triangulation $\Delta_{w}^{*}.$ Then the cluster variable $x_{w}$ can be written as 
$$x_{w} = \frac{1}{x_1 x_2 \dots x_n} \prod_{S \in \mathbb{S}_{w}^{*}} x_S.$$
\end{thm}

We now give the set $\mathbb{S}_w$ a poset structure.

A \textit{flip} of $S$ at diagonal $\delta_{i}$ is the local move that takes two edges in $S$ which are both incident to the same endpoint of the minimal quadrilateral $Q_i$ defined by the diagonal $\delta_{i}$, and replaces them with the other two edges of $Q_i$.

Directly below is the local picture of an $S$-path flip at a generic internal diagonal $\delta_{i}$.

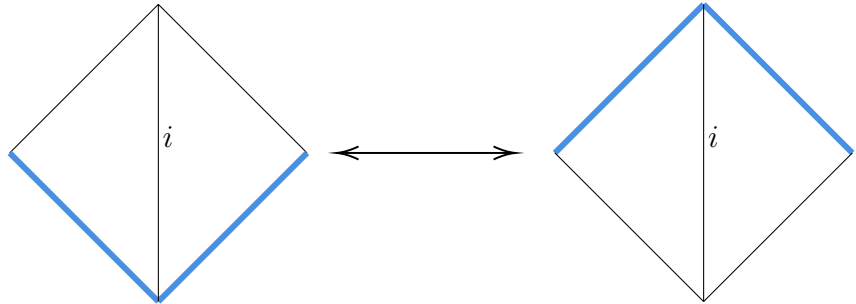
\begin {figure}[h!]
    \centering
    \caption{Flip of an $S$-path at diagonal $\delta_{i}$}
    \label{fig:flip3}
     \begin{tikzpicture}[x=0.75pt,y=0.75pt,yscale=-1,xscale=1]

\draw    (158,60) -- (233,135) ;
\draw    (83,135) -- (158,60) ;
\draw [color={rgb, 255:red, 74; green, 144; blue, 226 }  ,draw opacity=1 ][line width=2.25]    (83,135) -- (158,210) ;
\draw [color={rgb, 255:red, 74; green, 144; blue, 226 }  ,draw opacity=1 ][line width=2.25]    (158,210) -- (233,135) ;
\draw    (158,60) -- (158,210) ;
\draw [line width=0.75]    (250,135) -- (336,135) ;
\draw [shift={(338,135)}, rotate = 180] [color={rgb, 255:red, 0; green, 0; blue, 0 }  ][line width=0.75]    (10.93,-3.29) .. controls (6.95,-1.4) and (3.31,-0.3) .. (0,0) .. controls (3.31,0.3) and (6.95,1.4) .. (10.93,3.29)   ;
\draw [shift={(248,135)}, rotate = 0] [color={rgb, 255:red, 0; green, 0; blue, 0 }  ][line width=0.75]    (10.93,-3.29) .. controls (6.95,-1.4) and (3.31,-0.3) .. (0,0) .. controls (3.31,0.3) and (6.95,1.4) .. (10.93,3.29)   ;
\draw [color={rgb, 255:red, 74; green, 144; blue, 226 }  ,draw opacity=1 ][line width=2.25]    (433,60) -- (508,135) ;
\draw [color={rgb, 255:red, 74; green, 144; blue, 226 }  ,draw opacity=1 ][line width=2.25]    (358,135) -- (433,60) ;
\draw    (358,135) -- (433,210) ;
\draw    (433,210) -- (508,135) ;
\draw    (433,60) -- (433,210) ;

\draw (159,120) node [anchor=north west][inner sep=0.75pt]    {$i$};
\draw (434,120) node [anchor=north west][inner sep=0.75pt]    {$i$};

\end{tikzpicture}

\end {figure}

 Let $l_i$ be the endpoint of $\delta_{i}$ to the left of $\gamma_{a \rightarrow b}$, and $r_i$ the endpoint to the right. An \textit{up-flip} of $S$ at diagonal $\delta_{i}$ is a flip that meets either of the following two conditions:

  \begin{enumerate}[(1)]
        \item $i$ is odd, and the two edges in $S$ incident to $r_i$ are replaced with the two edges in $Q_i$ incident to $l_i$, or
        
        \item $i$ is even, and the two edges in $S$ incident to $l_i$ are replaced with the two edges in $Q_i$ incident to $r_i$.
    \end{enumerate}

\begin{defn}
       The \textit{minimal element} $S_{-}$ of $\mathbb{S}_w$ is the unique $S$-path that starts with the boundary edge $\delta_{2n+1}$, and only uses internal diagonals as edges, except for the first and last (boundary) edges. The \textit{maximal element} $S_{+}$ of $\mathbb{S}_w$ is the unique $S$-path that starts with the edge $\delta_{2n}$, and only uses internal diagonals as edges, except for the first and last (boundary) edges.
\end{defn}

\begin{defn}
       The \textit{poset structure on $\mathbb{S}_w$} is defined as follows. The unique minimal element of $\mathbb{S}_w$ is the minimal $S$-path $S_{-},$ and the unique maximal element is $S_{+}.$ An $S$-path $S_2$ covers another $S$-path $S_1$ if there exists a diagonal $\delta_{i}$ such that $S_2$ can be obtained from $S_1$ by performing a single up-flip of $S_1$ at $\delta_{i}$.
\end{defn}

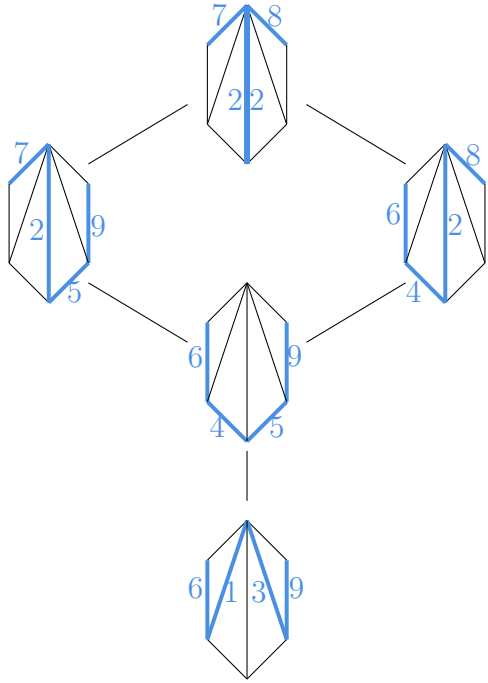
\begin {figure}[h!]
    \centering
    \caption{The poset $\mathbb{S}_{bb}$}
    \label{fig:S_bb}
    \begin{tikzpicture}[x=0.75pt,y=0.75pt,yscale=-1,xscale=1]

\draw    (340,324.97) -- (360,344.97) ;
\draw    (320,344.97) -- (340,324.97) ;
\draw [color={rgb, 255:red, 74; green, 144; blue, 226 }  ,draw opacity=1 ][line width=1.5]    (360,384.97) -- (360,344.97) ;
\draw [color={rgb, 255:red, 74; green, 144; blue, 226 }  ,draw opacity=1 ][line width=1.5]    (320,384.97) -- (320,344.97) ;
\draw    (320,384.97) -- (340,404.97) ;
\draw    (340,404.97) -- (360,384.97) ;
\draw [color={rgb, 255:red, 74; green, 144; blue, 226 }  ,draw opacity=1 ][line width=1.5]    (320,384.97) -- (340,324.97) ;
\draw [color={rgb, 255:red, 74; green, 144; blue, 226 }  ,draw opacity=1 ][line width=1.5]    (360,384.97) -- (340,324.97) ;
\draw    (340,404.97) -- (340,324.97) ;
\draw    (340,314.97) -- (340,289.84) ;
\draw    (370,234.97) -- (420,204.97) ;
\draw    (260,204.97) -- (310,234.97) ;
\draw    (370,114.97) -- (420,144.97) ;
\draw    (260,144.97) -- (310,114.97) ;
\draw    (340,204.97) -- (360,224.97) ;
\draw    (320,224.97) -- (340,204.97) ;
\draw [color={rgb, 255:red, 74; green, 144; blue, 226 }  ,draw opacity=1 ][line width=1.5]    (360,264.97) -- (360,224.97) ;
\draw [color={rgb, 255:red, 74; green, 144; blue, 226 }  ,draw opacity=1 ][line width=1.5]    (320,264.97) -- (320,224.97) ;
\draw [color={rgb, 255:red, 74; green, 144; blue, 226 }  ,draw opacity=1 ][line width=1.5]    (320,264.97) -- (340,284.97) ;
\draw [color={rgb, 255:red, 74; green, 144; blue, 226 }  ,draw opacity=1 ][line width=1.5]    (340,284.97) -- (360,264.97) ;
\draw    (320,264.97) -- (340,204.97) ;
\draw    (360,264.97) -- (340,204.97) ;
\draw    (340,284.97) -- (340,204.97) ;
\draw    (240,134.97) -- (260,154.97) ;
\draw [color={rgb, 255:red, 74; green, 144; blue, 226 }  ,draw opacity=1 ][line width=1.5]    (220,154.97) -- (240,134.97) ;
\draw [color={rgb, 255:red, 74; green, 144; blue, 226 }  ,draw opacity=1 ][line width=1.5]    (260,194.97) -- (260,154.97) ;
\draw    (220,194.97) -- (220,154.97) ;
\draw    (220,194.97) -- (240,214.97) ;
\draw [color={rgb, 255:red, 74; green, 144; blue, 226 }  ,draw opacity=1 ][line width=1.5]    (240,214.97) -- (260,194.97) ;
\draw    (220,194.97) -- (240,134.97) ;
\draw    (260,194.97) -- (240,134.97) ;
\draw [color={rgb, 255:red, 74; green, 144; blue, 226 }  ,draw opacity=1 ][line width=1.5]    (240,214.97) -- (240,134.97) ;
\draw [color={rgb, 255:red, 74; green, 144; blue, 226 }  ,draw opacity=1 ][line width=1.5]    (440,134.97) -- (460,154.97) ;
\draw    (420,154.97) -- (440,134.97) ;
\draw    (460,194.97) -- (460,154.97) ;
\draw [color={rgb, 255:red, 74; green, 144; blue, 226 }  ,draw opacity=1 ][line width=1.5]    (420,194.97) -- (420,154.97) ;
\draw [color={rgb, 255:red, 74; green, 144; blue, 226 }  ,draw opacity=1 ][line width=1.5]    (420,194.97) -- (440,214.97) ;
\draw    (440,214.97) -- (460,194.97) ;
\draw    (420,194.97) -- (440,134.97) ;
\draw    (460,194.97) -- (440,134.97) ;
\draw [color={rgb, 255:red, 74; green, 144; blue, 226 }  ,draw opacity=1 ][line width=1.5]    (440,214.97) -- (440,134.97) ;
\draw [color={rgb, 255:red, 74; green, 144; blue, 226 }  ,draw opacity=1 ][line width=1.5]    (340,64.97) -- (360,84.97) ;
\draw [color={rgb, 255:red, 74; green, 144; blue, 226 }  ,draw opacity=1 ][line width=1.5]    (320,84.97) -- (340,64.97) ;
\draw    (360,124.97) -- (360,84.97) ;
\draw    (320,124.97) -- (320,84.97) ;
\draw    (320,124.97) -- (340,144.97) ;
\draw    (340,144.97) -- (360,124.97) ;
\draw    (320,124.97) -- (340,64.97) ;
\draw    (360,124.97) -- (340,64.97) ;
\draw [color={rgb, 255:red, 74; green, 144; blue, 226 }  ,draw opacity=1 ][line width=2.25]    (340,144.97) -- (340,64.97) ;

\draw (314,359) node    {$\textcolor[rgb]{0.29,0.56,0.89}{6}$};
\draw (332,361) node    {$\textcolor[rgb]{0.29,0.56,0.89}{1}$};
\draw (346,361) node    {$\textcolor[rgb]{0.29,0.56,0.89}{3}$};
\draw (365,359) node    {$\textcolor[rgb]{0.29,0.56,0.89}{9}$};
\draw (314.06,242.44) node    {$\textcolor[rgb]{0.29,0.56,0.89}{6}$};
\draw (364,242) node    {$\textcolor[rgb]{0.29,0.56,0.89}{9}$};
\draw (325.06,277.44) node    {$\textcolor[rgb]{0.29,0.56,0.89}{4}$};
\draw (355.06,277.44) node    {$\textcolor[rgb]{0.29,0.56,0.89}{5}$};
\draw (265,176) node    {$\textcolor[rgb]{0.29,0.56,0.89}{9}$};
\draw (253.06,209.44) node    {$\textcolor[rgb]{0.29,0.56,0.89}{5}$};
\draw (234.06,178.44) node    {$\textcolor[rgb]{0.29,0.56,0.89}{2}$};
\draw (226.06,139.44) node    {$\textcolor[rgb]{0.29,0.56,0.89}{7}$};
\draw (414,170) node    {$\textcolor[rgb]{0.29,0.56,0.89}{6}$};
\draw (424.06,209.44) node    {$\textcolor[rgb]{0.29,0.56,0.89}{4}$};
\draw (445.06,175.44) node    {$\textcolor[rgb]{0.29,0.56,0.89}{2}$};
\draw (454.06,140.44) node    {$\textcolor[rgb]{0.29,0.56,0.89}{8}$};
\draw (345.06,112.44) node    {$\textcolor[rgb]{0.29,0.56,0.89}{2}$};
\draw (354.06,70.44) node    {$\textcolor[rgb]{0.29,0.56,0.89}{8}$};
\draw (334.06,112.44) node    {$\textcolor[rgb]{0.29,0.56,0.89}{2}$};
\draw (326.06,69.44) node    {$\textcolor[rgb]{0.29,0.56,0.89}{7}$};

\end{tikzpicture}
  
\end {figure}

\begin{ex}
    Fix $w=ab$. The poset of $S$-paths on the dual triangulated surface $\Sigma^{*}_{ab} = \Sigma_{bb}$ is displayed in Figure \ref{fig:S_bb}. Note that the ``middle'' internal diagonal in $S_{+}$ is used twice.
\end{ex}

\newpage 

We now give the analogue of Proposition \ref{bij_PAT}.

\begin{prop} \label{bij_LBS}
     Fix the word $w$. There are poset isomorphisms 

\begin{center}
\begin{tikzcd}
                                   & \mathbb{L}_w \arrow[ld] \arrow[rd] &                                    \\
\mathbb{B}_w \arrow[rr] \arrow[ru] &                                    & \mathbb{S}_w \arrow[lu] \arrow[ll]
\end{tikzcd}
\end{center}

which respect the additional node structure of each poset.
\end{prop}

\begin{proof} 
    Let $A(\widetilde{G_w})$ be the set of angles in $\widetilde{G_w}$ that are incident to a diagonal and a side of $\widetilde{G_w}$. Let $E(G_w)$ be the edges of $G_w \subset \widetilde{G_w}$. Let $A(\Sigma_w)$ be the set of angles in $\Sigma_w$ which are neither incident to $a$ nor $b$. Recall from \cite{yurikusa2019cluster} the canonical surjection $A(\widetilde{G_w}) \longrightarrow E(G_w),$ which induces the bijection $E(G_w) \longrightarrow A(\Sigma_w).$ The map $\mathbb{L}_w \longrightarrow \mathbb{B}_w$ is defined by taking the preimage of the lattice path $L$ under the map $A(\widetilde{G_w}) \longrightarrow E(G_w),$ and then folding the result to obtain an element in $A(\Sigma_w)$. That this map is well-defined can be seen by induction on the number of tiles $n$. That covering relations are preserved is straightforward to check. The inverse map is defined by unfolding and then applying the map $A(\widetilde{G_w}) \longrightarrow E(G_w).$ It follows that $\mathbb{L}_w \longrightarrow \mathbb{B}_w$ is a poset isomorphism respecting the additional node structure.

    The map $\mathbb{L}_w \longrightarrow \mathbb{S}_w$ is the folding map, additionally keeping track of the images of the edges selected by $L \in \mathbb{L}_w$ after each fold. One can check using induction that this map is well-defined. Again, it is straightforward to see that covering relations are preserved. The inverse map is induced by the unfolding map. It follows that $\mathbb{L}_w \longrightarrow \mathbb{S}_w$ is a poset isomorphism respecting the additional node structure.
    
    Finally, one can check that the map $\mathbb{B}_w \longrightarrow \mathbb{S}_w$ is the composition of the two poset isomorphisms above. 
\end{proof}

We illustrate Proposition \ref{bij_LBS} with the three dual posets $\mathbb{L}_{bb}, \mathbb{B}_{bb},$ and $\mathbb{S}_{bb}$ from the running example. By Lemma 3.2 in \cite{yurikusa2019cluster} there is a bijection between the angles in $\Sigma_w$ and the edges in $G_w$ induced by identifying certain pairs of angles in $\widetilde{G_w}.$ The pairs of angles that are identified are those that are opposite one another in the quadrilateral determined by two consecutive tiles of $\widetilde{G_w}.$ Any pair of angles in $\widetilde{G_w}$ which have been identified correspond to a single internal angle in $\Sigma_w$. 

Thus, given a lattice path $L \in \mathbb{L}_w,$ we can associate to it a collection of angles in $\widetilde{G}$ and fold the result to obtain a lattice path of angles $B$ in $\Sigma_w$. 

\begin {figure}[h!]
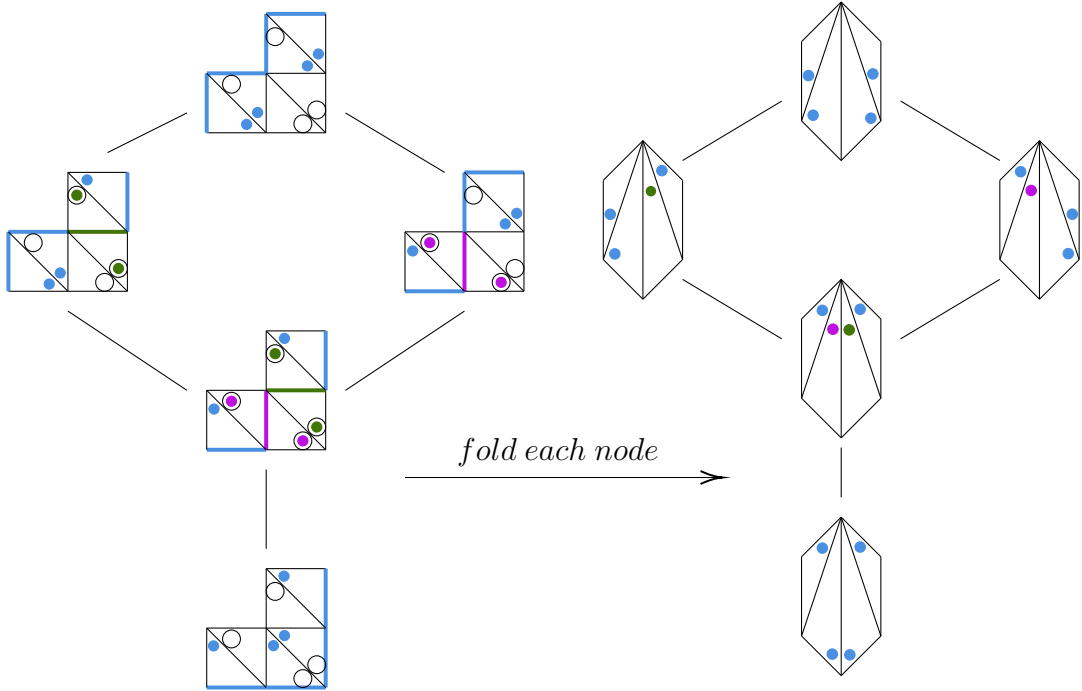

    \centering
    \caption{The map $\mathbb{L}_w \longrightarrow \mathbb{B}_{bb}$ via angle identification and folding}
    \label{fig:L_to_B}



\end {figure}
\newpage

The map $\mathbb{B}_w \longrightarrow \mathbb{S}_w$ is defined by taking the collection of edges which are opposite some vertex in $S \in \mathbb{S}_w$. See Figure \ref{fig:B_to_S} below.

\begin {figure}[h!]
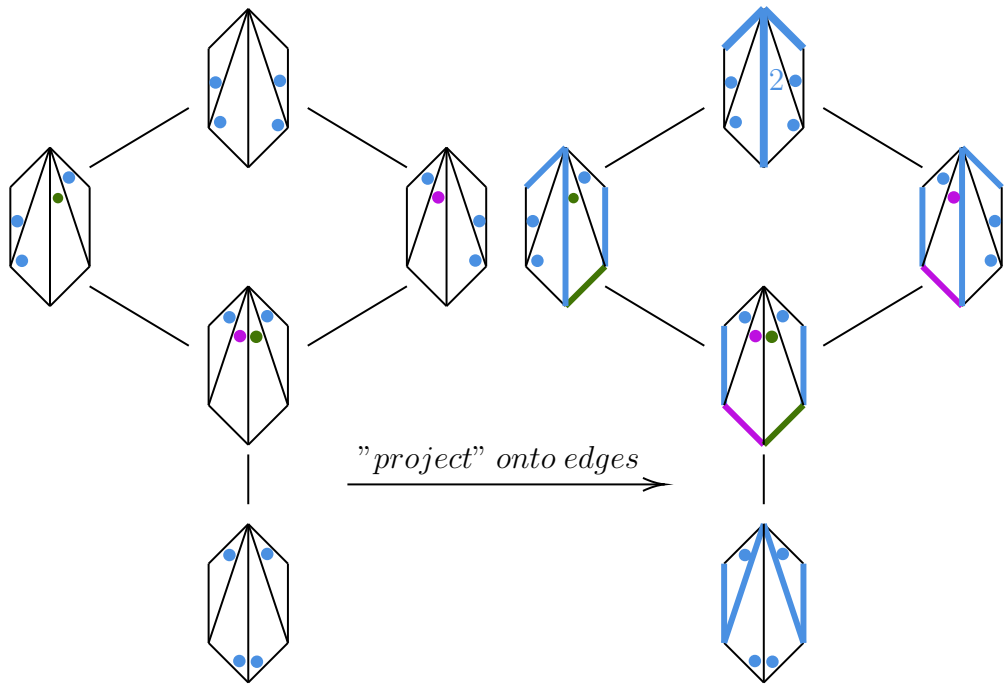

    \centering
    \caption{The map $\mathbb{B}_{bb} \longrightarrow \mathbb{S}_{bb}$ via associating edges to angles}
    \label{fig:B_to_S}

\tikzset{every picture/.style={line width=0.75pt}} 

%

\end {figure}

\newpage

The map $\mathbb{L}_w \longrightarrow \mathbb{S}_w$ is defined by folding and is the composition of the previous two.

\section{Expansion Duality}

\begin{thm} \label{duality_thm}
    Fix the word $w$. 

\begin{enumerate}[(a)]
    
    \item There is a explicit bijection $\text{Tree}(w)^{*} \longrightarrow \text{Tree}(w^{*})$ that preserves additional node structure, and weights of leaves. 
    
    \item There is a weight-preserving bijection $\text{Res}(w)^{*} \longrightarrow \text{Res}(w^{*})$ respecting additional node structure. 
    
    \item The Laurent polynomial $x_{w}^{*}$ is a cluster variable in $\mathcal{A}(\Sigma)_{w^{*}}$, and is equal to $x_{w}^{*} = x_{w^{*}}.$ 
    
     \item There are explicit isomorphisms of distributive lattices 
    $$\mathbb{P}_w \xrightarrow{\sim} \mathbb{L}_{w^{*}} , \mathbb{A}_w \xrightarrow{\sim} \mathbb{B}_{w^{*}}, \text{ and } \mathbb{T}_w \xrightarrow{\sim} \mathbb{S}_{w^{*}}$$
    
respecting the additional structure of each lattice, and making the following diagram commute.
    
    
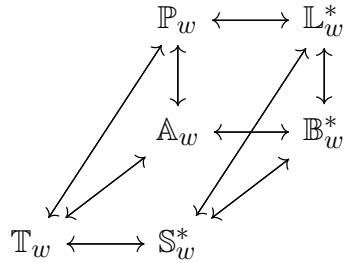
\begin{center}
\begin{tikzcd}
                                                & \mathbb{P}_w \arrow[ldd] \arrow[r] \arrow[d]        & \mathbb{L}_{w}^{*} \arrow[ldd] \arrow[l] \arrow[d] \\
                                                & \mathbb{A}_w \arrow[u] \arrow[r] \arrow[ld]         & \mathbb{B}_{w}^{*} \arrow[u] \arrow[l] \arrow[ld]  \\
\mathbb{T}_{w} \arrow[ru] \arrow[r] \arrow[ruu] & \mathbb{S}_{w}^{*} \arrow[ru] \arrow[l] \arrow[ruu] &                                                   
\end{tikzcd}

\end{center}

    Namely, corresponding nodes in the six posets have the same weight, except for the nodes of the $T$-path expansion poset. In this case, node weights have an additional factor of $\frac{1}{x_1 x_2 \dots x_n}$ that is not preset in any of the node weights for the other five expansion posets.

    \item The cluster variable $x_w$ can be written as
    $$x_{w} = \frac{1}{ x_1 x_2 \dots x_n} \sum_{L \in \mathbb{L}_{w^{*}} } x_L = \frac{1}{ x_1 x_2 \dots x_n} \sum_{B \in \mathbb{B}_{w^{*}} } x_B = \frac{1}{ x_1 x_2 \dots x_n} \sum_{S \in \mathbb{S}_{w^{*}} } x_S.$$

    \item There are isomorphisms of distributive lattices $D_w \cong \mathcal{I}(C_w) \cong \mathbb{P}_w$ and $D_{w^{*}} \cong \mathcal{I}(C_{w}^{*}) \cong \mathbb{L}_{w}.$ Thus, $\mathbb{P}_w$ and $\mathbb{L}_w$ are dual to one another in the sense of distributive lattices (see Definition \ref{dual_dist_def}).

\end{enumerate}

\end{thm}

\begin{proof}
    
    \begin{enumerate}[(a)]
        \item Let $t \in \text{Tree}(w)^{*}$. The map $\text{Tree}(w)^{*} \longrightarrow \text{Tree}(w^{*})$ is defined by the application of the triangle map to each node of the input tree $t$, except that we must additionally specify how to transform arcs $\delta \in \Sigma_w$ which are not contained in the triangulation $\Delta_w$ into arcs $\delta^{*} \in \Sigma_{w}^{*}$. 
        
        Suppose $\delta \notin \Delta_w$ is an arc in a diagram which is a node in $t$, and that $\delta$ crosses the triangles $\Delta_i , \dots , \Delta_j$. The obvious one-to-one correspondence induced by the triangle map between angles in $\Sigma_w$ and angles in $\Sigma_{w}^{*}$ gives a natural candidate for the image of the arc $\delta$.
        
        Namely, suppose the endpoints of $\delta$ are $v_i$ and $v_j$, corresponding to the angles $\alpha_i$ an $\alpha_j$ in $\Delta_i$ and $\Delta_j$, respectively. Let $\alpha_{i}^{*}$ and $\alpha_{j}^{*}$ be the angles which are the respective images of $\alpha_i$ and $\alpha_j$ under the triangle map. Let $v_{i}^{*}$ be the vertex that $\alpha_{i}^{*}$ is matched to, and let $v_{j}^{*}$ be the vertex that $\alpha_{j}^{*}$ is matched to.  Then $\delta^{*}$ is the arc that 
    \begin{enumerate}[(1)]
        \item starts at the vertex $v_{i}^{*},$
        
        \item ends at the vertex $v_{j}^{*},$
        
        \item passes through the center of precisely those triangles in $\Sigma_{w}^{*}$ which are the images of the triangles $\Delta_i , \dots \Delta_j$ under the triangle map, and 
        
        \item has as its only intersections (besides endpoints) the midpoint of each internal diagonal it crosses.
    \end{enumerate}

    For instance, suppose the arc $\delta$ starts at $a$, is contained in the triangulated subpolygon $[\Delta_0 , \Delta_1 , \dots , \Delta_{j}],$ and ends at vertex $v_j$ of $\Delta_{j}.$ The output $\delta^{*}$ for the two subcases $j$ even and $j$ odd are shown in Figure \ref{fig:arc_transform}. 
    
    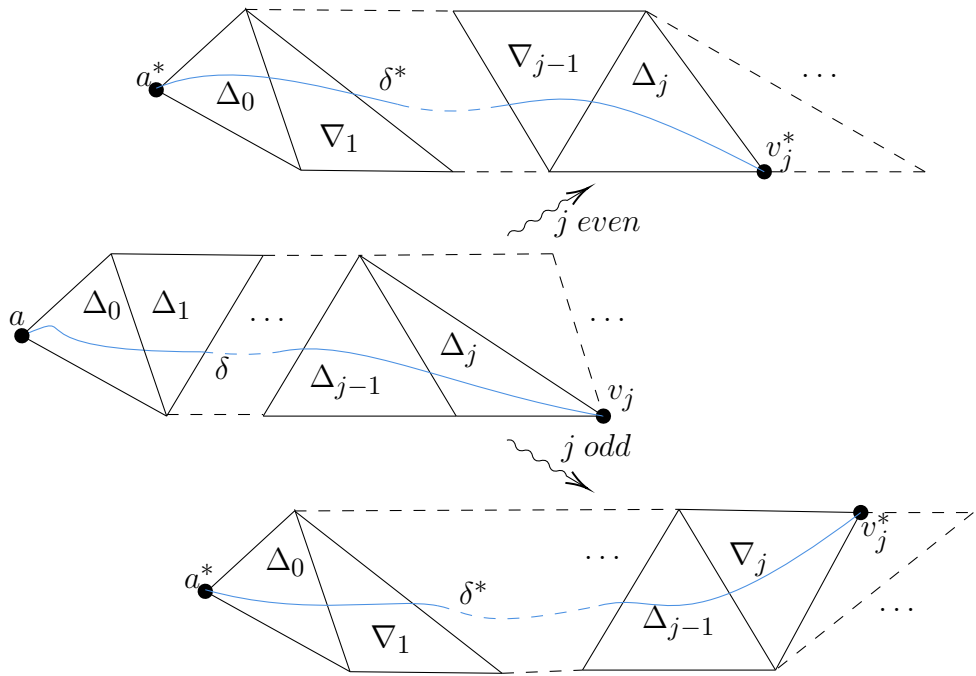
\begin {figure}[h!]
    \centering
    \caption{The arc $\delta$, and the output arc $\delta^{*}$ for both $j$ odd and $j$ even}
    \label{fig:arc_transform}
    \begin{tikzpicture}[x=0.75pt,y=0.75pt,yscale=-1,xscale=1]

\draw    (226.99,130.94) -- (275.67,212.07) ;
\draw    (178.31,212.07) -- (226.99,130.94) ;
\draw    (178.31,212.07) -- (275.67,212.07) ;
\draw [color={rgb, 255:red, 74; green, 144; blue, 226 }  ,draw opacity=1 ]   (56.61,171.51) .. controls (91.09,156.7) and (45.25,179.62) .. (147.07,179.62) ;
\draw    (301.63,121.21) .. controls (302.2,118.92) and (303.62,118.06) .. (305.91,118.63) .. controls (308.2,119.2) and (309.63,118.35) .. (310.2,116.06) .. controls (310.77,113.77) and (312.2,112.92) .. (314.49,113.49) .. controls (316.78,114.06) and (318.21,113.21) .. (318.78,110.92) .. controls (319.35,108.63) and (320.77,107.77) .. (323.06,108.34) .. controls (325.35,108.91) and (326.78,108.06) .. (327.35,105.77) .. controls (327.92,103.48) and (329.35,102.63) .. (331.64,103.2) -- (333.62,102.01) -- (340.48,97.9) ;
\draw [shift={(342.19,96.87)}, rotate = 509.04] [color={rgb, 255:red, 0; green, 0; blue, 0 }  ][line width=0.75]    (10.93,-3.29) .. controls (6.95,-1.4) and (3.31,-0.3) .. (0,0) .. controls (3.31,0.3) and (6.95,1.4) .. (10.93,3.29)   ;
\draw    (101.64,130.13) -- (178.31,130.94) ;
\draw    (101.64,130.13) -- (129.63,212.07) ;
\draw    (101.64,130.13) -- (56.61,171.51) ;
\draw    (129.63,212.07) -- (56.61,171.51) ;
\draw [shift={(56.61,171.51)}, rotate = 209.05] [color={rgb, 255:red, 0; green, 0; blue, 0 }  ][fill={rgb, 255:red, 0; green, 0; blue, 0 }  ][line width=0.75]      (0, 0) circle [x radius= 3.35, y radius= 3.35]   ;
\draw    (226.99,130.94) -- (349.9,212.07) ;
\draw [shift={(349.9,212.07)}, rotate = 33.43] [color={rgb, 255:red, 0; green, 0; blue, 0 }  ][fill={rgb, 255:red, 0; green, 0; blue, 0 }  ][line width=0.75]      (0, 0) circle [x radius= 3.35, y radius= 3.35]   ;
\draw    (275.67,212.07) -- (349.9,212.07) ;
\draw    (129.63,212.07) -- (178.31,130.94) ;
\draw [color={rgb, 255:red, 74; green, 144; blue, 226 }  ,draw opacity=1 ]   (186.42,179.62) .. controls (224.42,171.37) and (275.67,199.09) .. (349.9,212.07) ;
\draw [color={rgb, 255:red, 74; green, 144; blue, 226 }  ,draw opacity=1 ] [dash pattern={on 4.5pt off 4.5pt}]  (147.07,179.62) .. controls (160.32,180.84) and (173.04,181.65) .. (186.42,179.62) ;
\draw  [dash pattern={on 4.5pt off 4.5pt}]  (178.31,130.94) -- (226.99,130.94) ;
\draw  [dash pattern={on 4.5pt off 4.5pt}]  (129.63,211.26) -- (178.31,211.26) ;
\draw    (274.04,7.62) -- (322.72,88.75) ;
\draw    (322.72,88.75) -- (371.4,7.62) ;
\draw    (274.04,7.62) -- (371.4,7.62) ;
\draw    (197.37,87.94) -- (274.04,88.75) ;
\draw    (169.38,6) -- (197.37,87.94) ;
\draw    (169.38,6.81) -- (124.36,48.19) ;
\draw    (197.37,87.94) -- (124.36,47.38) ;
\draw [shift={(124.36,47.38)}, rotate = 209.05] [color={rgb, 255:red, 0; green, 0; blue, 0 }  ][fill={rgb, 255:red, 0; green, 0; blue, 0 }  ][line width=0.75]      (0, 0) circle [x radius= 3.35, y radius= 3.35]   ;
\draw    (274.04,88.75) -- (169.38,6.81) ;
\draw  [dash pattern={on 4.5pt off 4.5pt}]  (169.38,6.81) -- (274.04,7.62) ;
\draw  [dash pattern={on 4.5pt off 4.5pt}]  (274.04,88.75) -- (322.72,88.75) ;
\draw    (388.03,259.13) -- (436.71,340.26) ;
\draw    (339.35,340.26) -- (388.03,259.13) ;
\draw    (339.35,340.26) -- (436.71,340.26) ;
\draw    (388.03,259.13) -- (479.71,260.75) ;
\draw [shift={(479.71,260.75)}, rotate = 1.01] [color={rgb, 255:red, 0; green, 0; blue, 0 }  ][fill={rgb, 255:red, 0; green, 0; blue, 0 }  ][line width=0.75]      (0, 0) circle [x radius= 3.35, y radius= 3.35]   ;
\draw    (436.71,340.26) -- (479.71,260.75) ;
\draw [color={rgb, 255:red, 74; green, 144; blue, 226 }  ,draw opacity=1 ]   (347.47,307.81) .. controls (381,299.02) and (386.95,330.93) .. (479.71,260.75) ;
\draw [color={rgb, 255:red, 74; green, 144; blue, 226 }  ,draw opacity=1 ] [dash pattern={on 4.5pt off 4.5pt}]  (266.34,307.81) .. controls (283.1,311.19) and (286.62,320.65) .. (347.47,307.81) ;
\draw  [dash pattern={on 4.5pt off 4.5pt}]  (194.13,259.94) -- (388.03,259.13) ;
\draw  [dash pattern={on 4.5pt off 4.5pt}]  (298.79,341.88) -- (339.35,340.26) ;
\draw [color={rgb, 255:red, 74; green, 144; blue, 226 }  ,draw opacity=1 ]   (124.36,46.57) .. controls (158.84,31.76) and (203.59,44.81) .. (246.05,54.68) ;
\draw [color={rgb, 255:red, 74; green, 144; blue, 226 }  ,draw opacity=1 ] [dash pattern={on 4.5pt off 4.5pt}]  (246.05,54.68) .. controls (263.5,59.55) and (282.29,58.33) .. (293.11,56.3) ;
\draw    (222.12,341.07) -- (298.79,341.88) ;
\draw    (194.13,259.13) -- (222.12,341.07) ;
\draw    (194.13,259.94) -- (149.1,301.32) ;
\draw    (222.12,341.07) -- (149.1,300.5) ;
\draw [shift={(149.1,300.5)}, rotate = 209.05] [color={rgb, 255:red, 0; green, 0; blue, 0 }  ][fill={rgb, 255:red, 0; green, 0; blue, 0 }  ][line width=0.75]      (0, 0) circle [x radius= 3.35, y radius= 3.35]   ;
\draw    (298.79,341.88) -- (194.13,259.94) ;
\draw [color={rgb, 255:red, 74; green, 144; blue, 226 }  ,draw opacity=1 ]   (149.1,299.69) .. controls (199.81,314.16) and (243.89,303.34) .. (266.34,307.81) ;
\draw    (301.22,224.24) .. controls (303.51,223.67) and (304.94,224.52) .. (305.51,226.81) .. controls (306.08,229.1) and (307.51,229.96) .. (309.8,229.39) .. controls (312.09,228.82) and (313.51,229.67) .. (314.08,231.96) .. controls (314.65,234.25) and (316.08,235.1) .. (318.37,234.53) .. controls (320.66,233.96) and (322.09,234.81) .. (322.66,237.1) .. controls (323.23,239.39) and (324.66,240.25) .. (326.95,239.68) .. controls (329.24,239.11) and (330.66,239.96) .. (331.23,242.25) -- (333.21,243.44) -- (340.07,247.55) ;
\draw [shift={(341.79,248.58)}, rotate = 210.96] [color={rgb, 255:red, 0; green, 0; blue, 0 }  ][line width=0.75]    (10.93,-3.29) .. controls (6.95,-1.4) and (3.31,-0.3) .. (0,0) .. controls (3.31,0.3) and (6.95,1.4) .. (10.93,3.29)   ;
\draw  [dash pattern={on 4.5pt off 4.5pt}]  (324.34,130.13) -- (226.99,130.94) ;
\draw  [dash pattern={on 4.5pt off 4.5pt}]  (324.34,130.13) -- (349.9,212.07) ;
\draw    (371.4,7.62) -- (431.03,88.75) ;
\draw [shift={(431.03,88.75)}, rotate = 53.68] [color={rgb, 255:red, 0; green, 0; blue, 0 }  ][fill={rgb, 255:red, 0; green, 0; blue, 0 }  ][line width=0.75]      (0, 0) circle [x radius= 3.35, y radius= 3.35]   ;
\draw    (322.72,88.75) -- (431.03,88.75) ;
\draw [color={rgb, 255:red, 74; green, 144; blue, 226 }  ,draw opacity=1 ]   (293.11,56.3) .. controls (331.65,50.62) and (344.63,43.32) .. (431.03,88.75) ;
\draw  [dash pattern={on 4.5pt off 4.5pt}]  (512.16,88.75) -- (371.4,7.62) ;
\draw  [dash pattern={on 4.5pt off 4.5pt}]  (512.16,88.75) -- (431.03,88.75) ;
\draw  [dash pattern={on 4.5pt off 4.5pt}]  (536.5,260.75) -- (436.71,340.26) ;
\draw  [dash pattern={on 4.5pt off 4.5pt}]  (479.71,260.75) -- (536.5,260.75) ;

\draw (181.55,163.39) node    {$\dotsc $};
\draw (54.18,162.58) node    {$a$};
\draw (359.23,203.15) node    {$v_{j}$};
\draw (121.92,39.26) node    {$a^{*}$};
\draw (439.55,78.61) node    {$v^{*}_{j}$};
\draw (350.71,283.47) node    {$\dotsc $};
\draw (487.42,270.89) node    {$v^{*}_{j}$};
\draw (146.26,292.8) node    {$a^{*}$};
\draw (347.16,114.71) node    {$j\ even$};
\draw (346.35,228.3) node    {$j\ odd$};
\draw (373.83,42.51) node    {$\Delta _{j}$};
\draw (319.88,32.77) node    {$\nabla _{j-1}$};
\draw (164.52,50.62) node    {$\Delta _{0}$};
\draw (216.85,71.72) node    {$\nabla _{1}$};
\draw (189.67,285.9) node    {$\Delta _{0}$};
\draw (242.81,324.84) node    {$\nabla _{1}$};
\draw (423.32,284.28) node    {$\nabla _{j}$};
\draw (388.03,316.73) node    {$\Delta _{j-1}$};
\draw (353.15,163.39) node    {$\dotsc $};
\draw (460.24,40.07) node    {$\dotsc $};
\draw (499.18,309.43) node    {$\dotsc $};
\draw (97.18,156.09) node    {$\Delta _{0}$};
\draw (132.06,156.09) node    {$\Delta _{1}$};
\draw (220.09,195.03) node    {$\Delta _{j-1}$};
\draw (277.29,178.81) node    {$\Delta _{j}$};
\draw (244.02,41.29) node    {$\delta ^{*}$};
\draw (284.59,300.91) node    {$\delta ^{*}$};
\draw (157.62,187.73) node    {$\delta $};

\end{tikzpicture}

\end {figure}

        This map is well-defined by induction on $n$, and that edge weights of leaves are preserved is obvious. The inverse map $\text{Tree}(w)^{*} \longrightarrow \text{Tree}(w^{*})$ is defined similarly, and the composition is the identity. Hence the map in question is an invertible bijection respecting additional node structure as claimed.

         \item This follows directly from (a), since $\text{Res}(w)$ is equal to the union of the leaves in the trees in $\text{Tree}(w)$, and $\text{Res}(w)^{*}$ is equal to the union of the leaves in the trees in $\text{Tree}(w)^{*}.$
         
          \item From (b) we have 
          $$x_{w}^{*} = \frac{1}{x_1 x_2 \cdots x_n} \sum_{r^{*} \in \text{Res}(w)^{*}} x_{r^{*}} = \frac{1}{x_1 x_2 \cdots x_n} \sum_{r \in \text{Res}(w^{*})} x_{r} = x_{w^{*}}.$$
          
         \item The first poset isomorphism $\mathbb{P}_w \longrightarrow \mathbb{L}_{w^{*}}$ is defined by applying the map $G_w \mapsto G_{w}^{T_1 \circ T_2 \circ ... \circ T_n}$ to each node $P$ of $\mathbb{P}_w$, keeping track of the images under the maps $T_i$ of all the edges in $P$. That this map is well-defined is a straightforward induction on the number of tiles $n$ of $G_w$. That weights are preserved is obvious. If the perfect matching $P_2$ covers $P_1$, then $P_2$ can be obtained from $P_1$ by performing an up-twist at some tile $T_i$. Suppose $L_1$ is the image of $P_1$, and $L_2$ is the image of $P_2$. Then one can check that $L_2$ covers $L_1$ by using the twist-parity condition given directly after Figure \ref{fig:twist} and cases on the parity of $i$.  The inverse $\mathbb{L}_{w^{*}} \longrightarrow \mathbb{P}_{w}$ of this map is again application of the snake graph factorization $T_1 \circ \dots \circ T_n$ (the map $\mathbb{L}_{w^{*}} \longrightarrow \mathbb{P}_{w}$ is equal to the composition $T_{n} \circ \dots \circ T_{1}$, and these factors commute). Thus $\mathbb{P}_w \longrightarrow \mathbb{L}_{w^{*}}$ is indeed a poset isomorphism that respects additional node structure. 
    
    The second map $\mathbb{A}_w \longrightarrow \mathbb{B}_{w^{*}}$ takes any perfect matching of angles $A \in \mathbb{A}_w$ to a lattice path of angles in $\mathbb{B}_{w^{*}}$ via the triangle map. That this map is a well-defined poset isomorphism that preserves node structure now follows from the bijection $\mathbb{P}_w \longrightarrow \mathbb{L}_{w^{*}}$ and Propositions \ref{bij_PAT} and \ref{bij_LBS}.

    The third map $\mathbb{T}_{w} \longrightarrow \mathbb{S}_{w^{*}}$ takes any $T \in \mathbb{T}_w$ to an $S$-path on $\Sigma_{w}^{*}$ by first superimposing a blue edge onto each internal diagonal of $\Delta_w$ (this accounts for multiplying each $T$-path weight by $x_1 x_2 \dots x_n$), canceling any blue/red pairs that result, and then applying the triangle map. Again, that this gives a well-defined poset isomorphism preserving additional node structure follows from Propositions \ref{bij_PAT} and \ref{bij_LBS} and the map $\mathbb{P}_w \longrightarrow \mathbb{L}_{w^{*}}$. 
    
    That $\mathbb{P}_w$ is a distributive lattice given in Theorem 5.2 in \cite{musiker2013bases}. Thus, the rest are distributive lattices as well. 
        
      \item Follows directly from (d). 
      
      \item That $\mathcal{I}(C_w) \cong \mathbb{P}_w$ follows from Definition 5.3 and Theorem 5.4 in \cite{musiker2013bases}. That $\mathcal{I}(C_{w}^{*}) \cong \mathbb{L}_w$ follows from the construction on page 18 of \cite{knauer2018lattice}. Namely, the poset $C_{w}^{*}$ can be built from the minimal path $L_{-}$ in $G_w$ by deleting the first and last steps in $L_{-}$, and rotating the result $45^{\circ}$ clockwise. That $\mathbb{P}_w$ and $\mathbb{L}_w$ are dual as distributive lattices follows from the definition. 
    
    \end{enumerate}
    
\end{proof}

\begin{ex}
    We illustrate how applying the isomorphisms from  Theorem \ref{duality_thm} (d) to the three minimal elements in the running example from the previous chapter gives the three respective minimal elements from the running example in this chapter.
\end{ex}

\begin {figure}[h!]
    \centering
    \caption{Transforming the minimal element $P_{-}$ into the minimal element $L_{-}$}
    \label{fig:P_to_L}
    
   \begin{tikzpicture}[x=0.75pt,y=0.75pt,yscale=-1,xscale=1]

\draw [color={rgb, 255:red, 74; green, 144; blue, 226 }  ,draw opacity=1 ][line width=2.25]    (36.64,194.9) -- (86.56,194.9) ;
\draw    (36.64,144.98) -- (36.64,194.9) ;
\draw    (36.64,144.98) -- (86.56,144.98) ;
\draw    (86.56,144.98) -- (86.56,194.9) ;
\draw [color={rgb, 255:red, 74; green, 144; blue, 226 }  ,draw opacity=1 ][line width=2.25]    (86.56,95.06) -- (86.56,144.98) ;
\draw [color={rgb, 255:red, 74; green, 144; blue, 226 }  ,draw opacity=1 ][line width=2.25]    (36.64,95.06) -- (36.64,144.98) ;
\draw    (36.64,95.06) -- (86.56,95.06) ;
\draw    (36.64,45.14) -- (36.64,95.06) ;
\draw    (86.56,45.14) -- (86.56,95.06) ;
\draw [color={rgb, 255:red, 74; green, 144; blue, 226 }  ,draw opacity=1 ][line width=2.25]    (36.64,45.14) -- (86.56,45.14) ;
\draw [color={rgb, 255:red, 74; green, 144; blue, 226 }  ,draw opacity=1 ][line width=2.25]    (119.84,194.9) -- (169.76,194.9) ;
\draw    (119.84,144.98) -- (119.84,194.9) ;
\draw    (119.84,144.98) -- (169.76,144.98) ;
\draw    (169.76,144.98) -- (169.76,194.9) ;
\draw [color={rgb, 255:red, 74; green, 144; blue, 226 }  ,draw opacity=1 ][line width=2.25]    (169.76,194.9) -- (219.68,194.9) ;
\draw [color={rgb, 255:red, 74; green, 144; blue, 226 }  ,draw opacity=1 ][line width=2.25]    (169.76,144.98) -- (219.68,144.98) ;
\draw    (219.68,144.98) -- (219.68,194.9) ;
\draw    (219.68,144.98) -- (269.61,144.98) ;
\draw    (219.68,194.9) -- (269.61,194.9) ;
\draw [color={rgb, 255:red, 74; green, 144; blue, 226 }  ,draw opacity=1 ][line width=2.25]    (269.61,144.98) -- (269.61,194.9) ;
\draw [color={rgb, 255:red, 74; green, 144; blue, 226 }  ,draw opacity=1 ][line width=2.25]    (302.89,194.9) -- (352.81,194.9) ;
\draw    (302.89,144.98) -- (302.89,194.9) ;
\draw    (302.89,144.98) -- (352.81,144.98) ;
\draw    (352.81,144.98) -- (352.81,194.9) ;
\draw [color={rgb, 255:red, 74; green, 144; blue, 226 }  ,draw opacity=1 ][line width=2.25]    (352.81,194.9) -- (402.73,194.9) ;
\draw    (352.81,144.98) -- (402.73,144.98) ;
\draw [color={rgb, 255:red, 74; green, 144; blue, 226 }  ,draw opacity=1 ][line width=2.25]    (402.73,144.98) -- (402.73,194.9) ;
\draw    (352.81,95.06) -- (352.81,144.98) ;
\draw    (402.73,95.06) -- (402.73,144.98) ;
\draw [color={rgb, 255:red, 74; green, 144; blue, 226 }  ,draw opacity=1 ][line width=2.25]    (352.81,95.06) -- (402.73,95.06) ;
\draw [color={rgb, 255:red, 74; green, 144; blue, 226 }  ,draw opacity=1 ][line width=2.25]    (436.01,194.9) -- (485.93,194.9) ;
\draw    (436.01,144.98) -- (436.01,194.9) ;
\draw    (436.01,144.98) -- (485.93,144.98) ;
\draw    (485.93,144.98) -- (485.93,194.9) ;
\draw [color={rgb, 255:red, 74; green, 144; blue, 226 }  ,draw opacity=1 ][line width=2.25]    (485.93,194.9) -- (535.85,194.9) ;
\draw    (485.93,144.98) -- (535.85,144.98) ;
\draw [color={rgb, 255:red, 74; green, 144; blue, 226 }  ,draw opacity=1 ][line width=2.25]    (535.85,144.98) -- (535.85,194.9) ;
\draw    (485.93,95.06) -- (485.93,144.98) ;
\draw [color={rgb, 255:red, 74; green, 144; blue, 226 }  ,draw opacity=1 ][line width=2.25]    (535.85,95.06) -- (535.85,144.98) ;
\draw    (485.93,95.06) -- (535.85,95.06) ;
\draw  [dash pattern={on 0.84pt off 2.51pt}]  (20,128.34) -- (103.2,211.54) ;
\draw  [dash pattern={on 0.84pt off 2.51pt}]  (153.12,128.34) -- (236.32,211.54) ;
\draw  [dash pattern={on 0.84pt off 2.51pt}]  (336.17,78.42) -- (419.37,161.62) ;
\draw    (94.88,223.26) .. controls (118.12,239.73) and (149.75,241.66) .. (176.45,224.35) ;
\draw [shift={(178.08,223.26)}, rotate = 505.66] [color={rgb, 255:red, 0; green, 0; blue, 0 }  ][line width=0.75]    (10.93,-3.29) .. controls (6.95,-1.4) and (3.31,-0.3) .. (0,0) .. controls (3.31,0.3) and (6.95,1.4) .. (10.93,3.29)   ;
\draw    (244.65,223.26) .. controls (267.88,239.73) and (299.51,241.66) .. (326.22,224.35) ;
\draw [shift={(327.85,223.26)}, rotate = 505.66] [color={rgb, 255:red, 0; green, 0; blue, 0 }  ][line width=0.75]    (10.93,-3.29) .. controls (6.95,-1.4) and (3.31,-0.3) .. (0,0) .. controls (3.31,0.3) and (6.95,1.4) .. (10.93,3.29)   ;
\draw    (394.41,223.26) .. controls (417.65,239.73) and (449.28,241.66) .. (475.98,224.35) ;
\draw [shift={(477.61,223.26)}, rotate = 505.66] [color={rgb, 255:red, 0; green, 0; blue, 0 }  ][line width=0.75]    (10.93,-3.29) .. controls (6.95,-1.4) and (3.31,-0.3) .. (0,0) .. controls (3.31,0.3) and (6.95,1.4) .. (10.93,3.29)   ;

\draw (60.77,34.32) node    {$9$};
\draw (93.22,67.6) node    {$8$};
\draw (31.65,69.26) node    {$2$};
\draw (59.94,103.38) node    {$5$};
\draw (94.88,116.69) node    {$3$};
\draw (29.98,117.52) node    {$1$};
\draw (59.94,154.96) node    {$4$};
\draw (59.1,203.22) node    {$6$};
\draw (29.15,166.61) node    {$7$};
\draw (95.71,168.28) node    {$2$};
\draw (147.3,136.66) node    {$2$};
\draw (163.94,167.44) node    {$4$};
\draw (147.3,203.22) node    {$6$};
\draw (213.86,169.94) node    {$5$};
\draw (275.43,169.94) node    {$9$};
\draw (114.02,169.94) node    {$7$};
\draw (192.23,136.66) node    {$3$};
\draw (192.23,203.22) node    {$1$};
\draw (242.15,136.66) node    {$8$};
\draw (247.14,203.22) node    {$2$};
\draw (330.34,136.66) node    {$2$};
\draw (297.06,169.94) node    {$7$};
\draw (330.34,203.22) node    {$6$};
\draw (346.98,169.94) node    {$4$};
\draw (375.27,203.22) node    {$1$};
\draw (375.27,136.66) node    {$5$};
\draw (408.55,169.94) node    {$3$};
\draw (375.27,86.74) node    {$9$};
\draw (408.55,120.02) node    {$8$};
\draw (346.98,120.02) node    {$2$};
\draw (463.47,203.22) node    {$6$};
\draw (513.39,203.22) node    {$1$};
\draw (541.68,169.94) node    {$3$};
\draw (513.39,86.74) node    {$8$};
\draw (541.68,120.02) node    {$9$};
\draw (480.11,120.02) node    {$2$};
\draw (513.39,136.66) node    {$5$};
\draw (480.11,169.94) node    {$4$};
\draw (430.19,169.94) node    {$7$};
\draw (463.47,136.66) node    {$2$};
\draw (135.65,245.65) node    {$T_{1}$};
\draw (287.08,245.65) node    {$T_{2}$};
\draw (435.18,245.65) node    {$T_{3}$};

\end{tikzpicture}
\end {figure}
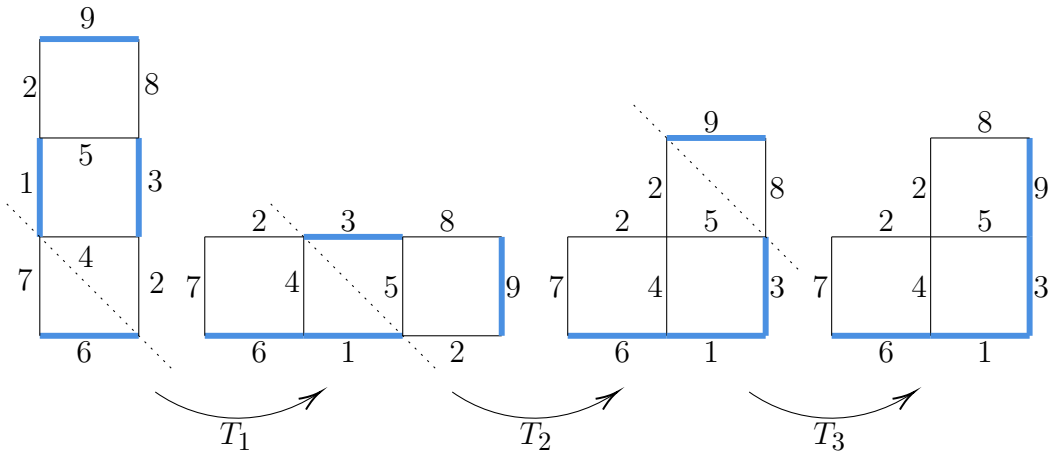

    Figure \ref{fig:P_to_L} shows the minimal matching $P_{-}$ in $\mathbb{P}_{ab}$ being sent to $L_{-}$ in $\mathbb{L}_{bb}$. 

\newpage

\begin {figure}[h!]
    \centering
    \caption{Transforming the minimal element $A_{-}$ into the minimal element $B_{-}$}
    \label{fig:A_to_B}
      \begin{tikzpicture}[x=0.75pt,y=0.75pt,yscale=-1,xscale=1]

\draw    (40,135.38) -- (65.13,110.25) ;
\draw    (65.13,160.5) -- (40,135.38) ;
\draw    (115.38,160.5) -- (65.13,160.5) ;
\draw    (115.38,110.25) -- (65.13,110.25) ;
\draw    (115.38,160.5) -- (140.5,135.38) ;
\draw    (140.5,135.38) -- (115.38,110.25) ;
\draw    (65.13,160.5) -- (65.13,110.25) ;
\draw    (115.38,160.5) -- (115.38,110.25) ;
\draw    (65.13,160.5) -- (115.38,110.25) ;
\draw  [color={rgb, 255:red, 74; green, 144; blue, 226 }  ,draw opacity=1 ][fill={rgb, 255:red, 74; green, 144; blue, 226 }  ,fill opacity=1 ] (57.87,120.87) .. controls (57.87,119.42) and (59.06,118.24) .. (60.51,118.24) .. controls (61.97,118.24) and (63.15,119.42) .. (63.15,120.87) .. controls (63.15,122.33) and (61.97,123.51) .. (60.51,123.51) .. controls (59.06,123.51) and (57.87,122.33) .. (57.87,120.87) -- cycle ;
\draw  [color={rgb, 255:red, 74; green, 144; blue, 226 }  ,draw opacity=1 ][fill={rgb, 255:red, 74; green, 144; blue, 226 }  ,fill opacity=1 ] (101.66,114.77) .. controls (101.66,113.32) and (102.85,112.13) .. (104.3,112.13) .. controls (105.76,112.13) and (106.94,113.32) .. (106.94,114.77) .. controls (106.94,116.23) and (105.76,117.41) .. (104.3,117.41) .. controls (102.85,117.41) and (101.66,116.23) .. (101.66,114.77) -- cycle ;
\draw  [color={rgb, 255:red, 74; green, 144; blue, 226 }  ,draw opacity=1 ][fill={rgb, 255:red, 74; green, 144; blue, 226 }  ,fill opacity=1 ] (73.67,155.69) .. controls (73.67,154.23) and (74.85,153.05) .. (76.31,153.05) .. controls (77.76,153.05) and (78.94,154.23) .. (78.94,155.69) .. controls (78.94,157.15) and (77.76,158.33) .. (76.31,158.33) .. controls (74.85,158.33) and (73.67,157.15) .. (73.67,155.69) -- cycle ;
\draw  [color={rgb, 255:red, 74; green, 144; blue, 226 }  ,draw opacity=1 ][fill={rgb, 255:red, 74; green, 144; blue, 226 }  ,fill opacity=1 ] (117.1,149.23) .. controls (117.1,147.77) and (118.28,146.59) .. (119.74,146.59) .. controls (121.19,146.59) and (122.37,147.77) .. (122.37,149.23) .. controls (122.37,150.69) and (121.19,151.87) .. (119.74,151.87) .. controls (118.28,151.87) and (117.1,150.69) .. (117.1,149.23) -- cycle ;
\draw    (630,96) -- (650,116) ;
\draw    (610,116) -- (630,96) ;
\draw    (650,156) -- (650,116) ;
\draw    (610,156) -- (610,116) ;
\draw    (610,156) -- (630,176) ;
\draw    (630,176) -- (650,156) ;
\draw    (610,156) -- (630,96) ;
\draw    (650,156) -- (630,96) ;
\draw    (630,176) -- (630,96) ;
\draw  [color={rgb, 255:red, 74; green, 144; blue, 226 }  ,draw opacity=1 ][fill={rgb, 255:red, 74; green, 144; blue, 226 }  ,fill opacity=1 ] (617.86,111.64) .. controls (617.86,110.18) and (619.04,109) .. (620.5,109) .. controls (621.96,109) and (623.14,110.18) .. (623.14,111.64) .. controls (623.14,113.1) and (621.96,114.28) .. (620.5,114.28) .. controls (619.04,114.28) and (617.86,113.1) .. (617.86,111.64) -- cycle ;
\draw  [color={rgb, 255:red, 74; green, 144; blue, 226 }  ,draw opacity=1 ][fill={rgb, 255:red, 74; green, 144; blue, 226 }  ,fill opacity=1 ] (622.86,165.14) .. controls (622.86,163.68) and (624.04,162.5) .. (625.5,162.5) .. controls (626.96,162.5) and (628.14,163.68) .. (628.14,165.14) .. controls (628.14,166.6) and (626.96,167.78) .. (625.5,167.78) .. controls (624.04,167.78) and (622.86,166.6) .. (622.86,165.14) -- cycle ;
\draw  [color={rgb, 255:red, 74; green, 144; blue, 226 }  ,draw opacity=1 ][fill={rgb, 255:red, 74; green, 144; blue, 226 }  ,fill opacity=1 ] (631.72,165.5) .. controls (631.72,164.04) and (632.9,162.86) .. (634.36,162.86) .. controls (635.82,162.86) and (637,164.04) .. (637,165.5) .. controls (637,166.96) and (635.82,168.14) .. (634.36,168.14) .. controls (632.9,168.14) and (631.72,166.96) .. (631.72,165.5) -- cycle ;
\draw  [color={rgb, 255:red, 74; green, 144; blue, 226 }  ,draw opacity=1 ][fill={rgb, 255:red, 74; green, 144; blue, 226 }  ,fill opacity=1 ] (636.92,111.1) .. controls (636.92,109.64) and (638.1,108.46) .. (639.56,108.46) .. controls (641.02,108.46) and (642.2,109.64) .. (642.2,111.1) .. controls (642.2,112.56) and (641.02,113.74) .. (639.56,113.74) .. controls (638.1,113.74) and (636.92,112.56) .. (636.92,111.1) -- cycle ;
\draw    (170,134.88) -- (195.13,109.75) ;
\draw    (195.13,160) -- (170,134.88) ;
\draw    (195.13,160) -- (195.13,109.75) ;
\draw  [color={rgb, 255:red, 74; green, 144; blue, 226 }  ,draw opacity=1 ][fill={rgb, 255:red, 74; green, 144; blue, 226 }  ,fill opacity=1 ] (187.87,120.37) .. controls (187.87,118.92) and (189.06,117.74) .. (190.51,117.74) .. controls (191.97,117.74) and (193.15,118.92) .. (193.15,120.37) .. controls (193.15,121.83) and (191.97,123.01) .. (190.51,123.01) .. controls (189.06,123.01) and (187.87,121.83) .. (187.87,120.37) -- cycle ;
\draw    (260.25,109.5) -- (210,109.5) ;
\draw    (210,159.75) -- (210,109.5) ;
\draw    (210,159.75) -- (260.25,109.5) ;
\draw  [color={rgb, 255:red, 74; green, 144; blue, 226 }  ,draw opacity=1 ][fill={rgb, 255:red, 74; green, 144; blue, 226 }  ,fill opacity=1 ] (246.54,114.02) .. controls (246.54,112.57) and (247.72,111.38) .. (249.18,111.38) .. controls (250.63,111.38) and (251.82,112.57) .. (251.82,114.02) .. controls (251.82,115.48) and (250.63,116.66) .. (249.18,116.66) .. controls (247.72,116.66) and (246.54,115.48) .. (246.54,114.02) -- cycle ;
\draw    (280.25,160) -- (230,160) ;
\draw    (280.25,160) -- (280.25,109.75) ;
\draw    (230,160) -- (280.25,109.75) ;
\draw  [color={rgb, 255:red, 74; green, 144; blue, 226 }  ,draw opacity=1 ][fill={rgb, 255:red, 74; green, 144; blue, 226 }  ,fill opacity=1 ] (238.54,155.19) .. controls (238.54,153.73) and (239.72,152.55) .. (241.18,152.55) .. controls (242.64,152.55) and (243.82,153.73) .. (243.82,155.19) .. controls (243.82,156.65) and (242.64,157.83) .. (241.18,157.83) .. controls (239.72,157.83) and (238.54,156.65) .. (238.54,155.19) -- cycle ;
\draw    (294.88,160) -- (320,134.88) ;
\draw    (320,134.88) -- (294.88,109.75) ;
\draw  [color={rgb, 255:red, 74; green, 144; blue, 226 }  ,draw opacity=1 ][fill={rgb, 255:red, 74; green, 144; blue, 226 }  ,fill opacity=1 ] (296.6,148.73) .. controls (296.6,147.27) and (297.78,146.09) .. (299.24,146.09) .. controls (300.69,146.09) and (301.87,147.27) .. (301.87,148.73) .. controls (301.87,150.19) and (300.69,151.37) .. (299.24,151.37) .. controls (297.78,151.37) and (296.6,150.19) .. (296.6,148.73) -- cycle ;
\draw    (294.88,160) -- (294.88,109.75) ;
\draw    (390,135.13) -- (415.13,110) ;
\draw    (415.13,160.25) -- (390,135.13) ;
\draw    (415.13,160.25) -- (415.13,110) ;
\draw  [color={rgb, 255:red, 74; green, 144; blue, 226 }  ,draw opacity=1 ][fill={rgb, 255:red, 74; green, 144; blue, 226 }  ,fill opacity=1 ] (407.87,120.62) .. controls (407.87,119.17) and (409.06,117.99) .. (410.51,117.99) .. controls (411.97,117.99) and (413.15,119.17) .. (413.15,120.62) .. controls (413.15,122.08) and (411.97,123.26) .. (410.51,123.26) .. controls (409.06,123.26) and (407.87,122.08) .. (407.87,120.62) -- cycle ;
\draw    (480,160.25) -- (429.75,160.25) ;
\draw    (480,160.25) -- (430,110) ;
\draw    (429.75,160.25) -- (430,110) ;
\draw  [color={rgb, 255:red, 74; green, 144; blue, 226 }  ,draw opacity=1 ][fill={rgb, 255:red, 74; green, 144; blue, 226 }  ,fill opacity=1 ] (464.72,155.36) .. controls (464.72,153.9) and (465.9,152.72) .. (467.36,152.72) .. controls (468.82,152.72) and (470,153.9) .. (470,155.36) .. controls (470,156.82) and (468.82,158) .. (467.36,158) .. controls (465.9,158) and (464.72,156.82) .. (464.72,155.36) -- cycle ;
\draw    (540,160.25) -- (489.75,160.25) ;
\draw    (540,160.25) -- (540,110) ;
\draw    (489.75,160.25) -- (540,110) ;
\draw  [color={rgb, 255:red, 74; green, 144; blue, 226 }  ,draw opacity=1 ][fill={rgb, 255:red, 74; green, 144; blue, 226 }  ,fill opacity=1 ] (498.29,155.44) .. controls (498.29,153.98) and (499.47,152.8) .. (500.93,152.8) .. controls (502.39,152.8) and (503.57,153.98) .. (503.57,155.44) .. controls (503.57,156.9) and (502.39,158.08) .. (500.93,158.08) .. controls (499.47,158.08) and (498.29,156.9) .. (498.29,155.44) -- cycle ;
\draw    (550,160.25) -- (575.13,135.13) ;
\draw    (575.13,135.13) -- (550,110) ;
\draw  [color={rgb, 255:red, 74; green, 144; blue, 226 }  ,draw opacity=1 ][fill={rgb, 255:red, 74; green, 144; blue, 226 }  ,fill opacity=1 ] (552,121.64) .. controls (552,120.18) and (553.18,119) .. (554.64,119) .. controls (556.1,119) and (557.28,120.18) .. (557.28,121.64) .. controls (557.28,123.1) and (556.1,124.28) .. (554.64,124.28) .. controls (553.18,124.28) and (552,123.1) .. (552,121.64) -- cycle ;
\draw    (550,160.25) -- (550,110) ;
\draw    (330,134) .. controls (331.67,132.33) and (333.33,132.33) .. (335,134) .. controls (336.67,135.67) and (338.33,135.67) .. (340,134) .. controls (341.67,132.33) and (343.33,132.33) .. (345,134) .. controls (346.67,135.67) and (348.33,135.67) .. (350,134) .. controls (351.67,132.33) and (353.33,132.33) .. (355,134) .. controls (356.67,135.67) and (358.33,135.67) .. (360,134) .. controls (361.67,132.33) and (363.33,132.33) .. (365,134) .. controls (366.67,135.67) and (368.33,135.67) .. (370,134) -- (378,134) ;
\draw [shift={(380,134)}, rotate = 180] [color={rgb, 255:red, 0; green, 0; blue, 0 }  ][line width=0.75]    (10.93,-3.29) .. controls (6.95,-1.4) and (3.31,-0.3) .. (0,0) .. controls (3.31,0.3) and (6.95,1.4) .. (10.93,3.29)   ;

\draw (153,130) node    {$=$};
\draw (590,132) node    {$=$};
\end{tikzpicture}
  
\end {figure}
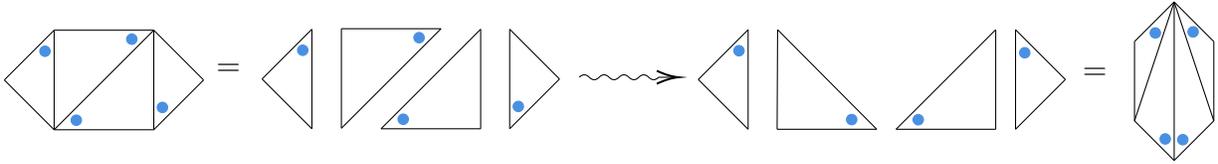

To compute the image of $A_{-} \in \mathbb{A}_{ab}$ under $\mathbb{A}_{ab} \longrightarrow \mathbb{B}_{bb},$ we apply the triangle map. This is shown in Figure \ref{fig:A_to_B}.

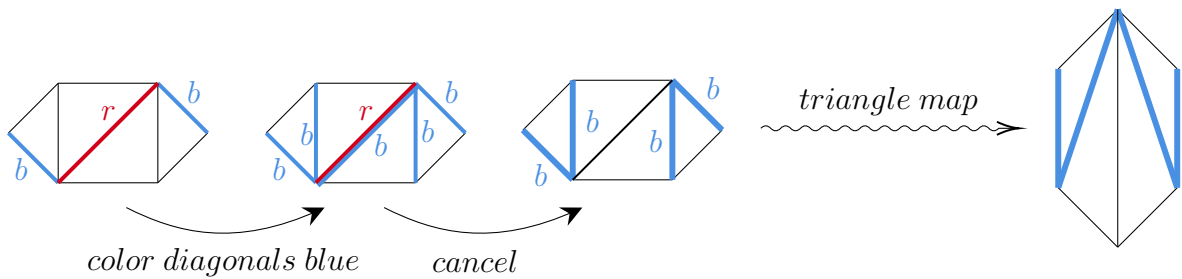
\begin {figure}[h!]
    \centering
    \caption{Transforming the minimal element $T_{-}$ into the minimal element $S_{-}$}
    \label{fig:T_to_S}
     \begin{tikzpicture}[x=0.75pt,y=0.75pt,yscale=-1,xscale=1]

\draw    (40.5,172.38) -- (65.63,147.25) ;
\draw [color={rgb, 255:red, 74; green, 144; blue, 226 }  ,draw opacity=1 ][line width=1.5]    (65.63,197.5) -- (40.5,172.38) ;
\draw    (115.88,197.5) -- (65.63,197.5) ;
\draw    (115.88,147.25) -- (65.63,147.25) ;
\draw    (115.88,197.5) -- (141,172.38) ;
\draw [color={rgb, 255:red, 74; green, 144; blue, 226 }  ,draw opacity=1 ][line width=1.5]    (141,172.38) -- (115.88,147.25) ;
\draw    (65.63,197.5) -- (65.63,147.25) ;
\draw    (115.88,197.5) -- (115.88,147.25) ;
\draw [color={rgb, 255:red, 208; green, 2; blue, 27 }  ,draw opacity=1 ][line width=1.5]    (65.63,197.5) -- (115.88,147.25) ;
\draw    (170.5,172.38) -- (195.63,147.25) ;
\draw [color={rgb, 255:red, 74; green, 144; blue, 226 }  ,draw opacity=1 ][line width=1.5]    (195.63,197.5) -- (170.5,172.38) ;
\draw    (245.88,197.5) -- (195.63,197.5) ;
\draw    (245.88,147.25) -- (195.63,147.25) ;
\draw    (245.88,197.5) -- (271,172.38) ;
\draw [color={rgb, 255:red, 74; green, 144; blue, 226 }  ,draw opacity=1 ][line width=1.5]    (271,172.38) -- (245.88,147.25) ;
\draw [color={rgb, 255:red, 74; green, 144; blue, 226 }  ,draw opacity=1 ][line width=1.5]    (195.63,197.5) -- (195.63,147.25) ;
\draw [color={rgb, 255:red, 74; green, 144; blue, 226 }  ,draw opacity=1 ][line width=1.5]    (245.88,197.5) -- (245.88,147.25) ;
\draw [color={rgb, 255:red, 208; green, 2; blue, 27 }  ,draw opacity=1 ][line width=1.5]    (195.63,197.5) -- (245.88,147.25) ;
\draw [color={rgb, 255:red, 74; green, 144; blue, 226 }  ,draw opacity=1 ][line width=1.5]    (196.75,199.08) -- (247,148.83) ;
\draw    (300,170.88) -- (325.13,145.75) ;
\draw [color={rgb, 255:red, 74; green, 144; blue, 226 }  ,draw opacity=1 ][line width=2.25]    (325.13,196) -- (300,170.88) ;
\draw    (375.38,196) -- (325.13,196) ;
\draw    (375.38,145.75) -- (325.13,145.75) ;
\draw    (375.38,196) -- (400.5,170.88) ;
\draw [color={rgb, 255:red, 74; green, 144; blue, 226 }  ,draw opacity=1 ][line width=2.25]    (400.5,170.88) -- (375.38,145.75) ;
\draw [color={rgb, 255:red, 74; green, 144; blue, 226 }  ,draw opacity=1 ][line width=2.25]    (325.13,196) -- (325.13,145.75) ;
\draw [color={rgb, 255:red, 74; green, 144; blue, 226 }  ,draw opacity=1 ][line width=2.25]    (375.38,196) -- (375.38,145.75) ;
\draw [color={rgb, 255:red, 0; green, 0; blue, 0 }  ,draw opacity=1 ][line width=0.75]    (325.13,196) -- (375.38,145.75) ;
\draw    (600,110) -- (630,140) ;
\draw    (570,140) -- (600,110) ;
\draw [color={rgb, 255:red, 74; green, 144; blue, 226 }  ,draw opacity=1 ][line width=2.25]    (630,200) -- (630,140) ;
\draw [color={rgb, 255:red, 74; green, 144; blue, 226 }  ,draw opacity=1 ][line width=2.25]    (570,200) -- (570,140) ;
\draw    (570,200) -- (600,230) ;
\draw    (600,230) -- (630,200) ;
\draw [color={rgb, 255:red, 74; green, 144; blue, 226 }  ,draw opacity=1 ][line width=2.25]    (570,200) -- (600,110) ;
\draw [color={rgb, 255:red, 74; green, 144; blue, 226 }  ,draw opacity=1 ][line width=2.25]    (630,200) -- (600,110) ;
\draw    (600,230) -- (600,110) ;
\draw    (420,170) .. controls (421.67,168.33) and (423.33,168.33) .. (425,170) .. controls (426.67,171.67) and (428.33,171.67) .. (430,170) .. controls (431.67,168.33) and (433.33,168.33) .. (435,170) .. controls (436.67,171.67) and (438.33,171.67) .. (440,170) .. controls (441.67,168.33) and (443.33,168.33) .. (445,170) .. controls (446.67,171.67) and (448.33,171.67) .. (450,170) .. controls (451.67,168.33) and (453.33,168.33) .. (455,170) .. controls (456.67,171.67) and (458.33,171.67) .. (460,170) .. controls (461.67,168.33) and (463.33,168.33) .. (465,170) .. controls (466.67,171.67) and (468.33,171.67) .. (470,170) .. controls (471.67,168.33) and (473.33,168.33) .. (475,170) .. controls (476.67,171.67) and (478.33,171.67) .. (480,170) .. controls (481.67,168.33) and (483.33,168.33) .. (485,170) .. controls (486.67,171.67) and (488.33,171.67) .. (490,170) .. controls (491.67,168.33) and (493.33,168.33) .. (495,170) .. controls (496.67,171.67) and (498.33,171.67) .. (500,170) .. controls (501.67,168.33) and (503.33,168.33) .. (505,170) .. controls (506.67,171.67) and (508.33,171.67) .. (510,170) .. controls (511.67,168.33) and (513.33,168.33) .. (515,170) .. controls (516.67,171.67) and (518.33,171.67) .. (520,170) .. controls (521.67,168.33) and (523.33,168.33) .. (525,170) .. controls (526.67,171.67) and (528.33,171.67) .. (530,170) .. controls (531.67,168.33) and (533.33,168.33) .. (535,170) .. controls (536.67,171.67) and (538.33,171.67) .. (540,170) -- (548,170) ;
\draw [shift={(550,170)}, rotate = 180] [color={rgb, 255:red, 0; green, 0; blue, 0 }  ][line width=0.75]    (10.93,-3.29) .. controls (6.95,-1.4) and (3.31,-0.3) .. (0,0) .. controls (3.31,0.3) and (6.95,1.4) .. (10.93,3.29)   ;
\draw    (100,210) .. controls (125.26,222.76) and (154.58,230.69) .. (197.36,211.23) ;
\draw [shift={(200,210)}, rotate = 514.5799999999999] [fill={rgb, 255:red, 0; green, 0; blue, 0 }  ][line width=0.08]  [draw opacity=0] (10.72,-5.15) -- (0,0) -- (10.72,5.15) -- (7.12,0) -- cycle    ;
\draw    (230,210) .. controls (255.26,222.76) and (284.58,230.69) .. (327.36,211.23) ;
\draw [shift={(330,210)}, rotate = 514.5799999999999] [fill={rgb, 255:red, 0; green, 0; blue, 0 }  ][line width=0.08]  [draw opacity=0] (10.72,-5.15) -- (0,0) -- (10.72,5.15) -- (7.12,0) -- cycle    ;

\draw (47,190.5) node    {$\textcolor[rgb]{0.29,0.56,0.89}{b}$};
\draw (134,151.5) node    {$\textcolor[rgb]{0.29,0.56,0.89}{b}$};
\draw (91,162.5) node    {$\textcolor[rgb]{0.82,0.01,0.11}{r}$};
\draw (177,190.5) node    {$\textcolor[rgb]{0.29,0.56,0.89}{b}$};
\draw (264,151.5) node    {$\textcolor[rgb]{0.29,0.56,0.89}{b}$};
\draw (221,162.5) node    {$\textcolor[rgb]{0.82,0.01,0.11}{r}$};
\draw (252,171.5) node    {$\textcolor[rgb]{0.29,0.56,0.89}{b}$};
\draw (228.33,178.5) node    {$\textcolor[rgb]{0.29,0.56,0.89}{b}$};
\draw (191,172.83) node    {$\textcolor[rgb]{0.29,0.56,0.89}{b}$};
\draw (309,193.5) node    {$\textcolor[rgb]{0.29,0.56,0.89}{b}$};
\draw (335,166.5) node    {$\textcolor[rgb]{0.29,0.56,0.89}{b}$};
\draw (367,176.5) node    {$\textcolor[rgb]{0.29,0.56,0.89}{b}$};
\draw (396,149.5) node    {$\textcolor[rgb]{0.29,0.56,0.89}{b}$};
\draw (79,229) node [anchor=north west][inner sep=0.75pt]    {$color\ diagonals\ blue$};
\draw (253,230) node [anchor=north west][inner sep=0.75pt]    {$cancel$};
\draw (438,148) node [anchor=north west][inner sep=0.75pt]    {$triangle\ map$};

\end{tikzpicture}

\end {figure}

The isomorphism $\mathbb{T}_{ab} \longrightarrow \mathbb{S}_{bb}$ is obtained by coloring all internal diagonals blue, canceling blue-red pairs if necessary, and then applying the triangle map. See Figure \ref{fig:T_to_S}.

\chapter {Expansion Posets as Intervals in Young's Lattice}

Loosely speaking, a graded poset is one whose elements can be arranged into ``horizontal ranks''. Each graded poset has an associated rank-generating function, which is a polynomial in one variable whose $j^{th}$ nonnegative integer coefficient records the number of elements sitting at rank $j$. The first objective of this chapter is to give some preliminaries on graded posets and their rank functions, and to note that each expansion poset we have studies thus far is graded.

Our second objective is to put a groupoid structure on the set of all snake graphs, and refine each orbit of this groupoid to a graded poset. The covering relation in each such poset resembles a flip of a lattice path inside a snake graph. We observe that each poset of snake graphs from our construction is isomorphic to one of the well-known lattices $L(m,n)$ whose rank generating functions are the classical $q$-binomial coefficients. For more on the posets $L(m,n)$ and their (symmetric and unimodal) rank generating functions (and much more on unimodality in general, and related concepts), see \cite{stanley1989log}, \cite{brenti178log}, and \cite{branden2015unimodality}.

 Our final objective in this chapter is to show how each poset $L(m,n)$ has a covering by intervals, each of which is isomorphic to one of the lattice path expansion posets considered above. This covering is such that two lattice path expansions embed into the same lattice $L(m,n)$ if and only if their underlying snake graphs are related by a morphism as mentioned above.  Finally, since each $L(m,n)$ is itself an interval in Young's lattice, this shows that each expansion poset is isomorphic to an interval in the latter.

\section{Graded Expansion Posets}

\begin{defn}
       Let $D$ be a finite poset. A \textit{chain in $D$} is a totally ordered subset of $D$. A \textit{maximal chain in $D$} is a chain that is not a proper subset of any other chain in $D$. The \textit{length} of a chain with $k$ elements is $k-1$. We say that $D$ is a \textit{graded poset} if all maximal chains in $D$ have the same finite length. If $D$ is graded then there exists a \textit{rank function} $\rho: D \longrightarrow \mathbb{N} = \{ 0,1,2,\dots \}$ that satisfies the following. 

\begin{enumerate}[(1)]
    \item The minimal elements of $D$ map to $0.$
    
    \item For every $x,y \in D,$ $x < y$ implies $\rho(x) < \rho(y).$
    
    \item If $x < y$ and there does not exist $z \in D$ such that $x < z < y$ (i.e., if $y$ \textit{strictly covers} $x$), then $\rho(y) = \rho(x) + 1.$
\end{enumerate}

 We say the element $x \in D$ has \textit{rank} $i$ if $\rho (x) = i$. The \textit{rank} of the finite graded poset $D$ is equal to the length of any maximal chain.

\end{defn}

It is an easy consequence of Birkhoff's Theorem \ref{birkhoff} above that every finite distributive lattice $D$ is graded. Indeed, for input the order ideal $I \in \mathcal{I}(C) \cong D$ the rank is $\rho(I) = \abs{I},$ the cardinality of $I$. Thus, $D_w \cong \mathcal{I}(C_w )$ is graded for each word $w$.

\begin{prop} The rank of any lattice path $L \in \mathbb{L}_w$ is the number of tiles enclosed by the symmetric difference $L \ominus L_{-}$ of $L$ with the minimal lattice path $L_{-}$ from Definition \ref{L_min_def}.
\end{prop}

\begin{proof}
    Follows from Theorem 5.1 in \cite{musiker2010cluster} and Proposition \ref{shape_duality}.
\end{proof}

\begin{defn}
       Suppose $w$ has length $n-1$ and consider the graded distributive lattice $\mathbb{L}_w$. The \textit{rank-generating function $\mathbb{L}_w (q)$ of the lattice $\mathbb{L}_w$} is the polynomial in $q$ of degree $n$ defined by $\mathbb{L}_w (q) = \sum_{i =0}^{n} r_i q^{i},$ where $r_i$ equals the number of lattice paths of rank $i$ in $\mathbb{L}_w$.
\end{defn}

\begin{defn} \label{unimodal_def}
       Let $\rho$ be the rank-generating function of a graded poset $D$ of rank $n$. The rank-generating function $\rho (q) = \sum_{i=0}^{n} r_i q^{i}$ is \textit{unimodal} if there exists some $m$ such that $r_0 \leq r_1 \leq \dots r_{m-1} \leq r_{m} \geq r_{m+1} \geq \dots \geq r_{n}.$ We say $\rho (q)$ is \textit{symmetric} if $r_{n - i} = r_{i}$ for each $i$. We call $D$ a \textit{rank-unimodal poset} (or just \textit{unimodal}) if its rank function is unimodal, and call $D$ a \textit{rank-symmetric poset} (or just \textit{symmetric}) if its rank function is symmetric.
\end{defn}

For instance, Fibonacci cubes are unimodal, and those of even order are symmetric.  \cite{munarini2002rank}.

\begin {figure}[h!]
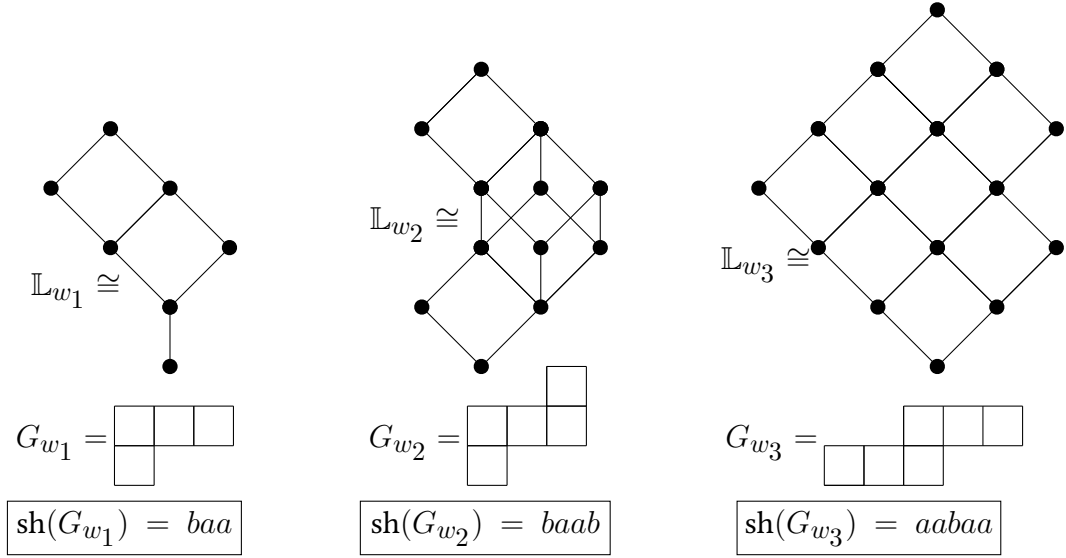

    \centering
    \caption{Posets, shapes, and rank functions for three snake graphs}
    \label{fig:3_snake_graphs}

  
\end {figure}

\begin{ex} \label{three_snakes_ex}
    Figure \ref{fig:3_snake_graphs} shows three snake graphs $G_{w_1} , G_{w_2} , $ and $G_{w_3},$ along with their respective shapes and lattices. Below these figures, we indicate the respective rank generating functions. Note the rank-generating functions $\mathbb{L}_{w_i} (q)$ are unimodal for $i = 1,2,3$ and symmetric for $i = 2$ or $3$. 
\end{ex}

\section{Posets of Snake Graphs and the $q$-binomial Coefficients}

\begin{defn} \label{definition_groupoid}     
    A \textit{groupoid} is a category such that every morphism is invertible.
\end{defn}

Let $\mathcal{L}^{n}$ be the set of snake graphs with $n \geq 1$ tiles. Let $\mathcal{L} = \cup \mathcal{L}^{n}.$ Consider $\mathcal{L}$ a groupoid by saying its elements are the objects, and declaring that there is a morphism between two snake graphs if they have the same number of tiles, and are related by the following local move: 

\begin {figure}[h!]
    \centering
    \caption{Local picture for a morphism between snake graphs}
    \label{fig:SG_morphism}
    \begin{tikzpicture}[x=0.75pt,y=0.75pt,yscale=-1,xscale=1]

\draw    (130,96.12) -- (130,142.99) ;
\draw    (176.87,96.12) -- (176.87,142.99) ;
\draw    (130,96.12) -- (176.87,96.12) ;
\draw    (130,142.99) -- (176.87,142.99) ;
\draw    (176.87,96.12) -- (176.87,142.99) ;
\draw    (176.87,96.12) -- (223.74,96.12) ;
\draw    (176.87,142.99) -- (223.74,142.99) ;
\draw    (223.74,96.12) -- (223.74,142.99) ;
\draw    (176.87,50) -- (176.87,96.87) ;
\draw    (176.87,50) -- (176.87,96.87) ;
\draw    (176.87,50) -- (223.74,50) ;
\draw    (223.74,50) -- (223.74,96.87) ;
\draw    (296.26,96.12) -- (296.26,142.99) ;
\draw    (343.13,96.12) -- (343.13,142.99) ;
\draw    (296.26,96.12) -- (343.13,96.12) ;
\draw    (296.26,142.99) -- (343.13,142.99) ;
\draw    (343.13,96.12) -- (343.13,142.99) ;
\draw    (343.13,96.12) -- (390,96.12) ;
\draw    (343.13,50) -- (343.13,96.87) ;
\draw    (343.13,50) -- (343.13,96.87) ;
\draw    (343.13,50) -- (390,50) ;
\draw    (390,50) -- (390,96.87) ;
\draw    (296.26,50) -- (296.26,96.87) ;
\draw    (296.26,50) -- (296.26,96.87) ;
\draw    (296.26,50) -- (343.13,50) ;
\draw    (230.38,99.19) -- (289.62,99.19) ;
\draw [shift={(291.62,99.19)}, rotate = 180] [color={rgb, 255:red, 0; green, 0; blue, 0 }  ][line width=0.75]    (10.93,-3.29) .. controls (6.95,-1.4) and (3.31,-0.3) .. (0,0) .. controls (3.31,0.3) and (6.95,1.4) .. (10.93,3.29)   ;
\draw [shift={(228.38,99.19)}, rotate = 0] [color={rgb, 255:red, 0; green, 0; blue, 0 }  ][line width=0.75]    (10.93,-3.29) .. controls (6.95,-1.4) and (3.31,-0.3) .. (0,0) .. controls (3.31,0.3) and (6.95,1.4) .. (10.93,3.29)   ;

\end{tikzpicture}
   
\end {figure}
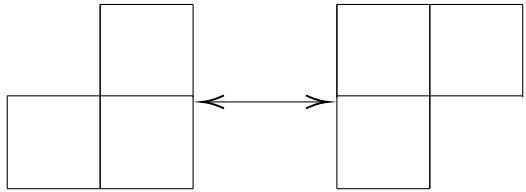

In particular, each straight snake graph (and each snake graph with less than three tiles) is the sole member of its orbit.

Each $\mathcal{L}^n$ is an induced subgroupoid of $\mathcal{L}$. For $n \geq 3$, we parameterize any orbit of $\mathcal{L}^{n}$ by the snake graph $G_{-}^{n,j}$ it contains which is of shape $a^{k_1-1} b^{k_2-1} = a^{k_1 - 1} b^{j}$, where $n = k_1 + k_2 - 1$. Use the notation $O_{j}^{n}$ for the $n-2$ orbits of the subgroupoid $\mathcal{L}^{n}$.

Refine each orbit $O_{j}^{n}$ of  $\mathcal{L}^{n} \hookrightarrow \mathcal{L}$ to a poset $\mathbb{O}_{j}^{n}$ by declaring that the snake graph $G^{n,j}_{-}$ is the minimal element, and by saying that performing the swap $ab \mapsto ba$ corresponds to going up in the poset.

\begin {figure}[h!]
    \centering
    \caption{The poset $\mathbb{O}_{2}^{5}$}
    \label{fig:poset_of_snake_graphs}
    \begin{tikzpicture}[x=0.75pt,y=0.75pt,yscale=-1,xscale=1]

\draw    (216.89,556.38) -- (229.97,556.38) ;
\draw    (216.89,543.3) -- (229.97,543.3) ;
\draw    (229.97,556.38) -- (229.97,543.3) ;
\draw    (216.89,556.38) -- (216.89,543.3) ;
\draw    (229.97,556.38) -- (243.06,556.38) ;
\draw    (229.97,543.3) -- (243.06,543.3) ;
\draw    (243.06,556.38) -- (243.06,543.3) ;
\draw    (229.97,556.38) -- (229.97,543.3) ;
\draw    (243.06,556.38) -- (256.14,556.38) ;
\draw    (243.06,543.3) -- (256.14,543.3) ;
\draw    (256.14,556.38) -- (256.14,543.3) ;
\draw    (243.06,556.38) -- (243.06,543.3) ;
\draw    (243.06,543.3) -- (256.14,543.3) ;
\draw    (243.06,530.21) -- (256.14,530.21) ;
\draw    (256.14,543.3) -- (256.14,530.21) ;
\draw    (243.06,543.3) -- (243.06,530.21) ;
\draw    (243.06,530.21) -- (256.14,530.21) ;
\draw    (243.06,517.13) -- (256.14,517.13) ;
\draw    (256.14,530.21) -- (256.14,517.13) ;
\draw    (243.06,530.21) -- (243.06,517.13) ;
\draw    (216.02,438.62) -- (229.1,438.62) ;
\draw    (216.02,425.54) -- (229.1,425.54) ;
\draw    (229.1,438.62) -- (229.1,425.54) ;
\draw    (216.02,438.62) -- (216.02,425.54) ;
\draw    (229.1,438.62) -- (242.18,438.62) ;
\draw    (229.1,425.54) -- (242.18,425.54) ;
\draw    (242.18,438.62) -- (242.18,425.54) ;
\draw    (229.1,438.62) -- (229.1,425.54) ;
\draw    (242.18,425.54) -- (255.27,425.54) ;
\draw    (242.18,438.62) -- (242.18,425.54) ;
\draw    (242.18,425.54) -- (255.27,425.54) ;
\draw    (242.18,412.45) -- (255.27,412.45) ;
\draw    (255.27,425.54) -- (255.27,412.45) ;
\draw    (242.18,425.54) -- (242.18,412.45) ;
\draw    (242.18,412.45) -- (255.27,412.45) ;
\draw    (242.18,399.37) -- (255.27,399.37) ;
\draw    (255.27,412.45) -- (255.27,399.37) ;
\draw    (242.18,412.45) -- (242.18,399.37) ;
\draw    (229.1,412.45) -- (242.18,412.45) ;
\draw    (229.1,425.54) -- (229.1,412.45) ;
\draw    (112.21,333.94) -- (125.3,333.94) ;
\draw    (112.21,320.86) -- (125.3,320.86) ;
\draw    (125.3,333.94) -- (125.3,320.86) ;
\draw    (112.21,333.94) -- (112.21,320.86) ;
\draw    (112.21,307.78) -- (125.3,307.78) ;
\draw    (125.3,320.86) -- (138.38,320.86) ;
\draw    (125.3,333.94) -- (125.3,320.86) ;
\draw    (138.38,320.86) -- (151.47,320.86) ;
\draw    (112.21,320.86) -- (112.21,307.78) ;
\draw    (138.38,320.86) -- (151.47,320.86) ;
\draw    (138.38,307.78) -- (151.47,307.78) ;
\draw    (151.47,320.86) -- (151.47,307.78) ;
\draw    (138.38,320.86) -- (138.38,307.78) ;
\draw    (138.38,307.78) -- (151.47,307.78) ;
\draw    (138.38,294.69) -- (151.47,294.69) ;
\draw    (151.47,307.78) -- (151.47,294.69) ;
\draw    (138.38,307.78) -- (138.38,294.69) ;
\draw    (125.3,307.78) -- (138.38,307.78) ;
\draw    (125.3,320.86) -- (125.3,307.78) ;
\draw    (321.56,335.69) -- (334.65,335.69) ;
\draw    (321.56,322.6) -- (334.65,322.6) ;
\draw    (334.65,335.69) -- (334.65,322.6) ;
\draw    (321.56,335.69) -- (321.56,322.6) ;
\draw    (334.65,335.69) -- (347.73,335.69) ;
\draw    (334.65,322.6) -- (347.73,322.6) ;
\draw    (347.73,335.69) -- (347.73,322.6) ;
\draw    (334.65,335.69) -- (334.65,322.6) ;
\draw    (347.73,335.69) -- (347.73,322.6) ;
\draw    (334.65,296.44) -- (347.73,296.44) ;
\draw    (347.73,309.52) -- (360.82,309.52) ;
\draw    (334.65,309.52) -- (334.65,296.44) ;
\draw    (347.73,322.6) -- (347.73,309.52) ;
\draw    (347.73,309.52) -- (360.82,309.52) ;
\draw    (347.73,296.44) -- (360.82,296.44) ;
\draw    (360.82,309.52) -- (360.82,296.44) ;
\draw    (347.73,309.52) -- (347.73,296.44) ;
\draw    (334.65,309.52) -- (347.73,309.52) ;
\draw    (334.65,322.6) -- (334.65,309.52) ;
\draw    (217.76,244.1) -- (230.84,244.1) ;
\draw    (217.76,231.01) -- (230.84,231.01) ;
\draw    (230.84,244.1) -- (230.84,231.01) ;
\draw    (217.76,244.1) -- (217.76,231.01) ;
\draw    (217.76,217.93) -- (230.84,217.93) ;
\draw    (230.84,231.01) -- (243.93,231.01) ;
\draw    (230.84,244.1) -- (230.84,231.01) ;
\draw    (217.76,231.01) -- (217.76,217.93) ;
\draw    (230.84,204.84) -- (243.93,204.84) ;
\draw    (243.93,217.93) -- (257.01,217.93) ;
\draw    (230.84,217.93) -- (230.84,204.84) ;
\draw    (243.93,231.01) -- (243.93,217.93) ;
\draw    (243.93,217.93) -- (257.01,217.93) ;
\draw    (243.93,204.84) -- (257.01,204.84) ;
\draw    (257.01,217.93) -- (257.01,204.84) ;
\draw    (243.93,217.93) -- (243.93,204.84) ;
\draw    (230.84,217.93) -- (243.93,217.93) ;
\draw    (230.84,231.01) -- (230.84,217.93) ;
\draw    (218.63,127.21) -- (231.72,127.21) ;
\draw    (218.63,114.13) -- (231.72,114.13) ;
\draw    (231.72,127.21) -- (231.72,114.13) ;
\draw    (218.63,127.21) -- (218.63,114.13) ;
\draw    (218.63,101.04) -- (231.72,101.04) ;
\draw    (218.63,87.96) -- (231.72,87.96) ;
\draw    (231.72,127.21) -- (231.72,114.13) ;
\draw    (218.63,114.13) -- (218.63,101.04) ;
\draw    (231.72,87.96) -- (244.8,87.96) ;
\draw    (244.8,101.04) -- (257.89,101.04) ;
\draw    (231.72,101.04) -- (231.72,87.96) ;
\draw    (218.63,101.04) -- (218.63,87.96) ;
\draw    (244.8,101.04) -- (257.89,101.04) ;
\draw    (244.8,87.96) -- (257.89,87.96) ;
\draw    (257.89,101.04) -- (257.89,87.96) ;
\draw    (244.8,101.04) -- (244.8,87.96) ;
\draw    (231.72,101.04) -- (244.8,101.04) ;
\draw    (231.72,114.13) -- (231.72,101.04) ;
\draw  [dash pattern={on 4.5pt off 4.5pt}] (204.68,538.5) .. controls (204.68,520.43) and (219.32,505.79) .. (237.39,505.79) .. controls (255.45,505.79) and (270.1,520.43) .. (270.1,538.5) .. controls (270.1,556.56) and (255.45,571.21) .. (237.39,571.21) .. controls (219.32,571.21) and (204.68,556.56) .. (204.68,538.5) -- cycle ;
\draw  [dash pattern={on 4.5pt off 4.5pt}] (204.68,420.74) .. controls (204.68,402.67) and (219.32,388.03) .. (237.39,388.03) .. controls (255.45,388.03) and (270.1,402.67) .. (270.1,420.74) .. controls (270.1,438.8) and (255.45,453.45) .. (237.39,453.45) .. controls (219.32,453.45) and (204.68,438.8) .. (204.68,420.74) -- cycle ;
\draw  [dash pattern={on 4.5pt off 4.5pt}] (100,316.06) .. controls (100,298) and (114.65,283.35) .. (132.71,283.35) .. controls (150.78,283.35) and (165.42,298) .. (165.42,316.06) .. controls (165.42,334.13) and (150.78,348.77) .. (132.71,348.77) .. controls (114.65,348.77) and (100,334.13) .. (100,316.06) -- cycle ;
\draw  [dash pattern={on 4.5pt off 4.5pt}] (309.35,316.06) .. controls (309.35,298) and (324,283.35) .. (342.06,283.35) .. controls (360.13,283.35) and (374.77,298) .. (374.77,316.06) .. controls (374.77,334.13) and (360.13,348.77) .. (342.06,348.77) .. controls (324,348.77) and (309.35,334.13) .. (309.35,316.06) -- cycle ;
\draw  [dash pattern={on 4.5pt off 4.5pt}] (204.68,224.47) .. controls (204.68,206.41) and (219.32,191.76) .. (237.39,191.76) .. controls (255.45,191.76) and (270.1,206.41) .. (270.1,224.47) .. controls (270.1,242.54) and (255.45,257.18) .. (237.39,257.18) .. controls (219.32,257.18) and (204.68,242.54) .. (204.68,224.47) -- cycle ;
\draw  [dash pattern={on 4.5pt off 4.5pt}] (204.68,106.71) .. controls (204.68,88.65) and (219.32,74) .. (237.39,74) .. controls (255.45,74) and (270.1,88.65) .. (270.1,106.71) .. controls (270.1,124.78) and (255.45,139.42) .. (237.39,139.42) .. controls (219.32,139.42) and (204.68,124.78) .. (204.68,106.71) -- cycle ;
\draw    (237.39,499.24) -- (237.39,458.68) ;
\draw    (237.39,186.53) -- (237.39,145.96) ;
\draw    (209.91,394.57) -- (160.19,343.54) ;
\draw    (313.28,288.58) -- (270.1,244.1) ;
\draw    (314.58,343.54) -- (263.56,393.26) ;
\draw    (204.68,244.1) -- (158.88,288.58) ;
\end{tikzpicture}
  
\end {figure}
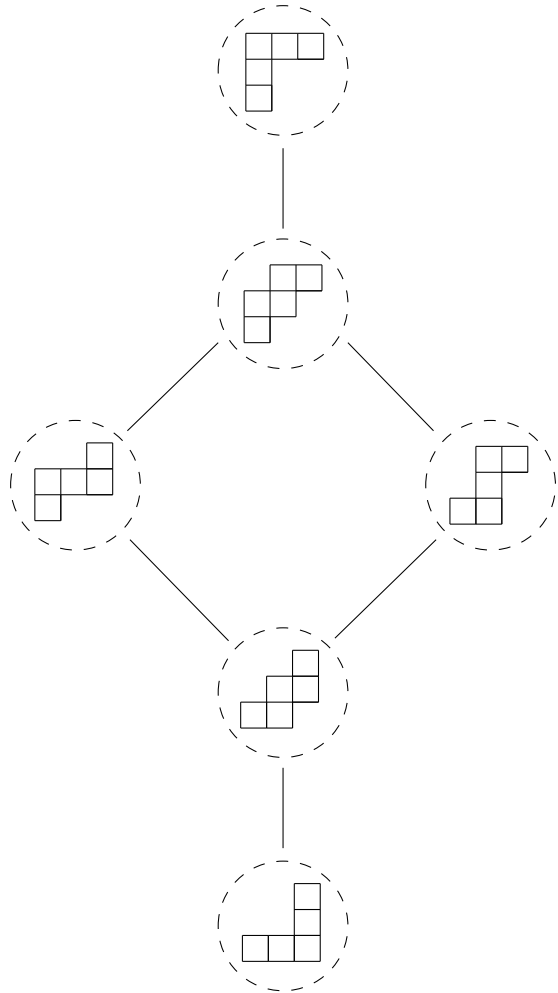

\begin{ex}
    Figure \ref{fig:poset_of_snake_graphs} show the poset $\mathbb{O}_{2}^{5}$ from the subgroupoid $\mathcal{L}^{5} \hookrightarrow \mathcal{L}$. 
\end{ex}

The rank functions of the posets $\mathbb{O}_{j}^{n}$ are well known.

\begin{defn}
       The \textit{$q$-binomial coefficients} are defined by $${n \brack k}_q = \frac{[n]_q!}{[k]_q ! [n-k]_q !},$$
       where $[k]_q! = (1+q)(1+q+q^2)\dots(1+q+ \dots + q^{k-1}).$ Each $q$-binomial coefficient is a rational function in the indeterminate $q$, and is in fact a polynomial function with positive coefficients. Note that taking the limit $q \rightarrow 1$ recovers the standard binomial coefficients.
\end{defn}

\begin{defn}
    The \textit{product} of two posets $(\mathcal{P} , \leq_{\mathcal{P}} )$ and $(\mathcal{Q} , \leq_{\mathcal{Q}})$ is the poset whose underlying set is equal to the Cartesian product $\mathcal{P} \times \mathcal{Q}$ with covering relations given by $$(p_1 , q_1) \leq_{\mathcal{P} \times \mathcal{Q}} (p_2 , q_2) \iff p_1 \leq_{\mathcal{P}} p_2 \text{ and } q_1 \leq_{\mathcal{Q}} q_2.$$ 
\end{defn}

Let $\mathbf{k}$ and $\mathbf{l}$ be the chain posets with $k$ and $l$ vertices, respectively. Consider the product poset $\mathbf{k \times l}.$ Then the rank generating function of the poset of order ideals of $\mathbf{l \times k}$ is equal to the $q$-binomial coefficient ${k+l \brack l}$ (see \cite{stanley1989log}). The $q$-binomial coefficients are also the rank generating functions for lattice paths in $\mathbb{Z}^{2}$ from the origin to the point $(k,l)$ with positive coordinates, or equivalently the Young diagrams that fit inside a $k$ by $l$ rectangular grid.

It is not hard to see that any $q$-binomial coefficient is symmetric (for instance see \cite{{stanley1989log}}). However, it is a nontrivial fact that these coefficients are unimodal. This was first proved by Sylvester in 1878 (see \cite{stanley1989log} and \cite{branden2015unimodality}). The first combinatorial proof of unimodality was given over one hundred years later by O'Hara in \cite{o1990unimodality}.

\begin{ex}
    For instance, the rank generating function of $\mathbb{O}_{3}^{7} \cong \mathcal{I}(\mathbf{2 \times 2})$ is equal to $$1 + q + 2q^2 + 3q^3 + 3q^4 + 3q^5 + 3q^6 + 2q^7 + q^8 + q^9.$$
    This is visualized in the next figure.
\end{ex}

\begin {figure}[h!]
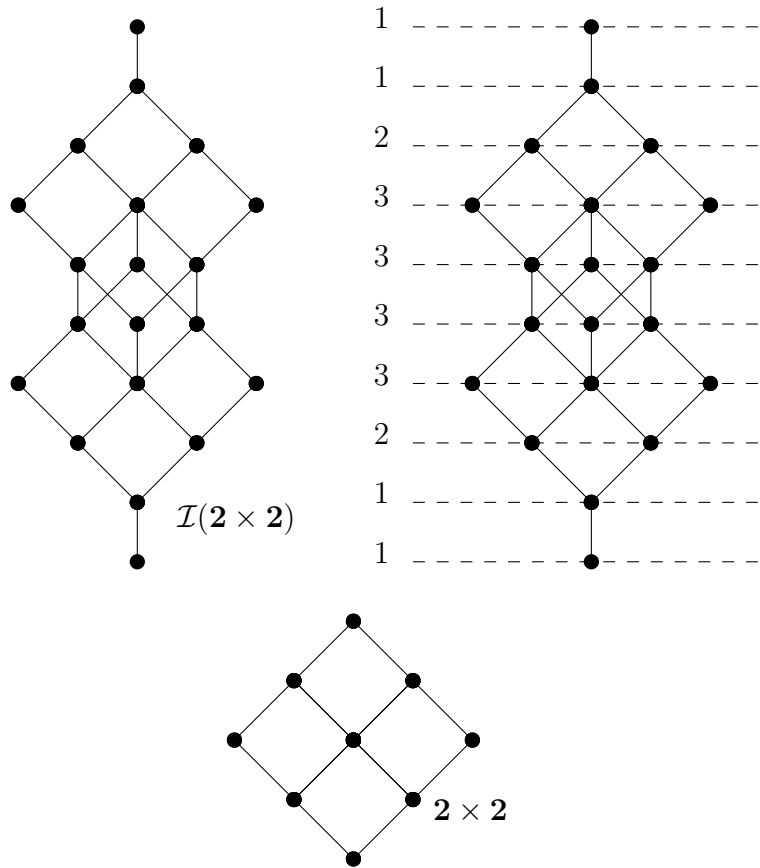

    \centering
    \caption{The poset $\mathbf{2 \times 2}$ (bottom), its lattice of order ideals $\mathcal{I} (\mathbf{2 \times 2})$ (top left), and the rank function of $\mathbb{O}_{3}^{7} \cong \mathcal{I}(\mathbf{2 \times 2})$ (top right)} .
    \label{fig:poset_lattice_order_ideals}

  
\end {figure}

\section {The Embeddings $\mathbb{L}_w \hookrightarrow \mathbb{O}_{j}^{n}$}

Recall that a \textit{partition} of a positive number $m$ is a weakly decreasing sequence $\lambda = (\lambda_1 , \lambda_2 , \dots , \lambda_{l} )$ such that $m = \lambda_1 + \lambda_2 + \dots + \lambda_l.$ For example, $\lambda = (4,3)$ is one partition of $7$. A \textit{Young diagram} is a way to visualize a partition as a left-justified collection of rows of boxes. For instance, to form the Young diagram associated to the partition $(4,3)$ of $7$, we draw an array of $7$ boxes, consisting of a row of $4$ boxes followed by a row of $3$ boxes. 

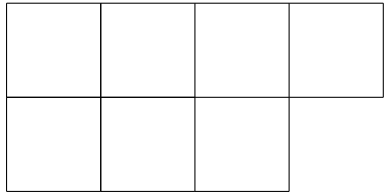
\begin {figure}[h!]
    \centering
    \caption{The Young diagram associated to the partition $(4,3)$ of $7$}
    \label{fig:Young_diagram}
    \begin{tikzpicture}[x=0.75pt,y=0.75pt,yscale=-1,xscale=1]

\draw    (570,317.5) -- (617.5,317.5) ;
\draw    (570,365) -- (617.5,365) ;
\draw    (617.5,317.5) -- (617.5,365) ;
\draw    (570,317.5) -- (570,365) ;
\draw    (617.5,317.5) -- (665,317.5) ;
\draw    (617.5,365) -- (665,365) ;
\draw    (665,317.5) -- (665,365) ;
\draw    (617.5,317.5) -- (617.5,365) ;
\draw    (570,270) -- (617.5,270) ;
\draw    (570,317.5) -- (617.5,317.5) ;
\draw    (617.5,270) -- (617.5,317.5) ;
\draw    (570,270) -- (570,317.5) ;
\draw    (617.5,270) -- (665,270) ;
\draw    (617.5,317.5) -- (665,317.5) ;
\draw    (665,270) -- (665,317.5) ;
\draw    (617.5,270) -- (617.5,317.5) ;
\draw    (665,270) -- (712.5,270) ;
\draw    (712.5,270) -- (712.5,317.5) ;
\draw    (665,317.5) -- (712.5,317.5) ;
\draw    (712.5,270) -- (760,270) ;
\draw    (760,270) -- (760,317.5) ;
\draw    (712.5,317.5) -- (760,317.5) ;
\draw    (712.5,317.5) -- (712.5,365) ;
\draw    (665,365) -- (712.5,365) ;

\end{tikzpicture}

\end {figure}

\begin{defn}
    Let $\mathcal{P}$ be a poset. Let $x,y \in \mathcal{P}$. The \textit{closed interval} $[x,y]$ is the subposet of $\mathcal{P}$ defined by $z \in [x,y]$ if and only if $x \leq z \leq y$.
\end{defn}

\textit{Young's lattice} is the infinite poset whose nodes are Young diagrams of partitions, ordered by inclusion (see \cite{sagan2013symmetric}). The minimal element is the empty set, and there is no maximal element. Each $\mathbb{O}_{j}^{n}$ is isomorphic to a finite closed interval $[ \varnothing , \lambda ]$ in Young's lattice (see Theorem \ref{young_emb} for a proof), whose minimal element is the empty set and with maximal element a rectangular array of boxes $\lambda$ (indeed, $\lambda$ is the smallest rectangular array of boxes containing the minimal snake graph $G^{n,j}_{-}$).

Consider the dual cluster variable $x_{w}^{*}$ on the snake graph $G_w$ with $n \geq 1$ tiles. Let $x_M$ be any Laurent monomial of $x_{w}^{*}$, represented by the lattice path $L$ on $G_w$. We can naturally associate to $L$ a word by labeling any E step in $L$ by $a,$ and any N step in $L$ by $b$. Let $G_{M}$ be the snake graph whose shape is determined by the assignment just described. Note that $G_{M}$ has two more tiles than $G_w$. 

\begin {figure}[h!]
    \centering
    \caption{Snake graphs from lattice paths}
    \label{fig:sg_to_lp}
     \begin{tikzpicture}[x=0.75pt,y=0.75pt,yscale=-1,xscale=1]

\draw [color={rgb, 255:red, 74; green, 144; blue, 226 }  ,draw opacity=1 ][line width=2.25]    (100,270) -- (140,270) ;
\draw    (100,230) -- (140,230) ;
\draw [color={rgb, 255:red, 74; green, 144; blue, 226 }  ,draw opacity=1 ][line width=2.25]    (140,230) -- (140,270) ;
\draw    (100,230) -- (100,270) ;
\draw    (140,270) -- (180,270) ;
\draw [color={rgb, 255:red, 74; green, 144; blue, 226 }  ,draw opacity=1 ][line width=2.25]    (140,230) -- (180,230) ;
\draw    (180,230) -- (180,270) ;
\draw    (140,190) -- (140,230) ;
\draw    (140,190) -- (180,190) ;
\draw [color={rgb, 255:red, 74; green, 144; blue, 226 }  ,draw opacity=1 ][line width=2.25]    (180,190) -- (180,230) ;
\draw [color={rgb, 255:red, 74; green, 144; blue, 226 }  ,draw opacity=1 ][line width=2.25]    (290,270) -- (330,270) ;
\draw [color={rgb, 255:red, 74; green, 144; blue, 226 }  ,draw opacity=1 ][line width=2.25]    (330,230) -- (330,270) ;
\draw [color={rgb, 255:red, 74; green, 144; blue, 226 }  ,draw opacity=1 ][line width=2.25]    (330,230) -- (370,230) ;
\draw [color={rgb, 255:red, 74; green, 144; blue, 226 }  ,draw opacity=1 ][line width=2.25]    (370,190) -- (370,230) ;
\draw    (270,290) -- (310,290) ;
\draw    (310,250) -- (310,290) ;
\draw    (270,250) -- (310,250) ;
\draw    (270,250) -- (270,290) ;
\draw    (310,290) -- (350,290) ;
\draw    (350,250) -- (350,290) ;
\draw    (310,250) -- (350,250) ;
\draw    (310,250) -- (310,290) ;
\draw    (310,250) -- (350,250) ;
\draw    (350,210) -- (350,250) ;
\draw    (310,210) -- (350,210) ;
\draw    (310,210) -- (310,250) ;
\draw    (350,250) -- (390,250) ;
\draw    (390,210) -- (390,250) ;
\draw    (350,210) -- (390,210) ;
\draw    (350,210) -- (350,250) ;
\draw    (350,210) -- (390,210) ;
\draw    (390,170) -- (390,210) ;
\draw    (350,170) -- (390,170) ;
\draw    (350,170) -- (350,210) ;
\draw    (200,220) .. controls (201.67,218.33) and (203.33,218.33) .. (205,220) .. controls (206.67,221.67) and (208.33,221.67) .. (210,220) .. controls (211.67,218.33) and (213.33,218.33) .. (215,220) .. controls (216.67,221.67) and (218.33,221.67) .. (220,220) .. controls (221.67,218.33) and (223.33,218.33) .. (225,220) .. controls (226.67,221.67) and (228.33,221.67) .. (230,220) .. controls (231.67,218.33) and (233.33,218.33) .. (235,220) .. controls (236.67,221.67) and (238.33,221.67) .. (240,220) .. controls (241.67,218.33) and (243.33,218.33) .. (245,220) .. controls (246.67,221.67) and (248.33,221.67) .. (250,220) .. controls (251.67,218.33) and (253.33,218.33) .. (255,220) .. controls (256.67,221.67) and (258.33,221.67) .. (260,220) .. controls (261.67,218.33) and (263.33,218.33) .. (265,220) .. controls (266.67,221.67) and (268.33,221.67) .. (270,220) -- (278,220) ;
\draw [shift={(280,220)}, rotate = 180] [color={rgb, 255:red, 0; green, 0; blue, 0 }  ][line width=0.75]    (10.93,-3.29) .. controls (6.95,-1.4) and (3.31,-0.3) .. (0,0) .. controls (3.31,0.3) and (6.95,1.4) .. (10.93,3.29)   ;

\end{tikzpicture}
  
\end {figure}
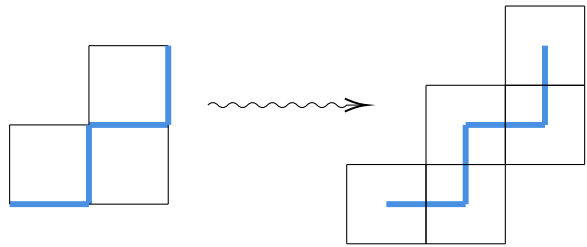

The above remarks give a map 
$$\phi: \text{Supp}(x_{w}^{*}) \longrightarrow \mathcal{L}^{n+2} , x_{M} \mapsto G_{M}.$$  

If we apply $\phi$ to the monomial weight of each node in $\mathbb{L}_{w}$ we obtain an isomorphic poset $\mathbb{I}_{w} \cong \mathbb{L}_{w},$ where each node of $\mathbb{I}_{w}$ is a snake graph in $\mathcal{L}^{n+2}$. The covering relation in $\mathbb{I}_{w}$ is the one shown in Figure \ref{fig:SG_morphism} above.

\begin {figure}[h!]
    \centering
    \caption{The posets $\mathbb{L}_{ab}$ and $\mathbb{I}_{ab} \cong \mathbb{L}_{ab}$}
    \label{fig:posets_L_and_I}
     \begin{tikzpicture}[x=0.75pt,y=0.75pt,yscale=-1,xscale=1]

\draw    (457.89,524.38) -- (470.97,524.38) ;
\draw    (457.89,511.3) -- (470.97,511.3) ;
\draw    (470.97,524.38) -- (470.97,511.3) ;
\draw    (457.89,524.38) -- (457.89,511.3) ;
\draw    (470.97,524.38) -- (484.06,524.38) ;
\draw    (470.97,511.3) -- (484.06,511.3) ;
\draw    (484.06,524.38) -- (484.06,511.3) ;
\draw    (470.97,524.38) -- (470.97,511.3) ;
\draw    (484.06,524.38) -- (497.14,524.38) ;
\draw    (484.06,511.3) -- (497.14,511.3) ;
\draw    (497.14,524.38) -- (497.14,511.3) ;
\draw    (484.06,524.38) -- (484.06,511.3) ;
\draw    (484.06,511.3) -- (497.14,511.3) ;
\draw    (484.06,498.21) -- (497.14,498.21) ;
\draw    (497.14,511.3) -- (497.14,498.21) ;
\draw    (484.06,511.3) -- (484.06,498.21) ;
\draw    (484.06,498.21) -- (497.14,498.21) ;
\draw    (484.06,485.13) -- (497.14,485.13) ;
\draw    (497.14,498.21) -- (497.14,485.13) ;
\draw    (484.06,498.21) -- (484.06,485.13) ;
\draw    (457.02,406.62) -- (470.1,406.62) ;
\draw    (457.02,393.54) -- (470.1,393.54) ;
\draw    (470.1,406.62) -- (470.1,393.54) ;
\draw    (457.02,406.62) -- (457.02,393.54) ;
\draw    (470.1,406.62) -- (483.18,406.62) ;
\draw    (470.1,393.54) -- (483.18,393.54) ;
\draw    (483.18,406.62) -- (483.18,393.54) ;
\draw    (470.1,406.62) -- (470.1,393.54) ;
\draw    (483.18,393.54) -- (496.27,393.54) ;
\draw    (483.18,406.62) -- (483.18,393.54) ;
\draw    (483.18,393.54) -- (496.27,393.54) ;
\draw    (483.18,380.45) -- (496.27,380.45) ;
\draw    (496.27,393.54) -- (496.27,380.45) ;
\draw    (483.18,393.54) -- (483.18,380.45) ;
\draw    (483.18,380.45) -- (496.27,380.45) ;
\draw    (483.18,367.37) -- (496.27,367.37) ;
\draw    (496.27,380.45) -- (496.27,367.37) ;
\draw    (483.18,380.45) -- (483.18,367.37) ;
\draw    (470.1,380.45) -- (483.18,380.45) ;
\draw    (470.1,393.54) -- (470.1,380.45) ;
\draw    (353.21,301.94) -- (366.3,301.94) ;
\draw    (353.21,288.86) -- (366.3,288.86) ;
\draw    (366.3,301.94) -- (366.3,288.86) ;
\draw    (353.21,301.94) -- (353.21,288.86) ;
\draw    (353.21,275.78) -- (366.3,275.78) ;
\draw    (366.3,288.86) -- (379.38,288.86) ;
\draw    (366.3,301.94) -- (366.3,288.86) ;
\draw    (379.38,288.86) -- (392.47,288.86) ;
\draw    (353.21,288.86) -- (353.21,275.78) ;
\draw    (379.38,288.86) -- (392.47,288.86) ;
\draw    (379.38,275.78) -- (392.47,275.78) ;
\draw    (392.47,288.86) -- (392.47,275.78) ;
\draw    (379.38,288.86) -- (379.38,275.78) ;
\draw    (379.38,275.78) -- (392.47,275.78) ;
\draw    (379.38,262.69) -- (392.47,262.69) ;
\draw    (392.47,275.78) -- (392.47,262.69) ;
\draw    (379.38,275.78) -- (379.38,262.69) ;
\draw    (366.3,275.78) -- (379.38,275.78) ;
\draw    (366.3,288.86) -- (366.3,275.78) ;
\draw    (562.56,303.69) -- (575.65,303.69) ;
\draw    (562.56,290.6) -- (575.65,290.6) ;
\draw    (575.65,303.69) -- (575.65,290.6) ;
\draw    (562.56,303.69) -- (562.56,290.6) ;
\draw    (575.65,303.69) -- (588.73,303.69) ;
\draw    (575.65,290.6) -- (588.73,290.6) ;
\draw    (588.73,303.69) -- (588.73,290.6) ;
\draw    (575.65,303.69) -- (575.65,290.6) ;
\draw    (588.73,303.69) -- (588.73,290.6) ;
\draw    (575.65,264.44) -- (588.73,264.44) ;
\draw    (588.73,277.52) -- (601.82,277.52) ;
\draw    (575.65,277.52) -- (575.65,264.44) ;
\draw    (588.73,290.6) -- (588.73,277.52) ;
\draw    (588.73,277.52) -- (601.82,277.52) ;
\draw    (588.73,264.44) -- (601.82,264.44) ;
\draw    (601.82,277.52) -- (601.82,264.44) ;
\draw    (588.73,277.52) -- (588.73,264.44) ;
\draw    (575.65,277.52) -- (588.73,277.52) ;
\draw    (575.65,290.6) -- (575.65,277.52) ;
\draw    (458.76,212.1) -- (471.84,212.1) ;
\draw    (458.76,199.01) -- (471.84,199.01) ;
\draw    (471.84,212.1) -- (471.84,199.01) ;
\draw    (458.76,212.1) -- (458.76,199.01) ;
\draw    (458.76,185.93) -- (471.84,185.93) ;
\draw    (471.84,199.01) -- (484.93,199.01) ;
\draw    (471.84,212.1) -- (471.84,199.01) ;
\draw    (458.76,199.01) -- (458.76,185.93) ;
\draw    (471.84,172.84) -- (484.93,172.84) ;
\draw    (484.93,185.93) -- (498.01,185.93) ;
\draw    (471.84,185.93) -- (471.84,172.84) ;
\draw    (484.93,199.01) -- (484.93,185.93) ;
\draw    (484.93,185.93) -- (498.01,185.93) ;
\draw    (484.93,172.84) -- (498.01,172.84) ;
\draw    (498.01,185.93) -- (498.01,172.84) ;
\draw    (484.93,185.93) -- (484.93,172.84) ;
\draw    (471.84,185.93) -- (484.93,185.93) ;
\draw    (471.84,199.01) -- (471.84,185.93) ;
\draw  [dash pattern={on 4.5pt off 4.5pt}] (445.68,506.5) .. controls (445.68,488.43) and (460.32,473.79) .. (478.39,473.79) .. controls (496.45,473.79) and (511.1,488.43) .. (511.1,506.5) .. controls (511.1,524.56) and (496.45,539.21) .. (478.39,539.21) .. controls (460.32,539.21) and (445.68,524.56) .. (445.68,506.5) -- cycle ;
\draw  [dash pattern={on 4.5pt off 4.5pt}] (445.68,388.74) .. controls (445.68,370.67) and (460.32,356.03) .. (478.39,356.03) .. controls (496.45,356.03) and (511.1,370.67) .. (511.1,388.74) .. controls (511.1,406.8) and (496.45,421.45) .. (478.39,421.45) .. controls (460.32,421.45) and (445.68,406.8) .. (445.68,388.74) -- cycle ;
\draw  [dash pattern={on 4.5pt off 4.5pt}] (341,284.06) .. controls (341,266) and (355.65,251.35) .. (373.71,251.35) .. controls (391.78,251.35) and (406.42,266) .. (406.42,284.06) .. controls (406.42,302.13) and (391.78,316.77) .. (373.71,316.77) .. controls (355.65,316.77) and (341,302.13) .. (341,284.06) -- cycle ;
\draw  [dash pattern={on 4.5pt off 4.5pt}] (550.35,284.06) .. controls (550.35,266) and (565,251.35) .. (583.06,251.35) .. controls (601.13,251.35) and (615.77,266) .. (615.77,284.06) .. controls (615.77,302.13) and (601.13,316.77) .. (583.06,316.77) .. controls (565,316.77) and (550.35,302.13) .. (550.35,284.06) -- cycle ;
\draw  [dash pattern={on 4.5pt off 4.5pt}] (445.68,192.47) .. controls (445.68,174.41) and (460.32,159.76) .. (478.39,159.76) .. controls (496.45,159.76) and (511.1,174.41) .. (511.1,192.47) .. controls (511.1,210.54) and (496.45,225.18) .. (478.39,225.18) .. controls (460.32,225.18) and (445.68,210.54) .. (445.68,192.47) -- cycle ;
\draw    (478.39,467.24) -- (478.39,426.68) ;
\draw    (450.91,362.57) -- (401.19,311.54) ;
\draw    (554.28,256.58) -- (511.1,212.1) ;
\draw    (555.58,311.54) -- (504.56,361.26) ;
\draw    (445.68,212.1) -- (399.88,256.58) ;
\draw [color={rgb, 255:red, 74; green, 144; blue, 226 }  ,draw opacity=1 ][line width=1.5]    (136,530) -- (166,530) ;
\draw    (136,500) -- (166,500) ;
\draw    (166,530) -- (166,500) ;
\draw    (136,530) -- (136,500) ;
\draw [color={rgb, 255:red, 74; green, 144; blue, 226 }  ,draw opacity=1 ][line width=1.5]    (196,530) -- (196,500) ;
\draw [color={rgb, 255:red, 74; green, 144; blue, 226 }  ,draw opacity=1 ][line width=1.5]    (166,530) -- (196,530) ;
\draw    (166,500) -- (196,500) ;
\draw    (166,470) -- (196,470) ;
\draw [color={rgb, 255:red, 74; green, 144; blue, 226 }  ,draw opacity=1 ][line width=1.5]    (196,500) -- (196,470) ;
\draw    (166,500) -- (166,470) ;
\draw [color={rgb, 255:red, 74; green, 144; blue, 226 }  ,draw opacity=1 ][line width=1.5]    (136,410) -- (166,410) ;
\draw    (136,380) -- (166,380) ;
\draw [color={rgb, 255:red, 74; green, 144; blue, 226 }  ,draw opacity=1 ][line width=1.5]    (166,410) -- (166,380) ;
\draw    (136,410) -- (136,380) ;
\draw    (196,410) -- (196,380) ;
\draw    (166,410) -- (196,410) ;
\draw [color={rgb, 255:red, 74; green, 144; blue, 226 }  ,draw opacity=1 ][line width=1.5]    (166,380) -- (196,380) ;
\draw    (166,350) -- (196,350) ;
\draw [color={rgb, 255:red, 74; green, 144; blue, 226 }  ,draw opacity=1 ][line width=1.5]    (196,380) -- (196,350) ;
\draw    (166,380) -- (166,350) ;
\draw [color={rgb, 255:red, 74; green, 144; blue, 226 }  ,draw opacity=1 ][line width=1.5]    (236,330) -- (266,330) ;
\draw    (236,300) -- (266,300) ;
\draw [color={rgb, 255:red, 74; green, 144; blue, 226 }  ,draw opacity=1 ][line width=1.5]    (266,330) -- (266,300) ;
\draw    (236,330) -- (236,300) ;
\draw    (296,330) -- (296,300) ;
\draw    (266,330) -- (296,330) ;
\draw    (266,300) -- (296,300) ;
\draw [color={rgb, 255:red, 74; green, 144; blue, 226 }  ,draw opacity=1 ][line width=1.5]    (266,270) -- (296,270) ;
\draw    (296,300) -- (296,270) ;
\draw [color={rgb, 255:red, 74; green, 144; blue, 226 }  ,draw opacity=1 ][line width=1.5]    (266,300) -- (266,270) ;
\draw    (36,330) -- (66,330) ;
\draw [color={rgb, 255:red, 74; green, 144; blue, 226 }  ,draw opacity=1 ][line width=1.5]    (36,300) -- (66,300) ;
\draw    (66,330) -- (66,300) ;
\draw [color={rgb, 255:red, 74; green, 144; blue, 226 }  ,draw opacity=1 ][line width=1.5]    (36,330) -- (36,300) ;
\draw    (96,330) -- (96,300) ;
\draw    (66,330) -- (96,330) ;
\draw [color={rgb, 255:red, 74; green, 144; blue, 226 }  ,draw opacity=1 ][line width=1.5]    (66,300) -- (96,300) ;
\draw    (66,270) -- (96,270) ;
\draw [color={rgb, 255:red, 74; green, 144; blue, 226 }  ,draw opacity=1 ][line width=1.5]    (96,300) -- (96,270) ;
\draw    (66,300) -- (66,270) ;
\draw    (136,250) -- (166,250) ;
\draw [color={rgb, 255:red, 74; green, 144; blue, 226 }  ,draw opacity=1 ][line width=1.5]    (136,220) -- (166,220) ;
\draw    (166,250) -- (166,220) ;
\draw [color={rgb, 255:red, 74; green, 144; blue, 226 }  ,draw opacity=1 ][line width=1.5]    (136,250) -- (136,220) ;
\draw    (196,250) -- (196,220) ;
\draw    (166,250) -- (196,250) ;
\draw    (166,220) -- (196,220) ;
\draw [color={rgb, 255:red, 74; green, 144; blue, 226 }  ,draw opacity=1 ][line width=1.5]    (166,190) -- (196,190) ;
\draw    (196,220) -- (196,190) ;
\draw [color={rgb, 255:red, 74; green, 144; blue, 226 }  ,draw opacity=1 ][line width=1.5]    (166,220) -- (166,190) ;
\draw    (66,340) -- (126,380) ;
\draw    (206,240) -- (256,270) ;
\draw    (86,260) -- (126,240) ;
\draw    (206,380) -- (266,340) ;
\draw    (166,420) -- (166,460) ;
\draw    (260,410) .. controls (261.67,408.33) and (263.33,408.33) .. (265,410) .. controls (266.67,411.67) and (268.33,411.67) .. (270,410) .. controls (271.67,408.33) and (273.33,408.33) .. (275,410) .. controls (276.67,411.67) and (278.33,411.67) .. (280,410) .. controls (281.67,408.33) and (283.33,408.33) .. (285,410) .. controls (286.67,411.67) and (288.33,411.67) .. (290,410) .. controls (291.67,408.33) and (293.33,408.33) .. (295,410) .. controls (296.67,411.67) and (298.33,411.67) .. (300,410) .. controls (301.67,408.33) and (303.33,408.33) .. (305,410) .. controls (306.67,411.67) and (308.33,411.67) .. (310,410) .. controls (311.67,408.33) and (313.33,408.33) .. (315,410) .. controls (316.67,411.67) and (318.33,411.67) .. (320,410) .. controls (321.67,408.33) and (323.33,408.33) .. (325,410) .. controls (326.67,411.67) and (328.33,411.67) .. (330,410) .. controls (331.67,408.33) and (333.33,408.33) .. (335,410) .. controls (336.67,411.67) and (338.33,411.67) .. (340,410) .. controls (341.67,408.33) and (343.33,408.33) .. (345,410) .. controls (346.67,411.67) and (348.33,411.67) .. (350,410) .. controls (351.67,408.33) and (353.33,408.33) .. (355,410) .. controls (356.67,411.67) and (358.33,411.67) .. (360,410) .. controls (361.67,408.33) and (363.33,408.33) .. (365,410) .. controls (366.67,411.67) and (368.33,411.67) .. (370,410) .. controls (371.67,408.33) and (373.33,408.33) .. (375,410) .. controls (376.67,411.67) and (378.33,411.67) .. (380,410) -- (388,410) ;
\draw [shift={(390,410)}, rotate = 180] [color={rgb, 255:red, 0; green, 0; blue, 0 }  ][line width=0.75]    (10.93,-3.29) .. controls (6.95,-1.4) and (3.31,-0.3) .. (0,0) .. controls (3.31,0.3) and (6.95,1.4) .. (10.93,3.29)   ;

\end{tikzpicture}
  
\end {figure}
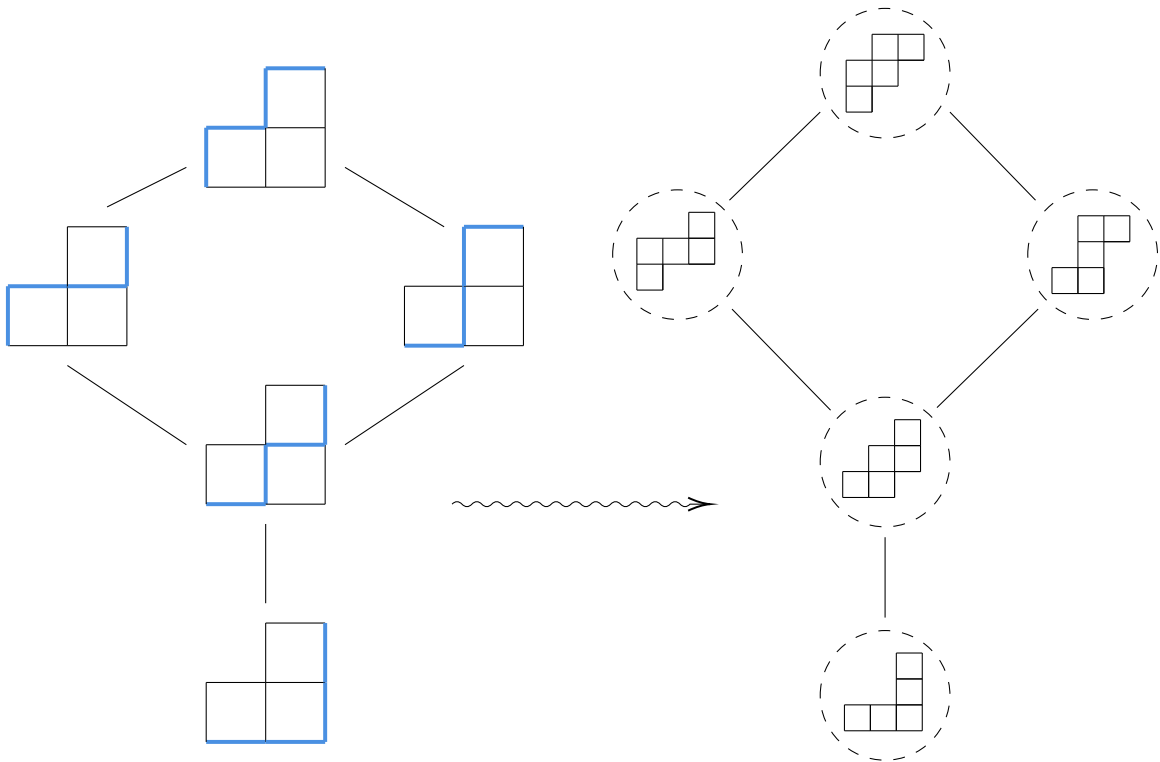

\begin{ex}
    Figure \ref{fig:posets_L_and_I} shows the poset $\mathbb{L}_{ab}$ on the left, and the isomorphic poset $\mathbb{I}_{ab}$ on the right. 
\end{ex}

The above discussion leads to the next result, which states that each poset $\mathbb{O}_{j}^{n+2}$ can be built from gluing together the lattice path posets corresponding to the snake graphs in some $\mathbb{O}_{k}^{n}$.

\begin{thm} \label{young_emb}
    Suppose the words $w , w_1 , $ and $w_2$ are each of length $n-1$. Let $\mathbb{L}_w$ be the poset of lattice paths on $G_w$. Let $\mathbb{L}_1 \doteq \mathbb{L}_{w_1}$ and $\mathbb{L}_2 \doteq \mathbb{L}_{w_{2}}$ be the posets of lattice paths on $G_1 \doteq G_{w_1}$ and $G_2 \doteq G_{w_2},$ respectively. Label their snake graph representations by $\mathbb{L}_1 \cong \mathbb{I}_1$ and $\mathbb{L}_2 \cong \mathbb{I}_2.$

\begin{enumerate}[(a)]

    \item Each poset $\mathbb{O}_{j}^{n+2}$ is isomorphic to the poset of lattice paths in a rectangular grid. 
    
    \item The rank generating function of the graded poset $\mathbb{O}_{j}^{n+2}$ is equal to a $q$-binomial coefficient.

    \item Each poset $\mathbb{O}_{j}^{n+2}$ is isomorphic to the closed interval $[\varnothing , \lambda]$ in Young's lattice. 

    \item Each lattice $\mathbb{L}_{w} \cong \mathbb{I}_w$ is isomorphic to a closed interval in one of the posets $\mathbb{O}_{j}^{n+2}$. 
    
    \item The two posets $\mathbb{L}_1 \cong \mathbb{I}_{1}$ and $\mathbb{L}_2 \cong \mathbb{I}_{2}$ are embedded as intervals into the same $\mathbb{O}_{j}^{n+2}$ if and only if $G_{1}$ and $G_{2}$ are both elements of the same poset $\mathbb{O}_{i}^{n}.$ Moreover, each $\mathbb{O}_{j}^{n+2}$ is covered by the collection of such embeddings.
\end{enumerate}
\end{thm}

\begin{proof} 

\begin{enumerate}[(a)]

    \item Suppose the minimal element of $\mathbb{O}_{j}^{n+2}$ has word $a^{k_1} b^{k_2},$ where $k_1 + k_2 - 1 = n$. Note that the covering relations in $\mathbb{O}_{j}^{n+2}$ preserve the number of $a$'s and $b$'s in the shape of a snake graph. That is, the nodes of $\mathbb{O}_{j}^{n+2}$ are precisely those snake graphs whose word has $k_1$ instances of the symbol $a$, and $k_2$ instances of the symbol $b$. This description makes clear the isomorphism between $\mathbb{O}_{j}^{n+2}$ and the lattice paths in a $k_1 \times k_2$ grid.
    
    \item As mentioned above, each poset of lattice paths in a rectangular grid is isomorphic to a closed interval $[\varnothing , \lambda]$ in Young's lattice. Now the claim follows from part $(a)$. 
    
    \item As mentioned above, the rank generating function of lattice paths in a grid is equal to some $q$-binomial coefficient. The result now follows from $(a)$. 
    
    \item Embed the underlying snake graph $G_w$ into the minimal rectangular grid containing it. This realizes the poset $\mathbb{L}_w$ as the interval $[L_{-} , L_{+}]$ inside the poset of lattice paths in this grid. As was indicated in part (a), the latter poset is isomorphic to $\mathbb{O}_{j}^{n+2}$ for some $j$. Thus the claim holds. 
    
    \item The snake graphs $G_1$ and $G_2$ are nodes in the same poset if and only both snake graphs are contained within the same minimal rectangular grid. This is true if and only if the lattice path posets $\mathbb{L}_1$ and $\mathbb{L}_2$ are embedded as intervals into the same poset of lattice paths in this minimal rectangular grid. By part $(a),$ this is true if and only if $\mathbb{L}_1$ and $\mathbb{L}_2$ are embedded into the same $\mathbb{O}_{j}^{n+2}$. That $\mathbb{O}_{j}^{n+2}$ is covered by the collection of these embeddings follows from the fact that any rectangular grid has a covering by snake graphs. 
\end{enumerate}
    
\end{proof}

\begin {figure}[h!]
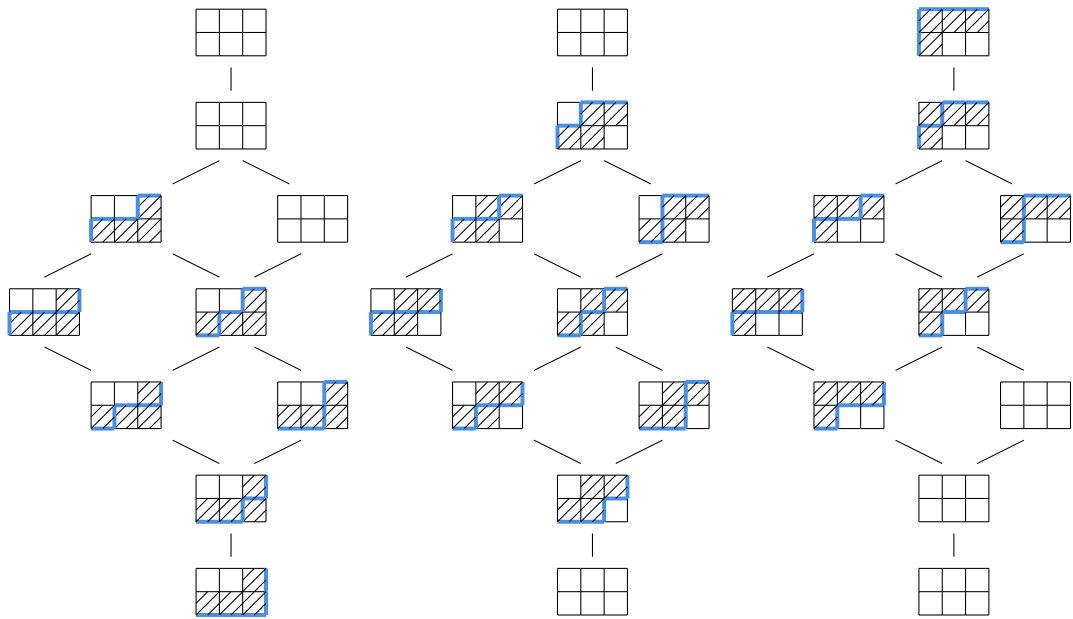

    \centering
    \caption{Three lattice path posets embedded into $\mathbb{O}_{2}^{6}$}
    \label{fig:lp_embedding_young}

  
\end {figure}

\begin{ex}
    Figure \ref{fig:lp_embedding_young} shows three copies of the same poset, each of which is isomorphic to $\mathbb{O}_{2}^{6}$. In each copy, we display one of the intervals (i.e., poset of lattice paths) in the cover from Theorem \ref{young_emb}. Note that the three distinct snake graphs corresponding to each embedded interval shown are precisely the elements of $\mathbb{O}_{1}^{4}$. 
\end{ex}

\begin{rmk}
    It is possible to dualize the above constructions (change the swap $ab \mapsto ba$ to $aa \mapsto bb$, reading off a snake graph from a perfect matching, etc.) to see explicitly that starting with perfect matchings on snake graphs gives analogous results.

    \begin {figure}[h!]
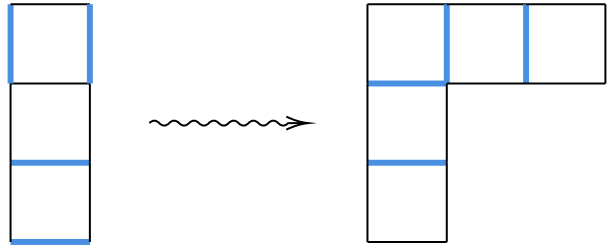

    \centering
    \caption{Snake graphs from perfect matchings}
    \label{fig:pm_sg}

\tikzset{every picture/.style={line width=0.75pt}} 

%

\end {figure}

\end{rmk}

\begin{ex}
    Below we show the ``dual'' of the poset $\mathbb{O}_{2}^{5}$. This poset can be obtained by gluing together the two posets $\mathbb{P}_{aa}$ and $\mathbb{P}_{bb}$, after realizing each perfect matching as a snake graph. Note that each node of this poset is dual (in the sense of snake graphs) to its respective node from Figure \ref{fig:poset_of_snake_graphs}.
\end{ex}

\begin {figure}[h!]
    \centering
    \caption{The poset dual to $\mathbb{O}_{2}^{5}$}
    \label{fig:snake_dual_poset}
\tikzset{every picture/.style={line width=0.75pt}} 

\begin{tikzpicture}[x=0.75pt,y=0.75pt,yscale=-1,xscale=1]

\draw    (692.11,527.38) -- (705.2,527.38) ;
\draw    (692.11,514.3) -- (705.2,514.3) ;
\draw    (705.2,527.38) -- (705.2,514.3) ;
\draw    (692.11,527.38) -- (692.11,514.3) ;
\draw    (705.2,501.21) -- (718.28,501.21) ;
\draw    (705.2,514.3) -- (718.28,514.3) ;
\draw    (718.28,514.3) -- (731.37,514.3) ;
\draw    (718.28,501.21) -- (731.37,501.21) ;
\draw    (731.37,514.3) -- (731.37,501.21) ;
\draw    (718.28,514.3) -- (718.28,501.21) ;
\draw    (718.28,501.21) -- (731.37,501.21) ;
\draw    (718.28,488.13) -- (731.37,488.13) ;
\draw    (731.37,501.21) -- (731.37,488.13) ;
\draw    (718.28,501.21) -- (718.28,488.13) ;
\draw    (706.33,419.7) -- (719.41,419.7) ;
\draw    (706.33,406.62) -- (706.33,393.54) ;
\draw    (719.41,419.7) -- (719.41,406.62) ;
\draw    (706.33,406.62) -- (719.41,406.62) ;
\draw    (706.33,393.54) -- (719.41,393.54) ;
\draw    (719.41,406.62) -- (719.41,393.54) ;
\draw    (706.33,354.28) -- (719.41,354.28) ;
\draw    (706.33,367.37) -- (719.41,367.37) ;
\draw    (706.33,380.45) -- (706.33,367.37) ;
\draw    (719.41,393.54) -- (719.41,380.45) ;
\draw    (719.41,367.37) -- (719.41,354.28) ;
\draw    (719.41,380.45) -- (719.41,367.37) ;
\draw    (706.33,380.45) -- (719.41,380.45) ;
\draw    (706.33,393.54) -- (706.33,380.45) ;
\draw    (587.44,299.94) -- (600.52,299.94) ;
\draw    (587.44,286.86) -- (600.52,286.86) ;
\draw    (600.52,299.94) -- (600.52,286.86) ;
\draw    (587.44,299.94) -- (587.44,286.86) ;
\draw    (600.52,299.94) -- (613.61,299.94) ;
\draw    (600.52,286.86) -- (613.61,286.86) ;
\draw    (600.52,299.94) -- (600.52,286.86) ;
\draw    (613.61,286.86) -- (626.69,286.86) ;
\draw    (613.61,299.94) -- (613.61,286.86) ;
\draw    (613.61,273.78) -- (626.69,273.78) ;
\draw    (626.69,286.86) -- (626.69,273.78) ;
\draw    (613.61,286.86) -- (613.61,273.78) ;
\draw    (613.61,260.69) -- (626.69,260.69) ;
\draw    (626.69,273.78) -- (626.69,260.69) ;
\draw    (613.61,273.78) -- (613.61,260.69) ;
\draw    (796.79,301.69) -- (809.87,301.69) ;
\draw    (796.79,288.6) -- (809.87,288.6) ;
\draw    (809.87,301.69) -- (809.87,288.6) ;
\draw    (796.79,301.69) -- (796.79,288.6) ;
\draw    (796.79,275.52) -- (809.87,275.52) ;
\draw    (796.79,262.44) -- (809.87,262.44) ;
\draw    (796.79,275.52) -- (796.79,262.44) ;
\draw    (796.79,288.6) -- (796.79,275.52) ;
\draw    (809.87,262.44) -- (822.96,262.44) ;
\draw    (822.96,275.52) -- (836.04,275.52) ;
\draw    (809.87,275.52) -- (809.87,262.44) ;
\draw    (822.96,275.52) -- (836.04,275.52) ;
\draw    (822.96,262.44) -- (836.04,262.44) ;
\draw    (836.04,275.52) -- (836.04,262.44) ;
\draw    (822.96,275.52) -- (822.96,262.44) ;
\draw    (809.87,275.52) -- (822.96,275.52) ;
\draw    (809.87,288.6) -- (809.87,275.52) ;
\draw    (679.9,181.93) -- (692.99,181.93) ;
\draw    (692.99,195.01) -- (706.07,195.01) ;
\draw    (692.99,181.93) -- (706.07,181.93) ;
\draw    (706.07,195.01) -- (719.16,195.01) ;
\draw    (679.9,195.01) -- (679.9,181.93) ;
\draw    (692.99,195.01) -- (692.99,181.93) ;
\draw    (732.24,181.93) -- (745.32,181.93) ;
\draw    (719.16,181.93) -- (732.24,181.93) ;
\draw    (719.16,195.01) -- (719.16,181.93) ;
\draw    (719.16,181.93) -- (732.24,181.93) ;
\draw    (719.16,195.01) -- (732.24,195.01) ;
\draw    (732.24,195.01) -- (732.24,181.93) ;
\draw    (745.32,195.01) -- (745.32,181.93) ;
\draw    (706.07,181.93) -- (719.16,181.93) ;
\draw    (706.07,195.01) -- (706.07,181.93) ;
\draw    (693.86,93.21) -- (706.94,93.21) ;
\draw    (693.86,80.13) -- (706.94,80.13) ;
\draw    (706.94,93.21) -- (706.94,80.13) ;
\draw    (693.86,93.21) -- (693.86,80.13) ;
\draw    (706.94,93.21) -- (720.03,93.21) ;
\draw    (706.94,80.13) -- (720.03,80.13) ;
\draw    (733.11,67.04) -- (733.11,53.96) ;
\draw    (720.03,67.04) -- (733.11,67.04) ;
\draw    (706.94,67.04) -- (720.03,67.04) ;
\draw    (706.94,80.13) -- (706.94,67.04) ;
\draw    (720.03,93.21) -- (720.03,80.13) ;
\draw    (706.94,53.96) -- (720.03,53.96) ;
\draw    (706.94,67.04) -- (706.94,53.96) ;
\draw    (720.03,67.04) -- (720.03,53.96) ;
\draw    (720.03,53.96) -- (733.11,53.96) ;
\draw    (720,80) -- (720,66.92) ;
\draw  [dash pattern={on 4.5pt off 4.5pt}] (674.58,386.29) .. controls (674.58,365.46) and (691.46,348.58) .. (712.29,348.58) .. controls (733.12,348.58) and (750,365.46) .. (750,386.29) .. controls (750,407.12) and (733.12,424) .. (712.29,424) .. controls (691.46,424) and (674.58,407.12) .. (674.58,386.29) -- cycle ;
\draw    (712.61,465.24) -- (712.61,427.68) ;
\draw    (710,147.56) -- (710,113) ;
\draw    (683.14,358.57) -- (636.42,310.54) ;
\draw    (787.5,253.58) -- (746.32,211.1) ;
\draw    (788.81,310.54) -- (740.78,357.26) ;
\draw    (678.9,211.1) -- (636.11,252.58) ;
\draw    (692.11,501.21) -- (705.2,501.21) ;
\draw    (705.2,514.3) -- (705.2,501.21) ;
\draw    (692.11,514.3) -- (692.11,501.21) ;
\draw    (706.33,419.7) -- (706.33,406.62) ;
\draw    (706.33,367.37) -- (706.33,354.28) ;
\draw    (613.61,299.94) -- (626.69,299.94) ;
\draw    (626.69,299.94) -- (626.69,286.86) ;
\draw    (679.9,195.01) -- (692.99,195.01) ;
\draw    (732.24,195.01) -- (745.32,195.01) ;
\draw  [dash pattern={on 4.5pt off 4.5pt}] (674.58,507.71) .. controls (674.58,486.88) and (691.46,470) .. (712.29,470) .. controls (733.12,470) and (750,486.88) .. (750,507.71) .. controls (750,528.54) and (733.12,545.42) .. (712.29,545.42) .. controls (691.46,545.42) and (674.58,528.54) .. (674.58,507.71) -- cycle ;
\draw  [dash pattern={on 4.5pt off 4.5pt}] (570,282.29) .. controls (570,261.46) and (586.88,244.58) .. (607.71,244.58) .. controls (628.54,244.58) and (645.42,261.46) .. (645.42,282.29) .. controls (645.42,303.12) and (628.54,320) .. (607.71,320) .. controls (586.88,320) and (570,303.12) .. (570,282.29) -- cycle ;
\draw  [dash pattern={on 4.5pt off 4.5pt}] (780,282.29) .. controls (780,261.46) and (796.88,244.58) .. (817.71,244.58) .. controls (838.54,244.58) and (855.42,261.46) .. (855.42,282.29) .. controls (855.42,303.12) and (838.54,320) .. (817.71,320) .. controls (796.88,320) and (780,303.12) .. (780,282.29) -- cycle ;
\draw  [dash pattern={on 4.5pt off 4.5pt}] (674.58,187.71) .. controls (674.58,166.88) and (691.46,150) .. (712.29,150) .. controls (733.12,150) and (750,166.88) .. (750,187.71) .. controls (750,208.54) and (733.12,225.42) .. (712.29,225.42) .. controls (691.46,225.42) and (674.58,208.54) .. (674.58,187.71) -- cycle ;
\draw  [dash pattern={on 4.5pt off 4.5pt}] (674.58,72.29) .. controls (674.58,51.46) and (691.46,34.58) .. (712.29,34.58) .. controls (733.12,34.58) and (750,51.46) .. (750,72.29) .. controls (750,93.12) and (733.12,110) .. (712.29,110) .. controls (691.46,110) and (674.58,93.12) .. (674.58,72.29) -- cycle ;

\end{tikzpicture}
\end {figure}
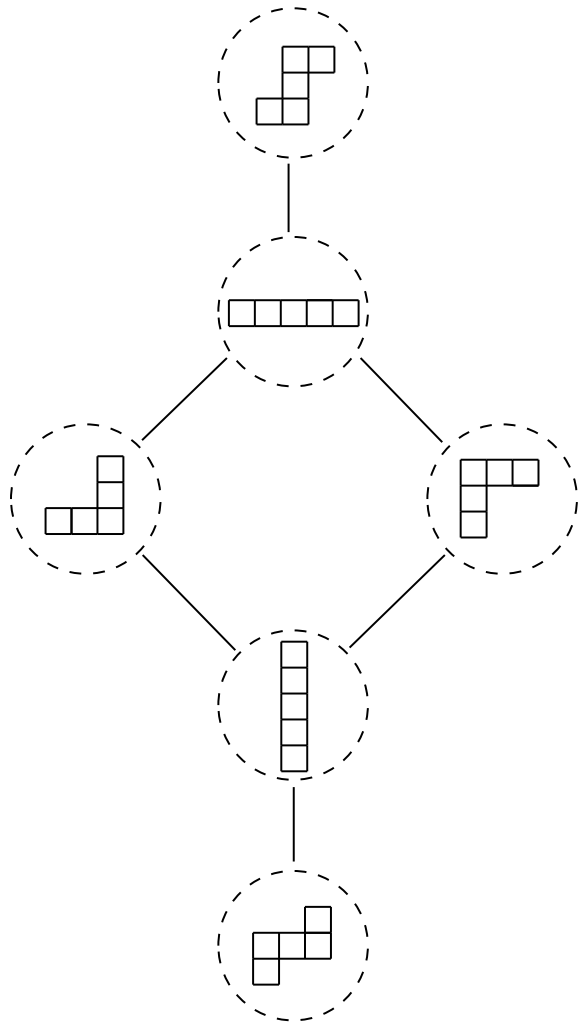

\chapter{A Recursion and Two Rank Formulas}

The first goal of this chapter is to give a recursive formula for the computation of $\mathbb{L}_w (q)$. Our second goal is to give a closed formula for $\mathbb{L}_w (q)$ in terms of products of \textit{hooks}, snake graphs whose shape is either $a^{k_1} b^{k_2}$ or $a^{k_2} b^{k_1}$ for $k_1 , k_2 \geq 2$. Thirdly, we combine this hook expansion with an interpretation of a snake graph as the central lattice path on a ``stretched'' zigzag snake graph to obtain a closed formula for $\mathbb{L}_w (q)$ in terms of the entries of the dual continued fraction $\text{CF}(w^{*})$.

\section{Lattice Path Recursion}

A straight segment is \textit{maximal} if it is not contained in any other straight segment. Decompose $G_w$ into a union of $d$ maximal overlapping straight segments as follows. Let $k_1$ be the number of tiles in the first row or column of $G_w$, $k_{d}$ the number of tiles in the last row or column of $G_w$ and $k_i + 1$ the number of tiles in the $i^{th}$ row or column of $G_w$ for $1 < i < d$. 

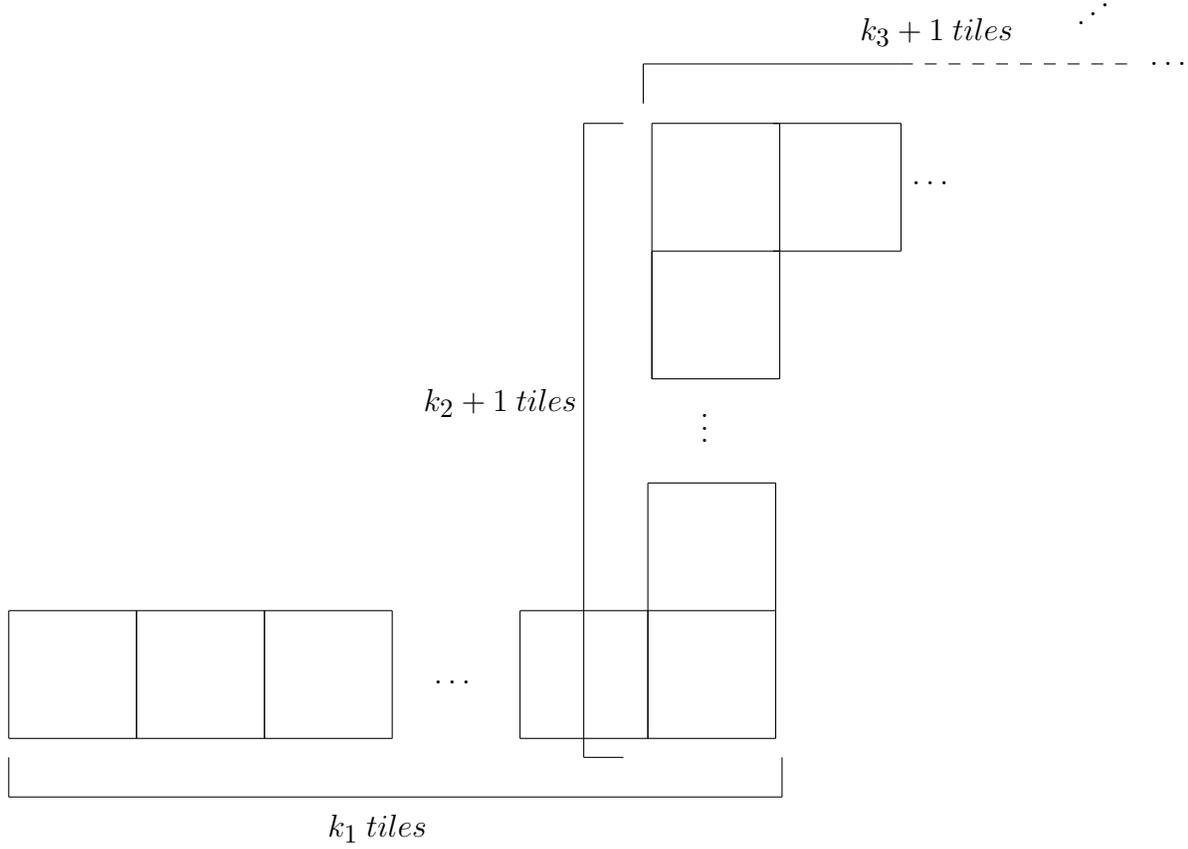
\begin {figure}[h!]
    \centering
    \caption{The maximal straight segments of $G_w$}
    \label{fig:max_straight_segments}
    \begin{tikzpicture}[x=0.75pt,y=0.75pt,yscale=-1,xscale=1]

\draw    (20,375.96) -- (20,440.43) ;
\draw    (84.47,375.96) -- (84.47,440.43) ;
\draw    (20,375.96) -- (84.47,375.96) ;
\draw    (20,440.43) -- (84.47,440.43) ;
\draw    (84.47,375.96) -- (84.47,440.43) ;
\draw    (148.94,375.96) -- (148.94,440.43) ;
\draw    (84.47,375.96) -- (148.94,375.96) ;
\draw    (84.47,440.43) -- (148.94,440.43) ;
\draw    (148.94,375.96) -- (148.94,440.43) ;
\draw    (213.4,375.96) -- (213.4,440.43) ;
\draw    (148.94,375.96) -- (213.4,375.96) ;
\draw    (148.94,440.43) -- (213.4,440.43) ;
\draw    (277.87,375.96) -- (277.87,440.43) ;
\draw    (342.34,375.96) -- (342.34,440.43) ;
\draw    (277.87,375.96) -- (342.34,375.96) ;
\draw    (277.87,440.43) -- (342.34,440.43) ;
\draw    (342.34,375.96) -- (342.34,440.43) ;
\draw    (406.81,375.96) -- (406.81,440.43) ;
\draw    (342.34,375.96) -- (406.81,375.96) ;
\draw    (342.34,440.43) -- (406.81,440.43) ;
\draw    (406.81,311.49) -- (406.81,375.96) ;
\draw    (342.34,311.49) -- (406.81,311.49) ;
\draw    (342.34,311.49) -- (342.34,375.96) ;
\draw    (344.37,194.47) -- (344.37,258.94) ;
\draw    (344.37,194.47) -- (344.37,258.94) ;
\draw    (408.84,194.47) -- (408.84,258.94) ;
\draw    (344.37,194.47) -- (408.84,194.47) ;
\draw    (344.37,258.94) -- (408.84,258.94) ;
\draw    (408.84,130) -- (408.84,194.47) ;
\draw    (344.37,130) -- (408.84,130) ;
\draw    (344.37,130) -- (344.37,194.47) ;
\draw    (470,130) -- (470,194.47) ;
\draw    (405.53,130) -- (470,130) ;
\draw    (405.53,194.47) -- (470,194.47) ;
\draw    (20,450) -- (20,470) ;
\draw    (410,450) -- (410,470) ;
\draw    (20,470) -- (410,470) ;
\draw    (310,450) -- (330,450) ;
\draw    (310,130) -- (330,130) ;
\draw    (310,450) -- (310,130) ;
\draw    (340,100) -- (340,120) ;
\draw    (340,100) -- (470,100) ;
\draw  [dash pattern={on 4.5pt off 4.5pt}]  (470,100) -- (590,100) ;

\draw (244.92,411.77) node    {$\dotsc $};
\draw (371,279) node    {$\vdots $};
\draw (486,160) node    {$\dotsc $};
\draw (206,486) node    {$k_{1} \ tiles$};
\draw (268,271) node    {$k_{2} +1\ tiles$};
\draw (488,84) node    {$k_{3} +1\ tiles$};
\draw (566.5,70) node    {$\iddots $};
\draw (606,100) node    {$\dotsc $};

\end{tikzpicture}
  
\end {figure}

Consider the dual continued fraction $\text{CF}(w^{*}) = [K_1 , K_2 , \dots , K_d],$ and assume $K_d \geq 2$ (this is no loss of generality, by the formula $[a_1 , a_2 , \dots a_m , 1] = [a_1 , a_2 , \dots a_m + 1]$). Form the continued fraction $\widehat{\text{CF}}(w^{*}) = [\widehat{K}_1 , \widehat{K}_2 , \dots , \widehat{K}_{\widehat{d}}]$ as follows: \newpage

\begin{enumerate}[(1)]
    \item If $K_1 = 1,$ then $\widehat{d} = d-1.$ In this case, the first entry of $\widehat{\text{CF}}(w^{*})$ is $\widehat{K}_{1} = K_1 + K_2 = 1 + K_2.$ The last entry of $\widehat{\text{CF}}(w^{*})$ is $\widehat{K}_{\widehat{d}-1} = K_{d}.$ The rest of the entries are $\widehat{K}_{i} = K_{i+1}+1$ for $i \neq 1, d-1.$
    
    \item If $K_1 \neq 1,$ then $\widehat{d} = d.$ In this case, the first entry of $\widehat{\text{CF}}(w^{*})$ is $\widehat{K}_1 = K_1,$ and the last entry is $\widehat{K}_d = K_d.$ The rest of the entries are $\widehat{K}_i = K_{i} + 1.$
\end{enumerate}

The next lemma follows from the previous definition, duality, and the fact that the entries in $\text{CF}(w)$ give a decomposition of $G_w$ into maximal zigzag segments. 

\begin{lem} \label{max_straight_cf}
    With the notation above, we have $k_i = \widehat{K}_i$ for each $i$.
\end{lem}

This result says that the entries in $\widehat{\text{CF}}(w^{*})$ are the lengths of maximal straight segments in the snake graph $G_w$.

We now recursively define a weight function $\rho$ on the vertices of the snake graph $G_w$ which assigns to each vertex of $G_w$ a polynomial in $q$ with positive integer coefficients. This weighting is such that if $G_{w^{\prime}}$ is a connected subsnake graph with tiles $T_1 , T_2 , \dots , T_j$, then $\rho(x,y) = \mathbb{L}_{w^{\prime}} (q),$ where $(x,y)$ is the coordinate of the NE corner of the tile $T_j$. By Lemma \ref{max_straight_cf}, these weights are functions of the entries in $\text{CF}(w^{*}).$

Suppose that $G_w$ begins by going right. Initialize the recurrence with the conditions $\rho(0,0) = 0$ and $ \rho (0,1) = \rho (x , 0) = 1, $ for $x > 1$.

\begin {figure}[h!]
    \centering
    \caption{Initial step in the recurrence}
    \label{fig:step1_rec}
        \begin{tikzpicture}[x=0.75pt,y=0.75pt,yscale=-1,xscale=1]

\draw    (46.03,473.12) -- (46.03,537.59) ;
\draw    (110.5,473.12) -- (110.5,537.59) ;
\draw    (46.03,473.12) -- (110.5,473.12) ;
\draw    (46.03,537.59) -- (110.5,537.59) ;
\draw    (110.5,473.12) -- (110.5,537.59) ;
\draw    (174.96,473.12) -- (174.96,537.59) ;
\draw    (110.5,473.12) -- (174.96,473.12) ;
\draw    (110.5,537.59) -- (174.96,537.59) ;
\draw    (174.96,473.12) -- (174.96,537.59) ;
\draw    (239.43,473.12) -- (239.43,537.59) ;
\draw    (174.96,473.12) -- (239.43,473.12) ;
\draw    (174.96,537.59) -- (239.43,537.59) ;
\draw    (303.9,473.12) -- (303.9,537.59) ;
\draw    (368.37,473.12) -- (368.37,537.59) ;
\draw    (303.9,473.12) -- (368.37,473.12) ;
\draw    (303.9,537.59) -- (368.37,537.59) ;
\draw    (368.37,473.12) -- (368.37,537.59) ;
\draw    (432.84,473.12) -- (432.84,537.59) ;
\draw    (368.37,473.12) -- (432.84,473.12) ;
\draw    (368.37,537.59) -- (432.84,537.59) ;
\draw    (432.84,408.65) -- (432.84,473.12) ;
\draw    (368.37,408.65) -- (432.84,408.65) ;
\draw    (368.37,408.65) -- (368.37,473.12) ;

\draw (270.95,508.94) node    {$\dotsc $};
\draw (48.18,544.75) node    {$0$};
\draw (112.65,544.75) node    {$1$};
\draw (177.11,544.75) node    {$1$};
\draw (241.58,544.75) node    {$1$};
\draw (370.52,544.75) node    {$1$};
\draw (434.99,544.75) node    {$1$};
\draw (306.05,544.75) node    {$1$};
\draw (41.01,480.28) node    {$1$};
\draw (397.74,387.16) node    {$\vdots $};

\end{tikzpicture}
  
\end {figure}
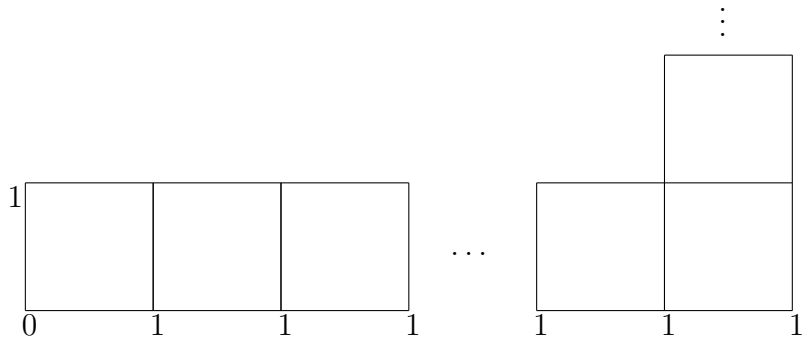

The remaining vertices in the first maximal straight segment are determined by the condition $$\rho (x , 1) = \rho (x , 0) + q \rho(x - 1 , 1).$$ The outputs of $\rho$ defined so far are displayed below in Figure \ref{fig:step2_rec}. Recall that each $[m]_q = 1 + q + \dots + q^{m-1}$ is the \textit{$q$-analog} of the integer $m$.

\begin {figure}[h!]
    \centering
    \caption{The recurrence for the first maximal straight segment of $G_w$}
    \label{fig:step2_rec}
    \begin{tikzpicture}[x=0.75pt,y=0.75pt,yscale=-1,xscale=1]

\draw    (50,150) -- (50,240) ;
\draw    (140,150) -- (140,240) ;
\draw    (50,150) -- (140,150) ;
\draw    (50,240) -- (140,240) ;
\draw    (140,150) -- (140,240) ;
\draw    (230,150) -- (230,240) ;
\draw    (140,150) -- (230,150) ;
\draw    (140,240) -- (230,240) ;
\draw    (230,150) -- (230,240) ;
\draw    (320,150) -- (320,240) ;
\draw    (230,150) -- (320,150) ;
\draw    (230,240) -- (320,240) ;
\draw    (410,150) -- (410,240) ;
\draw    (500,150) -- (500,240) ;
\draw    (410,150) -- (500,150) ;
\draw    (410,240) -- (500,240) ;
\draw    (500,150) -- (500,240) ;
\draw    (590,150) -- (590,240) ;
\draw    (500,150) -- (590,150) ;
\draw    (500,240) -- (590,240) ;
\draw    (590,60) -- (590,150) ;
\draw    (500,60) -- (590,60) ;
\draw    (500,60) -- (500,150) ;

\draw (364,200) node    {$\dotsc $};
\draw (53,250) node    {$0$};
\draw (143,250) node    {$1$};
\draw (233,250) node    {$1$};
\draw (323,250) node    {$1$};
\draw (503,250) node    {$1$};
\draw (593,250) node    {$1$};
\draw (413,250) node    {$1$};
\draw (124.5,161) node    {$[ 2]_{q}$};
\draw (214.5,161) node    {$[ 3]_{q}$};
\draw (305.5,161) node    {$[ 4]_{q}$};
\draw (558,161) node    {$[ k_{1} +1]_{q}$};
\draw (481,161) node    {$[ k_{1}]_{q}$};
\draw (44,160) node    {$1$};
\draw (549,40) node    {$\vdots $};

\end{tikzpicture}
  
\end {figure}
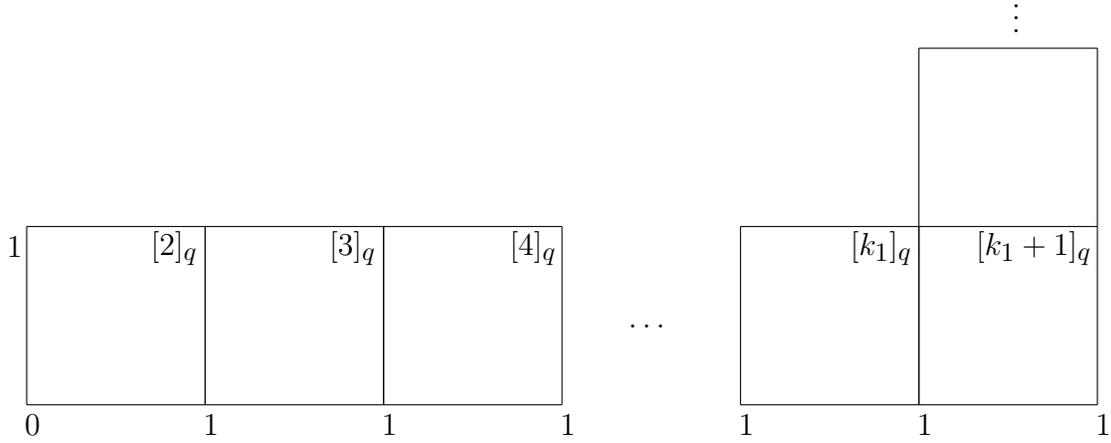

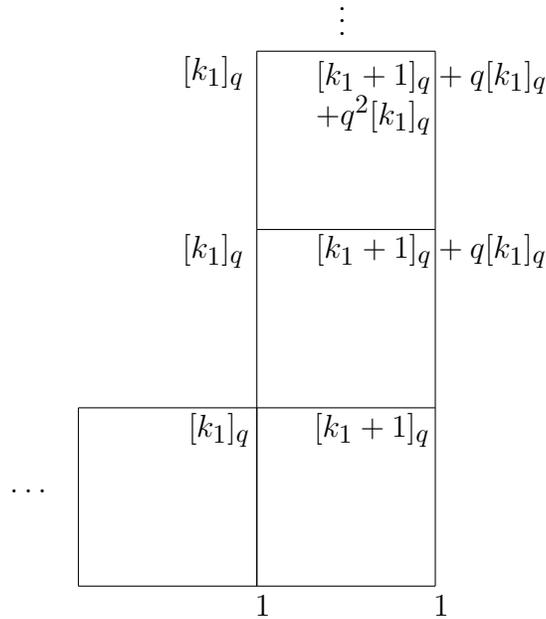
\begin {figure}[h!]
    \centering
    \caption{The recurrence for the second maximal straight segment of $G_w$}
    \label{fig:step3_rec}
    \begin{tikzpicture}[x=0.75pt,y=0.75pt,yscale=-1,xscale=1]

\draw    (120,400) -- (120,490) ;
\draw    (210,400) -- (210,490) ;
\draw    (120,400) -- (210,400) ;
\draw    (120,490) -- (210,490) ;
\draw    (210,400) -- (210,490) ;
\draw    (300,400) -- (300,490) ;
\draw    (210,400) -- (300,400) ;
\draw    (210,490) -- (300,490) ;
\draw    (300,310) -- (300,400) ;
\draw    (210,310) -- (300,310) ;
\draw    (210,310) -- (210,400) ;
\draw    (300,220) -- (300,310) ;
\draw    (210,220) -- (300,220) ;
\draw    (210,220) -- (210,310) ;

\draw (96,442) node    {$\dotsc $};
\draw (213,500) node    {$1$};
\draw (303,500) node    {$1$};
\draw (268,411) node    {$[ k_{1} +1]_{q}$};
\draw (191,411) node    {$[ k_{1}]_{q}$};
\draw (188,321) node    {$[ k_{1}]_{q}$};
\draw (188,231) node    {$[ k_{1}]_{q}$};
\draw (298,321) node    {$[ k_{1} +1]_{q} +q[ k_{1}]_{q}$};
\draw (253,200) node    {$\vdots $};
\draw (298,243) node    {$ \begin{array}{l}
[ k_{1} +1]_{q} +q[ k_{1}]_{q}\\
+q^{2}[ k_{1}]_{q}
\end{array}$};

\end{tikzpicture}
\end {figure}

Next, $\rho (k_1 - 1 , y) = \rho (k_1 - 1 , 1) =  [k_1]_q,$ for each $y > 1,$ and $$\rho(k_1 , y) = \rho(k_1 , y - 1 ) + q^{y} \rho(k_1 - 1 , y)$$
for $y > 1$. The first few outputs $\rho(k_1 - 1 , y) = [k_1]_q$ and $\rho(k_1 , y)$ for $y > 1$ are shown in Figure \ref{fig:step3_rec}.

To compute the weights of the vertices in the third (horizontal) straight segment, we use the same recurrence relation that was used to compute the weights in the first straight segment, only the initial values are different. Similarly, the vertices in the fourth (vertical) straight segment are computed using the second recurrence rule given above, except $q^y$ is replaced with $q^{y-k_2 + 1}$. We continue along in this fashion until all outputs have been computed. A similar recurrence holds when $G_w$ starts by going up.

\begin{ex}
    We use the recurrence just given to compute the lattice path rank function of the snake graph $G_{w_3}$ from \ref{fig:3_snake_graphs} above.
\end{ex}

\begin {figure}[h!]
    \centering
    \caption{Recursive computation of the rank function $\mathbb{L}_{w_3} (q)$}
    \label{fig:third_sg_rank_function}

\tikzset{every picture/.style={line width=0.75pt}} 

\begin{tikzpicture}[x=0.75pt,y=0.75pt,yscale=-1,xscale=1]

\draw    (188.26,767.71) -- (188.26,878.06) ;
\draw    (298.61,767.71) -- (298.61,878.06) ;
\draw    (188.26,767.71) -- (298.61,767.71) ;
\draw    (188.26,878.06) -- (298.61,878.06) ;
\draw    (298.61,767.71) -- (298.61,878.06) ;
\draw    (408.96,767.71) -- (408.96,878.06) ;
\draw    (298.61,767.71) -- (408.96,767.71) ;
\draw    (298.61,878.06) -- (408.96,878.06) ;
\draw    (408.96,767.71) -- (408.96,878.06) ;
\draw    (519.3,767.71) -- (519.3,878.06) ;
\draw    (408.96,767.71) -- (519.3,767.71) ;
\draw    (408.96,878.06) -- (519.3,878.06) ;
\draw    (408.96,657.36) -- (408.96,767.71) ;
\draw    (519.3,657.36) -- (519.3,767.71) ;
\draw    (408.96,657.36) -- (519.3,657.36) ;
\draw    (408.96,767.71) -- (519.3,767.71) ;
\draw    (519.3,657.36) -- (519.3,767.71) ;
\draw    (629.65,657.36) -- (629.65,767.71) ;
\draw    (519.3,657.36) -- (629.65,657.36) ;
\draw    (519.3,767.71) -- (629.65,767.71) ;
\draw    (629.65,657.36) -- (629.65,767.71) ;
\draw    (740,657.36) -- (740,767.71) ;
\draw    (629.65,657.36) -- (740,657.36) ;
\draw    (629.65,767.71) -- (740,767.71) ;
\draw  [dash pattern={on 0.84pt off 2.51pt}]  (538.57,624.08) .. controls (544.26,667) and (652.86,615.33) .. (740,657.36) ;

\draw (194.39,886.82) node    {$0$};
\draw (514.93,886.82) node    {$1$};
\draw (409.83,886.82) node    {$1$};
\draw (304.74,886.82) node    {$1$};
\draw (280.66,780.85) node    {$[ 2]_{q}$};
\draw (393.63,780.85) node    {$[ 3]_{q}$};
\draw (507.48,780.85) node    {$[ 4]_{q}$};
\draw (182.13,779.1) node    {$1$};
\draw (394.51,667) node    {$[ 3]_{q}$};
\draw (476.83,668.75) node    {$[ 4]_{q} +q[ 3]_{q}$};
\draw (618.7,666.56) node    {$[ 4]_{q} +q[ 4]_{q} +q^{3}[ 3]_{q}$};
\draw (461.07,610.07) node    {$\mathbb{L}_{w_{3}}( q) \ =\ [ 4]_{q} \ +q[ 4]_{q} \ +q^{2}[ 4]_{q} +q^{3}[ 3]_{q} =\ 1\ +\ 2q\ +\ 3q^{2} \ +\ 3q^{3} \ +\ 3q^{4} \ +\ 2q^{5} \ +\ q^{6}$};
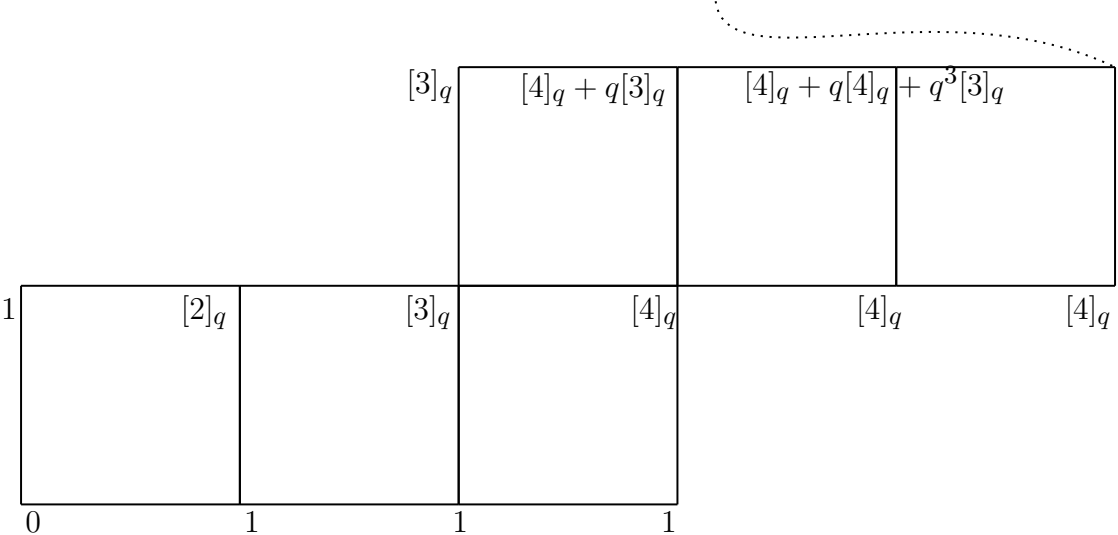
\draw (621.33,780.85) node    {$[ 4]_{q}$};
\draw (726.43,780.85) node    {$[ 4]_{q}$};

\end{tikzpicture}
\end {figure}

\section{Hook Rank Formula}

Now we provide a closed formula for any $\mathbb{L}_w (q)$. We continue to assume that $G_w$ starts by going right.

Suppose $G_w$ has $d$ straights segments. Say tile $T_i$ has an \textit{exposed NW corner} if there is one tile glued to its S edge, and one tile glued to its E edge (these three tiles correspond to a subword $ba$ in $\text{sh}(G_w)$). Similarly, we say tile $T_i$ has an \textit{exposed SE corner} if there is one tile glued to its N edge, and one tile glued to its W edge (corresponding to a subword $ab$ in $\text{sh}(G_w)$). Let $\text{NW}(G_w)$ be the set of tiles of $G_w$ with an exposed NW corner, and let $\text{SE}(G_w)$ be the set of tiles with an exposed SE corner.

If $G_w$ starts by going right, then $u \doteq \abs{\text{NW}(G_w)} = \left \lfloor{\frac{d-1}{2}}\right \rfloor.$ For each $T_i \in \text{NW}(G_w),$ let $T_{i}^{\text{SE}}$ and $T_{i}^{\text{NW}}$ respectively denote the SE and NW corner of $T_i$.

For each $i$, any lattice path in $\mathbb{L}_w$ must pass through either $T_{i}^{\text{SE}}$ or $T_{i}^{\text{NW}}$ (but not both). This allows us to partition $\mathbb{L}_w$ into $2^{u}$ sets of lattice paths by specifying which corner in each $T_i \in \text{NW}(G_w)$ a path must pass through.

Consider the natural ordering on the tiles $\text{NW}(G_w)$ inherited from the ordering of the tiles of $G_w$. If $t_i$ is one of the two corners of some $T_i \in \text{NW}(G_w)$, then the assignment $t_i \mapsto 0$ if $t_i = T_{i}^{\text{SE}}$ and $t_i \mapsto 1$ if $t_i = T_{i}^{\text{NW}}$ induces a poset isomorphism from $u$-tuples $(t_1 , t_2 , \dots , t_u)$ to $B_u$, the Boolean lattice of rank $u$. Here, one element $\sigma_2$ in $B_u$ covers another $\sigma_1$ if $\sigma_2$ can be obtained from $\sigma_1$ by switching one bit ``0'' in $\sigma_1$ to the bit ``1''.

We now build another poset $\mathcal{H}$, and give an explicit isomorphism from it to a Boolean lattice. Introduce the following notation:  

\begin{itemize}
    \item $H_{i,i+1} = 1 + q[k_i]_q [k_{i+1}]_q$
    
    \item $H_i = [k_i]_q, $
    
    \item Define \[
  H^{i,i+1} =
  \begin{cases}
                                   q^{k_i + k_{i+1} + 1} & \text{if $i \neq 1$ and $i+1 \neq l$} \\
                                   q^{k_i + k_{i+1}} & \text{if $i = 1$ or $i+1 = l.$}
  \end{cases}
\]
\end{itemize}

Define a multiplication $\circ$ on the symbols $H^{i,i+1}.$ For $i < j$, define
\[
  H^{i,i+1} \circ H^{j,j+1} =
  \begin{cases}
                                   H^{i,i+1}H^{j,j+1}  & \text{if $j \neq i+2$} \\
                                   q^{-1}(H^{i,i+1}H^{j,j+1}) & \text{if $j = i+2$.}
  \end{cases}
\]

This new multiplication can be extended to products of more than two symbols.. If no confusion will result, we omit the comma appearing in the subscripts and superscripts of these symbols (e.g., we write $H^{23}$ instead of $H^{2,3}$). We also omit the symbol ``$\circ$'' from the computations.

Starting with the minimal element $H_{12} H_{34} \dots H_{d-2 , d-1} H_{d}$ if $d$ is odd, or $H_{12} H_{34} \dots H_{d-3 , d-2} H_{d-1 , d}$ if $d$ is even, and interpreting the local procedures $$H_{i} H_{i+1} \mapsto H^{i} H^{i+1},$$ $$H_{i,i+1} H_{i+2} \mapsto H_i H^{i+1 , i+2},$$ $$H_i H_{i+1,i+2} \mapsto H^{i , i+1} H_{i+2},$$ and $$H_{i , i+1} H_{i+2 , i+3} \mapsto H_i H^{i+1 , i+2} H_{i+3}$$ as covering relations, gives a poset $\mathcal{H}$ with additional node structure.

We now give an explicit isomorphism between the poset $\mathcal{H}$ and the Boolean lattice $B_u$. Suppose for the moment that $d$ is even, so that the minimal element of the poset above is $H_{12} H_{34} \dots H_{d-1 , d}.$ Send this element $H_{12} H_{34} \dots H_{d-1 , d}$ to $(0,0,\dots , 0) \in B_n$. Now sending $H_{1} H^{23}H_{4}H_{56} \dots \mapsto (1,0,0,\dots),$ $H_{12}H_{3}H^{45}H_6 H_{78} \dots \mapsto (0,1,0,0,\dots),$ etc. induces the aforementioned isomorphism. In other words, each weight in $\mathcal{H}$ is given coordinates based on the symbols with upper subscripts that it contains. The assignment is similar when $i$ is odd.

Thus, the nodes $H_{\sigma}$ of $\mathcal{H}$ are indexed by $\sigma \in B_u$. The next result follows from the above constructions.

\begin{thm} \label{hook_expansion}
    Fix the word $w$ and consider the rank function $\mathbb{L}_w (q)$ of lattice paths on the snake graph $G_w$.  Let $B_u$ be the Boolean lattice of rank $u$. Recall the symbols $H_{\sigma}$, each parameterized by $\sigma \in B_{u}$ and representing a polynomial in $q$ with positive integer coefficients. Then we have  
$$\mathbb{L}_w (q) = \sum_{\sigma \in B_u} H_{\sigma}.$$
\end{thm}

\begin{ex}
    Consider the snake graph $G_w$ and the lattice $\mathbb{L}_w$ with its rank function $\mathbb{L}_w (q)$. 

\begin{enumerate}[(1)]
    \item If $\text{sh}(G_w) = a^{k_1 - 1} b^{k_2} a^{k_3 - 1},$ then $$\mathbb{L}_w (q) =  H_{12} H_3 +  H_1 H^{23}.$$
    
    \item If $\text{sh}(G_w) = a^{k_1 - 1} b^{k_2} a^{k_3} b^{k_4 -1},$ then $$\mathbb{L}_w (q) = H_{12} H_{34} + H_1 H^{23} H_4.$$
    
    \item If $\text{sh}(G_w) = a^{k_1 - 1} b^{k_2} a^{k_3} b^{k_4} a^{k_5-1},$ then $$\mathbb{L}_w (q) = H_{12}H_{34}H_5 + H_1 H^{23} H_4 H_5 + H_{12}H_{3} H^{45} + H_{1} H^{23}H^{45}.$$
    
    \item If $\text{sh}(G_w) = a^{k_1 - 1} b^{k_2} a^{k_3} b^{k_4} a^{k_5} b^{k_6 - 1},$ then 
    $$\mathbb{L}_w (q) = H_{12}H_{34}H_{56} + H_1 H^{23} H_4 H_{56} + H_{12}H_{3} H^{45}H_{6} + H_{1} H^{23} H^{45} H_6.$$
\end{enumerate}
\end{ex}

Similar formulas can be derived for when $G_w$ starts by going up, e.g., for $G_w$ such that $\text{sh}(G_w) = b^{k_1 - 1} a^{k_2} b^{k_3 - 1},$ we have $\mathbb{L}_w (q) = H_{1}H_{23} + H^{12}H_3.$

\section{Fibonacci Rank Formula}

The hook expansion formula from the previous section can be refined to an explicit closed formula for $\mathbb{L}_w (q),$ as a sum over ``face-weighted'' products of $q$-deformations of the entries $k_i$ from $\text{CF}(w)^{*}.$ Each term in this formula is the weight of a lattice path $L$ on a zigzag snake graph $\mathbb{G}_w$ with $d-1$ tiles.

By the weight of $L$, we mean the product of the weights attached to the edges of $L$, multiplied by the product of face weights in the symmetric difference $L \ominus L_{-}$ (recall $L_{-}$ is the minimal lattice path from Definition \ref{L_min_def} above). Each edge weight is either equal to $1$ or $[k_i]_q$ for some $i$, and each $[k_i]_q$ is used precisely once. Each face weight is either equal to $q$, or equal to a power of $q$ with exponent equal to the number of tiles in some hook of $G_w.$

The snake graph $\mathbb{G}_w$ is built from $G_w$ by first reflecting $G_w$ about the antidiagonal $a_1$, and then treating $G_w$ as the ``middle'' lattice path on $\mathbb{G}_w$. This construction in the case that $G_w$ is built from four maximal straight segments is illustrated below.

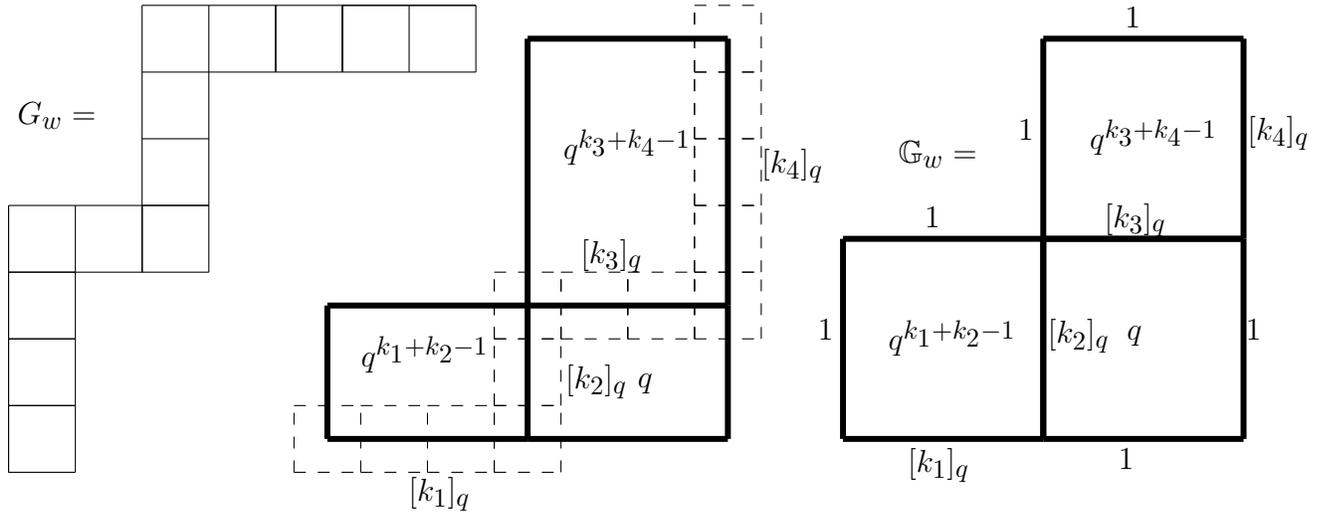
\begin {figure}[h!]
    \centering
    \caption{The construction of $\mathbb{G}_w$ from $G_w$}
    \label{fig:mathbb_G}
    \begin{tikzpicture}[x=0.75pt,y=0.75pt,yscale=-1,xscale=1]

\draw  [dash pattern={on 4.5pt off 4.5pt}]  (267.26,252.67) -- (300.93,252.67) ;
\draw  [dash pattern={on 4.5pt off 4.5pt}]  (267.26,219) -- (300.93,219) ;
\draw  [dash pattern={on 4.5pt off 4.5pt}]  (267.26,252.67) -- (267.26,219) ;
\draw  [dash pattern={on 4.5pt off 4.5pt}]  (300.93,252.67) -- (300.93,219) ;
\draw  [dash pattern={on 4.5pt off 4.5pt}]  (300.93,252.67) -- (334.59,252.67) ;
\draw  [dash pattern={on 4.5pt off 4.5pt}]  (300.93,219) -- (334.59,219) ;
\draw  [dash pattern={on 4.5pt off 4.5pt}]  (300.93,252.67) -- (300.93,219) ;
\draw  [dash pattern={on 4.5pt off 4.5pt}]  (334.59,252.67) -- (334.59,219) ;
\draw  [dash pattern={on 4.5pt off 4.5pt}]  (300.93,219) -- (300.93,185.33) ;
\draw  [dash pattern={on 4.5pt off 4.5pt}]  (300.93,219) -- (334.59,219) ;
\draw  [dash pattern={on 4.5pt off 4.5pt}]  (300.93,185.33) -- (334.59,185.33) ;
\draw  [dash pattern={on 4.5pt off 4.5pt}]  (300.93,219) -- (300.93,185.33) ;
\draw  [dash pattern={on 4.5pt off 4.5pt}]  (334.59,219) -- (334.59,185.33) ;
\draw  [dash pattern={on 4.5pt off 4.5pt}]  (300.93,185.33) -- (300.93,151.67) ;
\draw  [dash pattern={on 4.5pt off 4.5pt}]  (300.93,185.33) -- (334.59,185.33) ;
\draw  [dash pattern={on 4.5pt off 4.5pt}]  (300.93,151.67) -- (334.59,151.67) ;
\draw  [dash pattern={on 4.5pt off 4.5pt}]  (300.93,185.33) -- (300.93,151.67) ;
\draw  [dash pattern={on 4.5pt off 4.5pt}]  (334.59,185.33) -- (334.59,151.67) ;
\draw  [dash pattern={on 4.5pt off 4.5pt}]  (334.59,185.33) -- (334.59,151.67) ;
\draw  [dash pattern={on 4.5pt off 4.5pt}]  (334.59,185.33) -- (368.26,185.33) ;
\draw  [dash pattern={on 4.5pt off 4.5pt}]  (334.59,151.67) -- (368.26,151.67) ;
\draw  [dash pattern={on 4.5pt off 4.5pt}]  (334.59,185.33) -- (334.59,151.67) ;
\draw  [dash pattern={on 4.5pt off 4.5pt}]  (368.26,185.33) -- (368.26,151.67) ;
\draw  [dash pattern={on 4.5pt off 4.5pt}]  (368.26,185.33) -- (368.26,151.67) ;
\draw  [dash pattern={on 4.5pt off 4.5pt}]  (368.26,185.33) -- (401.93,185.33) ;
\draw  [dash pattern={on 4.5pt off 4.5pt}]  (368.26,151.67) -- (401.93,151.67) ;
\draw  [dash pattern={on 4.5pt off 4.5pt}]  (368.26,185.33) -- (368.26,151.67) ;
\draw  [dash pattern={on 4.5pt off 4.5pt}]  (401.93,185.33) -- (401.93,151.67) ;
\draw  [dash pattern={on 4.5pt off 4.5pt}]  (401.93,185.33) -- (401.93,151.67) ;
\draw  [dash pattern={on 4.5pt off 4.5pt}]  (401.93,185.33) -- (435.59,185.33) ;
\draw  [dash pattern={on 4.5pt off 4.5pt}]  (401.93,151.67) -- (435.59,151.67) ;
\draw  [dash pattern={on 4.5pt off 4.5pt}]  (401.93,185.33) -- (401.93,151.67) ;
\draw  [dash pattern={on 4.5pt off 4.5pt}]  (435.59,185.33) -- (435.59,151.67) ;
\draw  [dash pattern={on 4.5pt off 4.5pt}]  (401.93,151.67) -- (401.93,118) ;
\draw  [dash pattern={on 4.5pt off 4.5pt}]  (401.93,151.67) -- (435.59,151.67) ;
\draw  [dash pattern={on 4.5pt off 4.5pt}]  (401.93,118) -- (435.59,118) ;
\draw  [dash pattern={on 4.5pt off 4.5pt}]  (401.93,151.67) -- (401.93,118) ;
\draw  [dash pattern={on 4.5pt off 4.5pt}]  (435.59,151.67) -- (435.59,118) ;
\draw  [dash pattern={on 4.5pt off 4.5pt}]  (401.93,118) -- (401.93,84.33) ;
\draw  [dash pattern={on 4.5pt off 4.5pt}]  (401.93,118) -- (435.59,118) ;
\draw  [dash pattern={on 4.5pt off 4.5pt}]  (401.93,84.33) -- (435.59,84.33) ;
\draw  [dash pattern={on 4.5pt off 4.5pt}]  (401.93,118) -- (401.93,84.33) ;
\draw  [dash pattern={on 4.5pt off 4.5pt}]  (435.59,118) -- (435.59,84.33) ;
\draw  [dash pattern={on 4.5pt off 4.5pt}]  (401.93,84.33) -- (401.93,50.67) ;
\draw  [dash pattern={on 4.5pt off 4.5pt}]  (401.93,84.33) -- (435.59,84.33) ;
\draw  [dash pattern={on 4.5pt off 4.5pt}]  (401.93,50.67) -- (435.59,50.67) ;
\draw  [dash pattern={on 4.5pt off 4.5pt}]  (401.93,84.33) -- (401.93,50.67) ;
\draw  [dash pattern={on 4.5pt off 4.5pt}]  (435.59,84.33) -- (435.59,50.67) ;
\draw  [dash pattern={on 4.5pt off 4.5pt}]  (401.93,50.67) -- (401.93,17) ;
\draw  [dash pattern={on 4.5pt off 4.5pt}]  (401.93,50.67) -- (435.59,50.67) ;
\draw  [dash pattern={on 4.5pt off 4.5pt}]  (401.93,17) -- (435.59,17) ;
\draw  [dash pattern={on 4.5pt off 4.5pt}]  (401.93,50.67) -- (401.93,17) ;
\draw  [dash pattern={on 4.5pt off 4.5pt}]  (435.59,50.67) -- (435.59,17) ;
\draw    (56,252.67) -- (56,219) ;
\draw    (56,219) -- (89.67,219) ;
\draw    (56,252.67) -- (56,219) ;
\draw    (89.67,252.67) -- (89.67,219) ;
\draw    (56,219) -- (56,185.33) ;
\draw    (56,219) -- (89.67,219) ;
\draw    (56,185.33) -- (89.67,185.33) ;
\draw    (56,219) -- (56,185.33) ;
\draw    (89.67,219) -- (89.67,185.33) ;
\draw    (56,185.33) -- (56,151.67) ;
\draw    (56,185.33) -- (89.67,185.33) ;
\draw    (56,151.67) -- (89.67,151.67) ;
\draw    (56,185.33) -- (56,151.67) ;
\draw    (89.67,185.33) -- (89.67,151.67) ;
\draw    (56,151.67) -- (56,118) ;
\draw    (56,151.67) -- (89.67,151.67) ;
\draw    (56,118) -- (89.67,118) ;
\draw    (56,151.67) -- (56,118) ;
\draw    (89.67,151.67) -- (89.67,118) ;
\draw [line width=2.25]    (216.76,235.83) -- (317.76,235.83) ;
\draw [line width=2.25]    (317.76,168.5) -- (418.76,168.5) ;
\draw [line width=2.25]    (418.76,168.5) -- (418.76,33.83) ;
\draw [line width=2.25]    (317.76,235.83) -- (317.76,168.5) ;
\draw [line width=2.25]    (418.76,235.83) -- (317.76,235.83) ;
\draw [line width=2.25]    (418.76,33.83) -- (317.76,33.83) ;
\draw [line width=2.25]    (317.76,33.83) -- (317.76,168.5) ;
\draw [line width=2.25]    (216.76,235.83) -- (216.76,168.5) ;
\draw [line width=2.25]    (216.76,168.5) -- (317.76,168.5) ;
\draw [line width=2.25]    (418.76,235.83) -- (418.76,168.5) ;
\draw  [dash pattern={on 4.5pt off 4.5pt}]  (233.59,252.67) -- (267.26,252.67) ;
\draw  [dash pattern={on 4.5pt off 4.5pt}]  (233.59,219) -- (267.26,219) ;
\draw  [dash pattern={on 4.5pt off 4.5pt}]  (233.59,252.67) -- (233.59,219) ;
\draw  [dash pattern={on 4.5pt off 4.5pt}]  (199.93,252.67) -- (233.59,252.67) ;
\draw  [dash pattern={on 4.5pt off 4.5pt}]  (199.93,219) -- (233.59,219) ;
\draw  [dash pattern={on 4.5pt off 4.5pt}]  (199.93,252.67) -- (199.93,219) ;
\draw [line width=2.25]    (476.83,235.83) -- (577.83,235.83) ;
\draw [line width=2.25]    (476.83,134.83) -- (577.83,134.83) ;
\draw [line width=2.25]    (476.83,134.83) -- (476.83,235.83) ;
\draw [line width=2.25]    (577.83,134.83) -- (577.83,235.83) ;
\draw [line width=2.25]    (577.83,235.83) -- (678.83,235.83) ;
\draw [line width=2.25]    (577.83,134.83) -- (678.83,134.83) ;
\draw [line width=2.25]    (577.83,134.83) -- (577.83,235.83) ;
\draw [line width=2.25]    (678.83,134.83) -- (678.83,235.83) ;
\draw [line width=2.25]    (577.83,134.83) -- (678.83,134.83) ;
\draw [line width=2.25]    (577.83,33.83) -- (678.83,33.83) ;
\draw [line width=2.25]    (577.83,33.83) -- (577.83,134.83) ;
\draw [line width=2.25]    (678.83,33.83) -- (678.83,134.83) ;
\draw    (56,252.67) -- (89.67,252.67) ;
\draw    (89.67,151.67) -- (123.33,151.67) ;
\draw    (89.67,118) -- (123.33,118) ;
\draw    (123.33,151.67) -- (123.33,118) ;
\draw    (123.33,151.67) -- (123.33,118) ;
\draw    (123.33,151.67) -- (157,151.67) ;
\draw    (123.33,118) -- (157,118) ;
\draw    (123.33,151.67) -- (123.33,118) ;
\draw    (157,151.67) -- (157,118) ;
\draw    (123.33,118) -- (157,118) ;
\draw    (123.33,84.33) -- (157,84.33) ;
\draw    (123.33,118) -- (123.33,84.33) ;
\draw    (157,118) -- (157,84.33) ;
\draw    (123.33,50.67) -- (157,50.67) ;
\draw    (123.33,84.33) -- (123.33,50.67) ;
\draw    (157,84.33) -- (157,50.67) ;
\draw    (123.33,17) -- (157,17) ;
\draw    (123.33,50.67) -- (123.33,17) ;
\draw    (157,50.67) -- (157,17) ;
\draw    (157,50.67) -- (190.67,50.67) ;
\draw    (157,17) -- (190.67,17) ;
\draw    (190.67,50.67) -- (190.67,17) ;
\draw    (190.67,50.67) -- (190.67,17) ;
\draw    (190.67,50.67) -- (224.33,50.67) ;
\draw    (190.67,17) -- (224.33,17) ;
\draw    (190.67,50.67) -- (190.67,17) ;
\draw    (224.33,50.67) -- (224.33,17) ;
\draw    (224.33,50.67) -- (224.33,17) ;
\draw    (224.33,50.67) -- (258,50.67) ;
\draw    (224.33,17) -- (258,17) ;
\draw    (224.33,50.67) -- (224.33,17) ;
\draw    (258,50.67) -- (258,17) ;
\draw    (224.33,17) -- (258,17) ;
\draw    (258,50.67) -- (291.67,50.67) ;
\draw    (258,17) -- (291.67,17) ;
\draw    (291.67,50.67) -- (291.67,17) ;
\draw    (258,17) -- (291.67,17) ;

\draw (256.52,253.93) node [anchor=north west][inner sep=0.75pt]    {$[ k_{1}]_{q}$};
\draw (335.58,196.9) node [anchor=north west][inner sep=0.75pt]    {$[ k_{2}]_{q}$};
\draw (343.21,133.57) node [anchor=north west][inner sep=0.75pt]    {$[ k_{3}]_{q}$};
\draw (434.42,88.33) node [anchor=north west][inner sep=0.75pt]    {$[ k_{4}]_{q}$};
\draw (371.98,201.58) node [anchor=north west][inner sep=0.75pt]    {$q$};
\draw (334.6,78.18) node [anchor=north west][inner sep=0.75pt]    {$q^{k_{3} +k_{4} -1}$};
\draw (232.6,182.91) node [anchor=north west][inner sep=0.75pt]    {$q^{k_{1} +k_{2} -1}$};
\draw (679.51,71.49) node [anchor=north west][inner sep=0.75pt]    {$[ k_{4}]_{q}$};
\draw (508.18,239.09) node [anchor=north west][inner sep=0.75pt]    {$[ k_{1}]_{q}$};
\draw (578.82,174.18) node [anchor=north west][inner sep=0.75pt]    {$[ k_{2}]_{q}$};
\draw (607.49,115.58) node [anchor=north west][inner sep=0.75pt]    {$[ k_{3}]_{q}$};
\draw (614.44,239.25) node [anchor=north west][inner sep=0.75pt]    {$1$};
\draw (678.73,174.33) node [anchor=north west][inner sep=0.75pt]    {$1$};
\draw (462.94,174.33) node [anchor=north west][inner sep=0.75pt]    {$1$};
\draw (516.81,117.42) node [anchor=north west][inner sep=0.75pt]    {$1$};
\draw (563.94,73.33) node [anchor=north west][inner sep=0.75pt]    {$1$};
\draw (617.81,16.42) node [anchor=north west][inner sep=0.75pt]    {$1$};
\draw (600.1,73.18) node [anchor=north west][inner sep=0.75pt]    {$q^{k_{3} +k_{4} -1}$};
\draw (498.57,174.49) node [anchor=north west][inner sep=0.75pt]    {$q^{k_{1} +k_{2} -1}$};
\draw (618.86,177.33) node [anchor=north west][inner sep=0.75pt]    {$q$};
\draw (59.3,64.76) node [anchor=north west][inner sep=0.75pt]    {$G_{w} =\ $};
\draw (503.52,84.18) node [anchor=north west][inner sep=0.75pt]    {$\mathbb{G}_{w} =\ $};

\end{tikzpicture}
  
\end {figure}

\newpage

Summing over the weights of the paths on this zigzag snake graph gives a closed expression for the rank function $\mathbb{L}_w (q)$. 

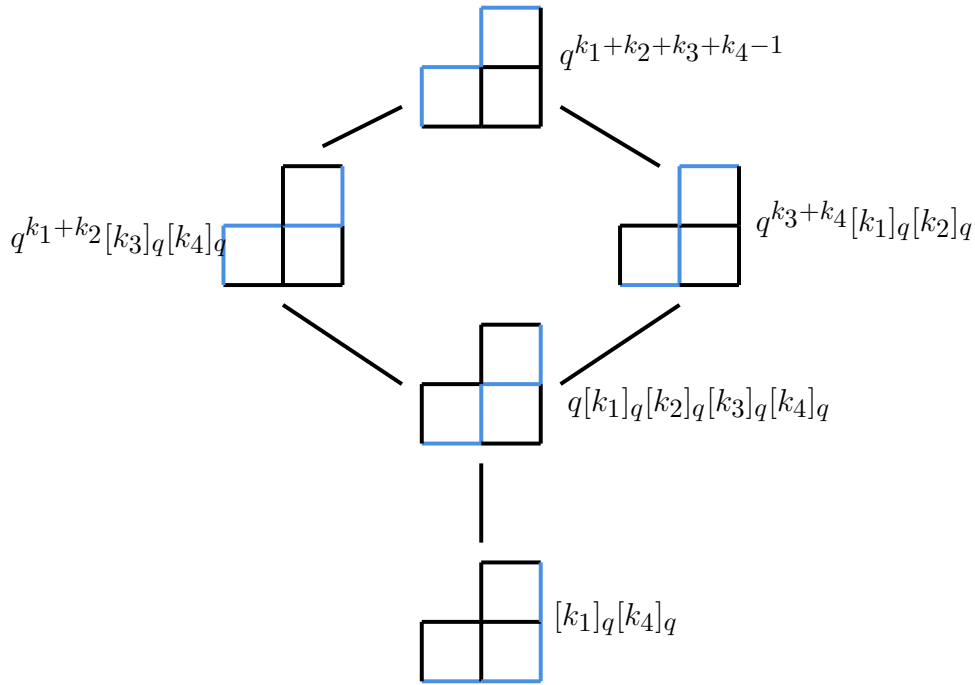
\begin {figure}[h!]
    \centering
    \caption{The poset $\mathbb{G}_{ab},$ its node weights, and the rank function $\mathbb{L}_{w} (q)$}
    \label{fig:mathbb_G_paths}

\tikzset{every picture/.style={line width=0.75pt}} 

\begin{tikzpicture}[x=0.75pt,y=0.75pt,yscale=-1,xscale=1]

\draw [color={rgb, 255:red, 74; green, 144; blue, 226 }  ,draw opacity=1 ][line width=1.5]    (261,383) -- (291,383) ;
\draw [line width=1.5]    (261,353) -- (291,353) ;
\draw [line width=1.5]    (291,383) -- (291,353) ;
\draw [line width=1.5]    (261,383) -- (261,353) ;
\draw [color={rgb, 255:red, 74; green, 144; blue, 226 }  ,draw opacity=1 ][line width=1.5]    (321,383) -- (321,353) ;
\draw [color={rgb, 255:red, 74; green, 144; blue, 226 }  ,draw opacity=1 ][line width=1.5]    (291,383) -- (321,383) ;
\draw [line width=1.5]    (291,353) -- (321,353) ;
\draw [line width=1.5]    (291,323) -- (321,323) ;
\draw [color={rgb, 255:red, 74; green, 144; blue, 226 }  ,draw opacity=1 ][line width=1.5]    (321,353) -- (321,323) ;
\draw [line width=1.5]    (291,353) -- (291,323) ;
\draw [color={rgb, 255:red, 74; green, 144; blue, 226 }  ,draw opacity=1 ][line width=1.5]    (261,263) -- (291,263) ;
\draw [line width=1.5]    (261,233) -- (291,233) ;
\draw [color={rgb, 255:red, 74; green, 144; blue, 226 }  ,draw opacity=1 ][line width=1.5]    (291,263) -- (291,233) ;
\draw [line width=1.5]    (261,263) -- (261,233) ;
\draw [line width=1.5]    (321,263) -- (321,233) ;
\draw [line width=1.5]    (291,263) -- (321,263) ;
\draw [color={rgb, 255:red, 74; green, 144; blue, 226 }  ,draw opacity=1 ][line width=1.5]    (291,233) -- (321,233) ;
\draw [line width=1.5]    (291,203) -- (321,203) ;
\draw [color={rgb, 255:red, 74; green, 144; blue, 226 }  ,draw opacity=1 ][line width=1.5]    (321,233) -- (321,203) ;
\draw [line width=1.5]    (291,233) -- (291,203) ;
\draw [color={rgb, 255:red, 74; green, 144; blue, 226 }  ,draw opacity=1 ][line width=1.5]    (361,183) -- (391,183) ;
\draw [line width=1.5]    (361,153) -- (391,153) ;
\draw [color={rgb, 255:red, 74; green, 144; blue, 226 }  ,draw opacity=1 ][line width=1.5]    (391,183) -- (391,153) ;
\draw [line width=1.5]    (361,183) -- (361,153) ;
\draw [line width=1.5]    (421,183) -- (421,153) ;
\draw [line width=1.5]    (391,183) -- (421,183) ;
\draw [line width=1.5]    (391,153) -- (421,153) ;
\draw [color={rgb, 255:red, 74; green, 144; blue, 226 }  ,draw opacity=1 ][line width=1.5]    (391,123) -- (421,123) ;
\draw [line width=1.5]    (421,153) -- (421,123) ;
\draw [color={rgb, 255:red, 74; green, 144; blue, 226 }  ,draw opacity=1 ][line width=1.5]    (391,153) -- (391,123) ;
\draw [line width=1.5]    (161,183) -- (191,183) ;
\draw [color={rgb, 255:red, 74; green, 144; blue, 226 }  ,draw opacity=1 ][line width=1.5]    (161,153) -- (191,153) ;
\draw [line width=1.5]    (191,183) -- (191,153) ;
\draw [color={rgb, 255:red, 74; green, 144; blue, 226 }  ,draw opacity=1 ][line width=1.5]    (161,183) -- (161,153) ;
\draw [line width=1.5]    (221,183) -- (221,153) ;
\draw [line width=1.5]    (191,183) -- (221,183) ;
\draw [color={rgb, 255:red, 74; green, 144; blue, 226 }  ,draw opacity=1 ][line width=1.5]    (191,153) -- (221,153) ;
\draw [line width=1.5]    (191,123) -- (221,123) ;
\draw [color={rgb, 255:red, 74; green, 144; blue, 226 }  ,draw opacity=1 ][line width=1.5]    (221,153) -- (221,123) ;
\draw [line width=1.5]    (191,153) -- (191,123) ;
\draw [line width=1.5]    (261,103) -- (291,103) ;
\draw [color={rgb, 255:red, 74; green, 144; blue, 226 }  ,draw opacity=1 ][line width=1.5]    (261,73) -- (291,73) ;
\draw [line width=1.5]    (291,103) -- (291,73) ;
\draw [color={rgb, 255:red, 74; green, 144; blue, 226 }  ,draw opacity=1 ][line width=1.5]    (261,103) -- (261,73) ;
\draw [line width=1.5]    (321,103) -- (321,73) ;
\draw [line width=1.5]    (291,103) -- (321,103) ;
\draw [line width=1.5]    (291,73) -- (321,73) ;
\draw [color={rgb, 255:red, 74; green, 144; blue, 226 }  ,draw opacity=1 ][line width=1.5]    (291,43) -- (321,43) ;
\draw [line width=1.5]    (321,73) -- (321,43) ;
\draw [color={rgb, 255:red, 74; green, 144; blue, 226 }  ,draw opacity=1 ][line width=1.5]    (291,73) -- (291,43) ;
\draw [line width=1.5]    (191,193) -- (251,233) ;
\draw [line width=1.5]    (331,93) -- (381,123) ;
\draw [line width=1.5]    (211,113) -- (251,93) ;
\draw [line width=1.5]    (331,233) -- (391,193) ;
\draw [line width=1.5]    (291,273) -- (291,313) ;

\draw (326,341) node [anchor=north west][inner sep=0.75pt]    {$[ k_{1}]_{q}[ k_{4}]_{q}$};
\draw (333,233) node [anchor=north west][inner sep=0.75pt]    {$q[ k_{1}]_{q}[ k_{2}]_{q}[ k_{3}]_{q}[ k_{4}]_{q}$};
\draw (428,139) node [anchor=north west][inner sep=0.75pt]    {$q^{k_{3} +k_{4}}[ k_{1}]_{q}[ k_{2}]_{q}$};
\draw (52,147) node [anchor=north west][inner sep=0.75pt]    {$q^{k_{1} +k_{2}}[ k_{3}]_{q}[ k_{4}]_{q}$};
\draw (329,55) node [anchor=north west][inner sep=0.75pt]    {$q^{k_{1} +k_{2} +k_{3} +k_{4} -1}$};
\draw (17,401) node [anchor=north west][inner sep=0.75pt]  [font=\footnotesize]  {$\mathbb{L}_{w} \ ( q) \ =[ k_{1}]_{q}[ k_{4}]_{q} +q[ k_{1}]_{q}[ k_{2}]_{q}[ k_{3}]_{q}[ k_{4}]_{q} +q^{k_{1} +k_{2}}[ k_{3}]_{q}[ k_{4}]_{q} +q^{k_{3} +k_{4}}[ k_{1}]_{q}[ k_{2}]_{q} +q^{k_{1} +k_{2} +k_{3} +k_{4} -1}$};

\end{tikzpicture}

\end {figure}

In particular, each such rank function expansion has a Fibonacci number of terms, which are arranged into a Fibonacci cube. As mentioned above, the snake graph $G_w$ is canonically associated to the node labeled by $000\dots 0$ in this Fibonacci cube (see Example \ref{fib_cube_000}).

Let $\mathbb{L}(\mathbb{G}_w)$ denote the poset of lattice paths on $\mathbb{G}_w.$ For any element $L \in \mathbb{L}(\mathbb{G}_w)$, let its weight as described above be denoted $q_{L}.$

\begin{thm}
    Let $w$ be any word. Suppose the snake graph $G_w$ is built from $d$ maximal straight segments. Then the rank function $\mathbb{L}_w (q)$ can be written as
    $$\mathbb{L}_w (q) = \sum_{L \in \mathbb{L} (\mathbb{G}_w)} q_{L}$$
    where $\mathbb{G}_w$ is the snake graph with $d-1$ tiles defined above, and $\mathbb{L} (\mathbb{G}_w)$ is its poset of lattice paths.
\end{thm}

\begin{proof} This is essentially the hook expansion formula above, except we decompose the set of lattice paths using the corners in $\text{NW}(G_w) \cup \text{SW}(G_w)$ instead of just $\text{NW}(G_w)$. 

\end{proof}

\chapter {Unimodality and Symmetry}

Recently, Morier-Genoud and Ovsienko gave a new notion of $q$-deformed continued fractions and rational numbers \cite{morier2018q}. The $q$-deformation $[\frac{r}{s} ]_q$ of a rational number $\frac{r}{s}$ is a rational function in $q$ defined by a continued fraction formula. It turns out that $[\frac{r}{s} ]_q$ is a rational function with positive integer coefficients.

The next result follows from Corollary B.4 in \cite{morier2018q} and Theorem \ref{duality_thm}. 

\begin{thm} Recall the notation introduced directly before Theorem \ref{CF_card_quotient_P}. 
Then for any $w$ we have
$$\big[ \text{CF}(w) \big ]_{q} = \frac{\mathbb{P}_w (q) }{ \mathbb{P}^{a_1}_{w}  (q)} = \frac{\mathbb{L}_{w^{*}} (q) }{ \mathbb{L}^{a_1}_{w^{*}}  (q)}$$
and
$$\big[ \text{CF}(w^{*}) \big ]_{q} = \frac{\mathbb{P}_{w^{*}} (q) }{ \mathbb{P}^{b_1}_{w^{*}}  (q)} = \frac{\mathbb{L}_{w} (q) }{ \mathbb{L}^{b_1}_{w}  (q)}.$$
\end{thm}

In \cite{morier2018q} it was conjectured that the numerator and denominator of any $q$-deformed rational is a unimodal polynomial (see Definition \ref{unimodal_def}). 

The unimodality of $\mathbb{L}_w$ is known to be true in some special cases. For instance, $\mathbb{L}_w$ is unimodal if 

\begin{enumerate}[(1)]
    \item $\text{sh}(G_w)$ is straight (trivial), or 
    
    \item $\text{sh}(G_w)$ is zigzag (Fibonacci cubes are unimodal, see \cite{munarini2002rank}), or 
    
    \item $C_{w^{*}}$ is an \textit{up-down poset} (see \cite{gansner1982lattice}, and \cite{gaebler2004alexander} for a relation to Alexander polynomials of $2$-bridge knots), a certain class of posets defined by the division algorithm.
\end{enumerate}

Our first goal in this chapter is to prove that $\mathbb{L}_w (q)$ is unimodal in the case that $G_w$ is built from at most four maximal straight segments (in fact, we prove something stronger than unimodality holds). Our second goal is to say when a poset $\mathbb{P}_w$ or $\mathbb{L}_w$ is symmetric, based on the shape of the underlying snake graph.

\section{Rank Unimodality}

We say the sequence $(r_j)$ has a \textit{plateau} if there exists $p \geq 1$ such that $r_{j_1} = r_{j_1 + 1} = \dots = r_{j_{1} + p}.$ A plateau is \textit{small} if $p=1$. The sequence $(r_j)$ is \textit{trapezoidal} if $(r_j)$ is symmetric and 
$$r_0 < r_1 < \dots < r_j = \dots = r_{n-j} >  \dots > r_n.$$

Say the sequence $(r_j)$ is \textit{weakly trapezoidal} if there exists some $t$ such that 
$$r_0 < r_1 < \dots < r_j = \dots = r_{j+t} >  \dots > r_n$$
such that if $n$ is odd then at least one of the middle two terms $r_{\left \lfloor{\frac{n}{2}}\right \rfloor}$ or $r_{\left \lceil{\frac{n}{2}}\right \rceil}$ is maximal, and if instead $n$ is even then at least the middle term $r_{\frac{n}{2}}$ is maximal. In particular, a weakly trapezoidal sequence is not assumed to be symmetric. A sequence is \textit{almost weakly trapezoidal} if
$$r_0 \leq r_1 < r_2 < \dots < r_j = \dots = r_{j+t} >  \dots > r_{n-2} > r_{n-1} \geq r_n.$$
and the subsequence $(r_j)_{j=1}^{n-1}$ is weakly trapezoidal. Say the sequence $(r_j)$ has \textit{unimodal growth} if the sequence $(\abs{r_{j+1} - r_{j}})$ is bimodal and nonconstant. The polynomial $\rho(q) = \sum_{i=0}^{n} r_i q^i$ is said to have a \textit{plateau} if its sequence of coefficients does, etc.

Consider the hook snake graph $G_w$ of shape $\text{sh}(G_w) = a^{k_1} b^{k_2}.$ Let $k = \text{min} \{ k_1 , k_2\},$ and let $G_{w_0}$ be the subsnake graph of $G_w$ of shape $\text{sh}(G_{w_0}) = a^{k-1} b^{k-1}$. Let $n_0 = 2k-1$ be the (odd) number of tiles of $G_{w_0}$.

\begin{prop} \label{hook_rank_function}
    Consider the hook snake graph $G_w$ of shape $\text{sh}(G_w) = a^{k_1 - 1} b^{k_2 - 1}.$ Then for $k_1 \leq k_2$ we have
    $$\mathbb{L}_w (q) = 1 + q + 2q^2 + \dots + k_1 q^{k_1} + \dots + k_1 q^{k_2} + \dots + 2q^{k_1 + k_2 -2} + q^{k_1 + k_2 - 1}.$$
\end{prop}

\begin{proof}
From the recursion given above we can compute $\mathbb{L}_w (q) = 1 + q[k_1]_q [k_2]_q.$ Now the claim follows from Proposition 1.5 (1) in \cite{cheng2016strict}.
\end{proof}

\begin{cor}
    Consider the hook snake graph $G_w$ of shape $\text{sh}(G_w) = a^{k_1 - 1} b^{k_2 - 1}.$ Let $k = \text{min}(k_1 , k_2)$. 

    \begin{enumerate}[(a)]

    \item The first and last coefficients of $\mathbb{L}_w (q)$ are both equal to $1$.
    
    \item The maximum value of the coefficients of $\mathbb{L}_w (q)$ is equal to $k.$
    
    \item The number of times the maximum value occurs is equal to $n-n_0+1 = n-2k+2,$ in degrees $k , k+1 , \dots , n-k+1.$
    
    \item The polynomial $\mathbb{L}_w (q)$ has at least one small plateau, consisting of the first two coefficients. If $\mathbb{L}_w (q)$ has a second plateau, then each entry of this plateau is equal to the maximum value. There are no other plateaus.   
    
    \item The polynomial $\mathbb{L}_w (q)$ is almost weakly trapezoidal. In particular, $\mathbb{L}_w (q)$ is unimodal.
\end{enumerate}
\end{cor}

\begin{proof}
    All parts (a)-(e) follow directly from Proposition \ref{hook_rank_function}.
\end{proof}

\begin{ex}
    Consider $G_w$ with $\text{sh}(G_w) = a^2 b^5,$ shown in the leftmost illustration in Figure below. The subsnake graph $G_0$ is shown shaded inside $G$. The middle figure is isomorphic to $\mathbb{L}_{w_0},$ and the rightmost lattice shown is isomorphic to $\mathbb{L}_w.$ The maximal coefficient in the polynomial $\mathbb{L}_w (q)$ is equal to $3,$ and this maximal coefficient occurs $4$ times. Here, $n=8$, $n_0 = 5$, $k= k_1 = 3,$ and $k_2 = 6$. 
\end{ex}

\begin {figure}[h!]
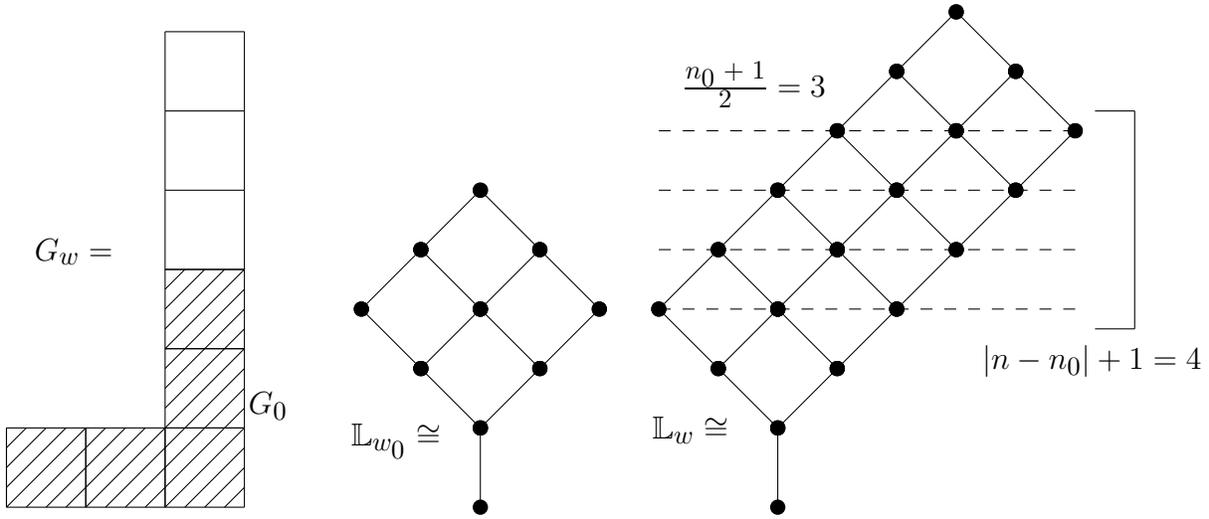

    \centering
    \caption{The snake graphs $G_{w}$ and $G_{0}$, and the associated posets $\mathbb{L}_{w_0}$ and $\mathbb{L}_w$}
    \label{fig:hook_rank_pic}

  
\end {figure}

\begin {figure}[h!]
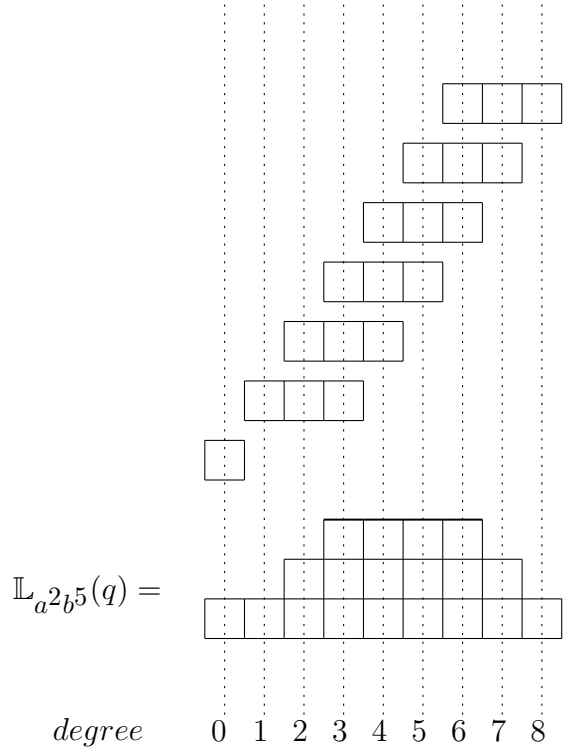

    \centering
    \caption{The rank function $\mathbb{L}_{a^2 b^5} (q)$ via block stacking}
    \label{fig:block_stacking}
    %

\end {figure}

It is useful to visualize any rank function from a snake graph in terms of stacking square tiles. Each coefficient $r_j$ of such a rank function is represented by $r_j$ tiles stacked vertically. Figure \ref{fig:block_stacking} shows how the rank function $\mathbb{L}_{a^2 b^5}(q)$ from the previous example is built, according to the expression $\mathbb{L}_{a^2 b^5} (q) = 1 + q[3]_q + q^2 [3]_q + \dots + q^{6}[3]_q.$

Let $n,m \geq 1,$ and suppose $n = k_1 + k_2 - 1$. Let $k = \text{min}(k_1 , k_2).$ Define the hook subsnake graphs $\mathcal{H} = \mathcal{H}_{k_1 , k_2 ,m}$ of $G_w$ to be of shape 

\[
  \text{sh}(\mathcal{H})  =
  \begin{cases}
                                   a^{k_{1}-1} b^{\text{min}(m-1,k_2 - 1)} & \text{if $k = k_1$} \\
                                    a^{\text{min}({m-1 , k_1 - 1})} b^{k_2 - 1} & \text{if $k = k_2$}
  \end{cases}
\]

Let the coefficients of the rank function of lattice paths on $\mathcal{H}_{k_1,k_2,,m}$ be called $\mathcal{H}_j(k_1,k_2,m).$ Note that when $m$ is equal to or exceeds the length of the longer of the two straight segments that $G_w$ is built from, then $\text{sh}(\mathcal{H}) = \text{sh}(G_w)$. In other words, writing $\mathbb{L}_w (q) = \sum_{j=0}^{n} r_j q^j$ we have that $\mathcal{H}_j (k_1 , k_2 , m) = r_j$ for $m$ sufficiently large.

Define $\chi = \chi(k_1 , k_2 , m)$ by 

\[
  \chi  =
  \begin{cases}
                                  k+m-1 & \text{if $m \leq n - 2k + 2$} \\
                                    \lfloor  \frac{n+m}{2}  \rfloor & \text{if $n-2k+2 < m < n$} \\
                                    n & \text{if $n < m$}
  \end{cases}
\]

and the coefficients

\[
  \widetilde{\mathcal{H}}_j (k_1 , k_2 , m)  =
  \begin{cases}
                                  1, & \text{if $n-2k+2 < m < n$ and $2 \mid (n+m)$ and $j = \frac{n+m}{2}$} \\
                                    \mathcal{H}_j(k_1 , k_2 , m), & \text{else}. 
  \end{cases}
\]

\begin{prop} \label{m_times_hook_prop}
    Consider the rank function $\mathbb{L}_w (q) = \sum_{j=0}^{n} r_j q^j$ of the hook snake graph $G_w$ of shape $\text{sh}(G_w) = a^{k_1 - 1} b^{k_2 - 1}$ built from two maximal straight segments of length $k_1 \geq 2$ and $k_2 \geq 2$. Let the number of tiles of $G_w$ be $n= k_1 + k_2 - 1$. Let $k = \text{min}(k_1 , k_2)$ and $m \geq 1$. Then we have 
    $$[m]_q \mathbb{L}_w (q) = [m]_q + \sum_{j=1}^{\chi} \widetilde{\mathcal{H}}_j(k_1,k_2,m) [n+m - 2j + 1]_q q^{j}.$$

In any case, $[m]_q \mathbb{L}_w (q)$ is weakly trapezoidal, and has unimodal growth for $m \geq 2$. 

Furthermore, 

\begin{enumerate}[(a)]

    \item If $m < n-2k+2$, then $[m]_q \mathbb{L}_w (q)$ has a unique plateau consisting of maximal values, which is of length $(n-2k+1)-(m-2)$ and is situated in degrees $k + (m-1) , k+m , \dots , n-k+1$

    \item If $m = n-2k+2$, then $[m]_q \mathbb{L}_w (q)$ has a maximum at degree $n-k+1$, and no plateau. 
    
    \item If $n-2k + 2 < m < n$ and $2 \nmid n+m$, then $[m]_q \mathbb{L}_w (q)$ has a unique (small) plateau, situated in degrees $\left \lfloor{\frac{n+m}{2}}\right \rfloor$ and $\left \lceil{\frac{n+m}{2}}\right \rceil$.
    
       \item If $n-2k+2 < m < n$ and $2 \mid n+m$, then $[m]_q \mathbb{L}_w (q)$ has a maximum at degree $\frac{n+m}{2}.$

    \item If $m = n$, then $[m]_q \mathbb{L}_w (q)$ has a unique (small) plateau, situated in degrees $n-1$ and $n$.

    \item If $m = n+1$, then $[m]_q \mathbb{L}_w (q)$ has a maximum at degree $n$, and no plateau.

    \item If $m > n+1$, then $[m]_q \mathbb{L}_w (q)$ has a unique plateau, which is of length $m-n$ and situated in degrees $n,n+1 , \dots , m-1.$

\end{enumerate}

\end{prop}

\begin{proof}
    Write the rank function $\mathbb{L}_w (q)$ as 
$$\mathbb{L}_w (q) = [1]_q + q[n]_q + q^{2} [n-2]_q + \dots + q^{k} [n-2k+2]_q,$$
so that the product $[m]_q \mathbb{L}_w (q)$ is 
\begin{equation} \label{eq1}
[m]_q \mathbb{L}_w (q) = [m]_q + q[m]_q [n]_q + q^{2} [m]_q [n-2]_q + \dots + q^{k} [m]_q [n-2k + 2]_q.
\end{equation}

From \cite{cheng2016strict} we have

\begin{equation} \label{eq2}
[m]_q [n]_q = \sum_{j=1}^{\text{min}(m,n)} [m+n-2j+1]_q q^{j-1}.
\end{equation}

Explicitly,

\begin{itemize}

\item If $m < n$, then 
$$[m]_q [n]_q = [n+(m-1)]_q + q[n+(m-3)]_q + \dots + q^{m-1} [n-(m-1)]_q.$$
\item  If $m=n,$ then  
$$[m]_q [n]_q = [n]_{q}^{2} =  [2n-1]_q + q[2n-3]_q + \dots + q^{n-1} [1]_q.$$

\item If $m > n$, then 
$$[m]_q [n]_q = [m+(n-1)]_q + q[m+(n-3)]_q + \dots + q^{n-1} [m-(n-1)]_q.$$

\end{itemize}

We consider seven cases, based on the value of $m$ relative $n-2k+2$ and $n$, and the parity of the sum $n+m$. In each case we use Equation \ref{eq2} above to expand each product of $q$-numbers in Equation \ref{eq1} and collect like terms to obtain the formula. This is shown explicitly in the first case only, as the rest are similar.

\begin{enumerate}[(a)]

    \item \underline{$m < n - 2k+ 2:$} Here, the statement is 
$$[m]_q \mathbb{L}_w (q) = [m]_q + \sum_{j=1}^{k+m-1} \mathcal{H}_j (k_1,k_2,m) [n+m -2j + 1]_q q^{j}.$$

 Using (2) to expand each term $q^j [m]_q [n - 2j]_q$ in (1) and collecting like terms gives the desired expression for $[m]_q \mathbb{L}_w (q)$: 

\begin{equation} \label{eq2}
\begin{split}
[m]_q \mathbb{L}_w (q) & = [m]_q \Big( 1 + q \sum_{j=1}^{k} [(k_1 + k_2) -2j + 1]_q q^{j-1} \Big)  \\
 & = [m]_q \Big( 1 + \sum_{j=1}^{k} [(k_1 + k_2 - 1) -2j + 2]_q q^{j} \Big) \\
 & = [m]_q \Big( 1 + \sum_{j=1}^{k} [n -2j + 2]_q q^{j} \Big) \\
  & = [m]_q + \sum_{j=1}^{k}  [m]_q [n -2j + 2]_q q^{j} \\
  & = [m]_q + \sum_{j=1}^{k} \sum_{i=1}^{m} [n -2(i+j-1) + 1]_q q^{i+j-1} \\
  & =   [m]_q + \sum_{j=1}^{k+m-1} \mathcal{H}_j (k_1,k_2,m) [n+m -2j + 1]_q q^{j}
\end{split}
\end{equation}

    Now let 
    $$\rho (q) \doteq \sum_{j=1}^{k+m-1} \mathcal{H}_j(k_1,k_2,m)  [n+m - 2j + 1]_q q^{j}.$$ Note that when we increment $j$ in this sum, the power of $q$ goes up by one, while the term inside the bracket decreases by two. This fact, along with the fact that the coefficients $\mathcal{H}_j (k_1,k_2,m)$ form a weakly trapezoidal sequence, shows that $\rho (q)$ is weakly trapezoidal with unimodal growth.

    In particular, the degrees of the terms whose coefficients make up the unique plateau of $\rho (q)$ are exactly the powers of $q$ occurring in the last term $q^{k+m-1} [n-m-2k+3].$ Thus we read off that $\rho(q)$ has a plateau of maximal coefficients of length $(n-2k+1)-(m-2)$ situated in degrees $k + (m-1) , k+m , \dots , n-k+1$.

    Since $k > 1$ we have $m-1 < m < n-2k + 2 < n-k+1$, so that adding the chain $[m]_q$ to $\rho (q)$ does not disturb the plateau of maximal coefficients of $\rho(q)$. In other words, $[m]_q + \rho (q)$ has the same plateau of maximal coefficients as $\rho(q)$. Since the coefficients $\mathcal{H}_{j}(k_1,k_2,m)$ form a weakly trapezoidal sequence, no other plateaus are created by adding the chain $[m]_q$ to $\rho(q)$. Hence $[m]_q \mathbb{L}_w (q)$ is weakly trapezoidal with unimodal growth, as claimed.

    \item \underline{$m = n - 2k + 2$:}
    
    In this case the formula can be computed to be 
    $$[m]_q \mathbb{L}_w (q) = [n-2k+2]_q + \sum_{j=1}^{n-k+1} \mathcal{H}_j (k_1,k_2,m) [2n-2k -2j + 3]_q q^{j},$$
    and the last term in the sum 
    $\rho(q) \doteq \sum_{j=1}^{n-k+1} \mathcal{H}_j (k_1,k_2,m) [2n-2k -2j + 3]_q q^{j}$ is $[1]_q q^{n-k+1}.$ Thus, $\rho (q)$ has a maximum at degree $n-k+1$. Just as before,  adding the chain $[n-2k+2]_q$ to $\rho (q)$ does not disturb the unique plateau of maximal coefficients of $\rho(q).$  Since the coefficients $\mathcal{H}_j (k_1,k_2,m)$ are weakly trapezoidal, no other plateaus are created by adding $[m]_q$ to $\rho(q)$. Thus the claim holds in this subcase as well.

    \item \underline{$n-2k+2 < m < n$ and $2 \nmid n+m:$} If $n$ and $m$ have opposite parity then there exists some $c$ such that $$n - 2c < m = n-2c+1 < n - 2c+2.$$ Using Equation \ref{eq2}, we see that 
$$[m]_q \mathbb{L}_w (q) = [m]_q + \sum_{j=1}^{n-c} \mathcal{H}_j(k_1,k_2,m) [n+m - 2j + 1]_q q^{j}.$$
This matches the form given in the statement, since $2 \nmid n+m$ implies $\lfloor \frac{n+m}{2} \rfloor = \frac{n+m-1}{2} = \frac{n+(n-2c+1) -1}{2} = n-c.$ Again by using a degree argument we see that $[m]_q \mathbb{L}_w (q)$ is weakly trapezoidal with unimodal growth.
    
    \item \underline{$n-2k+2 < m < n$ and $2 \mid n+m:$} If $2 \mid m+n$ then there exists some $c$ such that $m = n - 2c.$ As in the previous case we use \ref{eq2} to write
    $$[m]_q \mathbb{L}_w (q) = [m]_q + \sum_{j=1}^{n-c} \mathcal{H}_j(k_1,k_2,m) [n+m - 2j + 1]_q q^{j},$$
    except that now $2 \mid n+m$ so that $\frac{n+m}{2} = n-c.$ That the maximum of $\mathbb{L}_w (q)$ occurs in degree $\frac{n+m}{2}$ can be seen by computing the last term $[1]_q q^{\frac{n+m}{2}}$ in the sum.
    
    \item \underline{$m = n:$} Using Equation \ref{eq2} gives 

$$[m]_q \mathbb{L}_w (q) = [n]_q + \sum_{j=1}^{n} r_j [2n - 2j + 1]_q q^{j}.$$

The last term in the above sum $\sum_{j=1}^{n} r_j  [2n - 2j + 1]_q q^{j}$ is $[1]_q q^n,$ which shows that $\sum_{j=1}^{n} r_j  [2n - 2j + 1]_q q^{j}$ has a maximum in degree $n$. That the last two terms are  $[1]_q q^n$ and $2[3]_q q^{n-1}$ implies that whatever the maximum coefficient is, the previous coefficient is one less. Now adding $[n]_q = [n-1]_q + q^{n-1}$ to the sum $\sum_{j=1}^{n} r_j q^{j} [2n - 2j + 1]_q$ and again considering the term $[3]_q q^{n-1}$ shows that the small plateau in question exists, in degrees $n-1$ and $n$ as claimed. That it is unique follows easily, as does the fact that the polynomial $[m]_q \mathbb{L}_w (q)$ is weakly trapezoidal with unimodal growth.

    \item \underline{$m = n + 1:$} In this case we have 
$$[m]_q \mathbb{L}_w (q) = [n+1]_q + \sum_{j=1}^{n} r_j [2n - 2j + 2]_q q^{j}.$$
To see that $[m]_q \mathbb{L}_w (q)$ has a unique maximum coefficient, consider the last term $q^{n} [2]_q$ from the sum $\rho(q) \doteq \sum_{j=1}^{n} r_j q^{j} [2n - 2j + 2]_q.$ It follows that $\rho(q)$ has a unique maximal plateau in degrees $n$ and $n+1$. Now adding the chain $[n+1]_q = [n]_q + q^n$ to $\rho(q)$ shows that $[n+1]_q \mathbb{L}_w (q)$ has a unique maximum coefficient, in degree $n$. That $[m]_q \mathbb{L}_w (q)$ is weakly trapezoidal with unimodal growth follows immediately. 
    
    \item  The formula here is 
$$[m]_q \mathbb{L}_w (q) = [m]_q + \sum_{j=1}^{n} r_j [n + m - 2j + 1]_q q^{j}.$$
Again, the last coefficient $[m-n+1]_q q^{n}$ tells us the behavior of the plateau of the sum above, and the claim follows.

\end{enumerate}

\end{proof}

\begin{ex}
In Figure \ref{fig:times_m_blocks} below, we show the effect of multiplying the rank function $\mathbb{L}_{a^2 b^5} (q)$ by $[m]_q$ for $1 \leq m \leq 12.$ We have indicated each plateau of maximal coefficients of $[m]_q \mathbb{L}_{a^2 b^5} (q)$ for $1 \leq m \leq 9 = n+1$ with a bold line.

\begin {figure}[h!]
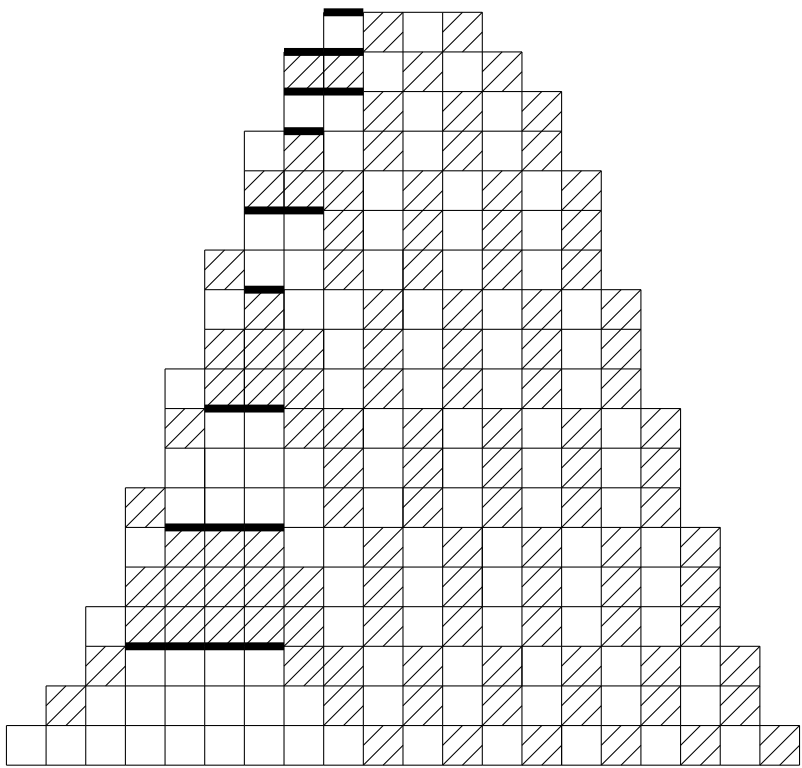

    \centering
    \caption{The rank function $[m]_q \mathbb{L}_{a^2 b^5} (q)$ for $1 \leq m \leq 12$}
    \label{fig:times_m_blocks}



\end {figure}

\end{ex}

\begin{thm}
    Consider the snake graph $G_w$ and the rank generating function $\mathbb{L}_w (q)$. 

\begin{enumerate}[(a)] 

\item Suppose the snake graph $G_w$ with $n=k_1 + k_2 + k_3 -1$ tiles has the shape $\text{sh}(G_w) = a^{k_1-1} b^{k_2} a^{k_3-1},$ where $k_1 \geq 2$, $k_2 \geq 1$, and $k_3 \geq 2$. Then the polynomial $\mathbb{L}_w (q)$ is weakly trapezoidal with unimodal growth.  In particular, $\mathbb{L}_w (q)$ is unimodal. 

\item Suppose the snake graph $G_w$ with $n=k_1+k_2+k_3+k_4 - 1$ tiles has the shape $\text{sh}(G_w) = a^{k_1-1} b^{k_2} a^{k_3} b^{k_4 - 1}.$ Then the polynomial $\mathbb{L}_w (q) = \sum_{j=1}^{n}$ is weakly trapezoidal with unimodal growth. In particular, $\mathbb{L}_w (q)$ is unimodal. 

\end{enumerate}
\end{thm}

\begin{proof}
    
    \begin{enumerate}[(a)]
    
        \item By Theorem \ref{hook_expansion} we can write 
        $$\mathbb{L}_{w} (q) = [k_3]_q \Big( 1+q[k_1]_q [k_2]_q \Big) + q^{k_2 + k_3} [k_1]_q,$$ By Proposition \ref{m_times_hook_prop} we have
        $$\mathbb{L}_{w} (q) = \rho(q) + q^{k_2 + k_3} [k_1]_q$$
        where $\rho(q)$ is the weakly trapezoidal polynomial 
        $$\rho(q) \doteq [k_3]_q + \sum_{j=1}^{\chi} \widetilde{\mathcal{H}}_j(k_1,k_2,k_3) [n - 2j + 1]_q q^{j}$$
        with unimodal growth.
        
        Note that the highest power occurring in $q^{k_2 + k_3} [k_1]_q$ is $q^{k_1 + k_2 + k_3 - 1}$, which is one higher than the degree of $\rho(q)$. This remark, along with the fact that the coefficients $\mathcal{H}_j (k_1 , k_2 , k_3)$, have unimodal growth, implies that adding the chain $q^{k_2 + k_3} [k_1]_q$ to $\rho(q)$ does not create any new plateaus. In particular, $\mathbb{L}_w$ is unimodal.
        
        \item Let $G_{w^{\prime}}$ be the subsnake graph of $G_w$ obtained by deleting the first $k_1$ tiles from $G_w$, and $G_{w^{\prime \prime}}$ the subsnake graph of $G_{w^{\prime}}$ obtained by deleting the first $k_1 + k_2 - 1$ tiles from $G_w$. Set $\rho^{\prime} (q) = \mathbb{L}_{w^{\prime}} (q)$ and $\rho^{\prime \prime} (q) = \mathbb{L}_{w^{\prime \prime}} (q).$ By the recurrence relation above we can write $$\mathbb{L}_w (q) = \bigg( [k_1]_q q \bigg) \rho^{\prime}(q) + \rho^{\prime \prime} (q).$$

        By part (a), $\rho^{\prime} (q)$ is weakly trapezoidal with unimodal growth. It is clear that multiplying $\rho^{\prime} (q)$ by $q [k_1]_q$ preserves these properties.

        Note that 
        $$\rho^{\prime}(q) =  [k_2]_q \bigg( 1 + q [k_3]_q [k_4]_q \bigg) + q^{k_2 + k_3} [k_4]_q.$$ In particular $\rho^{\prime}(q)$ is built from $\rho^{\prime \prime} (q) = 1 + q[k_3]_q [k_4]_q$ by first multiplying by $[k_2]_q$ and then adding to the result the chain $q^{k_2 + k_3} [k_4]_q.$ Since $k_2 \geq 2$, the degree of the first entry of the plateau of $\rho^{\prime}(q)$ is weakly larger than the degree of the first entry of the plateau of $\rho^{\prime \prime}(q).$ Multiplying $\rho^{\prime}(q)$ by $[k_1]_q q$ can only shift the start of the plateau of $\rho^{\prime}(q)$ further to the right, implying that the start of the plateau of $\big( [k_1]_q q \big) \rho^{\prime}(q) $ is to the right of the start of the plateau of $\rho^{\prime \prime}(q).$

        If the two plateaus in question overlap or are adjacent, then we are done. Otherwise, the lowest degree occurring in the plateau of $\big( [k_1]_q q \big) \rho^{\prime}(q) $ is strictly larger than $\text{max}(k_3 , k_4)$. In this case, the unimodal growth of the coefficients of $\big( [k_1]_q q \big) \rho^{\prime}(q) $ ensures that the addition of $\rho^{\prime \prime}(q) $ (whose consecutive coefficients are at most one apart) preserves the two properties in question. Hence, $\mathbb{L}_w$ is unimodal in this case.

    \end{enumerate}

\end{proof}

\section{Rank Symmetry}

\begin{defn}
    Fix a word $w = w_1 w_2 \dots w_{n-1}$ of length $l(w) = n-1$. Define the three words $$w^{T} = w_{1}^{*} w_{2}^{*} \dots w_{n-1}^{*},$$ $$w^{\circ} = w_{n-1} w_{n-2} \dots w_{1},$$ and $$\overline{w} = w_{n-1}^{*} w_{n-2}^{*} \dots w_{1}^{*}.$$ We say that $w$ is \textit{symmetric} if $w^{T} = w,$ and \textit{self-conjugate} if $\overline{w} = w.$ We say that $G_w$ is \textit{symmetric} if $\text{sh}(G_w)$ is symmetric, and \textit{self-conjugate} if $\text{sh}(G_w)$ is self-conjugate. 
\end{defn}

\begin{ex}
    The word $w_1 = \text{sh}(G_{w_1})^{*} = aab$ from Example \ref{three_snakes_ex} is neither symmetric nor self-conjugate, while $w_2 = \text{sh}(G_{w_2})^{*} = aabb$ is self-conjugate and $w_3 = \text{sh}(G_{w_3})^{*} = baaab$ is symmetric. The word  $\text{sh}(G_{w_{1}}) = w_{1}^{*} = baa$ is neither symmetric nor self-conjugate, while the words $\text{sh}(G_{w_{2}}) = w_{2}^{*} = baab$ and $\text{sh}(G_{w_{3}}) = w_{3}^{*} = aabaa$ are both symmetric.
\end{ex}

\begin{prop} \label{sym}
    Fix the word $w$ of length $l(w) = n-1$. Recall the internal edges $e_1 , e_2 , \dots e_{n-1}$ of the snake graph $G_w$ (see the discussion before Definition \ref{straight_zigzag_SG_def} above).

\begin{enumerate}[(a)]
    \item Suppose $l(w)$ is odd. Then the word $w = \text{sh}(G_{w})^{*}$ is symmetric if and only if the word $w^{*} = \text{sh}(G_{w})$ is symmetric. 
    
    \item If $l(w)$ is odd and $w$ is symmetric, then $G_w$ has $180^{\circ}$ rotational symmetry about the midpoint of its middle internal edge, and the same is true for $G_{w}^{*}.$
    
    \item Suppose $l(w)$ is even. Then the word $w = \text{sh}(G_{w})^{*}$ is symmetric if and only if the word $w^{*} = \text{sh}(G_{w})$ is self-conjugate, and $w = \text{sh}(G_{w})^{*}$ is self-conjugate if and only if $w^{*} = \text{sh}(G_{w})$ is symmetric.

    \item If $l(w)$ is even and $w$ is symmetric, then $G_w$ is symmetric about the diagonal of its center tile $T_{\frac{n+1}{2}},$ and $G_{w}^{*}$ has $180^{\circ}$ rotational symmetry about the center of its middle tile $T_{\frac{n+1}{2}}$.

    \item The snake graph $G_w$ has $180^{\circ}$ rotational symmetry if and only if $C_{w^{*}}$ has $180^{\circ}$ rotational symmetry. In this case, the graded poset $\mathbb{L}_w$ is order-theoretically self-dual.  
\end{enumerate}
\end{prop}

\begin{proof}

    \begin{enumerate}[(a)]
        \item If $w$ is symmetric of odd length, then setting $m = \frac{l(w)+1}{2}$ we can write $$w = w_1 w_2 \dots w_{m-1} w_{m} w_{m-1} \dots w_2 w_1.$$ If $m$ is even, then taking the dual of $w$ gives $$w^{*} = w_{1}^{*} w_{2} \dots w_{m-1}^{*} w_{m} w_{m-1}^{*} \dots w_2 w_{1}^{*},$$ which is symmetric. Similarly, if $m$ is odd then taking the dual gives $$w^{*} = w_{1}^{*} w_{2} \dots w_{m-1} w_{m}^{*} w_{m-1} \dots w_2 w_{1}^{*},$$ which is also symmetric. Taking dual words now gives (a).
        
        \item Follows directly from (a).

        \item If $l(w)$ is even and symmetric, then for $m = \frac{l(w)}{2}$ we can write 
        $$w = w_1 w_2 \dots w_{m-1} w_{m} w_{m} w_{m-1} \dots w_{2} w_{1}.$$ Suppose $m$ is even. Then taking the dual gives
        $$w^{*} = w_{1}^{*} w_2 \dots w_{m-1}^{*} w_{m} w_{m}^{*} w_{m-1} \dots w_{2}^{*} w_{1},$$ which is self-conjugate. Similarly, if $m$ is odd, then taking the dual of $w$ gives 
        $$w_{1}^{*} w_2 \dots w_{m-1} w_{m}^{*} w_{m} w_{m-1}^{*} \dots w_{2}^{*} w_{1},$$
        which is again self-conjugate.
        
        \item Follows directly from (c).
        
        \item By construction, $G_w$ has $180^{\circ}$ symmetry if and only if $C_{w}^{*}$ does. Let $G_{w}^{\circ}$ be $G_{w}$ rotated by $180^{\circ}$. The result now follows from noting that the order-theoretic dual of $\mathbb{L}_w$ is isomorphic to the poset of lattice paths on $G_{w}^{\circ}.$ 
    \end{enumerate}
\end{proof}

\begin{cor} \label{sym_cor}
    Consider the words $w$ and $w^{*}$ of length $n-1$, the associated snake graphs $G_w$ and $G_{w}^{*}$ with $n$ tiles, and the distributive lattices $D_w \cong \mathbb{P}_w \cong \mathbb{L}_{w^{*}}$ and $D_{w^{*}} \cong \mathbb{L}_w \cong \mathbb{P}_{w^{*}}$ of rank $n$. 

\begin{enumerate}[(a)] 
    
    \item If $\text{sh}(G_w)$ is symmetric, then the poset $\mathbb{L}_w$ is symmetric.  
    
    \item If either $l(w)$ is odd and $\text{sh}(G_w)$ is symmetric, or $l(w)$ is even and $\text{sh}(G_w)$ is self-conjugate, then the poset $\mathbb{P}_w$ is symmetric
    
\end{enumerate}

\end{cor}

\begin{proof}
    Follows directly from Proposition \ref{sym}.
\end{proof}

\begin{ex}
    
    \begin{enumerate}[(a)]
        \item The lattice path expansion on a zigzag snake graph with an even number of tiles is isomorphic to a Fibonacci cube of even order, and so is symmetric (this was shown in \cite{munarini2002rank}). Trivially, the perfect matching expansion poset on this snake graph is symmetric, since it is a chain.
        
        \item Consider the snake graph $G_w$ with the word $\text{CF}(w) = [a_1 , a_2 , \dots , a_k].$  Following \cite{ccanakcci2020snake}, the \textit{palindromification} of $G_w$ is the snake graph $G_{\leftrightarrow}$ with associated continued fraction $[a_n , a_{n-1} , \dots , a_2 , a_1 , a_2 , \dots , a_{n-1} , a_{n}].$ Every palindromification $G_{\leftrightarrow}$ has $180^{\circ}$ symmetry about its center tile (see Theorem A in \cite{ccanakcci2020snake}). Let $\mathbb{L}_{\leftrightarrow}$ be the poset of lattice paths on $G_{\leftrightarrow}.$ Then by Corollary \ref{sym_cor}, we know that $\mathbb{L}_{\leftrightarrow}$ are symmetric.  Since $l(w)$ is even, the poset $\mathbb{P}_{\leftrightarrow}$ of perfect matchings on the palindromification is not symmetric.
        
        \item A \textit{Markov snake graph} is a snake graph which can be built from the Christoffel path in a $p$ by $q$ grid, where $p$ and $q$ are relatively prime positive integers (see \cite{propp2005combinatorics}, \cite{ccanakcci2020snake}, and \cite{rabideau2018continued}). Figure \ref{fig:markov} shows one example of a Markov snake graph. The number of perfect matchings on any Markov snake graph is equal to a Markov number. Every Markov snake graph is a palindromification (see \cite{ccanakcci2020snake}), so the conclusions in (b) hold here as well.
        
        \begin {figure}[h!]
    \centering
    \caption{A Markov snake graph}
    \label{fig:markov}

\tikzset{every picture/.style={line width=0.75pt}} 

\begin{tikzpicture}[x=0.75pt,y=0.75pt,yscale=-1,xscale=1]

\draw    (540,260) -- (540,330) ;
\draw    (540,330) -- (610,330) ;
\draw    (540,260) -- (610,260) ;
\draw    (610,260) -- (610,330) ;
\draw    (610,330) -- (680,330) ;
\draw    (610,260) -- (680,260) ;
\draw    (680,260) -- (680,330) ;
\draw    (680,330) -- (750,330) ;
\draw    (680,260) -- (750,260) ;
\draw    (750,260) -- (750,330) ;
\draw    (540,190) -- (540,260) ;
\draw    (540,190) -- (610,190) ;
\draw    (610,190) -- (610,260) ;
\draw    (610,190) -- (680,190) ;
\draw    (680,190) -- (680,260) ;
\draw    (680,190) -- (750,190) ;
\draw    (750,190) -- (750,260) ;
\draw [color={rgb, 255:red, 208; green, 2; blue, 27 }  ,draw opacity=1 ]   (540,330) -- (750,190) ;
\draw [line width=1.5]    (575,295) -- (610,295) ;
\draw [line width=1.5]    (575,295) -- (575,330) ;
\draw [line width=1.5]    (610,295) -- (610,330) ;
\draw [line width=1.5]    (575,330) -- (610,330) ;
\draw [line width=1.5]    (610,330) -- (645,330) ;
\draw [line width=1.5]    (645,295) -- (645,330) ;
\draw [line width=1.5]    (610,295) -- (645,295) ;
\draw [line width=1.5]    (680,295) -- (680,330) ;
\draw [line width=1.5]    (645,295) -- (680,295) ;
\draw [line width=1.5]    (645,330) -- (680,330) ;
\draw [line width=1.5]    (645,260) -- (680,260) ;
\draw [line width=1.5]    (680,260) -- (680,295) ;
\draw [line width=1.5]    (645,260) -- (645,295) ;
\draw [line width=1.5]    (680,260) -- (715,260) ;
\draw [line width=1.5]    (645,225) -- (645,260) ;
\draw [line width=1.5]    (645,225) -- (680,225) ;
\draw [line width=1.5]    (680,225) -- (680,260) ;
\draw [line width=1.5]    (715,225) -- (715,260) ;
\draw [line width=1.5]    (680,225) -- (715,225) ;
\draw [line width=1.5]    (750,225) -- (750,260) ;
\draw [line width=1.5]    (715,225) -- (750,225) ;
\draw [line width=1.5]    (715,260) -- (750,260) ;

\end{tikzpicture}

\end {figure}
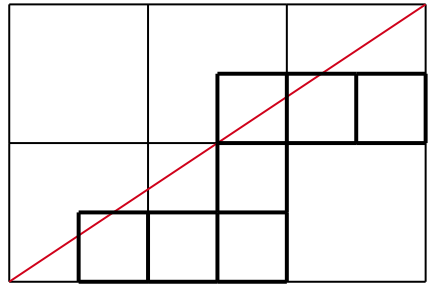

        \item It is well-known that the infinite continued fraction expansion of $\sqrt{2}$ is equal to $[1,2,2,2, \dots].$ Thus, we may write $1 + \sqrt{2} = [2,2,2,2, \dots].$ The latter infinite fraction is called the \textit{silver mean}. Consider any of the continued fractions obtained by truncating the continued fraction expansion of the silver mean after an odd number (larger than $1$) of $2$'s. An example of a snake graph associated to such a continued fraction is shown below.
        
        \begin {figure}[h!]
    \centering
    \caption{The snake graph with continued fraction $[2,2,2,2,2]$}
    \label{fig:pell}

\tikzset{every picture/.style={line width=0.75pt}} 

\begin{tikzpicture}[x=0.75pt,y=0.75pt,yscale=-1,xscale=1]

\draw    (260,390) -- (330,390) ;
\draw [line width=1.5]    (155,425) -- (190,425) ;
\draw [line width=1.5]    (155,425) -- (155,460) ;
\draw [line width=1.5]    (190,425) -- (190,460) ;
\draw [line width=1.5]    (155,460) -- (190,460) ;
\draw [line width=1.5]    (190,460) -- (225,460) ;
\draw [line width=1.5]    (225,425) -- (225,460) ;
\draw [line width=1.5]    (190,425) -- (225,425) ;
\draw [line width=1.5]    (260,425) -- (260,460) ;
\draw [line width=1.5]    (225,425) -- (260,425) ;
\draw [line width=1.5]    (225,390) -- (260,390) ;
\draw [line width=1.5]    (260,390) -- (260,425) ;
\draw [line width=1.5]    (225,390) -- (225,425) ;
\draw [line width=1.5]    (260,390) -- (295,390) ;
\draw [line width=1.5]    (225,460) -- (260,460) ;
\draw [line width=1.5]    (225,355) -- (225,390) ;
\draw [line width=1.5]    (225,355) -- (260,355) ;
\draw [line width=1.5]    (260,355) -- (260,390) ;
\draw [line width=1.5]    (295,355) -- (295,390) ;
\draw [line width=1.5]    (260,355) -- (295,355) ;
\draw [line width=1.5]    (330,355) -- (330,390) ;
\draw [line width=1.5]    (295,355) -- (330,355) ;
\draw [line width=1.5]    (295,390) -- (330,390) ;
\draw [line width=1.5]    (295,320) -- (295,355) ;
\draw [line width=1.5]    (330,320) -- (330,355) ;
\draw [line width=1.5]    (295,320) -- (330,320) ;
\draw [line width=1.5]    (295,285) -- (295,320) ;
\draw [line width=1.5]    (330,285) -- (330,320) ;
\draw [line width=1.5]    (295,285) -- (330,285) ;

\end{tikzpicture}

\end {figure}
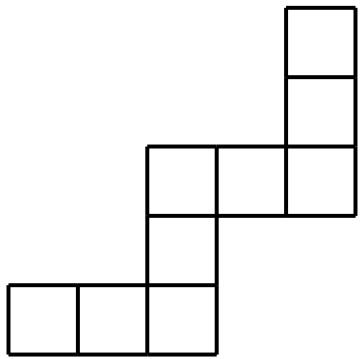

        By the Corollary \ref{sym_cor}, the perfect matching poset on the associated snake graph is symmetric, while the lattice path expansion poset on this snake graph is not symmetric.

    \end{enumerate}
    
\end{ex}

\chapter{Expansion Posets as Groupoid Orbits}

In this short chapter, we give an interpretation of the support of the cluster variable $x_w$ as an orbit of a groupoid. This implies that the support of any two distinct cluster variables written with respect to the same initial seed are disjoint, so that in particular any $x_w$ can be completely reconstructed from any one of its monomials.

Recall the definition of a groupoid, given as Definition  \ref{definition_groupoid} above. Consider the set of Laurent monomials from the extended cluster $(x_1 , x_2 , \dots , x_n , x_{n+1} , \dots , x_{2n+3}.)$ For each $k$ with $1 \leq k \leq n$, define the element 
\begin{equation}
    \hat{y}_k = \frac{\prod_{i \rightarrow k} x_i}{\prod_{k \rightarrow j} x_j}.
\end{equation}
(see \cite{fomin2007cluster}).

Define a groupoid $\mathcal{F}$ whose objects are the Laurent monomials from the extended cluster above. There is a morphism $x_{M} \longrightarrow x_{M^{\prime}}$ between two Laurent monomials  $x_M$ and $x_{M^{\prime}}$ if $x_{M^{\prime}} = \hat{y} x_M$ or $x_{M^{\prime}} = \hat{y}^{-1} x_M$ for some $\hat{y}$, such that the reduced fraction $x_{M^{'}}$ has the following properties:

\begin{enumerate}[(1)]

\item No frozen variable appears in the denominator of $x_{M^{\prime}}.$

\item No frozen variable which appears in the numerator of $x_{M^{\prime}}$ is squared. 

\end{enumerate}

The rest of the morphisms in $\mathcal{F}$ are compositions of such multiplications. 

\begin{thm}
    Consider the cluster variable $x_{w},$ and let $\text{Supp}(x_{w})$ be the set of Laurent monomials in the Laurent expansion of $x_w$. Let $x_M \in \text{Supp}(x_{w}),$ and let $\mathcal{O}(x_M)$ be the connected component of $\mathcal{F}$ containing $x_M$. Then $$\text{Supp}(x_{w}) = \mathcal{O}(x_M). $$
\end{thm}

\begin{proof}
Let $x_{M^{\prime}} \in \text{Supp}(x_w).$ If we represent both $x_{M}$ and $x_{M^{\prime}}$ as $T$-paths, then there is some sequence of $T$-path twists that takes $x_{M}$ to $x_{M^{\prime}}.$ Notice that if two $T$-paths are related by a twist then there is a morphism between their respective weights. Indeed, $T$-path twists algebraically are multiplication by some $\hat{y}$ as in (1), and furthermore a $T$-path has no red boundary edges nor does it use the same (blue) edge twice. Thus, there is a morphism $x_{M^{\prime}} \longrightarrow x_{M}$ and so $x_{M^{\prime}} \in \text{Supp}(x_M).$ The reverse inclusion is similar. 
\end{proof}

\begin{cor}
     Any cluster variable $x_{w}$ is completely determined by any one of its monomials. 
\end{cor}

\chapter {$T$-paths for Configurations of Flags}

The work in this chapter is joint with Nicholas Ovenhouse.

In \cite{fock2006moduli}, a generalization of the decorated Teichm\"uller spaces was introduced, called higher Teichm\"uller spaces. These moduli spaces take as input in their construction a Lie group $G$ and a marked surface $S$. It was shown in \cite{fock2006moduli} that when $G = \text{SL}_3$, the coordinate ring of this moduli spaces has a cluster structure. 

Here we focus on the case when $G = \text{SL}_3$ and $S$ is a disc with marked points on the boundary. In this case, the moduli space reduces to the moduli space of configurations of affine flags.

In this final chapter, our first goal is to give a Laurent expansion formula which generalizes the $T$-path formula from Section \ref{T_path_section}, in the special case when the initial seed is constructed from a fan triangulation. Our second goal is to describe the poset structure of some of these Laurent expansions.

\section{Decorated Teichm\"uller Spaces}

\begin{defn}
    Let $\Sigma$ be a marked surface of genus $g$ with $b$ boundary components and $m$ marked points located on boundary components, where each boundary component has at least one marked point. The $\textit{Teichm{\"u}ller space}$ is the space of marked complete metrics on $\Sigma$ having constant negative curvature $-1$ and a finite area, modulo the action of the connected component of the identity in the group of diffeomorphisms of $\Sigma.$ Denote this space by $\mathcal{T}^{m}_{g,b}.$
\end{defn}

\begin{defn}
    Let $\Sigma , g , b , m$ be as above. Let $p$ be a marked point on the boundary of $\Sigma$. A \textit{horocycle centered at $p$} is a circle orthogonal to any geodesic passing through $p$. Any horocycle centered at $p$ can be parameterized by a positive real number called the \textit{height} of the horocycle. The \textit{decorated Teichm{\"u}ller space} is denoted by $\widetilde{\mathcal{T}_{g,b}^{m}}$ and is defined as the total space of the trivial fiber bundle $\widetilde{\mathcal{T}_{g,b}^{m}} \twoheadrightarrow \mathcal{T}_{g,b}^{m}$ where the projection map is forgetting about horocycles. Each fiber is isomorphic to $\mathbb{R}^{m}_{+}$.
\end{defn}

We now recall the construction of \textit{Penner coordinates} on the decorated Teichm{\"u}ller space. Let $\Delta$ be an ideal triangulation of $\Sigma$ by geodesic arcs. Each edge $e$ of the ideal triangulation has infinite hyperbolic length. However, we can define the \textit{length} $l (e)$ of $e$ to be the signed finite length of the segment between the two horocycles centered at the endpoints of $e,$ where the sign of $l(e)$ is positive if the two horocycles don't intersect, and negative if they do.

Define $f(e) = \text{exp} \Big( \frac{l (e)}{2} \Big).$ The functions $f(e)$ where $e$ ranges over the edges in the triangulation $\Delta$ give a homeomorphism from the decorated Teichm{\"u}ller space to $\mathbb{R}^{6g+3b+2m-6}.$ The functions $f(e)$ are called \textit{Penner coordinates}, or \textit{lambda lengths}.

Any two Penner coordinates attached to the two arcs $p$ and $q$ involved in a flip of the underlying ideal triangulation are related by a  \textit{generalized Ptolemy relation}:
$$f(p) f(q) = f(a) f(c) + f(b) f(d),$$
where arcs $a$ and $c$ (respectively, $b$ and $d$) are opposite one another in the unique quadrilateral with edges from $\Delta$ determined by the arcs $p$ and $q$.

\begin{defn}
    The \textit{cluster algebra structure on $\widetilde{\mathcal{T}_{g,b}^{m}}$}, denoted $\mathcal{A} (\widetilde{\mathcal{T}_{g,b}^{m}})$, is given by building a quiver, depending on $\Delta$ as before, but with nodes now labeled by Penner coordinates.
\end{defn}

We have the following correspondences:
$$\text{cluster variables in $\mathcal{A} (\widetilde{\mathcal{T}_{g,b}^{m}})$ } \longleftrightarrow \text{ Penner coordinates of arcs in $\Delta$}$$
$$\text{ seeds of $\mathcal{A} (\widetilde{\mathcal{T}_{g,b}^{m}})$} \longleftrightarrow \text{ triangulations of $\Delta$}$$
$$\text{ seed mutations in $\mathcal{A} (\widetilde{\mathcal{T}_{g,b}^{m}})$} \longleftrightarrow \text{ flips in $\Delta$}$$

\section{Higher Teichm\"uller Spaces}

We recall in this section a generalization of Teichm\"uller spaces introduced in \cite{fock2006moduli}. 

\begin{defn}
    Let $G$ be a Lie group, and $S$ a marked surface. The space of $\textit{$G$-local systems}$ on $S$ is the character variety $\text{Hom}(\pi_1 (S) , G)/G,$ where the quotient is by conjugation.  
\end{defn}

The moduli space $\mathcal{M}$ defined in \cite{fock2006moduli} parameterizes local systems on $S$ with some additional structure, namely a choice of element in the quotient $G/U$ for each marked point. Here, $U = [B,B]$ is the unipotent radical of a Borel subgroup $B \subset G$.

When $G = \text{SL}_n$, the quotient $G/U$ parameterizes \textit{affine flags}. 

\begin{defn}
    An \textit{affine flag} in a vector space $V$ is a saturated chain of subspaces 
    $$V_1 \subset V_2 \subset \dots \subset V_{n-1} \subset V_n = \mathbb{R}^{n},$$
    along with a choice of nonzero vector $v_i \in V_{i+1} / V_{i}$ in each successive quotient for $1 < i < n-1$.  
\end{defn}

When $G = \text{SL}_2$, the space $\mathcal{M}$ parameterizes configurations of affine flags in $\mathbb{R}^{2}$. Note that an affine flag in $\mathbb{R}^{2}$ is simply a choice of non-zero vector. The surface type cluster algebras we have considered in the previous chapters may be interpreted as rings of functions on $\mathcal{M}$. 

Fix $G = \text{SL}_3$ for the remainder of this chapter. In \cite{fock2006moduli}, it was shown that in this case the coordinate ring of $\mathcal{M}$ has a cluster algebra structure. We now describe how to construct a seed for this cluster algebra. 

As we have seen, when $G = \text{SL}_2$ the nodes of the quiver associated to a triangulation are in one-to-one correspondence with the arcs (and boundary segments) of the triangulation, and the arrows of this quiver form clockwise $3$-cycles inside the triangles of the triangulation. Instead, when $G = \text{SL}_3$, each arc has two quiver nodes associated to it, and furthermore there is a vertex of the quiver for each triangle cut out by the triangulation. Now, each ideal triangle cut out by the triangulation has three $3$-cycles contained within it, each having the internal face vertex in common. This quiver is called a \textit{$3$-triangulation}. 

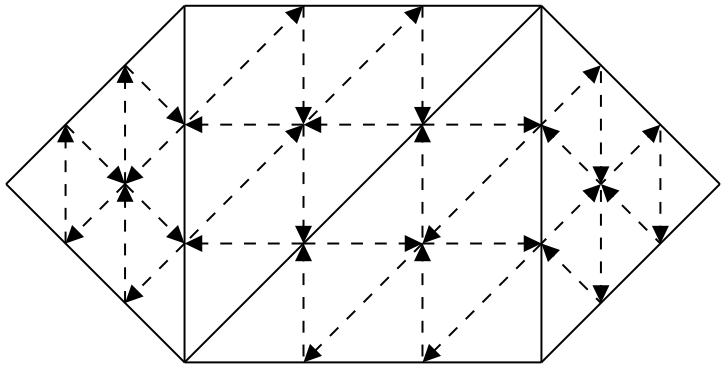
\begin {figure}[h!]
    \centering
    \caption{A $3$-triangulation}
    \label{fig:3_triangulation}

\tikzset{every picture/.style={line width=0.75pt}} 

\begin{tikzpicture}[x=0.75pt,y=0.75pt,yscale=-1,xscale=1]

\draw    (530,180) -- (620,90) ;
\draw    (620,270) -- (530,180) ;
\draw    (620,90) -- (800,90) ;
\draw    (620,270) -- (800,270) ;
\draw    (800,270) -- (890,180) ;
\draw    (890,180) -- (800,90) ;
\draw    (620,270) -- (620,90) ;
\draw    (620,270) -- (800,90) ;
\draw    (800,270) -- (800,90) ;
\draw  [dash pattern={on 4.5pt off 4.5pt}]  (560,150) -- (587.88,177.88) ;
\draw [shift={(590,180)}, rotate = 225] [fill={rgb, 255:red, 0; green, 0; blue, 0 }  ][line width=0.08]  [draw opacity=0] (8.93,-4.29) -- (0,0) -- (8.93,4.29) -- cycle    ;
\draw  [dash pattern={on 4.5pt off 4.5pt}]  (590,180) -- (562.12,207.88) ;
\draw [shift={(560,210)}, rotate = 315] [fill={rgb, 255:red, 0; green, 0; blue, 0 }  ][line width=0.08]  [draw opacity=0] (8.93,-4.29) -- (0,0) -- (8.93,4.29) -- cycle    ;
\draw  [dash pattern={on 4.5pt off 4.5pt}]  (560,210) -- (560,153) ;
\draw [shift={(560,150)}, rotate = 450] [fill={rgb, 255:red, 0; green, 0; blue, 0 }  ][line width=0.08]  [draw opacity=0] (8.93,-4.29) -- (0,0) -- (8.93,4.29) -- cycle    ;
\draw  [dash pattern={on 4.5pt off 4.5pt}]  (590,120) -- (617.88,147.88) ;
\draw [shift={(620,150)}, rotate = 225] [fill={rgb, 255:red, 0; green, 0; blue, 0 }  ][line width=0.08]  [draw opacity=0] (8.93,-4.29) -- (0,0) -- (8.93,4.29) -- cycle    ;
\draw  [dash pattern={on 4.5pt off 4.5pt}]  (620,150) -- (592.12,177.88) ;
\draw [shift={(590,180)}, rotate = 315] [fill={rgb, 255:red, 0; green, 0; blue, 0 }  ][line width=0.08]  [draw opacity=0] (8.93,-4.29) -- (0,0) -- (8.93,4.29) -- cycle    ;
\draw  [dash pattern={on 4.5pt off 4.5pt}]  (590,180) -- (590,123) ;
\draw [shift={(590,120)}, rotate = 450] [fill={rgb, 255:red, 0; green, 0; blue, 0 }  ][line width=0.08]  [draw opacity=0] (8.93,-4.29) -- (0,0) -- (8.93,4.29) -- cycle    ;
\draw  [dash pattern={on 4.5pt off 4.5pt}]  (590,180) -- (617.88,207.88) ;
\draw [shift={(620,210)}, rotate = 225] [fill={rgb, 255:red, 0; green, 0; blue, 0 }  ][line width=0.08]  [draw opacity=0] (8.93,-4.29) -- (0,0) -- (8.93,4.29) -- cycle    ;
\draw  [dash pattern={on 4.5pt off 4.5pt}]  (620,210) -- (592.12,237.88) ;
\draw [shift={(590,240)}, rotate = 315] [fill={rgb, 255:red, 0; green, 0; blue, 0 }  ][line width=0.08]  [draw opacity=0] (8.93,-4.29) -- (0,0) -- (8.93,4.29) -- cycle    ;
\draw  [dash pattern={on 4.5pt off 4.5pt}]  (590,240) -- (590,183) ;
\draw [shift={(590,180)}, rotate = 450] [fill={rgb, 255:red, 0; green, 0; blue, 0 }  ][line width=0.08]  [draw opacity=0] (8.93,-4.29) -- (0,0) -- (8.93,4.29) -- cycle    ;
\draw  [dash pattern={on 4.5pt off 4.5pt}]  (832.12,182.12) -- (860,210) ;
\draw [shift={(830,180)}, rotate = 45] [fill={rgb, 255:red, 0; green, 0; blue, 0 }  ][line width=0.08]  [draw opacity=0] (8.93,-4.29) -- (0,0) -- (8.93,4.29) -- cycle    ;
\draw  [dash pattern={on 4.5pt off 4.5pt}]  (857.88,152.12) -- (830,180) ;
\draw [shift={(860,150)}, rotate = 135] [fill={rgb, 255:red, 0; green, 0; blue, 0 }  ][line width=0.08]  [draw opacity=0] (8.93,-4.29) -- (0,0) -- (8.93,4.29) -- cycle    ;
\draw  [dash pattern={on 4.5pt off 4.5pt}]  (860,207) -- (860,150) ;
\draw [shift={(860,210)}, rotate = 270] [fill={rgb, 255:red, 0; green, 0; blue, 0 }  ][line width=0.08]  [draw opacity=0] (8.93,-4.29) -- (0,0) -- (8.93,4.29) -- cycle    ;
\draw  [dash pattern={on 4.5pt off 4.5pt}]  (827.88,182.12) -- (800,210) ;
\draw [shift={(830,180)}, rotate = 135] [fill={rgb, 255:red, 0; green, 0; blue, 0 }  ][line width=0.08]  [draw opacity=0] (8.93,-4.29) -- (0,0) -- (8.93,4.29) -- cycle    ;
\draw  [dash pattern={on 4.5pt off 4.5pt}]  (830,177) -- (830,120) ;
\draw [shift={(830,180)}, rotate = 270] [fill={rgb, 255:red, 0; green, 0; blue, 0 }  ][line width=0.08]  [draw opacity=0] (8.93,-4.29) -- (0,0) -- (8.93,4.29) -- cycle    ;
\draw  [dash pattern={on 4.5pt off 4.5pt}]  (830,237) -- (830,180) ;
\draw [shift={(830,240)}, rotate = 270] [fill={rgb, 255:red, 0; green, 0; blue, 0 }  ][line width=0.08]  [draw opacity=0] (8.93,-4.29) -- (0,0) -- (8.93,4.29) -- cycle    ;
\draw  [dash pattern={on 4.5pt off 4.5pt}]  (802.12,212.12) -- (830,240) ;
\draw [shift={(800,210)}, rotate = 45] [fill={rgb, 255:red, 0; green, 0; blue, 0 }  ][line width=0.08]  [draw opacity=0] (8.93,-4.29) -- (0,0) -- (8.93,4.29) -- cycle    ;
\draw  [dash pattern={on 4.5pt off 4.5pt}]  (802.12,152.12) -- (830,180) ;
\draw [shift={(800,150)}, rotate = 45] [fill={rgb, 255:red, 0; green, 0; blue, 0 }  ][line width=0.08]  [draw opacity=0] (8.93,-4.29) -- (0,0) -- (8.93,4.29) -- cycle    ;
\draw  [dash pattern={on 4.5pt off 4.5pt}]  (827.88,122.12) -- (800,150) ;
\draw [shift={(830,120)}, rotate = 135] [fill={rgb, 255:red, 0; green, 0; blue, 0 }  ][line width=0.08]  [draw opacity=0] (8.93,-4.29) -- (0,0) -- (8.93,4.29) -- cycle    ;
\draw  [dash pattern={on 4.5pt off 4.5pt}]  (677.88,92.12) -- (620,150) ;
\draw [shift={(680,90)}, rotate = 135] [fill={rgb, 255:red, 0; green, 0; blue, 0 }  ][line width=0.08]  [draw opacity=0] (8.93,-4.29) -- (0,0) -- (8.93,4.29) -- cycle    ;
\draw  [dash pattern={on 4.5pt off 4.5pt}]  (677.88,152.12) -- (620,210) ;
\draw [shift={(680,150)}, rotate = 135] [fill={rgb, 255:red, 0; green, 0; blue, 0 }  ][line width=0.08]  [draw opacity=0] (8.93,-4.29) -- (0,0) -- (8.93,4.29) -- cycle    ;
\draw  [dash pattern={on 4.5pt off 4.5pt}]  (737.88,92.12) -- (680,150) ;
\draw [shift={(740,90)}, rotate = 135] [fill={rgb, 255:red, 0; green, 0; blue, 0 }  ][line width=0.08]  [draw opacity=0] (8.93,-4.29) -- (0,0) -- (8.93,4.29) -- cycle    ;
\draw  [dash pattern={on 4.5pt off 4.5pt}]  (680,90) -- (680,147) ;
\draw [shift={(680,150)}, rotate = 270] [fill={rgb, 255:red, 0; green, 0; blue, 0 }  ][line width=0.08]  [draw opacity=0] (8.93,-4.29) -- (0,0) -- (8.93,4.29) -- cycle    ;
\draw  [dash pattern={on 4.5pt off 4.5pt}]  (680,150) -- (623,150) ;
\draw [shift={(620,150)}, rotate = 360] [fill={rgb, 255:red, 0; green, 0; blue, 0 }  ][line width=0.08]  [draw opacity=0] (8.93,-4.29) -- (0,0) -- (8.93,4.29) -- cycle    ;
\draw  [dash pattern={on 4.5pt off 4.5pt}]  (680,150) -- (680,207) ;
\draw [shift={(680,210)}, rotate = 270] [fill={rgb, 255:red, 0; green, 0; blue, 0 }  ][line width=0.08]  [draw opacity=0] (8.93,-4.29) -- (0,0) -- (8.93,4.29) -- cycle    ;
\draw  [dash pattern={on 4.5pt off 4.5pt}]  (680,210) -- (623,210) ;
\draw [shift={(620,210)}, rotate = 360] [fill={rgb, 255:red, 0; green, 0; blue, 0 }  ][line width=0.08]  [draw opacity=0] (8.93,-4.29) -- (0,0) -- (8.93,4.29) -- cycle    ;
\draw  [dash pattern={on 4.5pt off 4.5pt}]  (740,90) -- (740,147) ;
\draw [shift={(740,150)}, rotate = 270] [fill={rgb, 255:red, 0; green, 0; blue, 0 }  ][line width=0.08]  [draw opacity=0] (8.93,-4.29) -- (0,0) -- (8.93,4.29) -- cycle    ;
\draw  [dash pattern={on 4.5pt off 4.5pt}]  (740,150) -- (683,150) ;
\draw [shift={(680,150)}, rotate = 360] [fill={rgb, 255:red, 0; green, 0; blue, 0 }  ][line width=0.08]  [draw opacity=0] (8.93,-4.29) -- (0,0) -- (8.93,4.29) -- cycle    ;
\draw  [dash pattern={on 4.5pt off 4.5pt}]  (800,210) -- (742.12,267.88) ;
\draw [shift={(740,270)}, rotate = 315] [fill={rgb, 255:red, 0; green, 0; blue, 0 }  ][line width=0.08]  [draw opacity=0] (8.93,-4.29) -- (0,0) -- (8.93,4.29) -- cycle    ;
\draw  [dash pattern={on 4.5pt off 4.5pt}]  (800,150) -- (742.12,207.88) ;
\draw [shift={(740,210)}, rotate = 315] [fill={rgb, 255:red, 0; green, 0; blue, 0 }  ][line width=0.08]  [draw opacity=0] (8.93,-4.29) -- (0,0) -- (8.93,4.29) -- cycle    ;
\draw  [dash pattern={on 4.5pt off 4.5pt}]  (740,210) -- (682.12,267.88) ;
\draw [shift={(680,270)}, rotate = 315] [fill={rgb, 255:red, 0; green, 0; blue, 0 }  ][line width=0.08]  [draw opacity=0] (8.93,-4.29) -- (0,0) -- (8.93,4.29) -- cycle    ;
\draw  [dash pattern={on 4.5pt off 4.5pt}]  (680,210) -- (737,210) ;
\draw [shift={(740,210)}, rotate = 180] [fill={rgb, 255:red, 0; green, 0; blue, 0 }  ][line width=0.08]  [draw opacity=0] (8.93,-4.29) -- (0,0) -- (8.93,4.29) -- cycle    ;
\draw  [dash pattern={on 4.5pt off 4.5pt}]  (680,213) -- (680,270) ;
\draw [shift={(680,210)}, rotate = 90] [fill={rgb, 255:red, 0; green, 0; blue, 0 }  ][line width=0.08]  [draw opacity=0] (8.93,-4.29) -- (0,0) -- (8.93,4.29) -- cycle    ;
\draw  [dash pattern={on 4.5pt off 4.5pt}]  (740,153) -- (740,210) ;
\draw [shift={(740,150)}, rotate = 90] [fill={rgb, 255:red, 0; green, 0; blue, 0 }  ][line width=0.08]  [draw opacity=0] (8.93,-4.29) -- (0,0) -- (8.93,4.29) -- cycle    ;
\draw  [dash pattern={on 4.5pt off 4.5pt}]  (740,150) -- (797,150) ;
\draw [shift={(800,150)}, rotate = 180] [fill={rgb, 255:red, 0; green, 0; blue, 0 }  ][line width=0.08]  [draw opacity=0] (8.93,-4.29) -- (0,0) -- (8.93,4.29) -- cycle    ;
\draw  [dash pattern={on 4.5pt off 4.5pt}]  (740,210) -- (797,210) ;
\draw [shift={(800,210)}, rotate = 180] [fill={rgb, 255:red, 0; green, 0; blue, 0 }  ][line width=0.08]  [draw opacity=0] (8.93,-4.29) -- (0,0) -- (8.93,4.29) -- cycle    ;
\draw  [dash pattern={on 4.5pt off 4.5pt}]  (740,213) -- (740,270) ;
\draw [shift={(740,210)}, rotate = 90] [fill={rgb, 255:red, 0; green, 0; blue, 0 }  ][line width=0.08]  [draw opacity=0] (8.93,-4.29) -- (0,0) -- (8.93,4.29) -- cycle    ;

\end{tikzpicture}

\end {figure}

Flips of the triangulation are now governed not by a single quiver mutation, but by a sequence of four mutations. More precisely, to perform a flip of the triangulation at the diagonal $\delta_i$, one must first mutate at both quiver vertices along the diagonal $\delta_i$, and then subsequently mutate at the two vertices inside the triangles on either side of the diagonal $\delta_i$. 

As for traditional surface cluster algebras, any cluster variable we are now considering may be expressed in terms of the initial cluster variables attached to this initial triangulation, via the flips and their governing mutation sequences just described.

\section{Colored $\text{SL}_3$ Diagrams}

Refer to the two initial cluster variables attached to any edge $\delta_i$ of the triangulation as \textit{directed edges}, and call any cluster variable attached to the interior of a triangle a \textit{face}. Each edge variable is visualized as a directed edge which is directed away from the endpoint of $\delta_i$ that it is closest to. Each face variable is visualized as a ``filling'' of the triangle containing it.

We visualize Laurent monomials in the initial cluster variables by representing each variable in such a monomial as just described, and superimposing each such representation onto the same diagram. Variables occurring in the numerator of a Laurent monomial will be pictured as blue, and those in the denominator will be pictured as red.  

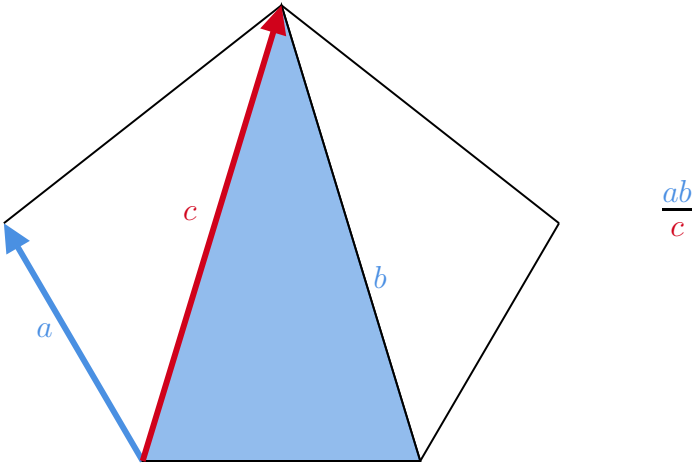
\begin {figure}[h!]
    \centering
    \caption{Visualization of a Laurent monomial}
    \label{fig:visualization2}

\tikzset{every picture/.style={line width=0.75pt}} 

\begin{tikzpicture}[x=0.75pt,y=0.75pt,yscale=-1,xscale=1]

\draw    (392,1100) -- (532,1100) ;
\draw [color={rgb, 255:red, 74; green, 144; blue, 226 }  ,draw opacity=1 ][line width=2.25]    (324.52,984.32) -- (392,1100) ;
\draw [shift={(322,980)}, rotate = 59.74] [fill={rgb, 255:red, 74; green, 144; blue, 226 }  ,fill opacity=1 ][line width=0.08]  [draw opacity=0] (14.29,-6.86) -- (0,0) -- (14.29,6.86) -- cycle    ;
\draw    (532,1100) -- (602,980) ;
\draw    (322,980) -- (462,870) ;
\draw    (462,870) -- (602,980) ;
\draw    (462,870) -- (532,1100) ;
\draw    (462,870) -- (392,1100) ;
\draw  [fill={rgb, 255:red, 74; green, 144; blue, 226 }  ,fill opacity=0.6 ] (462,870) -- (532,1100) -- (392,1100) -- cycle ;
\draw [color={rgb, 255:red, 208; green, 2; blue, 27 }  ,draw opacity=1 ][line width=2.25]    (460.54,874.78) -- (392,1100) ;
\draw [shift={(462,870)}, rotate = 106.93] [fill={rgb, 255:red, 208; green, 2; blue, 27 }  ,fill opacity=1 ][line width=0.08]  [draw opacity=0] (14.29,-6.86) -- (0,0) -- (14.29,6.86) -- cycle    ;

\draw (337,1029) node [anchor=north west][inner sep=0.75pt]    {$\textcolor[rgb]{0.29,0.56,0.89}{a}$};
\draw (411,970) node [anchor=north west][inner sep=0.75pt]    {$\textcolor[rgb]{0.82,0.01,0.11}{c}$};
\draw (507,1000) node [anchor=north west][inner sep=0.75pt]    {$\textcolor[rgb]{0.29,0.56,0.89}{b}$};
\draw (651,957) node [anchor=north west][inner sep=0.75pt]  [font=\Large]  {$\frac{\textcolor[rgb]{0.29,0.56,0.89}{ab}}{\textcolor[rgb]{0.82,0.01,0.11}{c}}$};

\end{tikzpicture}

\end {figure}

\begin {figure}[h!] 
    \centering
    \caption{Face resolution}
    \label{fig:face_resolution}

\tikzset{every picture/.style={line width=0.75pt}} 

\begin{tikzpicture}[x=0.75pt,y=0.75pt,yscale=-1,xscale=1]

\draw    (121.58,355.97) -- (200.81,355.97) ;
\draw    (121.58,276.75) -- (200.81,276.75) ;
\draw    (200.81,276.75) -- (200.81,355.97) ;
\draw    (121.58,276.75) -- (121.58,355.97) ;
\draw    (121.58,355.97) -- (200.81,276.75) ;
\draw  [fill={rgb, 255:red, 74; green, 144; blue, 226 }  ,fill opacity=0.6 ] (200.81,355.97) -- (121.58,276.75) -- (200.81,276.75) -- cycle ;
\draw  [fill={rgb, 255:red, 74; green, 144; blue, 226 }  ,fill opacity=0.6 ] (478.08,276.75) -- (398.86,355.97) -- (398.86,276.75) -- cycle ;
\draw    (272.1,355.97) -- (351.32,355.97) ;
\draw    (272.1,276.75) -- (351.32,276.75) ;
\draw    (351.32,276.75) -- (351.32,355.97) ;
\draw    (272.1,276.75) -- (272.1,355.97) ;
\draw    (272.1,355.97) -- (351.32,276.75) ;
\draw    (398.86,355.97) -- (478.08,355.97) ;
\draw    (398.86,276.75) -- (478.08,276.75) ;
\draw    (478.08,276.75) -- (478.08,355.97) ;
\draw    (398.86,276.75) -- (398.86,355.97) ;
\draw    (398.86,355.97) -- (478.08,276.75) ;
\draw  [fill={rgb, 255:red, 74; green, 144; blue, 226 }  ,fill opacity=0.6 ] (272.1,355.97) -- (351.32,276.75) -- (351.32,355.97) -- cycle ;
\draw [color={rgb, 255:red, 208; green, 2; blue, 27 }  ,draw opacity=1 ][line width=2.25]    (275.64,352.44) -- (351.32,276.75) ;
\draw [shift={(272.1,355.97)}, rotate = 315] [fill={rgb, 255:red, 208; green, 2; blue, 27 }  ,fill opacity=1 ][line width=0.08]  [draw opacity=0] (14.29,-6.86) -- (0,0) -- (14.29,6.86) -- cycle    ;
\draw [color={rgb, 255:red, 208; green, 2; blue, 27 }  ,draw opacity=1 ][line width=2.25]    (402.39,352.44) -- (478.08,276.75) ;
\draw [shift={(398.86,355.97)}, rotate = 315] [fill={rgb, 255:red, 208; green, 2; blue, 27 }  ,fill opacity=1 ][line width=0.08]  [draw opacity=0] (14.29,-6.86) -- (0,0) -- (14.29,6.86) -- cycle    ;
\draw [color={rgb, 255:red, 74; green, 144; blue, 226 }  ,draw opacity=1 ][line width=2.25]    (277.1,276.75) -- (351.32,276.75) ;
\draw [shift={(272.1,276.75)}, rotate = 0] [fill={rgb, 255:red, 74; green, 144; blue, 226 }  ,fill opacity=1 ][line width=0.08]  [draw opacity=0] (14.29,-6.86) -- (0,0) -- (14.29,6.86) -- cycle    ;
\draw [color={rgb, 255:red, 74; green, 144; blue, 226 }  ,draw opacity=1 ][line width=2.25]    (478.08,350.97) -- (478.08,276.75) ;
\draw [shift={(478.08,355.97)}, rotate = 270] [fill={rgb, 255:red, 74; green, 144; blue, 226 }  ,fill opacity=1 ][line width=0.08]  [draw opacity=0] (14.29,-6.86) -- (0,0) -- (14.29,6.86) -- cycle    ;

\draw (220.87,306.26) node [anchor=north west][inner sep=0.75pt]    {$=$};
\draw (364.51,306.26) node [anchor=north west][inner sep=0.75pt]    {$+$};

\end{tikzpicture}

\end {figure}
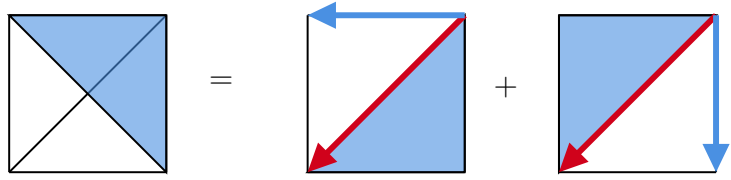

Any non-initial cluster variable may now be \textit{resolved} with respect to the initial triangulation. That is, it can be expressed as a sum over products of directed edges and faces from the initial triangulation. Figure \ref{fig:face_resolution} shows the result of resolving a non-initial face in a triangulated square. 

\newpage

\begin {figure}[h!]
    \centering
    \caption{Edge resolution}
    \label{fig:edge_resolution}

\tikzset{every picture/.style={line width=0.75pt}} 

\begin{tikzpicture}[x=0.75pt,y=0.75pt,yscale=-1,xscale=1]

\draw    (120,229.22) -- (199.22,229.22) ;
\draw    (120,150) -- (199.22,150) ;
\draw    (199.22,150) -- (199.22,229.22) ;
\draw    (120,150) -- (120,229.22) ;
\draw    (120,229.22) -- (199.22,150) ;
\draw    (270.52,229.22) -- (349.74,229.22) ;
\draw    (270.52,150) -- (349.74,150) ;
\draw    (349.74,150) -- (349.74,229.22) ;
\draw    (270.52,150) -- (270.52,229.22) ;
\draw    (270.52,229.22) -- (349.74,150) ;
\draw    (397.27,229.22) -- (476.49,229.22) ;
\draw    (397.27,150) -- (476.49,150) ;
\draw    (476.49,150) -- (476.49,229.22) ;
\draw    (397.27,150) -- (397.27,229.22) ;
\draw    (397.27,229.22) -- (476.49,150) ;
\draw [color={rgb, 255:red, 208; green, 2; blue, 27 }  ,draw opacity=1 ][line width=2.25]    (346.2,153.54) -- (270.52,229.22) ;
\draw [shift={(349.74,150)}, rotate = 135] [fill={rgb, 255:red, 208; green, 2; blue, 27 }  ,fill opacity=1 ][line width=0.08]  [draw opacity=0] (14.29,-6.86) -- (0,0) -- (14.29,6.86) -- cycle    ;
\draw [color={rgb, 255:red, 74; green, 144; blue, 226 }  ,draw opacity=1 ][line width=2.25]    (344.74,150) -- (270.52,150) ;
\draw [shift={(349.74,150)}, rotate = 180] [fill={rgb, 255:red, 74; green, 144; blue, 226 }  ,fill opacity=1 ][line width=0.08]  [draw opacity=0] (14.29,-6.86) -- (0,0) -- (14.29,6.86) -- cycle    ;
\draw    (524.03,229.22) -- (603.25,229.22) ;
\draw    (524.03,150) -- (603.25,150) ;
\draw    (603.25,150) -- (603.25,229.22) ;
\draw    (524.03,150) -- (524.03,229.22) ;
\draw    (524.03,229.22) -- (603.25,150) ;
\draw    (650.78,229.22) -- (730,229.22) ;
\draw    (650.78,150) -- (730,150) ;
\draw    (730,150) -- (730,229.22) ;
\draw    (650.78,150) -- (650.78,229.22) ;
\draw    (650.78,229.22) -- (730,150) ;
\draw [color={rgb, 255:red, 208; green, 2; blue, 27 }  ,draw opacity=1 ][line width=2.25]    (654.31,225.69) -- (730,150) ;
\draw [shift={(650.78,229.22)}, rotate = 315] [fill={rgb, 255:red, 208; green, 2; blue, 27 }  ,fill opacity=1 ][line width=0.08]  [draw opacity=0] (14.29,-6.86) -- (0,0) -- (14.29,6.86) -- cycle    ;
\draw [color={rgb, 255:red, 74; green, 144; blue, 226 }  ,draw opacity=1 ][line width=2.25]    (650.78,224.22) -- (650.78,150) ;
\draw [shift={(650.78,229.22)}, rotate = 270] [fill={rgb, 255:red, 74; green, 144; blue, 226 }  ,fill opacity=1 ][line width=0.08]  [draw opacity=0] (14.29,-6.86) -- (0,0) -- (14.29,6.86) -- cycle    ;
\draw [color={rgb, 255:red, 74; green, 144; blue, 226 }  ,draw opacity=1 ][line width=2.25]    (195.69,225.69) -- (120,150) ;
\draw [shift={(199.22,229.22)}, rotate = 225] [fill={rgb, 255:red, 74; green, 144; blue, 226 }  ,fill opacity=1 ][line width=0.08]  [draw opacity=0] (14.29,-6.86) -- (0,0) -- (14.29,6.86) -- cycle    ;
\draw [color={rgb, 255:red, 74; green, 144; blue, 226 }  ,draw opacity=1 ][line width=2.25]    (344.74,229.22) -- (270.52,229.22) ;
\draw [shift={(349.74,229.22)}, rotate = 180] [fill={rgb, 255:red, 74; green, 144; blue, 226 }  ,fill opacity=1 ][line width=0.08]  [draw opacity=0] (14.29,-6.86) -- (0,0) -- (14.29,6.86) -- cycle    ;
\draw  [fill={rgb, 255:red, 74; green, 144; blue, 226 }  ,fill opacity=0.6 ] (397.27,229.22) -- (476.49,150) -- (476.49,229.22) -- cycle ;
\draw  [fill={rgb, 255:red, 208; green, 2; blue, 27 }  ,fill opacity=0.5 ] (476.49,150) -- (397.27,229.22) -- (397.27,150) -- cycle ;
\draw [color={rgb, 255:red, 74; green, 144; blue, 226 }  ,draw opacity=1 ][line width=2.25]    (397.27,155) -- (397.27,229.22) ;
\draw [shift={(397.27,150)}, rotate = 90] [fill={rgb, 255:red, 74; green, 144; blue, 226 }  ,fill opacity=1 ][line width=0.08]  [draw opacity=0] (14.29,-6.86) -- (0,0) -- (14.29,6.86) -- cycle    ;
\draw [color={rgb, 255:red, 208; green, 2; blue, 27 }  ,draw opacity=1 ][line width=2.25]    (472.96,153.54) -- (397.27,229.22) ;
\draw [shift={(476.49,150)}, rotate = 135] [fill={rgb, 255:red, 208; green, 2; blue, 27 }  ,fill opacity=1 ][line width=0.08]  [draw opacity=0] (14.29,-6.86) -- (0,0) -- (14.29,6.86) -- cycle    ;
\draw [color={rgb, 255:red, 74; green, 144; blue, 226 }  ,draw opacity=1 ][line width=2.25]    (471.49,150) -- (397.27,150) ;
\draw [shift={(476.49,150)}, rotate = 180] [fill={rgb, 255:red, 74; green, 144; blue, 226 }  ,fill opacity=1 ][line width=0.08]  [draw opacity=0] (14.29,-6.86) -- (0,0) -- (14.29,6.86) -- cycle    ;
\draw  [fill={rgb, 255:red, 74; green, 144; blue, 226 }  ,fill opacity=0.6 ] (524.03,229.22) -- (603.25,150) -- (603.25,229.22) -- cycle ;
\draw  [fill={rgb, 255:red, 208; green, 2; blue, 27 }  ,fill opacity=0.5 ] (603.25,150) -- (524.03,229.22) -- (524.03,150) -- cycle ;
\draw [color={rgb, 255:red, 208; green, 2; blue, 27 }  ,draw opacity=1 ][line width=2.25]    (527.56,225.69) -- (603.25,150) ;
\draw [shift={(524.03,229.22)}, rotate = 315] [fill={rgb, 255:red, 208; green, 2; blue, 27 }  ,fill opacity=1 ][line width=0.08]  [draw opacity=0] (14.29,-6.86) -- (0,0) -- (14.29,6.86) -- cycle    ;
\draw [color={rgb, 255:red, 74; green, 144; blue, 226 }  ,draw opacity=1 ][line width=2.25]    (524.03,224.22) -- (524.03,150) ;
\draw [shift={(524.03,229.22)}, rotate = 270] [fill={rgb, 255:red, 74; green, 144; blue, 226 }  ,fill opacity=1 ][line width=0.08]  [draw opacity=0] (14.29,-6.86) -- (0,0) -- (14.29,6.86) -- cycle    ;
\draw [color={rgb, 255:red, 74; green, 144; blue, 226 }  ,draw opacity=1 ][line width=2.25]    (529.03,150) -- (603.25,150) ;
\draw [shift={(524.03,150)}, rotate = 0] [fill={rgb, 255:red, 74; green, 144; blue, 226 }  ,fill opacity=1 ][line width=0.08]  [draw opacity=0] (14.29,-6.86) -- (0,0) -- (14.29,6.86) -- cycle    ;
\draw [color={rgb, 255:red, 74; green, 144; blue, 226 }  ,draw opacity=1 ][line width=2.25]    (730,224.22) -- (730,150) ;
\draw [shift={(730,229.22)}, rotate = 270] [fill={rgb, 255:red, 74; green, 144; blue, 226 }  ,fill opacity=1 ][line width=0.08]  [draw opacity=0] (14.29,-6.86) -- (0,0) -- (14.29,6.86) -- cycle    ;

\draw (219.29,179.51) node [anchor=north west][inner sep=0.75pt]    {$=$};
\draw (362.92,179.51) node [anchor=north west][inner sep=0.75pt]    {$+$};
\draw (489.68,179.51) node [anchor=north west][inner sep=0.75pt]    {$+$};
\draw (616.43,179.51) node [anchor=north west][inner sep=0.75pt]    {$+$};

\end{tikzpicture}

\end {figure}
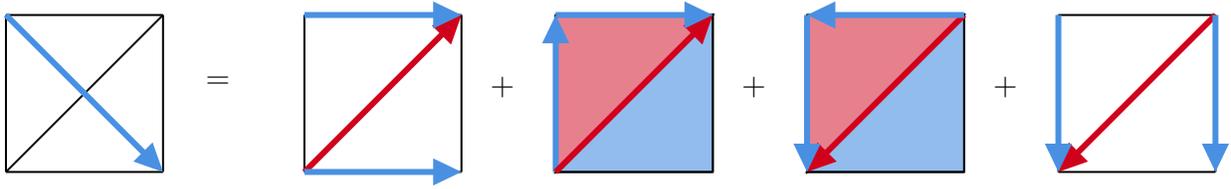

Figure \ref{fig:edge_resolution} shows the result of resolving a non-initial directed edge in a triangulated square.

    \begin{defn}
    Fix a triangulation $\Delta$ of the polygon $\Sigma$. A \textit{colored $\text{SL}_3$ diagram} is a multiset of directed edges and triangular faces from $\Delta$ drawn superimposed on the same polygon, where each element is colored either blue or red.
    \end{defn}
    
By the above discussion, colored $\text{SL}_3$ diagrams are in bijection with the set of Laurent monomials built from the initial cluster variables attached to the triangulation $\Delta$. The \textit{weight} of the $\text{SL}_3$ diagram $\pi$ is the associated Laurent monomial $x_{\pi}$.

\section{Crossings in Colored $\text{SL}_3$ Diagrams}

In this section, we say what it means for two colored $\text{SL}_3$ diagrams to be ``crossing''. From now on, we restrict to fan triangulations $\Delta$ (see Definition \ref{fan_def}). 

\begin{defn}
    Two directed arcs $i \rightarrow j$ and $k \rightarrow l$ are said to be \textit{crossing} if the underlying directed arcs are crossing in the usual sense. 
\end{defn}

\begin{defn} \label{cross_1_def}
    A directed arc $i \rightarrow j$ and a face $pqr$ are said to be \textit{crossing} if:
    
    \begin{enumerate}[(1)]
        \item $i=p$ and $q<j<r$ cyclically, i.e., $i \rightarrow j$ begins at one of the vertices of the triangle $pqr$, and has nontrivial intersection with the interior of $pqr$, or
        
        \item $p<i<q<r$ cyclically, i.e., the arc $i \rightarrow j$ begins and ends outside $pqr$, and intersects two sides of $pqr$. 
    \end{enumerate}
\end{defn}

See the next figure for an illustration of $(1)$ and $(2)$ in Definition \ref{cross_1_def}.

\begin {figure}[h!]
    \centering
    \caption{A directed edge and a face which cross}
    \label{fig:cross_1}

\begin{tikzpicture}[x=0.75pt,y=0.75pt,yscale=-1,xscale=1]

\draw [color={rgb, 255:red, 74; green, 144; blue, 226 }  ,draw opacity=1 ][line width=2.25]    (720,121.04) -- (615.78,121.78) ;
\draw [shift={(725,121)}, rotate = 179.59] [fill={rgb, 255:red, 74; green, 144; blue, 226 }  ,fill opacity=1 ][line width=0.08]  [draw opacity=0] (14.29,-6.86) -- (0,0) -- (14.29,6.86) -- cycle    ;
\draw  [fill={rgb, 255:red, 74; green, 144; blue, 226 }  ,fill opacity=0.5 ] (687.56,81.52) -- (671.8,177.8) -- (615.78,121.78) -- cycle ;
\draw  [fill={rgb, 255:red, 74; green, 144; blue, 226 }  ,fill opacity=0.6 ] (865,103.56) -- (785.78,159.22) -- (785.78,80) -- cycle ;
\draw [color={rgb, 255:red, 74; green, 144; blue, 226 }  ,draw opacity=1 ][line width=2.25]    (825,66) -- (825,171) ;
\draw [shift={(825,61)}, rotate = 90] [fill={rgb, 255:red, 74; green, 144; blue, 226 }  ,fill opacity=1 ][line width=0.08]  [draw opacity=0] (14.29,-6.86) -- (0,0) -- (14.29,6.86) -- cycle    ;

\draw (605,115) node [anchor=north west][inner sep=0.75pt]    {$i$};
\draw (724,111) node [anchor=north west][inner sep=0.75pt]    {$j$};
\draw (675,73) node [anchor=north west][inner sep=0.75pt]    {$r$};
\draw (665,179) node [anchor=north west][inner sep=0.75pt]    {$q$};
\draw (776,72) node [anchor=north west][inner sep=0.75pt]    {$r$};
\draw (775,160) node [anchor=north west][inner sep=0.75pt]    {$p$};
\draw (866,99) node [anchor=north west][inner sep=0.75pt]    {$q$};
\draw (820,173) node [anchor=north west][inner sep=0.75pt]    {$i$};
\draw (819,46) node [anchor=north west][inner sep=0.75pt]    {$j$};

\end{tikzpicture}

\end {figure}
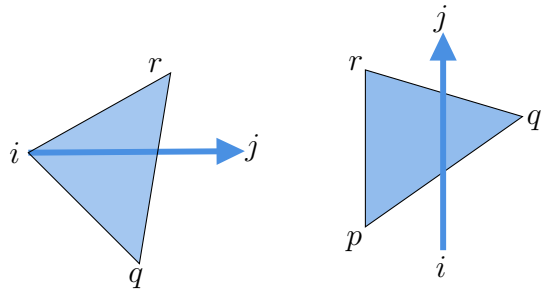

\begin{rmk}
    In (1) above, the directed arc $j \rightarrow i$ does \textit{not} cross the face $pqr$. In (2), the orientation of the arc is irrelevant. 
\end{rmk}

\begin{defn} \label{cross_2_def}
    Two faces $ijk$ and $pqr$ are considered to be \textit{crossing} if 
    
    \begin{enumerate}[(1)]
        \item $q=j$ and $i < j = q < r < k < p$ cyclically, or 
        
        \item $i < j < r < k < q < p$ cyclically. 
    \end{enumerate}
\end{defn}

See the next figure for an illustration of $(1)$ and $(2)$ in Definition \ref{cross_2_def}.

\begin {figure}[h!]
    \centering
    \caption{Two crossing faces}
    \label{fig:cross_2}
    \begin{tikzpicture}[x=0.75pt,y=0.75pt,yscale=-1,xscale=1]

\draw  [fill={rgb, 255:red, 74; green, 144; blue, 226 }  ,fill opacity=0.5 ] (280.78,1250) -- (70,1330) -- (240.78,1330) -- cycle ;
\draw  [fill={rgb, 255:red, 74; green, 144; blue, 226 }  ,fill opacity=0.5 ] (116.95,1375.39) -- (155.39,1244.61) -- (240.78,1330) -- cycle ;
\draw  [fill={rgb, 255:red, 74; green, 144; blue, 226 }  ,fill opacity=0.5 ] (410,1430) -- (490,1210) -- (490,1410) -- cycle ;
\draw  [fill={rgb, 255:red, 74; green, 144; blue, 226 }  ,fill opacity=0.5 ] (375.1,1258.13) -- (579.88,1371.56) -- (382.35,1340.28) -- cycle ;

\draw (55,1321) node [anchor=north west][inner sep=0.75pt]    {$p$};
\draw (145,1226) node [anchor=north west][inner sep=0.75pt]    {$k$};
\draw (275,1235) node [anchor=north west][inner sep=0.75pt]    {$r$};
\draw (245,1319) node [anchor=north west][inner sep=0.75pt]    {$j$};
\draw (105,1370) node [anchor=north west][inner sep=0.75pt]    {$i$};
\draw (365,1244) node [anchor=north west][inner sep=0.75pt]    {$p$};
\draw (580,1361) node [anchor=north west][inner sep=0.75pt]    {$r$};
\draw (369,1330) node [anchor=north west][inner sep=0.75pt]    {$q$};
\draw (483,1193) node [anchor=north west][inner sep=0.75pt]    {$k$};
\draw (399,1420) node [anchor=north west][inner sep=0.75pt]    {$i$};
\draw (492,1399) node [anchor=north west][inner sep=0.75pt]    {$j$};

\end{tikzpicture}

\end {figure}
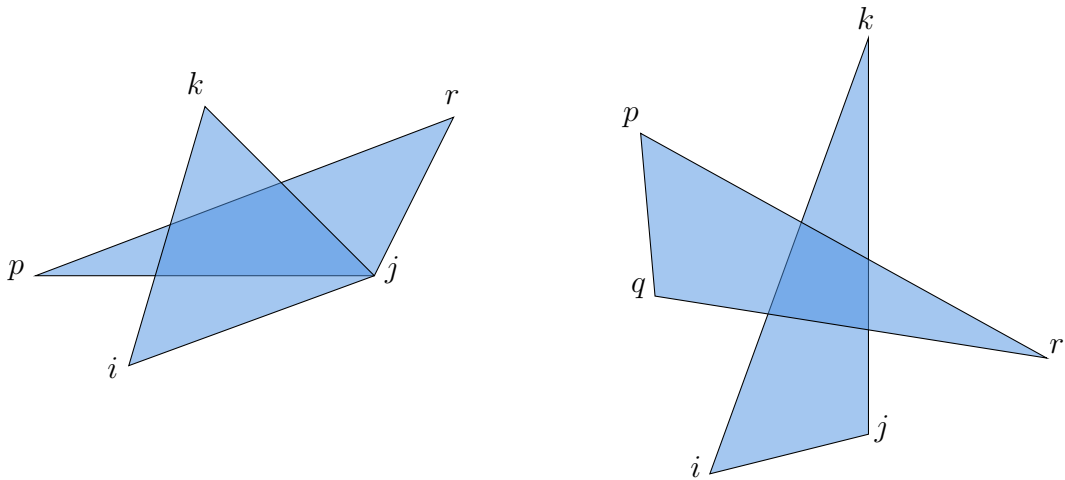

\section{Fork-Join Networks}

It is often convenient to model a simple process, algorithm, or computer program by 
a directed graph which encodes precedence. That is, the vertices are the states, and
the edges are directed so that $i \rightarrow j$ means that $i$ must happen before $j$.
In this way, a directed path $i_1 \rightarrow i_2 \rightarrow \cdots \rightarrow i_n$ models an algorithm
which is totally sequential (the steps have a definite order).

In an algorithm which utilizes parallelism, such as multithreading, the main procedure
may \textit{fork}, creating a \textit{child process}. After this, the main procedure and
the child process may execute in parallel. It may be necessary for the parent and child to
\textit{join} at some point before continuing. This means that the parent and child must
each finish their concurrent tasks before the main algorithm carries on in a sequential manner.

The next figure shows an example of a general fork-join network. 

\begin {figure}[h!]
    \centering
    \caption{A fork-join network}
    \label{fig:fork-join}
     \begin {tikzpicture}
        \coordinate (vstart) at (-1,0);
        \coordinate (v00) at (0,0);
        \coordinate (v10) at (1,0);
        \coordinate (v20) at (2,0);
        \coordinate (v30) at (3,0);
        \coordinate (v40) at (4,0);
        \coordinate (v50) at (5,0);
        \coordinate (v60) at (6,0);
        \coordinate (v70) at (7,0);
        \coordinate (v80) at (8,0);
        \coordinate (v90) at (9,0);
        \coordinate (v10_0) at (10,0);
        \coordinate (v11_0) at (11,0);
        \coordinate (v12_0) at (12,0);
        \coordinate (v13_0) at (13,0);

        \coordinate (v21) at (2,1);
        \coordinate (v31) at (3,1);
        \coordinate (v41) at (4,1);
        \coordinate (v51) at (5,1);
        \coordinate (v61) at (6,1);
        \coordinate (v71) at (7,1);
        \coordinate (v81) at (8,1);
        \coordinate (v91) at (9,1);
        \coordinate (v10_1) at (10,1);
        \coordinate (v11_1) at (11,1);

        \coordinate (v42) at (4,2);
        \coordinate (v52) at (5,2);
        \coordinate (v62) at (6,2);
        \coordinate (v72) at (7,2);
        \coordinate (v82) at (8,2);
        \coordinate (v92) at (9,2);

        \coordinate (v63) at (6,3);
        \coordinate (v73) at (7,3);

        \coordinate (v3n1) at (3,-1);
        \coordinate (v4n1) at (4,-1);

        \draw [fill=black] (vstart) circle (0.05);
        \draw [fill=black] (v00) circle (0.05);
        \draw [fill=black] (v10) circle (0.05);
        \draw [fill=black] (v20) circle (0.05);
        \draw [fill=black] (v30) circle (0.05);
        \draw [fill=black] (v40) circle (0.05);
        \draw [fill=black] (v50) circle (0.05);
        \draw [fill=black] (v60) circle (0.05);
        \draw [fill=black] (v70) circle (0.05);
        \draw [fill=black] (v80) circle (0.05);
        \draw [fill=black] (v90) circle (0.05);
        \draw [fill=black] (v10_0) circle (0.05);
        \draw [fill=black] (v11_0) circle (0.05);
        \draw [fill=black] (v12_0) circle (0.05);
        \draw [fill=black] (v13_0) circle (0.05);

        \draw [fill=black] (v21) circle (0.05);
        \draw [fill=black] (v31) circle (0.05);
        \draw [fill=black] (v41) circle (0.05);
        \draw [fill=black] (v51) circle (0.05);
        \draw [fill=black] (v61) circle (0.05);
        \draw [fill=black] (v71) circle (0.05);
        \draw [fill=black] (v81) circle (0.05);
        \draw [fill=black] (v91) circle (0.05);
        \draw [fill=black] (v10_1) circle (0.05);
        \draw [fill=black] (v11_1) circle (0.05);

        \draw [fill=black] (v42) circle (0.05);
        \draw [fill=black] (v52) circle (0.05);
        \draw [fill=black] (v62) circle (0.05);
        \draw [fill=black] (v72) circle (0.05);
        \draw [fill=black] (v82) circle (0.05);
        \draw [fill=black] (v92) circle (0.05);

        \draw [fill=black] (v63) circle (0.05);
        \draw [fill=black] (v73) circle (0.05);

        \draw [fill=black] (v3n1) circle (0.05);
        \draw [fill=black] (v4n1) circle (0.05);

        \draw [-latex] ($(vstart) + 0.05*(1,0)$) -- ($(v00) - 0.05*(1,0)$);
        \draw [-latex] ($(v00) + 0.05*(1,0)$) -- ($(v10) - 0.05*(1,0)$);
        \draw [-latex] ($(v10) + 0.05*(1,0)$) -- ($(v20) - 0.05*(1,0)$);
        \draw [-latex] ($(v20) + 0.05*(1,0)$) -- ($(v30) - 0.05*(1,0)$);
        \draw [-latex] ($(v30) + 0.05*(1,0)$) -- ($(v40) - 0.05*(1,0)$);
        \draw [-latex] ($(v40) + 0.05*(1,0)$) -- ($(v50) - 0.05*(1,0)$);
        \draw [-latex] ($(v50) + 0.05*(1,0)$) -- ($(v60) - 0.05*(1,0)$);
        \draw [-latex] ($(v60) + 0.05*(1,0)$) -- ($(v70) - 0.05*(1,0)$);
        \draw [-latex] ($(v70) + 0.05*(1,0)$) -- ($(v80) - 0.05*(1,0)$);
        \draw [-latex] ($(v80) + 0.05*(1,0)$) -- ($(v90) - 0.05*(1,0)$);
        \draw [-latex] ($(v90) + 0.05*(1,0)$) -- ($(v10_0) - 0.05*(1,0)$);
        \draw [-latex] ($(v10_0) + 0.05*(1,0)$) -- ($(v11_0) - 0.05*(1,0)$);
        \draw [-latex] ($(v11_0) + 0.05*(1,0)$) -- ($(v12_0) - 0.05*(1,0)$);
        \draw [-latex] ($(v12_0) + 0.05*(1,0)$) -- ($(v13_0) - 0.05*(1,0)$);

        \draw [-latex] ($(v21) + 0.05*(1,0)$) -- ($(v31) - 0.05*(1,0)$);
        \draw [-latex] ($(v31) + 0.05*(1,0)$) -- ($(v41) - 0.05*(1,0)$);
        \draw [-latex] ($(v41) + 0.05*(1,0)$) -- ($(v51) - 0.05*(1,0)$);
        \draw [-latex] ($(v51) + 0.05*(1,0)$) -- ($(v61) - 0.05*(1,0)$);
        \draw [-latex] ($(v61) + 0.05*(1,0)$) -- ($(v71) - 0.05*(1,0)$);
        \draw [-latex] ($(v71) + 0.05*(1,0)$) -- ($(v81) - 0.05*(1,0)$);
        \draw [-latex] ($(v81) + 0.05*(1,0)$) -- ($(v91) - 0.05*(1,0)$);
        \draw [-latex] ($(v91) + 0.05*(1,0)$) -- ($(v10_1) - 0.05*(1,0)$);
        \draw [-latex] ($(v10_1) + 0.05*(1,0)$) -- ($(v11_1) - 0.05*(1,0)$);

        \draw [-latex] ($(v42) + 0.05*(1,0)$) -- ($(v52) - 0.05*(1,0)$);
        \draw [-latex] ($(v52) + 0.05*(1,0)$) -- ($(v62) - 0.05*(1,0)$);
        \draw [-latex] ($(v62) + 0.05*(1,0)$) -- ($(v72) - 0.05*(1,0)$);
        \draw [-latex] ($(v72) + 0.05*(1,0)$) -- ($(v82) - 0.05*(1,0)$);
        \draw [-latex] ($(v82) + 0.05*(1,0)$) -- ($(v92) - 0.05*(1,0)$);

        \draw [-latex] ($(v63) + 0.05*(1,0)$) -- ($(v73) - 0.05*(1,0)$);

        \draw [-latex] ($(v3n1) + 0.05*(1,0)$) -- ($(v4n1) - 0.05*(1,0)$);

        \draw [-latex] ($(v10) + 0.05*(1,1)$) -- ($(v21) - 0.05*(1,1)$);
        \draw [-latex] ($(v31) + 0.05*(1,1)$) -- ($(v42) - 0.05*(1,1)$);
        \draw [-latex] ($(v52) + 0.05*(1,1)$) -- ($(v63) - 0.05*(1,1)$);
        \draw [-latex] ($(v4n1) + 0.05*(1,1)$) -- ($(v50) - 0.05*(1,1)$);

        \draw [-latex] ($(v20) + 0.05*(1,-1)$) -- ($(v3n1) - 0.05*(1,-1)$);
        \draw [-latex] ($(v73) + 0.05*(1,-1)$) -- ($(v82) - 0.05*(1,-1)$);
        \draw [-latex] ($(v92) + 0.05*(1,-1)$) -- ($(v10_1) - 0.05*(1,-1)$);
        \draw [-latex] ($(v11_1) + 0.05*(1,-1)$) -- ($(v12_0) - 0.05*(1,-1)$);
    \end {tikzpicture}
  
\end {figure}
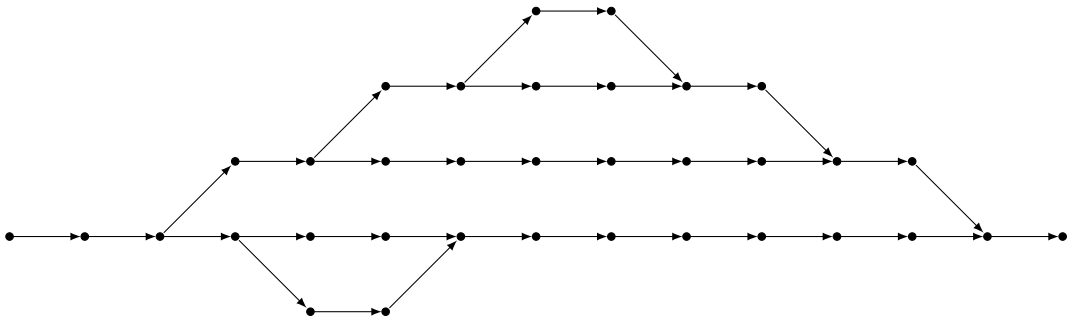

\newpage

We will only consider processes with the following two properties: 

\begin{enumerate}
    \item[(A1)] any fork must later be accompanied by a join, and
    
    \item[(A2)] there are no nested forks.
\end{enumerate}

From now on, whenever we say ``fork-join network'' we will mean a network satisfying these two restrictions. 

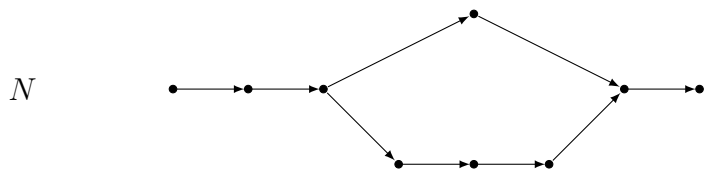
\begin {figure}[h!] \label{alt_fork_join}
    \centering
    \caption{A fork-join network satisfying conditions (A1) and (A2)}
    \label{fig:a1_a2}
    \begin{tikzpicture}
            \coordinate (v00) at (0,0);
        \coordinate (v10) at (1,0);
        \coordinate (v20) at (2,0);
        \coordinate (v3_1) at (3,-1);
        \coordinate (v4_1) at (4,-1);
        \coordinate (v41) at (4,1);
        \coordinate (v5_1) at (5,-1);
        \coordinate (v60) at (6,0);
        \coordinate (v70) at (7,0);

        \draw[fill=black] (v00) circle (0.05);
        \draw[fill=black] (v10) circle (0.05);
        \draw[fill=black] (v20) circle (0.05);
        \draw[fill=black] (v3_1) circle (0.05);
        \draw[fill=black] (v4_1) circle (0.05);
        \draw[fill=black] (v41) circle (0.05);
        \draw[fill=black] (v5_1) circle (0.05);
        \draw[fill=black] (v60) circle (0.05);
        \draw[fill=black] (v70) circle (0.05);

        \draw [-latex] ($(v00) + 0.05*(1,0)$) -- ($(v10) - 0.05*(1,0)$);
        \draw [-latex] ($(v10) + 0.05*(1,0)$) -- ($(v20) - 0.05*(1,0)$);
        \draw [-latex] ($(v20) + 0.05*(1,-1)$) -- ($(v3_1) - 0.05*(1,-1)$);
        \draw [-latex] ($(v3_1) + 0.05*(1,0)$) -- ($(v4_1) - 0.05*(1,0)$);
        \draw [-latex] ($(v4_1) + 0.05*(1,0)$) -- ($(v5_1) - 0.05*(1,0)$);
        \draw [-latex] ($(v5_1) + 0.05*(1,1)$) -- ($(v60) - 0.05*(1,1)$);
        \draw [-latex] ($(v60) + 0.05*(1,0)$) -- ($(v70) - 0.05*(1,0)$);

        \draw [-latex] ($(v20) + 0.03*(2,1)$) -- ($(v41) - 0.03*(2,1)$);
        \draw [-latex] ($(v41) + 0.03*(2,-1)$) -- ($(v60) - 0.03*(2,-1)$);

        \draw (-2,0) node {$N$};
  \end{tikzpicture}
\end {figure}

\begin{defn}
     An \textit{alternating fork-join network} is a directed graph $\Omega$ with the same underlying
    undirected graph as some fork-join network $N$ such that
    \begin{enumerate}[(1)]
        \item For any directed path through $N$ from beginning to end, the odd-numbered edges of $\Omega$ are oriented
              the same as the corresponding edge in $N$, and the even-numbered edges are oriented opposite.
        \item For any fork in $N$, the edge going into the fork must be inverted in $\Omega$.
        \item For any join in $N$, the edge leaving after the join must \textit{not} be inverted in $\Omega$.
    \end {enumerate}
\end{defn} 

We color each directed edge of an alternating fork-join network blue or red, according to the following rules. 

\begin{enumerate}[(1)]
    \item The first arrow is colored blue. 
    
    \item Away from forks and joins, colors alternate along paths. 
    
    \item At a fork or join, the three edges incident that vertex are all colored the same color.  
\end{enumerate}

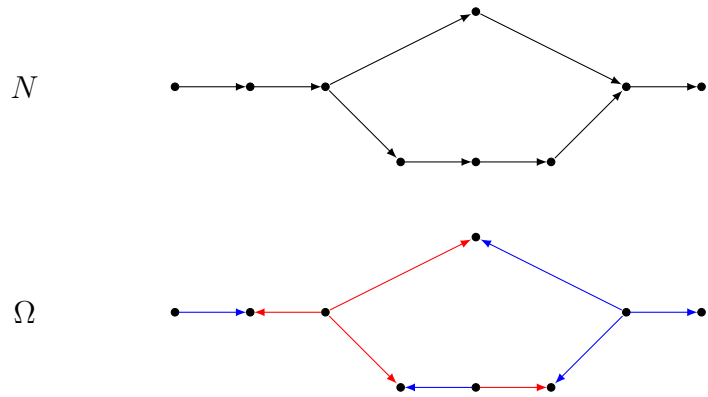
\begin {figure}[h!] 
    \centering
    \caption{Construction of an alternating fork-join network $\Omega$}
    \label{fig:alt_fork_join}
     \begin {tikzpicture}
        \coordinate (v00) at (0,0);
        \coordinate (v10) at (1,0);
        \coordinate (v20) at (2,0);
        \coordinate (v3_1) at (3,-1);
        \coordinate (v4_1) at (4,-1);
        \coordinate (v41) at (4,1);
        \coordinate (v5_1) at (5,-1);
        \coordinate (v60) at (6,0);
        \coordinate (v70) at (7,0);

        \draw[fill=black] (v00) circle (0.05);
        \draw[fill=black] (v10) circle (0.05);
        \draw[fill=black] (v20) circle (0.05);
        \draw[fill=black] (v3_1) circle (0.05);
        \draw[fill=black] (v4_1) circle (0.05);
        \draw[fill=black] (v41) circle (0.05);
        \draw[fill=black] (v5_1) circle (0.05);
        \draw[fill=black] (v60) circle (0.05);
        \draw[fill=black] (v70) circle (0.05);

        \draw [-latex] ($(v00) + 0.05*(1,0)$) -- ($(v10) - 0.05*(1,0)$);
        \draw [-latex] ($(v10) + 0.05*(1,0)$) -- ($(v20) - 0.05*(1,0)$);
        \draw [-latex] ($(v20) + 0.05*(1,-1)$) -- ($(v3_1) - 0.05*(1,-1)$);
        \draw [-latex] ($(v3_1) + 0.05*(1,0)$) -- ($(v4_1) - 0.05*(1,0)$);
        \draw [-latex] ($(v4_1) + 0.05*(1,0)$) -- ($(v5_1) - 0.05*(1,0)$);
        \draw [-latex] ($(v5_1) + 0.05*(1,1)$) -- ($(v60) - 0.05*(1,1)$);
        \draw [-latex] ($(v60) + 0.05*(1,0)$) -- ($(v70) - 0.05*(1,0)$);

        \draw [-latex] ($(v20) + 0.03*(2,1)$) -- ($(v41) - 0.03*(2,1)$);
        \draw [-latex] ($(v41) + 0.03*(2,-1)$) -- ($(v60) - 0.03*(2,-1)$);

        \draw (-2,0) node {$N$};


        \coordinate (V00)  at ($(v00) + (0,-3)$);
        \coordinate (V10)  at ($(v10) + (0,-3)$);
        \coordinate (V20)  at ($(v20) + (0,-3)$);
        \coordinate (V3_1) at ($(v3_1) + (0,-3)$);
        \coordinate (V4_1) at ($(v4_1) + (0,-3)$);
        \coordinate (V41)  at ($(v41) + (0,-3)$);
        \coordinate (V5_1) at ($(v5_1) + (0,-3)$);
        \coordinate (V60)  at ($(v60) + (0,-3)$);
        \coordinate (V70)  at ($(v70) + (0,-3)$);

        \draw[fill=black] (V00) circle (0.05);
        \draw[fill=black] (V10) circle (0.05);
        \draw[fill=black] (V20) circle (0.05);
        \draw[fill=black] (V3_1) circle (0.05);
        \draw[fill=black] (V4_1) circle (0.05);
        \draw[fill=black] (V41) circle (0.05);
        \draw[fill=black] (V5_1) circle (0.05);
        \draw[fill=black] (V60) circle (0.05);
        \draw[fill=black] (V70) circle (0.05);

        \draw [-latex, blue] ($(V00) + 0.05*(1,0)$)  -- ($(V10) - 0.05*(1,0)$);
        \draw [-latex, red]  ($(V20) - 0.05*(1,0)$)  -- ($(V10) + 0.05*(1,0)$);
        \draw [-latex, red]  ($(V20) + 0.05*(1,-1)$) -- ($(V3_1) - 0.05*(1,-1)$);
        \draw [-latex, blue] ($(V4_1) - 0.05*(1,0)$) -- ($(V3_1) + 0.05*(1,0)$);
        \draw [-latex, red]  ($(V4_1) + 0.05*(1,0)$) -- ($(V5_1) - 0.05*(1,0)$);
        \draw [-latex, blue] ($(V60) - 0.05*(1,1)$) -- ($(V5_1) + 0.05*(1,1)$);
        \draw [-latex, blue] ($(V60) + 0.05*(1,0)$)  -- ($(V70) - 0.05*(1,0)$);

        \draw [-latex, red]  ($(V20) + 0.03*(2,1)$)  -- ($(V41) - 0.03*(2,1)$);
        \draw [-latex, blue] ($(V60) - 0.03*(2,-1)$) -- ($(V41) + 0.03*(2,-1)$);

        \draw (-2,-3) node {$\Omega$};
    \end {tikzpicture}
  
\end {figure}

Let $\Omega_1$ and $\Omega_2$ be two directed graphs. A \textit{homomorphism of directed graphs} $g : \Omega_1 \longrightarrow \Omega_2$ is a mapping of vertex sets such that if $i \rightarrow j$ in $\Omega_1$, then $g(i) \rightarrow g(j)$ is an edge of $\Omega_2$.

\begin{defn}
    A homomorphism of directed graphs $g : \Omega_1 \longrightarrow \Omega_2$ is called an \textit{immersion} if for any vertex $v \in \Omega_1$, every vertex incident to $v$ maps to a distinct vertex under $g$. 
\end{defn}

\begin{defn}
    Let $\Omega_1$ be an alternating fork-join network, and $\Omega_2$ any directed graph. Suppose that $g: \Omega_1 \longrightarrow \Omega_2$ is an immersion. The \textit{image} of $g$ is the associated multiset of edges in $\Omega_2$. That is, if $i \rightarrow j$ is an edge in $\Omega_1$, then $g(i) \rightarrow g(j)$ is an edge in $\Omega_2$, and if more than one edge from $\Omega_1$ maps to $g(i) \rightarrow g(j)$, then we count the edge $g(i) \rightarrow g(j)$ with multiplicity equal to the cardinality of its preimage. 
\end{defn}

\section{Generalized $T$-paths} 

We now generalize to our current setup the $T$-paths from Section \ref{T_path_section}. To this end, we picture each face variable as a \textit{tripod}, as shown in the next figure. Each tripod resembles a web diagram from \cite{fomin2012tensor} and \cite{fomin2014webs}.

\begin {figure}[h!]
    \centering
    \caption{Faces as webs}
    \label{fig:face_web}

\tikzset{every picture/.style={line width=0.75pt}} 

\begin{tikzpicture}[x=0.75pt,y=0.75pt,yscale=-1,xscale=1]

\draw  [fill={rgb, 255:red, 74; green, 144; blue, 226 }  ,fill opacity=0.6 ] (520,1880) -- (440.78,1959.22) -- (440.78,1880) -- cycle ;
\draw    (440.78,1959.22) -- (520,1959.22) ;
\draw    (440.78,1880) -- (520,1880) ;
\draw    (520,1880) -- (520,1959.22) ;
\draw    (440.78,1880) -- (440.78,1959.22) ;
\draw    (440.78,1959.22) -- (520,1880) ;
\draw [color={rgb, 255:red, 0; green, 0; blue, 0 }  ,draw opacity=1 ][line width=0.75]    (440.78,1959.22) -- (520,1880) ;
\draw [color={rgb, 255:red, 0; green, 0; blue, 0 }  ,draw opacity=1 ][line width=0.75]    (520,1959.22) -- (520,1880) ;
\draw    (590,1960) -- (669.22,1960) ;
\draw    (590,1880.78) -- (669.22,1880.78) ;
\draw    (669.22,1880.78) -- (669.22,1960) ;
\draw    (590,1880.78) -- (590,1960) ;
\draw    (590,1960) -- (669.22,1880.78) ;
\draw [color={rgb, 255:red, 0; green, 0; blue, 0 }  ,draw opacity=1 ][line width=0.75]    (590,1960) -- (669.22,1880.78) ;
\draw [color={rgb, 255:red, 0; green, 0; blue, 0 }  ,draw opacity=1 ][line width=0.75]    (669.22,1960) -- (669.22,1880.78) ;
\draw [color={rgb, 255:red, 74; green, 144; blue, 226 }  ,draw opacity=1 ][line width=0.75]    (592.07,1882.95) -- (616,1908) ;
\draw [shift={(590,1880.78)}, rotate = 46.31] [fill={rgb, 255:red, 74; green, 144; blue, 226 }  ,fill opacity=1 ][line width=0.08]  [draw opacity=0] (8.93,-4.29) -- (0,0) -- (8.93,4.29) -- cycle    ;
\draw [color={rgb, 255:red, 74; green, 144; blue, 226 }  ,draw opacity=1 ][line width=0.75]    (666.55,1882.15) -- (616,1908) ;
\draw [shift={(669.22,1880.78)}, rotate = 152.91] [fill={rgb, 255:red, 74; green, 144; blue, 226 }  ,fill opacity=1 ][line width=0.08]  [draw opacity=0] (8.93,-4.29) -- (0,0) -- (8.93,4.29) -- cycle    ;
\draw [color={rgb, 255:red, 74; green, 144; blue, 226 }  ,draw opacity=1 ][line width=0.75]    (591.34,1957.32) -- (616,1908) ;
\draw [shift={(590,1960)}, rotate = 296.57] [fill={rgb, 255:red, 74; green, 144; blue, 226 }  ,fill opacity=1 ][line width=0.08]  [draw opacity=0] (8.93,-4.29) -- (0,0) -- (8.93,4.29) -- cycle    ;

\draw (539,1909) node [anchor=north west][inner sep=0.75pt]    {$=$};

\end{tikzpicture}

\end {figure}

Given the triangulation $\Delta$, define the directed graph $\Omega_{\Delta}$ as follows. There is one vertex of $\Omega_{\Delta}$ for each vertex of $\Delta$, and one vertex per triangle cut out by $\Delta$. For each edge in $\Delta$ (including boundary segments) connecting vertices $i$ and $j$, there is a pair of directed edges $i \rightarrow j$ and $j \rightarrow i$. For each triangle $ijk$ cut out by $\Delta$, there is a tripod consisting of three edges directed from the center vertex of the triangle to the three vertices $i,j,$ and $k$.

Note that if $g: \Omega \longrightarrow \Omega_{\Delta}$ is an immersion of a fork-join network, then the image is by definition a colored $\text{SL}_3$ diagram. 

\begin {defn}
    Suppose $\Omega$ is an alternating fork-join network, and $\Omega \to \Omega_\Delta$ is an immersion, which sends the source of $\Omega$ to $i$
    and the terminal vertex of $\Omega$ to $j$. We call the image of this immersion a 
    \textit{$T$-path of edge type} from $i$ to $j$ if:
    \begin {enumerate}
        \item[(T1)] the image of any path through $\Omega$ does not use any edge of $\Omega_\Delta$ more than once
        \item[(T2)] there are an odd number of elements in the diagram
        \item[(T3)] the red elements cross the directed arc $i \to j$
        \item[(T4)] once a path crosses an edge in $\Delta$, it never returns to the original half of the polygon
    \end {enumerate}
    The set of edge-type $T$-paths from $i$ to $j$ is denoted by $\mathbb{T}_{ij}$.
\end {defn}

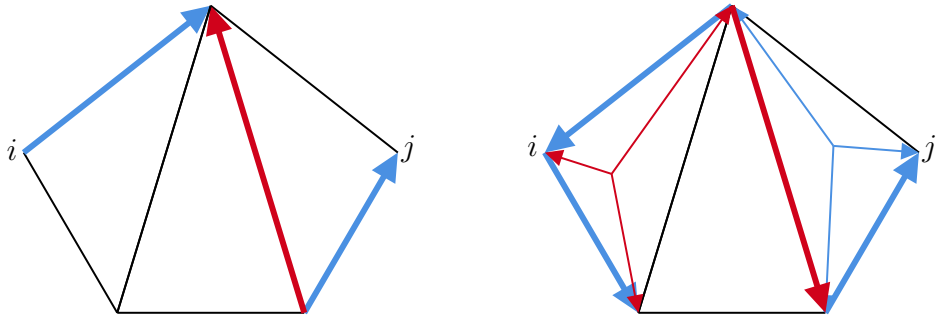
\begin{figure}[h!] 
    \centering
    \caption{Two edge-type $T$-paths}
    \label{fig:edge_T_paths}

\tikzset{every picture/.style={line width=0.75pt}} 

\begin{tikzpicture}[x=0.75pt,y=0.75pt,yscale=-1,xscale=1]

\draw    (193.69,2755) -- (288.03,2755) ;
\draw [color={rgb, 255:red, 0; green, 0; blue, 0 }  ,draw opacity=1 ][line width=0.75]    (146.52,2674.13) -- (193.69,2755) ;
\draw [color={rgb, 255:red, 74; green, 144; blue, 226 }  ,draw opacity=1 ][line width=2.25]    (288.03,2755) -- (332.69,2678.45) ;
\draw [shift={(335.21,2674.13)}, rotate = 480.26] [fill={rgb, 255:red, 74; green, 144; blue, 226 }  ,fill opacity=1 ][line width=0.08]  [draw opacity=0] (14.29,-6.86) -- (0,0) -- (14.29,6.86) -- cycle    ;
\draw [color={rgb, 255:red, 74; green, 144; blue, 226 }  ,draw opacity=1 ][line width=2.25]    (146.52,2674.13) -- (236.93,2603.09) ;
\draw [shift={(240.86,2600)}, rotate = 501.84] [fill={rgb, 255:red, 74; green, 144; blue, 226 }  ,fill opacity=1 ][line width=0.08]  [draw opacity=0] (14.29,-6.86) -- (0,0) -- (14.29,6.86) -- cycle    ;
\draw    (240.86,2600) -- (335.21,2674.13) ;
\draw [color={rgb, 255:red, 208; green, 2; blue, 27 }  ,draw opacity=1 ][line width=2.25]    (242.32,2604.78) -- (288.03,2755) ;
\draw [shift={(240.86,2600)}, rotate = 73.07] [fill={rgb, 255:red, 208; green, 2; blue, 27 }  ,fill opacity=1 ][line width=0.08]  [draw opacity=0] (14.29,-6.86) -- (0,0) -- (14.29,6.86) -- cycle    ;
\draw    (240.86,2600) -- (193.69,2755) ;
\draw [color={rgb, 255:red, 0; green, 0; blue, 0 }  ,draw opacity=1 ][line width=0.75]    (240.86,2600) -- (193.69,2755) ;
\draw    (456.51,2755) -- (550.85,2755) ;
\draw [color={rgb, 255:red, 74; green, 144; blue, 226 }  ,draw opacity=1 ][line width=2.25]    (409.34,2674.13) -- (453.99,2750.68) ;
\draw [shift={(456.51,2755)}, rotate = 239.74] [fill={rgb, 255:red, 74; green, 144; blue, 226 }  ,fill opacity=1 ][line width=0.08]  [draw opacity=0] (14.29,-6.86) -- (0,0) -- (14.29,6.86) -- cycle    ;
\draw [color={rgb, 255:red, 74; green, 144; blue, 226 }  ,draw opacity=1 ][line width=2.25]    (550.85,2755) -- (595.51,2678.45) ;
\draw [shift={(598.03,2674.13)}, rotate = 480.26] [fill={rgb, 255:red, 74; green, 144; blue, 226 }  ,fill opacity=1 ][line width=0.08]  [draw opacity=0] (14.29,-6.86) -- (0,0) -- (14.29,6.86) -- cycle    ;
\draw [color={rgb, 255:red, 74; green, 144; blue, 226 }  ,draw opacity=1 ][line width=2.25]    (413.27,2671.04) -- (503.68,2600) ;
\draw [shift={(409.34,2674.13)}, rotate = 321.84] [fill={rgb, 255:red, 74; green, 144; blue, 226 }  ,fill opacity=1 ][line width=0.08]  [draw opacity=0] (14.29,-6.86) -- (0,0) -- (14.29,6.86) -- cycle    ;
\draw    (503.68,2600) -- (598.03,2674.13) ;
\draw    (503.68,2600) -- (456.51,2755) ;
\draw [color={rgb, 255:red, 0; green, 0; blue, 0 }  ,draw opacity=1 ][line width=0.75]    (503.68,2600) -- (456.51,2755) ;
\draw [color={rgb, 255:red, 208; green, 2; blue, 27 }  ,draw opacity=1 ][line width=0.75]    (412.19,2675.04) -- (443.03,2684.91) ;
\draw [shift={(409.34,2674.13)}, rotate = 17.74] [fill={rgb, 255:red, 208; green, 2; blue, 27 }  ,fill opacity=1 ][line width=0.08]  [draw opacity=0] (8.93,-4.29) -- (0,0) -- (8.93,4.29) -- cycle    ;
\draw [color={rgb, 255:red, 208; green, 2; blue, 27 }  ,draw opacity=1 ][line width=0.75]    (501.94,2602.44) -- (443.03,2684.91) ;
\draw [shift={(503.68,2600)}, rotate = 125.54] [fill={rgb, 255:red, 208; green, 2; blue, 27 }  ,fill opacity=1 ][line width=0.08]  [draw opacity=0] (8.93,-4.29) -- (0,0) -- (8.93,4.29) -- cycle    ;
\draw [color={rgb, 255:red, 208; green, 2; blue, 27 }  ,draw opacity=1 ][line width=0.75]    (455.94,2752.05) -- (443.03,2684.91) ;
\draw [shift={(456.51,2755)}, rotate = 259.11] [fill={rgb, 255:red, 208; green, 2; blue, 27 }  ,fill opacity=1 ][line width=0.08]  [draw opacity=0] (8.93,-4.29) -- (0,0) -- (8.93,4.29) -- cycle    ;
\draw [color={rgb, 255:red, 74; green, 144; blue, 226 }  ,draw opacity=1 ][line width=0.75]    (595.03,2673.89) -- (554.9,2670.76) ;
\draw [shift={(598.03,2674.13)}, rotate = 184.47] [fill={rgb, 255:red, 74; green, 144; blue, 226 }  ,fill opacity=1 ][line width=0.08]  [draw opacity=0] (8.93,-4.29) -- (0,0) -- (8.93,4.29) -- cycle    ;
\draw [color={rgb, 255:red, 74; green, 144; blue, 226 }  ,draw opacity=1 ][line width=0.75]    (505.44,2602.43) -- (554.9,2670.76) ;
\draw [shift={(503.68,2600)}, rotate = 54.1] [fill={rgb, 255:red, 74; green, 144; blue, 226 }  ,fill opacity=1 ][line width=0.08]  [draw opacity=0] (8.93,-4.29) -- (0,0) -- (8.93,4.29) -- cycle    ;
\draw [color={rgb, 255:red, 74; green, 144; blue, 226 }  ,draw opacity=1 ][line width=0.75]    (551,2752) -- (554.9,2670.76) ;
\draw [shift={(550.85,2755)}, rotate = 272.75] [fill={rgb, 255:red, 74; green, 144; blue, 226 }  ,fill opacity=1 ][line width=0.08]  [draw opacity=0] (8.93,-4.29) -- (0,0) -- (8.93,4.29) -- cycle    ;
\draw [color={rgb, 255:red, 208; green, 2; blue, 27 }  ,draw opacity=1 ][line width=2.25]    (503.68,2600) -- (549.4,2750.21) ;
\draw [shift={(550.85,2755)}, rotate = 253.07] [fill={rgb, 255:red, 208; green, 2; blue, 27 }  ,fill opacity=1 ][line width=0.08]  [draw opacity=0] (14.29,-6.86) -- (0,0) -- (14.29,6.86) -- cycle    ;

\draw (136.49,2665.71) node [anchor=north west][inner sep=0.75pt]    {$i$};
\draw (399.04,2663.97) node [anchor=north west][inner sep=0.75pt]    {$i$};
\draw (598.36,2663.97) node [anchor=north west][inner sep=0.75pt]    {$j$};
\draw (335.54,2663.97) node [anchor=north west][inner sep=0.75pt]    {$j$};

\end{tikzpicture}

\end {figure}

Figure \ref{fig:edge_T_paths} shows two edge-type $T$-paths. Note that the $T$-path pictured on the right is an immersion of the alternating fork-join network $\Omega$ pictured in Figure \ref{fig:alt_fork_join}.

\begin{rmk}
    For any $T$-path defined in Section \ref{T_path_section}, there is a corresponding $\text{SL}_3$ edge-type $T$-path obtained by orienting edges (for instance, see the diagram on the left in Figure \ref{fig:edge_T_paths}).
\end{rmk}

Let $\gamma$ be a face which is not a triangle from the triangulation $\Delta$. Since we are assuming that $\Delta$ is a fan triangulation, there is at least one internal edge in $\Delta$ which terminates at one of the vertices of $\gamma$. We call any such edge a \textit{splitting edge}. The vertex that is shared by a splitting edge and the face $\gamma$ is called the \textit{near vertex}, and the other endpoint of this splitting edge is called the \textit{far vertex}.

We distinguish between two types of maximal (i.e., not contained in any subtriangulation) faces in a fan triangulation. There is the unique \textit{fan face}, which has every internal diagonal as a splitting edge; and there are $n$ \textit{split faces}, which have precisely one splitting edge. 

\begin {defn} Consider the face $ijk$. 
    A \textit{$T$-path of face type} (with vertices $i,j,k$) is a triple of immersed alternating fork-join networks, terminating at $i$, $j$, and $k$, respectively, and with common source a tripod vertex within a triangle of $\Delta$, such that
    \begin {enumerate}
        \item[(F1)] any of the three branches satisfies (T1), (T2), and (T4).
        \item[(F2)] their red elements cross the triangle $ijk$.
        \item[(F3)] if one of the three paths arrives at $i$, $j$, or $k$ using a forward blue edge, then that path must end
        \item[(F4)] when a branch crosses any splitting edge, it must 
                    cross either at the far end after using a blue edge or at the near end after
                    using a red edge.
    \end {enumerate} 
     The set of face-type $T$-paths with vertices $i,j,$ and $k$ is denoted by $\mathbb{T}_{ijk}$.
\end {defn}

See the next figure for two examples of face-type $T$-paths.

\begin {figure}[h!] \label{face_type_fig}
    \centering
    \caption{Two face-type $T$-paths}
    \label{fig:face-type}

\tikzset{every picture/.style={line width=0.75pt}} 

\begin{tikzpicture}[x=0.75pt,y=0.75pt,yscale=-1,xscale=1]

\draw    (326.5,310) -- (403,352.5) ;
\draw    (250,352.5) -- (326.5,310) ;
\draw    (403,437.5) -- (403,352.5) ;
\draw [color={rgb, 255:red, 74; green, 144; blue, 226 }  ,draw opacity=1 ][line width=2.25]    (250,437.5) -- (250,357.5) ;
\draw [shift={(250,352.5)}, rotate = 450] [fill={rgb, 255:red, 74; green, 144; blue, 226 }  ,fill opacity=1 ][line width=0.08]  [draw opacity=0] (14.29,-6.86) -- (0,0) -- (14.29,6.86) -- cycle    ;
\draw    (250,437.5) -- (326.5,480) ;
\draw    (326.5,480) -- (403,437.5) ;
\draw [color={rgb, 255:red, 208; green, 2; blue, 27 }  ,draw opacity=1 ][line width=2.25]    (564.5,480) -- (564.5,318.1) ;
\draw [shift={(564.5,313.1)}, rotate = 450] [fill={rgb, 255:red, 208; green, 2; blue, 27 }  ,fill opacity=1 ][line width=0.08]  [draw opacity=0] (14.29,-6.86) -- (0,0) -- (14.29,6.86) -- cycle    ;
\draw [color={rgb, 255:red, 208; green, 2; blue, 27 }  ,draw opacity=1 ][line width=2.25]    (400.43,433.21) -- (326.5,310) ;
\draw [shift={(403,437.5)}, rotate = 239.04] [fill={rgb, 255:red, 208; green, 2; blue, 27 }  ,fill opacity=1 ][line width=0.08]  [draw opacity=0] (14.29,-6.86) -- (0,0) -- (14.29,6.86) -- cycle    ;
\draw [color={rgb, 255:red, 74; green, 144; blue, 226 }  ,draw opacity=1 ][line width=2.25]    (326.5,475) -- (326.5,310) ;
\draw [shift={(326.5,480)}, rotate = 270] [fill={rgb, 255:red, 74; green, 144; blue, 226 }  ,fill opacity=1 ][line width=0.08]  [draw opacity=0] (14.29,-6.86) -- (0,0) -- (14.29,6.86) -- cycle    ;
\draw    (564.5,310) -- (641,352.5) ;
\draw    (488,352.5) -- (564.5,310) ;
\draw    (641,437.5) -- (641,352.5) ;
\draw    (488,437.5) -- (488,352.5) ;
\draw    (488,437.5) -- (564.5,480) ;
\draw    (564.5,480) -- (641,437.5) ;
\draw    (488,437.5) -- (564.5,310) ;
\draw    (641,437.5) -- (564.5,310) ;
\draw [color={rgb, 255:red, 74; green, 144; blue, 226 }  ,draw opacity=1 ][line width=0.75]    (379.2,369.5) -- (400.56,354.24) ;
\draw [shift={(403,352.5)}, rotate = 504.46] [fill={rgb, 255:red, 74; green, 144; blue, 226 }  ,fill opacity=1 ][line width=0.08]  [draw opacity=0] (8.93,-4.29) -- (0,0) -- (8.93,4.29) -- cycle    ;
\draw [color={rgb, 255:red, 74; green, 144; blue, 226 }  ,draw opacity=1 ][line width=0.75]    (379.2,369.5) -- (334.25,317.08) ;
\draw [shift={(332.3,314.8)}, rotate = 409.39] [fill={rgb, 255:red, 74; green, 144; blue, 226 }  ,fill opacity=1 ][line width=0.08]  [draw opacity=0] (8.93,-4.29) -- (0,0) -- (8.93,4.29) -- cycle    ;
\draw [color={rgb, 255:red, 74; green, 144; blue, 226 }  ,draw opacity=1 ][line width=0.75]    (379.2,369.5) -- (398.17,421.18) ;
\draw [shift={(399.2,424)}, rotate = 249.85] [fill={rgb, 255:red, 74; green, 144; blue, 226 }  ,fill opacity=1 ][line width=0.08]  [draw opacity=0] (8.93,-4.29) -- (0,0) -- (8.93,4.29) -- cycle    ;
\draw [color={rgb, 255:red, 208; green, 2; blue, 27 }  ,draw opacity=1 ][line width=2.25]    (250,437.5) -- (323.93,314.29) ;
\draw [shift={(326.5,310)}, rotate = 480.96] [fill={rgb, 255:red, 208; green, 2; blue, 27 }  ,fill opacity=1 ][line width=0.08]  [draw opacity=0] (14.29,-6.86) -- (0,0) -- (14.29,6.86) -- cycle    ;
\draw [color={rgb, 255:red, 74; green, 144; blue, 226 }  ,draw opacity=1 ][line width=2.25]    (564.5,310) -- (492.37,350.07) ;
\draw [shift={(488,352.5)}, rotate = 330.95] [fill={rgb, 255:red, 74; green, 144; blue, 226 }  ,fill opacity=1 ][line width=0.08]  [draw opacity=0] (14.29,-6.86) -- (0,0) -- (14.29,6.86) -- cycle    ;
\draw [color={rgb, 255:red, 74; green, 144; blue, 226 }  ,draw opacity=1 ][line width=2.25]    (564.5,310) -- (636.63,350.07) ;
\draw [shift={(641,352.5)}, rotate = 209.05] [fill={rgb, 255:red, 74; green, 144; blue, 226 }  ,fill opacity=1 ][line width=0.08]  [draw opacity=0] (14.29,-6.86) -- (0,0) -- (14.29,6.86) -- cycle    ;
\draw [color={rgb, 255:red, 74; green, 144; blue, 226 }  ,draw opacity=1 ][line width=2.25]    (564.5,480) -- (492.37,439.93) ;
\draw [shift={(488,437.5)}, rotate = 389.05] [fill={rgb, 255:red, 74; green, 144; blue, 226 }  ,fill opacity=1 ][line width=0.08]  [draw opacity=0] (14.29,-6.86) -- (0,0) -- (14.29,6.86) -- cycle    ;
\draw [color={rgb, 255:red, 208; green, 2; blue, 27 }  ,draw opacity=1 ][line width=2.25]    (638.43,433.21) -- (564.5,310) ;
\draw [shift={(641,437.5)}, rotate = 239.04] [fill={rgb, 255:red, 208; green, 2; blue, 27 }  ,fill opacity=1 ][line width=0.08]  [draw opacity=0] (14.29,-6.86) -- (0,0) -- (14.29,6.86) -- cycle    ;
\draw [color={rgb, 255:red, 208; green, 2; blue, 27 }  ,draw opacity=1 ][line width=2.25]    (490.57,433.21) -- (564.5,310) ;
\draw [shift={(488,437.5)}, rotate = 300.96] [fill={rgb, 255:red, 208; green, 2; blue, 27 }  ,fill opacity=1 ][line width=0.08]  [draw opacity=0] (14.29,-6.86) -- (0,0) -- (14.29,6.86) -- cycle    ;
\draw [color={rgb, 255:red, 74; green, 144; blue, 226 }  ,draw opacity=1 ][line width=0.75]    (590,406.9) -- (569.78,331.9) ;
\draw [shift={(569,329)}, rotate = 434.90999999999997] [fill={rgb, 255:red, 74; green, 144; blue, 226 }  ,fill opacity=1 ][line width=0.08]  [draw opacity=0] (8.93,-4.29) -- (0,0) -- (8.93,4.29) -- cycle    ;
\draw [color={rgb, 255:red, 74; green, 144; blue, 226 }  ,draw opacity=1 ][line width=0.75]    (590,406.9) -- (638.43,435.96) ;
\draw [shift={(641,437.5)}, rotate = 210.96] [fill={rgb, 255:red, 74; green, 144; blue, 226 }  ,fill opacity=1 ][line width=0.08]  [draw opacity=0] (8.93,-4.29) -- (0,0) -- (8.93,4.29) -- cycle    ;
\draw [color={rgb, 255:red, 74; green, 144; blue, 226 }  ,draw opacity=1 ][line width=0.75]    (590,406.9) -- (565.49,477.17) ;
\draw [shift={(564.5,480)}, rotate = 289.23] [fill={rgb, 255:red, 74; green, 144; blue, 226 }  ,fill opacity=1 ][line width=0.08]  [draw opacity=0] (8.93,-4.29) -- (0,0) -- (8.93,4.29) -- cycle    ;

\end{tikzpicture}

\end {figure}
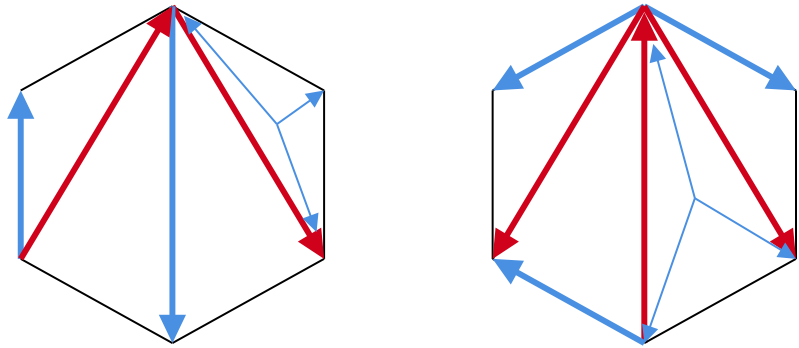

\section{The Recursion Lemmas}
In this section, we prove two lemmas which allow us to write any set of $T$-paths with respect to a fan triangulation in terms of shorter $T$-paths. These lemmas moreover show that the two resolution rules given in Figure \ref{fig:face_resolution} and Figure \ref{fig:edge_resolution} hold on the level of $T$-paths, as well. The first lemma concerns edge-type variables, and the second concerns face-type variables.

First, we show that the resolution of any cluster variable in a fan triangulation comes from immersing an alternating fork-join network, i.e., that (A1) and (A2) are satisfied for the preimage of any $T$-path in a fan triangulation. 

\begin{lem} \label{fork_join_lem}
    Let $\gamma$ be a directed edge or a face in a fan triangulation. If $\gamma$ is a directed edge, then every term in the resolution of $\gamma$ is the image of an alternating fork-join network. If $\gamma$ is a face, then every term in the resolution of $\gamma$ is the image of a triple of alternating fork-join networks. 
\end{lem}

\begin{proof}
    We proceed by induction on $n$, the number of internal diagonals in the triangulation. If $n=1$ then examining the resolution rules in Figures \ref{fig:face_resolution} and \ref{fig:edge_resolution} shows that the result is true in this case. So, suppose $n>1$. 
    
    First consider the case when $\gamma = i \rightarrow j$ is the longest arc. Resolve $\gamma$ with respect to the first edge from the triangulation that it crosses. This creates four diagrams, one of which having all its (directed edge) elements contained in the original triangulation. Of the remaining three diagrams, one will not contain any face variables. The first two colored edges of the latter diagram are contained within the triangulation, and the other edge is a blue edge contained in a subfan with one less triangle. If we resolve this blue edge, ignoring the first two edges, any child of this diagram will come from a fork-join network satisfying conditions (A1) and (A2). One can check that adding back to the leaves of this resolution the first two edges preserves the properties (A1) or (A2). The other two terms from the resolution each have five elements in their diagrams. In either case, the only element not in the original triangulation is a blue face; in fact, it is the same blue face in both diagrams. Now the result follows by an argument similar to the one used for the previous case.
    
    Suppose $\gamma$ is the fan face $ijk$, where $i$ is the fan vertex. Resolve this face over the first diagonal that its non-boundary edge crosses, where we consider this edge oriented from $k$ to $j$. Then the only element in each child diagram which is not contained in the triangulation is a blue fan face, contained in a subfan with one less internal diagonal. Now the result follows much as before. Namely, we first resolve the diagram which only contains the non-initial face just mentioned. By induction, properties (A1) and (A2) are satisfied for every child of this diagram. Now one can check that adding back in the elements from the resolution that aren't this blue face gives a set of $T$-paths that satisfy (A1) and (A2). 
    
    Suppose instead that $\gamma$ is one of the $n$ split faces $ijk$. Let $i$ be the near vertex of the splitting edge, and let $l$ be the other (fan) vertex. Then resolving with respect to the splitting edge gives two child diagrams whose non-initial elements are contained within some smaller fan triangulation. Each child diagram contains at most two elements which are not contained in the initial triangulation; one of these elements is the longest edge in a subfan (or, it is an edge in the triangulation), and the other is a fan face inside a subfan (or, it is a face of the triangulation). The result now follows by an argument similar to the one used in the previous two cases. 
\end{proof}

\begin{rmk} \label{faces_remark} Each of the three remarks below can be seen by induction and the resolution rules in Figure \ref{fig:face_resolution} and Figure \ref{fig:edge_resolution}.

\begin{enumerate}[(1)]
    \item Any diagram obtained from resolving longest edge in a fan triangulation either contains no faces, or contains at most one pair of faces, one red and one blue. 
    
    \item Any diagram obtained from resolving the fan face in a fan triangulation contains precisely one blue face, which is where the $T$-path begins. 
    
    \item Any diagram obtained from resolving a split face in a fan triangulation contains either one blue face (the face where the $T$-path begins), or this blue face along with a pair of faces, one red and one blue.
\end{enumerate}
    
\end{rmk}

Let $i \rightarrow j$ be the longest directed edge not in the triangulation $\Delta$. Suppose the first edge that $i \rightarrow j$ crosses has endpoints labeled $l$ and $k$, as shown in the next figure.

\begin {figure}[h!]
    \centering
    \caption{Set-up for Lemma \ref{rec_1}}
    \label{fig:set_up_1}

\tikzset{every picture/.style={line width=0.75pt}} 

\begin{tikzpicture}[x=0.75pt,y=0.75pt,yscale=-1,xscale=1]

\draw    (413.5,70) -- (490,112.5) ;
\draw    (337,112.5) -- (413.5,70) ;
\draw    (490,197.5) -- (490,112.5) ;
\draw [color={rgb, 255:red, 0; green, 0; blue, 0 }  ,draw opacity=1 ][line width=0.75]    (337,197.5) -- (337,112.5) ;
\draw    (337,197.5) -- (413.5,240) ;
\draw    (413.5,240) -- (490,197.5) ;
\draw [color={rgb, 255:red, 0; green, 0; blue, 0 }  ,draw opacity=1 ][line width=0.75]    (490,197.5) -- (413.5,70) ;
\draw [color={rgb, 255:red, 0; green, 0; blue, 0 }  ,draw opacity=1 ][line width=0.75]    (413.5,240) -- (413.5,70) ;
\draw [color={rgb, 255:red, 0; green, 0; blue, 0 }  ,draw opacity=1 ][line width=0.75]    (337,197.5) -- (413.5,70) ;
\draw [color={rgb, 255:red, 74; green, 144; blue, 226 }  ,draw opacity=1 ][line width=2.25]    (337,112.5) -- (485,112.5) ;
\draw [shift={(490,112.5)}, rotate = 180] [fill={rgb, 255:red, 74; green, 144; blue, 226 }  ,fill opacity=1 ][line width=0.08]  [draw opacity=0] (14.29,-6.86) -- (0,0) -- (14.29,6.86) -- cycle    ;

\draw (324,101) node [anchor=north west][inner sep=0.75pt]    {$i$};
\draw (492,102) node [anchor=north west][inner sep=0.75pt]    {$j$};
\draw (329,200) node [anchor=north west][inner sep=0.75pt]    {$l $};
\draw (409,52) node [anchor=north west][inner sep=0.75pt]    {$k$};

\end{tikzpicture}

\end {figure}

\begin{lem} \label{rec_1}
    \begin{enumerate}[(a)] Fix the fan triangulation $\Delta$ and consider the ``longest'' directed edge $\gamma = i \rightarrow j$ with $l$ and $k$ as in Figure \ref{fig:set_up_1}.  
        \item There is a bijection between the set $\mathbb{T}_{lj}$ and the subset of $\mathbb{T}_{ij}$ consisting of the $T$-paths which begin with the edge $ik$ and do not use the face $ikl.$ In terms of monomials, this bijection is realized as multiplication by $\frac{x_{ik}}{x_{lk}}.$
        
         \item There is a bijection between the set $\mathbb{T}_{kj}$ and the subset of $\mathbb{T}_{ij}$ consisting of the $T$-paths which begin with the edge $il$ and do not use the face $ikl.$ In terms of monomials, this bijection is realized as multiplication by $\frac{x_{il}}{x_{lk}}.$
         
         \item There is a bijection between $\mathbb{T}_{jkl} \bigsqcup \mathbb{T}_{jkl}$ and the subset of $\mathbb{T}_{ij}$ which includes the face $ikl$. 
    \end{enumerate}
\end{lem}

\begin{proof}
    $(a)$ Indeed, a path in $\mathbb{T}_{l j}$ either begins with $l \to k$ or it does not. If it does not, then prepending with $i \to k$
    and $(l \to k)^{-1}$ gives a $T$-path $i \to j$ which does not use the face $ikl$. If, on the other hand, it \textit{does} begin
    with $l \to k$, then removing this first edge and replacing it with $i \to k$ gives a path $i \to j$.

    The inverse of this map is described as follows. If a $T$-path from $i$ to $j$ begins with $i \to k$ and does not include the face $ikl$,
    then its second edge is either contained in the smaller polygon on the other side of edge $kl$, or its second edge is $(l \to k)^{-1}$
    and its third edge is contained in the smaller polygon. In the first case, we replace the first edge $i \to k$ with $l \to k$ to obtain
    a $T$-path from $k$ to $j$. In the second case, we remove the first two edges $(i \to k)(l \to k)^{-1}$, and what remains is a $T$-path
    from $k$ to $j$.

    $(b)$ Similar to $(a)$. 

    $(c)$ If a $T$-path uses face $ikl$, then it begins either $(i \to k)(l \to k)^{-1}(l \to i)(ikl)^{-1}$, or $(i \to l)(k \to l)^{-1}(k \to i) (ikl)^{-1}$.
    The remainder of the $T$-path is a face-path in $\mathbb{T}_{jkl}$. Conversely, any face $T$-path between $j$, $k$, and $l$ can be extended to two
    $T$-paths from $i$ to $j$ which use face $ikl$, corresponding to the two ways to go around the triangle.
\end{proof}

We now give the analogous lemma for face variables. Consider the face $ijk$ which is not cut out by the triangulation $\Delta$. Since we are only considering type A triangulations, there exists an internal diagonal from $\Delta$ which has as one of its endpoints one of the vertices $i$, $j$, or $k$. There are two cases; either $ijk$ is the fan face, or one of the $n$ split faces. We now fix notation in each case. 

Suppose $ijk$ is the fan face, where $i$ is the fan vertex. In this case, every internal diagonal splits this face. Fix one of these splitting edges. This edge will have $i$ as one of its endpoints; call the other endpoint $l$. 

If instead $ijk$ is split, then there is a unique splitting edge; call this internal diagonal $i l$, where $i$ is the near vertex and $l$ is the far (fan) vertex. 

The setup in either case is pictured in Figure \ref{fig:set_up_2} below. 

\newpage

\begin {figure}[h!] 
    \centering
    \caption{Set-up for Lemma \ref{rec_2}}
    \label{fig:set_up_2}

\tikzset{every picture/.style={line width=0.75pt}} 

\begin{tikzpicture}[x=0.75pt,y=0.75pt,yscale=-1,xscale=1]

\draw    (390.5,567) -- (467,609.5) ;
\draw    (314,609.5) -- (390.5,567) ;
\draw    (467,694.5) -- (467,609.5) ;
\draw [color={rgb, 255:red, 0; green, 0; blue, 0 }  ,draw opacity=1 ][line width=0.75]    (314,694.5) -- (314,609.5) ;
\draw    (314,694.5) -- (390.5,737) ;
\draw    (390.5,737) -- (467,694.5) ;
\draw [color={rgb, 255:red, 0; green, 0; blue, 0 }  ,draw opacity=1 ][line width=0.75]    (467,694.5) -- (390.5,567) ;
\draw [color={rgb, 255:red, 0; green, 0; blue, 0 }  ,draw opacity=1 ][line width=0.75]    (390.5,737) -- (390.5,567) ;
\draw [color={rgb, 255:red, 0; green, 0; blue, 0 }  ,draw opacity=1 ][line width=0.75]    (314,694.5) -- (390.5,567) ;
\draw  [fill={rgb, 255:red, 74; green, 144; blue, 226 }  ,fill opacity=0.6 ] (390.5,737) -- (314,609.5) -- (467,609.5) -- cycle ;
\draw    (172.5,568) -- (249,610.5) ;
\draw    (96,610.5) -- (172.5,568) ;
\draw    (249,695.5) -- (249,610.5) ;
\draw [color={rgb, 255:red, 0; green, 0; blue, 0 }  ,draw opacity=1 ][line width=0.75]    (96,695.5) -- (96,610.5) ;
\draw    (96,695.5) -- (172.5,738) ;
\draw    (172.5,738) -- (249,695.5) ;
\draw [color={rgb, 255:red, 0; green, 0; blue, 0 }  ,draw opacity=1 ][line width=0.75]    (249,695.5) -- (172.5,568) ;
\draw [color={rgb, 255:red, 0; green, 0; blue, 0 }  ,draw opacity=1 ][line width=0.75]    (172.5,738) -- (172.5,568) ;
\draw [color={rgb, 255:red, 0; green, 0; blue, 0 }  ,draw opacity=1 ][line width=0.75]    (96,695.5) -- (172.5,568) ;
\draw  [fill={rgb, 255:red, 74; green, 144; blue, 226 }  ,fill opacity=0.6 ] (172.5,568) -- (249,610.5) -- (96,610.5) -- cycle ;

\draw (384,738) node [anchor=north west][inner sep=0.75pt]    {$i$};
\draw (469,599) node [anchor=north west][inner sep=0.75pt]    {$j$};
\draw (385,552) node [anchor=north west][inner sep=0.75pt]    {$\ell $};
\draw (299,600) node [anchor=north west][inner sep=0.75pt]    {$k$};
\draw (167,740) node [anchor=north west][inner sep=0.75pt]    {$l$};
\draw (251,600) node [anchor=north west][inner sep=0.75pt]    {$j$};
\draw (167,553) node [anchor=north west][inner sep=0.75pt]    {$i$};
\draw (81,601) node [anchor=north west][inner sep=0.75pt]    {$k$};

\end{tikzpicture}

\end {figure}

\begin {lem} \label{rec_2}
    Consider the face $ijk$, which is either the fan face or a split face. If $ijk$ is the fan face, let $i,j,k,l$ be as in the left of
    Figure \ref{fig:set_up_2}; if $ijk$ is a split face, let $i,j,k,l$ be as in the right of
    Figure \ref{fig:set_up_2}.
    Then there is a bijection
    $$\mathbb{T}_{ijk} \cong ( \mathbb{T}_{ikl} \times \mathbb{T}_{ij}) \cup ( \mathbb{T}_{ijl} \times \mathbb{T}_{ik})$$
    In terms of monomials, this bijection is realized as multiplication by $x_{i l}$.
\end {lem}

\begin{proof}
   Define a map $( \mathbb{T}_{ikl} \times \mathbb{T}_{ij}) \cup ( \mathbb{T}_{ijl} \times \mathbb{T}_{ik}) \to \mathbb{T}_{ijk}$ by adding
    the red edge $i \to l$. More specifically,
    given a pair $(\varphi,\pi)$, where $\varphi \in \mathbb{T}_{ikl}$ and $\pi \in \mathbb{T}_{ij}$,
    adding the red edge $i \to l$ connects the branch of $\varphi$ ending at $l$ with $\pi$,
    creating a face-path in $\mathbb{T}_{ijk}$. Similarly, if we start with a pair $(\theta,\kappa)$, with
    $\theta \in \mathbb{T}_{ijl}$ and $\kappa \in \mathbb{T}_{ik}$, then adding the red edge $i \to l$ will
    connect the branch of $\theta$ ending at $l$ with $\kappa$, creating a face-path in $\mathbb{T}_{ijk}$.

    We want to show that this map is a bijection. If this is true, then its inverse must be given by
    adding the blue edge $i \to l$. Starting with something in $( \mathbb{T}_{ikl} \times \mathbb{T}_{ij}) \cup ( \mathbb{T}_{ijl} \times \mathbb{T}_{ik})$,
    adding the red edge $i \to l$ and then adding the blue edge $i \to l$ will cancel and the composition is
    the identity. All that needs to be checked is that for any face-path in $\mathbb{T}_{ijk}$, adding the blue edge $i \to l$
    gives a colored diagram which can be interpreted as the superposition of a face path in $\mathbb{T}_{ikl}$ with an
    edge-path in $\mathbb{T}_{ij}$ (or a face in $\mathbb{T}_{ijl}$ and an edge in $\mathbb{T}_{ik}$).

    The edge $il$ cuts the polygon in half, and any face-type $T$-path for triangle $ijk$ must
    ``begin'' with a triangle of $\Delta$ on one of the two halves. Two of the three vertices $ijk$ will already be on this
    half (either $i,k$ or $i,j$), and the third will be on the other half. Let's say we are in the first case, where
    the triangle is on the left half, containing $k$ (the other case will be similar). 

    Since the three branches of the $T$-path start in the left half of the polygon, the two branches going to $i$ and $k$ will
    not leave the left half (by property (T4)). The third branch, which ends at $j$, must cross to the right half along the way. 
    There are two cases we must consider: 
    when this path crosses diagonal $il$ and moves into the right half, it must do so either through vertex $i$ or $l$.

    \underline{$ijk$ is the fan face:} Part (b) of Remark \ref{faces_remark} implies that the branch that we are considering now cannot cross over into the right half of the polygon through vertex $l$. Indeed, if this branch were to cross over at $l$, it would need to do so after using a blue edge by (F4) above. By the remark, the only possibility for the next edge used in the path is a red edge in the triangulation (since no red faces are produced in the resolution of a fan face). But, since we are in a fan triangulation, we are forced to next use the red edge $l \rightarrow i$. This is a contradiction, since we will now leave the left half of the polygon through the fan vertex $i$. 
    
   Therefore we must cross over to the right side through the fan vertex $i$. By (F4), the branch in question arrives at vertex $i$ via a red edge.
    
    There are two possibilities to consider; either the last edge of the path before it arrives at $i$ is the red edge $i \rightarrow l$ , or it is some other red edge which is not the splitting edge we have arbitrarily chosen.  In the latter case, adding the blue edge $i \to l$ makes the beginning part (the part ending at $l$ after the added blue edge)
    a face-type path in $\mathbb{T}_{ikl}$.
    The remainder of the path is an edge-path in $\mathbb{T}_{ij}$. Suppose instead that the last edge of the path
    before it arrives at $i$ is the red edge $i \to l$ (which would again be immediately preceded by a blue edge ending at $l$).  Then our suggested map, which adds the blue edge $i \to l$, will cancel this last edge.
    This beginning part of the path is a face-path in $\mathbb{T}_{ikl}$, and the part of the path after the cancelled edge is an edge-path in $\mathbb{T}_{ij}$.
    Thus, the map indeed gives an element of $\mathbb{T}_{ikl} \times \mathbb{T}_{ij}$. Thus, the claim holds in this case.
    
    \underline{$ijk$ is a split face:} Consider the case where $ijk$ is split, and when we leave the left half at vertex $i$.  Condition (F4) above guarantees that this branch of the path must arrive at $i$ with a red edge. It is easy to see using the resolution rule in Figure \ref{fig:face_resolution} and part (3) in Remark \ref{faces_remark}, that the red edge used to arrive at vertex $i$ must be an edge in the triangulation (as opposed to an edge in a fork). Since we are in a fan triangulation, the only possibility is that the last edge of the path
    before it arrives at $i$ is the red edge $i \to l$ (which would be immediately preceded by a blue edge ending at $l$). As in the case of the fan face, superimposing the blue edge $i \rightarrow j$ cancels this red edge, and the claim again follows much as before. 
    
    Now suppose that the path crosses into the right half at vertex $l$. In this case, the path cannot possibly use the red edge $i \to l$
    (or else the path would cross over at $i$ instead). So adding the blue edge $i \to l$ does not cancel anything, and the blue edge will remain. By (F4), the path arrives at $l$ using a blue edge, and continues after $l$ with a red edge. Thus, the first part of the path (ending at $l$)
    is a face-path in $\mathbb{T}_{ikl}$, and after adding the blue edge $i \to l$, the remaining part of the path becomes an edge-path in $\mathbb{T}_{ij}$. Hence the claim holds for split faces as well. 
    
    This completes the proof.
\end{proof}




    
        


\section{Expansion Formulas in Fan Triangulations}

In this section, we give an expansion formula valid for any arc or face variable in a fan triangulation. 

\begin{thm} Consider the fan triangulation $\Delta$. 

\begin{enumerate}[(a)]
    \item Let $\gamma \colon i \to j$ be a directed edge not in the initial triangulation, with associated cluster variable $x_\gamma$.
            The Laurent expansion of $x_\gamma$ in the initial cluster is
            $$x_\gamma = \sum_{\pi \in \mathbb{T}_{ij}} x_\pi$$
            
            \item Let $\psi$ be a face (with vertices $i,j,k$) not in $\Delta$, with variable $x_\psi$. The Laurent expansion of $x_\psi$ is
            $$x_\psi = \sum_{\varphi \in \mathbb{T}_{ijk}} x_\varphi$$
\end{enumerate}
\end{thm}

\begin{proof}
    \begin{enumerate}[(a)]
        \item We are in the situation pictured in Figure \ref{fig:set_up_1}. Consider the quadrilateral with $ij$ and $k l$ as its two diagonals.
    By construction, triangle $ik l$ is in $\Delta$, but triangle $jk l$ may not be.
    The sequence of 4 mutations which achieves a flip of the diagonal of this quadrilateral gives
    $$x_\gamma = \frac{x_{ik}}{x_{l k}} x_{l j} + \frac{x_{ik} x_{l i}}{x_{ikl} x_{l k}} x_{jkl}
                + \frac{x_{i l}}{x_{kl}} x_{kj} + \frac{x_{i l} x_{ki}}{x_{ikl} x_{kl}} x_{jkl} $$
    The edges $l \to j$ and $k \to j$, and the triangle $jkl$, are within a smaller polygon, and so we may use induction:
    \[ x_\gamma = \frac{x_{ik}}{x_{l k}} \sum_{\pi \in \mathbb{T}_{l j}} x_\pi   
                + \frac{x_{i l}}{x_{kl}} \sum_{\pi \in \mathbb{T}_{kj}} x_\pi 
                + \left( \frac{x_{i l} x_{ki}}{x_{ikl} x_{kl}} + \frac{x_{ik} x_{l i}}{x_{ikl} x_{l k}} \right) \sum_{\pi \in \mathbb{T}_{jkl}} x_\pi\]
    The result now follows from Lemma \ref{rec_1}.
    
    \item The face $\psi = ijk$ is either the fan face shown on the left in Figure \ref{fig:set_up_2}, or a split face, as shown on the firght of Figure \ref{fig:set_up_2}.

    In either case, the flip mutation sequence gives
    \[ x_\psi = \frac{x_{ij}}{x_{il}} x_{ikl} + \frac{x_{ik}}{x_{il}} x_{ijl} \]
    In the first term, triangle $ikl$ and edge $i \to j$ are within smaller polygons. In the second term, triangle $ijl$ and edge $i \to k$
    are within smaller polygons. So we may use induction:
    \begin {align*} 
        x_\psi &= \frac{1}{x_{il}} \left( x_{ij} x_{ikl} + x_{ik} x_{ijl} \right) \\
                 &= \frac{1}{x_{il}} \left( \sum_{\substack{\pi \in \mathbb{T}_{ij}\\\varphi \in \mathbb{T}_{ikl}}} x_\pi x_\varphi 
                                             + \sum_{\substack{\kappa \in \mathbb{T}_{ik}\\\theta \in \mathbb{T}_{ijl}}} x_\kappa x_\theta \right)
    \end {align*}

    The result now follows from Lemma \ref{rec_2}.
    
    \end{enumerate}
\end{proof}

\section{Some Expansion Posets}

In this final section, we describe the poset structure on the Laurent monomials in the expansion of the fan face, and for the ``longest'' edge, in a fan triangulation. Covering relations are again multiplication by the $\widehat{y}$ variables defined in  \cite{fomin2007cluster}. The corresponding colored $\text{SL}_3$ diagrams are shown below. 

\begin {figure}[h!]
    \centering
    \caption{Covering relations}
    \label{fig:cov_rel}
   \begin{tikzpicture}[x=0.75pt,y=0.75pt,yscale=-1,xscale=1]

\draw [color={rgb, 255:red, 74; green, 144; blue, 226 }  ,draw opacity=1 ]   (531,299) -- (568.35,242.5) ;
\draw [shift={(570,240)}, rotate = 483.47] [fill={rgb, 255:red, 74; green, 144; blue, 226 }  ,fill opacity=1 ][line width=0.08]  [draw opacity=0] (8.93,-4.29) -- (0,0) -- (8.93,4.29) -- cycle    ;
\draw [color={rgb, 255:red, 74; green, 144; blue, 226 }  ,draw opacity=1 ]   (571,241) -- (608.35,297.5) ;
\draw [shift={(610,300)}, rotate = 236.53] [fill={rgb, 255:red, 74; green, 144; blue, 226 }  ,fill opacity=1 ][line width=0.08]  [draw opacity=0] (8.93,-4.29) -- (0,0) -- (8.93,4.29) -- cycle    ;
\draw [color={rgb, 255:red, 74; green, 144; blue, 226 }  ,draw opacity=1 ]   (608,300) -- (533,300) ;
\draw [shift={(530,300)}, rotate = 360] [fill={rgb, 255:red, 74; green, 144; blue, 226 }  ,fill opacity=1 ][line width=0.08]  [draw opacity=0] (8.93,-4.29) -- (0,0) -- (8.93,4.29) -- cycle    ;
\draw [color={rgb, 255:red, 208; green, 2; blue, 27 }  ,draw opacity=1 ]   (576.66,242.5) -- (615,300) ;
\draw [shift={(575,240)}, rotate = 56.31] [fill={rgb, 255:red, 208; green, 2; blue, 27 }  ,fill opacity=1 ][line width=0.08]  [draw opacity=0] (8.93,-4.29) -- (0,0) -- (8.93,4.29) -- cycle    ;
\draw [color={rgb, 255:red, 208; green, 2; blue, 27 }  ,draw opacity=1 ]   (526.66,297.5) -- (565,240) ;
\draw [shift={(525,300)}, rotate = 303.69] [fill={rgb, 255:red, 208; green, 2; blue, 27 }  ,fill opacity=1 ][line width=0.08]  [draw opacity=0] (8.93,-4.29) -- (0,0) -- (8.93,4.29) -- cycle    ;
\draw [color={rgb, 255:red, 208; green, 2; blue, 27 }  ,draw opacity=1 ]   (607,304) -- (530,304) ;
\draw [shift={(610,304)}, rotate = 180] [fill={rgb, 255:red, 208; green, 2; blue, 27 }  ,fill opacity=1 ][line width=0.08]  [draw opacity=0] (8.93,-4.29) -- (0,0) -- (8.93,4.29) -- cycle    ;
\draw    (730,230) -- (770,270) ;
\draw    (770,270) -- (730,310) ;
\draw    (730,230) -- (690,270) ;
\draw    (690,270) -- (730,310) ;
\draw    (730,230) -- (730,310) ;
\draw  [fill={rgb, 255:red, 208; green, 2; blue, 27 }  ,fill opacity=0.6 ] (730,230) -- (730,310) -- (690,270) -- cycle ;
\draw  [fill={rgb, 255:red, 74; green, 144; blue, 226 }  ,fill opacity=0.6 ] (730,310) -- (730,230) -- (770,270) -- cycle ;
\draw [color={rgb, 255:red, 74; green, 144; blue, 226 }  ,draw opacity=1 ][line width=1.5]    (692.83,272.83) -- (730,310) ;
\draw [shift={(690,270)}, rotate = 45] [fill={rgb, 255:red, 74; green, 144; blue, 226 }  ,fill opacity=1 ][line width=0.08]  [draw opacity=0] (11.61,-5.58) -- (0,0) -- (11.61,5.58) -- cycle    ;
\draw [color={rgb, 255:red, 208; green, 2; blue, 27 }  ,draw opacity=1 ][line width=1.5]    (730,310) -- (767.17,272.83) ;
\draw [shift={(770,270)}, rotate = 495] [fill={rgb, 255:red, 208; green, 2; blue, 27 }  ,fill opacity=1 ][line width=0.08]  [draw opacity=0] (11.61,-5.58) -- (0,0) -- (11.61,5.58) -- cycle    ;

\draw (517,298) node [anchor=north west][inner sep=0.75pt]    {$i$};
\draw (610,297) node [anchor=north west][inner sep=0.75pt]    {$j$};
\draw (565,223) node [anchor=north west][inner sep=0.75pt]    {$k$};
\draw (535,320) node [anchor=north west][inner sep=0.75pt]    {$face\ type$};
\draw (697,320) node [anchor=north west][inner sep=0.75pt]    {$edge\ type$};

\end{tikzpicture}
\end {figure}
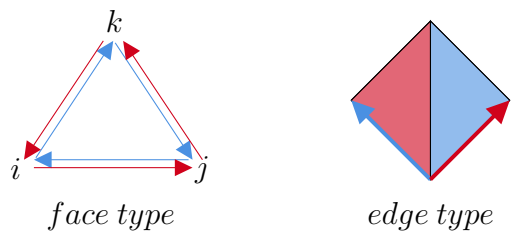

Suppose the triangulation $\Delta$ has $n$ internal diagonals.

\begin{prop} Consider the fan face $\gamma$. Then the poset $\mathbb{T}_{\gamma}$ is isomorphic to a chain with $n+1$ vertices. 
\end{prop}

\begin{proof}
    We induct on the number $n$ of internal diagonals. By the resolution rule given in Figure \ref{fig:face_resolution} the result is true for $n=1$. So suppose $n>1$. Label the vertices of the face $\gamma$ by $i,j,$ and $k$, where $k$ is the fan vertex of $\Delta$. Let $kl$ be the first internal diagonal of $\Delta$ that the directed edge $i \rightarrow j$ crosses. 
    
     Resolve $\gamma$ with respect to the arc $i l$ to obtain two diagrams (see Figure \ref{fig:hex_fan_face} for an illustration of the two resulting diagrams for the fan face in a hexagon).
     
       \begin {figure}[h!]
    \centering
    \caption{The first step in the resolution of the fan face in a hexagon}
    \label{fig:hex_fan_face}

\tikzset{every picture/.style={line width=0.75pt}} 

\begin{tikzpicture}[x=0.75pt,y=0.75pt,yscale=-1,xscale=1]

\draw    (226.58,2160) -- (303.08,2202.5) ;
\draw    (150.08,2202.5) -- (226.58,2160) ;
\draw    (303.08,2287.5) -- (303.08,2202.5) ;
\draw [color={rgb, 255:red, 0; green, 0; blue, 0 }  ,draw opacity=1 ][line width=0.75]    (150.08,2287.5) -- (150.08,2202.5) ;
\draw    (150.08,2287.5) -- (226.58,2330) ;
\draw    (226.58,2330) -- (303.08,2287.5) ;
\draw [color={rgb, 255:red, 0; green, 0; blue, 0 }  ,draw opacity=1 ][line width=0.75]    (303.08,2287.5) -- (226.58,2160) ;
\draw [color={rgb, 255:red, 0; green, 0; blue, 0 }  ,draw opacity=1 ][line width=0.75]    (226.58,2330) -- (226.58,2160) ;
\draw [color={rgb, 255:red, 0; green, 0; blue, 0 }  ,draw opacity=1 ][line width=0.75]    (150.08,2287.5) -- (226.58,2160) ;
\draw  [fill={rgb, 255:red, 74; green, 144; blue, 226 }  ,fill opacity=0.6 ] (150.06,2202.54) -- (226.58,2160) -- (150,2287.45) -- cycle ;
\draw    (429.66,2160) -- (506.16,2202.5) ;
\draw    (353.16,2202.5) -- (429.66,2160) ;
\draw    (506.16,2287.5) -- (506.16,2202.5) ;
\draw [color={rgb, 255:red, 0; green, 0; blue, 0 }  ,draw opacity=1 ][line width=0.75]    (353.16,2287.5) -- (353.16,2202.5) ;
\draw    (353.16,2287.5) -- (429.66,2330) ;
\draw    (429.66,2330) -- (506.16,2287.5) ;
\draw [color={rgb, 255:red, 0; green, 0; blue, 0 }  ,draw opacity=1 ][line width=0.75]    (506.16,2287.5) -- (429.66,2160) ;
\draw [color={rgb, 255:red, 0; green, 0; blue, 0 }  ,draw opacity=1 ][line width=0.75]    (429.66,2330) -- (429.66,2160) ;
\draw [color={rgb, 255:red, 0; green, 0; blue, 0 }  ,draw opacity=1 ][line width=0.75]    (353.16,2287.5) -- (429.66,2160) ;
\draw  [fill={rgb, 255:red, 74; green, 144; blue, 226 }  ,fill opacity=0.6 ] (354.89,2289.38) -- (427.15,2159.03) -- (504.87,2202.12) -- cycle ;
\draw [color={rgb, 255:red, 208; green, 2; blue, 27 }  ,draw opacity=1 ][line width=2.25]    (226.58,2160) -- (152.58,2283.17) ;
\draw [shift={(150,2287.45)}, rotate = 301] [fill={rgb, 255:red, 208; green, 2; blue, 27 }  ,fill opacity=1 ][line width=0.08]  [draw opacity=0] (14.29,-6.86) -- (0,0) -- (14.29,6.86) -- cycle    ;
\draw [color={rgb, 255:red, 74; green, 144; blue, 226 }  ,draw opacity=1 ][line width=2.25]    (226.58,2160) -- (298.71,2200.07) ;
\draw [shift={(303.08,2202.5)}, rotate = 209.05] [fill={rgb, 255:red, 74; green, 144; blue, 226 }  ,fill opacity=1 ][line width=0.08]  [draw opacity=0] (14.29,-6.86) -- (0,0) -- (14.29,6.86) -- cycle    ;
\draw [color={rgb, 255:red, 74; green, 144; blue, 226 }  ,draw opacity=1 ][line width=2.25]    (429.66,2160) -- (357.53,2200.07) ;
\draw [shift={(353.16,2202.5)}, rotate = 330.95] [fill={rgb, 255:red, 74; green, 144; blue, 226 }  ,fill opacity=1 ][line width=0.08]  [draw opacity=0] (14.29,-6.86) -- (0,0) -- (14.29,6.86) -- cycle    ;
\draw [color={rgb, 255:red, 208; green, 2; blue, 27 }  ,draw opacity=1 ][line width=2.25]    (428.66,2159) -- (355.57,2286.12) ;
\draw [shift={(353.08,2290.45)}, rotate = 299.9] [fill={rgb, 255:red, 208; green, 2; blue, 27 }  ,fill opacity=1 ][line width=0.08]  [draw opacity=0] (14.29,-6.86) -- (0,0) -- (14.29,6.86) -- cycle    ;

\draw (319,2229) node [anchor=north west][inner sep=0.75pt]    {$+$};

\end{tikzpicture}

\end {figure}
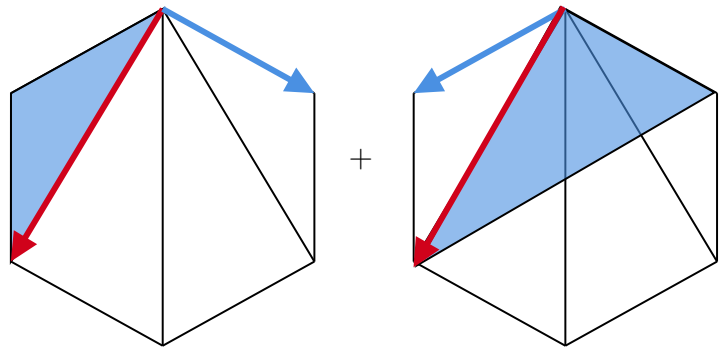
     Note that one of these diagrams has all of its elements contained in the triangulation $\Delta$. Call this diagram $D$. All of the elements in the other diagram are part of the fan triangulation, except one blue face, which is contained in a subfan triangulation with $n-1$ internal diagonals.

     Resolve the latter diagram with respect to the arc whose non-fan endpoint is ``closest'' to $l$. This produces two terms, with one of these terms having every element contained in the triangulation. This latter term is related to $D$ by the first covering relation shown in Figure \ref{fig:cov_rel}. Now the result follows by induction.
\end{proof}

\begin{ex}
    Figure \ref{fig:poset_fan_face} shows the expansion poset of the fan face in a hexagon. The minimal element is the leftmost element shown, and the maximal element is the rightmost element shown.

    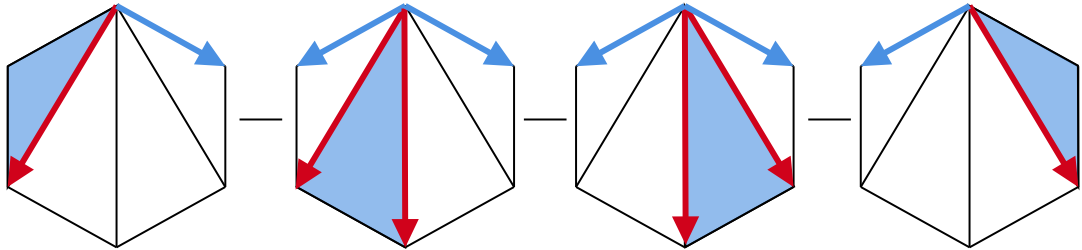
\begin {figure}[h!]
    \centering
    \caption{Expansion poset of a fan face}
    \label{fig:poset_fan_face}

\tikzset{every picture/.style={line width=0.75pt}} 

\begin{tikzpicture}[x=0.75pt,y=0.75pt,yscale=-1,xscale=1]

\draw    (124.91,1410) -- (179.77,1440.47) ;
\draw    (70.06,1440.47) -- (124.91,1410) ;
\draw    (179.77,1501.42) -- (179.77,1440.47) ;
\draw [color={rgb, 255:red, 0; green, 0; blue, 0 }  ,draw opacity=1 ][line width=0.75]    (70.06,1501.42) -- (70.06,1440.47) ;
\draw    (70.06,1501.42) -- (124.91,1531.9) ;
\draw    (124.91,1531.9) -- (179.77,1501.42) ;
\draw [color={rgb, 255:red, 0; green, 0; blue, 0 }  ,draw opacity=1 ][line width=0.75]    (179.77,1501.42) -- (124.91,1410) ;
\draw [color={rgb, 255:red, 0; green, 0; blue, 0 }  ,draw opacity=1 ][line width=0.75]    (124.91,1531.9) -- (124.91,1410) ;
\draw [color={rgb, 255:red, 0; green, 0; blue, 0 }  ,draw opacity=1 ][line width=0.75]    (70.06,1501.42) -- (124.91,1410) ;
\draw  [fill={rgb, 255:red, 74; green, 144; blue, 226 }  ,fill opacity=0.6 ] (70.04,1440.5) -- (124.91,1410) -- (70,1501.39) -- cycle ;
\draw    (555.11,1410) -- (609.96,1440.47) ;
\draw    (500.25,1440.47) -- (555.11,1410) ;
\draw    (609.96,1501.42) -- (609.96,1440.47) ;
\draw [color={rgb, 255:red, 0; green, 0; blue, 0 }  ,draw opacity=1 ][line width=0.75]    (500.25,1501.42) -- (500.25,1440.47) ;
\draw    (500.25,1501.42) -- (555.11,1531.9) ;
\draw    (555.11,1531.9) -- (609.96,1501.42) ;
\draw [color={rgb, 255:red, 0; green, 0; blue, 0 }  ,draw opacity=1 ][line width=0.75]    (609.96,1501.42) -- (555.11,1410) ;
\draw [color={rgb, 255:red, 0; green, 0; blue, 0 }  ,draw opacity=1 ][line width=0.75]    (555.11,1531.9) -- (555.11,1410) ;
\draw [color={rgb, 255:red, 0; green, 0; blue, 0 }  ,draw opacity=1 ][line width=0.75]    (500.25,1501.42) -- (555.11,1410) ;
\draw  [fill={rgb, 255:red, 74; green, 144; blue, 226 }  ,fill opacity=0.6 ] (609.76,1440.14) -- (610,1501.36) -- (555.11,1410) -- cycle ;
\draw    (270.53,1410) -- (325.39,1440.47) ;
\draw    (215.68,1440.47) -- (270.53,1410) ;
\draw    (325.39,1501.42) -- (325.39,1440.47) ;
\draw [color={rgb, 255:red, 0; green, 0; blue, 0 }  ,draw opacity=1 ][line width=0.75]    (215.68,1501.42) -- (215.68,1440.47) ;
\draw    (215.68,1501.42) -- (270.53,1531.9) ;
\draw    (270.53,1531.9) -- (325.39,1501.42) ;
\draw [color={rgb, 255:red, 0; green, 0; blue, 0 }  ,draw opacity=1 ][line width=0.75]    (325.39,1501.42) -- (270.53,1410) ;
\draw [color={rgb, 255:red, 0; green, 0; blue, 0 }  ,draw opacity=1 ][line width=0.75]    (270.53,1531.9) -- (270.53,1410) ;
\draw [color={rgb, 255:red, 0; green, 0; blue, 0 }  ,draw opacity=1 ][line width=0.75]    (215.68,1501.42) -- (270.53,1410) ;
\draw  [fill={rgb, 255:red, 74; green, 144; blue, 226 }  ,fill opacity=0.6 ] (270.1,1411.34) -- (270.53,1531.9) -- (215.71,1501.39) -- cycle ;
\draw    (411.51,1410) -- (466.37,1440.47) ;
\draw    (356.66,1440.47) -- (411.51,1410) ;
\draw    (466.37,1501.42) -- (466.37,1440.47) ;
\draw [color={rgb, 255:red, 0; green, 0; blue, 0 }  ,draw opacity=1 ][line width=0.75]    (356.66,1501.42) -- (356.66,1440.47) ;
\draw    (356.66,1501.42) -- (411.51,1531.9) ;
\draw    (411.51,1531.9) -- (466.37,1501.42) ;
\draw [color={rgb, 255:red, 0; green, 0; blue, 0 }  ,draw opacity=1 ][line width=0.75]    (466.37,1501.42) -- (411.51,1410) ;
\draw [color={rgb, 255:red, 0; green, 0; blue, 0 }  ,draw opacity=1 ][line width=0.75]    (411.51,1531.9) -- (411.51,1410) ;
\draw [color={rgb, 255:red, 0; green, 0; blue, 0 }  ,draw opacity=1 ][line width=0.75]    (356.66,1501.42) -- (411.51,1410) ;
\draw  [fill={rgb, 255:red, 74; green, 144; blue, 226 }  ,fill opacity=0.6 ] (411.51,1410) -- (411.77,1531.79) -- (466.59,1501.27) -- cycle ;
\draw [color={rgb, 255:red, 208; green, 2; blue, 27 }  ,draw opacity=1 ][line width=2.25]    (124.91,1410) -- (72.58,1497.1) ;
\draw [shift={(70,1501.39)}, rotate = 301] [fill={rgb, 255:red, 208; green, 2; blue, 27 }  ,fill opacity=1 ][line width=0.08]  [draw opacity=0] (14.29,-6.86) -- (0,0) -- (14.29,6.86) -- cycle    ;
\draw [color={rgb, 255:red, 74; green, 144; blue, 226 }  ,draw opacity=1 ][line width=2.25]    (124.91,1410) -- (175.4,1438.05) ;
\draw [shift={(179.77,1440.47)}, rotate = 209.05] [fill={rgb, 255:red, 74; green, 144; blue, 226 }  ,fill opacity=1 ][line width=0.08]  [draw opacity=0] (14.29,-6.86) -- (0,0) -- (14.29,6.86) -- cycle    ;
\draw [color={rgb, 255:red, 208; green, 2; blue, 27 }  ,draw opacity=1 ][line width=2.25]    (270.1,1411.34) -- (217.76,1498.45) ;
\draw [shift={(215.19,1502.73)}, rotate = 301] [fill={rgb, 255:red, 208; green, 2; blue, 27 }  ,fill opacity=1 ][line width=0.08]  [draw opacity=0] (14.29,-6.86) -- (0,0) -- (14.29,6.86) -- cycle    ;
\draw [color={rgb, 255:red, 74; green, 144; blue, 226 }  ,draw opacity=1 ][line width=2.25]    (270.53,1410) -- (220.05,1438.05) ;
\draw [shift={(215.68,1440.47)}, rotate = 330.95] [fill={rgb, 255:red, 74; green, 144; blue, 226 }  ,fill opacity=1 ][line width=0.08]  [draw opacity=0] (14.29,-6.86) -- (0,0) -- (14.29,6.86) -- cycle    ;
\draw [color={rgb, 255:red, 208; green, 2; blue, 27 }  ,draw opacity=1 ][line width=2.25]    (270.1,1411.34) -- (270.51,1526.9) ;
\draw [shift={(270.53,1531.9)}, rotate = 269.8] [fill={rgb, 255:red, 208; green, 2; blue, 27 }  ,fill opacity=1 ][line width=0.08]  [draw opacity=0] (14.29,-6.86) -- (0,0) -- (14.29,6.86) -- cycle    ;
\draw [color={rgb, 255:red, 74; green, 144; blue, 226 }  ,draw opacity=1 ][line width=2.25]    (270.53,1410) -- (321.02,1438.05) ;
\draw [shift={(325.39,1440.47)}, rotate = 209.05] [fill={rgb, 255:red, 74; green, 144; blue, 226 }  ,fill opacity=1 ][line width=0.08]  [draw opacity=0] (14.29,-6.86) -- (0,0) -- (14.29,6.86) -- cycle    ;
\draw [color={rgb, 255:red, 208; green, 2; blue, 27 }  ,draw opacity=1 ][line width=2.25]    (411.51,1410) -- (463.79,1497.14) ;
\draw [shift={(466.37,1501.42)}, rotate = 239.04] [fill={rgb, 255:red, 208; green, 2; blue, 27 }  ,fill opacity=1 ][line width=0.08]  [draw opacity=0] (14.29,-6.86) -- (0,0) -- (14.29,6.86) -- cycle    ;
\draw [color={rgb, 255:red, 208; green, 2; blue, 27 }  ,draw opacity=1 ][line width=2.25]    (411.51,1410) -- (411.93,1525.56) ;
\draw [shift={(411.94,1530.56)}, rotate = 269.8] [fill={rgb, 255:red, 208; green, 2; blue, 27 }  ,fill opacity=1 ][line width=0.08]  [draw opacity=0] (14.29,-6.86) -- (0,0) -- (14.29,6.86) -- cycle    ;
\draw [color={rgb, 255:red, 74; green, 144; blue, 226 }  ,draw opacity=1 ][line width=2.25]    (411.51,1410) -- (361.03,1438.05) ;
\draw [shift={(356.66,1440.47)}, rotate = 330.95] [fill={rgb, 255:red, 74; green, 144; blue, 226 }  ,fill opacity=1 ][line width=0.08]  [draw opacity=0] (14.29,-6.86) -- (0,0) -- (14.29,6.86) -- cycle    ;
\draw [color={rgb, 255:red, 74; green, 144; blue, 226 }  ,draw opacity=1 ][line width=2.25]    (411.51,1410) -- (462,1438.05) ;
\draw [shift={(466.37,1440.47)}, rotate = 209.05] [fill={rgb, 255:red, 74; green, 144; blue, 226 }  ,fill opacity=1 ][line width=0.08]  [draw opacity=0] (14.29,-6.86) -- (0,0) -- (14.29,6.86) -- cycle    ;
\draw [color={rgb, 255:red, 208; green, 2; blue, 27 }  ,draw opacity=1 ][line width=2.25]    (555.11,1410) -- (607.39,1497.14) ;
\draw [shift={(609.96,1501.42)}, rotate = 239.04] [fill={rgb, 255:red, 208; green, 2; blue, 27 }  ,fill opacity=1 ][line width=0.08]  [draw opacity=0] (14.29,-6.86) -- (0,0) -- (14.29,6.86) -- cycle    ;
\draw [color={rgb, 255:red, 74; green, 144; blue, 226 }  ,draw opacity=1 ][line width=2.25]    (555.11,1410) -- (504.62,1438.05) ;
\draw [shift={(500.25,1440.47)}, rotate = 330.95] [fill={rgb, 255:red, 74; green, 144; blue, 226 }  ,fill opacity=1 ][line width=0.08]  [draw opacity=0] (14.29,-6.86) -- (0,0) -- (14.29,6.86) -- cycle    ;
\draw    (186.94,1467.36) -- (208.45,1467.36) ;
\draw    (330.35,1467.36) -- (351.86,1467.36) ;
\draw    (473.76,1467.36) -- (495.27,1467.36) ;

\end{tikzpicture}

\end {figure}

\end{ex}

\begin{rmk}
There is a bijection between the terms in the expansion poset of a fan face, and the terms in the ``$\text{SL}_2$'' $T$-path poset obtained from expanding the longest (undirected) edge in the same triangulation. Namely, given a term in the fan face poset, deleting the unique blue face and replacing it with the unique blue boundary edge boarding the deleted face, which is not adjacent to the fan vertex, gives the bijection. An example of this is indicated below.

\begin {figure}[h!]
    \centering
    \caption{The bijection between a fan face poset and the corresponding $\text{SL}_2$ $T$-path poset}
    \label{fig:SL_3_2_bij}

\tikzset{every picture/.style={line width=0.75pt}} 

\begin{tikzpicture}[x=0.75pt,y=0.75pt,yscale=-1,xscale=1]

\draw    (213.51,1410) -- (250.43,1430.51) ;
\draw    (176.59,1430.51) -- (213.51,1410) ;
\draw    (250.43,1471.53) -- (250.43,1430.51) ;
\draw [color={rgb, 255:red, 0; green, 0; blue, 0 }  ,draw opacity=1 ][line width=0.75]    (176.59,1471.53) -- (176.59,1430.51) ;
\draw    (176.59,1471.53) -- (213.51,1492.05) ;
\draw    (213.51,1492.05) -- (250.43,1471.53) ;
\draw [color={rgb, 255:red, 0; green, 0; blue, 0 }  ,draw opacity=1 ][line width=0.75]    (250.43,1471.53) -- (213.51,1410) ;
\draw [color={rgb, 255:red, 0; green, 0; blue, 0 }  ,draw opacity=1 ][line width=0.75]    (213.51,1492.05) -- (213.51,1410) ;
\draw [color={rgb, 255:red, 0; green, 0; blue, 0 }  ,draw opacity=1 ][line width=0.75]    (176.59,1471.53) -- (213.51,1410) ;
\draw  [fill={rgb, 255:red, 74; green, 144; blue, 226 }  ,fill opacity=0.6 ] (176.58,1430.53) -- (213.51,1410) -- (176.55,1471.51) -- cycle ;
\draw    (503.05,1410) -- (539.97,1430.51) ;
\draw    (466.13,1430.51) -- (503.05,1410) ;
\draw    (539.97,1471.53) -- (539.97,1430.51) ;
\draw [color={rgb, 255:red, 0; green, 0; blue, 0 }  ,draw opacity=1 ][line width=0.75]    (466.13,1471.53) -- (466.13,1430.51) ;
\draw    (466.13,1471.53) -- (503.05,1492.05) ;
\draw    (503.05,1492.05) -- (539.97,1471.53) ;
\draw [color={rgb, 255:red, 0; green, 0; blue, 0 }  ,draw opacity=1 ][line width=0.75]    (539.97,1471.53) -- (503.05,1410) ;
\draw [color={rgb, 255:red, 0; green, 0; blue, 0 }  ,draw opacity=1 ][line width=0.75]    (503.05,1492.05) -- (503.05,1410) ;
\draw [color={rgb, 255:red, 0; green, 0; blue, 0 }  ,draw opacity=1 ][line width=0.75]    (466.13,1471.53) -- (503.05,1410) ;
\draw  [fill={rgb, 255:red, 74; green, 144; blue, 226 }  ,fill opacity=0.6 ] (539.84,1430.29) -- (540,1471.49) -- (503.05,1410) -- cycle ;
\draw    (311.52,1410) -- (348.44,1430.51) ;
\draw    (274.6,1430.51) -- (311.52,1410) ;
\draw    (348.44,1471.53) -- (348.44,1430.51) ;
\draw [color={rgb, 255:red, 0; green, 0; blue, 0 }  ,draw opacity=1 ][line width=0.75]    (274.6,1471.53) -- (274.6,1430.51) ;
\draw    (274.6,1471.53) -- (311.52,1492.05) ;
\draw    (311.52,1492.05) -- (348.44,1471.53) ;
\draw [color={rgb, 255:red, 0; green, 0; blue, 0 }  ,draw opacity=1 ][line width=0.75]    (348.44,1471.53) -- (311.52,1410) ;
\draw [color={rgb, 255:red, 0; green, 0; blue, 0 }  ,draw opacity=1 ][line width=0.75]    (311.52,1492.05) -- (311.52,1410) ;
\draw [color={rgb, 255:red, 0; green, 0; blue, 0 }  ,draw opacity=1 ][line width=0.75]    (274.6,1471.53) -- (311.52,1410) ;
\draw  [fill={rgb, 255:red, 74; green, 144; blue, 226 }  ,fill opacity=0.6 ] (311.23,1410.9) -- (311.52,1492.05) -- (274.62,1471.51) -- cycle ;
\draw    (406.41,1410) -- (443.33,1430.51) ;
\draw    (369.49,1430.51) -- (406.41,1410) ;
\draw    (443.33,1471.53) -- (443.33,1430.51) ;
\draw [color={rgb, 255:red, 0; green, 0; blue, 0 }  ,draw opacity=1 ][line width=0.75]    (369.49,1471.53) -- (369.49,1430.51) ;
\draw    (369.49,1471.53) -- (406.41,1492.05) ;
\draw    (406.41,1492.05) -- (443.33,1471.53) ;
\draw [color={rgb, 255:red, 0; green, 0; blue, 0 }  ,draw opacity=1 ][line width=0.75]    (443.33,1471.53) -- (406.41,1410) ;
\draw [color={rgb, 255:red, 0; green, 0; blue, 0 }  ,draw opacity=1 ][line width=0.75]    (406.41,1492.05) -- (406.41,1410) ;
\draw [color={rgb, 255:red, 0; green, 0; blue, 0 }  ,draw opacity=1 ][line width=0.75]    (369.49,1471.53) -- (406.41,1410) ;
\draw  [fill={rgb, 255:red, 74; green, 144; blue, 226 }  ,fill opacity=0.6 ] (406.41,1410) -- (406.58,1491.97) -- (443.48,1471.43) -- cycle ;
\draw [color={rgb, 255:red, 208; green, 2; blue, 27 }  ,draw opacity=1 ][line width=2.25]    (213.51,1410) -- (179.12,1467.22) ;
\draw [shift={(176.55,1471.51)}, rotate = 301] [fill={rgb, 255:red, 208; green, 2; blue, 27 }  ,fill opacity=1 ][line width=0.08]  [draw opacity=0] (14.29,-6.86) -- (0,0) -- (14.29,6.86) -- cycle    ;
\draw [color={rgb, 255:red, 74; green, 144; blue, 226 }  ,draw opacity=1 ][line width=2.25]    (213.51,1410) -- (246.06,1428.08) ;
\draw [shift={(250.43,1430.51)}, rotate = 209.05] [fill={rgb, 255:red, 74; green, 144; blue, 226 }  ,fill opacity=1 ][line width=0.08]  [draw opacity=0] (14.29,-6.86) -- (0,0) -- (14.29,6.86) -- cycle    ;
\draw [color={rgb, 255:red, 208; green, 2; blue, 27 }  ,draw opacity=1 ][line width=2.25]    (311.23,1410.9) -- (276.84,1468.13) ;
\draw [shift={(274.27,1472.41)}, rotate = 301] [fill={rgb, 255:red, 208; green, 2; blue, 27 }  ,fill opacity=1 ][line width=0.08]  [draw opacity=0] (14.29,-6.86) -- (0,0) -- (14.29,6.86) -- cycle    ;
\draw [color={rgb, 255:red, 74; green, 144; blue, 226 }  ,draw opacity=1 ][line width=2.25]    (311.52,1410) -- (278.97,1428.08) ;
\draw [shift={(274.6,1430.51)}, rotate = 330.95] [fill={rgb, 255:red, 74; green, 144; blue, 226 }  ,fill opacity=1 ][line width=0.08]  [draw opacity=0] (14.29,-6.86) -- (0,0) -- (14.29,6.86) -- cycle    ;
\draw [color={rgb, 255:red, 208; green, 2; blue, 27 }  ,draw opacity=1 ][line width=2.25]    (311.23,1410.9) -- (311.5,1487.05) ;
\draw [shift={(311.52,1492.05)}, rotate = 269.8] [fill={rgb, 255:red, 208; green, 2; blue, 27 }  ,fill opacity=1 ][line width=0.08]  [draw opacity=0] (14.29,-6.86) -- (0,0) -- (14.29,6.86) -- cycle    ;
\draw [color={rgb, 255:red, 74; green, 144; blue, 226 }  ,draw opacity=1 ][line width=2.25]    (311.52,1410) -- (344.07,1428.08) ;
\draw [shift={(348.44,1430.51)}, rotate = 209.05] [fill={rgb, 255:red, 74; green, 144; blue, 226 }  ,fill opacity=1 ][line width=0.08]  [draw opacity=0] (14.29,-6.86) -- (0,0) -- (14.29,6.86) -- cycle    ;
\draw [color={rgb, 255:red, 208; green, 2; blue, 27 }  ,draw opacity=1 ][line width=2.25]    (406.41,1410) -- (440.75,1467.25) ;
\draw [shift={(443.33,1471.53)}, rotate = 239.04] [fill={rgb, 255:red, 208; green, 2; blue, 27 }  ,fill opacity=1 ][line width=0.08]  [draw opacity=0] (14.29,-6.86) -- (0,0) -- (14.29,6.86) -- cycle    ;
\draw [color={rgb, 255:red, 208; green, 2; blue, 27 }  ,draw opacity=1 ][line width=2.25]    (406.41,1410) -- (406.68,1486.14) ;
\draw [shift={(406.7,1491.14)}, rotate = 269.8] [fill={rgb, 255:red, 208; green, 2; blue, 27 }  ,fill opacity=1 ][line width=0.08]  [draw opacity=0] (14.29,-6.86) -- (0,0) -- (14.29,6.86) -- cycle    ;
\draw [color={rgb, 255:red, 74; green, 144; blue, 226 }  ,draw opacity=1 ][line width=2.25]    (406.41,1410) -- (373.86,1428.08) ;
\draw [shift={(369.49,1430.51)}, rotate = 330.95] [fill={rgb, 255:red, 74; green, 144; blue, 226 }  ,fill opacity=1 ][line width=0.08]  [draw opacity=0] (14.29,-6.86) -- (0,0) -- (14.29,6.86) -- cycle    ;
\draw [color={rgb, 255:red, 74; green, 144; blue, 226 }  ,draw opacity=1 ][line width=2.25]    (406.41,1410) -- (438.96,1428.08) ;
\draw [shift={(443.33,1430.51)}, rotate = 209.05] [fill={rgb, 255:red, 74; green, 144; blue, 226 }  ,fill opacity=1 ][line width=0.08]  [draw opacity=0] (14.29,-6.86) -- (0,0) -- (14.29,6.86) -- cycle    ;
\draw [color={rgb, 255:red, 208; green, 2; blue, 27 }  ,draw opacity=1 ][line width=2.25]    (503.05,1410) -- (537.4,1467.25) ;
\draw [shift={(539.97,1471.53)}, rotate = 239.04] [fill={rgb, 255:red, 208; green, 2; blue, 27 }  ,fill opacity=1 ][line width=0.08]  [draw opacity=0] (14.29,-6.86) -- (0,0) -- (14.29,6.86) -- cycle    ;
\draw [color={rgb, 255:red, 74; green, 144; blue, 226 }  ,draw opacity=1 ][line width=2.25]    (503.05,1410) -- (470.5,1428.08) ;
\draw [shift={(466.13,1430.51)}, rotate = 330.95] [fill={rgb, 255:red, 74; green, 144; blue, 226 }  ,fill opacity=1 ][line width=0.08]  [draw opacity=0] (14.29,-6.86) -- (0,0) -- (14.29,6.86) -- cycle    ;
\draw    (255.25,1448.61) -- (269.73,1448.61) ;
\draw    (351.78,1448.61) -- (366.26,1448.61) ;
\draw    (448.3,1448.61) -- (462.78,1448.61) ;
\draw    (213.51,1529.87) -- (250.43,1550.38) ;
\draw    (176.59,1550.38) -- (213.51,1529.87) ;
\draw    (250.43,1591.41) -- (250.43,1550.38) ;
\draw [color={rgb, 255:red, 0; green, 0; blue, 0 }  ,draw opacity=1 ][line width=0.75]    (176.59,1591.41) -- (176.59,1550.38) ;
\draw    (176.59,1591.41) -- (213.51,1611.92) ;
\draw    (213.51,1611.92) -- (250.43,1591.41) ;
\draw [color={rgb, 255:red, 0; green, 0; blue, 0 }  ,draw opacity=1 ][line width=0.75]    (250.43,1591.41) -- (213.51,1529.87) ;
\draw [color={rgb, 255:red, 0; green, 0; blue, 0 }  ,draw opacity=1 ][line width=0.75]    (213.51,1611.92) -- (213.51,1529.87) ;
\draw [color={rgb, 255:red, 0; green, 0; blue, 0 }  ,draw opacity=1 ][line width=0.75]    (176.59,1591.41) -- (213.51,1529.87) ;
\draw    (503.05,1529.87) -- (539.97,1550.38) ;
\draw    (466.13,1550.38) -- (503.05,1529.87) ;
\draw    (539.97,1591.41) -- (539.97,1550.38) ;
\draw [color={rgb, 255:red, 0; green, 0; blue, 0 }  ,draw opacity=1 ][line width=0.75]    (466.13,1591.41) -- (466.13,1550.38) ;
\draw    (466.13,1591.41) -- (503.05,1611.92) ;
\draw    (503.05,1611.92) -- (539.97,1591.41) ;
\draw [color={rgb, 255:red, 0; green, 0; blue, 0 }  ,draw opacity=1 ][line width=0.75]    (539.97,1591.41) -- (503.05,1529.87) ;
\draw [color={rgb, 255:red, 0; green, 0; blue, 0 }  ,draw opacity=1 ][line width=0.75]    (503.05,1611.92) -- (503.05,1529.87) ;
\draw [color={rgb, 255:red, 0; green, 0; blue, 0 }  ,draw opacity=1 ][line width=0.75]    (466.13,1591.41) -- (503.05,1529.87) ;
\draw    (311.52,1529.87) -- (348.44,1550.38) ;
\draw    (274.6,1550.38) -- (311.52,1529.87) ;
\draw    (348.44,1591.41) -- (348.44,1550.38) ;
\draw [color={rgb, 255:red, 0; green, 0; blue, 0 }  ,draw opacity=1 ][line width=0.75]    (274.6,1591.41) -- (274.6,1550.38) ;
\draw    (274.6,1591.41) -- (311.52,1611.92) ;
\draw    (311.52,1611.92) -- (348.44,1591.41) ;
\draw [color={rgb, 255:red, 0; green, 0; blue, 0 }  ,draw opacity=1 ][line width=0.75]    (348.44,1591.41) -- (311.52,1529.87) ;
\draw [color={rgb, 255:red, 0; green, 0; blue, 0 }  ,draw opacity=1 ][line width=0.75]    (311.52,1611.92) -- (311.52,1529.87) ;
\draw [color={rgb, 255:red, 0; green, 0; blue, 0 }  ,draw opacity=1 ][line width=0.75]    (274.6,1591.41) -- (311.52,1529.87) ;
\draw    (406.41,1529.87) -- (443.33,1550.38) ;
\draw    (369.49,1550.38) -- (406.41,1529.87) ;
\draw    (443.33,1591.41) -- (443.33,1550.38) ;
\draw [color={rgb, 255:red, 0; green, 0; blue, 0 }  ,draw opacity=1 ][line width=0.75]    (369.49,1591.41) -- (369.49,1550.38) ;
\draw    (369.49,1591.41) -- (406.41,1611.92) ;
\draw    (406.41,1611.92) -- (443.33,1591.41) ;
\draw [color={rgb, 255:red, 0; green, 0; blue, 0 }  ,draw opacity=1 ][line width=0.75]    (443.33,1591.41) -- (406.41,1529.87) ;
\draw [color={rgb, 255:red, 0; green, 0; blue, 0 }  ,draw opacity=1 ][line width=0.75]    (406.41,1611.92) -- (406.41,1529.87) ;
\draw [color={rgb, 255:red, 0; green, 0; blue, 0 }  ,draw opacity=1 ][line width=0.75]    (369.49,1591.41) -- (406.41,1529.87) ;
\draw [color={rgb, 255:red, 208; green, 2; blue, 27 }  ,draw opacity=1 ][line width=2.25]    (213.51,1529.87) -- (176.55,1591.38) ;
\draw [color={rgb, 255:red, 74; green, 144; blue, 226 }  ,draw opacity=1 ][line width=2.25]    (213.51,1529.87) -- (250.43,1550.38) ;
\draw [color={rgb, 255:red, 208; green, 2; blue, 27 }  ,draw opacity=1 ][line width=2.25]    (311.23,1530.78) -- (274.27,1592.29) ;
\draw [color={rgb, 255:red, 74; green, 144; blue, 226 }  ,draw opacity=1 ][line width=2.25]    (311.52,1529.87) -- (274.6,1550.38) ;
\draw [color={rgb, 255:red, 208; green, 2; blue, 27 }  ,draw opacity=1 ][line width=2.25]    (311.23,1530.78) -- (311.52,1611.92) ;
\draw [color={rgb, 255:red, 74; green, 144; blue, 226 }  ,draw opacity=1 ][line width=2.25]    (311.52,1529.87) -- (348.44,1550.38) ;
\draw [color={rgb, 255:red, 208; green, 2; blue, 27 }  ,draw opacity=1 ][line width=2.25]    (406.41,1529.87) -- (443.33,1591.41) ;
\draw [color={rgb, 255:red, 208; green, 2; blue, 27 }  ,draw opacity=1 ][line width=2.25]    (406.41,1529.87) -- (406.7,1611.01) ;
\draw [color={rgb, 255:red, 74; green, 144; blue, 226 }  ,draw opacity=1 ][line width=2.25]    (406.41,1529.87) -- (369.49,1550.38) ;
\draw [color={rgb, 255:red, 74; green, 144; blue, 226 }  ,draw opacity=1 ][line width=2.25]    (406.41,1529.87) -- (443.33,1550.38) ;
\draw [color={rgb, 255:red, 208; green, 2; blue, 27 }  ,draw opacity=1 ][line width=2.25]    (503.05,1529.87) -- (539.97,1591.41) ;
\draw [color={rgb, 255:red, 74; green, 144; blue, 226 }  ,draw opacity=1 ][line width=2.25]    (503.05,1529.87) -- (466.13,1550.38) ;
\draw    (255.25,1568.48) -- (269.73,1568.48) ;
\draw    (351.78,1568.48) -- (366.26,1568.48) ;
\draw    (448.3,1568.48) -- (462.78,1568.48) ;
\draw [color={rgb, 255:red, 74; green, 144; blue, 226 }  ,draw opacity=1 ][line width=2.25]    (176.59,1550.38) -- (176.59,1591.41) ;
\draw [color={rgb, 255:red, 74; green, 144; blue, 226 }  ,draw opacity=1 ][line width=2.25]    (274.6,1591.41) -- (311.52,1611.92) ;
\draw [color={rgb, 255:red, 74; green, 144; blue, 226 }  ,draw opacity=1 ][line width=2.25]    (406.7,1611.01) -- (443.33,1591.41) ;
\draw [color={rgb, 255:red, 74; green, 144; blue, 226 }  ,draw opacity=1 ][line width=2.25]    (539.97,1591.41) -- (539.97,1550.38) ;

\draw (122.06,1445.95) node [anchor=north west][inner sep=0.75pt]    {$\text{SL}_{3}$};
\draw (122.73,1552.97) node [anchor=north west][inner sep=0.75pt]    {$\text{SL}_{2}$};

\end{tikzpicture}

\end {figure}
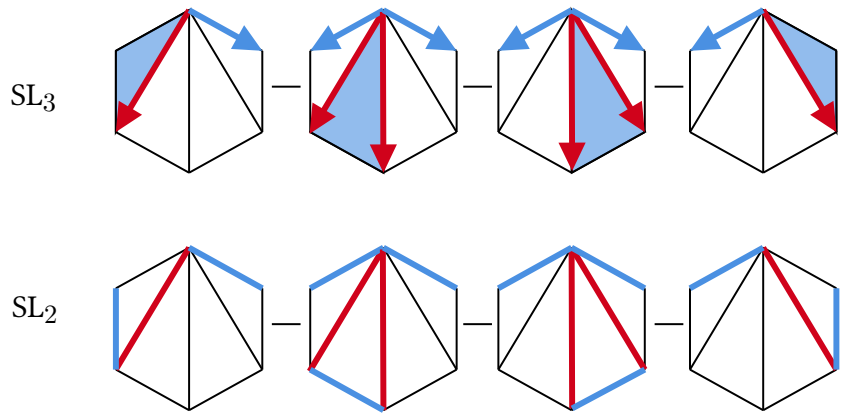

\end{rmk}

Finally, we consider the longest edge in a fan triangulation $\Delta$. Consider the lattice path $L$ in $\mathbb{Z}^{2}$ which starts at the origin and has associated word $(a b^2)^n = a b^2 a b^2 \dots$. Let $L_{\gamma}$ be the set of points in $\mathbb{Z}^{2}$ that are on or below $L$ and weakly above the horizontal axis of $\mathbb{Z}^2$. Make $L_{\gamma}$ into (the Hasse diagram of) a poset by declaring a node $(x_2 , y_2)$ covers another $(x_1 , y_1)$ iff either $x_2 = x_1 + 1$ or $y_2 = y_1 + 1$ (but not both).

\begin{prop} Let $\gamma = i \rightarrow j$ be the longest edge in a fan triangulation, and let $x_{\gamma}$ be the associated cluster variable.  

\begin{enumerate}[(a)]
    \item The Laurent expansion of $x_{\gamma}$ contains $(n+1)^2$ terms. 
    
    \item The minimal element is the $T$-path consisting of $3$ edges; it first goes along the boundary away from the fan vertex, then goes along the diagonal closest to $i$, and then goes along the boundary to $j$. 
    
    \item The maximal element is the $T$-path consisting of $3$ edges; it first goes to the fan vertex, then goes along the diagonal furthest from $i$, and then goes along the boundary edge to $j$. 
        
    \item The poset $\mathbb{T}_{\gamma}$ is isomorphic to $L_{\gamma}$.
\end{enumerate}

\end{prop}

\begin{proof}
    We induct on $n$. Parts (a)-(d) are true for the base case by the resolution rule given in Figure \ref{fig:edge_resolution}, so suppose $n > 1$. Resolve the arc $\gamma$ over the first diagonal it crosses to obtain four diagrams (these diagrams are shown below for the longest edge in a hexagon).

    \begin {figure}[h!]
    \centering
    \caption{The first step in the resolution of the longest edge in a hexagon}
    \label{fig:hex_longest_edge}

\tikzset{every picture/.style={line width=0.75pt}} 

\begin{tikzpicture}[x=0.75pt,y=0.75pt,yscale=-1,xscale=1]

\draw    (130.21,2908.83) -- (190.37,2942.25) ;
\draw    (70.06,2942.25) -- (130.21,2908.83) ;
\draw    (190.37,3009.08) -- (190.37,2942.25) ;
\draw [color={rgb, 255:red, 0; green, 0; blue, 0 }  ,draw opacity=1 ][line width=0.75]    (70.06,3009.08) -- (70.06,2942.25) ;
\draw    (70.06,3009.08) -- (130.21,3042.5) ;
\draw    (130.21,3042.5) -- (190.37,3009.08) ;
\draw [color={rgb, 255:red, 0; green, 0; blue, 0 }  ,draw opacity=1 ][line width=0.75]    (190.37,3009.08) -- (130.21,2908.83) ;
\draw [color={rgb, 255:red, 0; green, 0; blue, 0 }  ,draw opacity=1 ][line width=0.75]    (130.21,3042.5) -- (130.21,2908.83) ;
\draw [color={rgb, 255:red, 0; green, 0; blue, 0 }  ,draw opacity=1 ][line width=0.75]    (70.06,3009.08) -- (130.21,2908.83) ;
\draw    (289.89,2908.83) -- (350.04,2942.25) ;
\draw    (229.74,2942.25) -- (289.89,2908.83) ;
\draw    (350.04,3009.08) -- (350.04,2942.25) ;
\draw [color={rgb, 255:red, 0; green, 0; blue, 0 }  ,draw opacity=1 ][line width=0.75]    (229.74,3009.08) -- (229.74,2942.25) ;
\draw    (229.74,3009.08) -- (289.89,3042.5) ;
\draw    (289.89,3042.5) -- (350.04,3009.08) ;
\draw [color={rgb, 255:red, 0; green, 0; blue, 0 }  ,draw opacity=1 ][line width=0.75]    (350.04,3009.08) -- (289.89,2908.83) ;
\draw [color={rgb, 255:red, 0; green, 0; blue, 0 }  ,draw opacity=1 ][line width=0.75]    (289.89,3042.5) -- (289.89,2908.83) ;
\draw [color={rgb, 255:red, 0; green, 0; blue, 0 }  ,draw opacity=1 ][line width=0.75]    (229.74,3009.08) -- (289.89,2908.83) ;
\draw  [fill={rgb, 255:red, 74; green, 144; blue, 226 }  ,fill opacity=0.6 ] (231.11,3010.56) -- (287.92,2908.07) -- (349.03,2941.95) -- cycle ;
\draw [color={rgb, 255:red, 74; green, 144; blue, 226 }  ,draw opacity=1 ][line width=2.25]    (130.21,2908.83) -- (185.99,2939.82) ;
\draw [shift={(190.37,2942.25)}, rotate = 209.05] [fill={rgb, 255:red, 74; green, 144; blue, 226 }  ,fill opacity=1 ][line width=0.08]  [draw opacity=0] (14.29,-6.86) -- (0,0) -- (14.29,6.86) -- cycle    ;
\draw [color={rgb, 255:red, 74; green, 144; blue, 226 }  ,draw opacity=1 ][line width=2.25]    (70.06,2942.25) -- (70,3004.04) ;
\draw [shift={(70,3009.04)}, rotate = 270.05] [fill={rgb, 255:red, 74; green, 144; blue, 226 }  ,fill opacity=1 ][line width=0.08]  [draw opacity=0] (14.29,-6.86) -- (0,0) -- (14.29,6.86) -- cycle    ;
\draw [color={rgb, 255:red, 208; green, 2; blue, 27 }  ,draw opacity=1 ][line width=2.25]    (130.21,2908.83) -- (72.64,3004.79) ;
\draw [shift={(70.06,3009.08)}, rotate = 300.96] [fill={rgb, 255:red, 208; green, 2; blue, 27 }  ,fill opacity=1 ][line width=0.08]  [draw opacity=0] (14.29,-6.86) -- (0,0) -- (14.29,6.86) -- cycle    ;
\draw  [fill={rgb, 255:red, 208; green, 2; blue, 27 }  ,fill opacity=0.6 ] (229.73,2942.27) -- (289.89,2908.83) -- (229.68,3009.04) -- cycle ;
\draw [color={rgb, 255:red, 74; green, 144; blue, 226 }  ,draw opacity=1 ][line width=2.25]    (289.89,2908.83) -- (234.11,2939.82) ;
\draw [shift={(229.74,2942.25)}, rotate = 330.95] [fill={rgb, 255:red, 74; green, 144; blue, 226 }  ,fill opacity=1 ][line width=0.08]  [draw opacity=0] (14.29,-6.86) -- (0,0) -- (14.29,6.86) -- cycle    ;
\draw [color={rgb, 255:red, 208; green, 2; blue, 27 }  ,draw opacity=1 ][line width=2.25]    (287.92,2909.64) -- (232.81,3005.44) ;
\draw [shift={(230.32,3009.77)}, rotate = 299.90999999999997] [fill={rgb, 255:red, 208; green, 2; blue, 27 }  ,fill opacity=1 ][line width=0.08]  [draw opacity=0] (14.29,-6.86) -- (0,0) -- (14.29,6.86) -- cycle    ;
\draw    (450.17,2909.59) -- (510.32,2943) ;
\draw    (390.02,2943) -- (450.17,2909.59) ;
\draw    (510.32,3009.84) -- (510.32,2943) ;
\draw [color={rgb, 255:red, 0; green, 0; blue, 0 }  ,draw opacity=1 ][line width=0.75]    (390.02,3009.84) -- (390.02,2943) ;
\draw    (390.02,3009.84) -- (450.17,3043.26) ;
\draw    (450.17,3043.26) -- (510.32,3009.84) ;
\draw [color={rgb, 255:red, 0; green, 0; blue, 0 }  ,draw opacity=1 ][line width=0.75]    (510.32,3009.84) -- (450.17,2909.59) ;
\draw [color={rgb, 255:red, 0; green, 0; blue, 0 }  ,draw opacity=1 ][line width=0.75]    (450.17,3043.26) -- (450.17,2909.59) ;
\draw [color={rgb, 255:red, 0; green, 0; blue, 0 }  ,draw opacity=1 ][line width=0.75]    (390.02,3009.84) -- (450.17,2909.59) ;
\draw  [fill={rgb, 255:red, 74; green, 144; blue, 226 }  ,fill opacity=0.6 ] (391.38,3011.32) -- (448.19,2908.83) -- (509.31,2942.7) -- cycle ;
\draw  [fill={rgb, 255:red, 208; green, 2; blue, 27 }  ,fill opacity=0.6 ] (390,2943.03) -- (450.17,2909.59) -- (389.96,3009.8) -- cycle ;
\draw [color={rgb, 255:red, 74; green, 144; blue, 226 }  ,draw opacity=1 ][line width=2.25]    (390.01,2948) -- (389.96,3009.8) ;
\draw [shift={(390.02,2943)}, rotate = 90.05] [fill={rgb, 255:red, 74; green, 144; blue, 226 }  ,fill opacity=1 ][line width=0.08]  [draw opacity=0] (14.29,-6.86) -- (0,0) -- (14.29,6.86) -- cycle    ;
\draw [color={rgb, 255:red, 208; green, 2; blue, 27 }  ,draw opacity=1 ][line width=2.25]    (446.89,2913.13) -- (390.74,3010.59) ;
\draw [shift={(449.38,2908.8)}, rotate = 119.95] [fill={rgb, 255:red, 208; green, 2; blue, 27 }  ,fill opacity=1 ][line width=0.08]  [draw opacity=0] (14.29,-6.86) -- (0,0) -- (14.29,6.86) -- cycle    ;
\draw    (609.85,2908.83) -- (670,2942.25) ;
\draw    (549.7,2942.25) -- (609.85,2908.83) ;
\draw    (670,3009.08) -- (670,2942.25) ;
\draw [color={rgb, 255:red, 0; green, 0; blue, 0 }  ,draw opacity=1 ][line width=0.75]    (549.7,3009.08) -- (549.7,2942.25) ;
\draw    (549.7,3009.08) -- (609.85,3042.5) ;
\draw    (609.85,3042.5) -- (670,3009.08) ;
\draw [color={rgb, 255:red, 0; green, 0; blue, 0 }  ,draw opacity=1 ][line width=0.75]    (670,3009.08) -- (609.85,2908.83) ;
\draw [color={rgb, 255:red, 0; green, 0; blue, 0 }  ,draw opacity=1 ][line width=0.75]    (609.85,3042.5) -- (609.85,2908.83) ;
\draw [color={rgb, 255:red, 0; green, 0; blue, 0 }  ,draw opacity=1 ][line width=0.75]    (549.7,3009.08) -- (609.85,2908.83) ;
\draw [color={rgb, 255:red, 74; green, 144; blue, 226 }  ,draw opacity=1 ][line width=2.25]    (605.48,2911.26) -- (549.7,2942.25) ;
\draw [shift={(609.85,2908.83)}, rotate = 150.95] [fill={rgb, 255:red, 74; green, 144; blue, 226 }  ,fill opacity=1 ][line width=0.08]  [draw opacity=0] (14.29,-6.86) -- (0,0) -- (14.29,6.86) -- cycle    ;
\draw [color={rgb, 255:red, 208; green, 2; blue, 27 }  ,draw opacity=1 ][line width=2.25]    (607.27,2913.11) -- (549.63,3009.04) ;
\draw [shift={(609.85,2908.83)}, rotate = 121] [fill={rgb, 255:red, 208; green, 2; blue, 27 }  ,fill opacity=1 ][line width=0.08]  [draw opacity=0] (14.29,-6.86) -- (0,0) -- (14.29,6.86) -- cycle    ;
\draw [color={rgb, 255:red, 74; green, 144; blue, 226 }  ,draw opacity=1 ][line width=2.25]    (445.43,2911.31) -- (390.02,2943) ;
\draw [shift={(449.77,2908.83)}, rotate = 150.23] [fill={rgb, 255:red, 74; green, 144; blue, 226 }  ,fill opacity=1 ][line width=0.08]  [draw opacity=0] (14.29,-6.86) -- (0,0) -- (14.29,6.86) -- cycle    ;
\draw [color={rgb, 255:red, 74; green, 144; blue, 226 }  ,draw opacity=1 ][line width=2.25]    (229.73,2942.27) -- (231,3005.56) ;
\draw [shift={(231.11,3010.56)}, rotate = 268.84000000000003] [fill={rgb, 255:red, 74; green, 144; blue, 226 }  ,fill opacity=1 ][line width=0.08]  [draw opacity=0] (14.29,-6.86) -- (0,0) -- (14.29,6.86) -- cycle    ;
\draw [color={rgb, 255:red, 74; green, 144; blue, 226 }  ,draw opacity=1 ][line width=2.25]    (665.63,2944.67) -- (549.63,3009.04) ;
\draw [shift={(670,2942.25)}, rotate = 150.97] [fill={rgb, 255:red, 74; green, 144; blue, 226 }  ,fill opacity=1 ][line width=0.08]  [draw opacity=0] (14.29,-6.86) -- (0,0) -- (14.29,6.86) -- cycle    ;

\draw (200.96,2960.84) node [anchor=north west][inner sep=0.75pt]    {$+$};
\draw (358.22,2960.08) node [anchor=north west][inner sep=0.75pt]    {$+$};
\draw (518.62,2960.87) node [anchor=north west][inner sep=0.75pt]    {$+$};

\end{tikzpicture}

\end {figure}
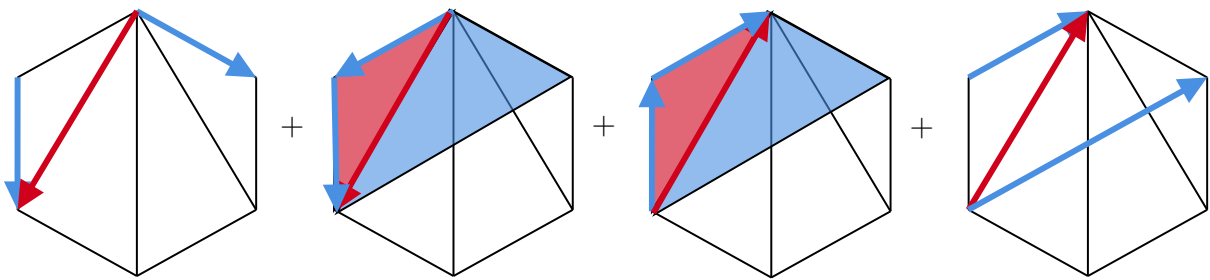

    One of these terms will consist of three directed edges that are contained within the triangulation; this is the minimal element. Call this minimal element $M$. 
    
    Two diagrams in the resolution consist of five terms. Both of these diagrams consist of a cycle around the first triangle, a red face inside the first triangle, and a blue face not contained in the triangulation. The only difference between these two diagrams is whether the initial $3$-cycle around the first triangle goes clockwise or counterclockwise. Call the term with the counterclockwise cycle $F_1$, and call the other $F_2$. 
    
    As we have seen, each of these two diagrams expands into a chain poset with $n$ terms. The terms from each of these posets are pairwise related by the face covering relation in the first triangle. Moreover, one can check that the minimal term in $F_1$ is connected to the minimal element $M$ by the edge covering relation about the first internal diagonal. Gluing together $M$ and the elements from $F_1$ and $F_2$ gives a poset. Call this poset $P_1$. Embed this poset into the discrete plane by sending $M$ to $(0,0),$ and situating $F_1$ and $F_2$ in the first quadrant, along the lines $y=0$ and $y=1$, respectively.  
    
    The remaining diagram from the resolution of $\gamma$ consists of three directed edges, two of which are elements in the triangulation. The remaining third element is the longest arc $\gamma^{\prime}$ in a subfan of $\Delta$ with $n-1$ diagonals. By induction, this poset, which we call $P_2$, is isomorphic to $L_{\gamma^{\prime}}.$ One can check that each of the $n$ elements along the bottom row of $P_2$ is related to the corresponding element in $F_2$ by an edge covering relation about the first diagonal that $\gamma$ crosses. Hence, $(d)$ is true. Part $(a)$ follows immediately from part $(d)$. That $(b)$ and $(c)$ are true follow by induction. This completes the proof. 
\end{proof}

\begin{ex}
    Figure \ref{fig:poset_longest_edge} shows how the expansion poset of the longest edge in a pentagon is constructed via inductive gluing. 
    
    \begin {figure}[h!]
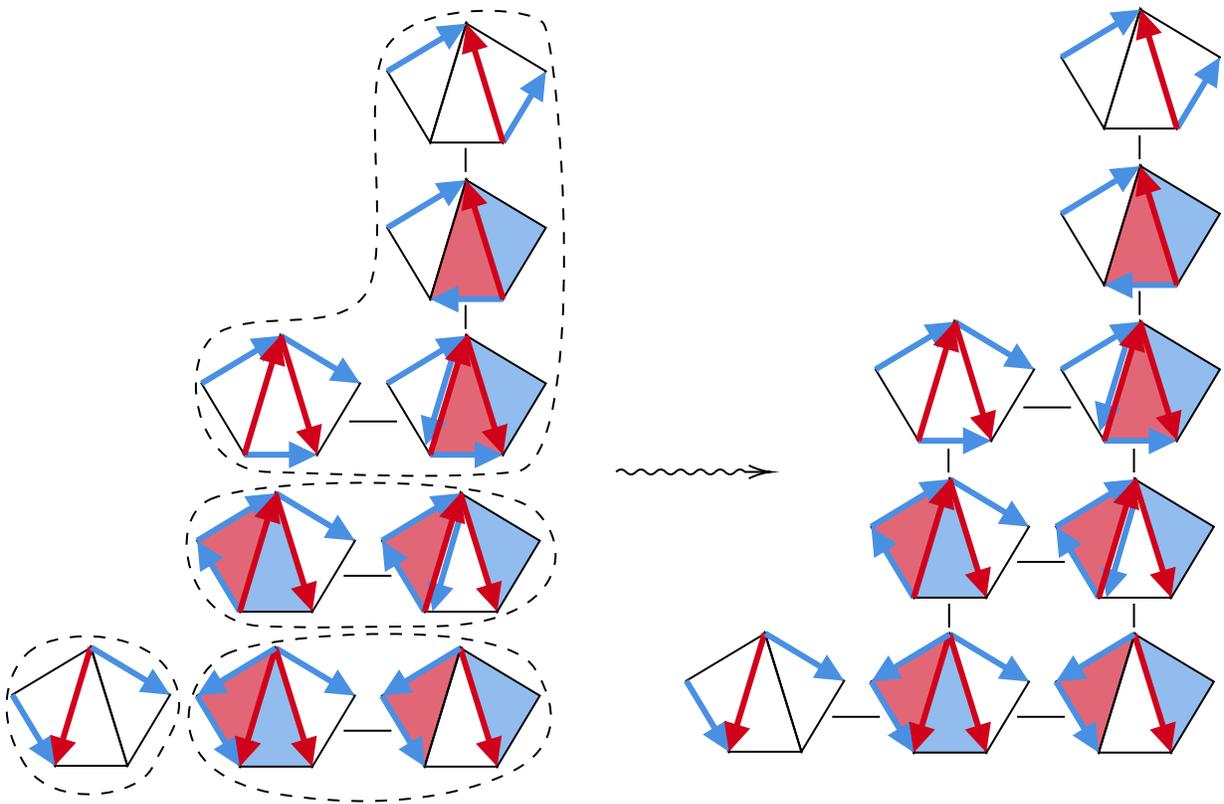

    \centering
    \caption{The gluing construction of the expansion poset of the longest edge in a pentagon}
    \label{fig:poset_longest_edge}

\tikzset{every picture/.style={line width=0.75pt}} 



\end {figure}

\end{ex}

\end {doublespace}

\renewcommand {\bibname} {REFERENCES}

\bibliographystyle {plain}
\bibliography {paper}

\begin{thebibliography}{10}

\bibitem{baileycluster}
Rachel Bailey and Emily Gunawan.
\newblock Cluster algebras and binary words.

\bibitem{birkhoff1937rings}
Garrett Birkhoff et~al.
\newblock Rings of sets.
\newblock {\em Duke Mathematical Journal}, 3(3):443--454, 1937.

\bibitem{branden2015unimodality}
Petter Br{\"a}nd{\'e}n.
\newblock Unimodality, log-concavity, real-rootedness and beyond.
\newblock {\em Handbook of enumerative combinatorics}, 87:437, 2015.

\bibitem{brenti178log}
Francesco Brenti.
\newblock Log-concave and unimodal sequences in algebra, combinatorics, and
  geometry: an update, jerusalem combinatorics’ 93, 71--89.
\newblock {\em Contemp. Math}, 178.

\bibitem{ccanakcci2018cluster}
{\.I}lke {\c{C}}anak{\c{c}}{\i} and Ralf Schiffler.
\newblock Cluster algebras and continued fractions.
\newblock {\em Compositio mathematica}, 154(3):565--593, 2018.

\bibitem{ccanakcci2020snake}
{\.I}lke {\c{C}}anak{\c{c}}{\i} and Ralf Schiffler.
\newblock Snake graphs and continued fractions.
\newblock {\em European Journal of Combinatorics}, 86:103081, 2020.

\bibitem{cheng2016strict}
Zhiyun Cheng, Sujoy Mukherjee, Jozef~H Przytycki, Xiao Wang, and Seung~Yeop
  Yang.
\newblock Strict unimodality of q-polynomials of rooted trees.
\newblock {\em arXiv preprint arXiv:1601.03465}, 2016.

\bibitem{felikson2012cluster}
Anna Felikson, Michael Shapiro, and Pavel Tumarkin.
\newblock Cluster algebras and triangulated orbifolds.
\newblock {\em Advances in Mathematics}, 231(5):2953--3002, 2012.

\bibitem{fock2006moduli}
Vladimir Fock and Alexander Goncharov.
\newblock Moduli spaces of local systems and higher teichm{\"u}ller theory.
\newblock {\em Publications Math{\'e}matiques de l'IH{\'E}S}, 103:1--211, 2006.

\bibitem{fomin2012tensor}
Sergey Fomin and Pavlo Pylyavskyy.
\newblock Tensor diagrams and cluster algebras.
\newblock {\em arXiv preprint arXiv:1210.1888}, 2012.

\bibitem{fomin2014webs}
Sergey Fomin and Pavlo Pylyavskyy.
\newblock Webs on surfaces, rings of invariants, and clusters.
\newblock {\em Proceedings of the National Academy of Sciences},
  111(27):9680--9687, 2014.

\bibitem{fomin2006cluster}
Sergey Fomin, Michael Shapiro, and Dylan Thurston.
\newblock Cluster algebras and triangulated surfaces. part i: Cluster
  complexes.
\newblock {\em arXiv preprint math/0608367}, 2006.

\bibitem{fomin2018cluster}
Sergey Fomin and Dylan Thurston.
\newblock {\em Cluster algebras and triangulated surfaces part II: Lambda
  lengths}, volume 255.
\newblock American Mathematical Society, 2018.

\bibitem{fomin2002cluster}
Sergey Fomin and Andrei Zelevinsky.
\newblock Cluster algebras i: foundations.
\newblock {\em Journal of the American Mathematical Society}, 15(2):497--529,
  2002.

\bibitem{fomin2007cluster}
Sergey Fomin and Andrei Zelevinsky.
\newblock Cluster algebras iv: coefficients.
\newblock {\em Compositio Mathematica}, 143(1):112--164, 2007.

\bibitem{gaebler2004alexander}
Rob Gaebler.
\newblock Alexander polynomials of two-bridge knots and links.
\newblock {\em Undergraduate Thesis}, 2004.

\bibitem{gansner1982lattice}
Emden~R Gansner.
\newblock On the lattice of order ideals of an up-down poset.
\newblock {\em Discrete Mathematics}, 39(2):113--122, 1982.

\bibitem{gekhtman2005cluster}
Michael Gekhtman, Michael Shapiro, Alek Vainshtein, et~al.
\newblock Cluster algebras and weil-petersson forms.
\newblock {\em Duke Mathematical Journal}, 127(2):291--311, 2005.

\bibitem{gekhtman2003cluster}
Michael Gekhtman, Michael~Zalmanovich Shapiro, and Alek~D Vainshtein.
\newblock Cluster algebras and poisson geometry.
\newblock {\em Moscow Mathematical Journal}, 3(3):899--934, 2003.

\bibitem{gunawan2019cluster}
Emily Gunawan, Gregg Musiker, and Hannah Vogel.
\newblock Cluster algebraic interpretation of infinite friezes.
\newblock {\em European Journal of Combinatorics}, 81:22--57, 2019.

\bibitem{hoft1985fibonacci}
H~H{\"o}ft and M~H{\"o}ft.
\newblock A fibonacci sequence of distributive lattices.
\newblock {\em Fibonacci Quart}, 23(3):232--237, 1985.

\bibitem{hsu1993fibonacci}
W-J Hsu.
\newblock Fibonacci cubes-a new interconnection topology.
\newblock {\em IEEE Transactions on Parallel and Distributed Systems},
  4(1):3--12, 1993.

\bibitem{klavvzar2013structure}
Sandi Klav{\v{z}}ar.
\newblock Structure of fibonacci cubes: a survey.
\newblock {\em Journal of Combinatorial Optimization}, 25(4):505--522, 2013.

\bibitem{knauer2018lattice}
Kolja Knauer, Leonardo Mart{\'\i}nez-Sandoval, and Jorge Luis~Ram{\'\i}rez
  Alfons{\'\i}n.
\newblock On lattice path matroid polytopes: integer points and ehrhart
  polynomial.
\newblock {\em Discrete \& Computational Geometry}, 60(3):698--719, 2018.

\bibitem{lee2019cluster}
Kyungyong Lee and Ralf Schiffler.
\newblock Cluster algebras and jones polynomials.
\newblock {\em Selecta Mathematica}, 25(4):58, 2019.

\bibitem{morier2018q}
Sophie Morier-Genoud and Valentin Ovsienko.
\newblock $ q $-deformed rationals and $ q $-continued fractions.
\newblock {\em arXiv preprint arXiv:1812.00170}, 2018.

\bibitem{munarini2002rank}
Emanuele Munarini and Norma~Zagaglia Salvi.
\newblock On the rank polynomial of the lattice of order ideals of fences and
  crowns.
\newblock {\em Discrete mathematics}, 259(1-3):163--177, 2002.

\bibitem{musiker2010cluster}
Gregg Musiker and Ralf Schiffler.
\newblock Cluster expansion formulas and perfect matchings.
\newblock {\em Journal of Algebraic Combinatorics}, 32(2):187--209, 2010.

\bibitem{musiker2011positivity}
Gregg Musiker, Ralf Schiffler, and Lauren Williams.
\newblock Positivity for cluster algebras from surfaces.
\newblock {\em Advances in Mathematics}, 227(6):2241--2308, 2011.

\bibitem{musiker2013bases}
Gregg Musiker, Ralf Schiffler, and Lauren Williams.
\newblock Bases for cluster algebras from surfaces.
\newblock {\em Compositio Mathematica}, 149(2):217--263, 2013.

\bibitem{o1990unimodality}
Kathleen~M O'Hara.
\newblock Unimodality of gaussian coefficients: a constructive proof.
\newblock {\em Journal of Combinatorial Theory, Series A}, 53(1):29--52, 1990.

\bibitem{propp2005combinatorics}
James Propp.
\newblock The combinatorics of frieze patterns and markoff numbers.
\newblock {\em arXiv preprint math/0511633}, 2005.

\bibitem{rabideau2018continued}
Michelle Rabideau and Ralf Schiffler.
\newblock Continued fractions and orderings on the markov numbers.
\newblock {\em arXiv preprint arXiv:1801.07155}, 2018.

\bibitem{sagan2013symmetric}
Bruce~E Sagan.
\newblock {\em The symmetric group: representations, combinatorial algorithms,
  and symmetric functions}, volume 203.
\newblock Springer Science \& Business Media, 2013.

\bibitem{schiffler2008cluster}
Ralf Schiffler.
\newblock A cluster expansion formula ($ a\_n $ case).
\newblock {\em the electronic journal of combinatorics}, 15(1):R64, 2008.

\bibitem{schiffler2010cluster}
Ralf Schiffler.
\newblock On cluster algebras arising from unpunctured surfaces ii.
\newblock {\em Advances in Mathematics}, 223(6):1885--1923, 2010.

\bibitem{schiffler2009cluster}
Ralf Schiffler and Hugh Thomas.
\newblock On cluster algebras arising from unpunctured surfaces.
\newblock {\em International Mathematics Research Notices},
  2009(17):3160--3189, 2009.

\bibitem{stanley1988differential}
Richard~P Stanley.
\newblock Differential posets.
\newblock {\em Journal of the American Mathematical Society}, 1(4):919--961,
  1988.

\bibitem{stanley1989log}
Richard~P Stanley.
\newblock Log-concave and unimodal sequences in algebra, combinatorics, and
  geometry.
\newblock {\em Ann. New York Acad. Sci}, 576(1):500--535, 1989.

\bibitem{uludaug2015jimm}
A~Muhammed Uluda{\u{g}} and Hakan Ayral.
\newblock Jimm, a fundamental involution.
\newblock {\em arXiv preprint arXiv:1501.03787}, 2015.

\bibitem{yurikusa2018combinatorial}
Toshiya Yurikusa.
\newblock Combinatorial cluster expansion formulas from triangulated surfaces.
\newblock {\em arXiv preprint arXiv:1808.01567}, 2018.

\bibitem{yurikusa2019cluster}
Toshiya Yurikusa.
\newblock Cluster expansion formulas in type a.
\newblock {\em Algebras and Representation Theory}, 22(1):1--19, 2019.

\end{thebibliography}

\end{document}